\documentclass[10pt,francais,leqno,makeidx]{smfart}

\usepackage[latin1]{inputenc}
\usepackage[francais]{babel}
\usepackage[T1]{fontenc}

\usepackage{amssymb,url,xspace,smfthm}
\usepackage{mathrsfs,euscript,color}

\usepackage[all]{xy}

\newcommand{\BibTeX}{{\scshape Bib}\kern-.08em\TeX}
\newcommand{\T}{\S\kern .15em\relax }
\newcommand{\AMS}{$\mathcal{A}$\kern-.1667em\lower.5ex\hbox
        {$\mathcal{M}$}\kern-.125em$\mathcal{S}$}

\theoremstyle{plain}

\newtheorem*{montheo}{\textsc{Théorème}}
\newtheorem*{montheo1}{\textsc{Théorème 1}}
\newtheorem*{montheo2}{\textsc{Théorème 2}}
\newtheorem*{mapropo}{\textsc{Proposition}}
\newtheorem*{mapropo1}{\textsc{Proposition 1}}
\newtheorem*{mapropo2}{\textsc{Proposition 2}}

\newtheorem*{monlem}{\textsc{Lemme}}
\newtheorem*{monlem1}{\textsc{Lemme 1}}
\newtheorem*{monlem2}{\textsc{Lemme 2}}
\newtheorem*{monlem3}{\textsc{Lemme 3}}
\newtheorem*{monlem4}{\textsc{Lemme 4}}

\newtheorem*{moncoro}{\textsc{Corollaire}}
\newtheorem*{moncoro1}{\textsc{Corollaire 1}}
\newtheorem*{moncoro2}{\textsc{Corollaire 2}}
\newtheorem*{moncoro3}{\textsc{Corollaire 3}}

\newtheorem*{marema}{\textsc{Remarque}}
\newtheorem*{marema1}{\textsc{Remarque 1}}
\newtheorem*{marema2}{\textsc{Remarque 2}}
\newtheorem*{marema3}{\textsc{Remarque 3}}
\newtheorem*{marema4}{\textsc{Remarque 4}}
\newtheorem*{mesrems}{\textsc{Remarques}}

\newtheorem*{madefi}{\textsc{Définition}}

\newtheorem*{mesdefi}{\textsc{Définitions}}

\newtheorem*{exemple}{\textsc{Exemple}}
\newtheorem*{exemples}{\textsc{Exemples}}

\newtheorem*{convention}{\textsc{Convention}}

\catcode`\@=11
\catcode`\=13
\catcode`\=13
\catcode`\=13
\catcode`\=13
\catcode`\=13
\catcode`\=13
\catcode`\=13
\catcode`\=13
\catcode`\=13
\catcode`\=13
\catcode`\=13
\def {\'e}
\def {\`e}
\def {\`a}
\def {\`u}
\def {\^e}
\def {\^a}
\def {\^o}
\def {\^{\i}}
\def {\^u}
\def {\c c}
\def {\"e}
\catcode`\:=13
\def :{~\string:}
\catcode`\;=13
\def ;{~\string;}
\catcode`\!=13
\def !{~\string!}

\def\cad{c'est--\`a--dire\ }
\def\ni{\noindent}

\def\v#1{\vskip#1mm}
\def\h#1{\hskip#1mm}
\def\ES#1{\EuScript{#1}}

\title[CARACT\`ERES TORDUS DES REPR\'ESENTATIONS ADMISSIBLES]
{CARACT\`ERES TORDUS DES REPR\'ESENTATIONS ADMISSIBLES}

\author{Bertrand Lemaire}
\address{Aix Marseille Universit, CNRS, Centrale Marseille, I2M, UMR 7373\\ 39 rue F. Joliot Curie, 
13453 Marseille, France}
\email{bertrand.lemaire@univ-amu.fr}

\makeindex
\begin{document}
\def\smfbyname{}
\small

\begin{abstract}
Soit $F$ un corps commutatif localement compact 
non archim\'edien de caractristique quelconque, ${\bf G}$ un groupe r\'eductif connexe d\'efini sur $F$, 
$\theta$ un $F$--automorphisme de ${\bf G}$, et $\omega$ un caract\`ere 
de ${\bf G}(F)$. On fixe une mesure de Haar $dg$ sur ${\bf G}(F)$. Si $\pi$ est une repr\'esentation 
complexe lisse irr\'eductible $(\theta,\omega)$--stable de ${\bf G}(F)$, \cad telle que 
$\pi\circ\theta\simeq \pi\otimes \omega$, 
le choix d'un isomorphisme $A$ de $\pi\otimes\omega$ 
sur $\pi\circ\theta$ d\'efinit une distribution $\Theta_\pi^A$ sur ${\bf G}(F)$, appel\'ee 
\og caract\`ere ($A$--)tordu de $\pi$\fg{}: pour toute fonction 
$f$ sur ${\bf G}(F)$, localement constante et \`a support compact, on pose 
$\Theta_\pi^A(f)={\rm trace}(\pi(fdg)\circ A)$. Dans cet article, on \'etudie ces distributions $\Theta_\pi^A$, 
sans hypoth\`ese restrictive sur $F$, ${\bf G}$ ou $\theta$. 
On prouve en particulier que la restriction de $\Theta_\pi^A$ \`a l'ouvert dense de ${\bf G}(F)$ 
form\'e des \'el\'ements $\theta$--quasi--r\'eguliers 
est donn\'ee par une fonction localement constante, et l'on 
d\'ecrit le comportement de cette fonction par rapport \`a l'induction parabolique et \`a la 
restriction de Jacquet. Cela nous am\`ene \`a reprendre la th\'eorie de Steinberg sur les 
automorphismes d'un groupe alg\'ebrique, d'un point de vue rationnel. 
\end{abstract}

\begin{altabstract}
Let $F$ be a non--Archimedean locally compact field (${\rm car}(F)\geq 0$), 
${\bf G}$ be a connected reductive group defined over $F$, $\theta$ be an 
$F$--automorphism 
of ${\bf G}$, and $\omega$ be a character of ${\bf G}(F)$. We fix a Haar measure $dg$ 
on ${\bf G}(F)$. For a smooth irreducible $(\theta,\omega)$--stable complex representation $\pi$ of ${\bf G}(F)$, 
that is such that $\pi\circ \theta\simeq \pi\otimes \omega$, the 
choice of an isomorphism $A$ from $\pi\otimes \omega$ to 
$\pi\circ \theta$ defines a distribution $\Theta_\pi^A$, called 
the \og ($A$--)twisted character of $\pi$\fg{}: for a compactly supported locally constant function $f$ on ${\bf G}(F)$, 
we put $\Theta_\pi^A(f)={\rm trace}(\pi(fdg)\circ A)$. 
In this paper, we study these distributions $\Theta_\pi^A$, 
without any restrictive hypothesis on $F$, ${\bf G}$ or $\theta$. 
We prove in particular that the restriction of $\Theta_\pi^A$ on the open dense subset of ${\bf G}(F)$ 
formed of those elements which are $\theta$--quasi--regular is given by a locally constant function, 
and we describe how this function behaves with 
respect to parabolic induction and Jacquet restriction. 
This leads us to take up again the Steinberg theory of 
automorphisms of an algebraic group, from a rationnal point of view.
\end{altabstract}



\altkeywords{non--Archimedean local field, reductive group, twisted space, quasi--semisimple element, quasi--regular element, 
admissible representation, distribution--character, function 
character, orbital integral, Weyl integration formula}
\maketitle
\tableofcontents

\section{Introduction}

\subsection{}\label{les objets} Soit $F$ un corps commutatif localement compact non archim\'edien, et soit 
${\bf G}$ un groupe alg\'ebrique r\'eductif connexe d\'efini sur $F$. On note $G={\bf G}(F)$ le groupe 
des points $F$--rationnels de ${\bf G}$. On munit $G$ de la topologie d\'efinie par $F$, ce 
qui en fait un groupe 
localement profini. Fixons un $F$--automorphisme 
$\theta$ de ${\bf G}$, un caract\`ere 
$\omega$ de ${\bf G}(F)$, et une mesure de Haar $dg$ sur $G$. Si $\pi$ est un repr\'esentation 
(complexe, lisse) admissible \textit{$(\theta,\omega)$--stable} de $G$, \cad telle que 
$\pi^\theta\simeq \omega\pi$ o\`u l'on a pos\'e $\pi^\theta=\pi\circ\theta$ et 
$\omega\pi=\pi\otimes\omega$, alors le choix d'un isomorphisme $A$ de $\omega\pi$ sur $\pi^\theta$ 
d\'efinit une distribution $\Theta_\pi^A$ sur $G$, appel\'ee \og caract\`ere ($A$)--tordu de $\pi$ \fg{}: 
pour toute fonction $f$ sur $G$, localement constante et \`a support compact, on pose
$$
\Theta_\pi^A(f)={\rm trace}(\pi(f)\circ A),
$$
o\`u $\pi(f)=\pi(fdg)$ d\'esigne l'op\'erateur $\int_Gf(g)\pi(g)dg$ sur l'espace de $\pi$. Puisque 
$\pi$ est admissible, l'op\'erateur $\pi(f)\circ A$ et de rang fini, et la distribution 
$\Theta_\pi^A$ sur $G$ est bien d\'efinie. Pour $x\in G$, elle v\'erifie la relation
$$
\Theta_\pi^A({^{x,\theta\!}f})=\omega(x)^{-1}\Theta_\pi^A(f),
$$
o\`u l'on a pos\'e ${^{x,\theta\!}f}(g)=f(x^{-1}g\theta(x))$, $g\in G$.

\subsection{}\label{survol}
Dans cet article, on tudie les principales 
propri\'et\'es des caract\`eres tordus de $G$ (voir \ref{description} pour une description dtaille), sans hypoth\`ese particuli\`ere sur $F$, ${\bf G}$ ou 
$\theta$. Ces propri\'et\'es sont souvent d\'ej\`a connues, mais ne sont en g\'en\'eral 
d\'emontr\'ees que dans le cas non tordu, \cad pour $(\theta,\omega)=({\rm id},1)$. 
D'autre part, dans le cas tordu, elles sont souvent \'enonc\'ees sous des hypoth\`eses 
restrictives: ${\rm car}(F)=0$; $\omega=1$; $\theta$ d'ordre 
fini; $\theta$ quasi--semisimple; $G={\rm GL}_n(F)$; etc --- citons en particulier l'article de Clozel \cite{Cl1}, 
dans lequel toutes ces hypoth\`eses sont v\'erifi\'ees. On donne ici des 
d\'emonstrations compl\`etes dans le cas tordu, valables pour $F$, ${\bf G}$ et $\theta$ 
quelconques.

Signalons qu'une version lgrement diffrente de cet article est prsente sur arXiv depuis le 21 juillet 2010 
(arXiv: 1007.3576v1 [math.RT]). Elle a dj t utilse par W.--W. Li \cite{Li} et par J.--L. Waldspurger \cite{W}.

\subsection{}\label{endoscopie} Ces deux types de torsion --- l'une donn\'ee par l'automorphisme $\theta$ et l'autre par le 
caract\`ere $\omega$ --- apparaissent dans de nombreuses applications du principe de fonctorialit\'e de 
Langlands; dj dans les premiers exemples, comme le changement de base cyclique ou la restriction des 
repr\'esentations \`a ${\bf G}_{\rm der}(F)$, on est amens  considrer des caractres tordus. La th\'eorie de l'endoscopie tordue 
\'etudie pr\'ecis\'ement les distributions $D$ sur $G$ v\'erifiant
$$
D({^{x,\theta\!}f})=\omega(x)^{-1}D(f)\leqno{(*)}
$$
pour toute fonction $f$ sur $G$, localement constante et  support compact, et tout $x\in G$. 
Cette th\'eorie permet en particulier de comparer les repr\'esentations irr\'eductibles de $G$ avec celles 
d'autres groupes $H$ qui lui sont associ\'es --- ceux faisant partie d'une donn\'ee endoscopique 
$(H,\ES{H},s,\xi)$ de $(G,\theta,\omega)$, cf. \cite{KS}. Dans les cas particuliers o\`u elles ont 
\'et\'e \'etablies (e.g. le changement de base et l'induction automorphe pour ${\rm GL}_n$), ces correspondances entre 
repr\'esentations s'expriment par des identit\'es de caract\`eres reliant les 
caract\`eres tordus de $G$ \`a des sommes pond\'er\'ees de caract\`eres de $H$. 
Ces identit\'es de caract\`eres sont le plus souvent obtenues dualement \`a partir d'identit\'es 
reliant, c\^ot\'e g\'eom\'etrique, les int\'egrales orbitales endoscopiques de $G$ --- qui sont des sommes 
d'int\'egrales orbitales $(\theta,\omega)$-tordues (voir plus loin, \ref{IO tordues}) pond\'er\'ees 
par des facteurs de transfert --- \`a des int\'egrales orbitales stables de $H$: c'est la \textit{conjecture de transfert}.

\subsection{}\label{analyse harmonique}
L'analyse harmonique en caract\'eristique non nulle est souvent plus compliqu\'ee qu'en caract\'eristique 
nulle, raison pour laquelle certains rsultats d'Harish--Chandra --- comme par exemple l'{\it intgrabilit 
locale des caractres} \cite{HC2} --- sont encore aujourd'hui valables en toute gnralit seulement 
en caractristique nulle (m\^eme dans le cas non tordu). La plupart des difficult\'es nouvelles sont li\'ees \`a 
des questions d'ins\'eparabilit\'e qui ne 
se posent pas en caract\'eristique nulle, ou en caract\'eristique r\'esiduelle \og suffisamment grande\fg. Certains estiment qu'il est aujourd'hui trop t\^ot pour aborder 
ces questions. On pense en revanche qu'il est important, lorsque c'est possible, de disposer 
d'\'enonc\'es vrais en toute caract\'eristique. Par ailleurs, depuis la d\'emonstration r\'ecente du lemme 
fondamental par Ng\^o et Waldspurger, l'analyse harmonique en caract\'eristique non nulle semble conna\^{\i}tre un 
regain d'int\'er\^et.

Le caractre tordu $\Theta_\pi^A$ est une fonction localement constante sur un certain ouvert dense 
de $G$; c'est l'une des principales proprits tablies dans cet article. 
En caractristique nulle, Clozel \cite{Cl2} a dmontr, pour 
$\theta$ d'ordre fini et $\omega=1$, que cette fonction est localement intgrable sur $G$, gnralisant ainsi le rsultat d'Harish--Chandra 
dj cit \cite{HC2}. 
Si $F$ est de caract\'eristique non nulle, cela n'est connu que dans quelques cas particuliers \cite{Le1, Le2}. 
Si la caractristique de $F$ est suffisamment grande (par rapport au rang du groupe),  
il est vraisemblable qu'avec des modifications mineures la m\'ethode d'Harish--Chandra \cite{HC2, Cl2} 
s'applique encore; cela m\'eriterait d'ailleurs d'\^etre r\'edig\'e.

Le r\'esultat dual, qui exprime les int\'egrales orbitales tordues en termes des caract\`eres --- et des variantes de 
ceux-ci, les {\it caract\`eres pond\'er\'es} ---, est une cons\'equence de la {\it formule des traces locale}, \'etablie par 
Arthur \cite{A} pour $F$ de caract\'eristique nulle, $\theta={\rm id}$ et $\omega=1$. Pour $F$ de caractristique nulle, 
la {\it formule des traces locale tordue} a rcemment t tablie par Waldspurger \cite{W}. Pour $F$ de caract\'eristique non nulle, 
l'affaire est loin d'\^etre r\'egl\'ee!

\subsection{}\label{espace tordu (intro)}
Pour tudier les distributions sur $G$ vrifiant la relation $(*)$ de \ref{endoscopie}, 
il est assez commode d'introduire l'espace tordu ${\bf G}^\natural={\bf G}\theta$ de Labesse. 
C'est une vari\'et\'e 
alg\'ebrique affine d\'efinie sur $F$, munie d'un isomorphisme de vari\'et\'es alg\'ebriques 
${\bf G}\rightarrow {\bf G}^\natural,\,g\mapsto g\theta$ lui aussi d\'efini sur $F$, et d'actions 
alg\'ebriques de ${\bf G}$ \`a gauche et \`a droite, commutant entre elles et 
v\'erifiant l'\'egalit\'e $g\cdot g'\theta\cdot g''= gg'\theta(g'')\theta$ pour tous $g,\,g'\!,\,g''\in {\bf G}$. 
Pour $g\in {\bf G}$, on note ${\rm Int}_{{\bf G}^\natural}(g\theta)$ l'automorphisme (alg\'ebrique) 
${\rm Int}_{\bf G}(g)\circ \theta$ de ${\bf G}$. On note $G^\natural=G\theta$ l'ensemble des points 
$F$--rationnels de ${\bf G}^\natural$, que l'on peut voir comme un $G$--espace (topologique) 
tordu. Pour $\gamma\in G^\natural$, l'automorphisme ${\rm Int}_{{\bf G}^\natural}(\gamma)$ de 
${\bf G}$ est d\'efini sur $F$, et induit un automorphisme du groupe topologique $G$, que 
l'on note ${\rm Int}_{G^\natural}(\gamma)$. Pour $g\in G$ et $\gamma\in G^\natural$, posant 
$\tau={\rm Int}_{G^\natural}(\gamma)$, on a
$$
g^{-1}\cdot\gamma\cdot g=g^{-1}\tau(g)\cdot\gamma.
$$

En fait, la donn\'ee de l'automorphisme $\theta$ n'est pas vraiment indispensable: il 
correspond au choix d'un point--base $\delta_1$ de ${\bf G}^\natural(F)$ tel que 
${\rm Int}_{{\bf G}^\natural}(\delta_1)=\theta$. La th\'eorie des groupes alg\'ebriques s'\'etend 
naturellement \`a celle des ${\bf G}$--espaces 
tordus. Par exemple, on appelle \textit{sous--espace parabolique} de $G^\natural$ un sous--espace 
topologique $P^\natural$ de $G^\natural$ de la forme $P^\natural=P\cdot\gamma$ pour 
un sous--groupe parabolique $P$ de $G$ et un \'el\'ement $\gamma$ de $G^\natural$ tels 
que ${\rm Int}_{G^\natural}(\gamma)(P)=P$.

Notons que les distributions $D$ sur $G$ vrifiant la condition $(*)$ de \ref{endoscopie} correspondent aux 
distributions $\ES{D}$ sur $G^\natural$ vrifiant
$$
\ES{D}({^x\phi})=\omega(x)^{-1}\ES{D}(\phi)\leqno{(**)}
$$
pour tout fonction $\phi$ sur $G^\natural$, localement constante et  support compact, et tout $x\in G$; o l'on a pos 
${^x\phi}(\delta)=\phi(x^{-1}\cdot\delta \cdot x)$, $\delta\in G^\natural$.

\subsection{}\label{caractres tordus (intro)}
On appelle \textit{$\omega$--repr\'esentation lisse de $G^\natural$} la donn\'ee 
d'une repr\'esentation lisse $(\pi,V)$ de $G$ et d'une application $\Pi:G^\natural\rightarrow {\rm Aut}_\Bbb{C}(V)$ 
v\'erifiant
$$
\Pi(x\cdot \gamma \cdot y)=\omega(y)\pi(x)\circ\Pi(\gamma)\circ \pi(y)
$$
pour tout $\gamma\in G^\natural$ et tous $x,\,y\in G$. 
La repr\'esentation $\pi$ \'etant d\'etermin\'ee par l'application $\Pi$, on la note aussi $\Pi^\circ$. 
La donn\'ee d'une $\omega$--repr\'esentation lisse $\Pi$ de $G^\natural$ 
\'equivaut \`a celle d'une repr\'esentation lisse $(\theta,\omega)$--stable $\pi$ de $G$ 
munie d'un isomorphisme $A$ de $\omega\pi$ sur $\pi^\theta$: on pose 
$\pi=\Pi^\circ$ et $A=\Pi(\delta_1)$. Si de plus $\pi$ est admissible, auquel cas on 
dit que $\Pi$ est admissible, alors le choix d'une mesure de Haar $d\gamma$ sur 
$G^\natural$ (cf. \ref{module d'un espace tordu}) d\'efinit 
une distribution $\Theta_\Pi$ sur $G^\natural$, appel\'ee \textit{caract\`ere de $\Pi$}: pour 
toute fonction $\phi$ sur $G^\natural$, localement constante et \`a support compact, on pose
$$
\Theta_\Pi(\phi)={\rm trace}(\Pi(\phi)),
$$
o\`u $\Pi(\phi)$ d\'esigne l'op\'erateur $\int_{G^\natural}\phi(\gamma)\Pi(\gamma)d\gamma$ sur 
l'espace de $\Pi$. Si $d\gamma$ est l'image de la mesure $dg$ par l'hom\'eomorphisme 
$G\rightarrow G^\natural,\,g\mapsto g\cdot \delta_1$, alors pour toute fonction $f$ sur $G$, 
localement constante et \`a support compact, notant $f^\natural$ la fonction $g\cdot \delta_1
\mapsto f(g)$ sur $G^\natural$, on a $\Pi(f^\natural)=\pi(f)\circ A$, d'o\`u
$$
\Theta_\Pi(f^\natural)=\Theta_\pi^A(f).
$$

\subsection{}\label{IO tordues}
Pour tout sous--groupe ferm\'e $H$ de $G$, on note $\mathfrak{h}$ (m\^eme lettre gothique) 
son alg\`ebre de Lie. Un \'el\'ement $\gamma$ de $G^\natural$ est dit \textit{quasi--r\'egulier} si pour tout sous--groupe 
parabolique $P$ de $G$, 
on a l'\'egalit\'e
$$\mathfrak{g}(1-\gamma)+\mathfrak{p}=\mathfrak{g},$$
o\`u l'on a pos\'e 
$\mathfrak{g}(1-\gamma)=\{X-{\rm Ad}_{G^\natural}(\gamma) (X):X\in \mathfrak{g}\}$, 
${\rm Ad}_{G^\natural}(\gamma)={\rm Lie}({\rm Int}_{G^\natural}(\gamma))$. Les \'el\'ements quasi--r\'eguliers forment 
un ensemble ouvert dense, disons $G_{\rm qr}^\natural$, dans $G^\natural$. La notion d'\'el\'ement 
quasi--r\'egulier g\'en\'eralise celle, plus classique, d'\'el\'ement \textit{(semisimple) r\'egulier}: les 
\'el\'ements r\'eguliers de $G^\natural$ sont par d\'efinition ceux qui n'annulent pas une certaine 
fonction r\'eguli\`ere $D_{{\bf G}^\natural}$ sur ${\bf G}^\natural$, d\'efinie comme le discriminant 
d'Harish--Chandra $D_{\bf G}$ sur ${\bf G}$ 
(rappelons que pour $g\in G$, on a $D_{\bf G}(g)\neq 0$ si et seulement si 
$g$ est semisimple r\'egulier). 
Les \'el\'ements r\'eguliers forment un ensemble ouvert dense, disons 
$G_{\rm reg}^\natural$, dans $G^\natural$, qui est contenu dans $G_{\rm qr}^\natural$. 
Pour $\gamma\in G_{\rm reg}^\natural$, la $G$--orbite 
$\ES{O}_{G}(\gamma)=\{g^{-1}\cdot\gamma\cdot g:g\in G\}$ est ferm\'ee dans $G^\natural$, et 
la composante connexe ${\bf G}_{\gamma}^\circ$ du centralisateur 
${\bf G}_{\gamma}=\{g\in {\bf G}:g^{-1}\cdot \gamma \cdot g=\gamma\}$ de $\gamma$ dans 
${\bf G}$ est un tore d\'efini sur $F$. 
Le choix d'une mesure de Haar $dg_\gamma$ sur 
le groupe des points $F$--rationnels $G_\gamma^\circ={\bf G}_\gamma^\circ (F)$ de ce tore, 
d\'efinit une distribution $\Lambda^G_{\omega}(\cdot,\gamma)$ sur $G^\natural$ --- appel\'ee 
\textit{$\omega$--int\'egrale orbitale}, ou \textit{int\'egrale orbitale $\omega$--tordue} de $\gamma$: 
pour toute fonction $\phi$ sur $G^\natural$, localement constante et \`a support compact, on pose
$$
\Lambda^G_{\omega}(\phi,\gamma)=\int_{G_\gamma^\circ\backslash G}\omega(g)\phi (g^{-1}\cdot \gamma \cdot g){dg\over dg_\gamma}.
$$
Puisque la $G$--orbite $\ES{O}_G(\gamma)$ est ferm\'ee dans $G^\natural$, l'int\'egrale est absolument convergente, 
et la distribution $\Lambda^G_{\omega}(\cdot,\gamma)$ est bien d\'efinie.

\subsection{}Les caractres tordus $\Theta_\Pi$ et les intgrales orbitales tordues $\Lambda^G_\omega(\cdot,\gamma)$ 
dfinis dans les numros prcdents sont les deux 
familles principales de distributions $\ES{D}$ sur $G^\natural$ vrifiant la condition $(**)$ de \ref{espace tordu (intro)}. On prouve dans cet 
article que pour toute $\omega$--reprsentation admissible $\Pi$ de $G^\natural$ telle que la reprsentation sous--jacente $\Pi^\circ$ de 
$G$ est de type fini (\cad de longueur finie, puisqu'elle est admissible), la restriction de $\Theta_\Pi$  $G_{\rm qr}$ est donne par une fonction localement constante, que l'on note encore $\Theta_\Pi$. En d'autres termes, pour toute fonction $\phi$ sur $G^\natural$, localement constante et  support compact {\it contenu dans $G_{\rm qr}$}, on a 
l'galit
$$
\Theta_\Pi(\phi)= \int_G\phi(\delta)\Theta_\Pi(\delta)d\delta;
$$
l'intgrale est absolument convergente (d'ailleurs c'est m\^eme une somme finie). La formule d'intgration de H.~Weyl permet alors 
de dvelopper l'intgrale ci--dessus en termes des intgrales orbitales tordues $\Lambda^G_\omega(\phi,\gamma)$, 
$\gamma\in G^\natural_{\rm reg}$. On subodore bien s\^ur que l'galit ci--dessus reste vraie m\^eme si le support de $\phi$ n'est pas 
contenu dans $G_{\rm qr}$ --- c'est la proprit d'intgrabilit locale des caractres, cf. \ref{analyse harmonique} --- mais cette 
question n'est pas aborde dans le prsent article.  

\subsection{}On l'a dit plus haut, le groupe $G={\bf G}(F)$ muni de la topologie dfinie par $F$, 
est localement profini. Mais comme groupe des points $F$--rationnels de ${\bf G}$, il h\'erite aussi de la structure 
du groupe alg\'ebrique ${\bf G}$. La division de l'article en six chapitres (voir plus haut la table des matires) est en grande 
partie commande par ce double point de vue: 
groupe topologique -- groupe alg\'ebrique. 

Les propri\'et\'es valables pour n'importe quel groupe topologique localement 
profini $G$ et n'importe quel automorphisme (topologique) $\theta$ de $G$, sont regroup\'ees 
dans le ch.~1. Dans les ch.~2 et 3, qui sont ind\'ependants du reste de l'article, 
on reprend\footnote{Nous ne pensions 
pas au d\'epart devoir reprendre en d\'etail cette th\'eorie, mais cela nous a vite sembl\'e indispensable, 
pour que nos r\'esultats soient valables en toute caract\'eristique (l'analyse harmonique en caract\'eristique 
non nulle est encore un terrain min\'e!).} la th\'eorie de Steinberg sur les automorphismes 
quasi--semisimples de ${\bf G}$, d'un point de vue rationnel. Dans les 
ch.~4, 5 et 6, on applique les r\'esultats \'etablis dans les deux premiers chapitres \`a 
l'\'etude des caract\`eres tordus de 
$G={\bf G}(F)$. 

L'article se conclut par trois annexes (cf. \ref{annexes} pour une description de leur contenu). 
La troisime (C) est  l'origine du laps de temps sparant ce texte de la version parue sur ArXiv en juillet 2010 
(cf. \ref{survol}), une assertion dans ladite version n'tant  
pas dmontre: la $G$--orbite $\ES{O}_G(\delta)$ d'un lment quasi--semisimple $\delta$ 
de $G$ est ferme dans $G$ pour la topologie dfinie par $F$. Pour dmontrer ce rsultat, il faut comprendre certains 
phnomnes d'insparabilit\footnote{La ${\bf G}$--orbite $\ES{O}_{\bf G}(\delta)$ de $\delta$ est ferme dans ${\bf G}^\natural$ pour la topologie de 
Zariski, mais l'application ${\bf G}\rightarrow \ES{O}_{\bf G}(\delta),\, g\mapsto g^{-1}\cdot \delta \cdot g$ peut ne pas \^etre sparable. En revanche si $s$ est un lment 
semisimple de ${\bf G}$, l'application ${\bf G}\rightarrow \ES{O}_{\bf G}(s),\, g\mapsto g^{-1}s g$ est toujours sparable.} qui n'existent pas dans le cas non tordu 
(ou si $F$ est de caractristique nulle). Nous avons pour cela repris les techniques de Bernstein--Zelevinski traitant des morphismes insparables 
\cite[Appendix]{BZ}. Entre--temps, nous nous sommes aper\c{c}us que Moret--Bailly avait r\'ecemment d\'emontr\'e dans un cadre plus 
g\'en\'eral\footnote{Signalons aussi le travail encore plus r\'ecent de Gabber, Gille et Moret--Bailly \cite{GGMB}.} le 
r\'esultat qui nous manquait \cite[theorem 1.3]{MB2} --- cf. la dmonstration de la proposition 1 de 
\ref{la topo p-adique}. Comme l'approche de Bernstein--Zelevinski conduit dans le cas trs particulier 
qui nous intresse  un rsultat explicite qui pourra nous servir dans un travail ultrieur, nous avons dcid de la conserver, tout en utilisant le rsultat de 
Moret--Bailly.

\subsection{}\label{description}D\'ecrivons bri\`evement le contenu des chapitres.

Dans le ch.~1, on d\'efinit les caract\`eres $(\theta,\omega)$--tordus 
d'un groupe localement profini quelconque $G$, o\`u $\theta$ et $\omega$ sont respectivement 
un automorphisme et un caract\`ere de $G$. On introduit 
l'espace (topologique) tordu $G^\natural=G\theta$, 
qui permet de \og voir\fg{}{} les caract\`eres tordus de $G$ comme des caract\`eres de 
$G^\natural$. Enfin on \'etablit une formule de descente 
pour les caract\`eres des $\omega$--repr\'esentations admissibles de $G^\natural$ 
qui sont induites \`a partir d'un sous--espace tordu 
$H^\natural$ de $G^\natural$ tel que $G$ soit compact modulo $H$. Il s'agit d'un 
analogue tordu de la formule bien connue de Van Dijk \cite{VD}, pour les caract\`eres--distributions.

Dans le ch.~2, on rappelle, pour un groupe r\'eductif connexe ${\bf H}$ et 
un $\overline{F}$--automorphisme \textit{quasi--semisimple} 
$\tau$ de ${\bf H}$, les principaux r\'esultats de la th\'eorie de Steinberg \cite{St, DM1}: description 
du groupe des points fixes ${\bf H}_\tau$, du groupe quotient ${\bf H}_\tau/{\bf H}_\tau^\circ$, etc. 
Pour $\tau$ quelconque, 
on introduit la notion d'\textit{\'el\'ement r\'egulier} de l'espace tordu ${\bf H}^\natural={\bf H}\tau$, 
qui g\'en\'eralise celle d'\'el\'ement semisimple r\'egulier 
de ${\bf H}$ (pour $h\in {\bf H}$, l'automorphisme int\'erieur ${\rm Int}_{{\bf H}}(h)$ 
est r\'egulier si et seulement si $h$ est semisimple r\'egulier). Les \'el\'ements r\'eguliers de ${\bf H}^\natural$ forment un 
ensemble ouvert dense dans ${\bf H}^\natural$, que l'on note ${\bf H}^\natural_{\rm reg}$. 
On montre qu'un \'el\'ement $\delta$ de ${\bf H}^\natural$ est r\'egulier si et seulement si les 
deux propri\'et\'es suivantes sont v\'erifi\'ees (thorme de \ref{automorphismes rguliers; le cas g}):
\begin{itemize}
\item le $\overline{F}$-automorphisme $\tau'={\rm Int}_{{\bf H}^\natural}(\delta)$ est quasi--semisimple,
\item le centralisateur connexe 
${\bf H}_\delta^\circ ={\bf H}_{\tau'}^\circ$ de $\delta$ dans ${\bf H}$ est un tore.
\end{itemize}
En ce cas, en posant ${\bf S}={\bf H}_\delta^\circ$, le centralisateur 
${\bf T}=Z_{{\bf H}}({\bf S})$ de ${\bf S}$ dans ${\bf H}$ est un tore maximal de ${\bf H}$ (cf. le lemme 2 de \ref{automorphismes rguliers; le cas g}). De plus, les ensembles 
${\bf S}^\natural= {\bf S}\cdot \delta$ et ${\bf T}^\natural={\bf T}\cdot \delta$ sont des sous--espaces 
tordus de ${\bf H}^\natural$: ${\bf S}^\natural$ est un espace tordu trivial, qu'on appelle 
\textit{tore maximal de ${\bf G}^\natural$}, et ${\bf T}^\natural$ est appel\'e \textit{sous--espace de 
Cartan de ${\bf G}^\natural$}. Notons que ${\bf S}^\natural$ d\'etermine le \textit{quadruplet de Cartan 
$({\bf S},{\bf S}^\natural,{\bf T},{\bf T}^\natural)$}, et que ${\bf T}^\natural$ d\'etermine 
\textit{le triplet de Cartan $({\bf S},{\bf T},{\bf T}^\natural)$}. Enfin pour tout $\delta'\in {\bf T}^\natural\cap {\bf G}^\natural_{\rm reg}$, 
on a ${\bf G}^\circ_{\delta'}={\bf S}$. 

Dans le ch.~3, on s'int\'eresse aux questions de rationalit\'e issues du ch.~2. On montre en particulier que si $\tau$ est un $F$--automorphisme 
quasi--semisimple de ${\bf H}$, alors le groupe ${\bf H}_\tau^\circ$ est d\'efini sur $F$, il existe un 
tore maximal $\tau$--stable ${\bf T}$ de ${\bf H}$ d\'efini sur $F$, et un sous--groupe 
de Borel $\tau$--stable de ${\bf H}$ contenant ${\bf T}$ et d\'efini sur une 
extension s\'eparable finie de $F$ dans $\overline{F}$.

Dans le ch.~4 (et dans les suivants), $G$ est le groupe des points 
$F$--rationnels d'un groupe r\'eductif connexe ${\bf G}$ d\'efini sur $F$, $\theta$ est un 
$F$--automorphisme de ${\bf G}$, et $\omega$ est un caract\`ere de $G$. On introduit 
la notion d'\'el\'ement \textit{quasi--r\'egulier} de l'espace tordu $G^\natural= 
G\theta$. L'ensemble des 
\'el\'ements quasi--r\'eguliers de $G^\natural$, que l'on note 
$G^\natural_{\rm qr}$, est ouvert dense dans $G^\natural$, et contient $G^\natural_{\rm reg}=G\cap {\bf G}^\natural_{\rm reg}$. 
On montre que les caract\`eres 
$\Theta_\Pi$ des $\omega$--repr\'esentations admissibles $\Pi$ 
de $G$ telles que $\Pi^\circ$ est de type fini, sont des fonctions localement constantes sur $G^\natural_{\rm qr}$. 
Ensuite on d\'ecrit comment le foncteur restriction de Jacquet op\`ere sur ces fonctions 
caract\`eres. La description de l'action du foncteur induction parabolique sur les caract\`eres--distributions 
est une simple application de la formule d\'emontr\'ee au ch.~1; son action 
sur les fonctions caract\`eres est d\'ecrite au ch.~6.

Dans le ch.~5, on explicite la fonction caract\`ere d'une $\omega$--repr\'esentation admissible 
$\Pi$ de $G^\natural$ telle que $\Pi^\circ$ est irr\'eductible et essentiellement de carr\'e int\'egrable modulo le 
centre de $G$, gr\^ace aux calculs effectu\'es dans le ch.~4. 
On applique ensuite ce r\'esultat au cas o\`u $\Pi^\circ$ est une 
repr\'esentation irr\'eductible cuspidale de $G$, induite compacte \`a partir d'un 
sous--groupe ouvert, compact modulo le centre, de $G$.

Dans le ch.~6, on introduit une famille de distributions \og duale\fg{} de 
celle des caract\`eres $(\theta,\omega)$--tordus: les int\'egrales orbitales $(\theta,\omega)$--tordues. Comme 
pour les caract\`eres, il est commode de les voir comme des distributions sur 
$G^\natural$: les \textit{$\omega$--int\'egrales 
orbitales}. On \'ecrit la formule 
de descente parabolique pour les $\omega$--int\'egrales orbitales, 
ainsi que la formule d'int\'egration de Weyl pour les fonctions int\'egrables sur 
$G^\natural$. On en d\'eduit la description de l'action du foncteur induction 
parabolique sur les fonctions caract\`eres. 

\subsection{}\label{annexes}
Dcrivons maintenant le contenu des annexes.

Dans l'annexe A, on caract\'erise les $\omega$--repr\'esentations lisses irr\'eductibles de $G^\natural$ en termes 
des foncteurs $V\mapsto V^K$, comme dans \cite{BZ}. On 
prouve aussi l'ind\'ependance lin\'eaires des caract\`eres--distributions $\Theta_\Pi$ pour 
les $\omega$--repr\'esentations admissibles $\Pi$ de $G^\natural$ telles que 
$\Pi^\circ$ est irr\'eductible. 

Dans l'annexe B, qui nous a t suggre par Guy Henniart, on remplace le corps des coefficients $\Bbb{C}$ par un corps $R$ de caract\'eristique $l$ diff\'erente 
de la caract\'eristique r\'esiduelle de $F$. 
On passe en revue les r\'esultats des ch.~1 et 4 dans ce nouveau cadre. On montre en particulier que, tout 
comme pour les repr\'esentations complexes, le caract\`ere--distribution d'une $(\omega,R)$--repr\'esentation 
admissible $\Pi:G^\natural\rightarrow {\rm Aut}_R(V)$ telle que $\Pi^\circ$ est de type fini, 
est donn\'e sur $G^\natural_{\rm qr}$ par une fonction localement 
constante.

Dans l'annexe C, on s'intresse aux proprits topologiques --- pour la topologie dfinie par $F$ --- de 
l'application $\alpha_F:{\bf Y}(F)\rightarrow {\bf X}(F)$ dduite par passage aux points $F$--rationnels 
d'un $F$--morphisme $\alpha:{\bf Y}\rightarrow {\bf X}$ de varits algbriques dfinies sur $F$. On sait dj d'aprs 
Bernstein--Zelevinski \cite[theorem A.2]{BZ} que l'image $\alpha_F({\bf Y}(F))$ est constructible dans ${\bf X}(F)$. 
Rcemment, des rsultats plus fins ont t obtenus par Moret--Bailly dans un cadre plus gnral \cite{MB1, MB2}. 
On sait en particulier d'aprs loc.~cit. que si le morphisme $\alpha$ est fini (resp. tale), alors l'application $\alpha_F$ 
est ferme (resp. un homomorphisme local). On prcise ce rsultat dans le cas particulier suivant (proposition de \ref{une consquence du lemme cl}): 
{\it soit ${\bf H}$ un groupe algbrique affine dfini sur $F$, et soit $\alpha: {\bf Y}\rightarrow {\bf X}$ un $F$--morphisme fini de varits algbriques 
affines irrductibles dfinies sur $F$. On suppose que les vari\'et\'es ${\bf Y}$ et ${\bf X}$ sont munies d'une action alg\'ebrique de ${\bf H}$ dfinie sur 
$F$, que $\alpha$ est ${\bf H}$--quivariant et que ${\bf Y}$ est un ${\bf H}$--espace homogne. 
Il existe un $F$--morphisme $\alpha_1: {\bf Y}_1\rightarrow {\bf X}_1$ de varits algbriques affines dfinies sur $F$, et des 
$F$--morphismes de varits algbriques $\gamma: {\bf Y}_1\rightarrow {\bf Y}$ et $\zeta: {\bf X}_1\rightarrow {\bf X}$, tels que:
\begin{enumerate}
\item[(1)]le morphisme $\alpha_1$ est fini et tale;
\item[(2)]l'application $\gamma_F:{\bf Y}_1(F)\rightarrow {\bf Y}(F)$ est un homomorphisme;
\item[(3)]l'application $\zeta_F:{\bf X}_1(F)\rightarrow {\bf X}(F)$ induit par restriction un homomorphisme
$$\alpha_{1,F}({\bf Y}_1(F))\rightarrow \alpha_F({\bf Y}(F));$$
\item[(4)]on a l'galit $\alpha_F= \zeta_F\circ \alpha_{1,F}\circ \gamma_F^{-1}$.
\end{enumerate}}
\ni En particulier, le morphisme (ferm) $\alpha_F:{\bf Y}(F)\rightarrow {\bf X}(F)$ est 
un homomorphisme local sur son image.

\subsection{}Dans tout l'article, on utilisera les notations et conventions d'\'ecriture 
suivantes.

Soit $H$ un groupe topologique. On 
appelle \textit{automorphisme} de $H$ un morphisme de groupes 
$H\rightarrow H$ qui est un hom\'eomorphisme, et \textit{caract\`ere} de $H$ (\`a ne pas 
confondre avec le caract\`ere d'une repr\'esentation de $H$) 
un morphisme de groupes continu $H\rightarrow 
\Bbb{C}^\times$. On note ${\rm Aut}(H)$ le groupe 
des automorphismes de $H$, et ${\rm Int}(H)$ le sous--groupe distingu\'e 
de ${\rm Aut}(H)$ form\'e des 
automorphismes int\'erieurs, \cad ceux de la forme
${\rm Int}_H(x):h\mapsto xhx^{-1}$ pour un $x\in H$. 
Si $\chi$ est un caract\`ere de $H$, on note $\chi^{-1}$ le 
caract\`ere $h\mapsto \chi(h)^{-1}=\chi(h^{-1})$ de $H$.

Soit $X$ un \textit{td--espace}, \cad un espace topologique s\'epar\'e 
tel que chaque point de $X$ poss\`ede une base de voisinages ouverts compacts. 
On note $C^\infty_{\rm c}(X)$ l'espace des fonctions complexes sur $X$ qui sont localement 
constantes et \`a support compact. Si 
$X$ est discret (resp. compact), l'espace $C^\infty_{\rm c}(X)$ 
est aussi not\'e $C_{\rm c}(X)$ (resp. $C^\infty(X)$). Soit 
$K$ et $K'$ deux groupes topologiques. 
Si $X$ est muni de deux actions continues 
$$
K\times X \rightarrow X,\,(k,x)\mapsto k\cdot x,\quad  
X\times K'\rightarrow X,\,(x,k')\mapsto x\cdot k'$$ commutant entre elles, 
on note $C^\infty_{\rm c}(K\backslash X/K')$ le sous--espace de 
$C^\infty_{\rm c}(X)$ form\'e des fonctions $f$ telles que $f(k\cdot x\cdot k')= f(x)$ pour 
tout $(k,x,k')\in K\times X\times K'$. Rappelons qu'on appelle \textit{distribution sur $X$} un 
\'el\'ement de l'espace dual $C^\infty_{\rm c}(X)^*={\rm Hom}_\Bbb{C}(C^\infty_{\rm c}(X),\Bbb{C})$.

Soit $G$ un groupe topologique localement profini. On note $\mathfrak{R}(G)$ la 
cat\'egorie des repr\'esentations complexes lisses de $G$, et l'on appelle simplement 
\textit{repr\'esentation lisse de $G$} un objet de $\mathfrak{R}(G)$. Si $(\pi,V)$ est une 
repr\'esentation lisse de $G$, pour tout sous--groupe 
ouvert compact $K$ de $G$, on note $V^K$ le sous--espace 
$\{v\in V:\pi(k)(v)=v,\,\forall k\in K\}$ de $V$. Rappelons qu'une repr\'esentation 
lisse $(\pi,V)$ de $G$ est dite \textit{admissible} si pour tout sous--groupe 
ouvert compact $K$ de $G$, l'espace $V^K$ est de dimension finie.

\vskip2mm
\begin{center}
\rule{45mm}{0.1mm}
\end{center}

Je remercie Guy Henniart pour ses nombreuses remarques et suggestions, qui m'ont permis peu à peu d'améliorer ce texte. C'est d'ailleurs lui qui au départ m'a incité à regarder les caractères tordus des représentations admissibles en caractéristique quelconque.

Je remercie vivement le rapporteur pour sa lecture minutieuse et critique du manuscrit, et pour m'avoir signalé un certain 
nombre d'erreurs vraiment gênantes. 

\section{Caractres tordus d'un groupe localement profini}

\v2 Dans ce chapitre, on fixe un groupe localement profini $G$ et une mesure de Haar {\it \`a gauche} $d_lg$ sur $G$.

\subsection{Module d'un automorphisme de $G$}\label{modules} Soit 
$\Delta_G:G\rightarrow \Bbb{R}_{>0}$ le {\it module de $G$}, i.e. le 
caract\`ere r\'eel de $G$ d\'efini (comme 
dans \cite[1.19]{BZ}) par
$$
\int_Gf(gx^{-1})d_lg=\Delta_G(x)\int_Gf(g)d_lg\quad (f\in C^\infty_{\rm c}(G),\,x\in G).
$$
En d'autres termes, pour $x\in G$, on a (abus d'\'ecriture) $d_l(gx)=\Delta_G(x)d_lg$.

Soit $\theta$ un automorphisme de $G$. La distribution $\mu$ sur $G$ d\'efinie par 
$$\mu(f)=\int_Gf(\theta^{-1}(g))d_lg\quad (f\in C^\infty_{\rm c}(G))
$$ est encore une mesure de Haar \`a gauche. Par cons\'equent 
il existe une constante $\Delta_G(\theta)\in \Bbb{R}_{>0}$, appel\'ee {\it module de $\theta$}, 
telle que
$$
\int_Gf(\theta^{-1}(g))d_lg=\Delta_G(\theta)\int_Gf(g)d_lg\quad (f\in C^\infty_{\rm c}(G)).
$$
En d'autres termes, on a (nouvel abus d'\'ecriture) $d_l(\theta(g))=\Delta_G(\theta)d_lg$.
L'application
$$
\Delta_G:{\rm Aut}(G)\rightarrow \Bbb{R}_{>0}
$$
ainsi d\'efinie ne d\'epend 
pas du choix de la mesure $d_lg$, et c'est un 
morphisme de groupes. Notons que pour $x\in G$, on a 
$$\Delta_G(x)=\Delta_G({\rm Int}_G(x^{-1})).
$$

\begin{monlem}
\begin{enumerate}
\item[(1)]On a $\Delta_G\circ\theta =\Delta_G$.
\item[(2)]S'il existe une partie ouverte compacte $\theta$--stable de $G$, 
alors $\Delta_G(\theta)=1$.
\end{enumerate}
\end{monlem}

\begin{proof}
Pour $x\in G$, on a
$$
\Delta_G(\theta(x))=\Delta_G({\rm Int}_G(\theta(x)^{-1}))=
\Delta_G(\theta\circ {\rm Int}_G(x^{-1})\circ \theta^{-1}).
$$
D'o\`u le point (1). Quant au point (2), il r\'esulte de ce que pour toute partie ouverte 
compacte $\Omega$ de $G$, on a ${\rm vol}(\theta(\Omega),d_lg)=\Delta_G(\theta){\rm vol}(\Omega,d_lg)$. 
\end{proof}

Soit $H$ un sous--groupe ferm\'e de $G$. Rappelons que pour toute partie ouverte $X$ de $G$ 
telle que $HX=X$, les conditions suivantes sont \'equivalentes:
\begin{itemize}
\item l'espace quotient $H\backslash X$ est compact;
\item il existe une partie compacte $\Omega$ de $G$, que l'on peut choisir ouverte compacte, 
telle que $X=H\Omega$.
\end{itemize}
Si ces deux conditions sont v\'erifi\'ees, on dit que $X$ est {\it compact modulo $H$}. 
Bien s\^ur, on peut d\'efinir la m\^eme notion pour une partie ouverte $X$ de $G$ 
telle que $XH=X$. Notons que si $X$ est un sous--groupe ouvert de $G$ tel que $HX=X$, 
alors $XH=X$ et les deux notions co\"{\i}ncident.

Soit $\boldsymbol{S}(H\backslash G)$ l'espace vectoriel des fonctions 
$f:G\rightarrow \Bbb{C}$ 
telle que:
\begin{itemize}
\item $f(hg)=\Delta_G(h)\Delta_H(h^{-1})f(g)$ ($h\in H$, $g\in G$);
\item il existe un sous--groupe ouvert compact $K_f$ de $G$ tel que $f(gk)=f(g)$ 
pour tout $g\in G$ et tout  $k\in K_f$;
\item il existe une partie ouverte compacte $\Omega_f$ de $G$ telle que pour 
$g\in G\smallsetminus H\Omega_f$, on a $f(g)=0$.
\end{itemize}
Notons que pour $f\in \boldsymbol{S}(H\backslash G)$, le support
$$
{\rm Supp}(f)=\{g\in G: f(g)\neq 0\}
$$ 
de $f$ est une partie ouverte ferm\'ee de $G$, compacte modulo $H$.

Soit $d_lh$ une mesure de Haar \`a gauche sur $H$, 
et $C^\infty_{\rm c}(G)\rightarrow \boldsymbol{S}(H\backslash G),\,f\mapsto \bar{f}$ 
l'application lin\'eaire d\'efinie par
$$
\bar{f}(g)=\Delta_G(g)\int_Hf(hg)d_lh\quad (f\in C^\infty_{\rm c}(G)).
$$
Soit $\bar{\mu}$ une mesure de Haar \`a droite sur l'espace quotient 
$H\backslash G$ \cite[1.21]{BZ}, et 
$\mu$ la distribution sur $G$ d\'efinie par
$$
\mu(f)=\bar{\mu}(\bar{f})\quad (f\in C^\infty_{\rm c}(G)).
$$
Pour $f\in C^\infty_{\rm c}(G)$ et $x\in G$, notant $\rho_x(f)\in C^\infty_{\rm c}(G)$ la 
fonction $g\mapsto f(gx)$, on a l'galit $\overline{\rho_x(f)}=\Delta_G(x^{-1})\bar{f}$, d'o\`u
$$
\mu(\rho_x(f))= \Delta_G(x^{-1})\bar{\mu}(\bar{f})=\Delta_G(x^{-1})\mu(f).
$$
On en d\'eduit que $\Delta_G^{-1}\mu$ est une mesure de Haar \`a droite sur 
$G$, i.e. que $\mu$ est une mesure de Haar \`a gauche sur $G$ \cite[1.19]{BZ}. 

Supposons de plus que le groupe $H$ est $\theta$--stable. D'apr\`es 
le lemme 1, pour $f\in \boldsymbol{S}(H\backslash G)$, on a $f\circ \theta^{-1}
\in \boldsymbol{S}(H\backslash G)$, 
et la forme lin\'eaire $\bar{\mu}'$ sur $\boldsymbol{S}(H\backslash G)$ d\'efinie par 
$\bar{\mu}'(f)=\bar{\mu}(f\circ\theta^{-1})$ est une mesure de Haar \`a droite 
sur $H\backslash G$. 
On a donc $\bar{\mu}'=\Delta_{H\backslash G}(\theta)\bar{\mu}$ pour une 
constante $\Delta_{H\backslash G}(\theta)>0$. Pour $f\in C^\infty_{\rm c}(G)$ et $g\in G$, 
on a
\begin{align*}
\overline{(f\circ \theta^{-1})}(g)&= \Delta_G(g)\int_Hf(\theta^{-1}(hg))d_lh\\
&= \Delta_G(\theta^{-1}(g))\Delta_H(\theta)\int_Hf(h\theta^{-1}(g))d_lh\\
&= \Delta_H(\theta)(\bar{f}\circ \theta^{-1})(g).
\end{align*}
Puisque $\mu(f\circ\theta^{-1})=\Delta_G(\theta)\mu(f)$ ($f\in C^\infty_{\rm c}(G)$), on 
en d\'eduit l'\'egalit\'e
$$
\Delta_G(\theta)=\Delta_H(\theta)\Delta_{H\backslash G}(\theta).\leqno{(*)}
$$
On peut bien sûr définir un caractère--module $\Delta_{H\backslash G}: H\rightarrow {\Bbb R}_{>0}$ en posant
$$
\Delta_{H\backslash G}(h) = \Delta_G(h)\Delta_H(h)^{-1}\quad (h\in H).
$$
Pour $h\in H$, l'automorphisme intérieur ${\rm Int}_G(h)$ de $G$ stabilise $H$, et d'après $(*)$, on a
$$
\Delta_{H\backslash G}(h) = \Delta_G({\rm Int}_G(h^{-1}))\Delta_H({\rm Int}_H(h))=\Delta_{H\backslash G}({\rm Int}_G(h^{-1})).
$$

\begin{marema}
{\rm 
Le caractère $\Delta_{H\backslash G}$ de $H$ vérifie encore $\Delta_{H\backslash G}\circ \theta= \Delta_{H\backslash G}$.  En revanche, l'existence d'une partie ouverte compacte $\theta$--stable de $H\backslash G$ n'implique pas que $\Delta_{H\backslash G}(\theta)=1$. En particulier, même si le groupe $G$ est compact modulo $H$, on a en général $\Delta_{H\backslash G}(\theta)\neq 1$. Prenons par exemple pour $G$ le groupe $GL(2,F)$ des points $F$--rationnels du groupe linéaire $GL(2)$, où $F$ est un corps localement compact non archimédien, et pour $H$ le sous--groupe de $GL(2,F)$ formé des matrices triangulaires supérieures. Alors $G$ est compact modulo $H$, et pour un élément $x$ dans le tore diagonal de $G$, l'automorphisme $\theta={\rm Int}_G(x)$ de $G$ stabilise $H$. On a $\Delta_G(\theta)=1$ mais $\Delta_H(\theta)= \Delta_H(x^{-1})$ peut ne pas être égal à $1$.}
\end{marema}

\subsection{Caract\`eres des repr\'esentations admissibles de $G$ (rappels)}\label{caractres} 
Soit $(\pi,V)$ une 
repr\'e\-sentation lisse de $G$. Pour $f\in C^\infty_{\rm c}(G)$, on note 
$\pi(f)=\pi(fd_lg)\in {\rm End}_\Bbb{C}(V)$ l'op\'erateur d\'efini par
$$
\pi(f)(v)= \int_G f(g)\pi(g)(v)d_lg\quad 
(v\in V);
$$
l'int\'egrale est en fait une somme finie. Si $K$ est un sous--groupe ouvert compact de 
$G$, alors pour toute fonction $f\in C_{\rm c}(K\backslash G)$, on a $\pi(f)(V)\subset V^K$. 
Supposons de plus que $\pi$ est admissible. Alors pour $f\in C^\infty_{\rm c}(G)$, 
l'op\'erateur $\pi(f)$ est de 
rang fini, 
et l'on peut d\'efinir sa trace:
$$
\Theta_\pi(f)={\rm tr}(\pi(f)).
$$
La distribution $\Theta_\pi$ sur $G$ ainsi d\'efinie, est appel\'ee le {\it caract\`ere de 
$\pi$} (elle d\'epend bien s\^ur du choix de la mesure $d_lg$). 
Elle est invariante par conjugaison: pour 
$f\in C^\infty_{\rm c}(G)$ et $x\in G$, notant ${^{x\!}f}=
{\rm Int}_G(x)(f)\in C^\infty_{\rm c}(G)$ la fonction 
$g\mapsto f(x^{-1}gx)$, 
on a l'\'egalit\'e
$$
\Theta_\pi({^{x\!}f})=\Delta_G(x^{-1})\Theta_\pi(f).
$$

\subsection{Caract\`eres tordus de $G$}\label{caractres tordus} 
Soit $(\pi,V)$ une repr\'esentation admissible de $G$, et soit 
un op\'erateur $A\in {\rm End}_\Bbb{C}(V)$. On 
note $\Theta_\pi^A$ la distribution sur $G$ d\'efinie par
$$
\Theta_\pi^A(f)={\rm tr}( \pi(f)\circ A) \quad (f\in C^\infty_{\rm c}(G)).
$$
On l'appelle le {\it caract\`ere $A$--tordu de $\pi$} (tout comme $\Theta_\pi$, la distribution 
$\Theta_\pi^A$ d\'epend du choix de la mesure $d_lg$). Notons que 
pour $f\in C^\infty_{\rm c}(G)$ et $x\in G$, 
on a l'\'egalit\'e
$$
\Theta_\pi^A({^{x\!}f})=\Delta_G(x^{-1})\Theta_\pi^{\pi(x^{-1})\circ A \circ \pi(x)}(f).\leqno{(*)}
$$
Si de plus $A\in {\rm Aut}_\Bbb{C}(V)$, alors notant $(\pi'\!,V)$ la repr\'esentation 
$g\mapsto A\circ \pi(g)\circ A^{-1}$ de $G$, on a $A\in {\rm Isom}_G(\pi,\pi')$.

\begin{marema}
{\rm Supposons que les sous--espace $A(V)$ et $\ker A$ de $V$ sont $G$--stables  
 (ce qui est toujours le cas si $A$ est un $\Bbb{C}$--automorphisme de $V$!). 
 Alors l'espace quotient $\overline{V}=V/\ker A$ d\'efinit une repr\'esentation 
 quotient 
(donc admissible)
 $\bar{\pi}$ de $\pi$. Notons $\bar{A}\in {\rm Isom}_\Bbb{C}(\overline{V},A(V))$ l'op\'erateur 
 d\'eduit de $A$ par passage au quotient. Il d\'efinit une repr\'esentation admissible  
 $(\pi'\!,A(V))$ de $G$: pour $g\in G$ et $v\in V$, on pose $\pi'(g)(v)=
 \bar{A}\circ \bar{\pi}(g)\circ \bar{A}^{-1}(v)$. En d'autres termes, on a 
 $\bar{A}\in {\rm Isom}_G(\bar{\pi},\pi')$, et 
 pour $f\in C^\infty_{\rm c}(G)$, l'op\'erateur 
 $A\circ \pi(f)$ appartient \`a ${\rm End}_\Bbb{C}(\overline{V},A(V))$.\hfill $\blacksquare$}
 \end{marema}
 
Soit maintenant $\theta$ un automorphisme de $G$, et $\omega$ un caract\`ere de 
$G$. Si $\pi$ est une repr\'esentation lisse de $G$, on note $\omega\pi$ et $\pi^\theta$ les 
repr\'esentations (lisses) de $G$ d\'efinies par
$$
\omega\pi = \pi\otimes \omega,\quad 
\pi^\theta=\pi\circ \theta.
$$
Soit $\pi$ une repr\'esentation admissible de $G$ telle que $\omega\pi\simeq\pi^\theta$, et 
soit $A\in {\rm Isom}_G(\omega\pi,\pi^\theta)$. 
Alors la distribution $\Theta_\pi^A$ est $\omega$--invariante par 
$\theta$--conjugaison: pour 
 $f\in C^\infty_{\rm c}(G)$ et $x\in G$, notant ${^{x,\theta\!}f}= {\rm Int}_{G,\theta}(x)(f)
 \in C^\infty_{\rm c}(G)$ la fonction $g\mapsto 
 f(x^{-1}g\theta(x))$, on a l'\'egalit\'e 
 $$
 \Theta_\pi^A({^{x,\theta\!}f})=\Delta_G(x^{-1})\omega(x^{-1})\Theta_\pi^A(f);
 \leqno{(**)}
 $$
 en utilisant que  $\Delta_G(\theta(x^{-1}))=\Delta_G(x^{-1})$ (lemme de \ref{modules}).

\subsection{Espaces topologiques tordus}\label{espaces topo tordus}
On peut modifier la situation pr\'ec\'edente 
de mani\`ere \`a \og voir \fg{} la $\theta$--conjugaison comme une 
conjugaison ordinaire. 
Soit $H$ un groupe topologique. 
\`A la suite de J.-P. Labesse [La], appelons $H$--{\it espace (topologique) tordu}\footnote{Labesse travaille 
en fait dans la cat\'egorie des groupes alg\'ebriques (cf. plus loin, \ref{éléments réguliers d'un espace tordu}), mais ses constructions se 
transposent ais\'ement \`a la cat\'egorie des groupes topologiques; 
c'est ce que nous faisons ici. D'ailleurs la notion d'espace tordu \'etant tout aussi 
g\'en\'erale que celle d'espace principal homog\`ene, elle s'adapte \`a 
n'importe quelle cat\'egorie de groupes.} la donn\'ee:
\begin{itemize}
\item d'un $H$--espace principal homog\`ene topologique $H^\natural$, \cad 
un espace topologique $H^\natural$ muni d'une op\'eration continue de 
$H$ \`a gauche, not\'ee 
$H\times H^\natural\rightarrow H^\natural,\,(h,\delta)\mapsto h\cdot\delta$, telle 
que pour un (et alors pour tout) $\delta\in H^\natural$, l'application 
$H\rightarrow H^\natural,\,h\mapsto h\cdot\delta$ est un hom\'eomorphisme;
\item et d'une application 
${\rm Int}_{H^\natural}:H^\natural\rightarrow {\rm Aut}(H)$ telle que
$$
{\rm Int}_{H^\natural}(h\cdot\delta)={\rm Int}_H(h)\circ {\rm Int}_{H^\natural}(\delta)\quad (h\in H,\, \delta\in H^\natural).
$$
\end{itemize}
On peut alors d\'efinir une op\'eration continue de $H$ sur $H^\natural$ \`a droite, 
$H^\natural\times H\rightarrow 
H^\natural,\, (\delta,h)\mapsto \delta\cdot h$, qui commute \`a l'action \`a gauche:
$$
 \delta\cdot h= {\rm Int}_{H^\natural}(\delta)(h)\cdot \delta.
$$
Pour $h\in H$, notons ${\rm Int}'_H(h):H^\natural\rightarrow H^\natural$ 
l'hom\'eomorphisme d\'efini par
$$
{\rm Int}'_H(h)(\delta)= h\cdot \delta \cdot h^{-1}= h\h{0.2}{\rm Int}_{H^\natural}(\delta)(h^{-1})\cdot \delta.
$$

L'image de $H^\natural$ 
dans le groupe ${\rm Out}(H)={\rm Aut}(H)/{\rm Int}(H)$ 
par l'application compos\'ee
$$
H^\natural \buildrel {\rm Int}_{H^\natural} \over{\hbox to 10mm{\rightarrowfill}} {\rm Aut}(H)\longrightarrow {\rm Out}(H),
$$
est r\'eduite \`a un point. Soit $\theta={\rm Int}_{H^\natural}(\delta_1)$ pour un 
$\delta_1\in H^\natural$. C'est un automorphisme de $H$ qui, par restriction, induit un automorphisme du 
centre $Z(H) $ de $H$. Ce dernier ne 
d\'epend pas du choix de $\delta_1$. On peut donc poser
$$
Z({H^\natural})=\{z\in Z(H):\theta(z)=z\}.
$$ 
C'est un sous--groupe ferm\'e de $Z(H)$, qui co\"{\i}ncide avec 
le centralisateur de $H^\natural$ dans $H$, i.e. on a
\begin{align*}
Z(H^\natural)&= \{h\in H: {\rm Int}_{H^\natural}(\delta)(h)=h,\,\forall \delta\in H^\natural\}\\ 
&= \{h\in H: h\cdot \delta\cdot h^{-1}=\delta,\,\forall \delta\in H^\natural\}.
\end{align*}

Si pour tout $\delta\in H^\natural$, on a ${\rm Int}_{H^\natural}(\delta)={\rm id}_H$, on dit que 
$H^\natural$ est un {\it $H$--espace tordu trivial}. En particulier, 
le groupe topologique $H$ est naturellement un $H$--espace tordu, et c'est un 
$H$--espace tordu trivial si et seulement s'il est commutatif.

\begin{marema1}{\rm Pour $\theta\in {\rm Aut}(H)$, on note $H\theta$ 
le sous-ensemble $ H \rtimes \theta$ 
de $H\rtimes {\rm Aut}(H)$, muni de la topologie rendant continue 
la bijection $H\rightarrow H\theta,\, h\mapsto h\theta =h\rtimes \theta$ (ainsi 
cette bijection est un hom\'eomorphisme). C'est un $H$--espace tordu. 
L'action $H\times H\theta\rightarrow H\theta$ est bien s\^ur 
donn\'ee par $x\cdot h\theta = (xh)\theta$ ($x,\,h\in H$), et l'application 
${\rm Int}_{H\theta}:H\theta\rightarrow {\rm Aut}(H)$ est donn\'ee par 
${\rm Int}_{H\theta}(h\theta)={\rm Int}_H(h)\circ \theta$ ($h\in H$).  Alors pour $x\in H$, 
l'application ${\rm Int}'_H(x^{-1}):H\theta\rightarrow H\theta$ est donn\'ee par
$$
{\rm Int}'_H(x^{-1})(h\theta)=(x^{-1}h\theta(x))\theta\quad (h\in H).\eqno{\blacksquare}
$$
}
\end{marema1}

\v1 Appelons {\it espace topologique tordu} 
la donn\'ee d'un groupe 
topologique $H$ et d'un $H$--espace tordu $H^\natural$; on note 
simplement $H^\natural$ l'espace topologique tordu d\'efini par la paire 
$(H,H^\natural)$, le groupe $H$ \'etant sous-entendu. 
Les espaces topologiques tordus s'organisent naturellement en une cat\'egorie: un {\it morphisme} entre 
deux espaces topologiques tordus $H_1^\natural$ et $H_2^\natural$ est la donn\'ee 
d'un morphisme de groupes continu $\phi:H_1\rightarrow H_2$ et 
d'une application 
$\Phi: H_1^\natural\rightarrow H_2^\natural$, v\'erifiant
$$
\Phi(x\cdot \delta\cdot y)=\phi(x)\cdot \Phi(\delta)\cdot \phi(y)\quad (\delta\in H_1;\;x,\,y\in H_1).
$$
Notons que $\Phi$ est une application continue, et que $\phi$ est enti\`erement d\'etermin\'e 
par $\Phi$. On pose donc 
$\Phi^\circ =\phi$, et l'on note simplement $\Phi$ le morphisme d\'efini par la paire $(\phi,\Phi)$. 
L'application $\Phi$ est un hom\'eomorphisme 
si et seulement si le morphisme $\Phi^\circ$ en est un, auquel cas on 
dit que le morphisme $\Phi $ est un isomorphisme (resp. un automorphisme si $H_1^\natural=
H_2^\natural$). 
Si $H^\natural$ un espace topologique tordu, 
on note ${\rm Aut}(H^\natural)$ le groupe des automorphismes de 
$H^\natural$.

On a des notions \'evidentes de sous--espace topologique tordu et d'espace topologique 
tordu quotient. Soit $H^\natural$ un espace topologique tordu. On appelle {\it sous--espace topologique 
tordu de $H^\natural$} la donn\'ee d'un sous--groupe topologique $J$ de $H$ et d'un sous--espace topologique $J^\natural$ de $H^\natural$, tels que les applications $H\times H^\natural\rightarrow 
H^\natural$ et ${\rm Int}_{H^\natural}:H^\natural\rightarrow {\rm Aut}(H)$ induisent sur $J^\natural$ une 
structure de $J$--espace tordu. On a donc $J^\natural = J\cdot \delta$ pour un \'el\'ement 
$\delta\in H^\natural$ tel que ${\rm Int}_{H^\natural}(\delta)(J)=J$; mais l'\'el\'ement $\delta$ n'est pas d\'etermin\'e 
de mani\`ere unique par $J^\natural$. 
De m\^eme, on appelle {\it espace topologique tordu quotient de $H^\natural$} 
la donn\'ee d'un sous--groupe topologique 
distingu\'e $I$ de $H$ tel que ${\rm Int}_{H^\natural}(\delta)(I)=I$ 
pour tout (i.e. pour un) $\delta\in H^\natural$; alors l'espace topologique quotient $H^\natural\!/I =I\backslash H^\natural$ des classes d'\'equivalence dans $H^\natural$ pour la 
translation \`a droite (resp. \`a gauche) par les \'el\'ements de $I$, est un 
($H/I$)--espace tordu. Plus g\'en\'eralement, pour tout sous--groupe topologique $I$ 
(non n\'ecessairement distingu\'e) de $H$ tel que ${\rm Int}_{H^\natural}(\delta)(I)=I$ pour un 
$\delta\in H^\natural$, on d\'efinit de la m\^eme mani\`ere les espaces topologiques quotients 
$H^\natural/I$ et $I\backslash H^\natural$, qui ne sont en g\'en\'eral pas des espaces tordus.

Soit $\Phi:H^\natural_1\rightarrow H^\natural_2$ un morphisme d'espaces topologiques tordus. 
L'image ${\rm Im}(\Phi)=\Phi(H^\natural_1)$ est naturellement munie d'une 
structure de $\Phi^\circ(H_1)$-espace tordu, mais il n'est en g\'en\'eral pas possible de d\'efinir 
le noyau de $\Phi$; il faudrait pour cela travailler avec des espaces topologiques tordus 
{\it point\'es}. 
En revanche, le noyau $\ker(\Phi^\circ)$ est bien d\'efini dans la 
cat\'egorie des groupes topologiques; et pour 
$\delta\in H_1^\natural$ et $x\in \ker(\Phi^\circ)$, puisque
$$
\Phi(\delta)= \Phi(\delta\cdot x)=\Phi({\rm Int}_{\smash{H_1^\natural}}(\delta)(x)\cdot \delta)= \Phi^\circ({\rm Int}_{\smash{H_1^\natural}}(\delta)(x))\cdot \Phi(\delta),
$$
on a ${\rm Int}_{\smash{H_1^\natural}}(\delta)(x)=x$. On peut donc d\'efinir l'espace topologique tordu quotient 
$H_1^\natural/\ker(\Phi^\circ)$ de $H_1^\natural$. Par construction, le morphisme 
$\Phi$ induit par passage au quotient un isomorphisme d'espaces topologiques 
tordus
$$
H_1^\natural/\ker(\Phi^\circ)\rightarrow \Phi(H_1^\natural).
$$

Notons que m\^eme si la cat\'egorie des espaces tordus n'est pas ab\'elienne, on 
a des notions naturelles de morphisme injectif (resp. surjectif), et de suite exacte courte.

\begin{marema2}{\rm 
Soit $H^\natural$ un $H$--espace tordu, et soit $\delta_1\in H^\natural$. Posons $\theta={\rm Int}_{H^\natural}(\delta_1)$.
Pour $z\in Z(H)$, la translation \`a gauche 
$H^\natural \rightarrow 
H^\natural,\,\delta\mapsto \Phi_z(\delta)=z\cdot 
\delta$ est un \'el\'ement de ${\rm Aut}(H^\natural)$, et l'application 
$Z(H)\rightarrow {\rm Aut}(H^\natural),\,z\mapsto \Phi_z$ est un 
morphisme de groupes injectif; on note $Z(H)'$ son image. Le groupe $Z(H)'$ 
co\"{\i}ncide avec le noyau du morphisme de groupes 
${\rm Aut}(H^\natural)\rightarrow {\rm Aut}(H),\,\Phi\mapsto \Phi^\circ$. En effet, 
pour $z\in Z(H)$, l'automorphisme $\Phi_z$ appartient \`a ce noyau. 
R\'eciproquement, soit $\Phi\in {\rm Aut}(H^\natural)$ tel que $\Phi^\circ ={\rm id}_H$. 
On a $\Phi(\delta_1)= h_1\cdot \delta_1$ pour un (unique) $h_1\in H$, et 
pour $x\in H$, on a
$$
h_1x\cdot \delta_1 = h_1\cdot \delta_1\cdot \theta^{-1}(x)= \Phi(\delta_1\cdot 
\theta^{-1}(x))=
\Phi(x\cdot \delta_1)= xh_1\cdot \delta_1.
$$
Par cons\'equent $h_1\in Z(H)$, et l'on a bien l'\'egalit\'e $Z(H)'= \ker\{{\rm Aut}(H^\natural)
\rightarrow {\rm Aut}(H)\}$. D'autre part,  pour 
$h\in H$, l'application ${\rm Int}'_H(h):H^\natural\rightarrow H^\natural$ 
est un \'el\'ement de ${\rm Aut}(H^\natural)$. Les ${\rm Int}'_H(h)$ 
pour $h\in H$, engendrent 
un sous--groupe de ${\rm Aut}(H^\natural)$, que l'on note ${\rm Int}(H)'$. Soit 
${\rm Int}(H^\natural)$ le sous--groupe de ${\rm Aut}(H^\natural)$ engendr\'e par 
$Z(H)'$ et ${\rm Int}(H)'$. Alors on a la suite exacte courte 
de groupes
$$
1\rightarrow Z(H)\rightarrow {\rm Int}(H^\natural)\rightarrow {\rm Int}(H)\rightarrow 1.
$$
Cette suite n'est en g\'en\'eral pas scind\'ee, sauf si la projection canonique 
$H/Z(H^\natural)\rightarrow 
H/Z(H)$ l'est. 
Notons que pour $z\in Z(H)$ et $h\in H$, on a
$$
{\rm Int}'_H(zh)= \Phi_{z\theta(z^{-1})}\circ {\rm Int}'_H(h).\eqno{\blacksquare}
$$
}
\end{marema2}

\begin{marema3}
{\rm Continuons avec les notations de la remarque prcdente: $\theta= {\rm Int}_{H^\natural}(\delta_1)$ pour un $\delta_1\in H^\natural$. 
L'application $\iota_{\delta_1}:H^\natural\rightarrow H\theta,\, h\cdot\delta_1\mapsto h\theta$ est 
un isomorphisme d'espaces topologiques tordus, tel que $\iota_{\delta_1}^\circ ={\rm id}_H$. 
Mais d'apr\`es la remarque 2, il d\'epend du choix de $\delta_1$, sauf si $Z(H)=\{1\}$. 
C'est pourquoi l'on pr\'ef\'erera travailler avec un $H$--espace tordu $H^\natural$ muni d'un point--base $\delta_1$, 
plut\^ot qu'avec un $H$--espace tordu $H\theta$.\hfill $\blacksquare$}
\end{marema3}

\subsection{Module d'un $G$--espace tordu}\label{module d'un espace tordu} Soit $G^\natural$ un $G$--espace tordu. 
On appelle {\it mesure de Haar \`a gauche (resp. \`a droite) sur $G^\natural$} une 
distribution $\mu$ sur $G^\natural$, invariante pour l'action de $G$ par translations 
\`a gauche (resp. \`a droite), et telle que pour toute fonction $\phi\in C^\infty_{\rm c}(G^\natural)$, $\phi\geq0$ 
et $\phi\neq 0$, on a $\mu(\phi)>0$. En d'autres termes, fix\'e un \'el\'ement 
$\delta_1\in G^\natural$, les mesures de Haar \`a gauche (resp. \`a droite) 
sur $G^\natural$ sont les images des mesures de Haar \`a gauche (resp. \`a droite) sur 
$G$ par l'hom\'eomorphisme $r_{\delta_1}:G\rightarrow G^\natural,\,g\mapsto g\cdot \delta_1$; 
de mani\`ere \'equivalente, ce sont les images des mesures de Haar \`a 
gauche (resp. \`a droite) sur $G$ par l'hom\'eomorphisme $l_{\delta_1}:G\rightarrow 
G^\natural,\,g\mapsto \delta_1\cdot g$. Si $\mu$ est une mesure de Haar \`a gauche 
ou \`a droite sur $G$, on note $\mu\cdot \delta_1$ son image par $r_{\delta_1}$, et 
$\delta_1\cdot \mu$ son image par $l_{\delta_1}$. 
On a
$$
\delta_1\cdot \mu = \Delta_G({\rm Int}_{G^\natural}(\delta_1))(\mu\cdot \delta_1).\leqno{(*)}
$$

Si $\mu$ est une mesure 
de Haar \`a gauche sur $G$, alors la mesure $\delta_1\cdot \mu$ sur $G^\natural$ 
ne d\'epend pas du choix de $\delta_1$; on l'appelle la {\it mesure de Haar \`a gauche sur 
$G^\natural$ associ\'ee \`a $\mu$}. De m\^eme, si $\mu$ est une mesure 
de Haar \`a droite sur $G$, alors la mesure $\mu\cdot \delta_1$ sur $G^\natural$ 
ne d\'epend pas du choix de $\delta_1$; on l'appelle la {\it mesure de Haar \`a droite sur 
$G^\natural$ associ\'ee \`a $\mu$}.

Soit $\Delta_{G^\natural}:G^\natural\rightarrow \Bbb{R}_{>0}$ 
le {\it module de $G^\natural$}, d\'efini par
$$
\Delta_{G^\natural}(\gamma)=\Delta_G({\rm Int}_{G^\natural}(\gamma)^{-1}).\leqno{(**)}
$$
Soit $d_l\gamma$ une mesure de Haar \`a gauche sur $G^\natural$. 

\begin{monlem}Soit $x\in G$ et $\gamma\in G^\natural$.
\begin{enumerate}
\item[(1)]On a (abus d'\'ecriture) $d_l(\gamma\cdot x)= \Delta_G(x)d_l\gamma$.
\item[(2)]On a $\Delta_{G^\natural}(\gamma\cdot x)=\Delta_{G^\natural}(x\cdot \gamma)= \Delta_G(x)
\Delta_{G^\natural}(\gamma)$.
\end{enumerate}
\end{monlem}

\begin{proof}
On peut supposer que $d_l\gamma= \delta_1\cdot d_lg$. 
Alors pour $\phi\in C^\infty_{\rm c}(G^\natural)$, on a
$$
\int_{G^\natural}\phi(\gamma\cdot x^{-1})d_l\gamma  = \int_G\phi(\delta_1\cdot gx^{-1})d_lg
= \Delta_G(x)\int_{G^\natural}\phi(\gamma)d_l\gamma.
$$
D'o\`u le point (1).

Posons $\theta={\rm Int}_{G^\natural}(\gamma)$. Puisque 
${\rm Int}_{G^\natural}(x\cdot \gamma)={\rm Int}_G(x)\circ \theta$, on a
$$
\Delta_{G^\natural}(x\cdot \gamma)=\Delta_G(\theta^{-1}\circ {\rm Int}_G(x^{-1}))=
\Delta_G(\theta^{-1})\Delta_G({\rm Int}_G(x^{-1})).
$$
D'o\`u l'\'egalit\'e $\Delta_{G^\natural}(x\cdot\gamma)= \Delta_G(x)\Delta_{G^\natural}(\gamma)$. 
Comme $\gamma\cdot x = \theta(x)\cdot \gamma$ et (1.1.1) $\Delta_G\circ\theta = \Delta_G$, 
on a aussi $\Delta_{G^\natural}(\gamma\cdot x)= \Delta_G(x)\Delta_{G^\natural}(\gamma)$. 
D'o\`u le point (2).
\end{proof}

D'apr\`es le lemme, pour $x\in G$, on a (abus d'\'ecriture) $\Delta_{G^\natural}
(\gamma\cdot x)^{-1}d_l(\gamma\cdot x)= \Delta_{G^\natural}(\gamma)^{-1}d_l\gamma$. On en d\'eduit que 
si $\mu$ 
est une mesure de Haar \`a gauche sur $G^\natural$, alors la distribution 
$\mu'$ sur $G^\natural$ d\'efinie par $d\mu'(\gamma)=
\Delta_{G^\natural}(\gamma)^{-1}d\mu(\gamma)$ 
est une mesure de Haar \`a droite sur $G^\natural$, et on les obtient 
toutes de cette mani\`ere. De plus, les trois conditions suivantes sont \'equivalentes:
\begin{itemize}
\item $G$ est unimodulaire;
\item le module $\Delta_{G^\natural}$ est constant;
\item les notions de mesures de 
Haar \`a gauche et \`a droite sur $G^\natural$ co\"{\i}ncident.
\end{itemize}
Si les conditions ci--dessus sont v\'erifi\'ees, on appelle simplement {\it mesure de Haar 
sur $G^\natural$}, une mesure de Haar \`a gauche (i.e \`a droite) 
sur $G^\natural$. Notons que si $G$ est unimodulaire et si 
$\mu$ est une mesure de Haar sur $G$, alors les mesures de Haar $\delta_1\cdot \mu$ et 
$\mu\cdot \delta_1$ sur $G^\natural$ sont reli\'ees par l'\'egalit\'e: 
$\mu\cdot \delta_1= \Delta_{G^\natural}(\delta_1)(\delta_1\cdot \mu)$. 

\v1 Si $G^\natural = G\theta$ pour un automorphisme $\theta$ de $G$, et si $\mu$ est une 
mesure de Haar \`a gauche ou \`a droite sur $G$, on pose 
$\mu\theta = \mu\cdot 1\theta$.

\subsection{Caract\`eres des $\omega$--repr\'esentations admissibles d'un $G$--espace tordu} \label{caractres tordus (bis)}
Soit $G^\natural$ un $G$--espace tordu, et $\omega$ un 
caract\`ere de $G$. On appelle {\it $\omega$--repr\'esentation lisse} 
{\it de $G^\natural$} la 
donn\'ee d'une repr\'esentation lisse 
$(\pi,V)$ de $G$ et d'une application 
$\Pi:G^\natural
\rightarrow {\rm Aut}_\Bbb{C}(V)$ v\'erifiant
$$
\Pi (x\cdot \gamma \cdot y)=\omega(y)\pi(x)\circ \Pi(\gamma)\circ \pi(y)\quad 
(\gamma\in G^\natural;\, x,\,y\in G).
$$
Notons que $\pi$ est enti\`erement dtermine par $\Pi$: on a
$$
\pi(x)= \Pi(x\cdot \gamma)\circ \Pi(\gamma)^{-1}\quad (x\in G,\, \gamma\in G^\natural).
$$
La $\omega$--representation lisse de $G$ d\'efinie par $(\pi,V)$ et 
$\Pi$ comme ci--dessus est not\'ee $(\Pi,V)$, ou plus 
simplement $\Pi$, et la repr\'esentation lisse $\pi$ de $G$ associ\'ee \`a $\Pi$ est 
not\'ee $\Pi^\circ$. On dit que $\Pi$ est {\it admissible} si $\Pi^\circ$ l'est.

Les $\omega$--repr\'esentations lisses de $G^\natural$ s'organisent naturellement en une 
cat\'egorie $\mathfrak{R}(G^\natural,\omega)$: si $(\Pi_1,V_1)$ et $(\Pi_2,V_2)$ sont 
deux $\omega$--repr\'esentations lisses de $G^\natural$, une fl\^eche entre 
$\Pi_1$ et $\Pi_2$ est par d\'efinition un {\it op\'erateur de $G^\natural$-entrelacement} 
entre $\Pi_1$ et $\Pi_2$, i.e. un morphisme $u\in {\rm End}_\Bbb{C}(V_1,V_2)$ tel que
$$
u\circ \Pi_1(\gamma)=\Pi_2(\gamma)\circ u\quad (\gamma\in G^\natural).
$$
On note ${\rm Hom}_{G^\natural}(\Pi_1,\Pi_2)$ l'espace des op\'erateurs d'entrelacement 
entre $\Pi_1$ et $\Pi_2$. Posons $\pi_i=\Pi_i^\circ$ ($i=1,\,2$). 
Tout op\'erateur $u\in {\rm Hom}_{G^\natural}(\Pi_1,\Pi_2)$ appartient 
\`a ${\rm Hom}_G(\pi_1,\pi_2)$. En effet, choisissons un $\delta\in G^\natural$ et posons 
$A_1=\Pi_1(\delta)$ et $A_2=\Pi_2(\delta)$. Comme $u\circ A_1= A_2\circ u$, 
pour tout $g\in G$, on a
\begin{align*}
u\circ \pi_1(g)\circ A_1&= u\circ \Pi_1(g\cdot \delta)\\
&= \Pi_2(g\cdot \delta)\circ u\\
&= \pi_2(g)\circ A_2\circ u\\
&=  \pi_2(g)\circ u\circ A_1
\end{align*}
d'o
$$u\circ \pi_1(g)=\pi_2(g)\circ u.
$$ En d'autres termes, l'application $\Pi\mapsto 
\Pi^\circ$ d\'efinit un foncteur d'oubli $\mathfrak{R}(G^\natural,\omega)\rightarrow \mathfrak{R}(G)$, et ce 
foncteur est fid\`ele. Pour plus de d\'etails sur les $\omega$--repr\'esentations lisses de $G^\natural$, 
en particulier celles qui sont \og irr\'eductibles\fg{}, on 
renvoie \`a l'annexe A.

\v1 Fixons un \'el\'ement $\delta_1\in G^\natural$ et posons $\theta={\rm Int}_{G^\natural}(\delta_1)$. 
Soit $\Pi$ une $\omega$--repr\'esentation lisse de $G^\natural$. Posons
$\pi=\Pi^\circ$ et $A=\Pi(\delta_1)$. Alors pour $g\in G$, on a $\pi(g)=\Pi(g\cdot\delta_1)\circ A^{-1}$ et 
$$
\omega(g) A\circ \pi(g)= \pi(\theta(g))\circ A;
$$
i.e. $A\in {\rm Isom}_G(\omega\pi,\pi^\theta)$. 
En d'autres termes, la donn\'ee d'une $\omega$--repr\'esentation lisse 
$\Pi$ de $G^\natural$ \'equivaut 
\`a celle d'une paire $(\pi,A)$ form\'ee d'une repr\'esentation lisse 
$\pi$ de $G$ 
telle que $\omega\pi\simeq \pi^\theta$, et 
d'un op\'erateur d'entrelacement $A\in {\rm Isom}_G(\omega\pi,\pi^\theta)$. De plus, si 
$\Pi_1$ et $\Pi_2$ sont deux $\omega$--repr\'esentations lisses 
de $G^\natural$, posant $\pi_i=\Pi_i^\circ$ et $A_i=\Pi_i(\delta_1)$ 
($i=1,\,2$), on a
$$
{\rm Hom}_{G^\natural}(\Pi_1,\Pi_2)=\{u\in {\rm Hom}_G(\pi_1,\pi_2): u\circ A_1 = A_2\circ u\}.
$$

\begin{marema}{\rm 
Pour qu'il existe une $\omega$--représentation $\Pi$ de $G^\natural$ telle que la représentation $\Pi^\circ$ de $G$ soit irréductible, 
il est n\'ecessaire que le 
caract\`ere $\omega$ de $G$ soit trivial sur $Z(G^\natural)$. En effet, 
l'ensemble $Z(G)(\theta-1)=\{\theta(z)z^{-1}:z\in Z(G)\}$ 
est un sous--groupe ferm\'e de $Z(G)$, qui ne d\'epend pas du choix de 
$\delta_1$, et l'on a la suite exacte courte de groupes topologiques
$$
1\rightarrow Z(G^\natural)\rightarrow Z(G)\rightarrow Z(G)(\theta-1)\rightarrow 1.
$$
Supposons qu'il existe une $\omega$--repr\'esentation 
lisse $(\Pi,V)$ de $G^\natural$ telle que la représentation $\pi=\Pi^\circ$ de $G$ est 
irréductible. D'après le lemme se Schur, pour $z\in Z(G)$, l'opérateur $\pi(z)$ est de la forme 
$\omega_\pi(z){\rm id}_V$ pour un élément $\omega_\pi(z)\in {\Bbb C}^\times$, et l'application $z\mapsto \omega_\pi(z)$ est un 
caractère de $Z(G)$. Pour $z\in Z(G)$ et $\gamma\in G^\natural$, posant $\pi=\Pi^\circ$, 
on a
$$
\omega_\pi(\theta(z))\Pi(\gamma)= \Pi(\theta(z)\cdot \gamma)= \Pi(\gamma \cdot z)= \omega(z)\omega_\pi(z)\Pi(\gamma),
$$
d'o\`u $\omega_\pi(z^{-1}\theta(z))= \omega(z)$. En particulier, on a 
$\omega\vert_{Z(G^\natural)}=1$.\hfill $\blacksquare$}
\end{marema}

Notons $d_l\gamma=\delta_1\cdot d_lg$ la mesure de 
Haar \`a gauche sur $G^\natural$ associ\'ee \`a $d_lg$. Rappelons (\ref{module d'un espace tordu}) que 
$\Delta_{G^\natural}(\delta_1)d_l\gamma = d_lg\cdot\delta_1$. 
Soit $(\Pi,V)$ une 
$\omega$--repr\'esentation lisse de $G^\natural$. 
Pour 
$\phi\in C^\infty_{\rm c}(G^\natural) $, on note 
$\Pi(\phi)=\Pi (\phi d_l\gamma)\in {\rm End}_\Bbb{C}(V)$ l'op\'erateur d\'efini (comme en 1.2) 
par 
$$
\Pi(\phi )(v) =\int_{G^\natural}\phi(\gamma)
\Pi(\gamma)(v) d_l\gamma \quad (v\in V).
$$
Supposons de plus que $\Pi$ est admissible. Alors cet op\'erateur $\Pi(\phi )$ est de rang fini. En effet, 
notant $C^\infty_{\rm c}(G)\rightarrow C^\infty_{\rm c}(G^\natural),\,
f\mapsto f^\natural$ l'isomorphisme de $\Bbb{C}$--espaces vectoriels d\'efini par 
$$
f^\natural(g\cdot \delta_1)=\Delta_{G^\natural}(\delta_1)^{-1}f(g)\quad 
(g\in G),
$$ et posant $\pi=\Pi^\circ$ et $A=\Pi(\delta_1)$, on a 
$$\Pi(f^\natural )= \pi(f)\circ A\quad (f\in C^\infty_{\rm c}(G)).
$$
On peut 
donc d\'efinir la trace de $\Pi(\phi )$:
$$
\Theta_\Pi(\phi)={\rm tr}(\Pi(\phi ))\quad (\phi\in C^\infty_{\rm c}(G^\natural)).
$$
La distribution $\Theta_\Pi$ sur $G^\natural$ ainsi d\'efinie, est appel\'ee le 
{\it caract\`ere de $\Pi$} (elle d\'epend bien s\^ur du choix de la mesure $d_l\gamma$). 
Pour $f\in C^\infty_{\rm c}(G)$, on a
$$
\Theta_\Pi(f^\natural)={\rm tr}(\pi(f)\circ A)=
 \Theta_\pi^A(f).
$$
Pour $\phi\in C^\infty_{\rm c}(G^\natural)$ et $x\in G$, notant 
${^x\phi}={\rm Int}_{G^\natural}(x)(\phi)\in C^\infty_{\rm c}(G^\natural)$ la fonction 
$\gamma\mapsto \phi(x^{-1}\cdot \gamma \cdot x)$, on a (d'aprs la relation $(**)$ de \ref{caractres tordus})
$$
\Theta_\Pi({^x\phi})=\Delta_G(x^{-1})\omega(x^{-1})\Theta_\Pi(\phi);\leqno{(*)}
$$
en utilisant (lemme de \ref{modules}) que 
$\Delta_G({\rm Int}_{G^\natural}(\delta_1)(x^{-1}))=\Delta_G(x^{-1})$.\goodbreak

\subsection{Induction compacte}\label{induction compacte} Soit $G^\natural$ un $G$--espace tordu, et 
$\omega$ un caract\`ere de $G$. Soit aussi $H$ un sous--groupe {\it ferm\'e} de $G$, et 
$H^\natural$ un $H$--espace tordu qui soit un sous--espace topologique tordu de $G^\natural$. 
Choisissons un \'el\'ement $\delta_1\in H^\natural$, et posons $\theta={\rm Int}_{G^\natural}(\delta_1)$. 
On a donc $\theta(H)=H$ et 
$\theta\vert_H={\rm Int}_{H^\natural}(\delta_1)$. On note encore $\omega$ 
le caract\`ere $\omega\vert_H$ de $H$, et $\theta$ l'automorphisme $\theta\vert_H$ de $H$. 

On d\'efinit comme suit un foncteur {\rm induction compacte (lisse)}
$$
{^\omega{\rm ind}}_{H^\natural}^{G^\natural}:\mathfrak{R}(H^\natural,\omega)\rightarrow \mathfrak{R}(G^\natural,\omega).
$$
Soit $(\Sigma,W)$ une $\omega$--repr\'esentation lisse de $H^\natural$. Posons $\sigma=\Sigma^\circ$ 
et $B=\Sigma(\delta_1)\in {\rm Isom}_H(\omega\sigma,\sigma^\theta)$. Notons 
$V^\natural={\rm ind}_H^{G^\natural}(W)$ 
l'espace des fonctions $\boldsymbol{\ES{F}}:G^\natural\rightarrow W$ telles que
\begin{itemize}
\item $\boldsymbol{\ES{F}}(h\cdot\gamma)=\sigma(h)(\boldsymbol{\ES{F}}(\gamma))$ pour tout $(h,\gamma)\in H\times G^\natural$,
\item il existe un sous--groupe 
ouvert compact $K_\ES{F}$ de $G$ tel que $\boldsymbol{\ES{F}}(\gamma\cdot k)=\boldsymbol{\ES{F}}(\gamma)$ 
pour tout $\gamma\in 
G$ et tout $k\in K_{\boldsymbol{\ES{F}}}$,
\item le support de $\boldsymbol{\ES{F}}$ est compact modulo $H$.
\end{itemize}

\v1 Soit $(\pi,V)$ la repr\'esentation induite compacte 
(lisse, non normalis\'ee) ${\rm ind}_H^G(\sigma,W)$. On rappelle la d\'efinition: 
$V$ est l'espace ${\rm ind}_H^G(W)$ obtenu en rempla\c{c}ant $G^\natural$ par $G$ 
dans la d\'efinition de $V^\natural$, 
et $\pi={\rm ind}_H^G(\sigma)$ est la repr\'esentation lisse de $G$ sur $V$ 
donn\'ee par
$$\pi(g)(\ES{F})(x)=\ES{F}(xg)\quad (\ES{F}\in V;\;g,\,x\in G).
$$
Pour 
$\ES{F}\in V$, notons $\ES{F}^\natural:G^\natural\rightarrow W$ la fonction d\'efinie par
$$
\ES{F}^\natural (g\cdot\delta_1)= \omega(\theta^{-1}(g))B(\ES{F}(\theta^{-1}(g)))\quad (g\in G),
$$
i.e. par
$$
\ES{F}^\natural(\delta_1\cdot g)= \omega(g)B(\ES{F}(g))\quad (g\in G).
$$
Pour $\ES{F}\in V$ et 
$(h,g)\in H\times G$, 
on a
\begin{align*}
\ES{F}^\natural (hg\cdot\delta_1)&= \omega(\theta^{-1}(hg))B(\ES{F}(\theta^{-1}(hg)))\\
& = \omega(\theta^{-1}(g))B\circ \omega\sigma(\theta^{-1}(h))(\ES{F}(\theta^{-1}(g)))\\
&= \omega(\theta^{-1}(g))\sigma(h)\circ B(\ES{F}(\theta^{-1}(g)))\\
& = \sigma(h)(\ES{F}^\natural(g\cdot\delta_1)).
\end{align*}
On en d\'eduit que $\ES{F}^\natural\in V^\natural$, et puisque $B\in {\rm Aut}_\Bbb{C}(W)$, 
l'application
$$V\rightarrow V^\natural,\,\ES{F}\mapsto \ES{F}^\natural
$$
est un isomorphisme de $\Bbb{C}$--espaces vectoriels. Cet isomorphisme 
ne d\'epend pas du choix de $\delta_1$: d'apr\`es le calcul ci--dessus, pour $\ES{F}\in V$, 
on a 
$$
\ES{F}^\natural(\delta\cdot g)= \omega(g)\Sigma(\delta)(\ES{F}(g))\quad 
(\delta\in H^\natural,\,g\in G).
\leqno{(1)}
$$
Pour $\gamma\in G^\natural$, soit $\Pi(\gamma)\in {\rm Aut}_\Bbb{C}(V)$ 
l'op\'erateur d\'efini par
$$
\Pi(\gamma)(\ES{F})(g)= \ES{F}^\natural(g\cdot \gamma)\quad (\ES{F}\in V).\leqno{(2)}
$$
Pour $\ES{F}\in V$, $\gamma=z\cdot\delta_1\in G^\natural$ et $x,\,y,\,g\in G$, on a
\begin{align*}
\Pi(x\cdot \gamma \cdot y)(\ES{F})(g) &= \ES{F}^\natural(gxz\theta(y)\cdot \delta_1)\\
& = \omega(\theta^{-1}(gxz)y)B(\ES{F}(\theta^{-1}(gxz)y)))\\
& = \omega(y) \omega(\theta^{-1}(gxz))B(\pi(y)(\ES{F})(\theta^{-1}(gxz)))\\
& = \omega(y) [\pi(y)(\ES{F})]^\natural(gx\cdot \gamma)\\
& = \omega(y)\Pi(\gamma)\circ \pi(y)(\ES{F})(gx)\\
&= \omega(y)\pi(x)\circ \Pi(\gamma)\circ \pi(y)(\ES{F})(g).
\end{align*}
L'application $G^\natural\rightarrow {\rm Aut}_\Bbb{C}(V),\,\gamma\mapsto \Pi(\gamma)$ 
est donc une $\omega$--repr\'esentation lisse de $G^\natural$, telle que $\Pi^\circ=\pi$. On pose
$$
{^\omega{\rm ind}}_{H^\natural}^{G^\natural}(\Sigma,W)=(\Pi,V).
$$
Notons que l'op\'erateur $A=\Pi(\delta_1)\in {\rm Isom}_G(\omega\pi,\pi^\theta)$ est donn\'e par
$$
A(\ES{F})(g)=\ES{F}^\natural(g\cdot \delta_1)\quad (\ES{F}\in V,\,g\in G).
$$
D'apr\`es les relations (1) et (2), pour $\ES{F}\in V$, $\delta\in H^\natural$ et $g,\,x\in G$, on a
$$
\Pi(g\cdot \delta)(\ES{F})(x)=\omega({\rm Int}_{H^\natural}(\delta)^{-1}(xg)) 
\Sigma(\delta)(\ES{F}({\rm Int}_{H^\natural}(\delta)^{-1}(xg));
$$
en particulier pour $x=1$, on a
$$
\Pi(\delta\cdot g)(\ES{F})(1)= \omega(g)\Sigma(\delta)(\ES{F}(g))
$$

Par construction, les foncteurs ${^\omega{\rm ind}}_{H^\natural}^{G^\natural}:\mathfrak{R}(H^\natural,\omega)\rightarrow \mathfrak{R}(G^\natural,\omega)$ et ${\rm ind}_H^G:\mathfrak{R}(H)\rightarrow \mathfrak{R}(G)$ commutent aux foncteurs d'oubli, au sens o\`u pour toute $\omega$--repr\'esentation lisse 
$\Sigma$ de $H^\natural$, on a
$$
{^\omega{\rm ind}}_{H^\natural}^{G^\natural}(\Sigma)^\circ = {\rm ind}_H^G(\Sigma^\circ).
$$
Notons que si $G$ est compact modulo $H$, la condition sur le support de $\ES{F}$ est automatiquement vrifie, 
et si de plus $\Sigma$ est admissible, alors $\Pi$ l'est aussi \cite[2.26]{BZ}.
\goodbreak

\subsection{Caract\`eres des induites compactes}\label{caractres des induites compactes}
Continuons avec les notations de \ref{induction compacte}. 
On fixe une mesure de Haar 
\`a gauche $d_lh$ sur $H$, et l'on note $d_rh$ la mesure de Haar \`a droite 
$\Delta_H(h^{-1})d_lh$ sur $H$. Pour toute repr\'esentation lisse $\sigma$ de $H$, on pose 
$\sigma(f)=\sigma(f d_lh)$ ($f\in C^\infty_{\rm c}(H)$), et si $\sigma$ est admissible, on note $\Theta_\sigma$ la distribution 
sur $H$ d\'efinie par $d_lh$ comme en 1.2.

Soit $d_l\gamma=\delta_1\cdot d_lg$ la mesure 
de Haar \`a gauche sur $G^\natural$ associ\'ee \`a $d_lg$, et soit $d_l\delta=
\delta_1\cdot d_lh$ la 
mesure de Haar \`a gauche sur $H^\natural$ associ\'ee \`a $d_lh$. 
Pour toute $\omega$--repr\'esentation 
lisse $\Pi$ de $G$, on pose $\Pi(\phi)=\pi(\phi d_l\gamma)$ 
($\phi\in C^\infty_{\rm c}(G^\natural)$), et si $\Pi$ est admissible, on note 
$\Theta_\Pi$ la distribution sur $G^\natural$ d\'efinie par $d_l\gamma$ 
comme en 1.6. De m\^eme, pour toute $\omega$--repr\'esentation lisse $\Sigma$ de $H^\natural$, 
on pose $\Sigma(\phi)=\Sigma(\phi d_l\delta)$ ($\phi\in C^\infty_{\rm c}(H^\natural)$), et si $\Sigma$ est admissible, on note $\Theta_\Sigma$ la 
distribution sur $H^\natural$ d\'efinie par $d_l\delta$.

Pour toute partie ouverte compacte $\Omega$ de $G$ telle que $H\cap \Omega\neq \emptyset$, on note
$$C^\infty_{\rm c}(G^\natural)\rightarrow C^\infty_{\rm c}(G^\natural),\phi\mapsto
\phi_\Omega = {^\omega\!\phi_\Omega}
$$ l'application 
lin\'eaire d\'efinie par
$$
\phi_\Omega(\gamma) =
{\rm vol}(H\cap \Omega,
d_rh)^{-1}
\int_\Omega\omega(g)\phi(g^{-1}\cdot \gamma \cdot g)d_lg
\quad (\gamma\in G^\natural).
$$
Notons que si $H\cap \Omega$ est un groupe (par exemple si $\Omega$ est un sous--groupe de $G$), 
alors puisque ce groupe est compact, on a 
$\Delta_H\vert_{H\cap \Omega}=1$ et 
${\rm vol}(H\cap \Omega,
d_rh)={\rm vol}(H\cap \Omega,
d_lh)$.

\v1 Consid\'erons les deux conditions suivantes:
\begin{enumerate}
\item[$(i)$]Le groupe $G$ est compact modulo $H$.
\item[$(ii)$]Il existe un sous--groupe ouvert compact $K$ de $G$ 
tel que $HK=KH$.
\end{enumerate}

\ni Si $G$ est compact modulo $H$, et si $K$ est un sous--groupe 
ouvert compact de $G$, alors il existe une partie ouverte compacte $\Omega$ de 
$G$ telle que $G=H\Omega$ et 
$K\Omega =\Omega$. Puisque $G= H\Omega$, l'ensemble $H\cap \Omega$ n'est pas vide. 
Si de plus $HK=KH$, alors l'ensemble $HK$ est un sous--groupe ouvert 
de $G$, et comme $G=H\Omega$, il est d'indice fini dans $G$. 
Choisissons 
un syst\`eme de repr\'esentants $\{x_1,\ldots ,x_n\}$ 
dans $G$ de l'espace quotient $HK\backslash G$, et posons 
$\Omega'=\bigcup_{i=1}^n Kx_i$ (l'union est disjointe). On peut supposer 
que $1\in \{x_1,\ldots ,x_n\}$. 
Alors $\Omega'$ est une partie ouverte compacte de $G$ v\'erifiant:
\begin{enumerate}
\item[$(iii)$] 
$G=H\Omega'$, $K\Omega'=\Omega'$ et $H \cap \Omega'= H\cap K$. 
\end{enumerate}

\begin{mapropo}On suppose que les conditions $(i)$ et $(ii)$ 
sont v\'erifi\'ees. Choisissons un sous--groupe ouvert compact $K$ de $G$ tel que 
$HK=KH$, et un système de représentants $\{x_1=1,x_2,\ldots ,x_n\}$ 
dans $G$ de l'espace quotient $HK\backslash G$. 
Posons $\Omega'=\bigcup_{i=1}^nKx_i$. 
Soit $\Sigma$ une 
$\omega$--repr\'esentation admissible de $H^\natural$, et soit 
$\Pi={^\omega{\rm ind}}_{H^\natural}^{G^\natural}(\Sigma)$. 
Pour toute fonction $\phi\in C^\infty_{\rm c}(G^\natural)$, 
on a la formule de 
descente
$$
\Theta_\Pi(\phi)=\Theta_\Sigma(\phi_{\Omega'}\vert_{H^\natural}).
$$
\end{mapropo}

\begin{proof} La d\'emonstration est longue et laborieuse, mais son principe est simple: 
fix\'ee la fonction $\phi\in C^\infty_{\rm c}(G^\natural)$, on construit un $\Bbb{C}$--espace vectoriel 
de dimension finie ${\boldsymbol{X}}$ et un op\'erateur $\boldsymbol{T}\in {\rm End}_\Bbb{C}({\boldsymbol{X}})$, 
tels que $\Theta_\Pi(\phi)={\rm tr}({\boldsymbol{T}})$; puis on calcule la trace de ${\boldsymbol{T}}$. Soit $W$ l'espace de $\Sigma$, et soit $V={\rm ind}_H^G(W)$ l'espace de $\Pi$. 
Posons $\sigma= \Sigma^\circ$ et $B= \Pi(\delta_1)$, $\pi=\Pi^\circ$ et $A=\Pi(\delta_1)$.

Par définition, l'union $\Omega'=\bigcup_{i=1}^nKx_i$ est disjointe, et $\Omega'$ est une partie ouverte compacte de $G$ v\'erifiant 
$(iii)$. Posons 
$K'=H\cap \Omega'$, et 
notons $\mathfrak{X}'$ le $\Bbb{C}$--espace vectoriel 
des fonctions localement constantes 
$\Omega'\rightarrow W$. On a une identification canonique 
$\mathfrak{X}'=C^\infty(\Omega')\otimes_\Bbb{C}W$. Puisque $G=H\Omega'$, l'application $V\rightarrow \mathfrak{X}',\,\ES{F}\mapsto \ES{F}\vert_{\Omega'}$ 
identifie $V$ au sous--espace $\mathfrak{X}$ de $\mathfrak{X}'$ 
form\'e des fonctions $\ES{F}'$ telles que $\ES{F}'(hx)=\sigma(h)(\ES{F}'(x))$ 
pour tout $(x,h)\in \Omega'\times  K'$. Soit $(\eta,C^\infty(\Omega'))$ la repr\'esentation 
admissible de $K$ d\'efinie par
$$
\eta(k)(\varphi)(x)=\varphi (k^{-1}x)\quad (k\in K,\,\varphi\in C^\infty(\Omega'),\,x\in \Omega').
$$
Posons $c={\rm vol}(K'\!,d_lh)$. 
Le $\Bbb{C}$--endomorphisme $u$ de $\mathfrak{X}'$ d\'efini par
$$
u(\varphi \otimes w)=c^{-1}\int_{H\cap K}\eta(h)(\varphi)\otimes \sigma(h)(w)d_lh
\quad (\varphi\in C^\infty(\Omega'),\,x\in \Omega') 
$$
est un projecteur sur $\mathfrak{X}$; i.e. on a 
$u(\mathfrak{X}')=\mathfrak{X}$ et $u\vert_\mathfrak{X}={\rm id}$.

Soit une fonction $\phi\in C^\infty_{\rm c}(G^\natural)$. Choisissons deux sous--groupes ouverts 
distingu\'es $K_1$ et $K_2$ de $K$ tels que:
\begin{itemize}
\item $\phi\in C_{\rm c}(K_1\backslash G^\natural /K_1)$ et $\omega\vert_{K_1}=1$;
\item $K_2\subset K_1$ et $x_i^{-1}K_2x_i \subset K_1\subset x_i^{-1}Kx_i$ 
($i=1,\ldots,n$).
\end{itemize}
Soit $e_{K^1}\in C^\infty_{\rm c}(G)$ la fonction ${\rm vol}(K_1,d_lg)^{-1}{\bf 1}_{K_1}$ o\`u 
${\bf 1}_{K_1}$ d\'esigne la fonction caract\'eristique de $K_1$. 
Puisque $\phi\in C_{\rm c}(K_1\backslash G^\natural /K_1)$ et 
$\omega\vert_{K_1}=1$, on a
$$
\pi(e_{K_1})\circ \Pi(\phi )\circ \pi(e_{K_1})= \Pi(\phi ).
$$
En particulier, on a
$$\Pi(\phi )(V)\subset \pi(e_{K_1})(V)=V^{K_1};
$$
et posant $V(K_1)=\langle \ES{F}-\pi(k)(\ES{F}):\ES{F}\in V,\,k\in K_1\rangle$, on a
$$
\Pi(\phi )(V(K_1))=0.
$$
D'autre part, puisque $x_i^{-1}K_2x_i \subset K_1\subset x_i^{-1}Kx_i$ 
($i=1,\ldots,n$), on a
$$K_2\Omega'K_1=\Omega'.
$$
Posons $K'_2=H\cap K_2$, et ${\boldsymbol{X}}=C(K_2\backslash \Omega')\otimes W^{K'_2}$; c'est un 
sous--espace vectoriel de dimension finie de $\mathfrak{X}'$. Alors via l'identification $V=\mathfrak{X}$, on a l'inclusion 
$V^{K_1}\subset {\boldsymbol{X}}$. Notons $T'$ 
l'op\'erateur $\Pi(\phi d_l\gamma)\circ u\in {\rm End}_\Bbb{C}(\mathfrak{X}')$. Il est de rang fini puisque $T'(\mathfrak{X}')\subset {\boldsymbol{X}}$, et l'on a
$$\Theta_\Pi(\phi)={\rm tr}(T').
$$ 
Posons ${\boldsymbol{T}}= T'\vert_{{\boldsymbol{X}}}\in {\rm End}_\Bbb{C}({\boldsymbol{X}})$. Par définition de la trace d'un opérateur de rang fini, on a ${\rm tr}(T')={\rm tr}(\boldsymbol{T})$. D'où le

\begin{monlem1}
On a $\Theta_\Pi(\phi)={\rm tr}({\boldsymbol{T}})$.
\end{monlem1}

\begin{marema1}
{\rm 
Notons $W(K'_2)$ le sous--espace vectoriel de $W$ engendr\'e par les lments 
$w- \sigma(k)(w)$ pour $w\in W$ et $k\in K'_2$. Alors on a la d\'ecomposition 
$W=W^{K'_2}\oplus W(K'_2)$. De m\^eme, en rempla\c{c}ant $(\sigma,W)$ par la reprsentation 
$(\eta,C^\infty(\Omega'))$, on a la d\'ecomposition 
$C^\infty(\Omega')=C^\infty(K_2\backslash \Omega')\oplus C^\infty(\Omega')(K_2)$. On a aussi la décomposition 
$V=V^{K_1}\oplus V(K_1)$. Par définition des groupes $K_1$, $K_2$ et $K'_2$, on a
$$
u(C(K_2\backslash \Omega')\otimes W(K'_2))=0,\quad u(C^\infty(\Omega')(K_2)\otimes W)\subset V(K_1).
$$ 
En effet, pour $\varphi\in C(K_2\backslash \Omega')$, 
$w\in W$ et $k'\in K'_2$, puisque $\omega(k')(\varphi)=\varphi$ et 
$\Delta_H\vert_{K'_2}=1$, on a
\begin{align*}
u(\varphi\otimes \sigma(k')(w))&= c^{-1}\int_{K'}\eta(h)(\varphi)\otimes \sigma(hk')(w)d_lh\\
&= c^{-1}\int_{K'}\eta(hk')(\varphi)\otimes \sigma(hk')(w)d_lh\\
&= u(\varphi\otimes w).
\end{align*}
D'où l'égalité $u(C(K_2\backslash \Omega')\otimes W(K'_2))=0$. D'autre part, via l'identification $V=\mathfrak{X}$, 
$V(K_1)$ 
s'identifie 
au sous--espace vectoriel de $\mathfrak{X}$ form\'e des fonctions 
$\ES{F}$ telles que $\int_{K_1}\ES{F}(xg)d_lg=0$ pour tout $x\in \Omega'$. Or 
pour $\varphi\in C^\infty(\Omega')$, $k\in K_2$, $w\in W$ et $x\in \Omega'$, on a
\begin{align*}
\int_{K_1}u(\eta(k)\varphi\otimes w)(xg)d_lg &= \int\!\!\!\!\int_{K_1\times K'} \eta(hk)(\varphi)(xg)\otimes \sigma(h)(w)d_lgd_lh\\
& = \int\!\!\!\!\int_{K_1\times K'}\varphi(k^{-1}h^{-1}xg)\otimes \sigma(h)(w)d_lgd_lh.
\end{align*}
Mais pour $h\in K'$, on a 
$k^{-1}h^{-1}= h^{-1}k'^{-1}$ avec $k'=hkh^{-1}\in K_2$, et l'on a $k'^{-1}x=xk''^{-1}$ 
avec $k''=x^{-1}k'x\in K_1$; on peut donc effectuer le changement de variables 
$g\mapsto k''g$. En d\'efinitive, on obtient l'\'egalit\'e
$$
\int_{K_1}u(\eta(k)\varphi\otimes w)(xg)d_lg=
\int_{K_1}u(\varphi\otimes w)(xg)d_lg.
$$
D'où l'inclusion $u(C^\infty(\Omega')(K_2)\otimes W)\subset V(K_1)$. \hfill $\blacksquare$
}
\end{marema1}

Calculons la trace de l'opérateur ${\boldsymbol{T}}$. Notons $(\check{\sigma},\check{W})$ la reprsentation de $H$ contragrdiente 
de $\sigma$. L'application canonique 
$W\times \check{W}\rightarrow \Bbb{C},\,(w,\check{w})\mapsto 
\langle w,\check{w}\rangle$ induit par restriction une 
identification $\check{W}^{K'_2}={\rm Hom}_\Bbb{C}(W^{K'_2},\Bbb{C})$. 
Posons $\widetilde{\boldsymbol{X}}=C(K_2\backslash \Omega')\otimes \check{W}^{K'_2}$, 
et notons
$${\boldsymbol{X}}\times \widetilde{\boldsymbol{X}}\rightarrow \Bbb{C},\,(\ES{F}, \widetilde{\ES{F}})\mapsto 
\langle \ES{F}, \widetilde{\ES{F}}\rangle$$
l'application bilin\'eaire non d\'eg\'en\'er\'ee 
d\'efinie par
$$
\langle \ES{F}, \widetilde{\ES{F}}\rangle
=\int_{\Omega'}\langle \ES{F}(g), \widetilde{\ES{F}}(g)\rangle d_lg
= \sum_{i=1}^n \Delta_G(x_i)\int_K\langle \ES{F}(gx_i), \widetilde{\ES{F}}(gx_i)\rangle d_lg.
$$
Elle induit une identification 
$\widetilde{\boldsymbol{X}}={\rm Hom}_\Bbb{C}({\boldsymbol{X}},\Bbb{C})$.

\begin{monlem2}
On a ${\rm tr}(\boldsymbol{T})=\Theta_\Sigma(\phi_{\Omega'}\vert_{H^\natural})$.
\end{monlem2}

\begin{proof} 
Fixons une base $\mathfrak{B}$ de $W^{K'_2}$ sur 
$\Bbb{C}$, et notons $\{\check{w}\}_{w\in\mathfrak{B}}$ la base de 
$\check{W}^{K'_2}$ 
duale de $\mathfrak{B}$. Choisissons un syst\`eme de 
repr\'esentants $\{k_1,\ldots ,k_s\}$ dans $K$ des \'el\'ements du groupe 
quotient $K_2\backslash K$, et pour $j=1,\ldots ,s$, posons 
$\varphi^{k_j}={\bf 1}_{K_2k_j}\in C(K_2\backslash 
K)$. Alors 
$\mathfrak{C}_\circ=\{\varphi^{k_j}\}_{j=1,\ldots ,s}$ est une base 
de $C(K_2\backslash K)$ sur $\Bbb{C}$. Pour $\varphi\in \mathfrak{C}_\circ$ et 
$i\in\{1,\ldots ,n\}$, notons $\varphi_i \in C(K_2\backslash\Omega')$ la fonction 
d\'efinie par ${\rm Supp}(\varphi_i)\subset Kx_i$ et $\varphi_i(kx_i)= 
\varphi(k)$ ($k\in K$). 
Alors $\mathfrak{C}=\{\varphi_i\}_{\varphi\in\mathfrak{C}_\circ,\,
i=1,\ldots ,n}$ est une base de $C(K_2\backslash\Omega')$ sur $\Bbb{C}$, et 
$\mathfrak{D}=\{\psi\otimes w\}_{ \psi\in\mathfrak{C},\,w\in\mathfrak{B}}$ est une base de 
${\boldsymbol{X}}$ sur $\Bbb{C}$. Posons $d={\rm vol}(K_2,d_lg)$. 
Puisque $\Delta_G\vert_K=1$, 
$\{d^{-1}(\Delta_G^{-1}\psi)\otimes \check{w})\}_{\psi\in \mathfrak{C},\,w\in\mathfrak{B}}$ 
est la base de $\widetilde{\boldsymbol{X}}$ duale de $\mathfrak{D}$. On a donc
\begin{align*}
{\rm tr}(\boldsymbol{T})&= d^{-1}\sum_{w\in\mathfrak{B}}\sum_{\psi\in\mathfrak{C}}
\langle \boldsymbol{T}(\psi\otimes w), (\Delta_G^{-1}\psi)\otimes \check{w}\rangle\\
&= d^{-1}\sum_{w\in\mathfrak{B}}\sum_{\psi\in\mathfrak{C}}\sum_{i=1}^n\int_K
 \psi(kx_i)\langle \boldsymbol{T}(\psi\otimes w)(kx_i),\check{w}\rangle dk\\
 &= d^{-1}\sum_{i=1}^n\sum_{w\in\mathfrak{B}}\sum_{\varphi\in\mathfrak{C}_\circ}
\int_K
 \varphi(k)\langle \boldsymbol{T}(\varphi_i\otimes w)(kx_i),\check{w}\rangle dk,
\end{align*}
o\`u l'on a pos\'e $dk=d_lg\vert_K$. 
Soit $\psi\in C(K_2\backslash \Omega')$, $w\in W$ et $x\in \Omega'$. Posons 
$\ES{F}=u(\psi\otimes w)\in V$ et notons $\ES{F}^\natural\in V^\natural$ la fonction 
$G^\natural\rightarrow W,\,g\cdot \delta_1\mapsto \omega(\theta^{-1}(g))B(\ES{F}(\theta^{-1}(g))$ 
--- cf. \ref{induction compacte}. 
Alors on a
\begin{align*}
\boldsymbol{T}(\psi\otimes w)(x)&=\int_{G^\natural}\phi(\gamma)\Pi(\gamma)(\ES{F})(x)d_l\gamma\\
&= \Delta_G(\theta)^{-1}\int_G\phi(g\cdot \delta_1) A(\ES{F})(xg)d_lg\\
&= \Delta_G(\theta)^{-1}\int_{G}\phi(x^{-1}g\cdot \delta_1)A(\ES{F})(g)d_lg\\
& = \int_{G^\natural}\phi(x^{-1}\cdot \gamma)\ES{F}^\natural(\gamma)d_l\gamma.
\end{align*}
Puisque $G=\coprod_{i=1}^nHKx_i$, on a les d\'ecompositions 
$$
G=\coprod_{i=1}^nH\theta(Kx_i),\quad 
G^\natural= \coprod_{i=1}^n H^\natural\cdot Kx_i=\coprod_{i=1}^n \delta_1\cdot HKx_i.
$$ D'apr\`es le lemme de \ref{module d'un espace tordu}, pour 
$i=1,\ldots ,n$, on a $d_l(\gamma\cdot x_i)=\Delta_G(x_i)d_l(\gamma)$. Comme $d_l\gamma=\delta_1\cdot d_lg$, 
$d_l\delta= \delta_1\cdot d_lh$, $H\cap K= K'$ et ${\rm vol}(K',d_lh)=c$, posant 
$\Delta_i= c^{-1}\Delta_G(x_i)$, on obtient
$$
\boldsymbol{T}(\psi\otimes w)(x)=\sum_{i=1}^n\Delta_i\int\!\!\!\!\int_{H^\natural\times K}
\phi(x^{-1}\cdot \delta\cdot kx_i)\ES{F}^\natural(\delta\cdot kx_i)d_l\delta dk.$$
Or pour $\delta\in H^\natural$, on a (relation (1) de \ref{induction compacte}) 
$\ES{F}^\natural(\delta\cdot kx_i)=\omega(kx_i)\Sigma(\delta)(\ES{F}(kx_i))$ et
\begin{align*}
\Sigma(\delta)(\ES{F}(kx_i)) &= c^{-1}\int_{K'}\omega(kx_i)\psi(k'^{-1}kx_i) \Sigma(\delta)\circ \sigma(k')(w)dk'\\
& = c^{-1}\int_{K'}\omega\psi(k'^{-1}kx_i) \Sigma(\delta\cdot k')(w)dk',
\end{align*}
o\`u l'on a pos\'e 
$dk'=d_lh\vert_{K'}$. 
En effectuant 
les changements de variables $\delta\mapsto \delta\cdot k'^{-1}$ (notons que  
$d_l(\delta\cdot k'^{-1})=\Delta_H(k'^{-1})d_l\delta=d_l\delta$) puis $k\mapsto k'k$ dans la formule pour $\boldsymbol{T}(\psi\otimes w)(x)$, on 
obtient
\begin{align*}
\lefteqn{\boldsymbol{T}(\psi\otimes w)(x)=}\\
&= c^{-1}\sum_{i=1}^n
\Delta_i\int\!\!\!\!\int\!\!\!\!\int_{H^\natural\times K\times K'}
\phi(x^{-1}\cdot \delta\cdot kx_i)\omega\psi(k'^{-1}kx_i)\Sigma(\delta\cdot k')(w)d_l\delta dk dk'\\
&=\sum_{i=1}^n
\Delta_i\int\!\!\!\!\int_{H^\natural\times K}
\phi(x^{-1}\cdot \delta\cdot kx_i)\omega\psi(kx_i)\Sigma(\delta)(w)d_l\delta dk.
\end{align*}
En particulier, si $\psi=\varphi_i$ et $x=k_0x_i$ pour un $\varphi\in \mathfrak{C}_\circ$, 
un $i\in\{1,\ldots ,n\}$ et un $k_0\in K$, alors on a
$$
\boldsymbol{T}(\varphi_i\otimes w)(k_0x_i)=
\omega(x_i)\Delta_i\int\!\!\!\!\int_{H^\natural\times K}
\phi(x_i^{-1}k_0^{-1}\cdot \delta\cdot kx_i)\omega\varphi(k)\Sigma(\delta)(w)d_l\delta dk.
$$
Injectons l'\'egalit\'e ci--dessus dans la formule pour ${\rm tr}(\boldsymbol{T})$. On obtient
$$
{\rm tr}(\boldsymbol{T})  = d^{-1}\sum_{w\in\mathfrak{B}}
\int_H \Phi(\delta)\langle \Sigma(\delta)(w),\check{w}\rangle d_l\delta
$$
o\`u
$$
\Phi(\delta)=\sum_{i=1}^n\omega(x_i)\Delta_i
\sum_{\varphi\in\mathfrak{C}_\circ}\int\!\!\!\!\int_{K\times K}\varphi(k_0)
\varphi(k)\omega(k)\phi(x_i^{-1}k_0^{-1}\cdot \delta \cdot kx_i)dk_0dk.
$$
Pour $i=1,\ldots ,n$ et 
$\delta\in H^\natural$, puisque $x_i^{-1}K_2x_i\subset K_1$, 
$\phi\in C_{\rm c}(K_1\backslash G^\natural / K_1)$ et $\omega\vert_{K_2}=1$, la fonction 
$$K\times K\rightarrow \Bbb{C},\, (k_0,k)\mapsto 
\omega(k)\phi(x_i^{-1}k_0^{-1}\cdot \delta \cdot kx_i)
$$ se 
factorise \`a travers $K_2\backslash K\times K_2\backslash K$. Par d\'efinition de la 
base $\mathfrak{C}_\circ$, on a
\begin{align*}
\Phi(\delta)& = d^2\sum_{i=1}^n\omega(x_i)\Delta_i
\sum_{k\in K_2\backslash K}\omega(k)f(x_i^{-1}k^{-1}\cdot \delta\cdot kx_i)\\
&= dc^{-1}\sum_{i=1}^n\omega(x_i)\Delta_G(x_i) \int_K
\omega(k)\phi(x_i^{-1}k^{-1})\cdot \delta \cdot kx_i)dk\\
&= dc^{-1}\int_{\Omega'}\omega(g)\phi(g^{-1}\cdot \delta \cdot g)d_lg\\
&= d{{\rm vol}(K'\!,d_rh)\over 
{\rm vol}(K'\!,d_lh)}\phi_{\Omega'}(\delta).
\end{align*}
Or $K'=H\cap K$ est un groupe compact, par cons\'equent 
$d_rh\vert_{K'}= d_lh\vert_{K'}$ et 
$\Phi (\delta)= d\phi_{\Omega'}(\delta)$. On a donc 
\begin{align*}
{\rm tr}(\boldsymbol{T}) & = \sum_{w\in\mathfrak{B}}
\int_H \phi_{\Omega'}(\delta)\langle \Sigma(\delta)(w),\check{w}\rangle d_l\delta \\
&= \sum_{w\in\mathfrak{B}}
\langle \Sigma(\phi_{\Omega'}\vert_{H^\natural}d_l\delta)(w),\check{w}\rangle.
\end{align*}
Posons $e_{K'_2}={\rm vol}(K'_2,d_lh)^{-1}{\bf 1}_{K'_2}$. Puisque $\phi_{\Omega'}\vert_{H^\natural}\in C_{\rm c}(K'_2\backslash H^\natural /K'_2)$ et 
$\omega\vert_{K'_2}=1$, on a
$$
\sigma(e_{K'_2})\circ \Sigma(\phi_{\Omega'}\vert_{H^\natural})\circ \sigma(e_{K'_2})=
\Sigma(\phi_{\Omega'}\vert_{H^\natural}).
$$
En particulier, $\Sigma(\phi_{\Omega'}\vert_{H^\natural})(W)\subset W^{K'_2}$ et 
$\Sigma(\phi_{\Omega'}\vert_{H^\natural})(W(K'_2))=0$, par cons\'equent  
$\Theta_\Sigma(\phi_{\Omega'}\vert_{H^\natural})$ co\"{\i}ncide avec la trace de l'op\'erateur 
$\Sigma(\phi_{\Omega'}\vert_{H^\natural})$ sur $W^{K'_2}$. 
On a donc bien
$$
{\rm tr}(\boldsymbol{T}) =\Theta_\Sigma(\phi_{\Omega'}\vert_{H^\natural}),
$$
ce qu'il fallait dmontrer.\end{proof}

Cela ach\`eve la d\'emonstration de la proposition.
\end{proof}

\begin{moncoro}On suppose qu'il existe un 
sous--groupe ouvert compact $K$ de $G$ tel que $G=HK$. 
Soit $\Sigma$ une 
$\omega$--repr\'esentation admissible de $H^\natural$, et soit 
$\Pi={^\omega{\rm ind}}_{H^\natural}^{G^\natural}(\Sigma)$. 
Pour toute fonction $\phi\in C^\infty_{\rm c}(G^\natural)$, 
on a la formule de 
descente
$$
\Theta_\Pi(\phi)=\Theta_\Sigma(\phi_K\vert_{H^\natural}).
$$
\end{moncoro}

\begin{marema2}{\rm On suppose que les conditions $(i)$ et $(ii)$ 
sont v\'erifi\'ees. 
Soit $K$ un sous--groupe ouvert compact de $G$ tel que $HK=KH$. 
Choisissons un syst\`eme de repr\'esentants $\{x_1=1,x_2,\ldots ,x_n\}$ dans 
$G$ de l'espace quotient $HK\backslash G$, 
et posons $\Omega'=\bigcup_{i=1}^nKx_i$ (l'union est disjointe). Alors $\Omega'$ vérifie $(iii)$, et pour 
$\phi\in C^\infty_{\rm c}(G)$ et $\gamma\in G^\natural$, on a
\begin{align*}
\phi_{\Omega'}(\gamma)&= {\rm vol}(H\cap K,d_lh)^{-1}
\int_{\Omega'}\omega(g)\phi (g^{-1}\cdot \gamma \cdot g)d_lg\\
&= \sum_{i=1}^n\omega(x_i)\Delta_G(x_i)({^{x_i}\phi})_K(\gamma).
\end{align*}
Notons que si le groupe $HK$ est $\theta$--stable (ce qui n'implique pas que $K$ soit 
$\theta$--stable), alors on peut choisir les $x_i$ de telle mani\`ere que l'ensemble 
$\{x_1=1,\ldots ,x_n\}$ soit $\theta$--stable. En effet, le groupe 
$\langle \theta \rangle  = \langle \theta^i: i\in \Bbb{Z}\rangle$ op\`ere sur 
l'espace quotient $HK\backslash G$, et pour chaque 
$\langle \theta \rangle$--orbite $\ES{O}$ dans $HK\backslash G$, il existe 
un $x\in G$ tel que $\ES{O}=\bigcup_{i\in\Bbb{Z}}HK\theta^i(x)$ (l'union est finie).
\hfill $\blacksquare$}
\end{marema2}

\begin{marema3}
{\rm Continuons avec les notations de la remarque 2. On peut, dans la proposition, remplacer $\Omega'$ par 
n'importe quelle partie $\Omega$ de $G$ de la forme $\Omega = \bigcup_{j=1}^my_j \Omega'$ pour des éléments 
$y_1,\ldots ,y_m\in H$: pour toute fonction $\phi\in C^\infty_{\rm c}(G^\natural)$, on a encore
$$
\Theta_\Pi(\phi)= \Theta_{\Sigma}(\phi_\Omega\vert_{H^\natural}).
$$
En effet, on peut supposer que pour $j,\, k\in \{1,\ldots, m\}$ tels que $j\neq k$, on a $y_j\Omega'\neq y_k\Omega'$. Alors l'union $\bigcup_{j=1}^m y_j\Omega'$ est disjointe: si $y_j\Omega' \cap y_{k}\Omega'\neq \emptyset$ pour des entiers $j,\,k\in \{1,\ldots ,m\}$, alors $y_jKx_i\cap y_{k}Kx_{l}$ pour des entiers $i,\,l\in \{1,\ldots ,n\}$. On a donc $i=l$
et $y_jK= y_kK$, d'où $y_j\Omega'= y_k\Omega'$ puisque $K\Omega'=\Omega'$, et $j=k$. Posons $K'= H\cap \Omega'\;(= H\cap K)$. 
Alors pour $\delta\in H^\natural$, on a
\begin{align*}
\phi_\Omega(h)& ={\rm vol}(H\cap \Omega,d_rh)^{-1}\sum_{j=1}^m
\int_{y_j\Omega'}\omega(g)\phi (g^{-1}\cdot \delta \cdot g)d_lg\\
& = {{\rm vol}(K'\!,d_rh)\over {\rm vol}(H\cap \Omega,d_rh)}\sum_{j=1}^m\omega(y_j)\phi_{\Omega'}
(y_j^{-1}\cdot\delta\cdot y_j).
\end{align*}
Or pour $\xi\in C^\infty_{\rm c}(H^\natural)$ et $y\in H$, on a (relation $(*)$ de \ref{caractres tordus (bis)})
$$\Theta_\Sigma({^y\xi})=
\Delta_H(y^{-1})\omega(y^{-1})\Theta_\Sigma(\xi).
$$
Par suite, posant 
$c={{\rm vol}(K'\!,d_rh)\over {\rm vol}(H\cap \Omega,d_rh)}
\sum_{j=1}^m\Delta_H(y_j^{-1})$, on obtient
$$
\Theta_\Sigma(\phi_\Omega\vert_{H^\natural})=
c\,\Theta_\Sigma(\phi_{\Omega'}\vert_{H^\natural}).
$$
Mais comme
$$
H\cap \Omega = \coprod_{j=1}^m H\cap y_j\Omega'= \coprod_{j=1}^m y_jK',
$$
on a
$$
{\rm vol}(H\cap \Omega,d_rh)=\sum_{j=1}^m {\rm vol}(y_jK'\!,d_rh)
=\sum_{j=1}^m \Delta_G(y_j^{-1}){\rm vol}(K'\!,d_lh).
$$
Donc $c=1$.\hfill$\blacksquare$
}
\end{marema3}

\subsection{Commentaires} La proposition et son corollaire {\it ne 
d\'ependent pas} du choix du point--base $\delta_1\in G^\natural$ --- rappelons que l'approche 
classique consiste 
\`a fixer l'automorphisme $\theta$ de $G$ et \`a \'etudier les caract\`eres $\theta$--tordus, 
ou plus g\'en\'eralement $(\theta,\omega)$--tordus, de $G$ --- cf. \ref{caractres tordus}. Remplacer $\delta_1$ par 
$\delta'_1=x\cdot \delta_1$ pour un $x\in G$, revient \`a remplacer $\theta={\rm Int}_{G^\natural}(\delta_1)$ par $\theta'={\rm Int}_G(x)\circ \theta$. 
Nous n'avons pas fait d'hypoth\`ese sur l'ordre de $\theta$, ni sur l'existence d'une 
partie ouverte (non vide!) compacte de $G$ qui soit $\theta$--stable; mais l'on pourrait essayer de choisir 
$\delta_1$ de telle mani\`ere que soit v\'erifi\'ee l'une ou l'autre des 
conditions suivantes:
\begin{enumerate}
\item[$(a)$]L'automorphisme $\theta$ est d'ordre fini.
\item[$(b)$]Pour tout sous--groupe ouvert compact $K$ de $G$, le sous--groupe 
$\bigcap_{i\in\Bbb{Z}}\theta^i(K)$ de $G$ est ouvert.
\item[$(c)$]Il existe une base de voisinages de $1$ dans $G$ form\'ee de sous--groupes 
ouverts compacts $\theta$--stables.
\item[$(d)$]Il existe un sous--groupe ouvert compact $\theta$--stable de $G$.
\item[$(e)$]Il existe une partie ouverte compacte non vide $\theta$--stable de $G$.
\end{enumerate}
On a les implications $(a)\Rightarrow (b)\Rightarrow (c)\Rightarrow (d)
\Rightarrow (e)$. Signalons 
quelques propri\'et\'es:

Si $\theta$ est d'ordre fini, alors toute partie compacte $\Omega$ de $G$ est 
contenue dans une partie ouverte compacte $\theta$--stable $\widetilde{\Omega}$ de 
$G$. En effet, il suffit 
de prendre $\widetilde{\Omega}=(\bigcup_{i\in \Bbb{Z}}\theta^i(\Omega))K$ 
pour un sous--groupe ouvert compact 
$\theta$--stable $K$ de $G$.

Si $K$ est un sous--groupe ouvert compact de $G$, et si $K'$ est un sous--groupe ouvert 
de $K$, alors $K''=\bigcap_{k\in K}k^{-1}K'k$ est un sous--groupe ouvert 
distingu\'e de $K$. Si de plus $K$ et $K'$ sont $\theta$--stables, 
alors $K''$ l'est aussi. On en d\'eduit que si la condition $(c)$ 
est v\'erifi\'ee, alors pour tout sous--groupe ouvert compact $\theta$--stable 
de $G$, il existe une base de voisinages de $1$ dans $G$ 
form\'ee de sous--groupes 
ouverts $\theta$--stables distingu\'es de $K$.

Les conditions $(d)$ et $(e)$ sont équivalentes. En effet, il s'agit de montrer l'implication $(e)\Rightarrow (d)$. 
Soit $U$ une partie ouverte compacte non vide $\theta$--stable de $G$. Pour tout $x\in U$, il existe un 
sous--groupe ouvert compact $K_x$ de $G$ tel que $K_x x\subset U$. Par compacité, on peut écrire $U$ comme réunion finie 
d'ensembles $K_{x_i}x_i$ pour des éléments $x_1,\ldots , x_n\in U$. Posons $K= \bigcap_{i=1}^n K_{x_i}$. C'est encore un 
sous--groupe ouvert compact de $G$, et on a $KU=U$. Pour tout $j\in {\Bbb Z}$, puisque $\theta^j(U)=U$, on a aussi 
$\theta^j(K)U=U$. Mais alors, le sous--groupe $K'$ de $G$ engendré par les $\theta^j(K)$ pour $j\in {\Bbb Z}$ vérifie encore 
$K'U=U$. A fortiori on a $K'\subset UU^{-1}$, et comme $UU^{-1}$ est compact, $K'$ l'est aussi. Le groupe $K'$ est ouvert compact et $\theta$--stable.

\section{Automorphismes d'un groupe rductif connexe}

Dans ce chapitre, on fixe un corps commutatif 
$F$ et une cl\^oture 
alg\'ebrique $\overline{F}$ de $F$. On note $p\geq 1$ l'exposant 
caract\'eristique de $F$; i.e. $p=1$ si la caract\'eristique ${\rm car}(F)$ de $F$ 
est nulle, et $p={\rm car}(F)$ sinon. Si $p>1$, on note $F^{p^{-\infty}}$ la cl\^oture radicielle de 
$F$ dans $\overline{F}$. La r\'ef\'erence adopt\'ee 
est le livre de Borel \cite{Bor}. Toutes les vari\'et\'es 
alg\'ebriques consid\'er\'ees sont suppos\'ees d\'efinies 
sur $\overline{F}$, et sont identifi\'ees 
\`a leur ensemble de points $\overline{F}$--rationnels. Les 
notions topologiques se r\'ef\`erent 
toujours \`a la topologie de 
Zariski. On fixe un groupe alg\'ebrique affine ${\bf H}$, que l'on 
supposera r\'eductif connexe \`a partir de \ref{groupes rductifs}.

\subsection{Groupes alg\'ebriques affines; g\'en\'eralit\'es}\label{groupes algbriques} 
Soit ${\bf H}^\circ$ la composante neutre du groupe ${\bf H}$; c'est un sous--groupe 
ferm\'e distingu\'e d'indice fini de ${\bf H}$ \cite[ch.~I, 1.2]{Bor}. 
On note $R({\bf H})=R({\bf H}^\circ)$ et $R_{\rm u}({\bf H})=
R_{\rm u}({\bf H}^\circ)$ le radical et le radical unipotent de ${\bf H}^\circ$ 
\cite[ch.~IV, 11.21]{Bor}. Par d\'efinition, 
$R({\bf H})$ est un sous--groupe ferm\'e distingu\'e de ${\bf H}$, connexe et r\'esoluble. 
D'autre part, comme $R_{\rm u}({\bf H})$ est l'ensemble des \'el\'ements unipotents de $R({\bf H})$, 
d'apr\`es \cite[ch.~III, 10.6]{Bor}, $R_{\rm u}({\bf H})$ est un sous--groupe ferm\'e distingu\'e de 
${\bf H}$, connexe et unipotent. 

Soit ${\bf H}'$ un sous--groupe ferm\'e 
de ${\bf H}$. On note 
$Z_{\bf H}({\bf H}')$ et $N_{\bf H}({\bf H}')$ le centralisateur et le normali\-sateur 
de ${\bf H}'$ dans ${\bf H}$ \cite[ch.~I, 1.7]{Bor}; ce 
sont des sous--groupes ferm\'es de ${\bf H}$. On note $Z({\bf H})$ le centre 
$Z_{\bf H}({\bf H})$ de ${\bf H}$. 
Le quotient ${\bf H}/{\bf H}'$ existe \cite[ch.~II, 6.8]{Bor}, i.e. il existe un {\it morphisme 
quotient} ${\bf H}\rightarrow {\bf H}/{\bf H}'$ au sens de \cite[ch.~II, 6.1]{Bor}, 
dont les fibres sont les classes $h {\bf H}'$ pour $h\in {\bf H}$. De plus on a 
\cite[ch.~II, 6.8]{Bor}:
\begin{itemize}
\item le quotient ${\bf H}/{\bf H}'$ est une vari\'et\'e alg\'ebrique lisse quasi--projective;
\item si ${\bf H}$ et ${\bf H}'$ sont d\'efinis sur $F$, alors ${\bf H}/{\bf H}'$ l'est aussi, 
i.e. le morphisme quotient ${\bf H}\rightarrow {\bf H}/{\bf H}'$ est d\'efini sur $F$;
\item si ${\bf H}'$ est distingu\'e dans ${\bf H}$, alors 
${\bf H}/{\bf H}'$ est un groupe 
alg\'ebrique affine.
\end{itemize}
En particulier, le quotient ${\bf H}/{\bf H}^\circ$ est un 
groupe alg\'ebrique affine (fini).

On note ${\rm Lie}({\bf H})={\rm Lie}({\bf H}^\circ)$ l'alg\`ebre de 
Lie\footnote{Dans \cite[ch.~I, 3.5]{Bor}, l'alg\`ebre de Lie de ${\bf H}$ est not\'ee 
$L({\bf H})$.}  de ${\bf H}$, 
i.e. l'espace tangent ${\rm T}({\bf H})_1={\rm T}({\bf H}^\circ)_1$ de ${\bf H}$ au point $1$ muni 
de sa structure naturelle d'alg\`ebre de Lie restreinte \cite[ch.~I, 3.5]{Bor}, et l'on pose 
$\frak{H}={\rm Lie}({\bf H})$.

Soit  $\phi:{\bf H}_1\rightarrow {\bf H}_2$ un morphisme de groupes alg\'ebriques, avec 
${\bf H}_1$ et ${\bf H}_2$ affines. Il induit 
par restriction un morphisme de groupes alg\'ebriques $\phi^\circ:{\bf H}^\circ_1 \rightarrow 
{\bf H}^\circ_2$.  Le 
noyau $\ker(\phi)$ est un sous--groupe ferm de ${\bf H}_1$, et l'image $\phi({\bf H}_1)$ est 
un sous--groupe ferm de ${\bf H}_2$; de plus on a l'galit \cite[ch.~I, 1.4]{Bor}
$$\phi({\bf H}_1)^\circ = \phi({\bf H}_1^\circ).
$$ On note 
${\rm Lie}(\phi)={\rm Lie}(\phi^\circ)$ la diff\'erentielle 
${\rm d}(\phi)_1:{\rm Lie}({\bf H}_1)
\rightarrow {\rm Lie}({\bf H}_2)$ de $\phi$ au point $1$; c'est un morphisme d'alg\`ebres de Lie.

Soit ${\bf V}$ une vari\'et\'e alg\'ebrique non vide (a priori ni affine, ni lisse) suppose munie d'une action 
alg\'ebrique de ${\bf H}$ disons \`a gauche\footnote{Le m\^eme discours s'applique bien 
s\^ur aussi \`a une action \`a droite.} 
${\bf H}\times {\bf V}\rightarrow {\bf V},\,(h,v)\mapsto h\cdot v$, et soit 
$v\in {\bf V}$. D'apr\`es \cite[ch.~I, 1.8]{Bor}, l'orbite
$${\bf H}\cdot v= \{ h\cdot v:h\in 
{\bf H}\}$$
est une vari\'et\'e alg\'ebrique lisse, localement ferm\'ee dans ${\bf V}$. Notons
$${\bf H}_v=\{h\in {\bf H}:h\cdot v= v\}$$ 
le stabilisateur de $v$ dans ${\bf H}$; c'est un 
sous--groupe ferm\'e de ${\bf H}$. Le morphisme 
de vari\'et\'es alg\'ebriques $\pi_v:{\bf H}\rightarrow {\bf H}\cdot v,\, h \mapsto h\cdot v$ 
induit, par passage au quotient, un morphisme bijectif de vari\'et\'es alg\'ebriques 
$\bar{\pi}_v:{\bf H}/{\bf H}_v \rightarrow {\bf H}\cdot v$, qui {\it n'est en g\'en\'eral pas} un 
isomorphisme. Pr\'ecis\'ement, notant ${\rm T}({\bf H}\cdot v)_v$ l'espace tangent de 
${\bf H}\cdot v$ au point $v$, et ${\rm d}(\pi_v)_1:\frak{H}\rightarrow {\rm T}({\bf H}\cdot v)_v$ 
la diff\'erentielle de $\pi_v$ au point $1$, on a les inclusions
$$
{\rm Lie}({\bf H}_v)\subset \ker ({\rm d}(\pi_v)_1),\quad 
{\rm d}(\pi_v)_1(\frak{H})\subset {\rm T}({\bf H}\cdot v)_v.
$$
D'apr\`es \cite[ch.~AG, 10.1]{Bor}, on a l'\'egalit\'e
$$
\dim ({\bf H})= \dim ({\bf H}_v)+ \dim({\bf H}\cdot v).\leqno{(*)}
$$
On en d\'eduit que les conditions 
suivantes sont \'equivalentes \cite[ch.~II, 6.7]{Bor}:
\begin{itemize}
\item $\bar{\pi}_v$ est un isomorphisme de vari\'et\'es alg\'ebriques, i.e. l'orbite ${\bf H}\cdot v$ est 
\og le \fg{} quotient de ${\bf H}$ par ${\bf H}_v$;
\item $\pi_v$ est s\'eparable;
\item ${\rm Lie}({\bf H}_v)= \ker ({\rm d}(\pi_v)_1)$;
\item ${\rm d}(\pi_v)_1(\frak{H})= {\rm T}({\bf H}\cdot v)_v$.
\end{itemize}
Si ces conditions sont v\'erifies, on dit aussi\footnote{Contrairement  l'usage, qui est de rserver l'adjectif {\it sparable} 
aux morphismes dominants.} que 
le morphisme ${\bf H}\rightarrow {\bf V},\, h\mapsto h\cdot v$ est s\'eparable. 
Notons que si $p=1$, alors $\pi_v$ est toujours s\'eparable.

\begin{exemples}
{\rm \begin{enumerate}
\item[(1)]Soit ${\bf H}'$ un sous--groupe ferm\'e 
de ${\bf H}$. Alors (par d\'efinition du quotient ${\bf H}/{\bf H}'$) 
le morphisme quotient $\pi:{\bf H}\rightarrow {\bf H}/{\bf H}'$ est surjectif, 
ouvert et s\'eparable: on a l'galit ${\rm Lie}({\bf H}')=\ker({\rm d}(\pi)_1)$.
\item[(2)]Soit  $\phi:{\bf H}_1\rightarrow {\bf H}_2$ un morphisme de groupes alg\'ebriques affines. Par 
passage au quotient, $\phi$ induit un morphisme de groupes alg\'ebriques 
$\bar{\phi}:{\bf H}_1/\ker(\phi)\rightarrow \phi({\bf H}_1)$, qui est un 
isomorphisme de groupes abstraits mais {\it n'est en g\'en\'eral pas} un 
isomorphisme de groupes alg\'ebriques. Pr\'ecis\'ement, on a l'inclusion 
${\rm Lie}(\ker(\phi))\subset \ker({\rm Lie}(\phi))$ avec \'egalit\'e 
si et seulement si $\phi$ 
est s\'eparable, 
i.e. si et seulement si $\bar{\phi}$ 
est un isomorphisme de groupes alg\'ebriques.\hfill $\blacksquare$
\end{enumerate}}
\end{exemples}

Consid\'erons la repr\'esentation adjointe 
${\rm Ad}_{\bf H}:{\bf H}\rightarrow {\rm GL}(\frak{H})$, d\'efinie par 
$${\rm Ad}_{\bf H}(h)={\rm Lie}({\rm Int}_{\bf H}(h))\quad (h\in {\bf H}).$$
Posons 
${\bf H}_{\rm ad}={\rm Ad}_{\bf H}({\bf H})$ et 
${\rm ad}_{\bf H}={\rm Lie}({\rm Ad}_{\bf H}):\frak{H}
\rightarrow {\rm End}(\frak{H})$. Considrons aussi le groupe quotient $\overline{\bf H}={\bf H}/Z({\bf H})$, et notons 
$\pi:{\bf H}\rightarrow \overline{\bf H}$ le morphisme quotient. 
D'apr\`es la propri\'et\'e universelle du quotient \cite[ch.~II, 6.3]{Bor}, il existe un {\it unique} 
morphisme de vari\'et\'es alg\'ebriques $\beta: \overline{\bf H}\rightarrow {\bf H}_{\rm ad}$ 
tel que ${\rm Ad}_{\bf H}= \beta\circ \pi$. De plus, 
$\beta$ est un morphisme de groupes alg\'ebriques, et puisque 
$\pi$ est s\'eparable, ${\rm Ad}_{\bf H}$ est s\'eparable si 
et seulement si $\beta$ est s\'eparable.

\begin{exemples}
{\rm \begin{enumerate}
\item[(3)]
Soit ${\bf H}=\Bbb{SL}_2$. Si $p=2$, alors $Z({\bf H})=\{1\}$, $\overline{\bf H}={\bf H}$ et 
${\bf H}_{\rm ad}=\Bbb{PGL}_2$. De plus (toujours si $p=2$), le 
morphisme ${\rm Ad}_{\bf H}:{\bf H}\rightarrow {\bf H}_{\rm ad}$ est 
un isomorphisme de groupes abstraits, mais il n'est pas s\'eparable.

\item[(4)] Si $p=1$ ou si ${\bf H}$ est connexe et semisimple, on a toujours 
l'\'egalit\'e $Z({\bf H})=\ker({\rm Ad}_{\bf H})$.

\item[(5)]Supposons $p>1$, et reprenons 
l'exemple de Chevalley d\'ecrit dans \cite[ch.~I, 3.15]{Bor}. Soit ${\bf H}$ le sous--groupe 
ferm\'e de $\Bbb{GL}_3$ form\'e des matrices 
$$\varphi(x,y)=
\left(\begin{array}{ccc}
x & 0 & 0\\
0 & x^p & y \\ 
0 & 0 & 1\\
\end{array}\right)
\quad (x\in \smash{\overline{F}}^\times,\,y\in \overline{F}).$$
On a $\varphi(x,y)\varphi(x',y')= \varphi(xx'\!, x^py'+y)$ et 
$\varphi(x,y)^{-1}=\varphi(x^{-1}\!,-x^{-p}y)$, d'o\`u
$$
{\rm Int}_{\bf H}(\varphi(x,y))(\varphi(x'\!,y'))= (x'\!, (1-x'^p)y + x^py').\leqno{\indent\indent (*)}
$$
En particulier le groupe ${\bf H}_{\rm a}=\varphi(1,\overline{F})\simeq \Bbb{G}_{\rm a}$ est distingu\'e 
dans ${\bf H}$, et 
comme
$$
\varphi(x,y)=\varphi(x,0)\varphi(1,y)\quad (x\in \smash{\overline{F}}^\times,\,y\in \overline{F}),
$$
posant 
${\bf H}_{\rm m}= \varphi(\smash{\overline{F}}^\times,0)\simeq \Bbb{G}_{\rm m}$, on a la d\'ecomposition en produit semidirect
$${\bf H}= {\bf H}_{\rm m}\ltimes {\bf H}_{\rm a}.$$ 
Gr\^ace \`a l'\'egalit\'e $(*)$, on obtient que $Z({\bf H})=\{1\}$ et $\ker({\rm Ad}_{\bf H})= {\bf H}_{\rm a}$. 
D'autre part (\`a nouveau gr\^ace \`a $(*)$), 
l'application commutateur ${\bf H}\times {\bf H}\rightarrow {\bf H},\, (h,h')\mapsto 
c_h(h')=hh'h^{-1}h'^{-1}$ est donn\'ee par
$$
c_{\varphi(x,y)}(\varphi(x'\!,y'))= (1, (1-x'^p)y - (1-x^p)y').
$$
On en d\'eduit que l'alg\`ebre de Lie de ${\bf H}$ est commutative (bien 
que ${\bf H}$ ne soit pas commutatif). 
Par suite l'inclusion
$${\rm Lie}({\bf H}_{\rm a})={\rm Lie}(\ker({\rm Ad}_{\bf H}))\subset 
\ker({\rm ad}_{\bf H})=\frak{H}$$ est stricte, et ${\rm Ad}_{\bf H}$ 
n'est pas s\'eparable.\hfill$\blacksquare$
\end{enumerate}}
\end{exemples}

\subsection{Automorphismes}\label{automorphismes} Appelons {\it $\overline{F}$--automorphisme de 
${\bf H}$} un automorphisme du groupe alg\'ebrique ${\bf H}$. 
On note ${\rm Aut}_{\overline{F}}({\bf H})$ le groupe des $\overline{F}$--automorphismes 
de ${\bf H}$, et ${\rm Int}_{\overline{F}}({\bf H})$ le sous--groupe distingu\'e de 
${\rm Aut}_{\overline{F}}({\bf H})$ form\'e des automorphismes int\'erieurs, i.e. ceux 
qui sont de la forme ${\rm Int}_{\bf H}(h): x\mapsto h xh^{-1}$ pour un 
$h\in {\bf H}$. 
L'application ${\bf H}\rightarrow {\rm Int}_{\overline{F}}({\bf H}),\, 
h\mapsto {\rm Int}_{\bf H}(h)$ se factorise en 
un isomorphisme de groupes abstraits $\overline{\bf H}\rightarrow {\rm Int}_{\overline{F}}({\bf H})$. 
Notons que pour $\tau\in {\rm Aut}_{\overline{F}}({\bf H})$ 
et $h,\,x\in {\bf H}$, 
on a
$$
{\rm Int}_{\bf H}(x^{-1})\circ {\rm Int}_{\bf H}(h)\circ \tau \circ {\rm Int}_{\bf H}(x)= 
{\rm Int}_{\bf H}(x^{-1}h\tau(x))\circ \tau.\leqno{(*)}
$$On note ${\rm Out}_{\overline{F}}({\bf H})$ le groupe quotient 
${\rm Aut}_{\overline{F}}({\bf H})/{\rm Int}_{\overline{F}}({\bf H})$, et 
${\rm Aut}_{\overline{F}}^0({\bf H})$ le sous--groupe de 
${\rm Aut}_{\overline{F}}({\bf H})$ form\'e des automorphismes dont l'image dans ${\rm Out}_{\overline{F}}({\bf H})$ est d'ordre fini.

Si ${\bf H}$ est d\'efini sur $F$, on note ${\rm Aut}_F({\bf H})$ le groupe des $F$--automorphismes 
de ${\bf H}$, i.e. le sous--groupe de 
${\rm Aut}_{\overline{F}}({\bf H})$ form\'e des automorphismes qui sont d\'efinis sur $F$, 
et l'on pose
$${\rm Aut}_F^0({\bf H})={\rm Aut}_F({\bf H})\cap {\rm Aut}_{\overline{F}}^0({\bf H}),
$$
$$
{\rm Int}_F({\bf H})= {\rm Aut}_F({\bf H})\cap {\rm Int}_{\overline{F}}({\bf H}).
$$
Si ${\bf H}$ et $Z({\bf H})$ sont d\'efinis sur $F$, alors $\overline{\bf H}$ l'est aussi 
\cite[ch.~II, 6.8]{Bor}, et via l'identification 
${\rm Int}_{\overline{F}}({\bf H})=\overline{\bf H}$, 
${\rm Int}_F({\bf H})$ co\"{\i}ncide avec le groupe $\overline{\bf H}(F)$ des 
points $F$--rationnels de $\overline{\bf H}$. 

\begin{mesrems}
{\rm \begin{enumerate}
\item[(1)]Si ${\bf H}$ est r\'eductif connexe, on verra plus loin (\ref{groupes rductifs}) que  
la projection canonique 
${\rm Aut}_{\overline{F}}({\bf H})\rightarrow 
{\rm Out}_{\overline{F}}({\bf H})$ est scind\'ee.
\item[(2)]Si ${\bf H}$ est un tore, alors on a 
${\rm Aut}_{\overline{F}}({\bf H})={\rm Out}_{\overline{F}}({\bf H})\simeq GL({\rm X}^*({\bf H}))$ 
o\`u ${\rm X}^*({\bf H})$ d\'esigne le groupe des caract\`eres 
alg\'ebriques de ${\bf H}$. En particulier si ${\bf H}$ est un tore de dimension $1$, 
alors le seul $\overline{F}$--automorphisme non trivial de 
${\bf H}$ est le passage \`a l'inverse $t\mapsto t^{-1}$.
\item[(3)]Si ${\bf H}$ est un tore d\'efini et 
d\'eploy\'e sur $F$, alors on a ${\rm Aut}_{\overline{F}}({\bf H})={\rm Aut}_F({\bf H})$.\hfill $\blacksquare$
\end{enumerate}}
\end{mesrems}

Soit $\tau\in {\rm Aut}_{\overline{F}}({\bf H})$. Par restriction, $\tau$ induit 
un $\overline{F}$--automorphisme de ${\bf H}^\circ$, que l'on note $\tau^\circ$. 
L'ensemble ${\bf H}_\tau=\{h\in {\bf H}:\tau(h)=h\}$ est un sous--groupe ferm\'e de ${\bf H}$, et 
${\bf H}_\tau^\circ =({\bf H}_\tau)^\circ$ 
est un sous--groupe ferm\'e de ${\bf H}^\circ$ qui co\"{\i}ncide avec 
$({\bf H}^\circ)_{\tau^\circ}^\circ$. Notons $1-\tau:{\bf H}\rightarrow {\bf H}$ le 
morphisme de vari\'et\'es alg\'ebriques 
$h\mapsto h\tau(h)^{-1}$, et 
posons
$${\bf H}(1-\tau)=\{h\tau(h)^{-1}:h\in {\bf H}\}\;(=\{h^{-1}\tau(h):h\in {\bf H}\}).
$$ D'apr\`es \ref{groupes algbriques}, 
${\bf H}(1-\tau)$ est
une vari\'et\'e alg\'ebrique lisse, localement ferm\'ee dans ${\bf H}$. Remarquons que 
si ${\bf H}$ est commutatif, alors $1-\tau$ est un morphisme de groupes 
alg\'ebriques, ${\bf H}_\tau=\ker(1-\tau)$, et ${\bf H}(1-\tau)={\rm Im}(1-\tau)$ est un sous--groupe 
ferm\'e de ${\bf H}$. Revenons au cas g\'en\'eral. Par passage au quotient, $1-\tau$ 
induit un morphisme bijectif de vari\'et\'es alg\'ebriques
$$
\overline{1-\tau}:{\bf H}/{\bf H}_\tau\rightarrow 
{\bf H}(1-\tau),
$$
qui n'est en g\'en\'eral pas un isomorphisme. 
Posons $\frak{H}_\tau=\ker({\rm id}_\frak{H}-{\rm Lie}(\tau))$. \`A nouveau d'apr\`es 
ce qui pr\'ec\`ede, on a l'inclusion 
${\rm Lie}({\bf H}_\tau)\subset \frak{H}_\tau$ avec \'egalit\'e si et 
seulement si $1-\tau$ est s\'eparable, i.e. si et seulement si 
$\overline{1-\tau}$ est un isomorphisme de vari\'et\'es alg\'ebriques.

Si ${\bf J}$ est un sous--groupe ferm\'e de ${\bf H}$ tel que $\tau({\bf J})={\bf J}$, alors 
la restriction $\sigma =\tau\vert_{\bf J}$ appartient \`a ${\rm Aut}_{\overline{F}}({\bf J})$, 
et l'on pose ${\bf J}_\tau={\bf J}_\sigma$, ${\bf J}^\circ_\tau ={\bf J}^\circ_\sigma$, 
${\bf J}(1-\tau)={\bf J}(1-\sigma)$, (etc.). Pour $h\in {\bf H}$ et $\tau={\rm Int}_{\bf H}(h)$, on a 
$\tau^\circ ={\rm Int}_{{\bf H}^\circ}(h)$, et l'on pose ${\bf H}_h={\bf H}_\tau$ et 
${\bf H}_h^\circ = {\bf H}_\tau^\circ$. 

Puisque ${\bf H}$ est affine, pour $h\in {\bf H}$, on a 
la {\it d\'ecomposition de Jordan} \cite[ch.~I, 4.4]{Bor}
$$h=h_{\rm s}h_{\rm u}$$
avec $h_{\rm s}\in {\bf H}$ 
semisimple, $h_{\rm u}\in {\bf H}$ unipotent, et 
$h_{\rm u}\in {\bf H}_{h_{\rm s}}$; cette d\'ecomposition 
est {\it unique}. En particulier, pour $h\in {\bf H}$, l'automorphisme intérieur ${\rm Int}_{\bf H}(h)$ de ${\bf H}$ se décompose en 
$$
{\rm Int}_{\bf H}(h)={\rm Int}_{\bf H}(h_{\rm s})\circ {\rm Int}_{\bf H}(h_{\rm u})= {\rm Int}_{\bf H}(h_{\rm u})\circ {\rm Int}_{\bf H}(h_{\rm s}).
$$
Si de plus ${\bf H}$ est d\'efini sur $F$, alors d'apr\`es loc.~cit., pour 
$h\in {\bf H}(F)$, les lments $h_{\rm s}$ et $h_{\rm u}$ appartiennent \`a 
${\bf H}(F^{p^{-\infty}})$, par conséquent les automorphismes ${\rm Int}_{\bf H}(h_{\rm s})$ et 
${\rm Int}_{\bf H}(h_{\rm u})$ de ${\bf H}$ sont définis sur $F^{p^{-\infty}}$. 

\begin{mesrems}
{\rm 
\begin{enumerate}
\item[(4)] Pour $h\in {\bf H}^\circ$, 
on a $h_{\rm s},\,h\in {\bf H}_{h_{\rm s}}^\circ$ \cite[ch.~IV, 11.12]{Bor} 
donc $h_{\rm u}\in {\bf H}_{h_{\rm s}}^\circ$. 
\item[(5)] Supposons $p=1$. Puisque le groupe alg\'ebrique affine ${\bf H}/{\bf H}^\circ$ est fini, il 
est form\'e d'\'el\'ements semisimples. 
Par cons\'equent tout \'el\'ement unipotent de ${\bf H}$ appartient \`a ${\bf H}^\circ$. 
En particulier pour $h\in {\bf H}$, on a $h_{\rm u}\in {\bf H}_{h_{\rm s}}^\circ$.\hfill $\blacksquare$
\end{enumerate}}
\end{mesrems}

\subsection{Groupes diagonalisables et tores}\label{groupes diago et tores} On note $\Bbb{G}_{\rm m}$ le groupe 
alg\'ebrique affine $\smash{\overline{F}}^\times$, et pour $k\in \Bbb{Z}_{\geq 1}$, on pose 
$\Bbb{G}_{\rm m}^k= \Bbb{G}_{\rm m}\times \cdots \times \Bbb{G}_{\rm m}$ ($k$ fois). Soit 
${\bf D}$ un {\it groupe diagonalisable}, c'est--\`a--dire un groupe alg\'ebrique affine isomorphe \`a 
un sous--groupe ferm\'e de $\Bbb{G}_{\rm m}^k$ pour un entier $k\geq 1$. On appelle {\it caract\`ere alg\'ebrique de ${\bf D}$} un morphisme de groupes alg\'ebriques ${\bf D}\rightarrow \Bbb{G}_{\rm m}$. 
On note ${\rm X}^*({\bf D})$ le groupe des caract\`eres alg\'ebriques de ${\bf D}$. C'est un 
$\Bbb{Z}$--module de type fini, sans $p$--torsion si $p>1$, et l'on a une identification canonique:
$$
{\bf D}={\rm Hom}_\Bbb{Z}({\rm X}^*({\bf D}),\smash{\overline{F}}^\times).
$$
Tout sous--groupe ferm\'e de ${\bf D}$ est 
un groupe diagonalisable. Pour tout morphisme de groupes 
alg\'ebriques $\phi:{\bf D}\rightarrow {\bf H}$, 
l'image $\phi({\bf D})$ est un groupe diagonalisable. En particulier, 
pour tout sous--groupe ferm\'e ${\bf D}'$ de ${\bf D}$, le quotient 
${\bf D}/{\bf D}'$ est un groupe diagonalisable, 
de groupe des caract\`eres 
$${\rm X}^*({\bf D}/{\bf D}')=\{\chi\in {\rm X}^*({\bf D}):\chi\vert_{{\bf D}'}={\rm id}\}.
$$
La composante neutre 
${\bf D}^\circ$ de ${\bf D}$ 
est un tore, de groupe des caract\`eres algbriques
$$
{\rm X}^*({\bf D}^\circ)= {\rm X}^*({\bf D})/{\rm X}^*({\bf D})_{\rm tor},
$$ o\`u 
${\rm X}^*({\bf D})_{\rm tor}$ d\'esigne le sous--module de torsion de ${\rm X}^*({\bf D})$. 
On a 
${\rm X}^*({\bf D}/{\bf D}^\circ)={\rm X}^*({\bf D})_{\rm tor}$, et 
d'apr\`es \cite[ch.~III, 8.7]{Bor}, il existe un sous--groupe ferm\'e {\it fini} ${\bf \Omega}$ de 
${\bf D}$ tel que le morphisme produit ${\bf D}^\circ \times {\bf \Omega}\rightarrow {\bf D}$ est 
un isomorphisme de groupes alg\'ebriques. Notons que 
pour tout morphisme de groupes 
alg\'ebriques $\phi:{\bf D}\rightarrow {\bf H}$, le 
groupe $\phi({\bf D}^\circ)=\phi({\bf D})^\circ$ est un tore.

Soit $\phi:{\bf D}\rightarrow {\bf D}'$ un morphisme 
de groupes alg\'ebriques, avec ${\bf D}$ et ${\bf D}'$ diagonalisables. Il induit un morphisme 
de $\Bbb{Z}$--modules
$$\phi^\sharp:{\rm X}^*({\bf D}')\rightarrow {\rm X}^*({\bf D}),\,
\chi'\mapsto \chi\circ \phi.
$$ Cela d\'efinit un foncteur contravariant 
de la cat\'egorie des groupes diagonalisables dans celle des $\Bbb{Z}$--modules de type 
fini. D'apr\`es \cite[ch.~III, 8.3]{Bor}, ce foncteur est pleinement fid\`ele. On a 
${\rm X}^*(\phi({\bf D}))= {\rm X}^*({\bf D}')/\ker({\phi}^\sharp)$, et ${\rm X}^*(\ker(\phi))$ 
co\"{\i}ncide avec:
\begin{itemize}
\item ${\rm coker}(\phi^\sharp)$ si $p=1$;
\item le quotient de ${\rm coker}(\phi^\sharp)$ par sa $p$--torsion 
si $p>1$.
\end{itemize}
Si $p>1$, 
$\phi$ est s\'eparable si et seulement si le $\Bbb{Z}$--module 
${\rm coker}(\phi^\sharp)$ est 
sans $p$--torsion.

\begin{exemple}{\rm 
Soit ${\bf T}$ un tore, et soit $m$ un entier $\geq 1$. 
Consid\'erons le morphisme de groupes alg\'ebriques $\phi:{\bf T}\rightarrow {\bf T}$ d\'efini 
par $\phi(x)=x^m$. On a:
\begin{itemize}
\item $\phi({\bf T})={\bf T}$;
\item si $(m,p)=1$, alors $\phi$ est s\'eparable et le noyau $\ker(\phi)$ est isomorphe (comme 
groupe abstrait) \`a $(\Bbb{Z}/m\Bbb{Z})^{\dim({\bf T})}$;
\item si $p>1$ et $m=p^r$ pour un entier $r$, alors $\ker(\phi)=\{1\}$ et ${\rm Lie}(\phi)=0$.
\hfill $\blacksquare$
\end{itemize}}
\end{exemple}

Soit ${\bf T}$ un 
tore, et soit $\tau\in {\rm Aut}_{\overline{F}}({\bf T})$. Alors $1-\tau$ est un morphisme de 
groupes alg\'ebriques, et puisque le groupe ${\bf T}(1-\tau)={\rm Im}(1-\tau)$ est connexe, 
c'est un tore. 
Soit $q_\tau:{\bf T}\rightarrow {\bf T}/{\bf T}(1-\tau)$ 
le morphisme quotient. On a 
la suite exacte longue de groupes alg\'ebriques:
$$
1\rightarrow {\bf T}_\tau \rightarrow {\bf T}\buildrel{1-\tau}\over{\longrightarrow} {\bf T} 
\buildrel{q_\tau}\over{\longrightarrow} {\bf T}/{\bf T}(1-\tau)\rightarrow 1.
$$
Posons $\frak{T}={\rm Lie}({\bf T})$. D'apr\`es \ref{groupes algbriques}, on a 
l'inclusion ${\rm Lie}({\bf T}_\tau)\subset \frak{T}_\tau(= \ker({\rm id}_\frak{T}-{\rm Lie}(\tau)))$ avec 
\'egalit\'e si et seulement si $1-\tau$ est s\'eparable.  

Supposons de plus que $\tau$ est d'ordre fini $m$. Alors $(\tau^\sharp)^m$ est l'identit\'e de 
${\rm X}^*({\bf T})$. Soit  
$\ES{N}_\tau:{\bf T}\rightarrow{\bf T}$ le morphisme de groupes alg\'ebriques 
d\'efini par
$$\ES{N}_\tau(t)=t\tau(t)\cdots \tau^{m-1}(t).
$$ 
Puisque $(1-\tau)^\sharp = {\rm id}-\tau^\sharp$ et $\ES{N}_\tau^\sharp 
={\rm id}+\tau^\sharp +\cdots + (\tau^\sharp)^{m-1}$, on a 
${\rm Im}({\rm id}-\tau^\sharp)\subset \ker(\ES{N}_\tau^\sharp)$ et
$$
\ES{N}_\tau^\sharp(\chi)- m\chi\in {\rm Im}({\rm id}-\tau^\sharp)\quad (\chi\in {\rm X}^*({\bf T})).
$$
Par suite, le quotient $\ker(\ES{N}_\tau^\sharp)/ {\rm Im}({\rm id}-\tau^\sharp)$ est un $\Bbb{Z}$--module 
de torsion d'exposant divisant $m$. D'autre part, 
${\rm X}^*({\bf T}_\tau/{\bf T}_\tau^\circ) ={\rm X}^*({\bf T}_\tau)_{\rm tor}$ 
co\"{\i}ncide avec la $p'$--torsion 
de ${\rm coker}({\rm id}-\tau^\sharp)$; o\`u par {\it $p'$--torsion} 
on entend:
\begin{itemize}
\item la torsion si $p=1$;
\item le quotient de la torsion par la $p$--torsion si $p>1$.
\end{itemize}
Comme $\ES{N}_\tau({\bf T})$ est un tore, le quotient ${\rm X}^*({\bf T})/\ker(\ES{N}_\tau^\sharp)={\rm X}^*(\ES{N}_\tau({\bf T}))$ est sans torsion. 
Par suite ${\rm X}^*({\bf T}_\tau/{\bf T}_\tau^\circ)$ 
co\"{\i}ncide avec la $p'$--torsion de $\ker(\ES{N}_\tau^\sharp)/{\rm Im}({\rm id}-\tau^\sharp)$. 
D'o\`u le

\begin{monlem}
Soit ${\bf T}$ un tore, et  soit $\tau$ 
un $\overline{F}$--automorphisme de ${\bf T}$ d'ordre fini $m$. Alors 
l'exposant de ${\bf T}_\tau/{\bf T}_\tau^\circ={\rm Hom}_\Bbb{Z}({\rm X}({\bf T}_\tau)_{\rm tor},\smash{\overline{F}}^\times)$ divise $m$.
\end{monlem}

\begin{marema}
{\rm Soit $\tau\in {\rm Aut}_{\overline{F}}({\bf H})$ 
et soit ${\bf T}$ un tore $\tau$--stable de ${\bf H}$. Notons $\tau'$ la 
restriction de $\tau$ \`a ${\bf T}$. 
Puisque $N_{\bf H}({\bf T})^\circ =Z_{\bf H}({\bf T})^\circ$ 
\cite[ch.~III, 8.10, cor. 2]{Bor}, le 
groupe $N_{\bf H}({\bf T})/Z_{\bf H}({\bf T})$ est fini. Par cons\'equent si 
$\tau\in {\rm Aut}_{\overline{F}}^0({\bf H})$, alors $\tau'$ est d'ordre fini.\hfill $\blacksquare$}
\end{marema}

\subsection{Automorphismes semisimples et unipotents}\label{automorphismes ss et u} Un $\overline{F}$--automorphisme 
$\tau$ de ${\bf H}$ est dit {\it localement fini}\footnote{On peut v\'erifier que cette notion 
est \'equivalente \`a la suivante (cf. \cite[ch.~I, convention pp. 81/82]{Bor}): l'alg\`ebre affine $\overline{F}[{\bf H}]$ de ${\bf H}$ 
est r\'eunion de sous--espaces $\tau^\sharp$--stables de dimension finie, 
o\`u $\tau^\sharp$ d\'esigne l'automorphisme de $\overline{F}[{\bf H}]$ d\'efini par 
$\tau^\sharp(f)(h)=f(\tau(h))$ pour $f\in \overline{F}[{\bf H}]$ et $h\in {\bf H}$.} s'il existe un morphisme de groupes 
alg\'ebriques
$$
\iota:{\bf H}\rightarrow \Bbb{GL}_n=GL_n(\overline{F})
$$ 
qui soit un isomorphisme (de 
groupes alg\'ebriques) de ${\bf H}$ 
sur un sous--groupe ferm\'e de $\Bbb{GL}_n$, et un \'el\'ement 
$g\in \Bbb{GL}_n$, tels que pour tout $h\in {\bf H}$, on 
ait $\iota\circ \tau(h)= g\iota(h)g^{-1}$.

Soit $\tau$ un $\overline{F}$--automorphisme localement fini de ${\bf H}$. On d\'efinit comme suit la 
{\it d\'ecomposition 
de Jordan} $\tau=\tau_{\rm s}\circ \tau_{\rm u}$ de $\tau$. Choisissons $\iota:{\bf H}\rightarrow 
\Bbb{GL}_n$ et $g$ comme ci-dessus, et posons ${\bf H}^\iota=\iota({\bf H})\subset \Bbb{GL}_n$. 
Soit 
$g=g_{\rm s}g_{\rm u}$ la d\'ecomposition de Jordan de $g$ (dans $\Bbb{GL}_n$). 
Puisque 
$g\in N_{\Bbb{GL}_n}({\bf H}^\iota)$ et que 
$N_{\Bbb{GL}_n}({\bf H}^\iota)$ est un sous--groupe ferm\'e de $\Bbb{GL}_n$, on a
$g_{\rm s},\,g_{\rm u}\in N_{\Bbb{GL}_n}({\bf H}^\iota)$. Notons  
$\tau^\iota_{\rm s}$ et ${\rm \tau}^\iota_{\rm u}$ les \'el\'ements de ${\rm Aut}_{\overline{F}}({\bf H}^\iota)$ 
d\'eduit de ${\rm Int}_{\Bbb{GL}_n}(g_{\rm s})$ et ${\rm Int}_{\Bbb{GL}_n}(g_{\rm u})$ 
par restriction, et posons $\tau_{\rm s}= \iota^{-1}\circ \tau^\iota_{\rm s}\circ \iota$ et 
$\tau_{\rm u}= \iota^{-1}\circ \tau^\iota_{\rm u}\circ \iota$. Les $\overline{F}$--automorphismes 
$\tau_{\rm s}$ et 
${\rm \tau}_{\rm u}$ de ${\bf H}$ sont bien d\'efinis, i.e. ils ne 
d\'ependent pas des choix de $\iota$ et $g$. En effet, fix\'e $\iota$, l'image $\bar{g}$ de 
$g$ dans le groupe alg\'ebrique affine $\smash{\overline{\bf H}}^\iota=N_{\Bbb{GL}_n}({\bf H}^\iota)/Z_{\Bbb{GL}_n}({\bf H}^\iota)$ 
est d\'etermin\'ee de mani\`ere unique par $\tau$, et si 
$\bar{g}=\bar{g}_{\rm s}\bar{g}_{\rm u}$ 
est la d\'ecomposition de Jordan de $\bar{g}$ dans 
$\smash{\overline{\bf H}}^\iota$, alors ${\rm Int}_{\smash{\overline{\bf H}}^\iota}(\bar{g}_{\rm s})$ et ${\rm Int}_{\smash{\overline{\bf H}}^\iota}(\bar{g}_{\rm u})$ 
d\'efinissent des $\overline{F}$--automorphismes de ${\bf H}^\iota$, 
qui co\"{\i}ncident avec $\tau^\iota_{\rm s}$ et $\tau^\iota_{\rm u}$; d'o\`u 
l'ind\'ependance par rapport au choix de $g$. Si maintenant 
$\iota_1:{\bf H}\rightarrow \Bbb{GL}_m$ est un autre morphisme de 
groupes alg\'ebriques qui soit un isomorphisme sur un sous--groupe 
ferm\'e ${\bf H}^{\iota_1}$ de $\Bbb{GL}_m$, et si $g_1$ est un lment de $\Bbb{GL}_m$ tel que 
$\iota_1\circ \tau(h)= g_1\iota(h)g_1^{-1}$ pour tout $h\in {\bf H}$, alors en notant 
$\jmath:{\bf H}\rightarrow \Bbb{GL}_{n+m}$ le morphisme compos\'e 
du morphisme produit $\iota\times \iota_1:{\bf H}\rightarrow \Bbb{GL}_n\times \Bbb{GL}_m$ et du morphisme diagonal par blocs 
$\delta:\Bbb{GL}_n\times \Bbb{GL}_m\rightarrow {\Bbb G}_{n+m}$, on 
a $\jmath \circ \tau(h)= \delta(g \iota (h)g^{-1}, g_1 \iota_1(h)g_1^{-1})$ pour tout $h\in {\bf H}$. Par conséquent $\tau_{\rm s}^j= \delta\circ (\tau_{\rm s}^\iota\times \tau_{\rm s}^{\iota_1})\circ \delta^{-1}$ et $\tau_{\rm u}^j= \delta\circ (\tau_{\rm u}^\iota\times \tau_{\rm u}^{\iota_1})\circ \delta^{-1}$. On en déduit que $
\iota_1^{-1}\circ \tau_{\rm s}^{\iota_1} \circ \iota_1=  \iota^{-1}\circ \tau_{\rm s}^{\iota}\circ \iota$ et
$\iota_1^{-1}\circ \tau_{\rm u}^{\iota_1} \circ \iota_1= 
\iota^{-1}\circ \tau_{\rm u}^{\iota}\circ \iota$.

Par construction, 
on a la d\'ecomposition $\tau=\tau_{\rm s}\circ \tau_{\rm u}=\tau_{\rm u}\circ \tau_{\rm s}$, et $\tau_{\rm s}$ et $\tau_{\rm u}$ 
sont localement finis. Cette d\'ecomposition est {\it unique} (car 
la d\'ecomposition de Jordan de $\bar{g}$ est unique). En 
particulier, si ${\bf J}$ est un sous--groupe ferm\'e $\tau$--stable de 
${\bf H}$, alors $\tau_{\rm s}({\bf J})={\bf J}$ et $\tau_{\rm u}({\bf J})={\bf J}$. 
Enfin si ${\bf H}$ est d\'efini sur $F$ et $\tau\in {\rm Aut}_F({\bf H})$, 
alors on peut choisir $\iota$ d\'efini sur $F$, et $\tau_{\rm s}$ et $\tau_{\rm u}$ sont 
d\'efinis sur $F^{p^{-\infty}}$. 

Pour $h\in {\bf H}$, on a bien s\^ur ${\rm Int}_{\bf H}(h)_{\rm s}=
{\rm Int}_{\bf H}(h_{\rm s})$ et ${\rm Int}_{\bf H}(h)_{\rm u}=
{\rm Int}_{\bf H}(h_{\rm u})$.

Un $\overline{F}$--automorphisme localement fini $\tau$ de ${\bf H}$ est dit 
{\it semisimple} si $\tau=\tau_{\rm s}$, et il est dit 
${\it unipotent}$ si $\tau=\tau_{\rm u}$.

Si $\tau$ est un $\overline{F}$--automorphisme localement fini de ${\bf H}$, 
alors le $\overline{F}$--automorphisme $\tau^\circ$ 
de ${\bf H}^\circ$ est localement fini, et l'on a $(\tau^\circ)_{\rm s}=
(\tau_{\rm s})^\circ$ et $(\tau^\circ)_{\rm u}=
(\tau_{\rm u})^\circ$. On peut donc poser $\tau_{\rm s}^\circ = (\tau_{\rm s})^\circ$ et 
$\tau_{\rm u}^\circ = (\tau_{\rm u})^\circ$.

Si $\tau\in {\rm Aut}_{\overline{F}}^0({\bf H}^\circ)$, alors il existe un groupe 
alg\'ebrique affine ${\bf H}'$ de composante neutre ${\bf H}^\circ$, tel que 
$\tau={\rm Int}_{{\bf H}'}(h')\vert_{{\bf H}^\circ}$ pour un \'el\'ement $h'\in {\bf H}'$; 
en particulier, $\tau$ est localement fini. En effet (voir aussi plus loin, \ref{l'application N}), 
soit $l$ l'ordre de l'image de $\tau$ dans ${\rm Out}_{\overline{F}}({\bf H}^\circ)$. 
Choisissons un \'el\'ement $h\in {\bf H}^\circ$ tel que $\tau^l = {\rm Int}_{{\bf H}^\circ}(h)$. 
Notons ${\bf H}^\circ \rtimes \langle \tau\rangle$ le produit semidirect (dans la cat\'egorie des 
groupes) de ${\bf H}^\circ$ par le groupe abstrait engendr\'e par $\tau$, et 
${\bf C}$ le sous--groupe cyclique de ${\bf H}^\circ \rtimes \langle \tau\rangle$ engendr\'e 
par $h^{-1}\!\rtimes \tau^l$. Alors ${\bf C}$ est distingu\'e dans ${\bf H}^\circ \rtimes \langle \tau\rangle$, 
et le groupe quotient ${\bf H}'= ({\bf H}^\circ \rtimes \langle \tau\rangle)/{\bf C}$ est 
naturellement muni d'une structure de groupe alg\'ebrique affine de composante 
neutre ${\bf H}^\circ$. On prend alors pour $h'$ l'image de $1\rtimes \tau$ dans ${\bf H}'$. 
Soit maintenant $h'=h'_{\rm s}h'_{\rm u}$ la d\'ecomposition de Jordan de $h'$ dans ${\bf H}'$. On a 
$\tau_{\rm s}={\rm Int}_{{\bf H}'}(h'_{\rm s})\vert_{{\bf H}^\circ}$ et $\tau_{\rm s}={\rm Int}_{{\bf H}'}(h'_{\rm u})\vert_{{\bf H}^\circ}$.

\begin{marema1}
{\rm Soit $\tau\in {\rm Aut}_{\overline{F}}({\bf H})$. 
Si $p>1$, alors $\tau$ est unipotent si et seulement s'il existe un 
entier $k\geq 1$ tel que $\tau^{p^k}={\rm id}_{\bf H}$. Si $p=1$ et 
$\tau$ est d'ordre fini, alors $\tau$ est toujours semisimple.
\hfill $\blacksquare$}
\end{marema1}

\begin{marema2}{\rm 
Soit ${\bf T}$ un tore. D'apr\`es la remarque de \ref{groupes diago et tores}, 
un $\overline{F}$--automorphisme de ${\bf T}$ 
est localement fini si et seulement s'il est d'ordre fini. Par suite si $p=1$, tout 
automorphisme localement fini de ${\bf T}$ est semisimple.\hfill $\blacksquare$}
\end{marema2}

\begin{mapropo}
Soit ${\bf T}$ un tore, et soit $\tau$ un $\overline{F}$--automorphisme 
unipotent de ${\bf T}$. Alors ${\bf T}_\tau$ est connexe, et le morphisme produit 
${\bf T}_\tau\times {\bf T}(1-\tau)\rightarrow {\bf T}$ est bijectif.
\end{mapropo}

\begin{proof}
Si $p=1$, alors $\tau={\rm id}_{\bf H}$ (d'aprs la remarque 2) et il n'y rien \`a d\'emontrer. 
On suppose donc $p>1$ et $\tau\neq 
{\rm id}_{\bf H}$. Alors $\tau$ est d'ordre fini $p^k$ pour un entier $k\geq 1$. 
Puisque ${\bf T}_\tau/{\bf T}_\tau^\circ = {\rm Hom}_\Bbb{Z}({\rm X}^*({\bf T}_\tau)_{\rm tor},
\smash{\overline{F}}^\times)$ et ${\rm X}^*({\bf T}_\tau)$ est sans 
$p$--torsion, les \'el\'ements de ${\bf T}_\tau/{\bf T}_\tau^\circ $ sont 
d'ordre premier \`a $p$. Comme d'autre part l'exposant de ${\bf T}_\tau/{\bf T}_\tau^\circ$ 
divise $p^k$ (lemme de \ref{groupes diago et tores}), 
on a ${\bf T}_\tau/{\bf T}_\tau^\circ=\{1\}$, i.e. ${\bf T}_\tau$ est connexe. 
Notons $\ES{N}_\tau$ le morphisme de 
groupes alg\'ebriques ${\bf T}\rightarrow {\bf T},\,t\mapsto t\tau(t)\cdots \tau^{p^k-1}(t)$. 
Puisque $\ker(\ES{N}_\tau^\sharp)/{\rm Im}({\rm id}-\tau^\sharp)$ est un 
groupe de $p$--torsion, on a $\ES{N}_\tau({\bf T})=\ker(1-\tau)={\bf T}_\tau$. 
De la m\^eme mani\`ere, comme 
$$
\ES{N}_\tau(t)\equiv t^{p^k}\;({\rm mod}\,{\bf T}(1-\tau))\quad (t\in {\bf T}),
$$
$\ker(\ES{N}_\tau)/{\bf T}(1-\tau)$ 
est un groupe de torsion d'exposant divisant $p^k$, et donc 
$\ker(\ES{N}_\tau)={\bf T}(1-\tau)$. Puisque 
$ {\bf T}_\tau\cap {\bf T}(1-\tau)=\{t\in {\bf T}: t^{p^k}=1\}$, on a 
${\bf T}_\tau\cap {\bf T}(1-\tau)=\{1\}$, et le morphisme 
produit $\mu: {\bf T}_\tau\times {\bf T}(1-\tau)\rightarrow {\bf T}$ est 
injectif. L'image ${\bf T}_\tau{\bf T}(1-\tau)=\mu({\bf T}_\tau\times {\bf T}(1-\tau))$ 
est un sous--groupe ferm\'e connexe de ${\bf T}$; c'est donc un 
tore. Puisque ${\bf T}_\tau=\ES{N}_\tau({\bf T})$ et ${\bf T}(1-\tau)=\ker(\ES{N}_\tau)$, 
on a
$$\dim({\bf T}_\tau{\bf T}(1-\tau))=\dim({\bf T}_\tau)+\dim({\bf T}(1-\tau))= \dim({\bf T}),
$$ et 
$\mu$ est surjectif. D'o\`u 
l'\'egalit\'e ${\bf T}_\tau{\bf T}(1-\tau)={\bf T}$.
\end{proof}

\begin{marema}{\rm 
Soit ${\bf T}$ et $\tau$ comme dans la proposition. 
Le morphisme produit 
$$\mu:{\bf T}_\tau\times {\bf T}(1-\tau)\rightarrow {\bf T}$$ n'est en g\'en\'eral pas 
un isomorphisme de groupes alg\'ebriques. En effet, prenons par exemple $p=2$, 
${\bf T}=\Bbb{G}_{\rm m}\times \Bbb{G}_{\rm m}$ et $\tau(x,y)=(y,x)$. Puisque 
$\tau^2={\rm id}_{\bf T}$, $\tau$ est unipotent. Pour $(x,y)\in \smash{\overline{F}}^\times 
\times \smash{\overline{F}}^\times$, notant $z\in \smash{\overline{F}}^\times$ 
l'\'el\'ement tel que $xy=z^2$, on a
$$(x,y)=(z,z)(z^{-1}x,zx^{-1})\in {\bf T}_{\tau}{\bf T}(1-\tau).$$ 
Mais puisque ${\rm Lie}({\bf T}_\tau)= {\rm Lie}({\bf T}(1-\tau))=\ker (1+{\rm Lie}(\tau))$ 
est strictement contenu dans ${\rm Lie}({\bf T})$, ${\rm Lie}(\mu)$ n'est pas bijectif, et 
$\mu$ n'est pas s\'eparable.\hfill $\blacksquare$}
\end{marema}

\subsection{Groupes r\'eductifs connexes}\label{groupes rductifs} 
{\it On suppose d\'esormais, et jusqu'\`a la fin du ch.~3,  que le groupe ${\bf H}$ est 
connexe et r\'eductif}. On note 
${\bf H}_{\rm der}$ le groupe d\'eriv\'e de ${\bf H}$, et 
$C({\bf H})$ le cocentre ${\bf H}/{\bf H}_{\rm der}$ de ${\bf H}$. Rappelons que l'on 
a pos\'e  $\frak{H}={\rm Lie}({\bf H})$ et $\overline{\bf H}= {\bf H}/Z({\bf H})$. On pose aussi 
$\frak{H}_{\rm der}={\rm Lie}({\bf H}_{\rm der})$, 
$\frak{Z}={\rm Lie}(R({\bf H}))$, $\frak{C}={\rm Lie}(C({\bf H}))$ et 
$\smash{\overline{\bf H}}_{\rm der}={\bf H}_{\rm der}/Z({\bf H}_{\rm der})$. 

Le radical $R({\bf H})$ co\"{\i}ncide avec le tore $Z({\bf H})^\circ$ \cite[ch.~IV, 11.21]{Bor}. Si ${\bf T}$ est un 
tore maximal de ${\bf H}$, alors ${\bf T}'= ({\bf T}\cap {\bf H}_{\rm der})^\circ$ est un tore maximal de ${\bf H}_{\rm der}$, 
et le morphisme produit
$$
R({\bf H})\times {\bf T}'\rightarrow {\bf H}_{\rm der}
$$
est surjectif, de noyau fini \cite[ch.~V, 21.1]{Bor}. On en dduit que 
le morphisme produit
$$\rho: R({\bf H}) \times {\bf H}_{\rm der}\rightarrow {\bf H}$$
est lui aussi surjectif, de noyau fini. D'apr\`es 
\cite[ch.~V, 22,4]{Bor} $\rho$ est une 
{\it isog\'enie centrale}, i.e. il 
v\'erifie les deux propri\'et\'es suivantes:
\begin{itemize}
\item $\ker(\rho)$ est fini et central dans 
$R({\bf H}) \times {\bf H}_{\rm der}$;
\item $\ker({\rm Lie}(\rho))$ est central 
dans $\frak{Z}\times \frak{H}_{\rm der}\;(={\rm Lie}(R({\bf H}) \times {\bf H}_{\rm der}))$.
\end{itemize}
Notons qu'en g\'en\'eral, $\rho$ {\it n'est pas} s\'eparable. 

Puisque ${\bf H}= Z({\bf H}){\bf H}_{\rm der}$ et 
${\bf H}_{\rm der}\cap Z({\bf H})=Z({\bf H}_{\rm der})$, 
l'inclusion ${\bf H}_{\rm der}\subset {\bf H}$ induit par passage aux quotients une 
identification $\smash{\overline{\bf H}}_{\rm der}=\overline{\bf H}$. Notons aussi 
que la repr\'esentation adjointe ${\rm Ad}_{\bf H}$ induit par restriction une isog\'enie centrale 
$$
{\bf H}_{\rm der}\rightarrow {\bf H}_{\rm ad}.
$$

Si ${\bf T}$ est un tore maximal de ${\bf H}$, le morphisme produit
$$
{\bf T}\times {\bf H}_{\rm der}\rightarrow {\bf H}
$$
est surjectif, et {\it sparable} \cite[ch.~V, 22.5]{Bor}. On en dduit que l'inclusion 
l'inclusion ${\bf T}\subset {\bf H}$ induit par passage aux quotients un 
isomorphisme de groupes alg\'ebriques 
$$
{\bf T}/({\bf T}\cap {\bf H}_{\rm der})\buildrel\simeq\over{\longrightarrow} C({\bf H}).
$$
Par suite le cocentre $C({\bf H})$ est un groupe connexe et diagonalisable, \cad un tore. On a l'inclusion 
$\frak{Z}+\frak{H}_{\rm der}\subset \frak{H}$ 
avec \'egalit\'e si et seulement si $\rho$ est s\'eparable, auquel cas on a la 
d\'ecomposition $\frak{H}=\frak{Z}\oplus \frak{H}_{\rm der}$. 
En g\'en\'eral (i.e. que $\rho$ soit ou non s\'eparable), 
si ${\bf T}$ est un tore maximal de ${\bf H}$, d'apr\`es loc.~cit. 
on a l'\'egalit\'e
$$
\frak{H}={\rm Lie}({\bf T})+\frak{H}_{\rm der}.\leqno{(*)}
$$

Soit $\tau\in {\rm Aut}_{\overline{F}}({\bf H})$. Alors 
$\tau$ induit par restriction des 
$\overline{F}$--automorphismes 
$\tau_{\rm cent}$ de $R({\bf H})$ et 
$\tau_{\rm der}$ de ${\bf H}_{\rm der}$. 
Puisque ${\bf H}_{\rm der}$ est semisimple, 
le groupe ${\rm Out}_{\overline{F}}({\bf H}_{\rm der})$ est fini, 
par cons\'equent 
on a l'\'egalit\'e ${\rm Aut}_{\overline{F}}({\bf H}_{\rm der})
={\rm Aut}_{\overline{F}}^0({\bf H}_{\rm der})$. D'apr\`es la remarque 2 de \ref{automorphismes ss et u}, les trois conditions suivantes 
sont \'equivalentes:
\begin{itemize}
\item $\tau$ est localement fini;
\item $\tau\in {\rm Aut}_{\overline{F}}^0({\bf H})$;
\item $\tau_{\rm cent}$ est d'ordre fini.
\end{itemize}
Posons 
${\rm Lie}(\tau)_{\rm cent}={\rm Lie}(\tau_{\rm cent})\in {\rm GL}(\frak{Z})$ et 
${\rm Lie}(\tau)_{\rm der}={\rm Lie}(\tau_{\rm der})\in {\rm GL}(\frak{H}_{\rm der})$. Alors pour 
$z\in Z({\bf H})$ et $h'\in {\bf H}_{\rm der}$, 
on a
$$({\rm Ad}_{\bf H}(zh)\circ{\rm Lie}(\tau))_{\rm cent}= {\rm Lie}(\tau)_{\rm cent}
$$
et
$$
({\rm Ad}_{\bf H}(zh')\circ{\rm Lie}(\tau))_{\rm der}= 
{\rm Ad}_{{\bf H}_{\rm der}}(h')\circ{\rm Lie}(\tau)_{\rm der}.
$$
Notons que si 
$\rho$ est s\'eparable, alors on a la d\'ecomposition 
${\rm Lie}(\tau) = {\rm Lie}(\tau)_{\rm cent}+
{\rm Lie}(\tau)_{\rm der}$.

Puisque 
le morphisme quotient $\pi:{\bf H}\rightarrow C({\bf H})$ est s\'eparable, le morphisme 
${\rm Lie}(\pi)$ induit une identification $\frak{H}/\frak{H}_{\rm der}= \frak{C}$. Notons 
$\tau_{\rm cocent}$ le $\overline{F}$--automorphisme de $C({\bf H})$ d\'eduit 
de $\tau$ par passage au quotient, et posons ${\rm Lie}(\tau)_{\rm cocent}=
{\rm Lie}(\tau_{\rm cocent})$; c'est l'\'el\'ement de ${\rm GL}(\frak{C})$ 
d\'eduit de ${\rm Lie}(\tau)$ par passage au quotient. Pour $h\in {\bf H}$, 
notons $c_h:{\bf H}\rightarrow {\bf H}$ l'application commutateur 
$x\mapsto hxh^{-1}x^{-1}$; c'est un morphisme de vari\'et\'es alg\'ebriques, 
dont l'image est contenue dans ${\bf H}_{\rm der}$. Puisque sa diff\'erentielle 
au point $1$ est donn\'ee par ${\rm d}(c_h)_1={\rm Ad}_{\bf H}(h)-{\rm id}_\frak{H}$, on a l'inclusion $({\rm Ad}_{\bf H}(h)-{\rm id}_\frak{H})(\frak{H})
\subset \frak{H}_{\rm der}$, d'o\`u l'\'egalit\'e
$$
({\rm Ad}_{\bf H}(h)\circ{\rm Lie}(\tau))_{\rm cocent}= {\rm Lie}(\tau)_{\rm cocent}.
$$

Soit $\Bbb{P}$ le sous--groupe parabolique de ${\rm GL}(\frak{H})$ form\'e 
des \'el\'ements $\Phi$ tels que $\Phi(\frak{H}_{\rm der})= \frak{H}_{\rm der}$. 
Alors $\Bbb{U}=R_{\rm u}(\Bbb{P})$ est donn\'e par
$$
\Bbb{U}=\{\Phi\in \Bbb{P}: \hbox{$(\Phi-{\rm id}_\frak{H})(\frak{H})\subset \frak{H}_{\rm der}$ et 
$(\Phi-{\rm id}_\frak{H})(\frak{H}_{\rm der})=0$}\},
$$
et $\Bbb{P}/\Bbb{U}= {\rm GL}(\frak{H}_{\rm der})\times {\rm GL}(\frak{C})$. 
De plus, ${\rm Lie}(\tau)$ appartient \`a $\Bbb{P}$, et 
$({\rm Lie}(\tau)_{\rm der},{\rm Lie}(\tau)_{\rm cocent})$ 
co\"{\i}ncide avec l'image de ${\rm Lie}(\tau)$ par le morphisme quotient 
$\Bbb{P}\rightarrow \Bbb{P}/\Bbb{U}$.

\begin{monlem}
Soit $\tau\in {\rm Aut}_{\overline{F}}({\bf H})$. On a 
$({\bf H}_{\rm der})_{\tau_{\rm der}}^\circ =({\bf H}_\tau^\circ\cap {\bf H}_{\rm der})^\circ$ et 
${\bf H}_\tau^\circ = R({\bf H})_\tau^\circ({\bf H}_{\rm der})_{\tau}^\circ$.
\end{monlem}

\begin{proof}
La premi\`ere \'egalit\'e est claire, puisque ${\bf H}_{\tau_{\rm der}}= {\bf H}_{\tau}\cap {\bf H}_{\rm der}$. 
Quant \`a la seconde \'egalit\'e, posons 
${\bf R}'={\bf H}_{\rm der}\cap R({\bf H})\subset Z({\bf H}_{\rm der})$, et soit $h\in {\bf H}_\tau$. \'Ecrivons $h=zh'$ avec $z\in R({\bf H})$ et 
$h'\in {\bf H}_{\rm der}$. Puisque $\tau(h)=h$, on a 
$h'^{-1}\tau(h')= z\tau(z)^{-1}\in {\bf R}'$ et 
l'image de $z\tau(z)^{-1}$ dans le groupe quotient ${\bf R}'\!/ {\bf R}'(1-\tau)$ ne d\'epend 
pas de la d\'ecomposition $h=zh'$ choisie; on note $\bar{z}_h$ cette image. L'application
$$
{\bf H}_\tau\mapsto {\bf R}'\!/{\bf R}'(1-\tau),\,h\mapsto \bar{z}_h
$$
ainsi d\'efinie, est un morphisme de groupes alg\'ebriques de noyau 
$R({\bf H})_\tau ({\bf H}_{\rm der})_{\tau}$. Comme le groupe quotient 
${\bf R}'\!/{\bf R}'(1-\tau)$ est fini, 
$R({\bf H})_\tau({\bf H}_{\rm der})_{\tau}$ est un sous--groupe ferm\'e d'indice fini de 
${\bf H}_\tau$. D'o\`u l'inclusion \cite[ch.~I, 1.2]{Bor} 
${\bf H}_\tau^\circ \subset R({\bf H})_\tau({\bf H}_{\rm der})_{\tau}$, puis 
l'\'egalit\'e ${\bf H}_\tau^\circ = [R({\bf H})_\tau({\bf H}_{\rm der})_{\tau}]^\circ$. 
Comme $[R({\bf H})_\tau({\bf H}_{\rm der})_{\tau}]^\circ$ co\"{\i}ncide 
avec l'image de $[R({\bf H})_\tau\times ({\bf H}_{\rm der})_\tau]^\circ = R({\bf H})_\tau^\circ \times ({\bf H}_{\rm der})_\tau^\circ$ 
par le morphisme produit $R({\bf H})_\tau\times ({\bf H}_{\rm der})_\tau\rightarrow {\bf H}_\tau$, 
on a bien l'galit ${\bf H}_\tau^\circ = 
R({\bf H})_\tau^\circ ({\bf H}_{\rm der})_{\tau}^\circ$.
\end{proof}

Choisissons un sous--groupe 
de Borel ${\bf B}$ de ${\bf H}$ et un tore maximal ${\bf T}$ de ${\bf B}$. Notons 
$\Phi=\Phi({\bf T},{\bf H})$ l'ensemble des racines de ${\bf T}$ dans ${\bf H}$, 
$\Phi^+ = \Phi({\bf T},{\bf B})$ le sous--ensemble de $\Phi$ form\'e des racines de 
${\bf T}$ dans ${\bf B}$, et $\Delta=\Delta ({\bf T},{\bf B})$ la base de $\Phi$ form\'ee 
des racines simples de $\Phi^+$. Pour $\alpha\in \Phi$, on note 
${\bf U}_\alpha$ le sous--groupe unipotent de ${\bf H}$ associ\'e \`a la racine 
$\alpha$. Pour chaque $\alpha\in \Delta$, 
choisissons un \'el\'ement $u_\alpha\in {\bf U}_\alpha\smallsetminus \{1\}$. Alors ${\rm Out}_{\overline{F}}({\bf H})$ est isomorphe au sous--groupe 
$\frak{A}={\rm Aut}_{\overline{F}}({\bf H},{\bf B},{\bf T}, \{u_\alpha\}_{\alpha\in \Delta})$ de ${\rm Aut}_{\overline{F}}({\bf H})$ form\'e des 
$\overline{F}$--automorphismes qui stabilisent ${\bf B}$, ${\bf T}$ et la famille 
$\{u_\alpha\}_{\alpha\in \Delta}$ \cite[prop. 2.13]{Sp}; i.e. on a la d\'ecomposition en 
produit semidirect
$$
{\rm Aut}_{\overline{F}}({\bf H})={\rm Int}_{\overline{F}}({\bf H})\rtimes \frak{A}.\leqno{(**)}
$$
Pour $\alpha\in \Delta$, notons $e_\alpha:\Bbb{G}_{\rm a} \rightarrow {\bf U}_\alpha$ 
l'{\it \'epinglage 
de ${\bf U}_\alpha$ d\'efini par $u_\alpha$}, i.e.
l'unique isomorphisme de groupes alg\'ebriques tel que 
$$
\left\{
\begin{array}{ll}
te_\alpha(x)t^{-1}= e_\alpha(\alpha(t)x)\hfill & (x\in \overline{F},\,t\in {\bf T})\\
e_\alpha(1)=u_\alpha\hfill & \\
\end{array}\right..
$$
Si $\tau\in {\rm Aut}_{\overline{F}}({\bf H})$ stabilise la paire $({\bf B},{\bf T})$, 
alors $\tau$ op\`ere par permutation sur l'ensemble $\Delta$, et pour chaque 
$\alpha\in \Delta$, il existe un (unique) $y_{\tau,\alpha}\in \overline{F}^\times$ tel que 
$$\tau\circ e_\alpha (x) = e_{\tau(\alpha)}(y_{\tau,\alpha}x)\quad (x\in \overline{F});
$$
par suite $\tau\in \frak{A}$ si et seulement si $y_{\tau,\alpha}=1$ ($\alpha\in \Delta$). On a 
donc $\frak{A}= {\rm Aut}_{\overline{F}}({\bf H},{\bf B},{\bf T}, \{e_\alpha\}_{\alpha\in \Delta})$. Choisissons 
une autre famille $\{u_{1,\alpha}\in {\bf U}_\alpha\smallsetminus\{1\}:\alpha\in \Delta\}$, et 
notons $\frak{A}_1$ et $e_{1,\alpha}$ ($\alpha\in \Delta$) les objets d\'efinis 
comme ci-dessus en rempla\c{c}ant $\{u_\alpha\}_{\alpha\in \Delta}$ par 
$\{u_{1,\alpha}\}_{\alpha\in \Delta}$. Pour chaque $\alpha\in \Delta$, il existe un unique lment $y_\alpha\in \overline{F}^\times$ 
tel que $e_{1,\alpha} (t) = e_{\alpha}(y_\alpha t)$. Choisissons 
un $t_1\in {\bf T}$ tel que pour chaque $\alpha \in \Delta$, on ait 
$\alpha(t_1)=y_\alpha$. Alors on a ${\rm Int}_{\bf H}(t_1)\circ e_\alpha = e_{1,\alpha}$ ($\alpha\in \Delta$). 
D'o\`u un isomorphisme de groupes
$$\frak{A}\rightarrow \frak{A}_1,\,\tau \mapsto \tau_1={\rm Int}_{\bf H}(t_1)\circ \tau \circ {\rm Int}_{\bf H}(t_1)^{-1},
$$
qui permet de passer de la d\'ecomposition $(**)$ \`a la d\'ecomposition 
$${\rm Aut}_{\overline{F}}({\bf H})={\rm Int}_{\overline{F}}({\bf H})\rtimes \frak{A}_1.$$
Prcisment, pour 
$h\in {\bf H}$ et $\tau\in \frak{A}$, on a 
$$
{\rm Int}_{\bf H}(h)\circ \tau = {\rm Int}_{\bf H}(h \tau(t_1)t_1^{-1})\circ \tau_1.
$$

\subsection{Rev\^etement universel}\label{rev\^etement universel}
Soit
$$
\pi:{\bf H}_{\rm SC}\rightarrow \bf{H}'={\bf H}_{\rm der}
$$ le 
rev\^etement universel de ${\bf H}_{\rm der}$ \cite[2.6.1]{T}: le groupe ${\bf H}_{\rm SC}$ 
est semisimple et simplement connexe, 
et $\pi$ est une isog\'enie centrale. Posons ${\bf B}'= {\bf B}
\cap {\bf H}'$ et ${\bf T}'= {\bf T}
\cap {\bf H}'$; ce sont respec\-tivement un sous--groupe de Borel et un 
tore maximal de ${\bf H}'$. 
D'apr\`es  \cite[ch.~V, 22.6]{Bor}, les images r\'eciproques ${\bf B}_{\rm sc}$ de 
${\bf B}'$, et ${\bf T}_{\rm sc}$ de ${\bf T}'$, par $\pi$ sont respectivement 
un sous--groupe de Borel et un tore maximal de ${\bf H}_{\rm SC}$. De plus 
(loc. cit.), l'application
$${\rm X}^*({\bf T}')\rightarrow {\rm X}^*({\bf T}_{\rm sc}),\,
\chi\mapsto \pi^\sharp(\chi)= \chi\circ(\pi\vert_{{\bf T}_{\rm sc}})
$$
induit une application bijective
$$\Phi'=\Phi({\bf T}'\! ,{\bf H}')\rightarrow 
\Phi({\bf T}_{\rm sc},{\bf H}_{\rm SC})=\widetilde{\Phi},$$ 
et pour $\alpha\in \Phi'$, $\pi$ induit par restriction un isomorphisme de groupes 
alg\'ebriques
$${\bf U}_{\pi^\sharp(\alpha)}\buildrel \simeq\over{\longrightarrow} {\bf U}_\alpha.
$$
 D'autre part, l'application 
${\rm X}^*({\bf T})\rightarrow {\rm X}^*({\bf T}') ,\,\chi\mapsto\chi'=\chi\vert_{{\bf T}'}$ 
induit une application bijective $\Phi\rightarrow \Phi'$, et pour $\alpha\in \Phi$, 
on a ${\bf U}_\alpha = {\bf U}_{\alpha'}\subset {\bf H}'$. On peut donc identifier 
$\Phi$ et $\Phi'$. Pour chaque $\alpha\in \Delta$, posons 
$\tilde{u}_\alpha = \pi^{-1}(u_\alpha)\in {\bf U}_{\pi^\sharp(\alpha)}$. Comme 
l'ensemble $\{\pi^\sharp(\alpha):\alpha\in \Delta\}$ co\"{\i}ncide avec 
la base $\Delta({\bf T}_{\rm sc},{\bf B}_{\rm sc})$ de $\widetilde{\Phi}$ associ\'ee 
\`a ${\bf B}_{\rm sc}$, on peut poser 
$\widetilde{\frak{A}} ={\rm Aut}_{\overline{F}}({\bf H}_{\rm SC},{\bf B}_{\rm sc},
{\bf T}_{\rm sc}, \{\tilde{u}_\alpha\}_{\alpha\in \Delta})$. On pose aussi 
$\frak{A}' ={\rm Aut}_{\overline{F}}({\bf H}',{\bf B}',
{\bf T}', \{u_\alpha\}_{\alpha\in \Delta})$. 

Le lemme suivant est d\^u  Steinberg \cite[9.16]{St}.

\begin{monlem}
Soit $\tau\in {\rm Aut}_{\overline{F}}({\bf H}_{\rm der})$. Il existe un  
unique $\tilde\tau\in {\rm Aut}_{\overline{F}}({\bf H}_{\rm SC})$ relevant 
$\tau$, i.e. v\'erifiant $\pi\circ \tau= \tilde{\tau}\circ\pi$.
\end{monlem}

L'unicit\'e du rel\`evement entra\^{\i}ne que l'application
$${\rm Aut}_{\overline{F}}({\bf H}_{\rm der})\rightarrow {\rm Aut}_{\overline{F}}({\bf H}_{\rm SC}),\,
\tau\mapsto \tilde{\tau}\leqno{(*)}$$
est un morphisme de groupes. 
Si $\tau\in {\rm Aut}_{\overline{F}}({\bf H}_{\rm der})$ stabilise la paire $({\bf B}',{\bf T}')$, alors par construction, 
son rel\`evement $\tilde{\tau}$ \`a ${\bf H}_{\rm SC}$ stabilise la paire $({\bf B}_{\rm sc},{\bf T}_{\rm sc})$; 
par cons\'equent l'application $(*)$ induit un isomorphisme de groupes 
$\frak{A}'\rightarrow \widetilde{\frak{A}}$. D'autre part $\pi$ induit par passage aux quotients un morphisme bijectif de groupes alg\'ebriques 
$\overline{\pi}: {\bf H}_{\rm SC}/Z({\bf H}_{\rm SC})\rightarrow {\bf H}_{\rm der}/Z({\bf H}_{\rm der})$, 
qui n'est en g\'en\'eral pas un isomorphisme. On en d\'eduit que l'application 
$(*)$ induit aussi un isomorphisme de groupes
${\rm Int}_{\overline{F}}({\bf H})\rightarrow
 {\rm Int}_{\overline{F}}({\bf H}_{\rm SC})$. L'application $(*)$  
 est donc un isomorphisme de groupes, qui pr\'eserve la 
 d\'ecomposition $(**)$ de \ref{groupes rductifs} pour ${\bf H}_{\rm der}$ et pour ${\bf H}_{\rm SC}$. 

\begin{exemple}{\rm 
Soit ${\bf H}\;(={\bf H}_{\rm der})$ le groupe 
$\Bbb{PGL}_2=\Bbb{GL}_2/{\bf Z}$, o\`u ${\bf Z}=Z(\Bbb{GL}_2)\simeq \Bbb{G}_{\rm m}$. 
Alors le rev\^etement universel $\pi: {\bf H}_{\rm SC}=\Bbb{SL}_2\rightarrow {\bf H}$ est le compos\'e de l'inclusion $\Bbb{SL}_2\subset \Bbb{GL}_2$ et du morphisme quotient 
$\Bbb{GL}_2\rightarrow \Bbb{PGL}_2$. Si $F$ est de caract\'eristique $2$, le morphisme $\pi$ est 
bijectif (on a $\pi=\overline{\pi}$) mais il n'est pas s\'eparable --- cf. \cite[ch.~III, 10.8, rem. p.~144]{Bor}.
\hfill $\blacksquare$}
\end{exemple}

\subsection{Automorphismes quasi--semisimples}\label{automorphismes qss} 
Commen\c{c}ons par le r\'esultat bien connu suivant, d\^u \`a Steinberg \cite[7.2]{St}:

\begin{montheo1}
Tout $\overline{F}$--automorphisme de ${\bf H}$ stabilise 
un sous--groupe de Borel de ${\bf H}$.
\end{montheo1}

On appelle {\it paire de Borel de ${\bf H}$} une 
paire $({\bf B},{\bf T})$ form\'ee d'un sous--groupe de Borel ${\bf B}$ de ${\bf H}$ et 
d'un tore maximal ${\bf T}$ de ${\bf B}$. Le groupe 
${\rm Aut}_{\overline{F}}({\bf H})$ op\`ere naturellement sur l'ensemble des paires de Borel 
de ${\bf H}$: si $\tau\in {\rm Aut}_{\overline{F}}({\bf H})$ et 
$({\bf B},{\bf T})$ est une paire de Borel de ${\bf H}$, on pose
 $\tau({\bf B},{\bf T})=(\tau({\bf B}),\tau({\bf T}))$. Rappelons que toutes les 
paires de Borel de ${\bf H}$ sont dans la m\^eme orbite sous 
${\rm Int}_{\overline{F}}({\bf H})$.
 
Un $\overline{F}$--automorphisme 
de ${\bf H}$ est dit {\it quasi--semisimple} s'il stabilise une paire de 
Borel de ${\bf H}$. Si $\tau\in {\rm Aut}_{\overline{F}}({\bf H})$, 
alors $\tau$ est quasi--semisimple si et seulement 
si $\tau_{\rm der}$ est quasi--semisimple. 

Si $h\in {\bf H}$, 
alors ${\rm Int}_{\bf H}(h)$ est quasi--semisimple si et seulement si 
$h$ est semisimple.

Si $\tau\in {\rm Aut}_{\overline{F}}({\bf H})$ est quasi--semisimple, un tore maximal 
$\tau$--stable ${\bf T}$ de ${\bf H}$ est dit {\it $\tau$--admissible} s'il est contenu 
dans un sous--groupe de Borel $\tau$--stable de ${\bf H}$. Si $\tau$ stabilise une 
paire de Borel $({\bf B},{\bf T})$ de ${\bf H}$, alors $\tau$ stabilise aussi la paire de 
Borel $({\bf B}^-,{\bf T})$ de ${\bf H}$ {\it oppos\'ee \`a $({\bf B},{\bf T})$}, o\`u ${\bf B}^-$ 
est l'unique sous--groupe de Borel contenant ${\bf T}$ tel que ${\bf B}^-\cap {\bf B}={\bf T}$. 

Soit $\tau\in {\rm Aut}_{\overline{F}}({\bf H})$. Alors il existe un 
$h\in {\bf H}$ tel que le $\overline{F}$--automorphisme 
${\rm Int}_{\bf H}(h)\circ\tau$ de ${\bf H}$ est quasi--semisimple. Plus 
pr\'ecis\'ement, on a le

\begin{monlem}
Soit $\tau\in {\rm Aut}_{\overline{F}}({\bf H})$, et soit 
$({\bf B},{\bf T})$ une paire de Borel de ${\bf H}$ telle que $\tau({\bf B})={\bf B}$. 
Alors il existe un unique $u\in R_{\rm u}({\bf B})$ tel que ${\rm Int}_{\bf H}(u)\circ \tau
({\bf B},{\bf T})=({\bf B},{\bf T})$.
\end{monlem}

\begin{proof} Puisque $\tau({\bf T})$ est un tore maximal de  
${\bf B}$, il existe un unique $u\in R_{\rm u}({\bf B})$ 
tel que $\tau({\bf T})={\rm Int}_{\bf H}(u^{-1})({\bf T})$. D'o\`u le lemme.
\end{proof}

D'aprs le thorme 1.8 et les propositions 1.11 et 1.12 de \cite{DM1}, on a le

\begin{montheo2}
Soit $\tau$ un $\overline{F}$--automorphisme quasi--semi\-simple de 
${\bf H}$.
\begin{enumerate}
\item[(1)]Le groupe ${\bf H}_\tau^\circ$ est r\'eductif, et le groupe 
${\bf H}_\tau/{\bf H}_\tau^\circ$ est form\'e d'\'el\'ements semisimples.
\item[(2)]Soit ${\bf P}$ un sous--groupe parabolique $\tau$--stable de ${\bf H}$. Il existe une 
composante de Levi de ${\bf P}$ qui soit $\tau$--stable. 
\item[(3)]Soit ${\bf P}$ un sous--groupe parabolique 
$\tau$--stable de ${\bf H}$, et 
${\bf L}$ une composante de Levi $\tau$--stable de ${\bf P}$. Posons 
${\bf U}=R_{\rm u}({\bf P})$. Alors 
${\bf U}_\tau$ est connexe, 
${\bf P}^\sharp= {\bf P}\cap {\bf H}_\tau^\circ$ est un sous--groupe 
parabolique de ${\bf H}_\tau^\circ$ de radical unipotent ${\bf U}_\tau $, 
et 
${\bf L}^\sharp={\bf L}\cap {\bf H}_\tau^\circ$ est une composante 
de Levi de ${\bf P}^\sharp$.
\item[(4)]Soit $({\bf B},{\bf T})$ une paire de Borel $\tau$--stable 
de ${\bf H}$. Alors $({\bf B}^\sharp,{\bf T}^\sharp)=({\bf B}\cap {\bf H}_\tau^\circ, 
{\bf T}\cap {\bf H}_\tau^\circ)$ est une paire de Borel de ${\bf H}_\tau^\circ$. R\'eciproquement, si 
$({\bf B}^\sharp,{\bf T}^\sharp)$ est une paire de Borel de ${\bf H}_\tau^\circ$, 
alors ${\bf T}={\bf Z}_{\bf H}({\bf T}^\sharp)$ est un tore maximal $\tau$--stable 
de ${\bf H}$, et il existe un sous--groupe de Borel $\tau$--stable ${\bf B}$ de 
${\bf H}$ tel que ${\bf B}^\sharp\subset {\bf B}$ et ${\bf T}\subset {\bf B}$.
\item[(5)]Soit ${\bf H}'$ un sous--groupe ferm\'e r\'eductif connexe de rang maximal 
de ${\bf H}$. Si ${\bf H}'$ est $\tau$--stable et si ${\bf H}'^\circ_\tau$ est un sous--groupe 
de rang maximal de ${\bf H}_\tau^\circ$, alors la restriction de $\tau$ \`a ${\bf H}'$ est encore 
quasi--semisimple.
\end{enumerate}
\end{montheo2}

\begin{marema1}
{\rm D'apr\`es le point (1), tout \'el\'ement unipotent de ${\bf H}_\tau$ 
appartient \`a ${\bf H}_\tau^\circ$. En particulier si $p>1$, alors l'ordre du groupe fini 
${\bf H}_\tau/{\bf H}_\tau^\circ$ est premier \`a $p$.

D'apr\`es le point (4), 
l'application ${\bf T}^\sharp \mapsto {\bf Z}_{\bf H}({\bf T}^\sharp)$ est une bijection 
entre l'ensemble des tores maximaux de ${\bf H}_\tau^\circ$ et l'ensemble des 
tores maximaux $\tau$--admissibles de ${\bf H}$. En particulier, si ${\bf H}_\tau^\circ$ 
est central dans ${\bf H}$, alors ${\bf H}$ est un tore.

D'apr\`es le point (3), le point (5) s'applique en particulier \`a toute composante de Levi $\tau$--stable 
${\bf H}'={\bf L}$ d'un 
sous--groupe parabolique $\tau$--stable ${\bf P}$ de ${\bf H}$. 
Notons aussi que puisque ${\bf L}\cap {\bf H}_\tau^\circ$ est connexe, 
on a ${\bf L}_\tau^\circ = {\bf L}\cap {\bf H}_\tau^\circ$ et 
${\bf P}_\tau^\circ = {\bf P}\cap {\bf H}_\tau^\circ$.\hfill $\blacksquare$
}
\end{marema1}

\begin{marema2}{\rm Supposons le groupe ${\bf H}$ semisimple, et soit $\tau$ un 
$\overline{F}$--automorphisme quasi--semisimple de ${\bf H}$. D'apr\`es \cite[cor.~9.4]{St} le groupe ${\bf H}_\tau/{\bf H}_\tau^\circ$ est 
ab\'elien. Si de plus ${\bf H}$ est simplement connexe, 
alors d'apr\`es \cite[theo.~8.2]{St}, le groupe ${\bf H}_\tau$ est connexe. \hfill $\blacksquare$
}
\end{marema2}

\begin{mapropo}
Soit $\tau\in {\rm Aut}_{\overline{F}}({\bf H})$ quasi--semisimple, 
$({\bf B},{\bf T})$ une paire de Borel $\tau$--stable de ${\bf H}$, et ${\bf U}=R_{\rm u}({\bf B})$. 
Le morphisme ${\bf U}\rightarrow {\bf U},\,u\mapsto u\tau(u)^{-1}$ est s\'eparable.
\end{mapropo}

\begin{proof}
L'application ${\bf U}\times {\bf U}\rightarrow {\bf U},\,(u,v)\mapsto u\cdot v=uv\tau(u)^{-1}$ 
est une action alg\'ebrique de ${\bf U}$ sur lui-m\^eme. Posant 
$\pi_v(u)=u\cdot v$ ($u,\,v\in {\bf U}$), il s'agit de montrer que pour 
$v=1$, le morphisme $\pi_1=(1-\tau)\vert_{\bf U}:{\bf U}\rightarrow {\bf U}$ est s\'eparable. 
En d'autres termes (cf. \ref{groupes algbriques}), il s'agit d'\'etablir l'\'egalit\'e 
${\rm Lie}({\bf U}_\tau)= \ker ({\rm id}-{\rm Lie}(\tau);{\rm Lie}({\bf U}))$.

Soit $\pi:{\bf H}_{\rm SC}\rightarrow {\bf H}_{\rm der}$ le rev\^etement universel 
de ${\bf H}_{\rm der}$, et soit $\tilde{\tau}\in {\rm Aut}_{\overline{F}}({\bf H}_{\rm SC})$ 
le rel\`evement de $\tau$ \`a ${\bf H}_{\rm SC}$ (cf. \ref{rev\^etement universel}). Notons ${\bf B}_{\rm sc}$, ${\bf T}_{\rm sc}$ et 
${\bf U}_{\rm sc}$ les images r\'eciproques de ${\bf B}\cap {\bf H}_{\rm der}$, ${\bf T}\cap {\bf H}_{\rm der}$ et 
${\bf U}$ par $\pi$. Alors on a (cf. \ref{rev\^etement universel}):
\begin{itemize}
\item $\tilde{\tau}$ stabilise la paire 
de Borel $({\bf B}_{\rm sc},{\bf T}_{\rm sc})$ de ${\bf H}_{\rm SC}$; 
\item $\pi$ induit un isomorphisme de groupes alg\'ebriques 
${\bf U}_{\rm sc}\rightarrow {\bf U}$;
\item pour $\tilde{u}\in {\bf U}_{\rm sc}$ et $u=\pi(\tilde{u})\in {\bf U}$, on a $\pi(\tilde{u}\tilde{\tau}(\tilde{u})^{-1})=
u\tau(u)^{-1}$.
\end{itemize}
On peut donc supposer ${\bf H}$ semisimple et simplement connexe. Puisque $\tau$ 
est quasi--semisimple, le groupe ${\bf H}_\tau$ est connexe (remarque 2). En particulier, 
${\bf S}={\bf T}_\tau \;(={\bf T}\cap {\bf H}_\tau)$ est un tore maximal de ${\bf H}_\tau$. 
D'apr\`es \cite[3.1.2]{T}, ${\bf H}$ se d\'ecompose en un produit direct 
${\bf H}={\bf H}_1\times \cdots \times {\bf H}_n$ o\`u chaque ${\bf H}_i$ est 
un groupe semisimple simplement connexe et presque simple, et cette d\'ecomposition 
est unique \`a permutation des ${\bf H}_i$ pr\`es. L'unicit\'e de la d\'ecomposition 
implique que $\tau$ permute les facteurs ${\bf H}_i$. On peut donc supposer que $\tau$ permute transitivement les facteurs 
${\bf H}_i$, et quitte à r\'eordonner ces facteurs, que $\tau({\bf H}_{i+1})= {\bf H}_i$ pour $i=1, \ldots n-1$. Alors en identifiant 
${\bf H}_{i+1}$ à ${\bf H}_1$ via $\tau^i$ pour $i=1,\ldots , n-1$, l'automorphisme $\tau$ de ${\bf H}= {\bf H}_1\times \cdots \times {\bf H}_1$ est donné par
$$
\tau(x)= (x_2, \ldots , x_n , \tau_1(x_1)),\quad x=(x_1, \ldots ,x_n)\in {\bf H},
$$
où $\tau_1\in {\rm Aut}_{\overline{F}}({\bf H}_1)$  est donné par la restriction de $\tau^n$ à ${\bf H}_1$. On peut donc supposer 
${\bf H}$ presque simple (et toujours semisimple simplement connexe), \cad $n=1$ et ${\bf H}={\bf H}_1$.

Posons $\Phi=\Phi({\bf T},{\bf H})$ et $\Phi^+=\Phi({\bf T},{\bf B})$ 
comme en \ref{groupes rductifs}. Pour chaque 
racine $\alpha\in \Phi$, notons $\Phi(\alpha)$ le sous--ensemble 
de $\Phi$ form\'e des racines $\beta$ dont la restriction \`a 
${\bf S}$ est proportionnelle \`a $\alpha\vert_{\bf S}$; i.e. telles que 
$\beta\vert_{\bf S}=c \alpha\vert_{\bf S}$ pour un élément $c\in {\Bbb R}_{>0}$
(on a donc $\Phi(\alpha)\subset \Bbb{R}\otimes_\Bbb{Z}{\rm X}^*({\bf S})$). 
On distingue deux cas (cf. la d\'emonstration de \cite[theo. 8.2]{St}):
\begin{enumerate}
\item[{\it I}:]$\Phi(\alpha)$ est une $\tau$--orbite, i.e. il existe un entier $l\geq 1$ 
tel que $\Phi(\alpha)=\{\alpha,\tau(\alpha),\ldots ,\tau^{l-1}(\alpha)\}$, $\tau^l(\alpha)=\alpha$, 
$\tau^i(\alpha) \neq\alpha $ 
et $\alpha + \tau^i(\alpha)\not\in \Phi$ pour $i=1,\ldots ,l-1$;
\item[{\it II}:] $\Phi(\alpha)=\{\alpha,\tau(\alpha),\alpha+\tau(\alpha)\}$, $\tau(\alpha)\neq \alpha$, 
$\tau^2(\alpha)=\alpha$ et $\alpha+ \tau(\alpha)\in \Phi$ 
(ce qui n'est possible que si $\Phi$ est de type $A_{2n}$).
\end{enumerate}
Notons ${\bf U}_{\Phi(\alpha)}$ le sous--groupe de ${\bf H}$ engendr\'e par les ${\bf U}_\beta$ pour $\beta\in \Phi(\alpha)$. Le groupe ${\bf U}$ est le produit direct, pris dans n'importe quel ordre, des groupes ${\bf U}_\alpha$ pour $\alpha\in \Phi^+$ \cite[ch.~IV, 14.4]{Bor}. 
On peut donc choisir un sous--ensemble $\Psi$ de $\Phi^+$ tel que ${\bf U}$ soit le produit direct, pris dans n'importe quel ordre, des groupes ${\bf U}_{\Phi(\alpha)}$ pour $\alpha\in \Psi$. D'apr\`es le th\'eor\`eme 2, le groupe ${\bf U}_\tau={\bf U}\cap {\bf H}_\tau$ 
est connexe et c'est le radical unipotent du sous--groupe de Borel ${\bf B}_\tau={\bf B}\cap {\bf H}_\tau$ de ${\bf H}_\tau$. 
Pour chaque racine $\alpha\in \Psi$, le groupe ${\bf U}_{\Phi(\alpha)}$ est 
$\tau$--stable, et ${\bf U}_\tau$ est le produit direct (pris dans n'importe quel ordre) des ${\bf U}_{\Phi(\alpha),\tau}={\bf U}_{\Phi(\alpha)}\cap {\bf U}_\tau$ pour $\alpha\in \Psi$. Le morphisme $1-\tau: {\bf U}\rightarrow {\bf U}$ induit, pour chaque $\alpha\in \Psi$, un morphisme $1-\tau:{\bf U}_{\Phi(\alpha)}\rightarrow {\bf U}_{\Phi(\alpha)}$, et il 
suffit de montrer que chacun des morphismes $(1-\tau)\vert_{{\bf U}_{\Phi(\alpha)}}$ est s\'eparable. 
Fixons donc une racine $\alpha\in \Psi$ et posons $\rho = (1-\tau)\vert_{{\bf U}_{\Phi(\alpha)}}$. 
Fixons aussi un \'epinglage 
$e_\alpha:\Bbb{G}_{\rm a}\rightarrow {\bf U}_\alpha$ de ${\bf U}_\alpha$ (cf. \ref{groupes rductifs}).

\v1
Commen\c{c}ons par le cas {\it I}. Pour 
$i=1,\ldots,l-1$, posons $e_{\tau^i(\alpha)}=\tau^i\circ e_\alpha: \Bbb{G}_{\rm a}\rightarrow {\bf U}_{\tau^i(\alpha)}$; c'est un \'epinglage de 
${\bf U}_{\tau^i(\alpha)}$. Puisque $\tau^l(\alpha)=\alpha$, il existe un $y\in \smash{\overline{F}}^\times$ tel que $\tau^l\circ e_\alpha (x)=e_\alpha(yx)$ pour tout $x\in \overline{F}$. 
Les groupes ${\bf U}_\beta$ 
($\beta\in \Phi(\alpha)$) commutent deux--\`a--deux, par suite $\rho $ est un 
morphisme de groupes, et pour $u=\prod_{i=0}^{l-1} e_{\tau^i(\alpha)}(x_i)\in {\bf U}_{\Phi(\alpha)}$, posant $x_{-1}=yx_{l-1}$, on a
$$
\rho(u)= \prod_{i=0}^{l-1} e_{\tau^i(\alpha)} (x_i-x_{i-1}).
$$
Le noyau $\ker \rho={\bf U}_{\Phi(\alpha),\tau}$ 
est non trivial si et seulement si $y=1$, auquel cas c'est l'ensemble 
des $\prod_{i=0}^{l-1}e_{\tau^i(\alpha)}(x)$ pour $x\in \overline{F}$. Un calcul analogue 
pour ${\rm Lie}(\rho)$ montre que
$${\rm Lie}(\ker (\rho))=\ker({\rm Lie}(\rho)).
$$ Ainsi $\rho$ est 
s\'eparable, et c'est un automorphisme si (et seulement si) $y\neq 1$. 

\v1
Traitons maintenant le cas {\it II}. Posons $e_{\tau(\alpha)}=\tau\circ e_\alpha:
\Bbb{G}_{\rm a}\rightarrow {\bf U}_{\tau(\alpha)}$, $\beta=\alpha +\tau(\alpha)$, et soit 
$e_\beta:\Bbb{G}_{\rm a}\rightarrow {\bf U}_\beta$ un \'epinglage de 
${\bf U}_\beta$. Un calcul simple dans $SL(3,\overline{F})$ montre qu'il existe un $u\in \smash{\overline{F}}^\times$ 
tel que $(e_\alpha(x),e_{\tau(\alpha)}(x'))=e_\beta(2uxx')$ pour tous 
$x,\,y\in \overline{F}$, o\`u $(\cdot ,\cdot)$ d\'esigne l'application commutateur. En particulier si $p=2$, les éléments 
$e_\alpha(x)$ et $e_{\tau{\alpha}}(x_{\tau(\alpha)})$ commutent. Si $p\neq 2$, on choisit $e_\beta$ de telle manière que $u={1\over 2}$. 
Puisque $\tau^2(\alpha)=\alpha$, il existe des $y,\,z\in \smash{\overline{F}}^\times$ tels que 
$\tau\circ e_{\tau(\alpha)}(x)=e_\alpha(yx)$ et $\tau\circ e_\beta (x)=e_\beta(zx)$ 
pour tout $x\in \overline{F}$. Si $p\neq 2$, la relation
$$
\tau(e_\alpha(x),e_{\tau(\alpha)}(x'))=(e_{\tau(\alpha)}(x),e_\alpha(yx'))= e_\beta(-yxx'),
$$
entra\^{\i}ne que $z=-y$. Pour $u=e_\alpha(x_\alpha)e_\beta(x_\beta)e_{\tau(\alpha)}(x_{\tau(\alpha)})\in {\bf U}_{\Phi(\alpha)}$, on a
\begin{align*}
\rho(u)&=e_{\alpha}(x_\alpha)e_\beta(x_\beta)e_{\tau(\alpha)}(x_{\tau(\alpha)})
e_{\alpha}(-y x_{\tau(\alpha)})e_\beta (-zx_\beta)e_{\tau(\alpha)}(-x_\alpha)\\
&= e_\alpha(x_\alpha-yx_{\tau(\alpha)})e_\beta((1-z)x_\beta +2uyx_{\tau(\alpha)}^2)
e_{\tau(\alpha)}(x_{\tau(\alpha)}-x_\alpha).
\end{align*}
On distingue plusieurs cas:
\begin{itemize}
\item si $y\notin\{\pm 1\}$ ou $z\notin\{\pm 1\}$, alors ${\bf U}_{\Phi(\alpha),\tau}=\{1\}$;
\item si  et $p\neq 2$ et $y=-1$, alors 
${\bf U}_{\Phi(\alpha),\tau}=\{e_\beta (x):x\in \overline{F}\}$; 
\item si $p\neq 2$ et $y=1$, alors 
${\bf U}_{\Phi(\alpha),\tau}=\{e_\alpha(x)e_\beta(-{1\over 2}x^2)e_{\tau(\alpha)}(x):x\in \overline{F}\}$;
\item si $p=2$ et $y=z=1$, alors ${\bf U}_{\Phi(\alpha),\tau}=\{e_\alpha(x)e_\beta(x')e_{\tau(\alpha)}(x): x,\,x'\in \overline{F}\}$.
\end{itemize}
Posons $E_\gamma = {\rm Lie}(e_\gamma)(1)$ ($\gamma\in \Phi(\alpha)$). Pour 
$X=x_\alpha E_\alpha + x_\beta E_\beta + x_{\tau(\alpha)}E_{\tau(\alpha)}\in {\rm Lie}({\bf U}_{\Phi(\alpha)})$, on a
$$
X-{\rm Lie}(\tau)(X)= (x_\alpha-yx_{\tau(\alpha)})E_\alpha + (1-z)x_\beta E_\beta + (x_{\tau(\alpha)}-x_\alpha )E_{\tau(\alpha)}.
$$
Soit $\frak{K}= \ker({\rm id}-{\rm Lie}(\tau); {\rm Lie}({\bf U}_{\Phi(\alpha)}))$. On a:
\begin{itemize}
\item si $y\notin\{\pm 1\}$ ou $z\not\in \{\pm 1\}$, alors $\frak{K}=\{0\}$;
\item si $p\neq 2$ et $y=-1$, alors $\frak{K}=\{xE_\beta: x\in \overline{F}\}$;
\item si $p\neq 2$ et $y=1$, alors $\frak{K}=\{x(E_\alpha + E_{\tau(\alpha)}):x\in \overline{F}\}$.
\item si $p=2$ et $y=z=1$, alors $\frak{K}= \{x(E_\alpha + E_{\tau(\alpha)})+ x' E_\beta: x,\, x'\in \overline{F}\}$.
\end{itemize}
On a donc
$${\rm Lie}(\ker \rho)= \frak{K},
$$
i.e. $\rho$ est s\'eparable. Cela ach\`eve la d\'emonstration de la proposition.
\end{proof}

\begin{moncoro}
Supposons ${\bf H}$ semisimple et simplement connexe, et soit 
$\tau\in {\rm Aut}_{\overline{F}}({\bf H})$ quasi--semisimple. Le morphisme ${\bf H}\rightarrow {\bf H},\,h\mapsto h\tau(h)^{-1}$ est s\'eparable.
\end{moncoro}

\begin{proof}
Soit $({\bf B},{\bf T})$ une paire de Borel $\tau$--stable de ${\bf H}$. Posons ${\bf U}=R_{\rm u}({\bf B})$ et 
reprenons les notations de la d\'emonstration de la proposition. 
Pour 
 $\alpha\in \Delta=\Delta({\bf T},{\bf B})$, notons $\alpha^\vee\in {\rm X}_*({\bf T})$ la 
coracine associ\'ee \`a $\alpha$, o\`u ${\rm X}_*({\bf T})$ d\'esigne  
le groupe des cocaract\`eres de ${\bf T}$. Puisque ${\bf H}$ est simplement connexe, 
on a ${\bf X}_*({\bf T})=\oplus_{\alpha\in \Delta}\Bbb{Z}\alpha^\vee$. 
Soit $t\in {\bf T}$. \'Ecrivons  
$t=\prod_{\alpha\in \Delta_\circ}\alpha^\vee(x_\alpha)$,  
$x_\alpha \in \smash{\overline{F}}^\times$. Puisque 
$\tau(\alpha^\vee(x))= \tau(\alpha)^\vee(x)$ ($\alpha\in \Delta,\,
x\in \smash{\overline{F}}^\times$), on a
$$
t\tau (t)^{-1}= 
\prod_{\alpha\in \Delta} \alpha^\vee(x_\alpha x_{\tau^{-1}(\alpha)}^{-1}).
$$
Le groupe ${\bf S}={\bf T}_\tau$ est l'ensemble des $t=\prod_{\alpha\in \Delta}\alpha^\vee(x_\alpha)$ tels que 
$x_{\tau(\alpha)}=x_\alpha$ pour tout $\alpha\in \Delta$. De la m\^eme mani\`ere, on obtient 
que $\frak{T}_\tau ={\rm ker}({\rm id}_\frak{T}-{\rm Lie}(\tau\vert_{\bf T}))$ est l'ensemble des 
$X=\sum_{\alpha\in \Delta}\alpha^\vee(X_\alpha)$ tels que $X_{\tau(\alpha)}=X_\alpha$ pour tout $\alpha\in \Delta$; 
o\`u l'on a pos\'e $\alpha^\vee(Y)={\rm Lie}(\alpha^\vee)(Y)$ pour $Y\in \overline{F}$. On en d\'eduit que 
${\rm Lie}({\bf S})=\frak{T}_\tau$. Le morphisme de groupes $(1-\tau)\vert_{\bf T}$ est donc s\'eparable, et 
d'apr\`es le lemme, le morphisme $(1-\tau)\vert_{\bf B}$ est lui aussi s\'eparable. Soit ${\bf U}^-$ le radical unipotent 
du sous--groupe de Borel ${\bf B}^-$ de ${\bf H}$ oppos\'e \`a ${\bf B}$ par rapport \`a ${\bf T}$. 
Puisque la paire de Borel $({\bf B}^-,{\bf T})$ de ${\bf H}$ est $\tau$--admissible, \`a nouveau d'apr\`es le lemme, 
le morphisme $(1-\tau)\vert_{{\bf U}-}$ est s\'eparable. D'apr\`es \cite[ch.~IV, 14.14]{Bor}, l'application 
produit ${\bf U}^-\!\times {\bf B}\rightarrow {\bf H}$ est un isomorphisme de ${\bf U}^-\!\times {\bf B}$ sur 
un ouvert de ${\bf H}$, et puisque les restrictions de $1-\tau$ \`a ${\bf U}^-$ et ${\bf B}$ sont 
s\'eparables, on a
$$
{\rm Lie}({\bf U}^-_\tau)+ {\rm Lie}({\bf B}_\tau)= {\rm Lie}({\bf U}^-)_\tau + {\rm Lie}({\bf B})_\tau= \frak{H}_\tau.
$$
Donc ${\rm Lie}({\bf H}_\tau)=\frak{H}_\tau$ et le morphisme ${\bf H}\rightarrow {\bf H},\,h\mapsto \tau(h)h^{-1}$ est s\'eparable.
\end{proof}
 
 \begin{marema3}
 {\rm Bien s\^ur, le corollaire n'est en g\'en\'eral plus vrai si ${\bf H}$ 
 n'est pas semisimple simplement connexe: si $p=2$ et $\tau$ est le passage \`a l'inverse dans $\Bbb{G}_{\rm m}$, le 
morphisme $1-\tau:\Bbb{G}_{\rm m}\rightarrow \Bbb{G}_{\rm m},\,x\mapsto x^{2}$ est 
bijectif mais n'est pas s\'eparable.\hfill $\blacksquare$
}\end{marema3}

\begin{marema4}{\rm Soit $\tau$ un automorphisme quasi--semisimple de ${\bf H}$, 
et $({\bf B},{\bf T})$ une paire de Borel $\tau$--stable de ${\bf H}$. Pour $\alpha\in \Phi=\Phi({\bf T},{\bf H})$, 
on note $l=l_{\tau,\alpha}$ le plus 
petit entier $\geq 1$ tel que $\tau^l(\alpha)=\alpha$, et $y=y_{\tau,\alpha}$ 
l'\'el\'ement de $ \smash{\overline{F}}^\times$ d\'efini par $\tau^l\circ e_\alpha(x)=e_\alpha(yx)$ 
($x\in \overline{F}$), o\`u $e_\alpha: \Bbb{G}_{\rm a}\rightarrow {\bf U}_\alpha$ est un 
\'epinglage de ${\bf U}_\alpha$ dont le choix n'a pas d'importance. D'ailleurs l'entier $l$ et l'\'el\'ement 
$y$ ne d\'ependent pas vraiment de $\alpha$ mais seulement de la $\tau$--orbite 
$\ES{O}=\{\alpha,\tau(\alpha),\ldots , \tau^{l-1}(\alpha)\}\subset \Phi$. On les note donc aussi 
$l_{\tau,\ES{O}}$ et $y_{\tau,\ES{O}}$. Soit $\Phi_\tau$ l'ensemble des $\tau$--orbites $\ES{O}$ dans 
$\Phi$ telles que $y_{\tau,\ES{O}}=1$, et $\Phi({\bf T}_\tau^\circ, {\bf H}_\tau^\circ)$ l'ensemble des racines de ${\bf T}_\tau^\circ$ 
dans ${\bf H}_\tau^\circ$. Pour une $\tau$--orbite $\ES{O}$ dans $\Phi$, on pose
$$
\ES{O}'= \vert \ES{O}\vert^{-1} \sum_{\alpha\in \ES{O}} \alpha \in {\rm X}^*({\bf T})\otimes_{\Bbb Z}{\Bbb R}.
$$
D'après \cite[th\'eo.~1.8]{DM1}, l'application $\Phi/\langle \tau \rangle \rightarrow {\rm X}^*({\bf T}_\tau^\circ),\, \ES{O}\mapsto \ES{O}'\vert_{{\rm X}^*({\bf T}_\tau^\circ)}$ induit par restriction une application surjective
$$
\Phi_\tau \rightarrow \Phi({\bf T}_\tau^\circ, {\bf H}_\tau^\circ).
$$
Cette surjection est bijective, sauf si $p=2$ et s'il existe une orbite $\ES{O}\in \Phi_\tau$ contenant deux racines $\alpha$ et $\alpha'$ telles que $\alpha + \alpha'\in \Phi$, auquel cas pour avoir une bijection il faut exclure ces orbites $\ES{O}$ de l'ensemble $\Phi_\tau$. Notons que, pour tout $p$, l'existence d'une telle orbite $\ES{O}\in \Phi_\tau$ n'est possible que si le 
diagramme de Dynkin de $G$ possède $k$ composantes de type $A_{2n}$, permutées transitivement par $\tau$ et telles que $\tau^k$ opère sur chacune d'elles par ``retournement'', \cad de la forme
$$
A_{2n}: 
\xymatrix{{\circ} \ar@{-}[r] ^<{\alpha_n}& \cdots\ar@{-}[r]  &  {\circ}\ar@{-}[r] ^<{\alpha_{2}}&{\circ} \ar@{-}[r] 
^<{\alpha_1}&{\circ} \ar@{-}[r] ^<{\tau^k(\alpha_1)}&   {\circ} \ar@{-}[r] ^<{\tau^k(\alpha_{2})} & 
\cdots \ar@{-}[r] & {\circ} \ar@{} ^<{\tau^k(\alpha_n)}
}.
$$
Alors pour toute racine $\alpha$ (dans l'une de ces $k$ composantes) telle que 
$\alpha + \tau^k(\alpha)$ soit une racine, les $\tau$--orbites de $\alpha$ et de $\alpha + \tau^k(\alpha)$ ont même image dans ${\rm X}^*({\bf T}_\tau^\circ)$, et
$$y_{\tau,\alpha}= - y_{\tau, \alpha + \tau^k(\alpha)}\;(\hbox{$= y_{\tau, \alpha + \tau^k(\alpha)}$ si $p=2$}).
$$

Pour une description explicite du système de racines $\Phi({\bf T}_\tau^\circ,{\bf H}_\tau^\circ)$, on 
renvoie \`a \cite{DM2}. \hfill $\blacksquare$}
\end{marema4}

\subsection{Automorphismes quasi--centraux}\label{automorphismes qc} 
Un automorphisme quasi--semisimple $\tau$ de ${\bf H}$ est 
dit {\it quasi--central} si pour tout automorphis\-me quasi--semisimple de 
${\bf H}$ de la forme ${\rm Int}_{\bf H}(h)\circ \tau $ avec $h\in {\bf H}$, 
on a $\dim({\bf H}_{{\rm Int}_{\bf H}(h)\circ \tau}^\circ)\leq  \dim({\bf H}_\tau^\circ)$. 

\begin{monlem}
Soit $\tau\in {\rm Aut}_{\overline{F}}({\bf H})$ quasi--semisimple, et 
soit $({\bf B},{\bf T})$ une paire de Borel $\tau$--stable de ${\bf H}$. Alors il existe un 
$t\in {\bf T}$ tel que le $\overline{F}$--automorphisme ${\rm Int}_{\bf H}(t)\circ \tau$ de 
${\bf H}$ est quasi--central.
\end{monlem}

\begin{proof}
Soit $h\in {\bf H}$ tel que le $\overline{F}$--automorphisme 
$\tau'={\rm Int}_{\bf H}(h)\circ \tau$ de ${\bf H}$ est quasi--central. Puisque 
$\tau'$ est quasi--semisimple, il existe un $x\in {\bf H}$ tel que $\tau'$ stabilise 
la paire de Borel ${\rm Int}_{\bf H}(x)\circ ({\bf B},{\bf T})$ de ${\bf H}$. D'apr\`es 
la relation $(*)$ de \ref{automorphismes}, quitte \`a remplacer $h$ par $x^{-1}h\tau(x)$, on peut supposer 
que $\tau'({\bf B},{\bf T})=({\bf B},{\bf T})$. Alors on a ${\rm Int}_{\bf H}(h)({\bf B},{\bf T})=
({\bf B},{\bf T})$, donc $h\in {\bf T}$.
\end{proof}

D'aprs \cite[prop.~1.21, cor.~1.25]{DM1}, on a la
\begin{mapropo}
Soit $\tau\in {\rm Aut}_{\overline{F}}({\bf H})$ quasi--central.
\begin{enumerate}
\item[(1)] Si $({\bf B},{\bf T})$ et $({\bf B}'\!,{\bf T}')$ sont deux paires de Borel $\tau$--stables de 
${\bf H}$, alors il existe un $h\in {\bf H}_\tau^\circ$ tel que $({\bf B}'\!,{\bf T}')=
{\rm Int}_{\bf H}(h)({\bf B},{\bf T})$.
\item[(2)] Soit ${\bf P}$ un sous--groupe parabolique $\tau$--stable de ${\bf H}$, et 
${\bf L}$ une composante de Levi $\tau$--stable de ${\bf P}$. Si $h\in {\bf H}$ v\'erifie 
$h^{-1}\tau(h)\in {\bf L}$, alors $h\in {\bf H}_\tau^\circ\, {\bf L}$. 
\item[(3)] L'application ${\bf P}\mapsto 
{\bf P}_\tau^\circ$ est une bijection entre l'ensemble des sous--groupes 
paraboliques $\tau$--stables de ${\bf H}$ et l'ensemble des 
sous--groupes paraboliques de ${\bf H}_\tau^\circ$. Soit ${\bf P}$ un 
sous--groupe parabolique $\tau$--stable de ${\bf H}$. 
L'application ${\bf L}\mapsto  {\bf L}_\tau^\circ$ est une 
bijection entre l'ensemble des composantes de Levi $\tau$--stables de 
${\bf P}$ et l'ensemble des composantes de Levi de ${\bf P}_\tau^\circ$, de 
bijection inverse ${\bf L}^\sharp \mapsto Z_{\bf H}(R({\bf L}^\sharp))$.
\end{enumerate}
\end{mapropo}

\begin{mesrems}{\rm 
\begin{enumerate}
\item[(1)]Un automorphisme quasi--semi\-simple $\tau$ de ${\bf H}$ est quasi--central 
si et seulement s'il v\'erifie les propri\'et\'es \'equivalentes suivantes \cite[d\'ef.--th\'eo.~1.15]{DM1}:
\begin{itemize}
\item[(i)]Tout sous--groupe de Borel de ${\bf H}_\tau^\circ$ est contenu dans un {\it unique} sous--groupe de Borel $\tau$--stable de ${\bf H}$.
\item[(ii)]Pour toute (resp. pour une) paire de Borel $\tau$--stable $({\bf B},{\bf T})$ de ${\bf H}$, posant ${\bf N}=N_{\bf H}({\bf T})$ et 
${\bf W}={\bf N}/{\bf T}$, tout \'el\'ement $\tau$--stable de ${\bf W}$ a un repr\'esentant dans 
${\bf N}\cap {\bf H}_\tau^\circ$.
\item[(iii)]Pour toute (resp. pour une) paire de Borel $\tau$--stable $({\bf B},{\bf T})$ de ${\bf H}$, posant $\Phi=\Phi({\bf T},{\bf H})$ 
et $\Delta=\Delta({\bf T},{\bf B})$, les éléments $y_{\tau,\alpha}\in \smash{\overline{F}}^\times$ (cf. 
la remarque 4 de \ref{automorphismes qss}) pour $\alpha\in \Delta$ vérifient la condition\footnote{D'apr\`es les relations de Chevalley, si cette condition est vérifiée pour toute racine simple $\alpha\in \Delta$, alors elle l'est pour toute racine $\alpha\in \Phi$ --- cf. \cite{DM2}. 
Par ailleurs, si $y_{\tau,\alpha}=1$ pour toute racine simple $\alpha\in \Delta$, alors on peut choisir des épinglages 
$e_\alpha: {\Bbb G}_{\rm a}\rightarrow {\bf U}_\alpha$ de ${\bf U}_\alpha$ pour $\alpha\in \Delta$ de telle manière 
que la paire de Borel épinglée $({\bf B},{\bf T}, \{e_\alpha\}_{\alpha\in \Delta})$ soit $\tau$--stable. Réciproquement, 
si $\tau$ stabilise une paire de Borel épinglée $({\bf B},{\bf T}, \{e_\alpha\}_{\alpha\in \Delta})$ de ${\bf H}$, alors $y_{\tau,\alpha}=1$ pour toute racine simple $\alpha\in \Delta$, et $\tau$ est quasi--central.} 
$$
y_{\tau,\alpha}\in \{\pm 1\},
$$
où $-1$ n'est autorisé que s'il existe deux racines de la $\tau$--orbite de $\alpha$ dont la somme appartient à $\Phi$.
\end{itemize}
\item[(2)]Soit un \'el\'ement $h\in {\bf H}$ tel que l'automorphisme int\'erieur $\tau={\rm Int}_{\bf H}(h)$ de ${\bf H}$ 
est quasi--central. Alors $h$ appartient au centre $Z({\bf H})$ de ${\bf H}$. 
En effet, puisque ${\rm Int}_{\bf H}(h)$ est 
quasi--semisimple, $h$ est semisimple, donc appartient \`a un tore maximal ${\bf T}$ de ${\bf H}$. D'apr\`es 
la condition (ii), tout \'el\'ement de $N_{\bf H}({\bf T})/{\bf T}$ a un repr\'esentant dans 
$N_{\bf H}({\bf T})\cap {\bf H}_\tau^\circ$. Ainsi $h$ est un \'el\'ement de ${\bf T}$ centralis\'e par 
$N_{\bf H}({\bf T})$, donc par ${\bf H}$ tout entier.
\item[(3)]Soit $\tau$ un automorphisme quasi--semisimple de ${\bf H}$, et $({\bf B},{\bf T})$ une paire de 
Borel $\tau$--stable de ${\bf H}$. Soit $\alpha\in \Phi({\bf T},{\bf H})$, et soit $e_\alpha:\Bbb{G}_{\rm a}\rightarrow {\bf U}_\alpha$ 
un \'epinglage de ${\bf U}_\alpha$. Posons $l=l_{\tau,\alpha}$ comme dans la remarque 4 de \ref{automorphismes qss}. Pour $t\in {\bf T}$, posant 
$\tau'={\rm Int}_{\bf H}(t)\circ \tau$ et $\ES{N}_{\tau,\alpha}(t)=t \tau(t)\cdots \tau^{l-1}(t)$, on a $l_{\tau'\!,\alpha}=l$ et
$$
\tau'^l\circ e_\alpha(x)= {\rm Int}_{\bf H}(\ES{N}_{\tau,\alpha}(t))\circ \tau^l\circ e_\alpha(x)= 
e_\alpha(\alpha(\ES{N}_{\tau,\alpha}(t))y_{\tau,\alpha}x)\quad (x\in \overline{F}).
$$
L'\'el\'ement $\alpha(\ES{N}_{\tau,\alpha}(t))\in \smash{\overline{F}}^\times$ ne d\'epend pas vraiment de 
$\alpha$, mais seulement de la $\tau$--orbite $\ES{O}$ de $\alpha$ dans $\Phi({\bf T},{\bf H})$. On le note donc 
aussi $a_{\tau,\ES{O}}(t)$. D'apr\`es la condition (iii), 
si l'on choisit $t$ de telle mani\`ere que pour chaque $\tau$--orbite 
$\ES{O}$ dans $\Delta({\bf T},{\bf B})$, on ait $a_{\tau,\ES{O}}(t)=y_{\tau,\ES{O}}^{-1}$, alors $\tau'$ est quasi--central.
\hfill $\blacksquare$
\end{enumerate}}
\end{mesrems}

\subsection{Automorphismes quasi--semisimples localement finis}\label{automorphismes qss lf} 
La notion d'automor\-phisme quasi--semisimple est 
justifi\'ee par le r\'esultat de Steinberg suivant \cite[7.5]{St}:

\begin{montheo}
Tout $\overline{F}$--automorphisme semisimple de ${\bf H}$ est 
quasi--semisimple.
\end{montheo}

\begin{marema1}{\rm 
D'apr\`es \cite[9]{St}, si $\tau'\in{\rm Aut}_{\overline{F}}({\bf H}_{\rm der})$ 
est quasi--semisimple, alors $\tau'$ est semisimple 
si et seulement si l'image de $\tau'$ dans ${\rm Out}_{\overline{F}}({\bf H}_{\rm der})$ 
est d'ordre premier \`a $p$.\hfill $\blacksquare$}
\end{marema1}

\begin{marema2}
{\rm Supposons $p=1$, et soit 
$\tau\in {\rm Aut}_{\overline{F}}({\bf H})$. Alors $\tau$ est quasi--semisimple 
si et seulement si $\tau_{\rm der}$ est semisimple. D'autre part si 
$\tau$ est localement fini, alors $\tau$ est quasi--semisimple si et seulement si $\tau$ est semisimple. 
En effet, supposons $\tau_{\rm der}$ semisimple et montrons que 
$\tau$ l'est aussi. \'Ecrivons la d\'ecomposition de Jordan 
$\tau=\tau_{\rm s}\circ \tau_{\rm u}$. Puisque $\tau_{\rm cent}$ 
est d'ordre fini et $\tau_{\rm der}$ est semisimple, il 
existe un entier $n\geq 1$ et un \'el\'ement $h\in {\bf H}$ 
semisimple, tels que 
$\tau^n={\rm Int}_{\bf H}(h)$. Alors $(\tau^n)_{\rm s}=(\tau_{\rm s})^n$ et 
$(\tau^n)_{\rm u}=(\tau_{\rm u})^n$. Mais puisque 
$\tau^n$ est semisimple, on a $(\tau_{\rm u})^n={\rm id}_{\bf H}$ et 
$\tau_{\rm u}$ est semisimple. Donc $\tau_{\rm u}=1$, i.e. $\tau$ est semisimple.
\hfill $\blacksquare$}
\end{marema2}

Soit $\tau\in {\rm Aut}_{\overline{F}}^0({\bf H})$. \'Ecrivons la d\'ecomposition 
de Jordan $\tau=\tau_{\rm s}\circ î\tau_{\rm u}$. D'apr\`es le thorme 2 de \ref{automorphismes qss} (et le 
thorme ci--dessus), le groupe 
${\bf H}_{\tau_{\rm s}}^\circ$ est r\'eductif, et puisque $\tau_{\rm s}\circ \tau_{\rm u}=
\tau_{\rm u}\circ \tau_{\rm s}$, $\tau$ induit par restriction un $\overline{F}$--automorphisme 
{\it unipotent} de ${\bf H}_{\tau_{\rm s}}^\circ$, que l'on note $\tau^*$. On a donc l'inclusion 
$({\bf H}_{\tau_{\rm s}}^\circ)_{\tau^*}={\bf H}_{\tau_{\rm s}}^\circ\cap {\bf H}_{\tau_{\rm u}}
\subset {\bf H}_\tau$.

D'aprs \cite[lemme~1.14 et cor.~1.33]{DM1}, on a la
\begin{mapropo}
Soit $\tau\in {\rm Aut}_{\overline{F}}^0({\bf H})$.
\begin{enumerate}
\item[(1)]$\tau$ est quasi--semisimple si et seulement si 
$\tau^*$ est quasi--semisimple.
\item[(2)]Si $\tau$ est quasi--semisimple et unipotent, alors $\tau$ est 
quasi--central et ${\bf H}_\tau$ est connexe.
\end{enumerate}
\end{mapropo}

\begin{moncoro}
Soit $\tau\in {\rm Aut}_{\overline{F}}^0({\bf H})$ 
quasi--semisimple. On a l'\'egalit\'e ${\bf H}_\tau^\circ = ({\bf H}_{\tau_{\rm s}}^\circ)_{\tau^*}$.
\end{moncoro}

\begin{proof} Puisque $({\bf H}_{\tau_{\rm s}})_{\tau^*}$ est connexe et contenu dans 
${\bf H}_\tau$, on a l'inclusion $({\bf H}_{\tau_{\rm s}}^\circ)_{\tau^*}\subset {\bf H}_\tau^\circ$. 
Pour l'inclusion inverse, identifions ${\bf H}$ \`a la composante neutre 
du groupe alg\'ebrique affine ${\bf H}'={\bf H}\rtimes \langle \tau \rangle /{\bf C}$ comme en \ref{automorphismes ss et u} et 
notons $\delta$ l'image de $1\rtimes \tau$ dans ${\bf H}'$.  Alors on a ${\bf H}_\tau^\circ= 
({\bf H}'_\delta)^\circ$. \'Ecrivons la d\'ecomposition de Jordan $\delta=\delta_{\rm s}\delta_{\rm u}$. 
L'unicit\'e de cette d\'ecomposition entra\^{\i}ne 
l'inclusion ${\bf H}'_\delta\subset {\bf H}'_{\delta_{\rm s}}\cap {\bf H}'_{\delta_{\rm u}}$. 
Or on a $({\bf H}'_{\delta_{\rm s}})^\circ = {\bf H}_{\tau_{\rm s}}$ et ${\bf H}'_{\delta_{\rm u}}\cap {\bf H}
={\bf H}_{\tau_{\rm u}}$, d'o\`u ${\bf H}_\tau^\circ\subset {\bf H}_{\tau_{\rm s}}^\circ \cap {\bf H}_{\tau_{\rm u}}=
({\bf H}_{\tau_{\rm s}}^\circ)_{\tau^*}$.
\end{proof}

\begin{marema3}
{\rm 
Supposons $p=1$. D'apr\`es la remarque (5) de \ref{automorphismes}, tout $\overline{F}$--automorphis\-me 
unipotent de ${\bf H}$ est de la 
forme ${\rm Int}_{\bf H}(u)$ pour un \'el\'ement unipotent $u\in {\bf H}$. En particulier, 
l'identit\'e de ${\bf H}$ est le seul $\overline{F}$--automorphisme quasi--semisimple 
unipotent de ${\bf H}$.
\hfill $\blacksquare$}
\end{marema3}

\begin{monlem}Soit $\tau\in {\rm Aut}_{\overline{F}}({\bf H})$ 
quasi--semisimple et unipotent, et soit $({\bf B},{\bf T})$ une paire de Borel 
$\tau$--stable de ${\bf H}$. Soit $h\in {\bf H}$ tel que $\tau'={\rm Int}_{\bf H}(h)\circ \tau$ 
est quasi--semisimple. Alors il existe un $t\in {\bf T}_\tau$ et un $x\in {\bf H}$ tels 
que $\tau' = {\rm Int}_{\bf H}(x^{-1})\circ {\rm Int}_{\bf H}(t)\circ \tau
\circ {\rm Int}_{\bf H}(x)$. De plus $\tau'$ est unipotent si et seulement si $t=1$.
\end{monlem}

\begin{proof}
Puisque $\tau'$ est quasi--semisimple, il existe un $x'\in {\bf H}$ tel que $\tau'$ 
stabilise la paire de Borel ${\rm Int}_{\bf H}(x')({\bf B},{\bf T})$ de ${\bf H}$. Alors
${\rm Int}_{\bf H}(x'^{-1}h\tau(x'))
\circ \tau({\bf B},{\bf T})=({\bf B},{\bf T})$, d'o\`u $x'^{-1}h\tau(x')\in {\bf T}$. D'apr\`es la proposition de \ref{automorphismes ss et u}, 
on a la dcomposition ${\bf T}={\bf T}_\tau{\bf T}(1-\tau)$. \'Ecrivons $x'^{-1}h\tau(x')= ty^{-1}\tau(y)$ 
avec $t\in {\bf T}_\tau$ et $y\in {\bf T}$. Puisque $ty^{-1}\tau(y)=y^{-1}t\tau(y)$, 
on obtient
$$
{\rm Int}_{\bf H}(x'^{-1}h\tau(x'))\circ \tau = {\rm Int}_{\bf H}(y^{-1})\circ {\rm Int}_{\bf H}(t)\circ \tau
\circ {\rm Int}_{\bf H}(y),
$$
d'o\`u $\tau'= {\rm Int}_{\bf H}(x'y^{-1})\circ {\rm Int}_{\bf H}(t)\circ \tau
\circ {\rm Int}_{\bf H}(yx'^{-1})$.

Si $t=1$, alors $\tau'$ est unipotent. R\'eciproquement, 
supposons $\tau'$ unipotent et posons $\tau''= {\rm Int}_{\bf H}(t)\circ \tau$. Alors 
$\tau''$ est 
quasi--semisimple et unipotent, et comme $t$ appartient  $ {\bf H}_\tau$, le $\overline{F}$--automorphisme 
$\tau''\circ \tau^{-1}= \tau^{-1}\circ 
\tau''$ de ${\bf H}$ est unipotent. Donc ${\rm Int}_{\bf H}(t)$ est unipotent, et $t=1$.
\end{proof}

\subsection{Automorphismes r\'eguliers; les automorphismes int\'erieurs}\label{automorphismes rguliers; le cas i} 
Pour $\tau\in {\rm Aut}_{\overline{F}}({\bf H})$, 
on note $r_\tau^1({\bf H})$ le plus petit entier $k\geq 1$ tel que la fonction 
${\bf H} \rightarrow \overline{F},\,\gamma\mapsto a_k(h,\tau)$ 
d\'efinie par le polyn\^ome en l'indetermin\'ee $t$ 
$$
P(h,\tau)(t)={\rm det}_{\overline{F}}(t- {\rm Ad}_{\bf H}(h)\circ{\rm Lie}(\tau)+{\rm id}_\frak{H};\frak{H})
=\sum_{i=0}^\infty a_i
(h,\tau)t^i,
$$
est non nulle; et l'on pose
$$
D_{\bf H}(\tau)=a_{r_\tau^1({\bf H})}(1,\tau)\in \overline{F}.
$$
On note aussi $r_\tau({\bf H})$ le plus grand entier $k' \geq 1$ v\'erifiant 
la propri\'et\'e: pour tout $h\in {\bf H}$ tel que ${\rm Int}_{\bf H}(h)\circ \tau$ est quasi--semisimple, 
on a $k'\leq \dim({\bf H}_{{\rm Int}_{\bf H}(h)\circ\tau}^\circ)$.

Soit $\tau\in {\rm Aut}_{\overline{F}}({\bf H})$. Pour 
$x\in \smash{\overline{F}}^\times$, notons $\frak{H}_\tau^x$ le sous--espace caract\'eristique de 
${\rm Lie}(\tau)$ associ\'e 
\`a la valeur propre $x$; i.e. posons
$$
\frak{H}_\tau^x=\{X\in \frak{H}:(x{\rm id}_\frak{H}-{\rm Lie}(\tau))^k(X)=0,\,\exists k\in \Bbb{Z}_{\geq 1}\}.
$$
On a la d\'ecomposition $\frak{H}=\oplus_{x\in \smash{\overline{F}}^\times}\frak{H}_\tau^x$. 
Par d\'efinition, on a $\dim_{\overline{F}}(\frak{H}_\tau^1)\geq r_\tau^1({\bf H})$ avec \'egalit\'e si et seulement si 
$D_{\bf H}(\tau)\neq 0$, auquel cas on a
$$
D_{\bf H}(\tau)={\rm \det}_{\overline{F}}({\rm id} -{\rm Lie}(\tau);\frak{H}/\frak{H}_\tau^1).
$$
Rappelons que $
\frak{H}_\tau= \ker\{{\rm id}_\frak{H}-{\rm Lie}(\tau)\}$. 
On a donc les inclusions
$${\rm Lie}({\bf H}_\tau^\circ)\subset \frak{H}_\tau\subset \frak{H}_\tau^1
$$ avec \'egalit\'es 
si $\tau$ est semisimple \cite[ch.~III, 9.1]{Bor}.

Soit
$\tau\in {\rm Aut}_{\overline{F}}^0({\bf H})$. \'Ecrivons la d\'ecomposition de 
Jordan $\tau=\tau_{\rm s}\circ \tau_{\rm u}$. Alors pour $x\in \smash{\overline{F}}^\times$, 
on a $\frak{H}_\tau^x=\frak{H}_{\tau_{\rm s}}^x=
\ker\{x{\rm id}_\frak{H}-{\rm Lie}(\tau_{\rm s})\}$. 
En particulier, on a 
l'\'egalit\'e $\frak{H}_\tau^1= \frak{H}_{\tau_{\rm s}}$, et
$$
D_{\bf H}(\tau)\neq 0\Leftrightarrow \dim({\bf H}_{\tau_{\rm s}}^\circ)=r_\tau^1({\bf H}).
\leqno{(*)}
$$

\begin{madefi}
{\rm Un $\overline{F}$--automorphisme 
$\tau$ de ${\bf H}$ est dit {\it r\'egulier} 
si $D_{\bf H}(\tau)\neq 0$.
}\end{madefi}

Soit $\tau\in {\rm Aut}_{\overline{F}}({\bf H})$. D'apr\`es \ref{automorphismes ss et u}, 
pour $z\in R({\bf H})$ et $h'\in {\bf H}_{\rm der}$, on a
$$
\dim_{\overline{F}}(\frak{H}_{{\rm Int}_{\bf H}(zh')\circ \tau}^1)=
\dim_{\overline{F}}((\frak{H}_{\rm der})_{{\rm Int}_{{\bf H}_{\rm der}}(h')\circ \tau_{\rm der}}^1)+
\dim_{\overline{F}}(\frak{C}_{\tau_{\rm cocent}}^1).
$$
Par suite on a
$$
r_\tau^1({\bf H})=r_{\tau_{\rm der}}^1({\bf H}_{\rm der})+\dim_{\overline{F}}(\frak{C}_{\tau_{\rm cocent}}^1),
$$
et $\tau$ est r\'egulier 
si et seulement si $\tau_{\rm der}$ est r\'egulier.

Soit $\tau\in {\rm Int}_{\overline{F}}({\bf H})$. Par d\'efinition, les entiers $r_\tau({\bf H})$ et 
$r_\tau^1({\bf H})$ ne d\'ependent pas de $\tau$. On les note respectivement 
$r({\bf H})$ et $r^1({\bf H})$.

\begin{mapropo}
Soit $\tau={\rm Int}_{\bf H}(h)$ pour un 
$h\in {\bf H}$. Alors $D_{\bf H}(\tau)\neq 0$ 
si et seulement si $h$ est semisimple et ${\bf H}_h^\circ$ est un tore 
(i.e. si et seulement si $h$ est \og semisimple r\'egulier\fg{} au sens de Borel \cite[ch.~IV, 12.2]{Bor}). De 
plus, on a $r({\bf H})=r^1({\bf H})$, et 
cet entier co\"{\i}ncide avec le rang de ${\bf H}$.
\end{mapropo}

\begin{proof}
\'Ecrivons les 
d\'ecompositions de Jordan $h=h_{\rm s}h_{\rm u}$ et 
$\tau=\tau_{\rm s}\circ \tau_{\rm u}$. On a 
$\tau_{\rm s}={\rm Int}_{\bf H}(h_{\rm s})$ et 
$\tau_{\rm u}={\rm Int}_{\bf H}(h_{\rm u})$, et d'apr\`es la remarque (4) de \ref{automorphismes}, $h_{\rm u}$ appartient \`a 
${\bf H}_{h_{\rm s}}^\circ$. 
Puisque $\tau_{\rm s}$ est semisimple, 
on a ${\rm Lie}({\bf H}_{h_{\rm s}}^\circ) =
\frak{H}_{\tau_{\rm s}}=\frak{H}_{\tau_{\rm s}}^1$, et comme 
le $\overline{F}$--automorphisme $\tau^*$ de ${\bf H}_{h_{\rm s}}^\circ$ est unipotent, on a
$$\dim({\bf H}_{\tau_{\rm s}}^\circ)=
\dim_{\overline{F}}(\frak{H}_{\tau_{\rm s}})\leq \dim_{\overline{F}}(\frak{H}_\tau^1).
$$ 
Par cons\'equent si $D_{\bf H}(\tau)\neq 0$, alors 
$D_{\bf H}(\tau_{\rm s})\neq 0$. En particulier, il existe un \'el\'ement semisimple $h'\in {\bf H}$ 
tel que $D_{\bf H}({\rm Int}_{\bf H}(h'))\neq 0$.

Choisissons une paire de Borel $({\bf B},{\bf T})$ de ${\bf H}$ 
telle que $h_{\rm s}\in {\bf T}$, et posons $\Phi=\Phi({\bf T},{\bf H})$. Alors on a l'inclusion 
${\bf T}\subset {\bf H}_{h_{\rm s}}^\circ$, avec \'egalit\'e si et seulement si 
pour toute racine $\alpha\in \Phi$, on a $\alpha(t)\neq 1$ \cite[ch.~IV, 12.2]{Bor}. 
L'ensemble des $t\in {\bf T}$ tels que $\alpha(t)\neq 1$ pour toute racine 
$\alpha\in \Phi$, est ouvert dense dans ${\bf T}$ 
(en particulier il est non vide). On en d\'eduit que $D_{\bf H}(\tau_{\rm s})\neq 0$ 
si et seulement si ${\bf H}_{h_{\rm s}}^\circ={\bf T}$, auquel cas $h_{\rm u}=1$. 
D'o\`u le lemme.
\end{proof}

\begin{moncoro}
Soit ${\bf T}$ un 
tore maximal de ${\bf H}$. Pour $t\in {\bf T}$, on a $D_{\bf H}({\rm Int}_{\bf H}(t))\neq 0$ si et seulement 
si $\alpha(t)\neq 1$ pour toute racine 
$\alpha\in \Phi({\bf T},{\bf H})$.
\end{moncoro}

\subsection{Automorphismes r\'eguliers; le cas g\'en\'eral}\label{automorphismes rguliers; le cas g} 
L'\'etude des $\overline{F}$--automorphismes ext\'e\-rieurs r\'eguliers 
de ${\bf H}$ est plus compliqu\'ee, du moins si $p>1$.

\begin{monlem1}
Supposons qu'il existe un $\overline{F}$--automorphisme 
r\'egulier $\tau$ de ${\bf H}$ tel que $\frak{H}_\tau^1=\frak{H}$. Alors ${\bf H}$ est un tore.
\end{monlem1}

\begin{proof}
Soit $\tau\in {\rm Aut}_{\overline{F}}({\bf H})$ tel 
que $D_{\bf H}(\tau)\neq 0$ 
et $\frak{H}_\tau^1=\frak{H}$. 
Alors pour tout $h\in {\bf H}$, on a  
$\frak{H}_{{\rm Int}_{\bf H}(h)\circ \tau}^1=\frak{H}$ et 
$D_{\bf H}({\rm Int}_{\bf H}(h)\circ \tau)\neq 0$. Par suite, 
quitte \`a remplacer $\tau$ par ${\rm Int}_{\bf H}(h)\circ \tau$ pour un 
$h\in {\bf H}$, on peut supposer $\tau$ quasi--semisimple. Ensuite, quitte 
\`a remplacer ${\bf H}$ par ${\bf H}_{\rm der}$, on peut supposer 
$\tau$ localement fini. \'Ecrivons la d\'ecomposition 
de Jordan $\tau=\tau_{\rm s}\circ \tau_{\rm u}$. Puisque 
$\frak{H}_{\tau_{\rm s}}=\frak{H}_\tau^1=\frak{H}$ et 
${\rm Lie}({\bf H}_{\tau_{\rm s}}^\circ)=\frak{H}_{\tau_{\rm s}}$, on a 
${\bf H}_{\tau_{\rm s}}^\circ ={\bf H}$ et $\tau_{\rm s}={\rm id}_{\bf H}$. Ainsi 
$\tau=\tau_{\rm u}$ est r\'egulier, quasi--semisimple et unipotent. 
D'apr\`es le thorme 2 de \ref{automorphismes qss}, le groupe ${\bf H}_\tau^\circ$ est 
r\'eductif. Soit ${\bf T}^\sharp$ un tore maximal de ${\bf H}_\tau^\circ$. Alors d'aprs loc.~cit., 
${\bf T}=Z_{\bf H}({\bf T}^\sharp)$ est un tore maximal 
($\tau$--admissible) de ${\bf H}$. Notons ${\bf \Omega}^\sharp$ 
l'ensemble des $t\in {\bf T}^\sharp$ tels que pour toute racine 
$\alpha$ de ${\bf T}^\sharp$ dans ${\bf H}$, on a $\alpha(t)\neq 1$. 
D'apr\`es \cite[ch.~III, 9.5]{Bor}, ${\bf \Omega}^\sharp$ est non vide, donc ouvert dense dans 
${\bf T}^\sharp$, et pour $t\in {\bf T}^\sharp$, on a 
${\bf H}_t^\circ =Z_{\bf H}({\bf T}^\sharp)={\bf T}$. Soit 
$\tau' = {\rm Int}_{\bf H}(t)\circ \tau$ pour un $t\in {\bf \Omega}^\sharp$. Puisque 
${\rm Int}_{\bf H}(t)$ est semisimple 
et commute \`a $\tau$, on a $({\rm Int}_{\bf H}(h)\circ \tau)_{\rm s}={\rm Int}_{\bf H}(h)$ et 
$({\rm Int}_{\bf H}(h)\circ \tau)_{\rm u}=\tau$, d'o\`u $\frak{H}_{\tau'}^1=\frak{H}_{{\rm Int}_{\bf H}(t)}
={\rm Lie}({\bf T})$. Par consquent $\frak{H}={\rm Lie}({\bf T})$, ce qui n'est possible que si ${\bf H}$ est un tore.
\end{proof}

\begin{monlem2}
Soit $\tau$ un $\overline{F}$--automorphisme
r\'egulier de ${\bf H}$, et posons ${\bf T}=Z_{\bf H}({\bf H}_\tau^\circ)$.
\begin{enumerate}
\item[(1)]$\tau$ est quasi--semisimple, 
${\bf H}_\tau^\circ$ est un tore, ${\bf T}$ est 
l'unique tore maximal $\tau$--admissible de ${\bf H}$, et $\frak{H}_\tau^1\subset {\rm Lie}({\bf T})$.
\item[(2)] Soit ${\bf B}$ est un sous--groupe 
de Borel $\tau$--stable de ${\bf H}$ contenant ${\bf T}$, et posons 
${\bf U}=R_{\rm u}({\bf B})$. 
L'application ${\bf U}\rightarrow {\bf U},\,u\mapsto u\tau(u)^{-1}$ est un automorphisme de vari\'et\'es alg\'ebriques.
\item[(3)]Supposons $\tau$ localement fini. 
Alors $\tau_{\rm s}$ est r\'egulier, 
${\bf T}=Z_{\bf H}({\bf H}_{\tau_{\rm s}}^\circ)$, et $D_{\bf H}(\tau)=D_{\bf H}(\tau_{\rm s})$.
\end{enumerate}
\end{monlem2}

\begin{proof}
Commen\c{c}ons par montrer que $\tau$ est quasi--semisimple. Puisque $\tau_{\rm der}$ 
est r\'egulier, quitte \`a remplacer ${\bf H}$ par ${\bf H}_{\rm der}$, on peut 
supposer $\tau$ localement fini. 
\'Ecrivons la d\'ecomposition de Jordan $\tau=\tau_{\rm s}\circ \tau_{\rm u}$, et 
posons ${\bf H}'={\bf H}_{\tau_{\rm s}}^\circ$. D'apr\`es le thorme de \ref{automorphismes qss lf} et le thorme 2 de \ref{automorphismes qss}, 
${\bf H}'$ est r\'eductif. Posons $\frak{H}'={\rm Lie}({\bf H}')$. Pour 
$h'\in {\bf H}'$ et $x\in \smash{\overline{F}}^\times$, puisque ${\rm Int}_{\bf H}(h')\circ 
\tau_{\rm s}=\tau_{\rm s}\circ {\rm Int}_{\bf H}(h')$, 
les $\overline{F}$--automorphismes ${\rm Ad}_{\bf H}(h')$ et 
${\rm Ad}_{\bf H}(h')\circ {\rm Lie}(\tau)$ de $\frak{H}$ stabilisent 
$\frak{H}_\tau^x= \frak{H}_{\tau_{\rm s}}^x$, et si $h'$ est 
suffisamment proche de $1$ dans ${\bf H}'$, alors pour tout 
$x\in \smash{\overline{F}}^\times \smallsetminus\{1\}$ 
tel que $\frak{H}_\tau^x\neq 0$, les valeurs propres de 
${\rm Ad}_{\bf H}(h')\circ {\rm Lie}(\tau)\vert_{\frak{H}_\tau^x}\in 
{\rm GL}(\frak{H}_\tau^x)$ sont toutes diff\'erentes de $1$. 
Choisissons un tel $h'$, et posons $\tau_1={\rm Int}_{\bf H}(h')\circ\tau\in 
{\rm Aut}_{\overline{F}}({\bf H})$ et 
$\tau^*_1={\rm Int}_{{\bf H}'}(h')\circ \tau^*\in 
{\rm Aut}_{\overline{F}}({\bf H}')$. Alors on a 
$\frak{H}^1_{\tau_1}= 
\frak{H}'^1_{\tau^*_1}\subset \frak{H}'=\frak{H}_\tau^1$. Mais  
puisque $\tau$ est r\'egulier, on a $\dim_{\overline{F}}(\frak{H}_{\tau_1}^1)\geq 
\dim_{\overline{F}}(\frak{H}')$. Par 
cons\'equent $\frak{H}'^1_{\tau_1^*} = \frak{H}'$. D'autre part, 
comme l'ensemble
des $x\in {\bf H}'$ tels que $D_{{\bf H}'}({\rm Int}_{{\bf H}'}(x)\circ \tau^*)\neq 0$ 
est ouvert dense dans ${\bf H}'$, on peut supposer $\tau^*_1$ r\'egulier. 
Cela implique que $\tau^*$ lui-m\^eme est r\'egulier. Comme $\tau^*$ est unipotent, 
d'apr\`es le lemme 1, ${\bf H}'$ est un tore. Puisqu'un $\overline{F}$--automorphisme d'un tore 
est automatiquement quasi--semisimple, d'apr\`es la proposition de \ref{automorphismes qss lf}, 
on obtient que $\tau$ est quasi--semisimple. 

Revenons au cas g\'en\'eral: $\tau\in {\rm Aut}_{\overline{F}}({\bf H})$. 
Choisissons une paire 
de Borel $\tau$--stable $({\bf B},{\bf T})$ de ${\bf H}$. Posons 
${\bf H}^\sharp={\bf H}_\tau^\circ$ et 
${\bf T}^\sharp={\bf T}\cap {\bf H}^\sharp$. 
D'apr\`es le thorme 2 de \ref{automorphismes qss}, le groupe 
${\bf H}^\sharp$ est r\'eductif, ${\bf T}^\sharp$ est 
un tore maximal de ${\bf H}^\sharp$, et ${\bf T}=Z_{\bf H}({\bf T}^\sharp)$.

Montrons que 
${\bf H}_\tau^\circ={\bf T}^\sharp$ et que $\frak{H}_\tau^1$ est contenu dans ${\rm Lie}({\bf T})$. 
Comme plus haut, pour $h\in {\bf H}^\sharp$ et $x\in \smash{\overline{F}}^\times$, 
on a ${\rm Int}_{\bf H}(h)\circ \tau (\frak{H}_\tau^x)=\frak{H}_\tau^x$, et si $h$ est suffisamment 
proche de $1$ dans ${\bf H}^\sharp$, alors on a $\frak{H}^1_{{\rm Int}_{\bf H}(h)\circ \tau}\subset 
\frak{H}^1_\tau$. D'apr\`es la d\'emonstration du lemme 1, il existe un $t\in {\bf T}^\sharp$ tel que 
$\frak{H}^1_{{\rm Int}_{\bf H}(t)\circ \tau}\subset 
\frak{H}^1_\tau$ et ${\bf H}_t^\circ=Z_{\bf H}({\bf T}^\sharp)={\bf T}$. Alors on a
$$
\frak{H}^1_{{\rm Int}_{\bf H}(t)\circ \tau}= \frak{H}_\tau^1\cap \frak{H}_{{\rm Int}_{\bf H}(t)}^1
=\frak{H}_\tau^1 \cap {\rm Lie}({\bf T}), 
$$
et $\dim_{\overline{F}}(\frak{H}^1_{{\rm Int}_{\bf H}(t)\circ \tau})\leq 
\dim_{\overline{F}}(\frak{H}^1_\tau)$ avec \'egalit\'e si et seulement si $\frak{H}_\tau^1\subset 
{\rm Lie}({\bf T})$. Puisque $\tau$ est r\'egulier, on a donc $\frak{H}_\tau^1\subset 
{\rm Lie}({\bf T})$. Par cons\'equent 
$${\rm Lie}({\bf H}^\sharp)\subset {\rm Lie}({\bf T})\cap 
{\rm Lie}({\bf H}^\sharp)={\rm Lie}({\bf T}^\sharp),
$$ ce qui n'est possible que si ${\bf H}^\sharp={\bf T}^\sharp$. Cela ach\`eve la d\'emonstration du point (1).

Montrons (2). Puisque ${\bf H}_\tau^\circ\subset 
{\bf T}$ et ${\bf H}_\tau/{\bf H}_\tau^\circ$ est form\'e 
d'\'el\'ements semisimples (thorme 2 de \ref{automorphismes qss}), le morphisme $(1-\tau)\vert_{\bf U}$ est injectif. 
Comme d'apr\`es le point (1) --- ou d'apr\`es la proposition de \ref{automorphismes qss} ---, on a
$$
{\rm Lie}({\bf U})\cap \frak{H}_\tau=
\{0\}= {\rm Lie}({\bf U}_\tau),
$$ 
ce morphisme est un automorphisme.

Supposons $\tau$ localement fini, 
et montrons (3). 
\'Ecrivons la d\'ecomposition de Jordan $\tau=\tau_{\rm s}\circ \tau_{\rm u}$. Puisque 
$\tau({\bf B},{\bf T})=({\bf B},{\bf T})$, on a 
$\tau_{\rm s}({\bf B},{\bf T})=({\bf B},{\bf T})$ et  $\tau_{\rm u}({\bf B},{\bf T})=({\bf B},{\bf T})$. 
D'apr\`es le premier paragraphe de la d\'emonstration, 
${\bf H}_{\tau_{\rm s}}^\circ$ est un tore, et d'apr\`es le thorme 2 de \ref{automorphismes qss}, on a 
$Z_{\bf H}({\bf H}^\circ_{\tau_{\rm s}})={\bf T}=Z_{\bf H}({\bf H}_\tau^\circ)$. 
Montrons que $\tau_{\rm s}$ est r\'egulier. 
Soit un \'el\'ement $h\in {\bf H}$ tel que le 
$\overline{F}$--automorphisme $\sigma={\rm Int}_{\bf H}(h)\circ \tau_{\rm s}$ de ${\bf H}$ 
est r\'egulier. Puisque $\sigma$ est quasi--semisimple (d'aprs le point (1)), il existe un 
$x\in {\bf H}$ tel que $\sigma$ stabilise la paire de Borel ${\rm Int}_{\bf H}(x)({\bf B},{\bf T})$ de 
${\bf H}$. D'apr\`es la relation $(*)$ de \ref{automorphismes}, quitte \`a remplacer $h$ par $x^{-1}h\tau_{\rm s}(x)$, on peut 
supposer que $\sigma({\bf B},{\bf T})=({\bf B},{\bf T})$. Alors on a 
${\rm Int}_{\bf H}(h)({\bf B},{\bf T})=({\bf B},{\bf T})$, donc $h\in {\bf T}$. D'apr\`es le point (1), ${\bf H}_{\sigma}^\circ$ est un tore 
et $Z_{\bf H}({\bf H}_{\sigma}^\circ)={\bf T}$. On a donc
$$
{\bf H}_\sigma^\circ ={\bf H}_\sigma^\circ \cap {\bf T}= {\bf H}_{\tau_{\rm s}}^\circ \cap {\bf T}
={\bf H}_{\tau_{\rm s}}^\circ.
$$
Comme $\dim_{\overline{F}}(\frak{H}^1_\sigma)\geq \dim({\bf H}_\sigma^\circ)
=\dim({\bf H}_{\tau_{\rm s}}^\circ)= \dim_{\overline{F}}(\frak{H}^1_{\tau_{\rm s}})$, on obtient 
que $\tau_{\rm s}$ est r\'egulier. D'o\`u l'\'egalit\'e
$$D_{\bf H}(\tau)=D_{\bf H}(\tau_{\rm s}),
$$ puisque 
$\frak{H}_\tau^1=\frak{H}_{\tau_{\rm s}}$.
\end{proof}

\begin{monlem3}
Soit $\tau\in {\rm Aut}_{\overline{F}}({\bf H})$ quasi--semisimple, 
et soit $({\bf B},{\bf T})$ une paire de Borel $\tau$--stable de ${\bf H}$.
\begin{enumerate}
\item[(1)]Il existe un 
$t\in {\bf T}$ tel que ${\rm Int}_{\bf H}(t)\circ \tau$ est r\'egulier.
\item[(2)]Soit $t\in {\bf T}$ tel que ${\rm Int}_{\bf H}(t)\circ \tau$ est r\'egulier. Alors 
on a ${\bf H}_{{\rm Int}_{\bf H}(t)\circ \tau}^\circ ={\bf T}_\tau^\circ$, et pour 
tout $u\in R_{\rm u}({\bf B})$, ${\rm Int}_{\bf H}(tu)\circ \tau$ est r\'egulier.
\end{enumerate}
\end{monlem3}

\begin{proof} 
La d\'emonstration du point (1) est identique \`a celle du lemme de \ref{automorphismes qc}.

Montrons (2). Soit $t\in {\bf T}$ tel que ${\rm Int}_{\bf H}(t)\circ \tau$ est r\'egulier. 
Posons $\tau'={\rm Int}_{\bf H}(t)\circ \tau$. Puisque 
$\tau'({\bf B},{\bf T})=({\bf B},{\bf T})$, 
d'apr\`es le lemme 1, ${\bf H}_{\tau'}^\circ$ est un sous--tore de 
${\bf T}$. Par cons\'equent 
${\bf H}_{\tau'}^\circ={\bf T}_{\tau'}^\circ = {\bf T}_\tau^\circ$. Soit 
$u\in {\bf U}=R_{\rm u}({\bf B})$, et 
posons $u'= tut^{-1}\in {\bf U}$. Soit $v\in {\bf U}$ tel 
que $u' = v^{-1}\tau'(v)$. Alors on a
$$
{\rm Int}_{\bf H}(tu)\circ \tau = {\rm Int}_{\bf H}(v^{-1})\circ {\rm Int}_{\bf H}(\tau'(v))\circ \tau'
=  {\rm Int}_{\bf H}(v^{-1})\circ \tau' \circ  {\rm Int}_{\bf H}(v),
$$
et ${\rm Int}_{\bf H}(tu)\circ \tau$ est r\'egulier.
\end{proof}

\begin{montheo}
Soit $\tau\in {\rm Aut}_{\overline{F}}({\bf H})$. 
\begin{enumerate}
\item[(1)]$\tau$ est r\'egulier si et seulement si 
$\tau$ est quasi--semisimple et ${\bf H}_\tau^\circ$ est un tore.
\item[(2)]Supposons $\tau$ quasi--semisimple. Alors $\tau$ est r\'egulier si et seulement si 
$$\dim({\bf H}_\tau^\circ)=r_\tau({\bf H}).$$
\item[(3)]Supposons $\tau$ localement fini. Alors $\tau$ est r\'egulier si et seulement si $\tau_{\rm s}$ 
est r\'egulier, i.e. si et seulement si ${\bf H}_{\tau_{\rm s}}^\circ$ est un tore. En d'autres termes, on a
$$
D_{\bf H}(\tau)=D_{\bf H}(\tau_{\rm s}).
$$ 
De plus, si $\tau$ est r\'egulier, alors on a $r_{\tau_{\rm s}}({\bf H})=r_\tau^1({\bf H})$.
\end{enumerate}
\end{montheo}

\begin{proof} 
Si $\tau$ est quasi--semisimple, on choisit une paire de Borel $({\bf B},{\bf T})$ de 
${\bf G}$ telle que $\tau({\bf B},{\bf T})= ({\bf B},{\bf T})$, et l'on pose ${\bf U}=R_{\rm u}({\bf B})$.

Montrons (1). Si $\tau$ est r\'egulier, alors d'apr\`es le lemme 2, $\tau$ est quasi--semisimple et 
${\bf H}_\tau^\circ$ est un tore; de plus (loc. cit.) $Z_{\bf H}({\bf H}_\tau^\circ)={\bf T}$ et 
$\frak{H}_\tau^1$ est contenu 
dans ${\rm Lie}({\bf T})$. R\'eciproquement, supposons 
que $\tau$ est quasi--semisimple et que ${\bf H}_\tau^\circ$ est un tore. 
On a donc $\tau({\bf B},{\bf T})=({\bf B},{\bf T})$. 
Puisque ${\bf H}_\tau^\circ$ est 
un tore, d'apr\`es le thorme 2 de \ref{automorphismes qss}, on a ${\bf H}^\circ_\tau\subset {\bf T}$ et 
$Z_{\bf H}({\bf H}_\tau^\circ)$ est un tore maximal 
($\tau$--admissible) de ${\bf H}$, par cons\'equent $Z_{\bf H}({\bf H}_\tau^\circ)={\bf T}$. 
Le morphisme $(1-\tau)\vert_{\bf U}$ est injectif (cf. la d\'emonstration du point (2) du lemme 2), 
et puisqu'il est s\'eparable (\ref{automorphismes qss}, proposition), c'est un automorphisme de vari\'et\'e alg\'ebrique. 
On en d\'eduit que
$$\frak{H}_\tau^1\cap {\rm Lie}({\bf U})=\{0\}.
$$
Le m\^eme 
raisonnement appliqu\'e au radical unipotent ${\bf U}^-$ du sous--groupe de Borel ${\bf B}^-$ de 
${\bf H}$ oppos\'e \`a ${\bf B}$ par rapport \`a ${\bf T}$, entra\^{\i}ne que 
$$\frak{H}_\tau^1\cap {\rm Lie}({\bf U}^-)=\{0\}.
$$
Par cons\'equent $\frak{H}_\tau^1$ est contenu dans ${\rm Lie}({\bf T})$. D'apr\`es 
le lemme 3, il existe un $t\in {\bf T}$ tel que $\tau_1= {\rm Int}_{\bf H}(t)\circ\tau$ 
est r\'egulier. Puisque 
$\tau_1({\bf B},{\bf T})=({\bf B},{\bf T})$, on a $Z_{\bf H}({\bf H}_{\tau_1}^\circ)={\bf T}$ et 
$\frak{H}_{\tau_1}^1\subset {\rm Lie}({\bf T})$. On en d\'eduit que
$$
\frak{H}_{\tau}^1 = \frak{H}_\tau^1\cap {\rm Lie}({\bf T})= \frak{H}_{\tau_1}^1
\cap {\rm Lie}({\bf T})= \frak{H}_{\tau_1}^1.
$$
Par suite $\dim_{\overline{F}}(\frak{H}_{\tau}^1)=r^1_\tau({\bf H})$ et $\tau$ est r\'egulier.

Supposons $\tau$ quasi--semisimple, et montrons (2). D'apr\`es (1), il s'agit de montrer que ${\bf H}_\tau^\circ$ est 
un tore si et seulement si $\dim({\bf H}_\tau^\circ)=r_\tau({\bf H})$. Supposons tout d'abord 
que ${\bf H}_\tau^\circ$ est un tore. Comme ${\bf H}_\tau^\circ \subset {\bf T}$ (cf. la d\'emonstration 
du point (1)), on a 
${\bf H}_\tau^\circ = {\bf T}_\tau^\circ$. Soit un \'el\'ement $h'\in {\bf H}$ tel que $\tau'={\rm Int}_{\bf H}(h')\circ \tau$ est quasi--semisimple. On veut montrer que $\dim({\bf H}_{\tau'}^\circ)\geq \dim({\bf H}_\tau^\circ)$. Il existe un $x\in {\bf H}$ tel 
$\tau'$ stabilise la paire de Borel ${\rm Int}_{\bf H}(x)({\bf B},{\bf T})$ de ${\bf H}$. Puisque 
$\tau''={\rm Int}_{\bf H}(x^{-1})\circ \tau'\circ {\rm Int}_{\bf H}(x)$ stabilise $({\bf B},{\bf T})$ et $\dim({\bf H}_{\tau'}^\circ)=\dim({\bf H}_{\tau''}^\circ)$, 
quitte \`a remplacer $\tau'$ par $\tau''$, on peut supposer que 
$\tau'({\bf B},{\bf T})=({\bf B},{\bf T})$. Alors ${\rm Int}_{\bf H}(h)({\bf B},{\bf T})=
({\bf B},{\bf T})$, i.e. $h\in {\bf T}$. Comme
$$
{\bf H}_{\tau'}^\circ \supset {\bf T}_{\tau'}^\circ = {\bf T}^\circ_\tau={\bf H}^\circ_\tau,
$$
on a $\dim({\bf H}^\circ_{\tau'})\geq \dim({\bf H}^\circ_\tau)$. Supposons maintenant que $\dim({\bf H}_\tau^\circ)=r_\tau({\bf H})$, et montrons que le groupe r\'eductif connexe 
${\bf H}'={\bf H}_\tau^\circ$ est un tore. Soit ${\bf T}'={\bf T}\cap {\bf H}'$ (c'est un 
tore maximal de ${\bf H}'$). Soit $t'\in {\bf T}'$, et posons $\tau'={\rm Int}_{\bf H}(t')\circ \tau'$. 
D'apr\`es la d\'emonstration du point (1) du lemme 2, si $t'$ est suffisamment proche de $1$ 
dans ${\bf T}'$, ce que l'on suppose, alors $\frak{H}_{\tau'}=\frak{H}_{\tau'}^1$ est contenu dans 
${\rm Lie}({\bf H}')$. Puisque $\frak{H}_{\tau'}={\rm Lie}({\bf H}^\circ_{\tau'})$, 
cela entra\^{\i}ne que ${\bf H}_{\tau'}^\circ$ est contenu dans ${\bf H}'$, puis que ${\bf H}_{\tau'}^\circ={\bf H}'^\circ_{{\rm Int}_{{\bf H}'}(t')}$. Supposons de plus que $t'$ 
est r\'egulier dans 
${\bf H}'$, i.e. que $D_{{\bf H}'}({\rm Int}_{{\bf H}'}(t'))\neq 0$. Alors 
d'aprs la proposition de \ref{automorphismes rguliers; le cas i}, ${\bf H}'^\circ_{{\rm Int}_{{\bf H}'}(t')}$ est un tore. D'autre part comme 
$\dim({\bf H}_{\tau'}^\circ)\leq \dim({\bf H}_\tau^\circ)=r_\tau({\bf H})$, par d\'efinition 
de $r_\tau({\bf H})$, on a l'\'egalit\'e 
${\bf H}_{\tau'}^\circ ={\bf H}_\tau^\circ$. Cela qui ach\`eve la d\'emonstration du point (2).
 
Supposons $\tau$ localement fini, et montrons (3). 
\'Ecrivons la d\'ecomposition de Jordan $\tau=\tau_{\rm s}\circ \tau_{\rm u}$. D'apr\`es le lemme 2, si 
$\tau$ est r\'egulier, alors $\tau_{\rm s}$ l'est aussi et on a l'galit
$$
D_{\bf H}(\tau)=D_{\bf H}(\tau_{\rm s}).
$$
D'apr\`es le point (1), $\tau_{\rm s}$ est 
r\'egulier si et seulement si ${\bf H}_{\tau_{\rm s}}^\circ$ est un tore: $\tau_{\rm s}$ 
est semisimple donc a fortiori quasi--semisimple (thorme 1 de \ref{automorphismes qss lf}). Supposons $\tau_{\rm s}$ 
r\'egulier et montrons que $\tau$ est r\'egulier. Comme ${\bf H}_{\tau_{\rm s}}^\circ$ est un tore, $\tau^*\in {\rm Aut}_{\overline{F}}({\bf H}_{\tau_{\rm s}}^\circ)$ est quasi--semisimple, par cons\'equent $\tau$ est quasi--semisimple (proposition de \ref{automorphismes qss lf}). Puisque ${\bf H}_\tau^\circ\subset {\bf H}_{\tau_{\rm s}}^\circ$ et 
${\bf H}_{\tau_{\rm s}}^\circ$ est un tore, ${\bf H}_\tau^\circ$ est un tore. D'apr\`es la relation $(*)$ de \ref{automorphismes rguliers; le cas i}, on a $\dim({\bf H}_{\tau_{\rm s}}^\circ)=r_\tau^1({\bf H})$, et 
d'apr\`es le point (2), on a $\dim({\bf H}_{\tau_{\rm s}}^\circ)=r_{\tau_{\rm s}}({\bf H})$, d'o\`u la derni\`ere assertion du point (3). Cela ach\`eve la 
d\'emonstration du th\'eor\`eme.
\end{proof}

\begin{moncoro}
Soit $\tau\in {\rm Aut}_{\overline{F}}({\bf H})$. On a $r_\tau({\bf H})\leq  r^1_\tau({\bf H})$ avec \'egalit\'e si et seulement si 
$\dim({\bf H}_{\tau'}^\circ)=\dim_{\overline{F}}(\frak{H}_{\tau'}^1)$ pour un 
(i.e. pour tout) $h\in {\bf H}$ tel que $\tau'={\rm Int}_{\bf H}(h)\circ\tau$ est r\'egulier.
\end{moncoro}

\begin{proof} Pour $h\in {\bf H}$, posant $\tau'={\rm Int}_{\bf H}(h)\circ \tau$, on a 
$$\dim({\bf H}_{\tau'}^\circ)\leq \dim_{\overline{F}}(\frak{H}_{\tau'}^1)\leq r^1_{\tau'}({\bf H}),$$
d'o\`u 
l'in\'egalit\'e $r_\tau({\bf H})\leq r_\tau^1({\bf H})$ puisqu'on sait (lemme 2) que si $\tau'$ est r\'egulier alors $\tau'$ est quasi--semisimple. Si $\tau'$ est r\'egulier, alors $\dim({\bf H}_{\tau'}^\circ)=r_\tau({\bf H})$ et 
$\dim_{\overline{F}}(\frak{H}_{\tau'}^1)=r_\tau^1({\bf H})$. D'o\`u le corollaire.
\end{proof}

\begin{marema1}{\rm 
Pour $\tau\in 
 {\rm Aut}_{\overline{F}}({\bf H})$, puisque
 $${\rm Lie}({\bf H}_{\tau}^\circ)\subset \frak{H}_{\tau}\subset 
\frak{H}_{\tau}^1,
$$ on a $\dim({\bf H}_{\tau}^\circ)=\dim_{\overline{F}}(\frak{H}_{\tau}^1)$ si et seulement si 
les inclusions ci-dessus sont des \'egalit\'es. D'autre part si $\tau$ est semisimple, on a ${\rm Lie}({\bf H}_\tau^\circ)=\frak{H}_\tau=\frak{H}_\tau^1$. On en d\'eduit que 
s'il existe un \'el\'ement $h\in {\bf H}$ tel que ${\rm Int}_{\bf H}(h)\circ\tau$ est 
semisimple r\'egulier, alors $r_\tau({\bf H})=r_\tau^1({\bf H})$.\hfill $\blacksquare$
}\end{marema1}

\v2 Si $p=1$, l'absence de $\overline{F}$--automorphismes 
quasi--semisimples unipotents non triviaux simplifie notablement la situation, 
du moins pour les automorphismes localement finis: 

\begin{monlem4}Supposons $p=1$, et soit 
$\tau\in {\rm Aut}_{\overline{F}}^0({\bf H})$. On a 
$r_\tau({\bf H})=r^1_\tau({\bf H})$.
\end{monlem4}

\begin{proof}
Puisque $\tau$ est localement fini, 
pour $h\in {\bf H}$,  
$\tau'={\rm Int}_{\bf H}(h)\circ\tau$ est encore localement fini, par cons\'equent 
si $\tau'$ est r\'egulier, alors $\tau'$ est semisimple (d'aprs la remarque 2 de \ref{automorphismes qss lf}). 
On conclut gr\^ace \`a la remarque 1 ci--dessus.
\end{proof}

\begin{marema2}{\rm Si $\tau$ n'est pas localement fini, 
alors l'in\'egalit\'e $r_\tau({\bf H})\leq r^1_\tau({\bf H})$ est en g\'en\'eral stricte (m\^eme si $p=1$). En effet, soit ${\bf T}$ le tore $\Bbb{G}_{\rm m}\times \Bbb{G}_{\rm m}$, et soit 
$\tau$ le $\overline{F}$--automorphisme de ${\bf T}$ d\'efini par $\tau(a,b)=\tau(ab,b)$ 
($a,\,b\in \smash{\overline{F}}^\times$); il n'est pas localement fini. 
Posons $\frak{T}={\rm Lie}({\bf T})\;(=\overline{F}\times \overline{F})$.
On a ${\rm Lie}(\tau)(X,Y)=(X+Y,Y)$ ($X,\,Y\in \overline{F}$).  
Par suite pour 
$h\in {\bf T}$, posant $\tau'={\rm Int}_{\bf H}(h)\circ\tau$, on a ${\rm Lie}({\bf T}_{\tau'})= 
\frak{T}_{\tau'}=\{(X,0): X\in \overline{F}\}$ et $\frak{T}_{\tau'}^1= \frak{T}$. Donc 
$r_\tau({\bf T})=1$ et $r^1_\tau({\bf T})=2$.\hfill $\blacksquare$}
\end{marema2}

\subsection{\'El\'ements r\'eguliers d'un ${\bf H}$--espace tordu}\label{éléments réguliers d'un espace tordu}
Il est commode de 
reformuler les r\'esultats de \ref{automorphismes rguliers; le cas i} et \ref{automorphismes rguliers; le cas g} en termes d'espaces tordus (voir 
\ref{espaces topo tordus}; on remplace ici la 
cat\'egorie des groupes topologiques par celle des groupes alg\'ebriques).  Comme dans \cite[I.3]{La}, 
appelons {\it ${\bf H}$--espace (alg\'ebri\-que) tordu} la donn\'ee:
\begin{itemize}
\item d'un ${\bf H}$--espace principal homog\`ene ${\bf H}^\natural$, i.e. une vari\'et\'e alg\'ebrique affine ${\bf H}^\natural$ munie d'une action 
alg\'ebrique de ${\bf H}$ \`a gauche, not\'ee ${\bf H}\times {\bf H}^\natural\rightarrow {\bf H}^\natural,\, (h,\delta)\mapsto h\cdot \delta$, 
telle que pour tout $\delta\in {\bf H}^\natural$, l'application ${\bf H}\rightarrow {\bf H}^\natural,\,h\mapsto h\cdot \delta$ 
est un isomorphisme de vari\'et\'es alg\'ebriques;
\item et d'une application ${\rm Int}_{{\bf H}^\natural}:{\bf H}^\natural\rightarrow {\rm Aut}_{\overline{F}}({\bf H})$ 
telle que
$$
{\rm Int}_{{\bf H}^\natural}(h\cdot \delta)={\rm Int}_{\bf H}(h)\circ {\rm Int}_{{\bf H}^\natural}(\delta)\quad (h\in {\bf H},\delta\in {\bf H}^\natural).
$$
\end{itemize}
On peut alors d\'efinir une action alg\'ebrique ${\bf H}^\natural\times {\bf H}\rightarrow 
{\bf H}^\natural,\, (\delta,h)\mapsto \delta\cdot h$ de ${\bf H}$ sur 
${\bf H}^\natural$ \`a droite, qui commute \`a l'action \`a gauche:
$$
 \delta\cdot h= {\rm Int}_{{\bf H}^\natural}(\delta)(h)\cdot \delta.
$$
Pour $h\in {\bf H}$, on note ${\rm Int}'_{{\bf H}}(h):{\bf H}^\natural\rightarrow {\bf H}^\natural$ 
l'automorphisme de vari\'et\'e alg\'ebrique d\'efini par
$$
{\rm Int}'_{\bf H}(\delta)=h\cdot \delta \cdot h^{-1}.
$$

Pour $\tau\in {\rm Aut}_{\overline{F}}({\bf H})$, on peut comme en \ref{espaces topo tordus} d\'efinir le 
${\bf H}$--espace tordu ${\bf H}\tau$. Il s'agit d'une vari\'et\'e 
alg\'ebrique affine, munie:
\begin{itemize}
\item d'un isomorphisme de vari\'et\'es 
alg\'ebriques ${\bf H}\rightarrow {\bf H}\tau,\, h\mapsto h\tau$, et 
d'actions alg\'ebriques de ${\bf H}$ \`a gauche 
${\bf H}\times {\bf H}\tau\rightarrow {\bf H}\tau,\, (h, h'\tau)\mapsto h\cdot h'\tau$ et 
\`a droite 
${\bf H}\tau \times {\bf H} \rightarrow {\bf H}\tau,\, (h'\tau,h)\mapsto h'\tau\cdot h$, 
commutant entre elles et 
v\'erifiant l'\'egalit\'e
$$
h\cdot h'\tau\cdot h'' =hh'\tau(h'')\tau\quad (h,\,h',\,h''\in {\bf H});
$$
\item et d'une application ${\rm Int}_{{\bf H}\tau}:{\bf H}\tau\rightarrow {\rm Aut}_{\overline{F}}({\bf H})$ 
donn\'ee par
$$
{\rm Int}_{{\bf H}\tau}(h\tau)={\rm Int}_{\bf H}(h)\circ \tau\quad (h\in {\bf H}).
$$
\end{itemize}

Mutatis mutandis, le ${\rm n}^\circ$ \ref{espaces topo tordus} reste vrai dans ce contexte. En particulier, si 
${\bf H}^\natural$ est un ${\bf H}$--espace tordu, le choix d'un \'el\'ement $\delta_1\in {\bf H}^\natural$ 
permet d'identifier ${\bf H}^\natural$ et ${\bf H}\tau_1$ o\`u $\tau_1={\rm Int}_{{\bf H}^\natural}(\delta_1)$, mais 
d'apr\`es la remarque 2 de \ref{espaces topo tordus}, cette identification n'est pas canonique.

\v1 
Soit ${\bf H}^\natural$ un ${\bf H}$--espace tordu. 
On note ${\rm Ad}_{{\bf H}^\natural}:{\bf H}^\natural\rightarrow {\rm GL}(\frak{H})$ l'application $\delta\mapsto 
{\rm Lie}({\rm Int}_{{\bf H}^\natural}(\delta))$. On a donc
$$
{\rm Ad}_{{\bf H}^\natural}(h\cdot\delta)={\rm Ad}_{\bf H}(h)\circ {\rm Ad}_{{\bf H}^\natural}(\delta)\quad (h\in {\bf H},\,\delta\in {\bf H}^\natural).
$$
On note $D_{{\bf H}^\natural}:{\bf H}^\natural\rightarrow \overline{F}$ la fonction d\'efinie par
$$
D_{{\bf H}^\natural}(\delta)=D_{\bf H}({\rm Int}_{{\bf H}^\natural}(\delta)).
$$
Pour $\delta\in {\bf H}^\natural$ et $\tau={\rm Int}_{{\bf H}^\natural}(\delta)$, on remplace l'indice $\tau$ 
par un indice $\delta$ dans toutes les notations pr\'ec\'edentes --- i.e. on pose 
${\bf H}_\delta^\circ ={\bf H}_\tau^\circ$, 
$\frak{H}_\delta^1=\frak{H}^1_\tau$, (etc.). Puisque ${\rm Int}_{{\bf H}^\natural}(h\cdot \delta)= 
{\rm Int}_{\bf H}(h)\circ \tau$ 
($h\in {\bf H}$), on a $\tau\in {\rm Aut}_{\overline{F}}^0({\bf H})$ si et seulement si 
${\rm Int}_{{\bf H}^\natural}(\delta')\in {\rm Aut}_{\overline{F}}^0({\bf H})$ pour tout $\delta'\in {\bf H}^\natural$; 
auquel cas on dit que ${\bf H}^\natural$ est {\it localement fini}. Par d\'efinition, 
les entiers $r_\delta({\bf H})$ et $r_\delta^1({\bf H})$ ne d\'ependent pas 
de l'\'el\'ement $\delta\in {\bf H}^\natural$. On 
les note $r({\bf H}^\natural)$ et $r^1({\bf H}^\natural)$. L'entier  
$r({\bf H}^\natural)$ est appel\'e {\it rang (absolu) de ${\bf H}^\natural$}.

Un \'el\'ement $\delta\in {\bf H}^\natural$ est dit {\it quasi--semisimple} (resp. {\it quasi--central}, 
{\it semisimple}, {\it unipotent}, {\it r\'egulier}) si le $\overline{F}$--automorphisme 
${\rm Int}_{{\bf H}^\natural}(\delta)$ de ${\bf H}$ est quasi--semisimple (resp. quasi--central, semisimple, 
unipotent, r\'egulier). Ces notions sont stables par ${\bf H}$--conjugaison. Notons que si 
${\bf H}^\natural$ n'est pas localement fini, alors aucun \'el\'ement de ${\bf H}^\natural$ 
n'est semisimple (resp. unipotent).

On note ${\bf H}^\natural_{\rm reg}$ l'ensemble des 
\'el\'ements r\'eguliers de ${\bf H}^\natural$. Par d\'efinition, la fonction
$$
{\bf H}^\natural\rightarrow \overline{F},\,\delta\mapsto D_{{\bf H}^\natural}(\delta)
$$ 
appartient \`a l'alg\`ebre affine $\overline{F}[{\bf H}^\natural]$. Puisque 
l'ensemble ${\bf H}^\natural_{\rm reg}$ est non vide, il est 
ouvert dense dans ${\bf H}^\natural$. D'apr\`es le thorme de \ref{automorphismes rguliers; le cas g}, 
un \'el\'ement $\delta\in {\bf H}^\natural$ est 
r\'egulier si et seulement s'il est quasi--semisimple et ${\bf H}_\delta^\circ$ est un tore. 
D'apr\`es le corollaire de \ref{automorphismes rguliers; le cas g}, on a
$$
r({\bf H}^\natural)\leq r^1({\bf H}^\natural)\leqno{(*)}
$$
avec \'egalit\'e si et seulement s'il existe un \'el\'ement $\delta\in {\bf H}^\natural_{\rm reg}$ tel que 
$\dim ({\bf H}_\delta^\circ)= \dim_{\overline{F}}(\frak{H}_\delta^1)$.

\subsection{Tores maximaux et sous--espaces de Cartan d'un ${\bf H}$--espace tordu}\label{tm et sec d'un espace tordu} 
Soit ${\bf H}^\natural$ un ${\bf H}$--espace tordu. 
Soit $\delta\in {\bf H}^\natural$, et posons 
$\tau={\rm Int}_{{\bf H}^\natural}(\delta)$. Soit $({\bf B},{\bf T})$ 
une paire de Borel de ${\bf H}$ telle que $\tau({\bf B})={\bf B}$. 
Posons ${\bf U}=R_{\rm u}({\bf B})$, et soit $u_0$ l'unique \'el\'ement de 
${\bf U}$ tel que ${\rm Int}_{\bf H}(u_0)\circ \tau$ stabilise $({\bf B},{\bf T})$. 
Posons $\delta_0=u_0\cdot \delta$ et $\tau_0={\rm Int}_{\bf H}(\delta_0)$. 
Puisque $\tau_0({\bf B},{\bf T})=({\bf B},{\bf T})$, $\delta_0$ est quasi--semisimple. 
Posons ${\bf S} ={\bf T}_{\delta_0}^\circ\;( ={\bf T}\cap {\bf H}_{\delta_0}^\circ)$. 
D'apr\`es le lemme 3 de \ref{automorphismes rguliers; le cas g}, on a la

\begin{mapropo}
\begin{enumerate}
\item[(1)]Tout \'el\'ement r\'egulier de ${\bf H}^\natural$ est de la forme 
$x^{-1}\cdot (t\cdot \delta_0) \cdot x$ pour 
un $t\in {\bf T}$ et un $x\in {\bf H}$.
\item[(2)] L'ensemble $({\bf T}\cdot \delta_0)\cap {\bf H}^\natural_{\rm reg}$ est 
ouvert dense dans ${\bf T}\cdot \delta_0$.
\item[(3)] Pour $\delta\in ({\bf T}\cdot \delta_0)\cap {\bf H}^\natural_{\rm reg}$, on a  ${\bf H}_\delta^\circ 
={\bf S}$.
\end{enumerate}
\end{mapropo}

Partons maintenant d'un élément $\delta'_0\in {\bf H}^\natural_{\rm reg}$, et posons 
${\bf S}^\natural={\bf H}_{\delta'_0}^\circ\cdot \delta'_0$. Le centralisateur
$$Z_{\bf H}({\bf S}^\natural)= 
\{h\in {\bf H}: {\rm Int}_{{\bf H}^\natural}(\gamma)(h)=h,\,\forall \gamma\in {\bf S}^\natural\}
$$
est un tore de ${\bf H}$, disons ${\bf S}$. Pour $\delta\in {\bf S}^\natural$, on a ${\bf S}^\natural={\bf S}\cdot \delta$ et ${\rm Int}_{{\bf H}^\natural}(\delta)\vert_{\bf S}={\rm id}_{\bf S}$. Cela munit ${\bf S}^\natural$ d'une structure de ${\bf S}$--espace tordu trivial.  
Le groupe ${\bf T}=Z_{\bf H}({\bf S})$ est un tore maximal de ${\bf H}$: 
pour $\delta\in {\bf S}^\natural\cap {\bf H}^\natural_{\rm reg}$, on a 
${\bf H}_\delta^\circ = {\bf S}$ --- en particulier on a ${\bf H}_{\delta'_0}^\circ = {\bf S}$ --- et ${\bf T}$ est l'unique tore maximal ${\rm Int}_{{\bf H}^\natural}(\delta)$--admissible de 
${\bf H}$. Posons
$$
{\bf T}^\natural = {\bf T}\cdot {\bf S}^\natural = {\bf S}^\natural \cdot {\bf T}\subset {\bf H}^\natural.
$$
Pour $\delta\in {\bf T}^\natural $, on a ${\bf T}^\natural={\bf T}\cdot \delta$ et 
${\rm Int}_{{\bf T}^\natural}(\delta)={\rm Int}_{{\bf H}^\natural}(\delta)\vert_{\bf T}$ est un $\overline{F}$--automorphisme 
de ${\bf T}$. Cela munit 
${\bf T}^\natural$ d'une structure de 
${\bf T}$--espace tordu. Pour $\delta\in {\bf T}^\natural\cap {\bf H}_{\rm reg}^\natural$, on a 
${\bf H}_\delta^\circ ={\bf S}$. De plus, on a
$Z({\bf T}^\natural)={\bf S}$ o\`u 
(rappel) $Z({\bf T}^\natural)$ est le centralisateur de ${\bf T}^\natural$ dans ${\bf T}$, d\'efini par 
$$Z({\bf T}^\natural)=\{t\in {\bf T}:{\rm Int}_{{\bf T}^\natural}(\gamma)(t)=t,\,\forall \gamma\in {\bf T}^\natural\}.
$$
Puisque
$$Z_{\bf H}({\bf T}^\natural)\subset 
Z_{\bf H}({\bf S}^\natural)= 
{\bf S}=Z({\bf T}^\natural)\subset Z_{\bf H}({\bf T}^\natural),
$$
les deux inclusions ci--dessus sont des \'egalit\'es. 

\begin{mesdefi}
{\rm On appelle:
\begin{itemize}
\item {\rm tore maximal de 
${\bf H}^\natural$} une partie ${\bf S}^\natural$ de la forme ${\bf S}^\natural={\bf H}^\circ_\delta\cdot \delta$ pour un $\delta\in {\bf H}^\natural_{\rm reg}$;
\item {\rm sous--espace de Cartan de ${\bf H}^\natural$} une partie ${\bf T}^\natural$ de la forme 
${\bf T}^\natural=Z_{\bf H}({\bf H}^\circ_\delta)\cdot \delta$ pour un 
$\delta\in {\bf H}^\natural_{\rm reg}$.
\end{itemize}
}\end{mesdefi}

D'apr\`es la proposition, les tores maximaux (resp. sous--espaces de Cartan) de ${\bf H}^\natural$ sont 
deux--\`a--deux conjugu\'es dans ${\bf H}$. On a les propri\'et\'es:
\begin{itemize}
\item Tout tore maximal ${\bf S}^\natural$ de ${\bf H}^\natural$ 
d\'etermine un quadruplet $({\bf S},{\bf S}^\natural,{\bf T},{\bf T}^\natural)$, 
qu'on appelle le {\it quadru\-plet de Cartan de ${\bf H}^\natural$ associ\'e \`a 
${\bf S}^\natural$}: on a 
${\bf S}=Z_{\bf H}({\bf S}^\natural)$, ${\bf T}=Z_{\bf H}({\bf S})$ et ${\bf T}^\natural={\bf T}\cdot {\bf S}^\natural$;
\item Tout sous--espace de Cartan ${\bf T}^\natural$ de ${\bf H}^\natural$ 
d\'etermine un triplet $({\bf S},{\bf T},{\bf T}^\natural)$, qu'on appelle le {\it triplet de Cartan de 
${\bf H}^\natural$ associ\'e \`a ${\bf T}^\natural$}: on a ${\bf S}=Z_{\bf H}({\bf T}^\natural)$ et ${\bf T}=Z_{\bf H}({\bf S})$.
\item Fix\'e un sous--espace de Cartan ${\bf T}^\natural$ de ${\bf H}^\natural$, 
l'application qui \`a $\delta\in {\bf T}^\natural$ associe la partie ${\bf S}\cdot \delta$ de ${\bf H}^\natural$ 
est une bijection du ${\bf S}\backslash {\bf T}$--espace tordu ${\bf S}\backslash {\bf T}^\natural$ 
sur l'ensemble des tores maximaux ${\bf S}^\natural$ de ${\bf H}^\natural$ tels que 
$Z_{\bf H}(Z_{\bf H}({\bf S}^\natural))\cdot {\bf S}^\natural={\bf T}^\natural$.
\end{itemize} 

\subsection{Orbites dans un ${\bf H}$--espace tordu}\label{orbites}
Soit ${\bf H}^\natural$ un 
${\bf H}$--espace tordu. Consid\'erons l'action de 
${\bf H}$ (\`a droite) sur ${\bf H}^\natural$ par conjugaison:
$$
{\bf H}^\natural \times {\bf H} \rightarrow {\bf H}^\natural,\,(\delta,h)\mapsto 
h^{-1}\cdot \delta \cdot h= 
h^{-1}{\rm Int}_{{\bf H}^\natural}(\delta)(h)\cdot \delta.
$$
Pour $\delta\in {\bf H}^\natural$, on note $\ES{O}_{\bf H}(\delta)$ la ${\bf H}$--orbite
$\{h^{-1}\cdot \delta \cdot h:h\in {\bf H}\}$; posant $\tau={\rm Int}_{{\bf H}^\natural}(\delta)$, on a
$$\ES{O}_{\bf H}(\delta)= {\bf H}(1-\tau)\cdot \delta.
$$ 
D'apr\`es \cite[ch.~I, 1.8]{Bor}, $\ES{O}_{\bf H}(\delta)$  est une 
vari\'et\'e alg\'ebrique lisse, et sa fermeture $\overline{\ES{O}_{\bf H}(\delta)}$ dans ${\bf H}^\natural$ est 
r\'eunion de $\ES{O}_{\bf H}(\delta)$ et de ${\bf H}$--orbites de dimension strictement inf\'erieure 
\`a $\dim(\ES{O}_{\bf H}(\delta))$. En particulier, la ${\bf H}$--orbite $\ES{O}_{\bf H}(\delta)$ est ouverte dans 
$\overline{\ES{O}_{\bf H}(\delta)}$ et localement ferm\'ee dans ${\bf H}^\natural$. Le groupe 
${\bf H}_\delta\;(=\{h\in {\bf H}: {\rm Int}_{{\bf H}^\natural}(\delta)(h)=h\})$ co\"{\i}ncide avec le stabilisateur 
de $\delta$ dans ${\bf H}$, et le morphisme de vari\'et\'es alg\'ebriques 
$\pi_\delta:{\bf H}\rightarrow \ES{O}_{\bf H}(\delta),\,h\mapsto h^{-1}\cdot \delta \cdot h$ se 
factorise \`a travers le quotient ${\bf H}_\delta\backslash {\bf H}$. On obtient donc 
un morphisme bijectif $\bar{\pi}_\delta:{\bf H}_\delta\backslash {\bf H}\rightarrow \ES{O}_{\bf H}(\delta)$, 
qui n'est en g\'en\'eral pas un isomorphisme. On a l'\'egalit\'e $(*)$ de \ref{groupes algbriques}
$$
\dim ({\bf H})= \dim ({\bf H}_\delta) + \dim (\ES{O}_{\bf H}(\delta)),
$$
et $\bar{\pi}_\delta$ est un 
isomorphisme si et seulement si le morphisme $\pi_\delta$ est s\'eparable, i.e. si et seulement si 
${\rm Lie}({\bf H}_\delta)=\ker ({\rm d}(\pi_\delta)_1)$. 

\begin{monlem}
Soit deux éléments $\delta,\, \delta'\in {\bf H}^\natural$ quasi--semisimples et unipotents. On a
$$
\ES{O}_{\bf H}(\delta')= \ES{O}_{\bf H}(\delta).
$$
\end{monlem}

\begin{proof}
Posons $\tau ={\rm Int}_{{\bf H}^\natural}(\delta)$ et $\tau'={\rm Int}_{{\bf H}^\natural}(\delta')$. 
D'après le lemme de \ref{automorphismes qss lf}, il existe un $x\in {\bf H}$ tel que
$$
\tau' = {\rm Int}_{\bf H}(x^{-1})\circ \tau \circ {\rm Int}_{\bf H}(x)\;(={\rm Int}_{{\bf H}^\natural}(x^{-1}\cdot \delta \cdot x)).
$$
Par conséquent $\delta'= z x^{-1}\cdot \delta \cdot x$ pour un élément $z\in Z(G)$, et quitte à remplacer $\delta'$ par $x\cdot \delta \cdot x^{-1}$, on peut supposer que $\delta'= z\cdot \delta$, i.e. que $\tau'=\tau$. Notons que si $p=1$, alors $\tau={\rm id}_{\bf H}$, et si $p>1$, 
alors $\tau^{p^k}={\rm id}_{\bf H}$ pour un entier $k\geq 1$. Soit ${\bf T}$ un tore maximal $\tau$--stable de ${\bf H}$. Puisque la restriction de $\tau$ à ${\bf T}$ est un automorphisme unipotent de ${\bf T}$, d'après la proposition de \ref{automorphismes ss et u}, l'élément $z$ s'écrit $z= u t^{-1}\tau(t)$ pour des éléments $u\in {\bf T}_\tau$ et $t\in {\bf T}$. Comme $\tau$ est un automorphisme localement fini de ${\bf H}$, on peut comme en \ref{automorphismes ss et u} identifier ${\bf H}$ à la composante neutre d'un groupe algébrique affine 
${\bf H}'=({\bf H}\rtimes \langle \tau \rangle )/{\bf C}$. On identifie ${\bf H}^\natural$ à l'image de ${\bf H}\rtimes \langle \tau \rangle$ dans ${\bf H}'$ via l'application $h\cdot \delta \mapsto (h,\tau)$. Alors $\delta$ s'identifie à l'image de $(1,\tau)$ dans ${\bf H}'$, 
et $t\cdot \delta'\cdot t'^{-1}$ s'identifie à l'image de $(u,\tau)$ dans ${\bf H}'$. Les éléments $\delta^{-1}$ et $t\cdot \delta' \cdot t^{-1}$ de ${\bf H}'$ sont unipotents, et comme ils commutent, leur produit $u\in {\bf T}$ est encore unipotent. Donc $u=1$, et $\delta'= t^{-1}\cdot \delta \cdot t$.
\end{proof}

\begin{mapropo}
Soit $\delta\in {\bf H}^\natural$. La ${\bf H}$--orbite 
$\ES{O}_{\bf H}(\delta)$ est ferm\'ee dans ${\bf H}^\natural$ si et seulement si 
$\delta$ est quasi--semisimple.
\end{mapropo}

\begin{proof} Posons $\tau={\rm Int}_{{\bf H}^\natural}(\delta)$. 
Comme $\ES{O}_{\bf H}(\delta)= {\bf H}(1-\tau)\cdot \delta$, 
si ${\bf H}$ est un tore, l'\'enonc\'e est vide: 
${\bf H}(1-\tau)$ est un sous--groupe 
ferm\'e de ${\bf H}$, et $\tau$ est quasi--semisimple. 
On suppose donc que ${\bf H}$ n'est pas un tore. 

Supposons $\delta$ quasi--semisimple. Soit 
$({\bf B},{\bf T})$ une paire de Borel $\tau$--stable de ${\bf H}$, et soit ${\bf U}=R_{\rm u}({\bf B})$. 
Pour $t\in {\bf T}$ et $u\in {\bf U}$, posant $\delta'=t^{-1}\cdot \delta \cdot t$, on a
\begin{align*}
(tu)^{-1}\cdot \delta \cdot tu &=  u^{-1}{\rm Int}_{{\bf H}^\natural}(\delta')(u)\cdot \delta'\\
&=  u^{-1}{\rm Int}_{{\bf H}^\natural}(\delta')(u)t^{-1}\tau(t)\cdot \delta
\end{align*}
o $u^{-1}{\rm Int}_{{\bf H}^\natural}(\delta')(u)\in {\bf U}$ et $t^{-1}\tau(t)\in {\bf T}(1-\tau)$. Par cons\'equent 
la ${\bf B}$--orbite
$$\ES{O}_{\bf B}(\delta)=\{b^{-1}\cdot \delta \cdot b:b\in {\bf B}\} $$
v\'erifie la double inclusion
$$
{\bf T}(1-\tau)\cdot\delta \subset \ES{O}_{\bf B}(\delta)\subset {\bf U}{\bf T}(1-\tau)\cdot \delta.
$$
Or ${\bf U}{\bf T}(1-\tau)={\bf T}(1-\tau){\bf U}$ est un sous--groupe ferm\'e $\tau$--stable 
de ${\bf B}$. Pour $b\in {\bf B}$, $u'\in {\bf U}$ et $t'\in {\bf T}(1-\tau)$, on a 
$$b^{-1}\cdot (u't'\cdot \delta)\cdot b= (b^{-1}ub)(b^{-1}t'\tau(b))\cdot \delta$$
o $b^{-1}ub\in {\bf U}$ et $b^{-1}t'\tau(b)\in {\bf U}{\bf T}(1-\tau)$. Par cons\'equent 
${\bf U}{\bf T}(1-\tau)\cdot \delta$ est une sous--vari\'et\'e ferm\'ee de ${\bf H}^\natural$, stable 
par ${\bf B}$--conjugaison. Posons $\Phi=\Phi({\bf T},{\bf H})$ et $\Delta= \Delta({\bf T},{\bf B})$. 
Soit 
${\bf T}_*$ le sous--groupe ferm\'e de ${\bf T}\cap {\bf H}_{\rm der}$ d\'efini par
$$
{\bf T}_*=\{t\in {\bf T}\cap {\bf H}_{\rm der}: \alpha(t)=\beta(t),\;\forall \alpha,\,\beta\in \Delta\}.
$$
Sa composante neutre ${\bf T}_*^\circ$ est un tore de dimension 1. De plus, on a l'inclusion 
${\bf T}_*^\circ \subset {\bf H}_\tau^\circ$. En effet, puisque ${\bf T}$ est 
$\tau$--stable, $\tau$ op\`ere sur $\Phi$: pour $\alpha\in \Phi$, on pose $\tau(\alpha)= \alpha\circ \tau^{-1}$. 
Comme ${\bf B}$ est $\tau$--stable, on a $\tau(\Delta)=\Delta$. On en d\'eduit 
que ${\bf T}_*$ est $\tau$--stable, et que $\alpha(\tau(t))= \alpha(t)$ ($\alpha\in \Delta$, $t\in {\bf T}$). 
Par suite ${\bf T}_*^\circ$ est $\tau$--stable, et comme le seul $\overline{F}$--automorphisme non 
trivial de ${\bf T}_*^\circ$ est le passage \`a l'inverse $t\mapsto t^{-1}$, on obtient que $\tau\vert_{T_*^\circ}={\rm id}$. D'o\`u l'inclusion cherch\'ee. 
Soit $\delta'\in {\bf U}{\bf T}(1-\tau)\cdot \delta$. \'Ecrivons 
$\delta'= u't^{-1}\tau(t)\cdot \delta$ avec $u'\in {\bf U}$ et $t\in {\bf T}$, et posons 
$v=t'ut'^{-1}\in {\bf U}$. Alors on a $t\cdot \delta'\cdot  t^{-1}=v\cdot \delta$. D'apr\`es la proposition de 
\ref{automorphismes rguliers; le cas i}, 
l'ensemble ${\bf T}_*^\circ\cap {\bf H}_{\rm reg}$ est non vide, et la 
fermeture $\overline{\ES{O}_{{\bf T}_*^\circ}(v)}$ de la ${\bf T}_*^\circ$--orbite 
$\ES{O}_{{\bf T}_*^\circ}(v)=\{tvt^{-1}:t\in {\bf T}_*^\circ\}$ dans ${\bf U}$ contient $1$. 
Par suite $\delta\in \overline{\ES{O}_{\bf T}(\delta')}$ et 
$\ES{O}_{\bf B}(\delta)\subset \overline{\ES{O}_{\bf B}(\delta')}$, o les fermetures sont dans 
${\bf U}{\bf T}(1-\tau)\cdot \delta$). En 
choisissant $\delta'$ tel que la ${\bf B}$--orbite $\ES{O}_{\bf B}(\delta')$ soit 
ferm\'ee dans ${\bf U}{\bf T}(1-\tau)\cdot \delta$ \cite[ch.~I, cor.~1.8]{Bor}, on 
obtient que la ${\bf B}$--orbite $\ES{O}_{\bf B}(\delta)$ est ferm\'ee dans ${\bf H}^\natural$. 
Et comme la vari\'et\'e quotient ${\bf H}/{\bf B}$ est compl\`ete, on en d\'eduit que la 
${\bf H}$--orbite $\ES{O}_{\bf H}(\delta)$ est ferm\'ee dans ${\bf H}^\natural$.

R\'eciproquement, supposons $\overline{\ES{O}_{\bf H}(\delta)}=\ES{O}_{\bf H}(\delta)$. 
Soit $({\bf B},{\bf T})$ une paire de Borel de ${\bf H}$ telle que $\tau({\bf B})={\bf B}$, et soit 
$u$ l'unique \'el\'ement de ${\bf U}=R_{\rm u}({\bf B})$ tel que ${\rm Int}_{\bf H}(u)\circ \tau
({\bf B},{\bf T})=({\bf B},{\bf T})$. Posons $\delta'=u\cdot \delta$ et $\tau'={\rm Int}_{\bf H}(\delta')$. 
On d\'efinit 
comme plus haut le sous--tore ${\bf T}_*^\circ$ de ${\bf T}$. D'apr\`es le 
paragraphe pr\'ec\'edent, on a l'inclusion ${\bf T}_*^\circ \subset {\bf H}_{\tau'}^\circ$, et 
la fermeture $\overline{\ES{O}_{{\bf T}_*^\circ}(u^{-1}\cdot \delta')}$ de la 
${\bf T}_*^\circ$--orbite $\ES{O}_{{\bf T}_*^\circ}(u^{-1}\cdot \delta')$ contient $\delta'$. 
Or $u^{-1}\cdot \delta'=\delta$ et la ${\bf H}$--orbite $\ES{O}_{\bf H}(\delta)$ est ferm\'ee 
dans ${\bf H}^\natural$. Par cons\'equent $\delta'\in \ES{O}_{\bf H}(\delta)$ et $\delta$ est quasi--semisimple.
\end{proof}

\begin{moncoro1}
Pour $\delta\in {\bf H}^\natural_{\rm reg}$, la ${\bf H}$--orbite 
$\ES{O}_{\bf H}(\delta)$ est ferm\'ee dans ${\bf H}^\natural$.
\end{moncoro1}

\begin{marema}
{\rm On peut aussi prouver le corollaire 1 directement, gr\^ace 
\`a la relation $(*)$ de \ref{groupes algbriques}. En effet, soit 
$\delta\in {\bf H}^\natural_{\rm reg}$ et soit $\delta'\in \overline{\ES{O}_{\bf H}(\delta)}$.  
Pour $h,\,x\in {\bf H}$, puisque
$$
{\rm Ad}_{{\bf H}^\natural}(x\cdot \delta\cdot x^{-1})= 
{\rm Ad}_{\bf H}(x)\circ {\rm Ad}_{{\bf H}^\natural}(\delta)\circ 
{\rm Ad}_{\bf H}(x)^{-1},$$
posant $P(h,\delta)=P(h,{\rm Int}_{\bf H}(\delta))$ comme en \ref{automorphismes rguliers; le cas i}, 
on a
$$
P(xhx^{-1}, x\cdot \delta\cdot x^{-1})= P(h, \delta).
$$
Par continuit\'e, on en d\'eduit que $P(1,\delta')=P(1,\delta)$. En particulier, on a
$D_{\bf H}(\delta')=D_{\bf H}(\delta)\neq 0$. 
Puisque $\delta'\in {\bf H}^\natural_{\rm reg}$, d'apr\`es le thorme de \ref{automorphismes rguliers; le cas g}, 
on a $\dim({\bf H}_{\delta'}^\circ)=r({\bf H}^\natural)=\dim({\bf H}_\delta^\circ)$. 
Gr\^ace \`a la relation $(*)$ de \ref{groupes algbriques}, on obtient
$$
\dim(\ES{O}_{\bf H}(\delta'))= \dim({\bf H})-r({\bf H}^\natural)= \dim(\ES{O}_{\bf H}(\delta)).
$$
Comme $\ES{O}_{\bf H}(\delta')\subset \overline{\ES{O}_{\bf H}(\delta)}$, cela entra\^{\i}ne que 
$\delta'$ appartient  $\ES{O}_{\bf H}(\delta)$.\hfill $\blacksquare$
}\end{marema}

\begin{moncoro2}
Soit un \'el\'ement 
$\delta\in {\bf H}^\natural$ unipotent. 
La fermeture 
$\overline{\ES{O}_{\bf H}(\delta)}$ est r\'eunion de ${\bf H}$--orbites unipotentes, et contient une unique ${\bf H}$--orbite fermée. 
Cette dernière est l'unique ${\bf H}$--orbite quasi--semisimple unipotente de ${\bf H}^\natural$.
\end{moncoro2}

\begin{proof}
Puisque $\delta$ est unipotent, tout 
\'el\'ement de $\overline{\ES{O}_{\bf H}(\delta)}$ l'est aussi. Par suite $\overline{\ES{O}_{\bf H}(\delta)}$ est r\'eunion de 
${\bf H}$--orbites unipotentes. Soit $\ES{O}$ une ${\bf H}$--orbite dans 
$\overline{\ES{O}_{\bf H}(\delta)}$ de dimension minimale. Elle est ferm\'ee 
dans ${\bf H}^\natural$, donc quasi--semisimple (proposition), et d'après le lemme, $\ES{O}$ est l'unique ${\bf H}$--orbite quasi--semisimple unipotente de ${\bf H}^\natural$.
\end{proof}

\section{Questions de rationnalit}

Continuons avec les hypothses et les notations du ch.~3: $F$ est un corps commutatif d'exposant caractristique $p\geq 1$, 
$\overline{F}$ est une cl\^oture algbrique de $F$, et ${\bf H}={\bf H}(\overline{F})$ est un groupe algbrique. On suppose de plus 
que ${\bf H}$ est connexe, rductif et dfini sur $F$, et l'on note $H={\bf H}(F)$ le groupe de ses points $F$--rationnels. \`A l'exception du ${\rm n}^\circ$ 
\ref{la topo p-adique}, 
les notions topologiques se rfrent toujours  la topologie de Zariski.

\subsection{G\'en\'eralit\'es (rappels)}\label{gnralits} 
Soit $F^{\rm sep}\!/F$ la sous--extension s\'eparable maximale 
de $\overline{F}/F$, et soit $\Sigma=\Sigma(F^{\rm sep}\!/F)$ son groupe de Galois. 

Soit ${\bf V}$ une vari\'et\'e alg\'ebrique d\'efinie sur $F$. Pour toute extension $F'$ de $F$, on note 
${\bf V}(F')$ l'ensemble de ses points $F'$--rationnels. D'aprs \cite[ch.~AG, 13.3]{Bor}, ${\bf V}(F^{\rm sep})$ est dense 
dans ${\bf V}$ (pour la topologie de Zariski). Une partie 
ferm\'ee de ${\bf V}$ est dite 
{\it $F$--ferme\footnote{Soit ${\bf H}'$ un groupe 
alg\'ebrique affine d\'efini sur $F$. Notons $\ES{H}'$ le $F$-sch\'ema en groupes affine 
lisse d'alg\'ebre affine $F[{\bf H}']$. 
Les parties $F$--fermes (resp. ferm\'ees et d\'efinies sur $F$) de ${\bf H}'$ correspondent 
bijectivement aux sous--$F$-sch\'emas ferm\'es {\it r\'eduits} (resp. {\it g\'eom\'etriquement r\'eduits}) de 
$\ES{H}'$: les unes et les autres correspondent bijectivement aux id\'eaux $I$ de $F[{\bf H}']$ 
tels que la $F$--alg\`ebre $F[{\bf H}']/I$ est r\'eduite (resp. tels que la $\overline{F}$--alg\`ebre 
$\overline{F}\otimes_F F[{\bf H}']/I$ est r\'eduite). \`A tout sous groupe ferm\'e 
de ${\bf H}'$ d\'efini sur $F$ correspond ainsi un sous--$F$--sch\'ema en groupes 
ferm\'e lisse de $\ES{H}'$, et r\'eciproquement. 

Soit $\tau\in {\rm Aut}_F({\bf H}')$. Le groupe 
${\bf H}'_\tau$ est $F$--ferm dans ${\bf H}'_\tau$, donc correspond \`a un 
sous--$F$-sch\'ema ferm\'e r\'eduit $\ES{H}'_\tau$ de $\ES{H}'$. On a 
${\bf H}'_\tau=\ES{H}'_\tau(\overline{F})$, et 
${\bf H}'_\tau$ est d\'efini sur $F$ si et seulement si $\ES{H}'_\tau$ est 
g\'eom\'etriquement r\'eduit, i.e. est un 
$F$-sch\'ema en groupes lisse.}(dans ${\bf V}$)} si elle est d\'efinie sur une sous--extension purement 
ins\'eparable de $\overline{F}/F$; ou, ce qui revient au m\^eme, si elle est d\'efinie 
sur $F^{p^{-\infty}}$. Ainsi, toute partie ferm\'ee de ${\bf V}$ est 
$F^{\rm sep}$--ferme, et toute partie $F$--ferme de ${\bf V}$ d\'efinie sur $F^{\rm sep}$ est 
d\'efinie sur $F$. En particulier si $p=1$, toute partie $F$--ferme de ${\bf V}$ est d\'efinie 
sur $F$. 

Soit $\pi:{\bf V}\rightarrow {\bf W}$ un 
morphisme de vari\'et\'es alg\'ebriques. Supposons que ${\bf V}$, ${\bf W}$ et 
$\pi$ sont d\'efinis sur $F$. Alors l'image ${\bf W}'=\pi({\bf V})$ est d\'efinie sur $F$. On 
a l'inclusion $\pi({\bf V}(F))\subset {\bf W}'(F)$ mais cette inclusion est en g\'en\'eral stricte --- 
m\^eme si le morphisme ${\bf V}\buildrel \pi\over{\longrightarrow} {\bf W}'$ est s\'eparable, et que $F=F^{\rm sep}$. 

\begin{marema1}
{\rm 
Si $({\bf W}'\!,\pi)$ est \og le\fg{} quotient de ${\bf V}$ par un groupe alg\'ebrique 
affine ${\bf H}'$ d\'efini sur $F$, alors par homog\'en\'eit\'e d'apr\`es \cite[ch.~AG, 13.2]{Bor}, 
on a l'\'egalit\'e
$$\pi({\bf V}(F^{\rm sep}))={\bf W}'(F^{\rm sep}).
\eqno \blacksquare
$$
}
\end{marema1}

Soit ${\bf H}'$ un groupe alg\'ebrique affine d\'efini sur $F$. 
Sa composante neutre ${\bf H}'^\circ$ est elle 
aussi d\'efinie sur $F$ \cite[ch.~I, 1.2]{Bor}, et ses composantes 
connexes sont toutes d\'efinies sur $F^{\rm sep}$ \cite[ch.~AG, 12.13]{Bor}. On a 
donc ${\bf H}'={\bf H}'(F^{\rm sep}){\bf H}'^\circ$. D'apr\`es \cite[ch.~V, 18.2]{Bor}, 
il existe un 
tore maximal de ${\bf H}'^\circ$ d\'efini sur $F$. Rappelons que 
${\bf H}'$ est {\it d\'eploy\'e sur $F$} s'il existe une suite  
$$
\{1\}={\bf H}'_n\subset {\bf H}'_{n-1}\subset \cdots \subset {\bf H}'_1={\bf H}'
$$ de sous--groupes d\'efinis 
sur $F$, tels que pour $i=1,\ldots ,n-1$, ${\bf H}'_{i+1}$ est distingu\'e dans ${\bf H}'_i$ et le groupe 
quotient ${\bf H}'_i/{\bf H}'_{i+1}$ est $F$--isomorphe au groupe multiplicatif $\Bbb{G}_{\rm m}$ ou 
au groupe additif $\Bbb{G}_{\rm a}$. Si ${\bf H}'$ est diagonalisable (e.g. un tore), alors ${\bf H}'$ se d\'eploie 
sur une sous--extension galoisienne finie de $F^{\rm sep}\!/F$ \cite[ch.~III, 8.11]{Bor}. 

Soit ${\bf V}$ une vari\'et\'e alg\'ebrique d\'efinie sur $F$, munie d'une action 
alg\'ebrique de ${\bf H}'$ \`a gauche ${\bf H}'\times {\bf V}\rightarrow {\bf V},\,(h,v)
\mapsto h \cdot v$, elle aussi d\'efinie sur $F$. Pour $v\in {\bf V}(F)$, la 
vari\'et\'e (lisse) ${\bf H}'\cdot v$ est d\'efinie sur $F$, et le groupe ${\bf H}'_v$ est $F$--ferm dans 
${\bf H}'$. Si de plus le morphisme $\pi_v:{\bf H}'\rightarrow {\bf V},\,h\mapsto h\cdot v$ est s\'eparable, 
alors ${\bf H}'_v$ est d\'efini sur $F$ \cite[ch.~II, 6.7]{Bor}. D'autre part, si le groupe 
${\bf H}'$ est connexe, r\'esoluble et $F$--d\'eploy\'e, et s'il op\`ere transitivement sur ${\bf V}$, alors ${\bf V}$ 
est affine et poss\`ede un point $F$--rationnel \cite[ch.~V, 15.11]{Bor}.

\begin{marema2}
{\rm Soit $\pi:{\bf V}\rightarrow {\bf W}$ un 
morphisme de vari\'et\'es alg\'ebriques. Supposons que ${\bf V}$, ${\bf W}$ et 
$\pi$ sont d\'efinis sur $F$. Supposons aussi que ${\bf V}$ et ${\bf W}'=\pi({\bf V})$ sont 
irr\'eductibles, et que ${\bf V}$ est muni d'une action alg\'ebrique (\`a gauche) de ${\bf H}'$, d\'efinie sur $F$, 
telle que les fibres de $\pi$ sont les orbites sous ${\bf H}'$; o\`u ${\bf H}'$ est toujours un groupe alg\'ebrique affine d\'efini sur $F$. 
Si le morphisme ${\bf V}\buildrel \pi\over{\longrightarrow}{\bf W}'$ est s\'eparable 
--- en particulier si $({\bf W}'\!,\pi)$ est \og le\fg{} quotient de ${\bf V}$ par ${\bf H}'$ --- 
et si ${\bf H}'$ est connexe, 
r\'esoluble et $F$--d\'eploy\'e, alors d'aprs \cite[ch.~V, 15.12]{Bor}, on a 
l'\'egalit\'e
$$\pi({\bf V}(F))={\bf W}(F).\eqno{\blacksquare}
$$}
\end{marema2}

\subsection{Gnralits; suite}\label{gnralits; suite}
Les rsultats ci--dessous sont valables pour n'importe quel groupe algbrique ${\bf H}$ rductif connexe et dfini sur $F$:
\begin{itemize}
\item ${\bf H}_{\rm der}$ est d\'efini sur $F$ \cite[ch.~I, cor.~2.3]{Bor};
\item $Z({\bf H})$ est d\'efini 
sur $F$ \cite[ch.~V, 18.2]{Bor}, et donc aussi $R({\bf H})=Z({\bf H})^\circ$ \cite[ch.~I, 1.2]{Bor};
\item ${\bf H}$ est d\'eploy\'e sur $F$ si et seulement s'il existe un tore maximal de 
${\bf H}$ d\'efini et d\'eploy\'e sur $F$ \cite[ch.~V, 18.7]{Bor};
\item si ${\bf T}$ est un tore maximal de ${\bf H}$ d\'efini et d\'eploy\'e sur une sous--extension 
$F'\!/F$ de $\overline{F}/F$, alors tout sous--groupe de Borel ${\bf B}$ de ${\bf H}$ 
contenant ${\bf T}$ est d\'efini et d\'eploy\'e sur $F'$ (cf. la d\'emonstration de loc. cit.);
\item si ${\bf S}$ est un tore de ${\bf H}$ d\'efini sur $F$, alors $Z_{\bf H}({\bf S})$ et $N_{\bf H}({\bf S})$ sont d\'efinis sur $F$ 
\cite[ch.~III, cor.~9.2 et ch.~V, 20.3]{Bor}; si de plus 
${\bf S}$ est d\'eploy\'e sur $F$, alors 
$Z_{\bf H}({\bf S})$ est une composante de Levi d'un sous--groupe parabolique de ${\bf H}$ d\'efini 
sur $F$ \cite[ch.~V, 20.4]{Bor};
\item si ${\bf P}$ est un sous--groupe parabolique de ${\bf H}$ d\'efini sur $F$, alors 
$R({\bf P})$ et $R_{\rm u}({\bf P})$ sont d\'efinis sur $F$, et les composantes de Levi dfinies sur $F$ de ${\bf P}$ 
sont les centralisateurs dans ${\bf H}$ des tores maximaux dfinis sur $F$ de $R({\bf P})$ \cite[ch.~20, 20.5]{Bor};
\item si ${\bf M}$ est une composante de Levi dfinie sur $F$ d'un sous--groupe parabolique ${\bf P}$ de ${\bf H}$ 
dfini sur $F$, alors le sous--groupe parabolique de ${\bf H}$ oppos\'e \`a ${\bf P}$ par rapport \`a ${\bf M}$, 
est d\'efini sur $F$ \cite[ch.~V, 20.5]{Bor}, et notant ${\bf S}$ le sous--tore $F$--d\'eploy\'e maximal 
de $R({\bf M})={\bf M}\cap R({\bf P})$ \cite[ch.~IV, 11.23]{Bor}, on a l'galit $Z_{\bf H}({\bf S})={\bf M}$ \cite[ch.~V, 20.6]{Bor}.
\item d'aprs ce qui prcde, si ${\bf P}$ est un sous--groupe parabolique de ${\bf H}$ d\'efini sur $F$, alors 
les composantes de Levi dfinies sur $F$ de ${\bf P}$ 
sont les centralisateurs dans ${\bf H}$ des tores $F$--dploys maximaux de $R({\bf P})$.
\end{itemize}

\begin{convention}
{\rm On applique au groupe $H={\bf H}(F)$ 
le langage des groupes alg\'ebriques de la mani\`ere 
habituelle suivante. On appelle:
\begin{itemize}
\item {\it tore de $H$} le groupe des points $F$--rationnels d'un tore de ${\bf H}$ d\'efini sur $F$;
\item {\it tore d\'eploy\'e de $H$} le groupe des points $F$--rationnels d'un tore 
de ${\bf H}$ d\'efini et d\'eploy\'e sur $F$;
\item {\it sous--groupe parabolique de $H$} 
le groupe des points $F$--rationnels $P={\bf P}(F)$ 
d'un sous--groupe parabolique ${\bf P}$ de ${\bf H}$ d\'efini sur $F$, et {\it composante de Levi} 
de $P$ le groupe des points $F$--rationnels d'une composante de Levi de ${\bf P}$ d\'efinie sur 
$F$;
\item {\it radical} (resp. {\it radical unipotent}) 
d'un sous--groupe parabolique ${\bf P}(F)$ de $H$ le groupe $R({\bf P})(F)$ 
(resp. le groupe $R_{\rm u}({\bf P})(F)$) --- on le note $R(P)$ (resp. $R_{\rm u}(P)$.
\end{itemize}}
\end{convention}

Soit $P={\bf P}(F)$ un sous--groupe parabolique de $H$. D'aprs les rsultats rappels plus haut, 
toute composante de Levi $M$ de $P$ 
est le centralisateur dans $H$ d'un tore dploy maximal de $R(P)$, et l'on a 
la dcomposition en produit semi--direct, appele {\it dcomposition de Levi de $P$}:
$$
P=M\ltimes R_{\rm u}(P).
$$
D'aprs \cite[ch.~V, 20.5]{Bor}, si $M$ et $M'$ sont deux 
composantes de Levi de $P$, alors il existe une unique $u\in R_{\rm u}({\bf P})(F)$ tel que 
$M'=uMu^{-1}$. De plus (loc.~cit.), le morphisme quotient ${\bf H}\rightarrow {\bf H}/{\bf P}$ induit une 
application {\it surjective} $H\rightarrow ({\bf H}/{\bf P})(F)$ , 
i.e. on a
$$({\bf H}/{\bf P})(F) = H/P.$$

\subsection{Points rationnels d'un ${\bf H}$--espace tordu d\'efini sur $F$}\label{points rationnels d'un espace tordu} 
Soit ${\bf H}^\natural$ un ${\bf H}$--espace tordu {\it d\'efini sur $F$}, c'est--\`a--dire tel que:
\begin{itemize}
\item ${\bf H}^\natural$ est une vari\'et\'e d\'efinie sur $F$;
\item les actions \`a gauche et \`a droite de ${\bf H}$ sur ${\bf H}^\natural$ sont d\'efinies sur $F$.
\end{itemize}
On note $H^\natural={\bf H}^\natural(F)$ 
l'ensemble des points $F$--rationnels de ${\bf H}^\natural$. On suppose que $H^\natural$ est {\it non vide}. 
Alors pour $\delta\in H^\natural$, le $\overline{F}$--automorphisme ${\rm Int}_{{\bf H}^\natural}(\delta)$ 
de ${\bf H}$
est d\'efini sur $F$, i.e. on a ${\rm Int}_{{\bf H}^\natural}(\delta)\in {\rm Aut}_F({\bf H})$. Munissons 
$H$ et $H^\natural$ de la topologie de Zariski h\'erit\'ee de 
${\bf H}$ et de ${\bf H}^\natural$. Alors $H^\natural$ est muni d'une structure 
de $H$--espace topologique tordu au sens de \ref{espaces topo tordus}: pour $\delta\in H^\natural$, 
on a $H^\natural = H\cdot \delta\subset {\bf H}^\natural$, et l'automorphisme 
${\rm Int}_{H^\natural}(\delta)$ de $H$ est donn\'e par la restriction de 
${\rm Int}_{{\bf H}^\natural}(\delta)$ \`a $H$. 
Notons que si ${\bf H}^\natural={\bf H}\tau$ pour un $\tau\in {\rm Aut}_F({\bf H})$, alors 
$H^\natural$ co\"{\i}ncide avec le $H$--espace topologique tordu 
$H\tau$ d\'efini en \ref{espaces topo tordus}. 

Puisque ${\bf H}^\natural$ est d\'efini sur $F$, la fonction 
${\bf H}^\natural \rightarrow \overline{F},\,\delta\mapsto D_{{\bf H}^\natural}
(\delta)$ appartient \`a l'alg\`ebre affine $F[{\bf H}^\natural]$. Pour 
$\delta\in H^\natural$, on pose
$$
D_{H^\natural}(\delta)=D_{{\bf H}^\natural}(\delta)\in F.
$$
L'ensemble
$${\bf H}^\natural \smallsetminus {\bf H}^\natural_{\rm reg}=\{\delta\in {\bf H}^\natural:
D_{{\bf H}^\natural}(\delta)=0\}$$ est 
$F$--ferm dans ${\bf H}^\natural$, par suite 
la vari\'et\'e ${\bf H}_{\rm reg}^\natural$ est d\'efinie sur $F$. 
On note
$$H^\natural_{\rm reg}= {\bf H}_{\rm reg}^\natural(F)\;(=H^\natural\cap {\bf H}^\natural_{\rm reg})
$$ l'ensemble de ses 
points $F$--rationnels. On a donc
$$
H^\natural_{\rm reg}=\{\delta\in H^\natural: D_{H^\natural}(\delta)\neq 0\}.
$$

\begin{marema}
{\rm Supposons le corps $F$ infini. 
Puisque ${\bf H}$ est rductif connexe, $H$ est dense dans ${\bf H}$ \cite[ch.~V, 18.3]{Bor}, et comme 
$H^\natural$ est suppos non vide, $H^\natural$ est dense dans 
${\bf H}^\natural$. Comme ${\bf H}^\natural_{\rm reg}$ est ouvert et 
dense dans ${\bf H}^\natural$, on en d\'eduit que $H^\natural_{\rm reg}$ 
est non vide. Par cons\'equent, 
$H^\natural_{\rm reg}$ est ouvert et dense dans $H^\natural$.\hfill $\blacksquare$}
\end{marema}

\subsection{La d\'ecomposition ${\rm Aut}_{F'}({\bf H})={\rm Int}_{F'}({\bf H})\rtimes \frak{A}_\circ$}\label{la dcomposition} 
Fixons un tore maximal ${\bf T}_\circ$ de ${\bf H}$ d\'efini sur $F$, et un sous--groupe de 
Borel ${\bf B}_\circ$ de ${\bf H}$ contenant ${\bf T}_\circ$. Notons 
${\bf B}_\circ^-$ le sous--groupe de Borel ${\bf B}_\circ^-$ de ${\bf H}$ 
oppos\'e \`a ${\bf B}_\circ$ par rapport \`a ${\bf T}_\circ$. Choisissons une sous--extension 
galoisienne finie $F_\circ/F$ de $F^{\rm sep}\!/F$ d\'eployant ${\bf T}_\circ$. 
Les groupes ${\bf B}_\circ$ et ${\bf B}_\circ^-$ sont donc tous les deux d\'efinis et d\'eploy\'es 
sur $F_\circ$. 
Comme en \ref{groupes rductifs}, posons $\Phi_\circ = \Phi({\bf T}_\circ,{\bf H})$, $\Phi_\circ^+=
\Phi({\bf T}_\circ,{\bf B}_\circ)$ et $\Delta_\circ = \Delta({\bf T}_\circ,{\bf B}_\circ)$. 
Posons aussi 
$\Delta_\circ^-=\{-\alpha:\alpha\in \Delta_\circ\}$. 
Pour $\alpha\in \Phi_\circ$, le groupe ${\bf U}_\alpha$ est d\'efini et d\'eploy\'e sur $F_\circ$. Pour chaque 
$\alpha\in \Delta_\circ$, choisissons un \'el\'ement $u_\alpha\in 
{\bf U}_\alpha(F_\circ)\smallsetminus\{1\}$, et 
notons $e_\alpha$ l'\'epinglage de 
${\bf U}_\alpha$ d\'efini par $u_\alpha$ (cf. \ref{groupes rductifs}); puisque $u_\alpha$ 
est $F_\circ$--rationnel, le morphisme $e_\alpha$ est d\'efini sur $F_\circ$. 
Rappelons qu'un {\it syst\`eme de Chevalley (relativement \`a ${\bf T}_\circ$)} est la donn\'ee 
d'une famille $\{f_\alpha\}_{\alpha\in \Phi_\circ}$ d'\'epinglages des ${\bf U}_{\alpha}$ 
v\'erifiant:
\begin{itemize}
\item pour $\alpha\in \Phi_\circ$, les \'epinglages $f_\alpha$ et $f_{-\alpha}$ 
sont associ\'es, i.e. l'\'el\'ement $f_{\alpha}(1)f_{-\alpha}(1)f_\alpha(1)$ 
appartient \`a $N_{\bf H}({\bf T}_\circ)$;
\item pour $\alpha,\,\beta\in \Phi_\circ$ , il existe 
$\epsilon= \pm 1$ tel que
$$f_{r_\alpha(\beta)}(x)= m_\alpha f_\beta(\epsilon x)m_\alpha^{-1}\quad (x\in \overline{F}),
$$ 
o\`u $r_\alpha \in N_{\bf H}({\bf T}_\circ)/{\bf T}_\circ$ d\'esigne 
la r\'eflection associ\'ee \`a la racine 
$\alpha$, et où on a noté $m_\alpha$ l'élément $f_{\alpha}(1)f_{-\alpha}(1)f_\alpha(1)\in 
N_{\bf H}({\bf T}_\circ)$.
\end{itemize}
On sait que la famille $\{e_\alpha\}_{\alpha\in \Delta_\circ}$ se prolonge 
en un {\it $F_\circ$--syst\`eme de Chevalley} $\{e_\alpha\}_{\alpha\in \Phi_\circ}$, i.e. 
un syst\`eme de 
Chevalley tel que les \'epinglages $e_\alpha:\Bbb{G}_{\rm a}\rightarrow {\bf U}_\alpha$ 
($\alpha\in \Phi_\circ$) sont d\'efinis sur $F_\circ$. De plus ce prolongement est unique 
\og au signe pr\`es\fg{}, c'est--\`a--dire au remplacement \'eventuel, pour certaines 
racines $\beta \in \Phi_\circ\smallsetminus (\Delta_\circ\cup \Delta_\circ^-)$, de 
 $e_\beta$ par $x\mapsto \bar{e}_\beta(x)=
e_\beta(-x)$. 
Alors pour $\alpha\in \Phi_\circ$, 
les \'el\'ements $u_\alpha=e_\alpha(1)\in {\bf U}_\alpha$ et 
$m_\alpha = u_\alpha u_{-\alpha}u_\alpha\in N_{\bf H}({\bf T}_\circ)$ 
sont $F_\circ$--rationnels.

Posons
$$\frak{A}_\circ = {\rm Aut}_{\overline{F}}({\bf H},{\bf B}_\circ,{\bf T}_\circ,
\{u_\alpha\}_{\alpha\in \Delta_\circ}),
$$ et soit $\tau_\circ\in \frak{A}_\circ$. 
Puisque $\tau_\circ({\bf B}_\circ,{\bf T}_\circ)=({\bf B}_\circ,{\bf T}_\circ)$, $\tau_\circ$ 
op\`ere sur l'ensemble $\Phi_\circ$ en stabilisant $\Delta_\circ$: pour 
$\alpha\in \Phi_\circ$, la racine $\tau_\circ(\alpha)\in \Phi_\circ$ est donn\'ee par 
${\bf U}_{\tau_\circ(\alpha)}=\tau_\circ({\bf U}_\alpha)$, et comme $\tau_\circ\circ e_\alpha= e_{\tau_\circ(\alpha)}$ 
($\alpha\in \Delta_\circ$), l'unicit\'e du $F_\circ$--syst\`eme de Chevalley 
prolongeant $\{e_\alpha\}_{\alpha\in \Delta_\circ}$ implique 
que pour $\beta\in \Phi_\circ \smallsetminus (\Delta_\circ \cup \Delta_\circ^-)$, on a 
$\tau_\circ \circ e_\beta=e_{\tau_\circ(\beta)}$ ou bien $\tau_\circ \circ e_\beta=
\bar{e}_{\tau_\circ(\beta)}$. Par cons\'equent pour $\alpha\in \Phi_\circ$, l'isomorphisme de 
groupes alg\'ebriques $\tau_\circ :{\bf U}_\alpha \rightarrow {\bf U}_{\tau_\circ(\alpha)}$ 
est d\'efini sur $F_\circ$. 
D'autre part puisque le tore ${\bf T}_\circ$ est dfini et d\'eploy\'e sur $F_\circ$, 
la restriction de $\tau_\circ $ \`a ${\bf T}_\circ$ 
est d\'efinie sur $F_\circ$. Ordonnons (de mani\`ere arbitraire) 
l'ensemble $\Phi_\circ^+$, i.e. posons $\Phi_\circ^+ =\{\alpha_1,\ldots ,\alpha_n\}$. 
Posons ${\bf U}_\circ =R_{\rm u}({\bf B}_\circ)$ et ${\bf U}_\circ^- =R_{\rm u}({\bf B}_\circ^-)$. 
D'apr\`es \cite[ch.~IV, 14.4]{Bor}, l'application produit 
$\prod_{i=1}^n {\bf U}_{\alpha_i}\rightarrow {\bf U}_\circ$ 
est un $F_\circ$--isomorphisme de vari\'et\'es alg\'ebriques; de la m\^eme manire, 
l'application produit $\prod_{i=1}^n {\bf U}_{-\alpha_i}\rightarrow {\bf U}_\circ^-$ 
est un $F_\circ$--isomorphisme de vari\'et\'es alg\'ebriques. 
D'apr\`es \cite[ch.~IV, 14.14]{Bor}, l'application produit 
${\bf U}_\circ^-\times {\bf T}_\circ 
\times {\bf U}_\circ\rightarrow {\bf H}$ 
est un $F_\circ$--isomorphisme de vari\'etes alg\'ebriques sur un ouvert 
${\bf V}_\circ$ de ${\bf H}$ d\'efini sur $F_\circ$. D'apr\`es ce qui 
pr\'ec\`ede, on a $\tau_\circ({\bf V}_\circ)={\bf V}_\circ$ et la 
restriction de $\tau_\circ$ \`a ${\bf V}_\circ$ est d\'efinie sur $F_\circ$. 
Par cons\'equent $\tau_\circ$ est d\'efini sur $F_\circ$. On a donc montr\'e 
l'inclusion
$$
\frak{A}_\circ \subset {\rm Aut}_{F_\circ}({\bf H}).\leqno{(*)}
$$
Cette inclusion jointe  la relation $(**)$ de \ref{groupes rductifs}, 
entra\^{\i}nent que pour toute sous--extension 
$F'\!/F_\circ$ de $\overline{F}/F_\circ$, on a la d\'ecomposition en produit semidirect
$$
{\rm Aut}_{F'}({\bf H})= {\rm Int}_{F'}({\bf H})\rtimes \frak{A}_\circ.\leqno{(**)}
$$

\begin{mesrems}
{\rm Soit $\pi:{\bf H}_{\rm SC}\rightarrow {\bf H}_{\rm der}$ 
le rev\^etement universel de ${\bf H}_{\rm der}$ (cf. \ref{rev\^etement universel}).
\begin{enumerate}
\item[(1)]D'apr\`es \cite[2.6.1]{T}, le groupe ${\bf H}_{\rm SC}$ est d\'efini sur 
$F$ et d\'eploy\'e sur $F_\circ$, et le morphisme $\pi$ est d\'efini sur $F$. Les 
paires de Borel $({\bf B}'_\circ,{\bf T}'_\circ)$ de ${\bf H}_{\rm der}$ 
et $({\bf B}_{\circ,{\rm sc}},{\bf T}_{\circ,{\rm sc}})$ de ${\bf H}_{\rm SC}$ 
d\'efinies comme en \ref{groupes rductifs} en rempla\c{c}ant $({\bf B},{\bf T})$ par $({\bf B}_\circ,{\bf T}_\circ)$, 
sont d\'efinies et d\'eploy\'ees sur $F_\circ$. 
Pour $\alpha\in \Delta_\circ$, l'\'el\'ement $\tilde{u}_\alpha=\pi^{-1}(u_\alpha)$ 
est $F_\circ$--rationnel. \`A partir de $\frak{A}_\circ$, on d\'efinit comme 
en \ref{rev\^etement universel} les sous--ensemble 
 $\widetilde{\frak{A}}_{\circ}\subset {\rm Aut}_{\overline{F}}({\bf H}_{\rm SC})$ et 
 $\frak{A}'_\circ\subset {\rm Aut}_{\overline{F}}({\bf H}_{\rm der})$. D'apr\`es $(*)$, 
 on a les 
inclusions
$$\widetilde{\frak{A}}_\circ \subset {\rm Aut}_{F_\circ}({\bf H}_{\rm SC}),\quad 
\frak{A}'_\circ \subset {\rm Aut}_{F_\circ}({\bf H}_{\rm der}).$$ 
D'autre part le morphisme bijectif
$$\overline{\pi}: {\bf H}_{\rm SC}/Z({\bf H}_{\rm SC})\rightarrow {\bf H}_{\rm der}/Z({\bf H}_{\rm der})$$ 
est d\'efini sur $F$, mais n'est en g\'en\'eral pas un isomorphisme (cf. l'exemple de \ref{rev\^etement universel}). On en 
d\'eduit que si $\tau\in {\rm Aut}_F({\bf H}_{\rm der})$, le rel\`evement 
$\tilde{\tau}$ de $\tau$ \`a ${\bf H}_{\rm SC}$ n'est pas n\'ecessairement d\'efini sur 
$F$ (m\^eme si $F=F_\circ$).
\item[(2)]D'apr\`es \cite[3.1.2]{T}, ${\bf H}_{\rm SC}$ se d\'ecompose en un produit direct 
$${\bf H}_{\rm SC}={\bf H}_{{\rm SC},1}\times \cdots \times {\bf H}_{{\rm SC},n},$$ o\`u chaque 
${\bf H}_{{\rm SC},i}$ est un groupe semisimple simplement connexe d\'efini sur $F$ et presque 
$F$--simple; de plus cette d\'ecomposition est unique \`a permutation des ${\bf H}_{{\rm SC},i}$ pr\`es. 
L'unicit\'e de la d\'ecomposition implique que si $\tau\in {\rm Aut}_F({\bf H}_{\rm SC})$, alors 
$\tau$ permute les facteurs ${\bf H}_{{\rm SC},i}$. Supposons 
${\bf H}$ presque $F$--simple, i.e. supposons ${\bf H}_{\rm SC}={\bf H}_{{\rm SC},1}$. 
Alors d'aprs loc. cit., il existe une sous--extension finie $L/F$ de $F^{\rm sep}\!/F$ telle 
que ${\bf H}_{\rm SC}$ est $F$--isomorphe \`a ${\rm Res}_{L/F}({\bf H}^*)$ pour 
un groupe semisimple simplement connexe ${\bf H}^*$ d\'efini sur $L$ et (absolument) 
presque simple. Ici, ${\rm Res}_{L/F}$ d\'esigne le foncteur {\it restriction des scalaires} 
de la cat\'egorie des groupes alg\'ebriques affines d\'efinis sur $L$ dans celle 
des groupes alg\'ebriques affines d\'efinis sur $F$. Le groupe 
${\bf H}_{\rm SC}$ se d\'ecompose donc, {\it sur $L$}, en un produit direct 
$${\bf H}_{\rm SC}= {\bf H}_1^*\times\cdots \times {\bf H}_m^*,$$ o\`u chaque groupe 
${\bf H}_i^*$ est $L$--isomorphe \`a ${\bf H}^*$, et $m=[L:F]$; comme plus haut, cette 
d\'ecom\-position est unique \`a permutation des ${\bf H}_i^*$ pr\`es, et si $\tau\in {\rm Aut}_F({\bf H}_{\rm SC})$, alors $\tau$ 
permute les facteurs ${\bf H}_i^*$.\hfill $\blacksquare$
\end{enumerate}}
\end{mesrems}

\subsection{Automorphismes stabilisant un sous--groupe de Borel d\'efini sur $F^{\rm sep}$}\label{automorphismes stabilisant B} 
On s'int\'eresse 
dans ce ${\rm n}^\circ$ \`a la version $F^{\rm sep}$--rationnelle du th\'eor\`eme 1 de \ref{automorphismes qss}: 
pour un $F$--automorphisme $\tau$ de ${\bf H}$, on aimerait savoir s'il existe un sous--groupe de 
Borel $\tau$--stable de ${\bf H}$ qui soit d\'efini sur une sous--extension de 
$F^{\rm sep}\!/F$ (i.e. sur $F^{\rm sep}$). Si $p>1$, la r\'eponse est n\'egative en g\'en\'eral:

\begin{exemple}{\rm 
Soit $F$ le corps de s\'eries formelles 
$\Bbb{F}_2((\varpi))$, ${\bf H}$ le groupe 
$\Bbb{GL}_2$, et $x\in {\bf H}$ la matrice $\left(\begin{array}{cc} 0& 1\\ \varpi & 0\end{array}\right)$. Prenons pour $\tau$ le 
$F$--automorphisme int\'erieur 
${\rm Int}_{\bf H}(x)$. Un sous--groupe de Borel de ${\bf H}$ 
est $\tau$--stable si et seulement s'il contient $x$, or aucun sous--groupe de Borel de 
${\bf H}$ d\'efini sur $F^{\rm sep}$ ne contient $x$.\hfill $\blacksquare$}
\end{exemple}

Pour traiter la question 
qui nous int\'eresse, on peut bien s\^ur supposer $F=F^{\rm sep}$.

\begin{monlem}
On suppose $F=F^{\rm sep}$. Soit un \'el\'ement $h\in H$ dont 
la d\'ecomposition de Jordan $h=h_{\rm s}h_{\rm u}$ est $F$--rationnelle, i.e. 
tel que $h_{\rm s}$ et $h_{\rm u}$ appartiennent \`a $H$. Alors il existe un 
sous--groupe de Borel ${\bf B}$ de ${\bf H}$ d\'efini sur $F$, tel 
que $h\in {\bf B}(F)$.
\end{monlem}

\begin{proof}
Supposons tout d'abord $h$ unipotent, i.e. $h=h_{\rm u}\in H$. 
Via le choix d'un $F$--plongement 
$\iota:{\bf H}\rightarrow \Bbb{GL}_n$, identifions ${\bf H}$ \`a un sous--groupe ferm\'e 
de $\Bbb{GL}_n$ d\'efini sur $F$. D'apr\`es \cite[ch.~I, theo. 4.8]{Bor}, il existe un \'el\'ement 
$x\in \Bbb{GL}_n(F)$ tel que $xhx^{-1}$ appartient au sous--groupe 
$\Bbb{U}_n$ de $\Bbb{GL}_n$ form\'e des matrices strictement triangulaires sup\'erieures. 
Notons $\Bbb{B}_n$ le sous--groupe ferm\'e de $\Bbb{GL}_n$ 
form\'e des matrices triangulaires sup\'erieures, et posons $\Bbb{B}_n^x=x^{-1}\Bbb{B}_n x$. 
Alors $\Bbb{B}_n^x$ est un sous--groupe de Borel de $\Bbb{GL}_n$ d\'efini sur $F$, 
${\bf H}\cap \Bbb{B}_n^x$ est un sous--groupe ferm\'e de ${\bf H}$ d\'efini sur $F$, 
et $h$ appartient au groupe $({\bf H}\cap \Bbb{B}_n^x)(F)$ des points $F$--rationnels de ${\bf H}\cap \Bbb{B}_n^x$. 
D'apr\`es \cite[ch.~IV, 11.14]{Bor}, la 
composante neutre ${\bf B}=({\bf H}\cap \Bbb{B}_n^x)^\circ$ de ${\bf H}\cap \Bbb{B}_n^x$ est un 
sous--groupe de Borel de ${\bf H}$, et comme $N_{\bf H}({\bf B})={\bf B}$ \cite[ch.~IV, 
11.16]{Bor}, on a ${\bf B}={\bf H}\cap \Bbb{B}_n^x$.

Passons au cas g\'en\'eral. Soit ${\bf M}=Z_{\bf H}(h_{\rm s})^\circ$; c'est un groupe r\'eductif connexe 
d\'efini sur $F$, et $h_{\rm s}$ et $h_{\rm u}$ appartiennent  ${\bf M}$. D'apr\`es le paragraphe 
pr\'ec\'edent, il existe un sous--groupe de Borel ${\bf B}_{\bf M}$ de ${\bf M}$ 
d\'efini sur $F$, tel que $h_{\rm u}\in {\bf B}_{\bf M}$. Soit ${\bf T}$ un tore maximal de 
${\bf B}_{\bf M}$ d\'efini sur $F$. Comme 
$h_{\rm s}\in Z({\bf M})\subset {\bf T}$, $h$ appartient \`a ${\bf B}_{\bf M}(F)$. 
D'apr\`es \cite[ch.~IV, 11.14]{Bor}, il existe un sous--groupe de 
Borel ${\bf B}$ de 
${\bf H}$ tel que ${\bf B}_{\bf M}= ({\bf B}\cap {\bf M})^\circ$. 
Comme ${\bf T}$ est un tore maximal de ${\bf H}$ d\'efini 
(donc d\'eploy\'e) sur $F$ --- rappelons que $F=F^{\rm sep}$ ---, 
le groupe ${\bf B}$ est lui aussi d\'efini et d\'eploy\'e sur $F$. D'o\`u le lemme. 
\end{proof}

\begin{mapropo}
On suppose $F=F^{\rm sep}$. Soit 
$\tau\in {\rm Aut}_F({\bf H})$ tel que la d\'ecomposition de Jordan 
$\tau_{\rm der}=(\tau_{\rm der})_{\rm s}\circ (\tau_{\rm der})_{\rm u}$ de 
$\tau_{\rm der}\in {\rm Aut}_F({\bf H}_{\rm der})$ est d\'efinie sur $F$, i.e. tel que 
$(\tau_{\rm der})_{\rm s}$ et $(\tau_{\rm der})_{\rm u}$ appartiennent  ${\rm Aut}_F({\bf H}_{\rm der})$. Alors 
il existe un sous--groupe de Borel $\tau$--stable de ${\bf H}$ d\'efini sur $F$.
\end{mapropo}

\begin{proof}
D'apr\`es \cite[ch.~IV, 11.14]{Bor}, 
l'application ${\bf B}\mapsto {\bf B}\cap {\bf H}_{\rm der}$ est une bijection 
de l'ensemble des sous--groupes de Borel de ${\bf H}$ sur l'ensemble des sous--groupes 
de Borel de ${\bf H}_{\rm der}$, 
et ${\bf B}$ est $\tau$--stable (resp. d\'efini sur $F$) si et seulement si ${\bf B}\cap {\bf H}_{\rm der}$ 
est $\tau$--stable (resp. d\'efini sur $F$). On peut donc supposer ${\bf H}={\bf H}_{\rm der}$. 
Alors $\tau=\tau_{\rm der}$ induit par passage au quotient un $F$--automorphisme 
$\overline{\tau}$ de $\overline{\bf H}={\bf H}/Z({\bf H})$ dont la d\'ecomposition de Jordan 
$\overline{\tau}= \overline{\tau}_{\rm s}\circ \overline{\tau}_{\rm u}$ est d\'efinie 
sur $F$. D'apr\`es loc. cit., notant $\phi$ le morphisme quotient 
${\bf H} \rightarrow \overline{\bf H}$, l'application ${\bf B}\mapsto \phi({\bf B})$ est 
une bijection 
de l'ensemble des sous--groupes de Borel de ${\bf H}$ sur l'ensemble des sous--groupes 
de Borel de $\overline{\bf H}$; et ${\bf B}$ est $\tau$--stable (resp. d\'efini sur $F$) 
si et seulement si $\phi({\bf B})$ est $\overline{\tau}$--stable (resp. d\'efini sur $F$). 
On peut donc suppposer $Z({\bf H})=\{1\}$.

Choisissons comme en \ref{automorphismes ss et u} un morphisme de groupe alg\'ebrique 
$$\iota:{\bf H}\rightarrow \Bbb{GL}_n$$ qui soit un isomorphisme sur un sous--groupe ferm\'e de $\Bbb{GL}_n$, et un \'el\'ement 
$g\in \Bbb{GL}_n$, tels que pour tout $h\in {\bf H}$ on ait $\iota \circ \tau ={\rm Int}_{\Bbb{GL}_n}(g)\circ \iota$. 
Puisque ${\bf H}$ est d\'efini sur $F$, on peut choisir $\iota$ d\'efini sur $F$, et identifier 
${\bf H}$ au sous--groupe (ferm\'e, d\'efini sur $F$) $\iota({\bf H})$ de $\Bbb{GL}_n$. Notons 
$\Bbb{H}$ le groupe quotient $N_{\Bbb{GL}_n}({\bf H})/Z_{\Bbb{GL}_n}({\bf H})$, et 
$\pi:N_{\Bbb{GL}_n}({\bf H})\rightarrow \Bbb{H}$ le morphisme quotient; 
il est d\'efini sur $F$, et induit par restriction un morphisme (de groupes algbriques)
${\bf H}\rightarrow \Bbb{H}$ qui est un isomorphisme sur un sous--groupe 
ferm\'e de $\Bbb{H}$ d\'efini sur $F$. On identifie ${\bf H}$ \`a ce sous--groupe. 
Puisque l'automorphisme $\tau$ est d\'efini sur $F$, la projection 
$\bar{g}=\pi(g)\in \Bbb{H}$ est $F$--rationnelle. \'Ecrivons les d\'ecompositions de Jordan 
$g=g_{\rm s}g_{\rm u}$ de $g$ et $\bar{g}=\bar{g}_{\rm s}\bar{g}_{\rm u}$ de $g$ et de $\bar{g}$. 
Alors $\tau_{\rm s}= {\rm Int}_{\Bbb{GL}_n}(g_{\rm s})\vert_{\bf H}$, 
$\tau_{\rm u}= {\rm Int}_{\Bbb{GL}_n}(g_{\rm u})\vert_{\bf H}$, 
$\bar{g}_{\rm s}=\pi(g_{\rm s})$ et $\bar{g}_{\rm u}=
\pi(g_{\rm u})$, et puisque par hypoth\`ese $\tau_{\rm s}$ et $\tau_{\rm u}$ sont 
d\'efinis sur $F$, $\bar{g}_{\rm s}$ et $\bar{g}_{\rm u}$ appartiennent \`a $\Bbb{H}(F)$. 
On peut donc appliquer le lemme: il existe un sous--groupe de Borel $\Bbb{B}$ de $\Bbb{H}$ 
d\'efini sur $F$, tel que $\bar{g}\in \Bbb{B}(F)$. 
\`A nouveau d'apr\`es \cite[ch.~IV, 11.14]{Bor}, 
$B=\Bbb{B}\cap {\bf H}$ est un sous--groupe de Borel de 
${\bf H}$, et il est d\'efini sur $F$. Comme 
$$\tau(B)= \bar{g}(\Bbb{B}\cap {\bf H})\bar{g}^{-1}= \bar{g}\Bbb{B}\bar{g}^{-1}\cap {\bf H}=
B,
$$
la proposition est d\'emontr\'ee.
\end{proof}

\subsection{Automorphismes stabilisant une paire de Borel d\'efinie sur $F^{\rm sep}$}\label{automorphismes stabilisant (B,T)} 
On prouve dans ce ${\rm n}^\circ$ qu'un $F$--automorphisme quasi--semisimple 
$\tau$ de ${\bf H}$ stabilise une paire de Borel de ${\bf H}$ d\'efinie sur $F^{\rm sep}$ 
--- pr\'ecis\'ement, on prouve qu'il existe un tore maximal $\tau$--admissible de 
${\bf H}$ d\'efini sur $F$ ---, et que le groupe ${\bf H}_\tau^\circ$ 
est d\'efini sur $F$.

Pour $\tau\in {\rm Aut}_F({\bf H})$, on note ${\bf H}_\tau^{\rm sep}$ le sous--groupe $F$--ferm\'e 
de ${\bf H}_\tau$ form\'e des composantes connexes qui poss\`edent un point $F^{\rm sep}$--rationnel, 
i.e. on pose
$$
{\bf H}_\tau^{\rm sep}=({\bf H}_\tau\cap {\bf H}(F^{\rm sep})){\bf H}_\tau^\circ.
$$
On a
$$({\bf H}_\tau^{\rm sep})^\circ={\bf H}_\tau^\circ.
$$
Par cons\'equent si ${\bf H}_\tau^{\rm sep}$ est d\'efini sur $F$, 
alors ${\bf H}_\tau^\circ$ l'est aussi. R\'eciproquement si ${\bf H}_\tau^\circ$ est d\'efini sur $F$, alors ${\bf H}_\tau^{\rm sep}$ 
est d\'efini sur $F^{\rm sep}$, donc sur $F$ puisqu'il est $F$--ferm\'e dans ${\bf H}$.

\begin{mesrems}{\rm 
\begin{enumerate}
\item[(1)]Pour $\tau\in {\rm Aut}_F({\bf H})$, le groupe ${\bf H}_\tau$ 
est $F$--ferm\'e dans ${\bf H}$ donc d\'efini sur $F^{p^{-\infty}}$. En particulier si 
$p=1$, alors ${\bf H}_\tau$ est d\'efini sur $F$, donc ${\bf H}_\tau^\circ$ l'est aussi.
\item[(2)]Soit $\tau\in {\rm Aut}_F({\bf H})$ tel que ${\bf H}_\tau^\circ$ est d\'efini sur $F$. D'apr\`es 
\cite[ch.~AG, 12.3 et 13.3]{Bor} on a ${\bf H}_\tau^{\rm sep}={\bf H}_\tau$ si et seulement si 
${\bf H}_\tau$ est d\'efini sur $F^{\rm sep}$, i.e. sur $F$ puisque ${\bf H}_\tau$ est $F$--ferm\'e dans ${\bf H}$. 
D'autre part pour toute sous--extension $F'/F$ de $F^{\rm sep}/F$, on a
\begin{align*}
{\bf H}_\tau\cap {\bf H}(F')&= ({\bf H}_\tau\cap {\bf H}(F^{\rm sep}))^{{\rm Gal}(F^{\rm sep}/F')}\\
& ={\bf H}_\tau^{\rm sep}(F^{\rm sep})^{{\rm Gal}(F^{\rm sep}/F')}\\
&={\bf H}_\tau^{\rm sep}(F').
\end{align*}
\item[(3)]On suppose $F=F^{\rm sep}$. Soit $\tau\in {\rm Aut}_F({\bf H})$ quasi--semisimple 
tel que ${\bf H}_\tau^\circ$ est d\'efini sur $F$. 
Choisissons un tore maximal ${\bf S}$ 
de ${\bf H}_\tau^\circ$ d\'efini sur $F$. Les groupes ${\bf T}=Z_{\bf H}({\bf S})$ et 
${\bf N}=N_{\bf H}({\bf T})$ sont $\tau$--stables et 
d\'efinis sur $F$. Soit ${\bf B}^\sharp$ un sous--groupe de Borel de ${\bf H}_\tau^\circ$ 
contenant ${\bf S}$, et soit $h\in {\bf H}_\tau$. Alors $h({\bf B}^\sharp,{\bf S})h^{-1}=
x({\bf B}^\sharp,{\bf S})x^{-1}$ pour un $x\in {\bf H}_\tau^\circ$, et 
$x^{-1}h\in N_{\bf H}({\bf S})\subset {\bf N}$. On a donc 
$${\bf H}_\tau= {\bf H}_\tau^\circ {\bf N}_\tau={\bf N}_\tau{\bf H}_\tau^\circ.$$ Posant ${\bf N}_\tau^{\rm sep}=
{\bf N}\cap {\bf H}_\tau^{\rm sep}$, on a aussi
$${\bf H}_\tau^{\rm sep}={\bf H}_\tau^\circ {\bf N}_\tau^{\rm sep}=
{\bf N}_\tau^{\rm sep}{\bf H}_\tau^\circ.$$ Puisque ${\bf T}$ est d\'efini et d\'eploy\'e sur $F$, 
${\bf T}_\tau$ l'est aussi, et ${\bf T}_\tau= {\bf T}_\tau(F){\bf T}_\tau^\circ$. 
Par cons\'equent ${\bf T}_\tau{\bf H}_\tau^\circ ={\bf T}_\tau(F){\bf H}_\tau^\circ$ 
est un sous--groupe ferm\'e de ${\bf H}_\tau^{\rm sep}$, d\'efini sur $F$, et ${\bf T}_\tau$ est contenu 
dans ${\bf N}_\tau^{\rm sep}$. Puisque ${\bf H}_\tau^{\rm sep}$ est d\'efini sur $F$, 
${\bf N}_\tau^{\rm sep}$ l'est aussi, et ${\bf H}_\tau$ est d\'efini sur $F$ si et seulement si ${\bf N}_\tau^{\rm sep}=
{\bf N}_\tau$.
\item[(4)]On suppose ${\bf H}$ semisimple et simplement connexe. Soit $\tau\in {\rm Aut}_F({\bf H})$ quasi--semisimple. 
Le morphisme $1-\tau:{\bf H}\rightarrow {\bf H}$ est s\'eparable (corollaire de \ref{automorphismes qss}), par suite ${\bf H}_\tau={\bf H}_\tau^\circ$ 
est d\'efini sur $F$. \hfill $\blacksquare$
\end{enumerate}}
\end{mesrems}

\begin{monlem}
Soit $\tau\in {\rm Aut}_F({\bf H})$ quasi--semisimple. 
Alors les trois conditions suivantes sont \'equivalentes\footnote{On verra plus loin (thorme) qu'elles sont v\'erifi\'ees.}:
\begin{itemize}
\item le groupe ${\bf H}_\tau^{\rm sep}$ est d\'efini sur $F$;
\item le groupe ${\bf H}_\tau^\circ$ est d\'efini sur $F$;
\item il existe un tore maximal $\tau$--admissible de ${\bf H}$ d\'efini sur $F$.
\end{itemize}
\end{monlem}

\begin{proof}
L'\'equivalence des deux premi\`eres conditions a \'et\'e d\'emontr\'ee 
plus haut. Si ${\bf H}_\tau^\circ$ est d\'efini sur $F$, alors il existe un tore 
maximal ${\bf S}$ de ${\bf H}_\tau^\circ$ d\'efini sur $F$, et ${\bf T}=Z_{\bf H}({\bf S})$ 
est un tore maximal $\tau$--admissible de ${\bf H}$ d\'efini sur $F$.

Supposons qu'il existe un tore maximal $\tau$--admissible ${\bf T}$ 
de ${\bf H}$ d\'efini sur $F$, et montrons que le groupe ${\bf H}_\tau^\circ$ est d\'efini sur $F$. 
Comme ${\bf H}_\tau$ est $F$--ferm\'e dans ${\bf H}$, ${\bf H}_\tau^\circ$ l'est aussi et il 
suffit de montrer que ${\bf H}_\tau^\circ$ est d\'efini sur $F^{\rm sep}$. 
On peut donc 
supposer $F=F^{\rm sep}$. 
Soit ${\bf B}$ un sous--groupe de Borel $\tau$--stable de ${\bf H}$ contenant ${\bf T}$. 
Alors $({\bf B}^\sharp\!,{\bf T}^\sharp)=({\bf B}\cap {\bf H}_\tau^\circ, {\bf T}\cap {\bf H}_\tau^\circ)$ 
est une paire de Borel de ${\bf H}_\tau^\circ$ d\'efinie sur $F$ (puisque ${\bf T}$ est d\'efini et d\'eploy\'e sur $F$, 
${\bf T}^\sharp$ l'est aussi). 
Notons ${\bf B}'$ le 
sous--groupe de Borel de ${\bf H}$ oppos\'e \`a ${\bf B}$ 
par rapport \`a ${\bf T}$. Alors ${\bf B}'^\sharp= {\bf B}'\cap {\bf H}_\tau^\circ$ est 
le sous--groupe de Borel de ${\bf H}_\tau^\circ$ oppos\'e \`a ${\bf B}^\sharp$ 
par rapport \`a ${\bf T}^\sharp$, et tout comme ${\bf B}^\sharp$, 
${\bf B}'^\sharp$ est d\'efini sur $F$. Le groupe ${\bf U}'^\sharp=R_{\rm u}({\bf B}'^\sharp)$ 
est d\'efini sur $F$, et l'application produit 
${\bf U}'^\sharp\times {\bf B}^\sharp\rightarrow  {\bf H}_\tau^\circ$ est un isomophisme de vari\'et\'es alg\'ebriques 
sur un ouvert de ${\bf H}_\tau^\circ$ \cite[ch.~IV, 14.14]{Bor}. Cet ouvert est d\'efini sur $F$ 
et engendre ${\bf H}_\tau^\circ$, par cons\'equent ${\bf H}_\tau^\circ$ est d\'efini sur $F$.
\end{proof}

On peut maintenant d\'emontrer le r\'esultat principal de ce ${\rm n}^\circ$:

\begin{montheo}
Soit $\tau\in {\rm Aut}_F({\bf H})$ 
quasi--semisimple. Le groupe ${\bf H}_\tau^\circ$ est d\'efini sur $F$, et il existe un 
tore maximal $\tau$--admissible de ${\bf H}$ d\'efini sur $F$. En particulier, il existe une 
paire de Borel $\tau$--stable de ${\bf H}$ d\'efinie sur $F^{\rm sep}$.
\end{montheo}

\begin{proof}
Il suffit de montrer que ${\bf H}_\tau^\circ$ est d\'efini sur $F^{\rm sep}$. 
On peut supposer $p>1$ et $F=F^{\rm sep}$. La d\'emonstration 
s'organise comme suit: on commence par se ramener au cas o\`u $\tau$ est localement fini. 
Puis par descente de ${\bf H}$ \`a ${\bf H}_{\tau_{\rm s}}^\circ$, on se ram\`ene au cas o\`u 
$\tau$ est (quasi--semisimple) unipotent, ce qui permet d'utiliser le lemme de \ref{automorphismes qss lf}.

D'apr\`es le lemme de \ref{groupes rductifs}, 
on a ${\bf H}_\tau^\circ = R({\bf H})_\tau^\circ ({\bf H}_{\rm der})_{\tau_{\rm der}}^\circ$. Puisque $R({\bf H})$ est un tore d\'efini sur $F\;(=F^{\rm sep})$, il est d\'eploy\'e sur $F$ 
\cite[ch.~III, 8.11]{Bor}, et $R({\bf H})_\tau^\circ$ 
est lui aussi un tore d\'efini et d\'eploy\'e sur $F$ \cite[ch.~III, 8.4]{Bor}. Par 
cons\'equent si $({\bf H}_{\rm der})_{\tau}^\circ$ est d\'efini sur 
$F$, alors ${\bf H}_\tau^\circ$ l'est aussi (c'est l'image du morphisme produit
$R({\bf H})_\tau^\circ \times ({\bf H}_{\rm der})_{\tau_{\rm der}}^\circ\rightarrow {\bf H}$). Quitte \`a remplacer ${\bf H}$ 
par ${\bf H}_{\rm der}$ et $\tau$ par $\tau_{\rm der}$, on peut donc 
supposer $\tau\in {\rm Aut}_F^0({\bf H})$. Alors on \'ecrit la d\'ecomposition de 
Jordan $\tau=\tau_{\rm s}\circ \tau_{\rm u}$. Rappelons que $\tau_{\rm s}$ et 
$\tau_{\rm u}$ appartiennent \`a ${\rm Aut}_{F^{p^{-\infty}}}({\bf H})$.

Montrons que le groupe ${\bf H}_{\tau_{\rm s}}^\circ$ est d\'efini sur $F$. 
Identifions ${\bf H}$ \`a la composante neutre 
du groupe alg\'ebrique affine ${\bf H}'={\bf H}\rtimes \langle \tau \rangle /{\bf C}$ comme en \ref{automorphismes ss et u}, et 
notons $\delta$ l'image de $1\rtimes \tau$ dans ${\bf H}'$. Puisque ${\bf H}$ et $\tau$ sont d\'efinis sur 
$F$, ${\bf H}'$ l'est aussi et $\delta$ appartient  ${\bf H}'(F)$. \'Ecrivons la d\'ecomposition de Jordan 
$\delta=\delta_{\rm s}\delta_{\rm u}$. On a $\delta_{\rm s}\in {\bf H}(F^{p^{-\infty}})$, 
$\tau_{\rm s}={\rm Int}_{{\bf H}'}(\delta_{\rm s})^\circ$ et 
${\bf H}_{\tau_{\rm s}}^\circ= 
({\bf H}'_{\delta_{\rm s}})^\circ$. Puisque ${\bf H}'$ est affine et d\'efini sur $F$, 
il existe un morphisme de groupes 
alg\'ebriques
$$\iota:{\bf H'}\rightarrow \Bbb{GL}_n
$$ d\'efini sur $F$, qui soit un 
isomorphisme sur un sous--groupe ferm\'e de de $\Bbb{GL}_n$. Identifions 
${\bf H}'$ \`a $\iota({\bf H}')$. Notons $\Bbb{T}$ le tore maximal diagonal de $\Bbb{GL}_n$. 
Il existe un $g\in {\rm GL}_n(\overline{F})$ 
tel que
$$g^{-1}\delta_{\rm s}g=
{\rm diag}(x_1,\ldots,x_1; x_2,\ldots ,x_2; \ldots ; x_k,\ldots ,x_k)\in \Bbb{T}
$$ 
pour des $x_i\in \smash{\overline{F}}^\times$ deux--\`a--deux distincts. 
Posons $y=g^{-1}\delta_{\rm s}g$.  On a  
$$
(\Bbb{GL}_n)_y=\Bbb{GL}_{n_1}\times \cdots \times \Bbb{GL}_{n_k}, \quad 
n_1+\cdots + n_k =n,
$$
o\`u $n_i$ est 
la multiplicit\'e de la valeur propre $x_i$. Pour $m\in \Bbb{Z}_{\geq 1}$, 
puisque $p>1$, l'application
$$\smash{\overline{F}}^\times\rightarrow \smash{\overline{F}}^\times,\,x\mapsto x^{p^m}$$ 
est bijective, et on a l'galit $(\Bbb{GL}_n)_{y^{p^m}}=(\Bbb{GL}_n)_y$; par suite 
$(\Bbb{GL}_n)_{\delta_{\rm s}^{p^m}}=(\Bbb{GL}_n)_{\delta_{\rm s}}$, d'o\`u
$$
{\bf H}'_{\delta_{\rm s}^{p^m}}={\bf H}'\cap (\Bbb{GL}_n)_{\delta_{\rm s}^{p^m}}=
{\bf H}'\cap (\Bbb{GL}_n)_{\delta_s}={\bf H}'_{\delta_s}.
$$
Choisissons un 
entier $m\geq 1$ tel que $(\tau_{\rm u})^{p^m}={\rm id}_{\bf H}$, et 
posons $\sigma = \tau^{p^m}= (\tau_{\rm s})^{p^m}$. Puisque $\sigma$ est semisimple et 
d\'efini sur $F$, ${\bf H}_{\tau_{\rm s}}^\circ = {\bf H}_\sigma^\circ$ est d\'efini sur $F$.

\v1
Comme $\tau$ est d\'efini sur $F$, $\tau^*\;(=\tau\vert_{{\bf H}_{\tau_{\rm s}^\circ}})$ l'est aussi. 
Or d'apr\`es la proposition et le corollaire de \ref{automorphismes qss lf}, 
$\tau^*$ est quasi--semisimple et unipotent,  
et ${\bf H}_\tau^\circ = ({\bf H}_{\tau_{\rm s}}^\circ)_{\tau^*}$. 
Quitte \`a remplacer ${\bf H}$ par ${\bf H}_{\tau_{\rm s}}^\circ$ et $\tau$ par $\tau^*$, 
on peut donc supposer $\tau$ unipotent. 
D'apr\`es la proposition de \ref{automorphismes stabilisant B}, $\tau$ stabilise un sous--groupe de Borel ${\bf B}$ de 
${\bf H}$ d\'efini sur $F$. Puisque ${\bf B}_\circ$ (cf. \ref{la dcomposition}) et ${\bf B}$ sont d\'efinis 
sur $F$, il existe un $y\in {\bf H}(F)$ tel que ${\bf B}=y^{-1}{\bf B}_\circ y$. 
Quitte \`a remplacer $\tau$ par ${\rm Int}_{\bf H}(y)\circ \tau\circ {\rm Int}_{\bf H}(y^{-1})$, 
on peut supposer $\tau({\bf B}_\circ)={\bf B}_\circ$. 
\'Ecrivons $\tau = {\rm Int}_{\bf H}(h)\circ \tau_\circ$ avec $h\in {\bf H}$ 
et $\tau_\circ \in \frak{A}_\circ$ (relation $(**)$ de \ref{groupes rductifs}). Soit un entier $m\geq 1$ tel 
que $\tau^{p^m}={\rm id}_{\bf H}$. Puisque
$$
{\rm Int}_{\bf H}(h\tau_\circ(h)\cdots \tau_\circ^{p^m-1}(h))\circ \tau_\circ^{p^m}=\tau^{p^m}= {\rm id}_{\bf H},
$$
d'apr\`es la relation $(**)$ de \ref{la dcomposition}, on a $h\tau_\circ(h)\cdots \tau_\circ^{p^m-1}(h)\in Z({\bf H})$ et $\tau_\circ^{p^m}={\rm id}_{\bf H}$. 
Donc $\tau_\circ$ est unipotent, et d'apr\`es le lemme de \ref{automorphismes qss lf}, il existe un $x\in {\bf H}$ tel que 
$$
\tau={\rm Int}_{\bf H}(x^{-1})\circ \tau_\circ \circ {\rm Int}_{\bf H}(x)={\rm Int}_{\bf H}(x^{-1}\tau_\circ(x))
\circ \tau_\circ.
$$
Comme $\tau\circ\tau_\circ^{-1}({\bf B}_\circ)={\bf B}_\circ$, 
$x^{-1}\tau_\circ(x)$ appartient \`a ${\bf B}_\circ$.

Puisque ${\bf V}_\circ\;(={\bf U}_\circ^-{\bf B}_\circ)$ est ouvert dans ${\bf H}$ et que ${\bf H}(F)$ est 
dense dans ${\bf H}$, quitte \`a remplacer $\tau$ par 
${\rm Int}_{\bf H}(y')\circ \tau\circ {\rm Int}_{\bf H}(y'^{-1})$ pour un $y'\in {\bf H}(F)$, on peut 
supposer $x\in {\bf V}_\circ$. \'Ecrivons $x=\bar{u}b$ avec 
$\bar{u}\in {\bf U}_\circ^-$ et $b\in {\bf B}_\circ$. 
Comme $x^{-1}\tau_\circ(x)= b^{-1}\bar{u}^{-1}\tau_\circ(\bar{u})\tau_\circ(b)\in {\bf B}_\circ$, on a 
$\bar{u}^{-1}\tau_\circ(\bar{u})\in {\bf B}_\circ$; or ${\bf U}_\circ^-\cap {\bf B}_\circ =\{1\}$, donc 
$\bar{u}\in {\bf U}_\circ^-\cap {\bf H}_{\tau_\circ}$. On obtient
$$
\tau={\rm Int}_{\bf H}(b^{-1}\tau_\circ(b))
\circ \tau_\circ
={\rm Int}_{\bf H}(b^{-1})\circ \tau_\circ \circ {\rm Int}_{\bf H}(b).
$$
Puisque $\tau$ et $\tau_\circ$ sont d\'efinis sur $F$, 
l'image $\bar{b}_1$ de $b_1=b^{-1}\tau_\circ(b)$ dans ${\bf B}_\circ/Z({\bf H})$ est $F$--rationnelle. 
\'Ecrivons $b= tu$, $t\in {\bf T}_\circ$, $u\in {\bf U}_\circ$, et posons $t_1=t^{-1}\tau_\circ(t)$ et 
$\tau'_\circ = {\rm Int}_H(t_1)\circ \tau_\circ$. On a
$$
b_1= u^{-1}t^{-1}\tau_\circ(t) \tau_\circ(u)= u^{-1}\tau'_\circ(u) t_1
$$
o $u^{-1}\tau'_\circ(u)\in {\bf U}_\circ$ et 
$t_1\in {\bf T}_\circ$. Comme $\bar{b}_1$ est $F$--rationnel, $u^{-1}\tau'_\circ (u)$ et 
l'image de $t_1$ dans ${\bf T}/Z({\bf G})$ le sont aussi. En particulier, 
$\tau'_\circ$ est d\'efini sur $F$ et stabilise la paire de Borel $({\bf B}_\circ,{\bf T}_\circ)$ de ${\bf H}$. D'apr\`es la proposition 
de \ref{automorphismes qss}, le morphisme 
$(1-\tau'_\circ)\vert_{{\bf U}_\circ}$ est s\'eparable. D'autre part 
le groupe $({\bf U}_\circ)_{\tau'_\circ}$ est connexe (thorme de \ref{automorphismes qss}), 
r\'esoluble (car nilpotent), et 
d\'eploy\'e sur $F\;(=F^{\rm sep})$. Comme $\tau'_\circ(u^{-1})u\in {\bf U}_\circ(F)$, cela implique 
\cite[ch.~V, 15.12]{Bor} qu'il 
existe un \'el\'ement $u'\in {\bf U}_\circ(F)$ tel que $\tau'_\circ(u^{-1})u= \tau'_\circ(u'^{-1})u'$. 
Donc $u^{-1}\tau_\circ(u)=u'^{-1}\tau_\circ(u')$ et
$$
\tau= {\rm Int}_{\bf H}(u^{-1}\tau'_\circ(u)t_1)\circ \tau_\circ = {\rm Int}_{\bf H}(u'^{-1}\tau'_\circ(u')))
\circ \tau'_\circ ={\rm Int}_{\bf H}(u'^{-1})\circ \tau'_\circ \circ {\rm Int}_{\bf H}(u').
$$
Le tore maximal $u'^{-1}{\bf T}_\circ u'$ de ${\bf H}$ est d\'efini sur $F$ et $\tau$--admissible. 
Puisque ce tore est d\'eploy\'e sur $F$, tout sous--groupe de Borel de 
${\bf H}$ le contenant est d\'efini et d\'eploy\'e sur $F$, ce qui ach\`eve la d\'emonstration 
du th\'eor\`eme.
\end{proof}

\begin{moncoro}
Soit $\tau\in {\rm Aut}_F({\bf H})$ 
quasi--central. Le groupe ${\bf H}_\tau$ est d\'efini sur $F$.
\end{moncoro}

\begin{proof}
On peut supposer $F=F^{\rm sep}$. Choisissons un tore maximal ${\bf S}$ 
de ${\bf H}_\tau^\circ$ d\'efini sur $F$, et posons ${\bf T}=Z_{\bf H}({\bf S})$, 
${\bf N}=N_{\bf H}({\bf T})$ et ${\bf W}={\bf N}/{\bf T}$. Les groupes ${\bf N}$, ${\bf T}$ et 
${\bf W}$ sont d\'efinis sur $F$, et d'apr\`es la remarque (3), on a l'galit ${\bf H}_\tau={\bf N}_\tau{\bf H}_\tau^\circ$. 
Puisque $\tau$ est quasi--central, tout \'el\'ement $\tau$--stable de ${\bf W}$ se rel\`eve en un \'el\'ement de 
${\bf N}\cap {\bf H}_\tau^\circ$ (remarque (1) de \ref{automorphismes qc}). On a donc ${\bf N}_\tau= ({\bf N}\cap {\bf H}_\tau^\circ){\bf T}_\tau=
{\bf T}_\tau({\bf N}\cap {\bf H}_\tau^\circ)$. Par cons\'equent ${\bf H}_\tau=
{\bf T}_\tau{\bf H}_\tau^\circ= {\bf T}_\tau(F){\bf H}_\tau^\circ$, d'o\`u le r\'esultat.
\end{proof}

\begin{marema}{\rm 
Soit $\tau\in {\rm Aut}_F({\bf H})$ quasi--semisimple. On 
peut pr\'eciser la remarque (3), dont on reprend les hypoth\`eses et les notations (en 
particulier $F=F^{\rm sep}$, et $\tau\in {\rm Aut}_F({\bf H})$ est quasi--semisimple). 

\begin{enumerate}
\item[(5)]
Posons ${\bf W}={\bf N}/{\bf T}$, 
et notons ${\bf W}_\tau$ le sous--groupe de ${\bf W}$ form\'e des \'el\'ements qui sont 
$\tau$--stables. La projection canonique ${\bf N}\rightarrow {\bf W}$ induit un morphisme injectif 
de groupes (et m\^eme de groupes alg\'ebriques) ${\bf N}_\tau/{\bf T}_\tau \rightarrow {\bf W}_\tau$. 
On note ${\bf W}_\tau^*\subset {\bf W}_\tau$ son image.

Puisque ${\bf T}$ est d\'efini et d\'eploy\'e sur 
$F$, pour chaque $\tau$--orbite $\ES{O}$ dans $\Phi({\bf T},{\bf H})$, l'\'el\'ement $y_{\tau,\ES{O}}$ 
d\'efini dans la remarque 4 de \ref{automorphismes qss}, appartient \`a $F$. On en d\'eduit (remarque (3) de \ref{automorphismes qc}) 
qu'il existe un \'el\'ement 
$t\in {\bf T}(F)$ tel que le $F$--automorphisme $\tau'={\rm Int}_{\bf H}(t)\circ\tau$ de ${\bf H}$ est 
quasi--central. Le tore maximal ${\bf T}$ de ${\bf H}$ est $\tau'$-admissible, et l'on a 
${\bf W}_{\tau'}={\bf W}_\tau$ et ${\bf T}_{\tau'}={\bf T}_\tau$. Puisque $\tau'$ est d\'efini 
sur $F$, ${\bf H}_{\tau'}$ l'est aussi (corollaire), et ${\bf N}_{\tau'}={\bf N}\cap {\bf N}_{\tau'}$ est 
d\'efini sur $F$. Notons que d'apr\`es la remarque (1) de \ref{automorphismes qc}, on a l'galit 
${\bf W}_{\tau'}^*={\bf W}_{\tau'}$. 

Soit maintenant un \'el\'ement $n\in {\bf N}_\tau$. Alors $w=n{\bf T}$ appartient \`a 
${\bf W}_{\tau}^*\subset {\bf W}_\tau={\bf W}_{\tau'}^*$, donc se rel\`eve en un \'el\'ement 
$n'\in {\bf N}_{\tau'}$, que l'on peut choisir dans ${\bf N}_{\tau'}(F)$ puisque 
${\bf N}_{\tau'}={\bf N}_{\tau'}(F){\bf T}_\tau^\circ$. 
Soit $x=nn'^{-1}\in {\bf T}$. Posant $t^{w-1}=ntn^{-1}t^{-1}=n'tn'^{-1}t^{-1}$, on a 
$t^{w-1}=x\tau(x)^{-1}\in  {\bf T}(1-\tau)$. R\'eciproquement 
si $w\in {\bf W}_\tau$ v\'erifie $t^{w-1}\in {\bf T}(1-\tau)$, alors $w$ se rel\`eve \`a ${\bf N}_\tau$. On a 
donc
$$
{\bf W}_\tau^*=\{w\in {\bf W}_\tau: t^{w-1}\in {\bf T}(1-\tau)\}.
$$
Puisque $t\in {\bf T}(F)$ et ${\bf W}={\bf W}(F)$, l'\'el\'ement $t^{w-1}$ appartient au groupe 
${\bf T}(1-\tau)(F)$ 
des points $F$--rationnels de ${\bf T}(1-\tau)$. Notons ${\bf T}(F)(1-\tau)$ le sous--groupe de 
${\bf T}(1-\tau)(F)$ form\'e des $y\tau(y)^{-1}$ pour $y\in {\bf T}(F)$. L'isomorphisme de groupes 
${\bf N}_\tau/{\bf T}_\tau\rightarrow {\bf W}_\tau^*$ induit par restriction une application injective 
$${\bf N}_\tau^{\rm sep}(F)/{\bf T}_\tau(F)\rightarrow {\bf W}_\tau^*,$$ o\`u (rappel) 
${\bf N}_\tau^{\rm sep}={\bf N}\cap {\bf H}_\tau^{\rm sep}$. Soit ${\bf W}_\tau^{*,{\rm sep}}
\subset {\bf W}_\tau^*$ son image. Par d\'efinition, ${\bf W}_\tau^{*,{\rm sep}}$ est le sous--groupe 
de ${\bf W}_\tau$ form\'e des \'el\'ements qui se rel\`event \`a 
${\bf N}_\tau^{\rm sep}(F)={\bf N}_\tau\cap {\bf H}(F^{\rm sep})$, et de la m\^eme mani\`ere, on obtient l'galit
$$
{\bf W}_\tau^{*,{\rm sep}}=\{w\in {\bf W}_\tau: t^{w-1}\in {\bf T}(F)(1-\tau)\}.
$$
Puisque ${\bf T}_\tau{\bf H}_\tau^\circ={\bf H}_\tau^\circ{\bf T}_\tau$ est d\'efini sur $F$, le groupe ${\bf H}_\tau$ est d\'efini 
sur $F$ si et seulement si ${\bf W}_\tau^{*,{\rm sep}}={\bf W}_\tau^*$. \hfill $\blacksquare$
\end{enumerate}}
\end{marema}

\subsection{Tores maximaux et sous--espaces de Cartan de ${\bf H}^\natural(F)$}\label{tores maximaux et sous-espaces de Cartan} 
Pour $\delta\in H^\natural_{\rm reg}$, le groupe ${\bf H}_\delta^\circ$ est un tore, et comme 
$\tau={\rm Int}_{\bf H}(\delta)$ est quasi--semisimple, il est d\'efini sur $F$. Par cons\'equent le groupe 
${\bf T}=Z_{\bf H}({\bf H}_\delta^\circ)$ est lui aussi d\'efini sur $F$ \cite[ch.~III, 9.2]{Bor} --- notons 
que puisque ${\bf T}$ est l'unique tore maximal $\tau$--admissible de ${\bf H}$, 
cette assertion r\'esulte aussi du lemme de \ref{automorphismes stabilisant B}.

\begin{madefi}
{\rm 
On appelle:
\begin{itemize}
\item {\it tore maximal de $H^\natural$} l'ensemble des points $F$--rationnels ${\bf S}^\natural(F)$ 
d'un tore maximal ${\bf S}^\natural$ de ${\bf H}^\natural$ d\'efini sur $F$ et tel que ${\bf S}^\natural\cap 
H_{\rm reg}^\natural\neq \emptyset$.
\item {\it sous--espace de Cartan de $H^\natural$} l'ensemble des points $F$--rationnels ${\bf T}^\natural(F)$ d'un sous--espace de Cartan ${\bf T}^\natural$ de ${\bf H}^\natural$ d\'efini sur $F$ et tel que ${\bf T}^\natural\cap H_{\rm reg}^\natural\neq \emptyset$.
\end{itemize}
}\end{madefi}

Par d\'efinition, un tore maximal de $H^\natural$ est une partie de la forme 
${\bf H}^\circ_\delta(F)\cdot \delta$ pour un $\delta\in H^\natural_{\rm reg}$, et 
un sous--espace de Cartan de $H^\natural$ est une partie de la forme 
$Z_{\bf H}({\bf H}_\delta^\circ)(F)\cdot \delta$ pour un $\delta\in H^\natural_{\rm reg}$. 
Tout tore maximal $S^\natural$ de $H^\natural$ d\'efinit un 
{\it quadruplet de Cartan $(S,S^\natural,T,T^\natural)$ de $H^\natural$}:
\begin{itemize}
\item $S={\bf H}_\delta^\circ(F)$ pour un (resp. pour tout) $\delta\in S^\natural \cap H^\natural_{\rm reg}$ --- c'est un tore de $H$, et $S^\natural$ est un 
$S$--espace tordu trivial;
\item $T=Z_{\bf H}({\bf H}_\delta^\circ)(F)$ pour 
un (resp. pour tout) $\delta\in S^\natural \cap H^\natural_{\rm reg}$ --- c'est un tore maximal de $H$;
\item $T^\natural= T\cdot S^\natural = S^\natural\cdot T$ --- c'est un $T$--espace tordu.
\end{itemize}
De m\^eme, tout sous--espace de Cartan $T^\natural$ de $H^\natural$ d\'efinit un {\it triplet de Cartan $(S,T,T^\natural)$ de $H^\natural$}:
\begin{itemize}
\item $S={\bf H}_\delta^\circ(F)$ pour un (resp. pour tout) $\delta \in T^\natural\cap H^\natural_{\rm reg}$;
\item $T= Z_{\bf H}({\bf H}_\delta^\circ)(F)$ pour un (resp. pour tout) $\delta \in T^\natural\cap H^\natural_{\rm reg}$.
\end{itemize}
L'application qui a $\delta\in T^\natural$ associe le quadruplet 
$(S,S\cdot\delta, T,T^\natural)$ est une bijection du $S\backslash T$--espace tordu 
$S\backslash T^\natural$ sur l'ensemble des quadruplets de Cartan de $H^\natural$ prolongeant 
$(S,T,T^\natural)$. 

\begin{mapropo}
Supposons le corps $F$ infini, et soit $(S,S^\natural,T,T^\natural)$ 
un quadruplet de Cartan de $H^\natural$. Alors on a $S=Z_H(S^\natural)=Z_H(T^\natural)$ et $T=Z_H(S)$.
\end{mapropo}

\begin{proof}
Soit ${\bf S}^\natural$ le sous--espace de Cartan de ${\bf H}$ 
donn\'e par ${\bf S}^\natural={\bf H}_\delta^\circ\cdot \delta$ pour un (resp. pour tout) 
$\delta\in S^\natural\cap H^\natural_{\rm reg}$, et soit $({\bf S},{\bf S}^\natural,{\bf T}, {\bf T}^\natural)$ le 
quadruplet de Cartan de ${\bf H}^\natural$ associ\'e \`a ${\bf S}^\natural$. Alors ${\bf S}$, ${\bf S}^\natural$, 
${\bf T}$, ${\bf T}^\natural$ sont d\'efinis sur $F$ et leurs ensembles de points $F$--rationnels 
co\"{\i}ncident avec $S$, $S^\natural$, $T$, $T^\natural$. Posant 
$Z_{\bf H}(S^\natural)=\{h\in {\bf H}: \delta\cdot h= h\cdot \delta \}$, on a l'inclusion 
$Z_{\bf H}({\bf S}^\natural)\subset Z_{\bf H}(S^\natural)$. Soit $h\in {\bf H}$. Posons 
${\bf S}^\natural_h=\{\delta\in {\bf S}^\natural: \delta \cdot h =h\cdot \delta\}$. C'est une sous--vari\'et\'e 
ferm\'ee de ${\bf S}^\natural$, et l'on a ${\bf S}^\natural_h={\bf S}^\natural$ si et seulement 
si $h\in Z_{\bf H}({\bf S}^\natural)$. De m\^eme, posant $S^\natural_h=\{\delta\in S^\natural: \delta\cdot h=h\cdot \delta\}$, on a $S^\natural_h = S^\natural$ si et seulement si $h\in Z_{\bf H}(S^\natural)$. Puisque $F$ est infini, 
$S^\natural$ est dense dans ${\bf S}^\natural$ (cf. la remarque de \ref{points rationnels d'un espace tordu}), par cons\'equent si $S^\natural_h=S^\natural$ 
alors ${\bf S}^\natural_h={\bf S}^\natural$, d'o\`u l'inclusion $Z_{\bf H}(S^\natural)\subset Z_{\bf H}({\bf S}^\natural) 
\;(={\bf S})$, puis l'\'egalit\'e $Z_H(S^\natural)=S$. Le m\^eme raisonnement entra\^{\i}ne que 
$Z_H(T^\natural)=S$, et que $Z_H(S)=T$.
\end{proof}

\begin{moncoro}
\begin{enumerate}
\item[(1)]Deux quadruplets de Cartan 
$(S,S^\natural,T,T^\natural)$ et $(S'\!,S'^\natural,T'\!,T'^\natural)$ de 
$H^\natural$ sont conjugu\'es dans $H$ si et seulement si les tores maximaux 
$S^\natural$ et $S'^\natural$ de $H^\natural$ le sont.
\item[(2)]Deux triplets de Cartan 
$(S,T,T^\natural)$ et $(S'\!,T'\!,T'^\natural)$ de 
$H^\natural$ sont conjugu\'es dans $H$ si et seulement si les sous--espaces 
de Cartan $T^\natural$ et $T'^\natural$ de $H^\natural$ le sont.
\end{enumerate}
\end{moncoro}

\begin{proof}
Prouvons (1). Il s'agit de montrer que si $S'^\natural= h\cdot S^\natural\cdot h^{-1}$ pour un $h\in H$, alors on a
$$(S'\!,S'^\natural,T'\!,T'^\natural)= h\cdot (S,S^\natural,T,T^\natural)\cdot h^{-1}.
$$ Puisque 
$S=Z_H(S^\natural)$ et $S'=Z_H(S'^\natural)$, on a $S'=hSh^{-1}$, et puisque $T=Z_H(S)$ et 
$T'=Z_H(S')$, on a $T'=hTh^{-1}$. On en d\'eduit que
$$
T'^\natural = T'\cdot S'^\natural = (hTh^{-1})\cdot (h\cdot S^\natural\cdot h^{-1}) = h\cdot (T\cdot S^\natural)\cdot h^{-1}= h\cdot T^\natural \cdot h^{-1}.
$$
Le point (2) s'obtient de la m\^eme mani\`ere, en utilisant que $S=Z_H(T^\natural)$. 
\end{proof}

\subsection{${\bf H}(F)$--orbites dans ${\bf H}^\natural(F)$}\label{H(F)-orbites} 
Pour $\delta\in H^\natural$, 
la ${\bf H}$--orbite $\ES{O}_{\bf H}(\delta)$ est 
d\'efinie sur 
$F$, et l'ensem\-ble $\ES{O}_{\bf H}(\delta)(F)$ de ses 
points $F$--rationnels
est r\'eunion de $H$--orbites de la forme
$$
\ES{O}_H(\delta')=\{h^{-1}\cdot \delta' \cdot h: h\in H\}
$$
pour $\delta'\in \ES{O}_{\bf H}(\delta)(F)$. D'autre part, le morphisme
$\pi_\delta:{\bf H}\rightarrow \ES{O}_{\bf H}(\delta),\,h\mapsto h^{-1}\cdot \delta\cdot h$ 
est d\'efini sur $F$, et s'il est s\'eparable (e.g. si $p=1$, ou si $\delta$ est semisimple d'apr\`es 
\cite[ch.~III, 9.1]{Bor}), alors ${\bf H}_\delta$ est d\'efini sur $F$ et $\ES{O}_{\bf H}(\delta)$ est 
\og le\fg{} quotient de ${\bf H}$ par ${\bf H}_\delta$ \cite[ch.~II, 6.7]{Bor}; en ce cas 
$\pi_\delta$ induit une application 
surjective ${\bf H}(F^{\rm sep})\rightarrow \ES{O}_{\bf H}(\delta)(F^{\rm sep})$, et l'\'etude des $H$--orbites dans 
$\ES{O}_{\bf H}(\delta)(F)$ se ram\`ene \`a un probl\`eme de 
cohomologie galoisienne.

\begin{mesrems}{\rm 
\begin{enumerate}
\item[(1)]Soit un \'el\'ement $\delta\in H^\natural$ tel que le groupe ${\bf H}_\delta$ est d\'efini sur $F$. 
Alors le morphisme bijectif
$$\bar{\pi}_\delta:{\bf H}_\delta\backslash {\bf H}\rightarrow \ES{O}_{\bf H}(\delta),\,h\mapsto h^{-1}\cdot\delta\cdot h$$ est 
d\'efini sur $F$. Par passage aux points $F$--rationnels, il induit une application injective 
$$({\bf H}_\delta\backslash {\bf H})(F)\hookrightarrow \ES{O}_{\bf H}(\delta)(F)$$ qui n'est en g\'en\'eral 
pas surjective, m\^eme si $F=F^{\rm sep}$ --- elle l'est si le morphisme $\pi_\delta$ est s\'eparable, puisqu'en ce cas $\bar{\pi}_\delta$ est un 
isomorphisme. Comme le morphisme quotient ${\bf H}\rightarrow {\bf H}_\delta\backslash {\bf H}$ 
induit une application surjective ${\bf H}(F^{\rm sep})\rightarrow ({\bf H}_\delta\backslash {\bf H})(F^{\rm sep})$, 
l'application ${\bf H}(F^{\rm sep})\rightarrow \ES{O}_{\bf H}(\delta)(F^{\rm sep}),\,h\mapsto h^{-1}\cdot \delta\cdot h$ est 
surjective si et seulement si l'application $({\bf H}_\delta\backslash {\bf H})(F^{\rm sep})\rightarrow \ES{O}_{\bf H}(\delta)(F^{\rm sep})$ est bijective. 
Si maintenant on suppose seulement que ${\bf H}_\delta^\circ$ est d\'efini sur $F$ (e.g. si $\delta$ quasi--semisimple, d'aprs le thorme de 
\ref{automorphismes stabilisant (B,T)}), 
alors ${\bf H}_\delta^{\rm sep} =({\bf H}_\delta\cap {\bf H}(F^{\rm sep})){\bf H}_\delta^\circ$ est d\'efini sur $F$, 
et puisque ${\bf H}_\delta^{\rm sep}(F^{\rm sep})= {\bf H}_\delta\cap {\bf H}(F^{\rm sep})$, on obtient aussi 
que l'application ${\bf H}(F^{\rm sep})\rightarrow \ES{O}_{\bf H}(\delta)(F^{\rm sep}),\,h\mapsto h^{-1}\cdot \delta\cdot h$ est 
surjective si et seulement si l'application $({\bf H}_\delta^{\rm sep}\backslash {\bf H})(F^{\rm sep})\rightarrow \ES{O}_{\bf H}(\delta)(F^{\rm sep})$ 
est bijective.
\item[(2)]Soit $\delta\in H^\natural$ quasi--semisimple tel que le morphisme $\pi_\delta$ est 
s\'eparable. Supposons que le groupe ${\bf H}_\delta$ est connexe (e.g. si $\delta$ est unipotent, ou 
si ${\bf H}$ est semisimple simplement connexe). 
Si $\delta$ est r\'egulier, alors ${\bf H}_\delta= {\bf H}_\delta^\circ$ est un tore. Ce 
tore se d\'eploie sur une sous--extension 
finie $F_1/F$ de $F^{\rm sep}/F$, et d'apr\`es la remarque 2 de \ref{gnralits}, pour toute sous--extension $F'/F_1$ de $\overline{F}/F_1$, 
$\pi_\delta$ induit une application surjective 
${\bf H}(F')\rightarrow \ES{O}_{\bf H}(\delta)(F')$. 
\item[(3)]Supposons que $F$ est un corps \og de type (F)\fg{} au sens 
de \cite[ch.~III, \S4.2]{Se}; i.e. que $F$ est parfait et que pour chaque entier $n\geq 1$, il 
n'existe qu'un nombre fini de sous--extensions de $\overline{F}/F$ de degr\'e $n$. 
Alors d'apr\`es \cite[ch.~III, \S4.4, th\'eo. 5]{Se}, pour $\delta\in H^\natural$, 
l'ensemble $\ES{O}_{\bf H}(\delta)(F)$ est r\'eunion d'un nombre {\it fini} de $H$--orbites.
\hfill $\blacksquare$
\end{enumerate}}
\end{mesrems}

\subsection{La topologie $\varpi$--adique (cas d'un corps local non archim\'edien)}\label{la topo p-adique} 
On suppose dans ce ${\rm n}^\circ$ que 
$F$ est un corps commutatif localement compact 
non archim\'edien. On note $\frak{o}_F$ l'anneau des entiers de 
$F$, et l'on choisit une uniformisante $\varpi$ de $F$. 

Si ${\bf X}$ est une vari\'et\'e alg\'ebrique affine d\'efinie sur $F$, on 
peut munir l'ensemble $X={\bf X}(F)$ de ses points $F$--rationnels de la 
{\it topologie d\'efinie par $F$}, appel\'ee aussi {\it topologie $\varpi$--adique}: c'est la topologie 
la {\it moins fine} rendant continues les applications $X\rightarrow F$ induites 
par les \'el\'ements de l'alg\`ebre affine $F[{\bf X}]$. Elle est {\it plus fine} 
que la topologie de Zariski restreinte \`a $X$. Cela fait de $X$ un td--espace. 
Si de plus ${\bf X}$ est lisse, 
alors $X$ est une {\it vari\'et\'e (analytique) $\varpi$--adique} au sens de \cite[ch.~V, \S2]{HC1}, et 
sa dimension comme vari\'et\'e $\varpi$--adique, not\'ee 
$\dim(X)$, co\"{\i}ncide avec $\dim({\bf X})$.

Munissons les ensembles de points $F$--rationnels $H={\bf H}(F)$ et 
$H^\natural={\bf H}^\natural(F)$ de la topologie $\varpi$--adique. 
Cela fait de $H$ un groupe topologique localement profini, et 
de $H^\natural$ un td--espace et un $H$--espace topologique tordu. Les td--espaces $H$ et 
$H^\natural$ sont des vari\'et\'es $\varpi$--adiques, de 
m\^eme dimension $\dim({\bf H})=\dim({\bf H}^\natural)$, et pour 
$\delta\in H^\natural$, le td--espace $\ES{O}_{\bf H}(\delta)(F)$ est une 
vari\'et\'e $\varpi$--adique de dimension 
$\dim({\bf H})-\dim({\bf H}_\delta)$. Rappelons que pour $\delta\in H^\natural$, on a 
pos\'e
$${\bf H}_\delta^{\rm sep} ={\bf H}_\delta(F^{\rm sep}){\bf H}_\delta^\circ$$
et 
$$H_\delta=\{h\in H: {\rm Int}_{H^\natural}(\delta)(h)=h\}.$$
Pour $\delta\in H^\natural$, on a
$$\dim({\bf H}_\delta)=\dim({\bf H}_\delta^\circ)=\dim({\bf H}_\delta^{\rm sep}),
$$ 
et si ${\bf H}_\delta^\circ$ est d\'efini sur $F$, i.e. si ${\bf H}_\delta^{\rm sep}$ est d\'efini sur $F$, 
on a $H_\delta={\bf H}_\delta^{\rm sep}(F)$.

Si ${\bf H}=\Bbb{GL}_n$, l'ensemble 
$\{1+\varpi^kM(n,\frak{o}_F):k\in \Bbb{Z}_{\geq 1}\}$ 
est une base de voisinages de $1$ dans ${\rm GL}_n(F)$ form\'ee de 
sous--groupes ouverts compacts. Dans le cas g\'en\'eral, la 
topologie $\varpi$--adique sur $H$ co\"{\i}ncide avec celle d\'eduite de ${\rm GL}_n(F)$ par restriction, via 
le choix de n'importe quel $F$--plongement ${\bf H}\hookrightarrow \Bbb{GL}_n$. 
En particulier, $H$ est r\'eunion d\'enombrable d'ouverts compacts 
(cela r\'esulte par exemple de 
la d\'ecomposition de Cartan pour 
${\rm GL}_n(F)$).

\begin{mapropo1}
Pour $\delta\in H^\natural$ quasi--semisimple, la $H$--orbite 
$\ES{O}_H(\delta)$ est ferm\'ee dans $H$ (pour la topologie $\varpi$--adique), et l'application bijective 
$H_\delta\backslash H \rightarrow \ES{O}_H(\delta),\,h\mapsto h^{-1}\cdot \delta \cdot h$ est un hom\'eomorphisme.
\end{mapropo1}

\begin{proof}
Soit $\delta\in H^\natural$ quasi--semisimple. D'apr\`es le thorme de \ref{automorphismes stabilisant (B,T)}, le groupe 
${\bf H}_\delta^\circ$ est d\'efini sur $F$, et $H_\delta={\bf H}_\delta^{\rm sep}(F)$ est une vari\'et\'e 
$\varpi$--adique de dimension $\dim({\bf H}_\delta^\circ)$. Par suite la vari\'et\'e $\varpi$--adique 
$\ES{O}_{\bf H}(\delta)(F)$ a pour dimension $\dim(H)-\dim(H_\delta)$.

La premire assertion est une consquences de rsultats rappels dans l'Annexe C. 
Supposons le morphisme $\pi_\delta$ s\'eparable. En ce cas ${\bf H}_\delta$ est d\'efini sur $F$, 
et la $H$--orbite $\ES{O}_H(\delta)$ est une 
sous--vari\'et\'e $\varpi$--adique de $\ES{O}_{\bf H}(\delta)(F)$, de m\^eme dimension $\dim(H)-\dim(H_\delta)$. 
D'aprs le lemme 1 de \ref{cas particuliers}, elle est ouverte dans $\ES{O}_{\bf H}(\delta)(F)$, et toutes les $H$--orbites dans $\ES{O}_{\bf H}(\delta)(F)$ 
sont ouvertes et ferm\'ees dans $\ES{O}_{\bf H}(\delta)(F)$.

Supposons maintenant que le morphisme $\pi_\delta$ n'est pas s\'eparable. On a donc $p>1$, et $\pi_\delta$ 
se d\'ecompose en
$$\pi_\delta =\bar{\pi}'_\delta\circ q: {\bf H} \buildrel q\over{\longrightarrow} {\bf H}_\delta^{\rm sep}\backslash {\bf H} 
\buildrel \bar{\pi}'_\delta\over{\longrightarrow}\ES{O}_{\bf H}(\delta)$$ o\`u 
$q$ est le morphisme quotient --- il est s\'eparable, et d\'efini sur $F$ puisque ${\bf H}_\delta^{\rm sep}$ l'est --- 
et $\bar{\pi}'_\delta$ est un morphisme non s\'eparable d\'efini sur $F$. Ces morphismes 
sont surjectifs, et m\^eme dominants puisque toutes les vari\'et\'es sont irr\'eductibles. D'aprs le lemme de \ref{un rsultat bien connu}, il 
existe un ouvert affine non vide ${\bf U}\subset \ES{O}_{\bf H}(\delta)$ tel que, posant ${\bf U}'=(\bar{\pi}'_\delta)^{-1}({\bf U})$, 
le morphisme
$$
\eta= \bar{\pi}'_\delta\vert_{{\bf U}'}: {\bf U}'\rightarrow {\bf U}
$$
est {\it fini}, i.e. le comorphisme (injectif)
$$
\eta^\sharp: \overline{F}[{\bf U}]\rightarrow \overline{F}[{\bf U}']
$$
fait de la $\overline{F}$--algbre $\overline{F}[{\bf U}']$ un $\overline{F}[{\bf U}]$--module de type fini. Par homognit on en dduit que le 
morphisme $\bar{\pi}'_\delta$ lui--m\^eme est fini. Puisque le morphisme $q$ est sparable, l'application 
$q_F: H\rightarrow ({\bf H}_\delta^{\rm sep}\backslash{\bf H})(F)$ est ouverte (lemme 1 de \ref{cas particuliers}), et la $H$--orbite 
$q_F(H)$ est ouverte et ferme dans $({\bf H}_\delta^{\rm sep}\backslash{\bf H})(F)$. D'autre part puisque le morphisme $\bar{\pi}'_\delta$ 
est fini, l'application $(\bar{\pi}'_\delta)_F:({\bf H}_\delta^{\rm sep}\backslash{\bf H})(F)\rightarrow \ES{O}_{\bf H}(\delta)(F)$ est ferme 
(remarque (3) de \ref{cas particuliers}), et la $H$--orbite
$$
(\bar{\pi}'_\delta)_F(q_F(H))= \ES{O}_H(\delta)
$$
est ferme dans $\ES{O}_{\bf H}(\delta)(F)$. D'o\`u la premi\`ere assertion 
de la proposition, puisque d'apr\`es la proposition de \ref{orbites}, 
$\ES{O}_{\bf H}(\delta)(F)$ est une partie ferm\'ee 
de $H^\natural$. 

La seconde assertion r\'esulte de \cite[ch.~1, cor.~1.6]{BZ}.
\end{proof}

\begin{marema}
{\rm Soit $L=\overline{F}(\ES{O}_{\bf H}(\delta))$ et $M=\overline{F}({\bf H}_\delta^{\rm sep}\backslash {\bf H})$ 
les corps des fonctions rationnelles sur $\ES{O}_{\bf H}(\delta)$ et sur ${\bf H}_\delta^{\rm sep}\backslash {\bf H}$. 
Le comorphisme
$$
(\bar{\pi}'_\delta)^\sharp: \overline{F}[\ES{O}_{\bf H}(\delta)]\rightarrow \overline{F}[{\bf H}_\delta^{\rm sep}\backslash {\bf H}]
$$
induit par passage aux corps des fractions un morphisme injectif de corps, disons $\iota: L \rightarrow M$, qui fait de $M$ une extension 
finie de $L$. Soit $L'/L$ 
la sous--extension s\'eparable maximale de $M/L$. Son degr\'e $m$ co\"{\i}ncide avec le cardinal du noyau de $\bar{\pi}'_\delta$, 
c'est--\`a--dire avec le cardinal de ${\bf H}_\delta^{\rm sep}\backslash {\bf H}_\delta$, lequel est premier \`a $p$ 
(remarque 1 de \ref{automorphismes qss}). D'autre part si le morphisme $\pi_\delta$ n'est pas sparable, alors 
l'extension $L/L'$ est purement ins\'eparable de degr\'e 
$q=p^s$, $s\geq 1$.\hfill $\blacksquare$}
\end{marema}

Si $p=1$, alors $F$ est un corps \og de type (F) \fg{} au sens de \cite[ch.~III, \S4.2]{Se} (cf. la remarque (3) de \ref{H(F)-orbites}) 
et pour $\delta\in H^\natural$, 
l'ensemble $\ES{O}_{\bf H}(\delta)(F)$ est r\'eunion finie de $H$--orbites 
(loc.~cit.). En g\'en\'eral, on a le r\'esultat plus faible suivant:

\begin{mapropo2}
Soit $\delta\in H^\natural$. Supposons que le 
morphisme $\pi_\delta$ est s\'eparable, que le 
groupe ${\bf H}_\delta^\circ$ est 
r\'eductif, et que le groupe quotient ${\bf H}_\delta/{\bf H}_\delta^\circ$ est d'ordre 
premier \`a $p$. Alors 
l'ensemble $\ES{O}_{\bf H}(\delta)(F)$ est r\'eunion {\rm finie} de $H$--orbites.
\end{mapropo2}

\begin{proof}
On peut supposer $p>1$. Alors $F$ est isomorphe \`a un corps de 
s\'eries formelles $\Bbb{F}_q((\varpi))$ o\`u $\Bbb{F}_q$ d\'esigne le corps fini \`a $q\,(=p^r)$ 
\'el\'ements. Puisque le morphisme $\pi_\delta$ est s\'eparable, le groupe 
${\bf H}_\delta$ est d\'efini sur $F$, le quotient ${\bf H}/{\bf H}_\delta$ l'est aussi, et le 
morphisme ${\bf H} \rightarrow \ES{O}_{\bf H}(\delta),\, h\mapsto h\cdot \delta \cdot h^{-1}$ induit une 
application bijective $({\bf H}/{\bf H}_\delta)(F) \rightarrow \ES{O}_{\bf H}(\delta)(F)$. 
Pour toute sous--extension $F'/F$ de $F^{\rm sep}\!/F$ et toute vari\'et\'e alg\'ebrique affine 
${\bf V}$ d\'efini sur $F'$, on note ${\rm H}^1(F',{\bf V})$ l'ensemble point\'e 
${\rm H}^1(\Sigma(F^{\rm sep}\!/F'),{\bf V}(F^{\rm sep}))$ d\'efini dans \cite[ch.~III, \S1]{Se}.   
D'apr\`es \cite[ch.~I, \S5.4, cor.~1]{Se}, le quotient de $({\bf H}/{\bf H}_\delta)(F)$ par $H$ s'identifie au 
noyau (dans la cat\'egorie des ensembles point\'es) de l'application canonique 
${\rm H}^1(F,{\bf H}_\delta)\rightarrow {\rm H}^1(F,{\bf H})$. Il suffit donc de montrer que l'ensemble 
${\rm H}^1(F,{\bf H}_\delta)$ est fini. Puisque ${\bf H}_\delta^\circ$ est d\'efini sur 
$F$ et distingu\'e dans ${\bf H}_\delta$, on a 
la suite exacte longue d'ensembles point\'es \cite[ch.~I, \S5.5, prop.~38]{Se}
$$
{\rm H}^1(F,{\bf H}_\delta^\circ)\rightarrow {\rm H}^1(F,{\bf H}_\delta)\rightarrow 
{\rm H}^1(F,{\bf H}_\delta/{\bf H}_\delta^\circ).
$$
Puisque ${\bf H}_\delta^\circ$ est r\'eductif connexe, d'apr\`es \cite[ch.~III, th\'eo.~3.12]{BT3} 
(cf. \cite[ch.~III, \S4.3, rem.~2]{Se}), 
l'ensemble ${\rm H}^1(F,{\bf H}_\delta^\circ)$ est fini. Il suffit donc de montrer 
que l'ensemble $H^1(F,{\bf H}_\delta/{\bf H}_\delta^\circ)$ est fini. 
Soit $F^{\rm mod}\!/F$ la sous--extension mod\'er\'ement ramifi\'ee maximale 
de $F^{\rm sep}\!/F$, et soit $\Sigma^{\rm mod}={\rm Gal}(F^{\rm mod}\!/F)$ son groupe de Galois. 
Puisque ${\bf H}_\delta$ est d\'efini sur $F$, le groupe quotient 
${\bf H}_\delta/{\bf H}_\delta^\circ$ l'est aussi, et comme les composantes connexes de 
${\bf H}_\delta$ sont d\'efinies sur $F^{\rm sep}$ \cite[ch.~AG, 12.3]{Bor}, on a ${\bf H}_\delta/{\bf H}_\delta^\circ=
({\bf H}_\delta/{\bf H}_\delta^\circ)(F^{\rm sep})$. Comme le groupe de Galois ${\rm Gal}(F^{\rm sep}\!/F^{\rm mod})$ 
est un pro-$p$-groupe et que (par hypoth\`ese) 
le groupe ${\bf H}_\delta/{\bf H}_\delta^\circ$ est d'ordre premier \`a $p$, 
on a ${\rm H}^1(F^{\rm mod},{\bf H}_\delta
/{\bf H}_\delta^\circ)=0$. 
D'apr\`es \cite[ch.~I, \S2.7.b]{Se}, on a donc une identification canonique
$$
{\rm H}^1(F,{\bf H}_\delta/{\bf H}_\delta^\circ)= 
{\rm H}^1(\Sigma^{\rm mod},({\bf H}_\delta/{\bf H}_\delta^\circ)(F^{\rm mod})).
$$
Or pour chaque entier $n\geq 1$, il n'existe qu'un nombre fini 
de sous--extensions de $F^{\rm mod}\!/F$ de degr\'e $n$. Par 
suite le groupe $\Sigma^{\rm mod}$ est \og de type (F)\fg{} au sens de \cite[ch.~III, \S4.1]{Se}, et 
l'ensemble ${\rm H}^1(\Sigma^{\rm mod},({\bf H}_\delta/{\bf H}_\delta^\circ)(F^{\rm mod}))$ 
est fini \cite[ch.~III, \S4.1, prop.~8]{Se}. D'o\`u la proposition.
\end{proof}

D'apr\`es le thorme et la remarque 1 de \ref{automorphismes qss}, on a le

\begin{moncoro}
Pour $\delta\in H^\natural$ quasi--semisimple tel que le morphisme $\pi_\delta$ est s\'eparable 
(e.g. si ${\bf H}$ est semisimple et simplement connexe, d'apr\`es le corollaire de \ref{automorphismes qss}),  
l'ensemble $\ES{O}_{\bf H}(\delta)(F)$ est r\'eunion finie de $H$--orbites.
\end{moncoro}

\begin{exemples}
{\rm Rappelons que le groupe $F^\times/(F^\times)^2$ est fini si et seulement si $F$ est de caract\'eristique diff\'erente de $2$.
\begin{enumerate}
\item[(1)]Soit ${\bf H}=\Bbb{SL}_2$. Il existe une unique ${\bf H}$--orbite unipotente 
non triviale dans ${\bf H}$: celle de l'\'el\'ement $u=\left(\begin{array}{cc}1&1\\ 0&1\end{array}\right)\in H$. Pour 
$x\in F^\times$, notons $\ES{O}_x$ la $H$--orbite de l'\'el\'ement 
$u_x=\left(\begin{array}{cc}1&x\\ 0&1\end{array}\right)$. Alors $\ES{O}_x$ ne d\'epend que de l'image de $x$ 
dans $F^\times/(F^\times)^2$, et $\ES{O}_{\bf H}(u)(F)$ est l'union disjointe des $H$--orbites $\ES{O}_x$ pour $x\in F^\times/(F^\times)^2$.
\item[(2)]Soit ${\bf H}=\Bbb{GL}_2$, et soit $\gamma$ l'\'el\'ement $\left(\begin{array}{cc}0&1\\\varpi & 0\end{array}\right)$ de 
$H$. Son polyn\^ome caract\'eristique est $t^2+\varpi$. Il est irr\'eductible sur $F$ et s\'eparable 
si et seulement si la caract\'eristique de $F$ est diff\'erente de $2$, auquel cas $\gamma$ est semisimple r\'egulier. 
Si $F$ est de caract\'eristique $2$, alors $\gamma$ est conjugu\'e dans 
${\bf H}$ \`a l'\'el\'ement $\left(\begin{array}{cc}\varpi^{1\over 2}& 1\\ 0 & \varpi^{1\over 2}\end{array}\right)$, o\`u 
$\varpi^{1\over 2}$ est l'unique racine (double) de $t^2+\varpi$ dans $\overline{F}$.
\item[(3)]Soit $\tau$  le $F$--automorphisme $t\mapsto t^{-1}$ de ${\bf H}=\Bbb{G}_{\rm m}$, 
et posons ${\bf H}^\natural ={\bf H}\tau$. Pour $\delta\in H^\natural$ et $t\in {\bf H}$, on a 
$t^{-1}\cdot \delta \cdot t= t^{-2}\cdot \delta$, par cons\'equent 
$\ES{O}_{\bf H}(\delta)={\bf H}^\natural$ et $\ES{O}_{\bf H}(\delta)(F)$ est l'union disjointe 
des $H$--orbites $\ES{O}_H(x\cdot \delta)$ pour $x$ parcourant un syst\`eme de repr\'esentants 
dans $F^\times$ des classes de $F^\times/(F^\times)^2$. Notons que le morphisme $1-\tau$ 
de ${\bf H}$ est s\'eparable si et seulement si $F$ est de caract\'eristique diff\'erente de $2$.\hfill $\blacksquare$
\end{enumerate}}
\end{exemples}

\section{Caract\`eres tordus d'un groupe r\'eductif $\frak{p}$--adique}

Dans ce chapitre, on fixe un corps commutatif localement compact 
non archim\'edien $F$, et un groupe alg\'ebrique r\'eductif connexe 
${\bf G}$ d\'efini sur $F$. On note $G={\bf G}(F)$ le groupe des 
points $F$--rationnels de ${\bf G}$ {\it muni 
de la topologie $\varpi$--adique}, o $\varpi$ d\'esigne une uniformisante 
de $F$. Pour la thorie de base des groupe rductifs $\varpi$--adiques, on renvoie  \cite{BT1, BT2}. 
On fixe un ${\bf G}$--espace tordu 
${\bf G}^\natural$ d\'efini sur $F$ et poss\'edant un point $F$--rationnel 
$\delta_1$. On note $G^\natural$ le $G$--espace tordu 
${\bf G}^\natural(F)$, et $\theta$ le $F$--automorphisme 
${\rm Int}_{{\bf G}^\natural}(\delta_1)$ de ${\bf G}$. On a donc 
$G^\natural=G\cdot \delta_1\subset {\bf G}^\natural$. 
On fixe aussi un caract\`ere $\omega$ de $G$.

\subsection{Paires paraboliques de $G$}\label{paires paraboliques}
On appelle {\it paire parabolique de $G$} une paire $(P,A)$ form\'ee 
d'un sous--groupe parabolique $P$ de $G$ et 
d'un tore d\'eploy\'e maximal $A$ du radical $R(P)$ de $P$.

Si $P$ est un sous--groupe parabolique de $G$, on note $U_P=R_{\rm u}(P)$ 
son radical unipotent. Si $(P,A)$ 
est une paire parabolique de 
$G$, on note $M_A=Z_G(A)$ le centralisateur de $A$ dans $G$. Alors 
$M_A$ est une {\it composante de Levi de $P$}, 
i.e. on a la d\'ecomposition en produit 
semidirect
$$P=M_A\ltimes  U_P.
$$ De plus, la paire $(P,A)$ est $\theta$--stable 
si et seulement si les groupes $U_P$ et $M_A$ sont 
$\theta$--stables.

Fixons une paire parabolique minimale $(P_\circ, A_\circ)$ de $G$. 
On ne suppose pas qu'il existe une telle paire qui soit $\theta$--stable\footnote{M\^eme si l'on 
peut toujours s'arranger pour que ce soit le cas, cf. la remarque plus loin.}. 
Posons $U_\circ =R_{\rm u}(P_\circ)$ et $M_\circ=M_{A_\circ}$, et notons 
$\ES{P}_\circ$ l'ensemble des sous--groupes paraboliques de $G$ contenant 
$P_\circ$. Pour $P\in \ES{P}_\circ$, on note:
\begin{itemize}
\item  $M_P$ l'unique composante de Levi de $P$ contenant $M_\circ$; 
\item $P^-$ le sous--groupe parabolique de $G$ oppos\'e \`a $P$ par rapport \`a $M_P$;
\item $A_P$ le tore d\'eploy\'e maximal du centre $M_P$.
\end{itemize}
On a donc
$$
M_P=M_{A_P}=P\cap P^-\quad (P\in \ES{P}_\circ).
$$ 
De plus, l'application $P\mapsto (P,A_P)$ identifie $\ES{P}_\circ$ \`a l'ensemble des paires 
paraboliques $(P,A)$ de $G$ telles que $P\supset P_\circ$ et $A\subset A_\circ$, et  
$\ES{P}_\circ$ param\'etrise l'ensemble des classes de conjugaison de 
paires paraboliques de $G$.

Notons $\ES{P}_\circ^\theta$ le sous--ensemble de 
$\ES{P}_\circ$ form\'e des paires paraboliques $\theta$--stables, i.e. posons
$$
\ES{P}_\circ^\theta = \{P\in \ES{P}_\circ: \hbox{$\theta(P)=P$ et $\theta(A_P)=A_P$}\}.
$$
Notons que $\ES{P}_\circ^\theta$ est non vide (la paire parabolique maximale $(G,A_G)$ est 
$\theta$--stable). Notons aussi que --- contrairement \`a ce que la notation pourrait faire croire --- 
si la paire $(P_\circ,A_\circ)$ n'est pas $\theta$--stable, alors:
\begin{itemize}
\item $\theta$ {\it n'op\`ere pas} sur $\ES{P}_\circ$ (vu comme ensemble de paires 
paraboliques);
\item $\ES{P}_\circ^\theta$ {\it ne param\'etrise pas} l'ensemble des classes de conjugaison de 
paires paraboliques $\theta$--stables de $G$.
\end{itemize}
%

\begin{marema}
{\rm 
Puisque $(\theta(P_\circ),\theta(A_\circ))$ 
est encore une paire parabolique minimale de $G$, il existe un $x\in G$ 
tel que $\theta(P_\circ)=x^{-1}P_\circ x$ et $\theta(A_\circ)=x^{-1}A_\circ x$. Par 
suite, posant $\delta'_1=x\cdot \delta_1$ et $\theta'={\rm Int}_{\bf G}(\delta'_1)\in {\rm Aut}_F({\bf G})$, 
on a $\theta'(P_\circ)=P_\circ$ et $\theta'(A_\circ)=A_\circ$, i.e. 
la paire $(P_\circ,A_\circ)$ est $\theta'$--stable.
}
\end{marema}

\begin{monlem}
Supposons que la paire $(P_\circ,A_\circ)$ est $\theta$--stable. 
Alors $\theta$ op\`ere sur l'ensemble $\ES{P}_\circ$ de mani\`ere compatible avec l'identification 
\og $P=(P,A_P)$ \fg{}: pour $P\in \ES{P}_\circ$, on a $A_{\theta(P)}=\theta(A_P)$. De plus, 
$\ES{P}_\circ^\theta$ 
param\'etrise l'ensemble des classes de 
$G$-conjugaison de paires paraboliques $\theta$--stables de $G$.
\end{monlem}

\begin{proof}
Pour $P\in \ES{P}_\circ$, puisque $\theta(P_\circ)=P_\circ$, on a 
$\theta(P)\in \ES{P}_\circ$, et puisque $\theta(A_\circ)=A_\circ$, on a 
$\theta(A_P)\subset A_\circ$ et $\theta(A_P)= A_{\theta(P)}$. Par suite, 
$\theta$ op\`ere sur 
$\ES{P}_\circ$ de mani\`ere compatible avec l'identification \og $P=(P,A_P)$\fg. 
En particulier, on a  
$\ES{P}_\circ^\theta=\{P\in \ES{P}_\circ: \theta(P)=P\}$. 

Soit $(P,A)$ une paire parabolique $\theta$--stable de $G$. Il existe un $g\in G$ tel que $g^{-1}Pg\supset P_\circ$ et $g^{-1}Ag \subset A_\circ$, i.e. 
tel que $g^{-1}Pg\in \ES{P}_\circ$ et $A_{g^{-1}Pg}= g^{-1}Ag$. Posons 
$x=g^{-1}\theta(g)$. 
Puisque
$$P_\circ \subset \theta(g^{-1}Pg)=
\theta(g^{-1})P\theta(g)= x^{-1}(g^{-1}Pg)x,
$$
on a 
$\theta(g^{-1}Pg)=g^{-1}Pg$. Par suite, $\ES{P}_\circ^\theta$ 
param\'etrise l'ensemble des classes de 
$G$-conjugaison de paires paraboliques $\theta$--stables de $G$.
\end{proof}

\subsection{Mesures normalises}\label{mesures normalises}
Soit $\vert\;\vert_F$ la valeur absolue sur $F$ normalis\'ee par 
$$\vert \varpi\vert_F =q^{-1},
$$
o $q$ est le cardidal du corps r\'esiduel de $F$. 
Soit $K_\circ$ le stabilisateur dans 
$$G^1=\bigcap_{\psi\in {\rm X}^*_F({\bf G})}\ker \vert \psi\vert_F
$$ 
d'un sommet sp\'ecial de l'appartement 
$\ES{A}_\circ=\ES{A}(G,A_\circ)$ associ\'e \`a $A_\circ$ de l'immeuble (non \'etendu) de $G$; o ${\rm X}^*_F({\bf G})$ est 
le groupe des caractres algbriques de ${\bf G}$ dfinis sur $F$. 
Ainsi $K_\circ$ est un sous--groupe ouvert compact 
maximal (sp\'ecial) de $G$, et pour tout $P\in \ES{P}_\circ$, on a:
\begin{enumerate}
\item[(i)]$G=K_\circ P =K_\circ M_P U_P$;
\item[(ii)]$K_\circ\cap P =(K_\circ \cap M_P)(K_\circ \cap U_P)$;
\item[(iii)]$K_\circ\cap P^- =(K_\circ \cap M_P)(K_\circ \cap U_{P^-})$.
\end{enumerate}
%

\begin{mesrems}
{\rm 
\begin{enumerate}
\item[(1)]Pour $i\in \Bbb{Z}$, $\theta^i(K_\circ)$ est 
encore un sous--groupe ouvert compact maximal sp\'ecial de 
$G$, et si la paire $(P_\circ,A_\circ)$ est $\theta$--stable, alors 
$\theta^i(K_\circ)$ v\'erifie les propri\'et\'es (i), (ii), (iii) ci-dessus 
(pour tout $P\in \ES{P}_\circ$). Notons qu'il n'est en g\'en\'eral 
pas possible de choisir $\theta'$ et $K_\circ$ tels que $(P_\circ,A_\circ)$ et 
$K_\circ$ soient $\theta'$--stables, m\^eme si 
$G$ est {\it non ramifi\'e} (c'est--\`a--dire quasi--d\'eploy\'e sur $F$ 
et d\'eploy\'e sur une extension non ramifi\'ee de $F$), cf. la remarque suivante.
\item[(2)]Prenons pour 
${\bf G}$ le groupe $\Bbb{SL}_{2/F}$, et pour $\theta$ l'automorphisme donn\'e par la conjugaison 
par un \'el\'ement de ${\rm GL}_2(F)$ dont le d\'eterminant est une uniformisante 
de $F$. Le groupe ${\bf G}$ est d\'eploy\'e sur $F$ --- donc a fortiori non ramifi\'e ---, et 
pour tout $x\in {\bf G}(F)$, le $F$--automorphisme 
${\rm Int}_{\bf G}(x)\circ \theta$ de ${\bf G}$ 
change le type des sous--groupes hypersp\'eciaux de ${\bf G}(F)$. Il ne peut donc en 
stabiliser aucun.\hfill $\blacksquare$
\end{enumerate}
}\end{mesrems}

\begin{madefi}
{\rm Pour un sous--groupe ferm\'e $H$ de $G$, on appelle 
{\it mesure de Haar \`a gauche sur $H$ normalis\'ee par $K_\circ$} 
l'unique mesure de Haar \`a gauche $d_lh$ sur $H$ telle que 
$${\rm vol}(H\cap K_\circ,d_lh)=1.$$}
\end{madefi}

Soit $dg$ la mesure de Haar sur $G$ normalis\'ee par 
$K_\circ$. Si $K$ est un 
sous--groupe ouvert compact de $G$, on note aussi 
$dk$ la mesure $dg\vert_K$ sur $K$.

Pour $P\in \ES{P}_\circ$, on note $dm_P$ et 
$du_P$ les mesures de Haar sur $M_P$ et $U_P$ normalis\'ees par 
$K_\circ$, et l'on pose 
$d_lp_P=dm_Pdu_P$; c'est la mesure de Haar \`a gauche sur $P$ 
normalis\'ee par $K_\circ$. Pour $x\in P$, on 
a donc (abus d'\'ecriture)
$$\Delta_P(x)d_l(xpx^{-1})_P=d_lp_P.
$$ Gr\^ace \`a $dm_P$, 
on d\'efinit les caract\`eres de $M_P$ comme en \ref{caractres}. Lorsqu'il n'y aura pas de confusion possible, on omettra 
l'indice $P$ dans les notations $dm_P$, $du_P$ et $d_lp_P$. Pour toute 
repr\'esentation lisse $\sigma$ de $M_P$, on pose
$$\sigma(f)=\sigma(fdm_P)\quad (f\in C^\infty_{\rm c}(M_P)),
$$ et si $\sigma$ est admissible, on note $\Theta_\sigma$ et 
$\Theta_\sigma^B$ (pour $B\in {\rm End}_{\Bbb C}(V)$) le caract\`ere 
et le caract\`ere $B$--tordu de $\sigma$ d\'efinis gr\^ace \`a $dm_P$ 
comme en \ref{caractres} et \ref{caractres tordus}.

\subsection{sous--espaces paraboliques de $G^\natural$}\label{sous-espaces paraboliques}
On appelle {\it sous--espace 
parabolique de $G^\natural$} un sous--espace topologique tordu de $G^\natural$ (cf. \ref{espaces topo tordus}) qui 
est un $P$--espace tordu pour un sous--groupe parabolique $P$ de $G$; i.e. un 
sous--espace topologique de $G^\natural$ de la forme $P\cdot\gamma$ pour un sous--groupe 
parabolique $P$ de $G$ et un \'el\'ement $\gamma\in G^\natural$ tel que 
${\rm Int}_{G^\natural}(\gamma)(P)=P$. Si $P$ est un sous--groupe parabolique de $G$, 
puisque $N_G(P)=P$, il existe {\it au plus} un sous--espace parabolique de $G$ de 
$G^\natural$ qui est un $P$--espace tordu. En d'autres termes, l'application 
$P^\natural \mapsto P$ de l'ensemble des sous--espaces parabolique de 
$G^\natural$ dans l'ensemble des sous--groupes paraboliques de $G$, est injective. 

Notons que pour tout sous--espace parabolique $P^\natural= P\cdot \gamma$ de 
$G^\natural$, notant $N_G(P^\natural)$ 
le normalisateur $\{g\in G: g^{-1}\cdot P^\natural\cdot g=P^\natural\}$ 
de $P^\natural$ dans $G^\natural$, on a
$$N_G(P^\natural)=N_G(P)=P.
$$ 

Si $P^\natural$ est un sous--espace parabolique de $G^\natural$, on appelle 
{\it composante de Levi de $P^\natural$} un sous--espace topologique tordu de 
$P^\natural$ qui est un $M$--espace tordu pour une composante de Levi de $M$ de $P$; 
i.e. un sous--espace topologique $M^\natural$ de $P^\natural$ de la forme 
$M\cdot \gamma$ pour une composante de Levi $M$ de $P$ et un \'el\'ement $\gamma\in P^\natural$ tel que ${\rm Int}_{G^\natural}(\gamma)(M)=M$.

\begin{monlem1}
Soit $P^\natural$ un sous--espace parabolique de $G^\natural$. 
Il existe une composante de Levi $M^\natural$ de $P^\natural$. L'application $M^\natural\mapsto M$ de 
l'ensemble des composantes de Levi de $P^\natural$ dans l'ensemble 
des composantes de Levi de $P$, est bijective. En particulier si $M^\natural$ et 
$M'^\natural$ sont deux composantes de Levi de $P^\natural$, 
alors il existe un unique $u\in U_P$ tel que $M'^\natural=u\cdot M^\natural
\cdot u^{-1}$.
\end{monlem1}

\begin{proof}
\'Ecrivons $P^\natural=P\cdot \gamma$, et soit  
$\tau={\rm Int}_{G^\natural}(\gamma)$. Puisque $\tau(P)=P$, on a $\tau(U_P)=U_P$. 
Soit $M$ une composante de Levi de $P$. Alors $\tau(M)$ est une autre composante de Levi 
de $M$, et il existe un unique $u\in U_P$ tel que $\tau(M)=u^{-1}Mu$. Posons $\gamma'=u\cdot \gamma\in P^\natural$ et $\tau'={\rm Int}_{G^\natural}(\gamma')$. Alors $\tau'(M)=M$ et $M\cdot \gamma'$ est 
une composante de Levi de $P^\natural$.

D'apr\`es ce qui pr\'ec\`ede, l'application $M^\natural\mapsto M$ 
de l'ensemble des composantes de Levi de $P^\natural$ dans l'ensemble 
des composantes de Levi $P$, est surjective. Elle est aussi injective car 
pour toute composante de Levi $M$ de $P$, on a $N_G(M)\cap P= M$. 
D'o\`u la derni\`ere assertion du lemme.
\end{proof}

Soit $P^\natural$ un sous--espace parabolique de $G^\natural$. 
Puisque
$$
{\rm Int}_{P^\natural}(\gamma)(U_P)=U_P\quad (\gamma\in P^\natural),
$$ 
on peut consid\'erer l'espace topologique tordu quotient $P^\natural/U_P$ (cf. \ref{espaces topo tordus}). 
Pour toute composante de Levi $M^\natural$ de $P^\natural$, on a
la {\it d\'ecomposition de Levi} 
$$
P^\natural = M^\natural\cdot U_P;
$$
i.e. tout \'el\'ement $\gamma\in P^\natural$ se d\'ecompose de 
mani\`ere unique en $\gamma=\delta\cdot u$ o $\delta\in M^\natural$, $u\in U_P$ et 
${\rm Int}_{P^\natural}(\delta)(U_P)=U_P$. De plus, l'inclusion $M^\natural\subset P^\natural$ induit 
un isomorphisme d'espaces topologiques tordus $M^\natural \rightarrow P^\natural/U_P$.

Soit $\ES{P}_\circ^\natural$ l'ensemble des sous--espaces paraboliques $P^\natural$ de $G^\natural$ tel 
que $P=N_G(P^\natural)$ est un lment de $\ES{P}_\circ$. L'application 
$$
\ES{P}_\circ^\natural\rightarrow \ES{P}_\circ,\, P^\natural\mapsto P=N_G(P^\natural)
$$
est injective, d'image le sous--ensemble de $\ES{P}_\circ$, disons $\ES{P}_\circ^{(\natural)}$, form des $P\in \ES{P}_\circ$ 
tels qu'il existe un sous--espace parabolique $P^\natural$ de $G^\natural$ qui est un $P$--espace tordu. 
D'aprs le lemme, pour chaque 
$P^\natural\in \ES{P}_\circ^\natural$, il existe une unique composante 
de Levi de $P^\natural$ qui est un $M_P$--espace tordu; on la note $M_P^\natural$. 
On en d\'eduit que $\ES{P}_\circ^\natural$ param\'etrise:
\begin{itemize}
\item l'ensemble des classes de 
$G$--conjugaison de sous--espaces paraboliques de $G^\natural$;
\item l'ensemble des classes de 
$G$--conjugaison de sous--espaces paraboliques $P^\natural$ de $G^\natural$ munis d'une 
d\'ecomposition de Levi $P^\natural=M^\natural_{P}\cdot U_P$.
\end{itemize}
\ni Pour $P\in \ES{P}_\circ^\theta$, l'ensemble $P\cdot\delta_1$ est un 
sous--espace parabolique de $G^\natural$. D'o\`u l'inclusion 
$$
\ES{P}_\circ^\theta\subset \ES{P}_\circ^{(\natural)}.
$$

\begin{monlem2}
Supposons que la paire $(P_\circ,A_\circ)$ est 
$\theta$--stable. Alors on a l'\'egalit\'e $\ES{P}_\circ^\theta=\ES{P}_\circ^{(\natural)}$, et pour 
$P^\natural \in \ES{P}_\circ^{\natural}$ de groupe sous--jacent $P=N_G(P^\natural)$, on a 
$P^\natural = P\cdot \delta_1$ et $M_P^\natural = M_P\cdot \delta_1$.
\end{monlem2}

\begin{proof}
Il s'agit de montrer que $\ES{P}_\circ^{(\natural)}\subset \ES{P}_\circ^\theta$. 
Soit $P^\natural\in \ES{P}_\circ^\natural$. \'Ecrivons $P^\natural=P\cdot \gamma$,  
et $\gamma= g\cdot \delta_1$ avec $g\in G$. Soit $\tau={\rm Int}_{G^\natural}(\gamma)$. 
Puisque $\tau(P)=P$ et $\tau={\rm Int}_G(g)\circ\theta$, on a 
$\theta(P)={\rm Int}_G(g^{-1})(P)$. Or $\theta(P)\in \ES{P}_\circ$, d'o\`u 
$\theta(P)=P$ et $g\in N_G(P)=P$.

Soit $P^\natural\in \ES{P}_\circ^\natural$ de groupe sous--jacent $P=N_G(P^\natural)$. On vient de montrer que 
$P^\natural= P\cdot \delta_1$. D'après le lemme de \ref{paires paraboliques}, on a $\theta(A_P)= A_P$ d'où $\theta(M_P)=M_P$, 
et donc (d'après le lemme 1) $M_P^\natural = M_P\cdot \delta_1$.
\end{proof}

Pour $P^\natural \in \ES{P}_\circ^\natural$, soit $d\gamma_{M_P}=\delta \cdot dm_P$ ($\delta\in M_P^\natural$) 
la mesure de Haar \`a gauche 
sur $M_P^\natural$ associ\'ee \`a $dm_P$ (c'est aussi une mesure de Haar \`a 
droite, cf. \ref{module d'un espace tordu}), et soit $d_l\gamma_P=d\gamma_{M_P}\cdot du_P$ la mesure 
produit sur $P^\natural=M_P^\natural\cdot U_P$. Alors $d_l\gamma_P$ est 
la mesure de Haar \`a gauche $\delta\cdot d_lp_P$ ($\delta\in P^\natural$) associ\'ee \`a $d_lp_P$. 
Pour toute 
$\omega$--repr\'esentation lisse $\Sigma$ de $M_P^\natural$, 
on pose
$$
\Sigma(\phi)=\Sigma(\phi\, d\gamma_{M_P})\quad (\phi\in C_{\rm c}^\infty(M_P^\natural)),
$$
et si $\Sigma$ est admissible, on note $\Theta_\Sigma$ le caract\`ere de $\Sigma$ 
d\'efini gr\^ace \`a $d\gamma_{M_P}$ comme en \ref{caractres tordus}.

Pour $P=G$, on pose $d\gamma=d\gamma_G$.

\v1 Soit $P^\natural$ un sous--espace parabolique de $G^\natural$. 
D'apr\`es la relation $(**)$ de \ref{module d'un espace tordu}, le module $\Delta_{P^\natural}$ de $P^\natural$  
se factorise \`a travers $P^\natural/U_P$. Soit $\delta_{P^\natural}: P^\natural\rightarrow \Bbb{R}_{>0}$ 
l'application d\'efinie par
$$
\delta_{P^\natural}(\gamma)= \Delta_{P^\natural}(\gamma)^{-1}
\quad (\gamma\in P^\natural).
$$
Elle aussi se factorise \`a travers $P^\natural/U_P$, et d'apr\`es le lemme de \ref{module d'un espace tordu}, on a
$$
\delta_{P^\natural}(p\cdot\gamma)= \delta_P(p)\delta_{P^\natural}(\gamma)\quad 
(p\in P,\,\gamma\in P^\natural),\leqno{(*)}
$$
o\`u $\delta_P$ est le caract\`ere $\Delta_P^{-1}$ de $P$.  
Pour $r\in \Bbb{R}$, on note $\delta_{P^\natural}^r:P^\natural\rightarrow \Bbb{R}_{>0}$ l'application dfinie par 
$$\delta_{P^\natural}^r(\gamma)= \delta_{P^\natural}(\gamma)^r.
$$

\begin{marema}
{\rm Puisque le groupe $G$ est unimodulaire, d'après le lemme de \ref{module d'un espace tordu}, le module $\Delta_{G^\natural}:G^\natural\rightarrow {\Bbb R}_{>0}$ est constant. On verra en \ref{la condition P2} (Annexe A) que l'on peut choisir le point--base $\delta_1\in G^\natural$ de telle manière qu'il existe une base de voisinages de $1$ dans $G$ formée de sous--groupes ouverts compacts de $G$ normalisés par $\delta_1$ (i.e. $\theta$--stables). Cela implique en particulier (d'après le lemme 1 de \ref{modules}) que $\Delta_{G^\natural}=1$. On a donc $\delta_{G^\natural}=1$. De la même manière, pour un sous--groupe parabolique $P^\natural$ de $G^\natural$, 
puisque le groupe quotient $P/U_P$ est unimodulaire, le module $\Delta_{P^\natural/U_P}$ de $P^\natural/U_P$ est constant, et l'on a 
$\Delta_{P^\natural/U_P}=1$. On a donc
$$
\delta_{P^\natural}(\gamma)= \Delta_{P^\natural}(\gamma)^{-1}\Delta_{P^\natural/U_P}(\gamma\cdot U_P)
\quad (\gamma\in P^\natural).
$$
}
\end{marema}

\subsection{\'El\'ements r\'eguliers et quasi--r\'eguliers de $G^\natural$}\label{lments rguliers et quasi-rguliers} 
On note $\frak{g}={\rm Lie}(G)$ l'alg\`ebre de Lie de $G$, et $\frak{g}^*={\rm Hom}_F(\frak{g},F)$ son  dual alg\'ebrique. Pour tout sous--groupe 
parabolique $P$ de $G$, on pose $\frak{p}={\rm Lie}(P)$ et 
$\frak{u}_P={\rm Lie}(U_P)$.

Soit ${\rm Ad}_{G^\natural}: G^\natural \rightarrow {\rm GL}(\frak{g})$ 
l'application $\gamma\mapsto {\rm Lie}({\rm Int}_{{\bf G}^\natural}
(\gamma)):\frak{g}\rightarrow \frak{g}$. 
Rappelons que pour $g\in G$ et $\gamma\in G^\natural$, on a ${\rm Int}_{G^\natural}(g\cdot\gamma)=
{\rm Int}_G(g)\circ {\rm Int}_{G^\natural}(\gamma)$. On a donc
$$
{\rm Ad}_{G^\natural}(g\cdot\gamma)={\rm Ad}_G(g)\circ {\rm Ad}_{G^\natural}
(\gamma)\quad (g\in G,\,\gamma\in G^\natural).
$$
Soit ${\rm Ad}_{G^\natural}^*:G^\natural\rightarrow 
{\rm Aut}_F(\frak{g}^*)$ l'application d\'efinie par
$$
\langle X, {\rm Ad}^*_{G^\natural}(\gamma)(Y)\rangle =
\langle {\rm Ad}_{G^\natural}(\gamma)^{-1}(X), Y\rangle\quad (X\in\frak{g},\,Y\in \frak{g}^*,\,\gamma\in G^\natural).
$$
Pour toute partie $\Omega$ de $\frak{g}$, on note $\Omega^\perp$ le sous--espace vectoriel 
de $\frak{g}^*$ d\'efini par
$$
\Omega^\perp=\{Y\in \frak{g}^*: \langle X,Y\rangle =0,\,\forall X\in \Omega\}.
$$ 
Pour $\gamma\in G^\natural$, on pose
\begin{align*}
& \frak{g}_\gamma = \ker ({\rm id}_\frak{g}-{\rm Ad}_{G^\natural}(\gamma))\subset \frak{g},\\
&\frak{g}(1-\gamma)= {\rm Im} ({\rm id}_\frak{g}-{\rm Ad}_{G^\natural}(\gamma))
\subset \frak{g},\\ 
&\frak{g}_\gamma^*= \ker ({\rm id}_{\frak{g}^*}-{\rm Ad}_{G^\natural}^*(\gamma)^{-1})\subset \frak{g}^*.
\end{align*}
On a donc
$$
\frak{g}_\gamma^*=\frak{g}(1-\gamma)^\perp\quad (\gamma\in G^\natural).
$$

\begin{madefi}
{\rm 
Comme dans \cite[appendix, prop.~A.2]{BH1}, un \'el\'ement $\gamma\in G^\natural$ 
est dit {\it quasi--r\'egulier} si pour tout sous--groupe 
parabolique $P$ de $G$, on a $\frak{g}(1-\gamma) + \frak{p} =\frak{g}$; i.e. si l'on a $\frak{g}_\gamma^*\cap \frak{p}^\perp=\{0\}$.
}
\end{madefi}

Rappelons que l'on a pos\'e $G^\natural_{\rm reg}=\{\gamma\in G^\natural: D_{G^\natural}(\gamma)\neq 0\}$. 
On note 
$G^\natural_{\rm qr}$ le sous--ensemble de $G^\natural$ form\'e 
des \'el\'ements quasi--r\'eguliers.  
Puisque l'ensemble $G_{\rm reg}^\natural$ est 
ouvert dense dans $G^\natural$ pour la topologie de Zariski, il l'est \`a fortiori 
pour la topologie $\varpi$--adique. Quant \`a l'ensemble $G^\natural_{\rm qr}$, 
il n'est pas d\'efini g\'eom\'etriquement, mais on verra plus loin (corollaire) 
qu'il poss\`ede lui 
aussi ces deux propri\'et\'es.

\begin{mapropo}
On a l'inclusion $G^\natural_{\rm reg}\subset G^\natural_{\rm qr}$.
\end{mapropo}

\begin{proof}
Soit $\gamma\in G_{\rm reg}^\natural$. D'apr\`es le thorme de \ref{automorphismes stabilisant (B,T)}, 
le tore ${\bf T}=Z_{\bf G}({\bf H}_\gamma^\circ)$ est d\'efini sur $F$. 
Posons $\frak{t}={\rm Lie}(T)$ et 
$\frak{t}^*={\rm Hom}_F(\frak{t},F)$. Comme en \ref{automorphismes rguliers; le cas i}, notons $\frak{g}_\gamma^1$ 
le sous--espace caract\'eristique de ${\rm Ad}_{G^\natural}(\gamma)$ associ\'e 
\`a la valeur propre $1$. L'application 
${\rm id}_\frak{g}-{\rm Ad}_{G^\natural}(\gamma)$ induit, par passage au quotient, un automorphisme 
du $F$--espace vectoriel $\frak{g}/\frak{g}_\gamma^1$. Puisque 
$\frak{g}_\gamma^1\subset \frak{t}$ (lemme 2 de \ref{automorphismes rguliers; le cas i}), on a l'\'egalit\'e 
$\frak{g}= \frak{t} + \frak{g}(1-\gamma)$. Pour $t\in T$, puisque 
$$ \frak{t}\cap {\rm Ad}_G(t)(\frak{g}(1-\gamma))=
\frak{t}\cap \frak{g}(1-\gamma),
$$
on a l'galit
$${\rm Ad}_G(t)(\frak{g}(1-\gamma))=\frak{g}(1-\gamma).
$$ 
L'inclusion $\frak{t}\subset \frak{g}$ induit donc, 
par passage aux quotients, une identification $T$--quivariante 
(pour l'action adjointe de $T$ sur $\frak{t}$ et $\frak{g}$):
$$\frak{t}/(\frak{t}\cap \frak{g}(1-\gamma))=
 \frak{g}/\frak{g}(1-\gamma).
 $$
 D'o\`u une identification $T$--quivariante (pour l'action coadjointe de $T$ sur $\frak{g}^*$ et 
 sur $\frak{t}^*$):
 $$
 \frak{g}_\gamma^* \;(={\rm Hom}_F(\frak{g}/\frak{g}(1-\gamma),F))
 ={\rm Hom}_F(\frak{t}/(\frak{t}\cap \frak{g}(1-\gamma)),F)\subset \frak{t}^*. 
 $$
On en d\'eduit que pour 
$Y\in \frak{g}_\gamma^*$, le stabilisateur de $Y$ dans 
 $G$ pour l'action coadjointe de $G$ sur $\frak{g}^*$, contient $T$. 
 Puisque $T$ est dense dans ${\bf T}$ pour la topologie de Zariski, 
 pour $Y\in \frak{g}^*_\gamma$, le stabilisateur de $Y$ dans 
 ${\bf G}$ pour l'action coadjointe de ${\bf G}$ sur $\frak{g}^*\otimes_F\overline{F}$, contient 
 ${\bf T}$. Soit ${\bf B}$ un sous--groupe de Borel de ${\bf G}$ contenant ${\bf T}$, et 
 soit ${\bf U}=R_{\rm u}({\bf B})$. Pour $Y\in \frak{g}^*_\gamma$, la ${\bf B}$--orbite de $Y$ 
 co\"{\i}ncide avec sa ${\bf U}$--orbite, laquelle est ferm\'ee 
 dans $\frak{g}^*\otimes_F\overline{F}$ pour la topologie de Zariski, 
d'aprs \cite[ch.~I, 4.10]{Bor}. Comme la vari\'et\'e quotient ${\bf G}/{\bf B}$ 
 est compl\`ete, on en d\'eduit que la ${\bf G}$--orbite de $Y$ 
 est ferm\'ee dans $\frak{g}^*\otimes_F\overline{F}$. Par cons\'equent la $G$--orbite de $Y$ est ferm\'ee dans $\frak{g}^*$ {\it pour la topologie 
 de Zariski}.

Soit maintenant $P$ un sous--groupe 
parabolique de $G$. Choisissons une composante de Levi $L$ de $P$, et notons  
$P^-$ le sous--groupe parabolique de $G$ oppos\'e \`a $P$ par rapport 
\`a $L$. On a la d\'ecomposition $\frak{g}=\frak{u}_{P^-}\oplus \frak{p}$. 
Posons $\frak{u}_{P^-}^*={\rm Hom}_F(\frak{u}_{P^-},F)$. Alors on a l'identification 
$P^-$--quivariante (pour l'action coadjointe de $P^-$ sur $\frak{g}^*$ et sur $\frak{u}_{P^-}^*$):
$$
\frak{p}^\perp =\frak{u}_{P^-}^*.
$$
Or pour tout 
$Y\in \frak{p}^\perp=\frak{u}_{P^-}^*$, la 
fermeture de la $P^-$--orbite de $Y$ dans 
$\frak{g}^*$ pour la topologie de Zariski, contient $0$. D'o\`u 
l'\'egalit\'e $\frak{g}_\gamma^*\cap \frak{p}^\perp =\{0\}$.
\end{proof}

\begin{moncoro}
L'ensemble $G^\natural_{\rm qr}$ est ouvert 
dense dans $G^\natural$ pour la topologie $\varpi$--adique.
\end{moncoro}

\begin{proof}
Puisque $G^\natural_{\rm reg}$ est dense dans $G^\natural$, 
l'ensemble $G^\natural_{\rm qr}$ l'est aussi.

On proc\`ede ensuite comme 
dans \cite[appendix, prop. A.3]{BH1}. Posons $\frak{N}=\bigcup_P \frak{p}^\perp \subset \frak{g}^*$ 
o\`u $P$ parcourt l'ensemble des sous--groupes paraboliques de $G$. Pour $P\in \ES{P}_\circ$, on a 
$\frak{p}^\perp \subset \frak{p}_\circ^\perp$. Pour tout sous--groupe 
parbolique $P$ de $G$, il existe un $g\in G$ tel que $g^{-1}Pg\in \ES{P}_\circ$, et puisque 
$G=K_\circ P_\circ$, on peut choisir $g$ dans $K_\circ$. Or pour $P=kP'k^{-1}$ o 
$k\in K_\circ$ et $P'\in \ES{P}_\circ$, on a
\begin{align*}
& \frak{p}={\rm Ad}_G(k)(\frak{p}'),\\
&\frak{p}^\perp= {\rm Ad}_G^*(k)(\frak{p}'^\perp)\subset {\rm Ad}_G^*(k)(\frak{p}_\circ^\perp).
\end{align*}
Par cons\'equent
$$\frak{N}=\bigcup_{k\in K_\circ}{\rm Ad}_G^*(k)(\frak{p}_\circ^\perp).
$$ 
Montrons que $\frak{N}$ est ferm\'e dans $\frak{g}^*$. 
Soit une suite $\{X_n:n\in \Bbb{Z}_{\geq 1}\}$ dans $\frak{N}$ qui converge vers un \'el\'ement 
$X\in \frak{g}^*$. Pour chaque entier $n\geq 1$, on \'ecrit $X_n= {\rm Ad}_G^*(k_n)(Y_n)$ o 
$k_n\in K_\circ $ et $Y_n\in \frak{p}_\circ^\perp$. Puisque $K_\circ$ est compact, 
quitte \`a remplacer $\{X_n\}$ par une sous--suite, on peut supposer que la suite $\{k_n:n\in \Bbb{Z}_{\geq 1}\}$ converge vers un \'el\'ement $k\in K_\circ$. Alors $Y_n={\rm Ad}_G^*(k_n^{-1})(X_n)$ tend vers ${\rm Ad}_G^*(k^{-1})(X)$ quand $n$ tend vers $+\infty$. Mais puisque $\frak{p}_\circ^\perp$ est ferm\'e dans $\frak{g}^*$, on a ${\rm Ad}_G^*(k^{-1})(X)\in \frak{p}_\circ^\perp$ et $X\in \frak{N}$. Donc 
$\frak{N}$ est ferm\'e dans $\frak{g}^*$.

Choisissons une base $\{e_1,\ldots ,e_d\}$ de $\frak{g}^*$ sur $F$. 
Notons $\frak{g}^*_0$ le sous--ensemble de $\frak{g}^*$ form\'e des $X=\sum_{i=1}^d x_ie_i$ ($x_i\in F$) tels que $\max\{\vert x_i\vert_F :i=1,\ldots ,d\}=1$. 
Alors $\frak{g}^*_0$ est compact dans $\frak{g}^*$. Posons $\frak{N}_0=\frak{N}\cap \frak{g}_0^*$. Alors $\frak{N}\smallsetminus \{0\}=F^\times \frak{N}_0$, et $\frak{N}_0$ est 
compact dans $\frak{g}^*$. 

Montrons que $G^\natural \smallsetminus G^\natural_{\rm qr}$ est ferm\'e dans $G^\natural$. 
Soit $\{\gamma_n:n\in \Bbb{Z}_{\geq 1}\}$ une suite dans $G^\natural 
\smallsetminus G^\natural_{\rm qr}$ qui converge vers un \'el\'ement $\gamma\in G^\natural$. 
Pour chaque entier $n\geq 1$, puisque $\frak{g}_{\gamma_n}^*\cap \frak{N}\neq \{0\}$ et 
$\frak{N}\smallsetminus \{0\}=F^\times\frak{N}_0$, 
il existe un \'el\'ement $X_n\in \frak{N}_0$ tel que 
${\rm Ad}_{G^\natural}(\gamma_n)(X_n)=X_n$. Puisque $\frak{N}_0$ est compact, 
quitte \`a remplacer la suite $\{\gamma_n\}$ par une sous--suite, 
on peut supposer que la suite $\{X_n\}$ converge vers un \'el\'ement $X\in \frak{N}_0$, et 
que la suite $\{k_n\}$ converge vers un \'el\'ement $k\in K_\circ$. 
Alors ${\rm Ad}_{G^\natural}^*(\gamma_n)(X_n)$ tend vers ${\rm Ad}_{G^\natural}^*(\gamma)(X)=X$ 
quand $n$ tend vers $+\infty$, et $\gamma\not\in G^\natural_{\rm qr}$. Donc $G^\natural_{\rm qr}$ est 
ouvert dans $G^\natural$.
\end{proof}

\begin{marema}
{\rm Pour $G^\natural=G={\rm GL}_n(F)$, d'apr\`es 
\cite[appendix, A2]{BH1}, un \'el\'ement $g\in G$ est dans $G_{\rm qr}$ si et seulement si son polyn\^ome 
caract\'eristique $P_{\rm car}(g)$ est produit de polyn\^omes irr\'eductibles sur $F$ deux----deux distincts, 
et il est dans $G_{\rm reg}$ si et seulement s'il est dans $G_{\rm qr}$ et 
si chaque facteur irr\'eductible de $P_{\rm car}(g)$ est s\'eparable. De mani\`ere 
\'equivalente, $G_{\rm reg}$ est l'ensemble des $g\in G$ tels que le polyn\^ome $P_{\rm car}(g)$ 
a $n$ racines distinctes dans $\overline{F}$.\hfill $\blacksquare$}
\end{marema}

\subsection{L'application $\ES{N}_{\theta,g_0}:G\rightarrow G$ pour $\theta$ localement 
fini}\label{l'application N}
Dans ce ${\rm n}^\circ$, on suppose de plus que le ${\bf G}$--espace tordu ${\bf G}^\natural$ est localement fini (cf. 
\ref{lments rguliers d'un espace tordu}); i.e. 
que $\theta\in {\rm Aut}_F^0({\bf G})$. 

Puisque le morphisme quotient 
${\bf G}_{\rm der}\rightarrow {\bf G}_{\rm ad}={\bf G}_{\rm der}/Z({\bf G}_{\rm der})$ est une 
$F$--isog\'enie centrale, 
tout \'el\'ement $\bar{g}\in {\bf G}_{\rm ad}(F)$ d\'efinit un $F$--automorphisme 
${\rm Int}_{\bf G}(\bar{g})$ de ${\bf G}$: on choisit un $g\in {\bf G}_{\rm der}$ 
qui rel\`eve $\bar{g}$, 
et l'on pose ${\rm Int}_{\bf G}(\bar{g})(x)={\rm Int}_{\bf G}(g)(x)$ ($x\in {\bf G}$); 
l'\'el\'ement ${\rm Int}_{\bf G}(\bar{g})(x)$ ne d\'epend pas du choix de 
$g$, et l'application ${\rm Int}_{\bf G}(\bar{g}):{\bf G}\rightarrow {\bf G}$ ainsi d\'efinie 
est un $F$--automorphisme. Soit $l_0$ le plus petit entier $k\geq 1$ tel que $\theta^k={\rm Int}_{\bf G}(\bar{g})$ pour 
un $\bar{g}\in {\bf G}_{\rm ad}(F)$, et soit $\bar{g}_0\in 
{\bf G}_{\rm ad}(F)$ tel que $\theta^{l_0}={\rm Int}_{\bf G}(\bar{g}_0)$. Puisque ${\bf G}_{\rm ad}(F)$ quotienté par 
l'image du morphisme canonique ${\bf G}_{\rm der}(F)\rightarrow {\bf G}_{\rm ad}(F)$ est de torsion, il existe un 
\'el\'ement $g_0\in {\bf G}_{\rm der}(F)$ et un entier $r\geq 1$ tels que 
$\bar{g}_0^r$ est l'image de $g_0$  
par ce morphisme canonique. Puisque
$$
\theta(g_0)\theta(x)\theta(g_0^{-1})= \theta^{rl_0+1}(x)= g_0\theta(x)g_0^{-1}\quad 
(x\in G),
$$
l'\'el\'ement $g_0^{-1}\theta(g_0)$ appartient au centre de ${\bf G}_{\rm der}(F)$. Par suite, 
quitte \`a remplacer $g_0$ par $g_0^k$ (et $r$ par $kr$) pour un entier $k\geq 1$, on peut 
supposer que $\theta(g_0)=g_0$. Fixons 
un tel $g_0$ minimisant $r$, et posons 
$l=rl_0$. Notons que $\bar{g}_0$ et $l$ sont 
d\'etermin\'es de mani\`ere unique par $\theta$, et que 
$g_0$ est d\'etermin\'e de 
mani\`ere unique modulo $Z(G^\natural)\cap {\bf G}_{\rm der}(F)$. Soit 
$$\ES{N}=\ES{N}_{\theta,g_0}:G\rightarrow G$$ l'application d\'efinie par
$$
\ES{N}(g)=g\theta(g)\cdots \theta^{l-1}(g)g_0\quad (g\in G).
$$
Elle d\'epend du choix de $g_0$, mais l'application 
$$\ES{N}'=\ES{N}'_\theta:G\rightarrow G/Z(G^\natural)$$ obtenue en composant 
$\ES{N}$ avec la projection canonique $G\rightarrow G/Z(G^\natural)$, n'en d\'epend pas. 
Pour $g,\,x\in G$, on a
$$
\ES{N}(x^{-1}g\theta(x))=x^{-1}\ES{N}(g)x.\leqno{(*)}
$$

Soit 
${\bf G}\rtimes \langle \theta\rangle $ le produit semidirect (dans la cat\'egorie des groupes) 
de ${\bf G}$ par le groupe abstrait engendr\'e par $\theta$, 
et soit ${\bf C}$ le sous--groupe cyclique de ${\bf G}\rtimes \langle \theta\rangle$ 
engendr\'e par $g_0^{-1}\rtimes \theta^l$. 
On peut d\'efinir 
le groupe quotient
$${\bf G}^+={\bf G}\rtimes \langle \theta\rangle/ {\bf C}.$$
Notons $\tilde{\theta}$ 
l'image de $1\rtimes \theta$ dans ${\bf G}^+$, et identifions ${\bf G}$ \`a l'image de 
${\bf G}\rtimes 1$ dans ${\bf G}^+$. Posons ${\bf G}_0= {\bf G}$ et ${\bf G}_i={\bf G}\tilde{\theta}^i$ ($i=1,\ldots ,l-1$). 
Alors on a la d\'ecomposition en union union disjointe:
$${\bf G}^+=\coprod_{i=0}^{l-1}{\bf G}_i.
$$ 
La multiplication et le passage \`a l'inverse dans ${\bf G}^+$ sont donn\'es par les 
relations suivantes, 
pour $x,\,y\in {\bf G}$ et $i,\,j\in \{0,\ldots , l-1\}$:
\begin{align*}
& \tilde{\theta}^l=g_0,\\
&(x\tilde{\theta}^i)(y\tilde{\theta}^j)=x\theta^i(y)\tilde{\theta}^{i+j},\\
&(x\tilde{\theta}^i)^{-1}=\tilde{\theta}^{-i}x^{-1}= \theta^{-i}(x^{-1})\tilde{\theta}^{-i}.
\end{align*}
Cela munit ${\bf G}^+$ d'une structure de groupe alg\'ebrique affine 
d\'efini sur $F$, de composante neutre ${\bf G}$, 
tel que $\tilde{\theta}$ appartient au groupe $G^+={\bf G}(F)$ des 
points $F$--rationnels de ${\bf G}^+$. Les ${\bf G}_i$ pour $i=0,\ldots ,l-1$, sont 
les composantes connexes de ${\bf G}^+$. Elles sont toutes d\'efinies sur $F$, et 
pour $i=0,\ldots ,l-1$, l'ensemble ${\bf G}_i(F)=G\tilde{\theta}^i$ des points $F$--rationnels de ${\bf G}_i$ 
muni de la topologie $\varpi$--adique, est un $G$--espace tordu. 
Identifions $G^\natural$ \`a ${\bf G}_1(F)$ via l'application $g\cdot \delta_1\mapsto g\tilde{\theta}$ 
($g\in G$). Alors pour 
$\gamma =g\cdot\delta_1\in G^\natural$, 
l'inverse de $\gamma$ dans $G^+$ est donn\'e par
$$
\gamma^{-1}= \tilde{\theta}^{-1}g^{-1}=
g_0^{-1}\tilde{\theta}^{l-1}g^{-1}=
g_0^{-1} \theta^{l-1}(g^{-1})\tilde{\theta}^{l-1},
$$
et l'on a
\begin{align*}
&\gamma^l=\ES{N}(g),\\
&g^{-1}\gamma^lg = \theta(\gamma^l).
\end{align*}
%

\begin{marema}
{\rm 
Supposons que la paire $(P_\circ,A_\circ)$ est $\theta$--stable. Alors 
puisque $\theta^l={\rm Int}_G(g_0)$, 
on a $g_0P_\circ g_0^{-1}=P_\circ$ et $g_0A_\circ g_0^{-1}=A_\circ$. En 
particulier, on a
$$
g_0\in N_G(A_\circ)\cap P_\circ= M_\circ.\eqno{\blacksquare}
$$}
\end{marema}

\subsection{La paire parabolique $(P_{[\gamma]},A_{[\gamma]})$ de $G$ associ\'ee \`a 
$\gamma\in G^\natural$}\label{la paire parabolique associe  gamma}
Pour $\tau\in {\rm Aut}(G)$ tel que la restriction de $\tau$ à $Z(G)$ est d'ordre fini, on note:
\begin{itemize}
\item  $P_{[\tau]}$ l'ensemble des $g\in G$ tels que $\{\tau^n(g):
n\in \Bbb{Z}_{\geq 1}\}$ est une partie born\'ee de $G$; 
\item $U_{[\tau]}$ le sous--ensemble de $P_{[\tau]}$ form\'e 
des $g\in G$ tels que $\lim_{n\rightarrow +\infty}\tau^n(g)=1$; 
\item $P^-_{\smash{[\tau]}}$ l'ensemble des $g\in G$ tels que $\{\tau^{-n}(g):
n\in \Bbb{Z}_{\geq 1}\}$ est une partie born\'ee de $G$;
\item $U^-_{\smash{[\tau]}}$ le sous--ensemble de $P^-_{\smash{[\tau]}}$ form\'e 
des $g\in G$ tels que $\lim_{n\rightarrow +\infty}\tau^{-n}(g)=1$;
\item $M_{[\tau]} = P_{[\tau]} \cap P^-_{\smash{[\tau]}}$. 
\end{itemize}
Ces cinq ensembles sont par d\'efinition des sous--groupes de $G$. 

Si le ${\bf G}$--espace tordu 
${\bf G}^\natural$ est localement fini, pour $\gamma\in G^\natural$ et $\tau={\rm Int}_{G^\natural}(\gamma)$, on remplace 
$\tau$ par $\gamma$ dans les notations ci-dessus; i.e. on pose 
$P_{[\gamma]} =P_{[\tau]}$, $U_{[\gamma]}=U_{[\tau]}$, (etc.). Pour $\gamma=g\cdot \delta_1\in G^\natural$, 
on a
$${\rm Int}_{G^\natural}(\gamma)^{-1}={\rm Int}_G(\theta^{-1}(g^{-1}))\circ \theta^{-1}.
$$ 
Par suite, posant $\bar{\gamma}=\theta^{-1}(g^{-1})\theta^{-1}\in G\theta^{-1}$, on a 
(si le ${\bf G}$--espace tordu 
${\bf G}^\natural$ est localement fini)
$$
P^-_{\smash{[\gamma]}}= P_{[\bar{\gamma}]},\quad U^-_{\smash{[\gamma]}}= U_{[\bar{\gamma}]}.
$$

\begin{marema1}
{\rm Supposons que le ${\bf G}$--espace tordu 
${\bf G}^\natural$ est localement fini, et identifions  $G^\natural$ 
\`a la composante connexe $G_1=G\tilde{\theta}$ de 
$G^+=\coprod_{i=0}^{l-1}G\tilde{\theta}^i$ comme en \ref{l'application N}. Pour $\gamma\in G^+$, l'automorphisme 
$g\mapsto \gamma g \gamma^{-1}$ de $G^+$ induit par restriction un 
automorphisme de $G$, que l'on note ${\rm Int}_{G^+}(\gamma)$. Si 
$\gamma=g\cdot \delta_1\in G^\natural$, cet automorphisme co\"{\i}ncide avec ${\rm Int}_{G^\natural}(\gamma)={\rm Int}_G(g)\circ\theta$. 
Pour $\gamma\in G^+$, la restriction de ${\rm Int}_{G^+}(\gamma)$ à $Z(G)$ est d'ordre fini, et l'on peut d\'efinir comme ci-dessus les sous--groupes $P_{[\gamma]}$,  
$U_{[\gamma]}$, $P^-_{\smash{[\gamma]}}$,  
$U^-_{\smash{[\gamma]}}$, $M_{[\gamma]}$ de $G$. Pour $n\in \Bbb{Z}$, on a l'égalité
${\rm Int}_{G^+}(\gamma)^n= {\rm Int}_{G^+}(\gamma^n)$. On en d\'eduit que pour 
$n\in \Bbb{Z}_{\geq 1}$, on a:
\begin{itemize}
\item $P_{[\gamma]}= P_{[\gamma^n]}$,
\item $U_{[\gamma]}$ et $U_{[\gamma^n]}$;
\end{itemize}
on a aussi:
\begin{itemize}
\item $P^-_{\smash{[\gamma]}}=P_{[\gamma^{-n}]}$,
\item $U^-_{\smash{[\gamma]}}= U_{[\gamma^{-n}]}$,
\item $M_{[\gamma]} =M_{[\gamma^n]}=M_{[\gamma^{-n}]}$.
\end{itemize}
En particulier, pour $g\in G$, 
on a $P_{[g\cdot \delta_1]} =P_{[\ES{N}(g)]}$, $U_{[g\cdot \delta_1]} =U_{[\ES{N}(g)]}$, (etc.).\hfill $\blacksquare$}
\end{marema1}

Soit $\Phi_\circ$ l'ensemble 
des racines de $A_\circ$ dans $G$, et soit 
$\Delta_\circ\subset \Phi_\circ$ l'ensemble des racines simples 
associ\'ees \`a $P_\circ$. Soit $v_F$ la valuation sur $\overline{F}$ normalis\'ee par $v_F(F^\times)=\Bbb{Z}$. 
Notons:
\begin{itemize}
\item $A_\circ^-$ l'ensemble des 
$t\in A_\circ$ tels que $v_F(\alpha(t))\geq 0$ pour toute racine $\alpha\in \Delta_\circ$;
\item $A_\circ^{-,\bullet}$ l'ensemble des 
$t\in A_\circ$ tels que $v_F(\alpha(t))> 0$ pour toute racine $\alpha\in \Delta_\circ$.
\end{itemize}
Pour tout $a\in A_\circ$, il existe un 
$g\in N_G(A_\circ)$ tel que $gag^{-1}\in A_\circ^-$, et si $a\in A_\circ^{-,\bullet}$, 
alors pour tout sous--groupe 
ouvert compact $J$ de $U_{P_\circ}$, on a 
$$\bigcap_{i\in \Bbb{Z}}a^i Ja^{-i}= \{1\},\quad 
\bigcup_{i\in \Bbb{Z}}a^i Ja^{-i}= U_{P_\circ}.
$$

\begin{mapropo}
Soit $\tau\in {\rm Aut}_F^0({\bf G})$. Alors:
\begin{itemize}
\item $P_{[\tau]}$ est un sous--groupe 
parabolique de $G$, et l'on a $U_{[\tau]}= U_{P_{[\tau]}}$;
\item $M_{[\tau]}$ est une composante de 
Levi de $P_{[\tau]}$;
\item $P^-_{\smash{[\tau]}}$ est le sous--groupe parabolique 
de $G$ oppos\'e \`a $P_{[\tau]}$ par rapport \`a $M_{[\tau]}$, et l'on a
$U^-_{\smash{[\tau]}}=U_{P^-_{\smash{[\tau]}}}$.
\end{itemize}
\end{mapropo}

\begin{proof}Quitte \`a remplacer ${\bf G}^\natural$ 
par ${\bf G}\tau$, on peut supposer que $\tau={\rm Int}_{{\bf G}^\natural}(\gamma)$ pour un $\gamma\in G^\natural$. 
Le ${\bf G}$--espace tordu ${\bf G}^\natural$ est alors localement fini. Identifions 
$G^\natural$ \`a la composante connexe $G_1=G\tilde{\theta}$ de 
$G^+=\coprod_{i=0}^{l-1}G\tilde{\theta}^i$, comme en \ref{l'application N}. Soit $x= \gamma^l\in G$, 
et notons $x=x_{\rm s}x_{\rm u}$ la d\'ecomposition de Jordan de $x$ (dans ${\bf G}$). 
Si ${\rm car}(F)=0$, 
alors $x_{\rm s}$ et $x_{\rm u}$ appartiennent \`a $G$. Si ${\rm car}(F)= p>1$, alors il existe un 
entier $m\geq 1$ tel que $(x_{\rm u})^{p^m}=1$. 
Puisque $x^n=(x_{\rm s})^n(x_{\rm u})^n$ ($n\in \Bbb{Z}$), 
quitte \`a remplacer $x$ par $x^n$ pour un entier $n\geq 1$, on peut 
supposer que $x_{\rm s}$ et $x_{\rm u}$ appartiennent \`a $G$. D'apr\`es le thorme de \ref{automorphismes qss lf} et le thorme 2 de \ref{automorphismes qss}, 
le groupe ${\bf H}={\bf G}_{x_{\rm s}}^\circ$ est r\'eductif, et d'apr\`es 
\cite[ch.~III, 9.1]{Bor} il est d\'efini sur $F$. Le radical ${\bf S}=R({\bf H})$ est un tore d\'efini $F$. 
Soit ${\bf A}$ le tore $F$--d\'eploy\'e maximal de ${\bf S}$, 
et soit ${\bf S}'$ le tore $F$--anisotrope maximal de ${\bf S}$. Notons $H$, $S$, $A$, $S'$ les groupes des points $F$--rationnels de ${\bf H}$, ${\bf S}$, ${\bf A}$, ${\bf S}'$. Rappelons que $x_{\rm s}$ et $x_{\rm u}$ appartiennent  ${\bf H}$ (\ref{automorphismes}, remarque (4)); en particulier $x_{\rm s}$ appartient  
$Z({\bf H})$. Puisque $Z({\bf H})^\circ =R({\bf H})$ et 
que le morphisme produit ${\bf A}\times 
{\bf S}'\rightarrow {\bf S}$ est une $F$--isog\'enie \cite[ch.~III, 8.15]{Bor}, quitte \`a remplacer une nouvelle fois 
$x$ par $x^n$ pour un entier $n\in \Bbb{Z}_{\geq 1}$, on peut supposer 
que $x_{\rm s}\in S$ et $s=as'$ pour des \'el\'ements $a\in A$ et $s'\in S'$. 
On a donc $x=as'x_{\rm u}$ avec $ax_{\rm u}=x_{\rm u}a$ et $s'x_{\rm u}=x_{\rm u}s'$. Puisque $S'$ est compact et que le 
sous--groupe de $H$ engendr\'e par $x_{\rm u}$ est born\'e, le sous--groupe de $H$ engendr\'e par 
$s'x_{\rm u}$ est lui aussi born\'e. On en d\'eduit que les groupes $P_{[\gamma]}=P_{[x]}$ et $P_{[a]}$ co\"{\i}ncident. De 
m\^eme, on a $U_{[\gamma]}= U_{[x]}=U_{[a]}$, (etc.).

On proc\`ede ensuite comme dans \cite[\S2]{Ca2}. Soit $y\in G$ tel que $yay^{-1}\in A_\circ^-$, et 
soit $I_{yay^{-1}}$ l'ensemble 
des racines $\alpha\in \Delta_\circ$ tels que $v_F(\alpha(yay^{-1}))=0$. Soit 
$P\in \ES{P}_\circ$ le sous--groupe parabolique de $G$ engendr\'e par $P_\circ$ et 
les sous--groupes radiciels de $G$ associ\'es aux 
racines $-\alpha$ pour $\alpha\in I_{yay^{-1}}$. Posons $M=M_P$. Alors on a:
\begin{itemize}
\item $P_{[a]}=y^{-1}P y$ et 
$U_{[a]}= y^{-1}U_Py$;
\item $P^-_{\smash{[a]}}= y^{-1}P^-y$ et $U^-_{\smash{[a]}}=y^{-1}U_{P^-}y$;
\item $M_{[a]}=y^{-1}My$.
\end{itemize}
D'o\`u la proposition.
\end{proof}

Notons que si ${\bf G}$ est un tore, disons ${\bf T}$, alors pour $\tau\in {\rm Aut}_F^0({\bf T})$, on a:
\begin{itemize}
\item $P_{[\tau]}=P^-_{\smash{[\tau]}} =M_{[\tau]}= T$,
\item $U_{[\tau]} = U^-_{\smash{[\tau]}}=\{1\}$.
\end{itemize}

On voudrait étendre les définitions des sous--groupes $P_{[\tau]}$, $U_{[\tau]}$, (etc.) de $G$ au cas d'un $F$--automor\-phisme $\tau$ de ${\bf G}$ qui n'est pas localement fini, de manière à ce que la proposition reste vraie. Soit donc $\tau\in {\rm Aut}_F({\bf G})$. Puisque $\tau_{\rm der}\in {\rm Aut}_F({\bf G}_{\rm der})= {\rm Aut}_F^0({\bf G}_{\rm der})$, on peut dfinir comme plus haut les sous--groupes $P_{[\tau_{\rm der}]}$, $U_{[\tau_{\rm der}]}$, (etc.) de $G_{\rm der}$. D'après la proposition, il 
existe un sous--groupe parabolique ${\bf P}'$ de 
${\bf G}_{\rm der}$ dfini sur $F$ et une composante de Levi ${\bf M}'$ de ${\bf P}'$ dfinie sur $F$ tels que, 
notant ${\bf P}'^-$ le sous--groupe parabolique de ${\bf G}_{\rm der}$ oppos  ${\bf P}'$ par rapport  ${\bf M}'$, on a:
\begin{itemize}
\item $P_{[\tau_{\rm der}]}={\bf P}'(F)$ et $U_{[\tau_{\rm der}]}=R_{\rm u}({\bf P}')(F)$;
\item $M_{[\tau_{\rm der}]}={\bf M}'(F)$;
\item $P^-_{\smash{[\tau_{\rm der}]}}={\bf P}'^-(F)$ et 
$U^-_{\smash{[\tau_{\rm der}]}}=R_{\rm u}({\bf P}'^-)(F)$.
\end{itemize}
Notons ${\bf P}$, ${\bf M}$, ${\bf P}^-$ les sous--groupes 
$R({\bf G}){\bf P}'$, $R({\bf G}){\bf M}'$, $R({\bf G}){\bf P}'^-$ de ${\bf G}$. Ils sont d\'efinis sur $F$, 
${\bf P}$ est un sous--groupe parabolique de ${\bf G}$, ${\bf M}$ est une composante de Levi de ${\bf P}$, 
et ${\bf P}^-$ est le sous--groupe parabolique de ${\bf G}$ oppos\'e \`a ${\bf P}$ par rapport \`a ${\bf M}$. 
Les groupes $R_{\rm u}({\bf P})=R_{\rm u}({\bf P}')$ et $R_{\rm u}({\bf P}^-)=R_{\rm u}({\bf P}'^-)$ sont 
eux aussi d\'efinis sur $F$; on les note ${\bf U}$ et ${\bf U}^-$. On pose:
\begin{itemize}
\item $P_{[\tau]}= {\bf P}(F)$; 
\item $U_{[\tau]}= {\bf U}(F)\;(= U_{[\tau_{\rm der}]})$;
\item $P^-_{[\tau]}= {\bf P}^-(F)$ ; 
\item $U^-_{[\tau]}= {\bf U}^-(F)\;(= U^-_{[\tau_{\rm der}]})$;
\item $M_{[\tau]}= {\bf M}(F)\;(=P_{[\tau]}\cap P^-_{[\tau]})$.
\end{itemize}
Si ${\bf G}= {\bf G}_{\rm der}$, ces définitions co\"{\i}ncident avec les précédentes, et pour ${\bf G}$ quelconque, on a:

\begin{itemize}
\item $P_{[\tau]}\cap {\bf G}_{\rm der}(F)= P_{[\tau_{\rm der}]}$;
\item $P^-_{[\tau]}\cap {\bf G}_{\rm der}(F)= P^-_{[\tau_{\rm der}]}$;
\item $M_{[\tau]}\cap {\bf G}_{\rm der}(F)= M_{[\tau_{\rm der}]}$.
\end{itemize}
On en déduit (pour ${\bf G}$ quelconque) que si $\tau\in {\rm Aut}_F^0({\bf G})$, ces définitions co\"{\i}ncident avec les précédentes, 
\cad que les notations sont cohérentes. De plus, par construction, la proposition reste vraie pour $\tau\in {\rm Aut}_F({\bf G})$. 
En particulier si ${\bf G}$ est un tore, disons ${\bf T}$, alors pour $\tau\in {\rm Aut}_F({\bf T})$, on a encore:
\begin{itemize}
\item $P_{[\tau]}=P^-_{\smash{[\tau]}} =M_{[\tau]}= T$,
\item $U_{[\tau]} = U^-_{\smash{[\tau]}}=\{1\}$.
\end{itemize}
%

\begin{marema2}
{\rm (On ne suppose plus que ${\bf G}$ est un tore.) Continuons avec les notations précédentes, et 
choisissons un tore maximal 
${\bf T}$ de ${\bf M}$ d\'efini sur $F$, et une sous--extension finie $E/F$ 
de $F^{\rm sep}/F$ d\'eployant ${\bf T}$. Notons 
$P_{E,[\tau]}$, $U_{E,[\tau]}$, (etc.) les sous--groupes de ${\bf G}(E)$ d\'efinis comme ci-dessus en rempla\c{c}ant $F$ 
par $E$. D'apr\`es la relation $(*)$ de \ref{groupes rductifs}, le morphisme produit 
${\bf T}\times {\bf G}_{\rm der}\rightarrow {\bf G}$ 
est s\'eparable, et comme ${\bf T}'={\bf T}\cap {\bf G}_{\rm der}$ est un tore d\'efini sur $F$ et 
d\'eploy\'e sur $E$, d'apr\`es la remarque 2 de \ref{gnralits}, on a l'galit
$$
{\bf G}(E)={\bf T}(E){\bf G}_{\rm der}(E).
$$
On en d\'eduit que
$$
P_{E,[\tau]}= {\bf P}(E)={\bf T}(E){\bf P}'(E)={\bf T}(E)P_{E,[\tau_{\rm der}]},
$$
et
$$
U_{E,[\tau]}={\bf U}(E)= {\bf U}'(E)= U_{E,[\tau_{\rm der}]}.
$$
On en déduit aussi que
$$
P^-_{\smash{E,[\tau]}}={\bf T}(E)P_{E,[\tau_{\rm der}]},\quad U^-_{\smash{E,[\tau]}}=U^-_{E,[\tau_{\rm der}]},
$$
et
$$
M_{E,[\tau]}={\bf M}(E)= {\bf T}(E)M_{E,[\tau_{\rm der}]}.
$$
Enfin notons que par définition, on a $P_{[\tau]}=P_{E,[\tau]}\cap G$, $U_{[\tau]}=U_{E,[\tau]}\cap G$, (etc.). 
}
\end{marema2}

\v1
Comme plus haut, pour $\gamma\in G^\natural$ et $\tau={\rm Int}_{G^\natural}(\gamma)$, on remplace 
$\tau$ par $\gamma$ dans les notations ci-dessus; i.e. on pose 
$P_{[\gamma]} =P_{[\tau]}$, $U_{[\gamma]}=U_{[\tau]}$, (etc.). Pour $\gamma=g\cdot \delta_1\in G^\natural$, 
avec la définition de $\bar{\gamma}$ donnée plus haut, on a encore (même si le ${\bf G}$--espace tordu 
${\bf G}^\natural$ n'est pas localement fini)
$$
P^-_{\smash{[\gamma]}}= P_{[\bar{\gamma}]},\quad U^-_{\smash{[\gamma]}}= U_{[\bar{\gamma}]}.
$$

\begin{madefi}
{\rm 
Pour $\gamma\in G^\natural$, on note $A_{[\gamma]}$ le tore d\'eploy\'e maximal du 
centre de $M_{[\gamma]}$. La paire $(P_{[\gamma]},A_{[\gamma]})$ est appel\'ee 
{\it la paire parabolique de $G$ associ\'ee \`a $\gamma$}. Un \'el\'ement 
$\gamma\in G^\natural$ tel que $P_{[\gamma]}\in \ES{P}_\circ$ et $A_{[\gamma]} =A_{P_{[\gamma]}}$ 
est dit {\it en position standard}.}
\end{madefi}

Tout \'el\'ement de $G^\natural$ est $G$--conjugu\'e \`a 
un \'el\'ement en position standard. En effet, 
pour $\gamma\in G^\natural$ et $x\in G$, on a
$$
(P_{[x^{-1}\cdot \gamma \cdot x]},A_{[x^{-1}\cdot \gamma \cdot x]})= (x^{-1}P_{[\gamma]} x,x^{-1}
A_{[\gamma]} x),$$
et fix\'e $\gamma\in G^\natural$, on peut choisir $x\in G$ tel que 
$P_{[x^{-1}\cdot\gamma x]}\supset P_\circ$ et $A_{[x^{-1}\cdot\gamma x]}\subset A_\circ$; autrement dit  tel que 
$P=P_{[x^{-1}\cdot\gamma x]}\in\ES{P}_\circ$ et $A_P=A_{[x^{-1}\cdot\gamma \cdot x]}$. 

Pour $\gamma\in G^\natural$, les sous--groupes paraboliques 
$P_{[\gamma]}$ et $P^-_{[\gamma]}$ de $G$ sont ${\rm Int}_{G^\natural}(\gamma)$--stables (par 
d\'efinition). Les sous--groupes $M_{[\gamma]}$ et $A_{[\gamma]}$ de $G$ sont donc eux aussi 
${\rm Int}_{G^\natural}(\gamma)$--stables. En particulier, $P_{[\gamma]}^\natural=P_{[\gamma]}\cdot\gamma$ et 
$P^{-,\natural}_{[\gamma]}=P^-_{[\gamma]}\cdot \gamma$ sont deux sous--espaces 
paraboliques de $G^\natural$, et $M^\natural_{[\gamma]} = 
M_{[\gamma]}\cdot \gamma$ est 
une composante de Levi de $P_{[\gamma]}^\natural$ (resp. de $P^{-,\natural}_{[\gamma]}$), 
qui vrifie
$$
M^\natural_{\smash{[\gamma]}} =P_{[\gamma]}^\natural\cap P^{-,\natural}_{[\gamma]}.
$$ Notons que 
$\gamma$ est en position standard si et seulement si 
$P_{[\gamma]}\in \ES{P}_\circ^{(\natural)}$ et $M_{[\gamma]}= M_{P_{[\gamma]}}$, i.e. si et seulement si 
$P_{[\gamma]}^\natural\in \ES{P}_\circ^\natural$ et $M_{[\gamma]}^\natural =M_{P_{[\gamma]}}^\natural$.

\begin{marema3}
{\rm Soit $\gamma = g\cdot \delta_1\in G^\natural$. 
Si $P_{[\gamma]}$ est $\theta$--stable, alors on a $g\in N_G(P_{[\gamma]})=P_{[\gamma]}$. 
Si $P_{[\gamma]}$ et $P^-_{\smash{[\gamma]}}$ sont $\theta$--stables, alors 
on a $g\in P_{[\gamma]}\cap P^-_{\smash{[\gamma]}}=M_{[\gamma]}$. R\'eciproquement, si 
$g\in P_{[\gamma]}$ (resp. si $g\in M_{[\gamma]}$), alors $P_{[\gamma]}$ est $\theta$--stable 
(resp. $P_{[\gamma]}$ et $P^-_{\smash{[\gamma]}}$ sont $\theta$--stables).

Supposons de plus que la paire $(P_\circ,A_\circ)$ est $\theta$--stable. 
Si $\gamma$ est en position standard, alors $g$ appartient  $M_{[\gamma]}$. En effet, 
supposons que $P=P_{[\gamma]}$ appartient  $\ES{P}_\circ$ et que $A_{[\gamma]}= A_P$. 
Puisque ${\rm Int}_{G^\natural}(\gamma)(P)=P$, on a 
$\theta(P)=g^{-1}Pg\in \ES{P}_\circ$ d'o\`u $g\in N_G(P)=P$, et 
puisque ${\rm Int}_{G^\natural}(\gamma)(M_{[\gamma]})=M_{[\gamma]} =M_P$, on a 
$\theta(M_P)=g^{-1}M_P g$. Or d'apr\`es le lemme de 
\ref{paires paraboliques}, on a $\theta(A_P)=A_P$ d'o\`u $\theta(M_P)=M_P$. Donc 
$g\in N_G(M_P)\cap P=M_P$.\hfill $\blacksquare$}
\end{marema3}

\begin{marema4}
{\rm Soit $P\in \ES{P}_\circ^{(\natural)}$, et soit $\gamma\in M_P^\natural$ 
tel que $P_{[\gamma]}=P$ et $A_{[\gamma]}=A_P$. Posons 
$M=M_P$, $U=U_P$, $U^-=U_P^-$, et $\tau={\rm Int}_{G^\natural}(\gamma)$. 
D'apr\`es \cite{D}, 
il existe une base $\{K_i:i\in \Bbb{Z}_{\geq 0}\}$ de voisinages de $1$ dans $G$ form\'ee de sous--groupes ouverts compacts 
$K_i$ de $G$ tels que, posant $K_i^-=K_i\cap U^-$, $K_i^0=K_i\cap M$ et 
$K_i^+= K_i\cap U$, on a
\begin{enumerate}
\item[-] $K_i= K_i^-K_i^0K_i^+$,
\item[-] $\tau^{-1}(K_i^-)\subset K_i^-$, $\tau(K_i^0)= K_i^0$ et $\tau(K_i^+)\subset K_i^+)$.
\end{enumerate}
Rappelons la construction. Fixons un syst\`eme de 
coordonn\'ees locales $\alpha$ au voisinage de $1$ dans $G$, i.e. 
un $\frak{o}_F$--r\'eseau $\Lambda$ dans $\frak{g}$  
et une application $\alpha:\Lambda \rightarrow G$ qui soit un isomorphisme 
de vari\'et\'es $\varpi$--adiques sur un voisinage de $1$ dans $G$. 
On peut choisir $\alpha$ tel que, posant 
$\Lambda^-=\Lambda \cap \frak{u}_{P^-}$, $\Lambda^0=\Lambda \cap \frak{m}_P$ et 
$\Lambda^+=\Lambda \cap \frak{u}_P$, on ait:
\begin{enumerate}
\item[-] $\alpha(0)=1$ et ${\rm d}\alpha_0={\rm id}_\frak{g}$;
\item[-] $\Lambda =\Lambda^-\oplus \Lambda^0
\oplus \Lambda^+$;
\item[-] ${\rm Ad}_{G^\natural}(\gamma)^{-1}(\Lambda^-)
\subset \Lambda^-$, ${\rm Ad}_{G^\natural}(\gamma)(\Lambda^0)
= \Lambda^0$ et ${\rm Ad}_{G^\natural}(\gamma)(\Lambda^+)
\subset \Lambda^+$.
\end{enumerate}
Alors il existe un entier $i_0\geq 0$ tel que, posant 
$K_i= \alpha(\varpi^{i+i_0}\Lambda)$ pour $i\in \Bbb{Z}_{\geq 0}$, l'ensemble 
$\{K_i:i\in \Bbb{Z}_{\geq 0}\}$ v\'erifie les conditions demand\'ees.\hfill $\blacksquare$}
\end{marema4}

\subsection{Le principe de submersion d'Harish--Chandra}\label{le principe de submersion}
Soit $P\in \ES{P}_\circ$. Pour $\gamma\in G^\natural$, consid\'erons 
l'application
$$
\psi_{P,\gamma}:G\times P\rightarrow G^\natural,\, (g,p)\mapsto 
g^{-1}\cdot \gamma \cdot gp
\quad (g\in G,\,p\in P). 
$$
Pour $G^\natural=G$ et $\gamma$ r\'egulier (\cad semisimple rgulier, cf. la proposition de \ref{automorphismes rguliers; le cas i}), le r\'esultat suivant est 
d\^u \`a Harish--Chandra \cite[theo. 1]{HC3}.

\begin{mapropo}
Pour $\gamma\in G^\natural_{\rm qr}$, l'application 
$\psi_{P,\gamma}$ est partout submersive.
\end{mapropo}

\begin{proof}
Pour $x\in G$ et $y\in P$, on a
$$
\psi_{P,\gamma}(xg,py)= \psi_{P,x^{-1}\cdot \gamma\cdot x}(g,p)\cdot y
$$
Il suffit donc de montrer que $\psi_{P,\gamma}$ est submersive en $(g,p)=(1,1)$. 
En calculant ${\rm d}(\psi_{P,\gamma})_{1,1}$, on obtient que $\Psi_{P,\gamma}$ est 
submersive en $(1,1)$ si et seulement si on a l'\'egalit\'e
$$
\frak{g}(1-\gamma) + \frak{p}'=\frak{g};
$$
o\`u l'on a pos\'e $\frak{p}'={\rm Ad}_{G^\natural}(\gamma)(\frak{p})$. 
D'o\`u la proposition.
\end{proof}

On peut donc appliquer ici le {\it principe de submersion d'Harish--Chandra}: d'aprs \cite[theo. 11, p.49]{HC1}, il 
existe une unique application lin\'eaire
$$
C^\infty_{\rm c}(G\times P)\rightarrow C^\infty_{\rm c}(G^\natural),\,\alpha\mapsto f_{\alpha,\gamma}
$$
telle que pour toute fonction $\Phi\in C^\infty_{\rm c}(G^\natural)$, on a l'\'egalit\'e
$$
\int\!\!\!\!\int_{G\times P}\alpha(g,p)\Phi(g^{-1}\cdot\gamma\cdot gp)
dgd_lp =\int_{G^\natural}f_{\alpha,\gamma}(\gamma')\Phi(\gamma')d\gamma'.
$$
D'ailleurs, l'\'egalit\'e ci-dessus reste vraie pour toute fonction $\Phi$ localement 
int\'egrable sur $G^\natural$ (par rapport \`a une mesure de Haar sur $G^\natural$).

\begin{monlem}
Soit une fonction $\alpha\in C^\infty_{\rm c}(G\times P)$. L'application
$$
G^\natural_{\rm qr}\rightarrow C^\infty_{\rm c}(G^\natural),\,\gamma\mapsto f_{\alpha,\gamma}
$$
est localement constante.
\end{monlem}

\begin{proof}
L'application
$$
G_{\rm qr}^\natural\times G\times P\rightarrow G^\natural_{\rm qr}
\times G^\natural,\,(\gamma,g,p)\mapsto
(\gamma, \psi_{P,\gamma}(g,p))
$$
est partout submersive. Par suite il existe une unique application 
lin\'eaire
$$
C^\infty_{\rm c}(G_{\rm qr}^\natural\times G\times P)\rightarrow 
C^\infty_{\rm c}(G^\natural_{\rm qr}
\times G^\natural),\,\beta\mapsto f_\beta
$$
telle que pour toute fonction $\Psi\in C^\infty_{\rm c}(G^\natural_{\rm qr}\times 
G^\natural)$, on a l'\'egalit\'e
$$
\int\!\!\!\!\int\!\!\!\!\int_{G^\natural_{\rm qr}\times G\times P}
\beta(\gamma,g,p)\Psi(\gamma,g^{-1}\cdot \gamma \cdot gp)d\gamma dg d_lp=
\int\!\!\!\!\int_{G^\natural_{\rm qr}\times G^\natural}f_\beta(\gamma,\gamma')
\Psi(\gamma,\gamma')d\gamma d\gamma'.
$$
Soit $\Psi=\lambda\otimes \Phi\in C^\infty_{\rm c}(G^\natural_{\rm qr}\times G^\natural)$ 
pour des fonctions $\lambda\in C^\infty_{\rm c}(G^\natural_{\rm qr})$ et 
$\Phi\in C^\infty_{\rm c}(G^\natural)$. Alors posant
$$
\Lambda_\beta(\Phi,\gamma)=\int\!\!\!\!\int_{G\times P}
\beta(\gamma,g,p)\Phi(g^{-1}\cdot \gamma \cdot gp)dg d_lp\quad (\gamma\in G^\natural_{\rm qr}),
$$
on a l'galit
$$
\int_{G^\natural_{\rm qr}}\lambda(\gamma)\Lambda_\beta(\Phi,\gamma)d\gamma
=\int_{G^\natural_{\rm qr}}\lambda(\gamma)\left\{\int_{G^\natural}f_\beta(\gamma,\gamma')
\Phi(\gamma') d\gamma'\right\} d\gamma.
$$
Puisque l'\'egalit\'e ci-dessus est vraie pour toute fonction 
$\lambda\in C^\infty_{\rm c}(G^\natural_{\rm qr})$, on en d\'eduit 
que pour toute fonction 
$\Phi\in C^\infty_{\rm c}(G^\natural)$ et tout $\gamma\in G^\natural_{\rm qr}$, on a l'\'egalit\'e
$$
\Lambda_\beta(\Phi,\gamma)=
\int_{G^\natural}f_\beta(\gamma,\gamma')\Phi(\gamma') d\gamma'.
$$

Fixons un \'el\'ement $\gamma_0\in G^\natural_{\rm qr}$ et une fonction 
$\mu\in C^\infty_{\rm c}(G^\natural_{\rm qr})$ telle que $\mu(\gamma_0)=1$. 
Soit aussi une fonction $\alpha\in C^\infty_{\rm c}(G\times P)$. Posons 
$\beta = \mu\otimes \alpha\in C^\infty_{\rm c}(G^\natural_{\rm qr}\times G\times P)$. 
Soit $\Omega$ un voisinage de $\gamma_0$ dans $G^\natural_{\rm qr}$ tel que 
$\mu\vert_\Omega=1$. Alors pour $(\gamma,g,p)\in \Omega\times G\times P$ et 
$\Phi\in C^\infty_{\rm c}(G^\natural)$, on a 
$\beta(\gamma,g,p)=\alpha(g,p)$; par cons\'equent 
\begin{align*}
\Lambda_\beta(\Phi,\gamma)
&=\int\!\!\!\!\int_{G\times P}\alpha(g,p)\Phi(g^{-1}\cdot \gamma \cdot gp)dgd_lp\\
&=\int_{G^\natural}f_{\alpha,\gamma}(\gamma')\Phi(\gamma')d\gamma'
\end{align*}
et
$$
\int_{G^\natural}f_{\alpha,\gamma}(\gamma')\Phi(\gamma')d\gamma'=
\int_{G^\natural}f_\beta(\gamma,\gamma')\Phi(\gamma') d\gamma'.
$$
L'\'egalit\'e ci--dessus \'etant vraie pour tout $\gamma\in \Omega$ et toute 
fonction $\Phi\in C^\infty_{\rm c}(G^\natural)$, on en d\'eduit l'galit
$$
f_{\alpha,\gamma}(\gamma')=f_\beta(\gamma,\gamma')\quad (\gamma\in \Omega,\,\gamma'\in G^\natural).
$$
D'o\`u le lemme, puisque $f_\beta\in C^\infty_{\rm c}(G^\natural_{\rm qr}\times G^\natural)$.
\end{proof}

\subsection{Les op\'erateurs ${\rm T}_\gamma$ pour $\gamma\in G^\natural$ quasi--r\'egulier}
\label{les oprateurs T} 
Soit $(\Pi, V)$ une 
$\omega$--repr\'esen\-ta\-tion admissible de $G^\natural$. Posons $\pi=\Pi^\circ$. Notons ${\rm End}_0(V)$ 
l'image canonique de $V\otimes_\Bbb{C}\check{V}$ dans ${\rm End}_\Bbb{C}(V)$, 
o\`u $\check{V}$ d\'esigne l'espace de la contragr\'ediente $\check{\pi}$ de $\pi$. De mani\`ere 
\'equivalente, ${\rm End}_0(V)$ est l'espace des applications lin\'eaires $u:V\rightarrow V$ 
telles que les deux applications $G\rightarrow {\rm End}_\Bbb{C}(V),\,g\mapsto \pi(g)\circ u$ et 
$G\rightarrow {\rm End}_\Bbb{C}(V),\,g\mapsto u\circ \pi(g)$, 
sont localement constantes.

Soit $M_\circ^+$ l'ensemble 
des $m\in M_\circ$ tels que pour toute racine $\alpha$ de $A_\circ$ dans $U_\circ$, on a
$\langle H_{M_\circ}(m),\alpha\rangle \leq 0$, o
$$H_{M_\circ}:
M_\circ\rightarrow {\rm X}_*(A_\circ)\otimes_\Bbb{Z}\Bbb{R}$$ est l'application d'Harish--Chandra, 
d\'efinie comme suit: pour $m\in M$ et $\lambda\in {\rm X}^*(A_\circ)$, il existe un entier 
$d\geq 1$ tel que $\lambda^d=\mu\vert_{A_\circ}$ pour un (unique) caract\`ere rationnel
$\mu\in {\rm X}^*_F(M_\circ)$, et l'on pose $
\langle H_M(m),\lambda\rangle= -{1\over d}v_F(\mu(m))$. On a 
la {\it d\'ecomposition de Cartan} \cite[prop.~4.4.3]{BT1}:
$$
G=K_\circ M_\circ^+ K_\circ.\leqno{(*)}
$$

\begin{marema}
{\rm 
Posons
$A_\circ^+=\{t^{-1}: t\in A_\circ^-\}$ (cf. \ref{la paire parabolique associe  gamma}). L'ensemble 
$A_\circ \cap M_\circ^+$ est contenu dans $A_\circ^+$, et si ${\bf G}$ est semisimple, alors on a 
$A_\circ\cap M_\circ^+=A_\circ^+$.\hfill $\blacksquare$}
\end{marema}

Fixons un sous--groupe ouvert compact $K$ de $G$. Pour $\gamma\in G^\natural$, notons 
${\rm T}_\gamma = {\rm T}_{\Pi,\gamma}^K$ le ${\Bbb C}$--endomorphisme de $V$ d\'efini 
par
\begin{align*}
{\rm T}_\gamma &= {\rm vol}(K,dk)^{-1}\int_K\omega(k^{-1})\Pi(k^{-1}\cdot \gamma \cdot k)dk\\
& = {\rm vol}(K,dk)^{-1}\int_K\pi(k^{-1})\circ \Pi(\gamma)\circ \pi(k)dk.
\end{align*}

Pour $G^\natural=G$, $\omega=1$, $\gamma$ r\'egulier et $K=K_\circ$, le th\'eor\`eme suivant est 
d\^u \`a Harish--Chandra \cite[theo. 2]{HC3}. L'hypoth\`ese $K=K_\circ$ a ensuite 
\'et\'e supprim\'ee par Rader et Silberger \cite{RS}.

\begin{montheo}
Supposons que $\pi$ est de type fini (i.e. de longueur finie). Pour $\gamma\in G^\natural_{\rm qr}$, l'op\'erateur 
${\rm T}_\gamma$ appartient \`a l'espace ${\rm End}_0(V)$, et l'application 
$G^\natural_{\rm qr}\rightarrow {\rm End}_0(V),\,\gamma\mapsto {\rm T}_\gamma$ est localement constante.
\end{montheo}

\begin{proof}
Notons que si le th\'eor\`eme est vrai, alors pour tout sous--groupe 
compact $K'$ de $G$ contenant $K$, le th\'eor\`eme reste vrai si l'on remplace 
$T_\gamma$ par $T_{\Pi,\gamma}^{K'}$. Quitte \`a remplacer $K$ par un groupe plus petit, on 
peut donc supposer que $K$ est un sous--groupe distingu\'e de $K_\circ$ tel que $\omega\vert_K=1$.

Notons $K_1=K,\, K_2,\,\ldots, \,K_n$ les classes de $K_\circ/K$ ($K_i= x_iK_\circ = K_\circ x_i$ pour un $x_i\in K_\circ$), et pour $\gamma\in 
G^\natural_{\rm qr}$, 
posons  $c={\rm vol}(K,dg)$ et
$$
{\rm T}_{\gamma,i}=c^{-1}\int_{K_i}\omega(g^{-1})\Pi(g^{-1}\cdot \gamma \cdot g)dg.\quad (i=1,\ldots ,n).
$$
Pour $i=1,\ldots ,n$ et $x\in K_i$, on a donc 
$$
{\rm T}_{\gamma,i}= \omega(x^{-1})T_{x^{-1}\cdot \gamma \cdot x}= \pi(x^{-1})\circ T_\gamma \circ \pi(x).
$$

Puisque $\pi$ est de type fini, il existe un sous--groupe 
ouvert distingu\'e $K'$ de $K_\circ$ tel que $V$ est engendr\'e sur $G$ par 
le sous--espace $V^{K'}=\{v\in V: \pi(k')(v)=v,\, \forall k'\in K'\}$.

Puisque le groupe $M_\circ/A_\circ$ 
est compact, il existe une partie compacte $\Omega$ de $M_\circ$ 
telle que $M_\circ^+\subset \Omega (A_\circ \cap M_\circ^+)$. On en d\'eduit qu'il existe un sous--groupe 
ouvert compact $J_{P_\circ}$ de $P_\circ$ tel que pour tout $m\in M_\circ^+$, on a 
l'inclusion $m^{-1}J_{P_\circ}m\subset K'$. Pour $i=1,\ldots , n$, notons 
$\alpha_i\in C^\infty_{\rm c}(G\times P_\circ)$ la 
fonction d\'efinie par
$$
\alpha_i(g,p)=\left\{
\begin{array}{ll}c^{-1}\omega(g^{-1})\hfill & \hbox{si $(g,p)\in K_i\times J_{P_\circ}$}\\
0\hfill & \hbox{sinon}
\end{array}\right..
$$
D'apr\`es \ref{le principe de submersion}, pour $i=1,\ldots ,n$ et $\gamma\in G^\natural_{\rm qr}$, 
il existe une unique fonction $f_{\alpha_i,\gamma}\in C^\infty_{\rm c}(G^\natural)$ telle que pour toute fonction 
$\Phi$ localement int\'egrable sur $G^\natural$, on a l'\'egalit\'e
$$
c^{-1}\int\!\!\!\!\int_{K_i\times J_{P_\circ}}\omega(g^{-1})\Phi(g^{-1}\cdot \gamma\cdot gp)dgd_lp
=\int_{G^\natural}f_{\alpha_i,\gamma}(\gamma')\Phi(\gamma')d\gamma'.
$$
Puisque pour $(v,\check{v})\in V\times \check{V}$, l'application
$$
G^\natural\rightarrow \Bbb{C},\, \gamma' \mapsto \langle \Pi(\gamma')(v),\check{v}\rangle
$$
est localement constante (donc localement int\'egrable), pour 
$i=1,\ldots ,n$ et $\gamma\in G^\natural_{\rm qr}$, on a
\begin{align*}
\int_{J_{P_\circ}}{\rm T}_{\gamma,i}\circ \pi(p)d_lp &=
c^{-1}\int\!\!\!\!\int_{K_i\times J_{P_\circ}}\omega(g^{-1})\Pi(g^{-1}\cdot \gamma\cdot gp)dgd_lp \\
&=\Pi(f_{\alpha_i,\gamma}).
\end{align*}
Fixons un 
\'el\'ement $\gamma_0\in G^\natural_{\rm qr}$. Pour 
chaque $i$, la fonction $G^\natural_{\rm qr}\rightarrow  C^\infty_{\rm c}(G^\natural),\, 
\gamma\mapsto f_{\alpha_i,\gamma}$ est localement constante (lemme de \ref{le principe de submersion}). 
Par suite il existe un 
voisinage $\ES{V}$ de $\gamma_0$ dans $G^\natural_{\rm qr}$ et un 
sous--groupe ouvert $K''$ de $K$ distingu\'e dans $K_\circ$, tels que pour $i=1,\ldots ,n$ et 
$\gamma\in \ES{V}$, la fonction 
$f_{\alpha_i,\gamma}$ est bi--invariante (i.e.  gauche et  droite) par $K''$. 

Puisque $K'$ est distingu\'e dans $K_\circ$, d'apr\`es la dcomposition de Cartan $(*)$,
$V$ est engendr\'e (sur $\Bbb{C}$) par les vecteurs 
$\pi(km)(v)$ pour $k\in K_\circ$, $m\in M_\circ^+$ et $v\in V^{K'}$. Fixons de  
tels $k$, $m$, $v$, et notons 
$i\in \{1,\ldots ,n\}$ l'indice tel que $k\in K_i$. Soit $\gamma\in \ES{V}$. 
On a
$$
{\rm T}_\gamma \circ \pi(km)(v)= \pi(k)\circ {\rm T}_{\gamma,i}\circ \pi(m)(v),
$$
et l'on a aussi
\begin{align*}
\Pi(f_{\alpha_i,\gamma})\circ \pi(m)(v) & = 
\int_{J_{P_\circ}}{\rm T}_{\gamma,i}\circ \pi(p)\circ \pi(m)(v)d_lp \\
& = \int_{J_{P_\circ}}{\rm T}_{\gamma,i}\circ \pi(m)\circ \pi(m^{-1}pm)(v)d_lp\\
& = {\rm T}_{\gamma,i}\circ \pi(m) \circ \left\{\int_{J_{P_\circ}} \pi(m^{-1}pm)(v)d_lp\right\}\\
& = {\rm vol}(J_{P_\circ},d_lp) {\rm T}_{\gamma,i}\circ \pi(m)(v).
\end{align*}
Soit maintenant $e_{K''}\in C^\infty_{\rm c}(G)$ la fonction 
caract\'eristique de $K''$ divis\'ee par ${\rm vol}(K''\!,dg)$. 
D'apr\`es les calculs ci-dessus, on a
\begin{align*}
\lefteqn{\pi(e_{K''})\circ {\rm T}_\gamma \circ \pi(km)(v)}\\
& = {\rm vol}(J_{P_\circ},d_lp)^{-1}
\pi(e_{K''})\circ \pi(k)\circ \Pi(f_{\alpha_i,\gamma})\circ \pi(m)(v).
\end{align*}
Or $K_\circ$ normalise $K''$, par cons\'equent les op\'erateurs $\pi(e_{K''})$ et $\pi(k)$ 
commutent. Comme la fonction $f_{\alpha_i,\gamma}$ est $K''$--invariante \`a gauche, 
on a $\pi(e_{K''})\circ \Pi(f_{\alpha_i,\gamma})= \Pi(f_{\alpha_i,\gamma})$. On en d\'eduit que
\begin{align*}
\pi(e_{K''})\circ {\rm T}_\gamma \circ \pi(km)(v) &= {\rm vol}(J_{P_\circ},d_lp)^{-1}\pi(k)\circ 
\Pi(f_{\alpha_i,\gamma})\circ \pi(m)(v)\\
&=\pi(k)\circ {\rm T}_{\gamma,i}\circ \pi(m)(v)\\
& = {\rm T}_\gamma \circ \pi(km)(v).
\end{align*}
Cela \'etant vrai pour tous $k,\,m,\,v$ et tout $\gamma$, on a montr\'e que 
$$
\pi(e_{K''})\circ {\rm T}_\gamma = {\rm T}_\gamma\quad (\gamma\in \ES{V}).
$$

Pour $y\in G^\natural$, d'apr\`es la d\'efinition de ${\rm T}_y$, on a $\pi(e_{K''})\circ {\rm T}_y = 
{\rm T}_y\circ \pi(e_{K''})$. 
Pour $\gamma\in \ES{V}$, on a donc $
\pi(e_{K''})\circ {\rm T}_\gamma = {\rm T}_\gamma\circ \pi(e_{K''})={\rm T}_\gamma$, d'o\`u 
${\rm T}_\gamma\in {\rm End}_0(V)$. D'autre part, puisque 
$K''$ est distingu\'e dans $K$, on a aussi
$$
{\rm T}_{x\cdot \gamma \cdot y}= {\rm T}_\gamma \quad (x,\, y\in K'',\,\gamma\in \ES{V}).
$$
L'application $G^\natural_{\rm qr}\rightarrow {\rm End}_0(V),\,\gamma\mapsto {\rm T}_\gamma$ est 
donc bien localement constante.
\end{proof}

Puisque $\pi$ est admissible, 
les \'el\'ements de ${\rm End}_0(V)$ sont des op\'erateurs de rang fini sur $V$. Pour 
$g\in G_{\rm qr}^\natural$, on peut donc d\'efinir la trace de ${\rm T}_\gamma$, que l'on note 
${\rm tr}({\rm T}_\gamma)$.

\begin{moncoro}
(On suppose toujours $\pi$ admissible et de longueur finie.) 
Le caract\`ere $\Theta_\Pi$ est repr\'esent\'e sur 
$G^\natural_{\rm qr}$ par la fonction localement constante
$$\gamma\mapsto \Theta_\Pi(\gamma)= {\rm tr}({\rm T}_\gamma).
$$ 
En d'autres termes, pour toute fonction $\phi\in C^\infty_{\rm c}(G^\natural_{\rm qr})$, on a 
l'\'egalit\'e
$$
\Theta_\Pi(\phi)=\int_{G^\natural_{\rm qr}}\phi(\gamma)\Theta_\Pi(\gamma)d\gamma.
$$
\end{moncoro}

\begin{proof}
Pour $\phi\in C^\infty_{\rm c}(G^\natural)$, notons $\phi_0 \in 
C^\infty_{\rm c}(G^\natural)$ la fonction d\'efinie par
$$
\phi_0(\gamma)={\rm vol}(K,dk)^{-1}\int_K\omega(k^{-1})\phi(k\cdot \gamma \cdot k^{-1})dk.
$$
Pour $x\in G$, on a (relation $(*)$ de \ref{caractres tordus (bis)})
$$\Theta_\Pi({^x\phi})=\omega(x^{-1})\Theta_\Pi(\phi),$$
par cons\'equent $\Theta_\Pi(\phi)=\Theta_\Pi(\phi_0)$.

Soit une fonction $\phi\in C^\infty_{\rm c}(G^\natural_{\rm qr})$. Alors 
$\phi_0\in C^\infty_{\rm c}(G^\natural_{\rm qr})$, et d'apr\`es le thorme, on a
$$
\int_{G^\natural}\phi(\gamma){\rm tr}({\rm T}_\gamma)d\gamma = {\rm tr}\left(\int_{G^\natural}
\phi(\gamma){\rm T}_\gamma d\gamma \right).
$$
On en d\'eduit que
\begin{align*}
\int_{G^\natural}\phi(\gamma){\rm tr}(T_\gamma)d\gamma &={\rm tr}\left(\int_G\left\{\int_K
\phi(\gamma)\omega(k^{-1})\Pi(k^{-1}\cdot \gamma \cdot k)dk\right\} d\gamma\right)\\
&= {\rm tr}\left(\int_G\phi_0(\gamma)\Pi(\gamma)d\gamma\right)\\
& = \Theta_\Pi(\phi_0).
\end{align*}
D'o\`u le corollaire.
\end{proof}

Puisque pour $\phi\in C^\infty_{\rm c}(G^\natural)$ et $x\in G$, on a 
$\Theta_\Pi({^x\phi})=\omega(x^{-1})\Theta_\Pi(\phi)$, on a l'\'egalit\'e
$$
\Theta_\Pi(x^{-1}\cdot\gamma \cdot x)=\omega(x)\Theta_\Pi(\gamma)\quad (\gamma\in G^\natural_{\rm qr},\, x\in G).\leqno{(**)}
$$

\subsection{Induction parabolique et caract\`eres}\label{induction parabolique et caractres}
Pour $P\in \ES{P}_\circ$, on note
$$
\iota_P^G:\frak{R}(M_P)\rightarrow \frak{R}(G)
$$
le foncteur {\it induction parabolique normalis\'ee}. On rappelle la 
d\'efinition: pour toute repr\'esentation lisse $(\sigma,W)$ de $M_P$, 
on note encore $\sigma$ la repr\'esentation lisse $\sigma\otimes 1:M_P\ltimes U_P\rightarrow {\rm Aut}_\Bbb{C}(W)$ de $P$, et l'on pose 
$\iota_P^G(\sigma)={\rm ind}_P^G(\delta_P^{1/2}\sigma)$ o\`u 
(rappel) $\delta_P$ est le caract\`ere $\Delta_P^{-1}$ de $P$.

\begin{monlem}
Pour tout sous--groupe 
parabolique $P$ de $G$, on a $\omega\vert_{U_P}=1$.
\end{monlem}

\begin{proof}
Soit $P$ un sous--groupe parabolique de $G$. Soit 
$x\in G$ tel que $x^{-1}Px\supset P_\circ$. Alors $U_{x^{-1}Px}= x^{-1}U_Px \subset U_{P_\circ}$, 
et comme $\omega(x^{-1}ux)=\omega(u)$ ($u\in U_P$), il suffit de montrer que $\omega\vert_{U_{P_\circ}}=1$. Soit $J$ un sous--groupe ouvert compact 
de $G$ tel que $\omega\vert_J=1$, et soit $a\in A_\circ^{-,\bullet}$. Puisque 
$\bigcup_{i\in \Bbb{Z}}a^i(J\cap U_{P_\circ})a^{-i}= U_{P_\circ}$ et $\omega\vert_{J\cap U_{P_\circ}}=1$, on en déduit que 
$\omega\vert_{U_{P_\circ}}=1$.
\end{proof} 

Pour $P^\natural \in \ES{P}_\circ^\natural$, on d\'efinit comme suit un foncteur
$$
{^\omega{\iota}}_{P^\natural}^{G^\natural}:\frak{R}(M_P^\natural,\omega)\rightarrow  \frak{R}(G^\natural,\omega).
$$
Soit $(\Sigma,W)$ une $\omega$--repr\'esentation lisse de $M_P^\natural$. Posons 
$\sigma=\Sigma^\circ$. 
Pour $\gamma\in P^\natural$, on \'ecrit $\gamma =\delta\cdot u$ avec $\delta\in M_P^\natural$ et $u\in U_P$ 
(rappelons que 
l'\'ecriture est unique), et l'on pose $\Sigma(\gamma)= \Sigma(\delta)$. Pour $\gamma\in P^\natural$ et 
$p,\,p'\in P$, en \'ecrivant 
$p=mu$ et $p'=m'u'$ avec $m,\,m'\in M_P$ et $u,\,u'\in U_P$, on a
$$
p\cdot \gamma \cdot p'= (m\cdot \gamma\cdot m')\cdot {\rm Int}_{P^\natural}(\gamma\cdot m')^{-1}(u)u'
$$
avec ${\rm Int}_{P^\natural}(\gamma\cdot m')^{-1}(u)u'\in U_P$, d'o\`u (en utilisant le lemme)
$$
\Sigma(p\cdot\gamma\cdot p')= \omega(m')\sigma(m)\circ \Sigma(\gamma)\circ \sigma(m')=\omega(p')
\sigma(p)\circ \Sigma(\gamma)\circ \sigma(p').
$$
En d'autres termes, 
l'application $\gamma\mapsto \Sigma(\gamma)$ est une 
$\omega$--repr\'esentation lisse de $P^\natural$. 
D'apr\`es la relation $(*)$ de \ref{sous-espaces paraboliques} et le lemme de \ref{modules}, $\delta_{P^\natural}^{1/2}\Sigma$ est encore 
une $\omega$--repr\'esentation de 
$P^\natural$, telle que $(\delta_{P^\natural}^{1/2}\Sigma)^\circ= \delta_P^{1/2}\sigma$. On peut donc 
poser
$${^\omega{\iota}}_{P^\natural}^{G^\natural}(\Sigma)= 
{^\omega{\rm ind}}_{P^\natural}^{G^\natural}(\delta_{P^\natural}^{1/2}\Sigma).
$$

Par construction, les foncteurs ${^\omega{\iota}}_{P^\natural}^{G^\natural}:\frak{R}(M_P^\natural,\omega)\rightarrow  \frak{R}(G^\natural,\omega)$ et 
$\iota_P^G:\frak{R}(M_P)\rightarrow \frak{R}(G)$ commutent aux foncteurs d'oubli: pour toute 
$\omega$--repr\'esentation lisse $\Sigma$ de $M_P^\natural$, on a
$$
{^\omega{\iota}}_{P^\natural}^{G^\natural}(\Sigma)^\circ = \iota_P^G(\Sigma^\circ).
$$

Pour $P^\natural\in \ES{P}_\circ^\natural$, on note
$$
C^\infty_{\rm c}(G^\natural)\rightarrow C^\infty_{\rm c}(M_P^\natural),\,\phi\mapsto
{^\omega{\phi}_{P^\natural,K_\circ}}
$$ l'application 
lin\'eaire d\'efinie par
$$
{^\omega\phi_{P^\natural,K_\circ}}(\delta)=\delta_{P^\natural}^{1/2}(\delta)
\int\!\!\!\!\int_{U_P\times K_\circ} \omega(k)\phi(k^{-1}\cdot \delta \cdot uk)dudk
\quad (\delta\in M_P^\natural).
$$
Le th\'eor\`eme suivant est une simple reformulation du corollaire  
de \ref{caractres des induites compactes}. Il g\'en\'eralise la formule bien connue de Van Dijk \cite{VD, Cl1}.

\begin{montheo}
Soit $P^\natural\in \ES{P}_\circ^\natural$, $\Sigma$ une 
$\omega$--repr\'esentation admissible de $M_P^\natural$, 
et $\Pi={^\omega{\iota}}_{P^\natural}^{G^\natural}(\Sigma)$. Alors 
pour toute fonction $\phi\in C^\infty_{\rm c}(G^\natural)$, 
on a la formule de descente
$$
\Theta_\Pi(\phi)=\Theta_\Sigma({^\omega\phi_{P^\natural,K_\circ}}).
$$
\end{montheo}

\begin{proof}
Notons 
$\Sigma'$ la $\omega$--repr\'esentation $\delta_{P^\natural}^{1/2}(\Sigma\otimes 1)$ de $P^\natural=M^\natural_P\cdot U_P$, 
et soit $\Theta_{\Sigma'}$ le caract\`ere de $\Sigma'$ d\'efini gr\^ace \`a la mesure de Haar 
\`a gauche $d_l\gamma_P$ sur 
$P^\natural$ comme en \ref{caractres tordus (bis)}. Soit une fonction $\phi\in C^\infty_{\rm c}(G)$. 
Puisque $G=PK_\circ = K_\circ P$ et 
$${\rm vol}(P\cap K_\circ, \Delta_P^{-1}(p)d_lp)=
{\rm vol}(P\cap K_\circ, d_lp)=1,$$ 
d'apr\`es le corollaire de \ref{caractres des induites compactes}, on a
$$
\Theta_\Pi(\phi)= \Theta_{\Sigma'}(\phi_{K_\circ}\vert_{P^\natural}),
$$
o\`u
$$
\phi_{K_0}(\gamma)= \int_{K_\circ}\omega(k)\phi(k^{-1}\cdot\gamma\cdot k)dk\quad 
(\gamma\in G^\natural).
$$
Or posant $d\delta = d\gamma_{M_P}$, on a
\begin{align*}
\Sigma'(\phi_{K_\circ}\vert_{P^\natural}d_l\gamma_P)&=\int\!\!\!\!\int_{M_P^\natural\times U_P}\phi_{K_\circ}(\delta\cdot u)
\Sigma'(\delta\cdot u)d\delta du\\
&=\int_{M_P^\natural}\delta_{P^\natural}^{1/2}(\delta)\left\{\int_{U_P}\phi_{K_\circ}(\delta\cdot u)du\right\}\Sigma(\delta)d\delta\\
& = \int_{M_P^\natural}{^\omega\phi_{P^\natural,K_\circ}}(\delta)\Sigma(\delta)d\delta.
\end{align*}
D'o\`u le th\'eor\`eme.
\end{proof}

\begin{marema}
{\rm 
Si de plus la repr\'esentation $\Sigma^\circ$ est de type fini (i.e. 
de longueur finie, puisqu'elle est admissible), alors $\Pi^\circ$ l'est aussi, et gr\^ace \`a 
la formule d'int\'egration de H.Weyl \'etablie au ch.~5, on peut 
\'ecrire la formule de descente du thorme en termes de fonctions caract\`eres 
(voir \ref{formule d'intgration de Weyl}, corollaire 3).\hfill $\blacksquare$}
\end{marema}

\subsection{Restriction de Jacquet et caract\`eres}\label{restriction de Jacquet et caractres} 
Pour $P\in \ES{P}_\circ$, on note
$$
r_G^P:\frak{R}(G)\rightarrow \frak{R}(M_P)
$$
le foncteur {\it restriction de Jacquet normalis\'ee}. On rappelle la d\'efinition: pour toute repr\'esentation 
lisse $(\pi,V)$ de $G$, on note $V(U_P)$ le sous--espace de $V$ engendr\'e par les vecteurs 
$\pi(u)(v)-v$ pour $u\in U_P$ et $v\in V$, et l'on pose $V_P=V/V(U_P)$. L'espace 
$V(U_P)$ est $P$--stable, et l'on note $(\pi_P,V_P)$  la repr\'esentation (lisse) de $M_P$ 
d\'eduite de $\pi$ par restriction et passage aux quotients. Enfin on pose 
$r_G^P(\pi)= \delta_P^{-1/2}\otimes \pi_P$.

Pour $P\in P_\circ^\natural$, on d\'efinit comme suit un foncteur
$$
{^\omega{r}}_{G^\natural}^{P^\natural}:\frak{R}(G^\natural,\omega)\rightarrow \frak{R}(M_P^\natural,\omega).
$$
Soit $(\Pi,V)$ une $\omega$--repr\'esentation lisse de $G^\natural$. Posons $\pi=\Pi^\circ$. 
Puisque $\omega\vert_{U_P}=1$ (lemme de \ref{induction parabolique et caractres}), pour $\gamma\in P^\natural$, $u\in U_P$ et $v\in V$, on a
$$
\Pi(\gamma)(\pi(u)(v)-v)= \Pi({\rm Int}_{G^\natural}(\gamma)(u))(\Pi(\gamma)(v))-\Pi(\gamma)(v)\in V(U_P).
$$
Par restriction et passage au quotient, $\Pi$ induit donc une application 
$\Pi_{P^\natural}:P^\natural\rightarrow {\rm Aut}_\Bbb{C}(V_P)$, qui se factorise \`a travers 
$M^\natural_P$. Pour 
$\delta\in M_P^\natural$, $m,\,m'\in M_P$ et $v\in V$, on a
\begin{align*}
\Pi_{P^\natural}(m\cdot \delta\cdot m')(v+V(U_P))&= \Pi(m\cdot\delta\cdot m')(v)+V(U_P)\\
&= \omega(m')\pi(m)\circ \Pi(\delta)\circ \pi(m')(v)+V(U_P)\\
&=\omega(m')\pi_P(m)\circ \Pi_{P^\natural}(\delta)\circ \pi_P(m')(v+V(U_P)).
\end{align*}
En d'autres termes, $\Pi_{P^\natural}$ est une $\omega$--repr\'esentation lisse de $M_P^\natural$. 
D'apr\`es la relation $(*)$ de \ref{sous-espaces paraboliques} et le lemme de \ref{modules}, 
${\delta}_{P^\natural}^{-1/2}\Pi_{P^\natural}$ est encore une $\omega$--repr\'esentation lisse 
de $M_P^\natural$, telle que $({\delta}_{P^\natural}^{-1/2}\Pi_{P^\natural})^\circ = \delta_P^{-1/2}\pi_P$. On 
peut donc poser
$${^\omega{r}}_{G^\natural}^{P^\natural}(\Pi)= 
{\delta}_{P^\natural}^{-1/2}\Pi_{P^\natural}.
$$

Par construction, les foncteurs ${^\omega{r}}_{G^\natural}^{P^\natural}:\frak{R}(G^\natural,\omega)\rightarrow 
\frak{R}(M_P^\natural,\omega)$ et 
$r_G^P:\frak{R}(G)\rightarrow \frak{R}(M_P)$ commutent aux foncteurs d'oubli: pour toute 
$\omega$--repr\'esentation lisse $\Sigma$ de $M_P^\natural$, on a
$$
{^\omega{r}}_{G^\natural}^{P^\natural}(\Pi)^\circ = r_G^P(\Pi^\circ).
$$

\begin{marema}
{\rm 
Pour $P^\natural\in \ES{P}_\circ^\natural$, 
le foncteur ${^\omega{r}}_{G^\natural}^{P^\natural}$ est un adjoint \`a gauche du foncteur ${^\omega\iota}_{P^\natural}^{G^\natural}$: 
pour toute $\omega$--repr\'esentation lisse $\Sigma$ de $M_P^\natural$ et toute 
$\omega$--repr\'esentation lisse $\Pi$ de $G^\natural$, on a un $\Bbb{C}$--isomorphisme 
canonique
$$
{\rm Hom}_{G^\natural}(\Pi,{^\omega\iota}_{P^\natural}^{G^\natural}(\Sigma))\simeq 
{\rm Hom}_{M_P^\natural}({^\omega{r}}_{G^\natural}^{P^\natural}(\Pi),\Sigma),
$$
fonctoriel en $\Pi$ et en $\Sigma$. Pr\'ecis\'ement, posant 
$\pi=\Pi^\circ$ et $\sigma=\Sigma^\circ$, il se d\'eduit par restriction du 
$\Bbb{C}$--isomorphisme canonique (fonctoriel en $\pi$ et en $\sigma$)
$$
{\rm Hom}_G(\pi,\iota_P^G(\sigma))\simeq {\rm Hom}_{M_P}(r_G^P(\pi),\sigma)\eqno{\blacksquare}
$$}
\end{marema}

Notons $\frak{F}$ le $\Bbb{C}$--espace vectoriel des fonctions 
$f:\Bbb{Z}\rightarrow \Bbb{C}$, et $\tau$ le $\Bbb{C}$--automorphisme de 
$\frak{F}$ d\'efini par $\tau(f)(n)=f(n+1)$ ($n\in \Bbb{Z}$). 
Une fonction $f\in \frak{F}$ est dite {\it $\tau$--finie} si les $\tau^i(f)$, $i\in \Bbb{Z}$, 
engendrent un sous--espace vectoriel de $\frak{F}$ 
de dimension finie. Une fonction $f\in \frak{F}$ est $\tau$--finie 
si et seulement s'il existe un polyn\^ome $P(t)\in \Bbb{C}[t]$, $P(0)\neq 0$, tel que $P(\tau)(f)=0$; o $P(\tau)$ 
est le $\Bbb{C}$--endomorphisme de $\frak{F}$ donn par $P(\tau)=\sum_{i=0}^ka_i\tau^i$ si $P(t)=\sum_{i=0}^ka_it^i$. Les 
fonctions $\tau$--finies $f:\Bbb{Z}\rightarrow \Bbb{C}$ forment un sous--espace vectoriel de $\mathfrak{F}$.
 
\begin{monlem1}
Soit $f_1,\,f_2\in \frak{F}$ deux fonctions $\tau$--finies. Supposons 
qu'il existe un entier $n_0\geq 1$ tel que pour tout $n\in \Bbb{Z}_{\geq n_0}$, on a $f_1(n)=f_2(n)$. 
Alors $f_1=f_2$.
\end{monlem1}

\begin{proof}
Quitte  remplacer $f_i$ par $\tau^{n_0}(f_i)$, on peut supposer $n_0=0$. Il s'agit de montrer qu'une 
fonction $\tau$--finie $f\in \mathfrak{F}$ telle que $f(n)=0$ pour tout entier $n\geq 0$, ne peut \^etre que la fonction nulle. 
Supposons par l'absurde que $f\neq 0$, et soit $j$ le plus petit entier $>1$ tel que $f(-j)\neq 0$. Pour $k\in {\Bbb Z}$, 
posons $f_k= \tau^{-j-k}(f)$. On a $f_k(k)=f(-j)\neq 0$ et $f_k(k+n)=0$ pour tout entier $n\geq 1$. Les fonctions $f_k$ sont 
linairement indpendantes, ce qui contredit le fait que $f$ est $\tau$--finie. Donc $f=0$.
\end{proof}

\begin{exemple}
{\rm Soit $u\in {\rm End}_\Bbb{C}(X)$ pour un $\Bbb{C}$--espace vectoriel 
$X$ de dimension finie, et soit $\lambda_1,\ldots ,\lambda_r\in \Bbb{C}^\times$ 
($\lambda_i\neq \lambda_j$ pour $i\neq j$)
les valeurs propres non nulles de $u$. Pour $i=1,\ldots ,r$, notons $m_i$ la multiplicit\'e de 
$\lambda_i$ dans le polyn\^ome caract\'eristique de $u$. 
Alors la fonction $f:\Bbb{Z}_{\geq 1}\rightarrow \Bbb{C}$ 
d\'efinie par
$$
f(n)={\rm tr}(u^n)=m_1\lambda_1^n+\cdots + m_r\lambda_r^n,
$$
se prolonge de manire unique en une fonction $\tau$--finie $f\in \frak{F}$.\hfill$\blacksquare$}
\end{exemple}

Pour $G^\natural=G$, $\omega=1$ et $g$ r\'egulier, le rsultat suivant est d\^u 
\`a Casselman \cite[theo. 5.2]{Ca2}.

\begin{montheo}Soit $P^\natural\in \ES{P}_\circ^\natural$, 
$\Pi$ une $\omega$--repr\'esentation admissible de $G^\natural$, 
et $\Sigma ={^\omega{r}_{G^\natural}^{P^\natural}}(\Pi)$. On 
suppose que $\Pi^\circ$ est de type fini (i.e. de longueur finie). Pour tout 
$\gamma\in M_P^\natural\cap G_{\rm qr}^\natural$ tel que 
$P_{[\gamma]}=P$ et $A_{[\gamma]}=A_P$, on a l'\'egalit\'e
$$
\delta_{P^\natural}^{-1/ 2}(\gamma)\Theta_\Pi(\gamma)=\Theta_\Sigma(\gamma).
$$
\end{montheo}

\begin{proof}
Posons $M=M_P$ et $M^\natural=M_P^\natural$. 
On a clairement l'inclusion 
$$M^\natural \cap G^\natural_{\rm qr}\subset M^\natural_{\rm qr}, $$
par cons\'equent l'\'enonc\'e a un sens. Posons $\pi=\Pi^\circ$. 

Soit un \'el\'ement $\gamma \in M^\natural
\cap G^\natural_{\rm qr}$ tel que $P_{[\gamma]} =P$ et $A_{[\gamma]} =A_P$. Posons 
$\tau={\rm Int}_{G^\natural}(\gamma)$. D'apr\`es la remarque 4 de \ref{la paire parabolique associe  gamma}, il 
existe un sous--groupe ouvert compact $K$ de $G$ tel que, posant 
$K^-=K\cap U_{P^-}$, $K^0=K\cap M$, $K^+=K\cap U_P$, on a:
\begin{itemize}
\item $K=K^-K^0K^+ = K^+K^0K^-$;
\item $\tau^{-1}(K^-)\subset K^-$, $\tau(K^0)=K^0$, $\tau(K^+)\subset K^+$;
\item $\omega\vert_K=1$;
\item $K\cdot \gamma\cdot K\subset G^\natural_{\rm qr}$;
\item $\Theta_\Pi\vert_{K\cdot \gamma\cdot  K}=\Theta_\Pi(\gamma)$;
\item $\Theta_{\Pi_{P^\natural}}\vert_{K^0\cdot \gamma}=\Theta_{\Pi_{P^\natural}}(\gamma)$.
\end{itemize}
Pour tout sous--groupe compact $J$ de $G$, on note $e_J$ la mesure de 
Haar normalis\'ee sur $J$ (i.e. telle que ${\rm vol}(J,e_J)=1$), identifi\'ee \`a une 
distribution \`a support compact sur $G$. Puisque $K=K^-K^0K^+=K^+K^0K^-$, 
on a les \'egalit\'es (dans l'alg\`ebre des 
distributions \`a support compact sur $G$)
$$
e_K = e_{K^-}* e_{K^0}*e_{K^+}= e_{K^+}* e_{K^0}*e_{K^-}.
$$

Pour $\delta\in G^\natural$, on note $\phi^K_\delta$ 
la fonction caract\'eristique de $K\cdot \delta \cdot K\subset G^\natural $ divis\'e par 
${\rm vol}(K\cdot \delta\cdot K,d\gamma_G)$, o\`u (rappel) 
$d\gamma_G = \delta_1\cdot dg$ (cf. \ref{module d'un espace tordu}). On a donc
$$
\Pi(\phi^K_\delta)= \pi(e_K)\circ \Pi(\delta)\circ \pi(e_K).
$$
Pour $n\in \Bbb{Z}_{\geq 1}$, puisque $\tau(K^0K^+)\subset K^0K^+$ et $\tau^{-n}(K^-)\subset K^-$, on a
\begin{align*}
 &\Pi(\phi^K_\gamma)\circ \pi(e_K)\circ \Pi(\gamma)^n\circ \pi(e_K)\\
 &\h{7}=\pi(e_K)\circ \Pi(\gamma)\circ
 \pi(e_K)\circ \Pi(\gamma)^n\circ \pi(e_K)\\
 &\h{7}=\pi(e_K)\circ \Pi(\gamma)\circ \pi(e_{K^0K^+})\circ \pi(e_{K^-})\circ \Pi(\gamma)^n\circ \pi(e_K)\\
 &\h{7}= \pi(e_K)\circ \pi(e_{\tau(K^0K^+)})\circ \Pi(\gamma)\circ \Pi(\gamma)^n\circ 
 \pi(e_{\tau^{-n}(e_{K^-})})\circ \pi(e_K)\\
 &\h{7}=\pi(e_K)\circ \Pi(\gamma)^{n+1}\circ \pi(e_K).
\end{align*}
Par r\'ecurrence sur $n$, on a  donc
$$
\Pi(\phi^K_\gamma)^n = \pi(e_K)\circ \Pi(\gamma)^n\circ \pi(e_K)\quad (n\in \Bbb{Z}_{\geq 1}).
$$
Pour $n\in \Bbb{Z}_{\geq 1}$, notons $V^{K,n}_\gamma$ le sous--espace 
$$\Pi(\phi^K_\gamma)^n(V^K)=\pi(e_K)\circ \Pi(\gamma)^n(V^K)$$ de $V^K=\pi(e_K)(V)$. 
Pour 
$m,\,n\in \Bbb{Z}_{\geq 1}$, on a donc
$$
V^{K,m+n}_\gamma = \Pi(\phi^K_\gamma)^m(V^{K,n}_\gamma)\subset 
\Pi(\phi^K_\gamma)^m(V^K)=
V^{K,m}_\gamma.
$$

Rappelons que l'espace $V(U_P)$ co\"{\i}ncide avec l'ensemble des $v\in V$ tels que 
$$\int_{\Omega_v}\pi(u)(v)du=0$$ pour un sous--groupe ouvert compact $\Omega_v$ de $U_P$. 
Puisque la reprsentation $\pi$ est admissible, l'espace $V(U_P)\cap V^K$ est de dimension finie, 
et il existe un sous--groupe 
ouvert compact $\Omega$ de $U_P$ tel que  pour tout 
$v\in V(U_P)\cap V^K$, on a $\int_\Omega\pi(u)(v)du=0$. 
Quitte \`a remplacer $\Omega$ par un sous--groupe plus 
gros, on peut supposer que $K^+\subset \Omega$. Choisissons un entier 
$n_0\geq 1$ tel que $\tau^{n_0}(\Omega)\subset K^+$.

Notons $p:V\rightarrow V_P$ la projection canonique, et posons 
$$V_P^{K^0}=\pi_P(e_{K^0})(V_P)= r_G^P(e_{K^0})(V_P).$$

\begin{monlem2}
Soit $m\in \Bbb{Z}_{\geq 1}$ et $n\in \Bbb{Z}_{\geq n_0}$.
\begin{enumerate}
\item[(1)]Pour $v\in V^K$, on a 
$$\Pi(\phi^K_\gamma)^m(v)=\pi(e_{K^+})\circ \Pi(\gamma)^m(v),\quad p\circ \Pi(\phi^K_\gamma)^m(v)= \Pi_{P^\natural}(\gamma)^m\circ p(v).$$
\item[(2)]On a 
$\Pi(\phi^{K,n}_\gamma)(V(U_P)\cap V^K)=0$.
\item[(3)]L'application
$p:V\rightarrow V_P$ induit par restriction un 
$\Bbb{C}$--isomorphisme $V^{K,n}_\gamma \rightarrow V_P^{K^0}$.
\item[(4)]On a $\Pi(\phi^K_\gamma)^m(V^{K,n}_\gamma)=
V^{K,n}_\gamma$.
\end{enumerate}
\end{monlem2}

\begin{proof}
Montrons (1). Soit $v\in V^K$. Puisque 
$\tau^{-m}(K^0K^-)\subset K$, on a
\begin{align*}
\Pi(\phi^K_\gamma)^m(v) &=\pi(e_K)\circ \Pi(\gamma)^m(v)\\
&=\pi(e_{K^+})\circ \pi(e_{K^0K^-})\circ \Pi(\gamma)^m(v)\\
&=\pi(e_{K^+})\circ \Pi(\gamma)^m\circ \pi(e_{\tau^{-m}(K^0K^-)})(v)\\
&= \pi(e_{K^+})\circ \Pi(\gamma)^m(v).
\end{align*}
Or $\pi(e_{K^+})\circ \Pi(\gamma)^m(v)$ est contenu dans 
$\Pi(\gamma)^m(v) +V(U_P)$, d'o\`u le point (1).

Montrons (2). Soit
$v\in V(U_P)\cap V^K$. D'apr\`es (1), on a
\begin{align*}
\Pi(\phi^{K,n}_\gamma)(v)&=
\pi(e_{K^+})\circ \Pi(\gamma)^n(v)\\
&= {\rm vol}(K^+\!,du)^{-1}\int_{K^+}\pi(u)\circ \Pi(\gamma)^n(v)du\\
&={\rm vol}(K^+\!,du)^{-1}\Pi(\gamma)^n\circ 
\int_{K^+}\pi(\tau^{-n}(u))(v)du\\
&=0
\end{align*}
car $\tau^{-n}(K^+)\supset \tau^{-n_0}(K^+)
\supset \Omega$.

Montrons (3). Fixons un entier $n\geq n_0$ et posons $W=V^{K,n}_\gamma$. 
Montrons la surjectivit\'e. Soit $\bar{v}\in V_P^{K^0}$. Puisque 
${\rm Int}_{M^\natural}(\gamma)^{-n}(K^0)={\rm Int}_{G^\natural}(\gamma)^{-n}(K^0)= K^0$, 
on a $\Pi_{P^\natural}(\gamma)^{-n}(\bar{v})\in V_P^{K^0}$. D'apr\`es le 
\og premier lemme de Jacquet\fg{}\footnote{Dans le cas admissible, le rsultat est prouv 
dans \cite[theo.~3.3.3]{Ca1} --- cf. aussi \cite[3.2]{Be}. Il a ensuite 
t tendu au cas non admissible par Bernstein (loc.~cit).}, 
la projection 
canonique $p:V\rightarrow V_P$ induit une application surjective 
$V^K\rightarrow V_P^{K^0}$. 
Choisissons un $v'\in V^K$ tel que $p(v')= \Pi_{P^\natural}(\gamma)^{-n}(\bar{v}) $, et 
posons $v=\Pi(f^K_\gamma)^{n}(v')\in W$. D'apr\`es (1), on a
$p(v)= \Pi_{P^\natural}(\gamma)^n\circ p(v')= \bar{v}$. 
Montrons l'injectivit\'e. Soit $w\in V(U_P)\cap W$. \'Ecrivons 
$w= \pi(e_K)\circ \Pi(\gamma)^n(v'')$ avec $v''\in V^K$. 
D'apr\`es (1), on a $w=\pi(e_{K^+})\circ \Pi(\gamma)^n(v'')$, 
et puisque $w\in V(U_P)\cap V^K$ et 
$K^+\subset \Omega$, on a
\begin{align*}
0& =\int_\Omega \pi(u)(w)du\\
& = \int_\Omega \pi(u) \pi(e_{K^+})\circ \Pi(\gamma)^n(v'')du\\
& =\int_\Omega \pi(u)\circ \Pi(\gamma)^{n}(v'')du\\
& = \Pi(\gamma)^n\circ \int_\Omega\pi(\tau^{-n}(u))(v'')du.
\end{align*}
Donc $v''\in V(U_P)$. D'apr\`es (2), on a donc $w=0$.

Montrons (4). On a
$$
\Pi(\phi^K_\gamma)^m(V^{K,n}_\gamma)=V^{K,n+m}_\gamma \subset V^{K,n}_\gamma.
$$
Or d'apr\`es (3), on a $\dim_\Bbb{C}(V^{K,n+m}_\gamma)=\dim_\Bbb{C}(V_P^{K^0})
=\dim_\Bbb{C}(V^{K,n}_\gamma)$, par cons\'equent l'inclusion ci-dessus est une 
\'egalit\'e.
\end{proof}

Posons $W=V^{K,n_0}_\gamma$ et $\overline{W}=V_P^{K^0}$, et pour $\delta\in M^\natural$, notons 
$\phi^{K^0}_\delta$ la fonction caract\'eristique de 
$K^0\cdot \gamma \cdot K^0$ divis\'ee par ${\rm vol}(K^0\cdot \delta \cdot K^0\!,d\gamma_M)$. 
D'apr\`es le lemme 2, pour
$n\in \Bbb{Z}_{\geq n_0}$, on a
$$
{\rm tr}(\Pi(\phi^K_\gamma)^n)= {\rm tr}(\Pi(\phi^K_\gamma)^n\vert_{W})=
{\rm tr}(\Pi_{P^\natural}(\gamma)^n\vert_{\overline{W}})= {\rm tr}(\Pi_{P^\natural}(\phi^{K^0}_\gamma)^n).
$$
Gr\^ace au lemme 1, on en d\'eduit que
$$
\Theta_\Pi(\gamma)=\Theta_\Pi(\phi^K_\gamma)
=\Theta_{\Pi_{P^\natural}}(\phi^{K^0}_\gamma)=\Theta_{\Pi_{P^\natural}}(\gamma).
$$
D'o\`u le th\'eor\`eme puisque 
$\Theta_{\Pi_{P^\natural}}(\gamma)=\delta_{P^\natural}^{1/2}(\gamma)\Theta_\Sigma(\gamma)$.
\end{proof}

Pour $G^\natural=G$, $\omega=1$ et $g$ r\'egulier, le corollaire suivant est d\^u \`a Deligne \cite{D}.

\begin{moncoro}
Soit $\Pi$ une $\omega$--repr\'esentation admissible 
de $G^\natural$. On suppose que $\Pi^\circ$ est de type fini, et que pour 
tout $P^\natural\in \ES{P}_\circ^\natural\smallsetminus \{G^\natural\}$, on a 
${^\omega{r}_{G^\natural}^{P^\natural}}(\Pi)=0$. Alors pour tout 
$\gamma \in G_{\rm qr}^\natural$ tel que 
$P_{[\gamma]} \neq G$, on a $\Theta_\Pi(\gamma)=0$.
\end{moncoro}

\begin{proof}
Soit $\gamma \in G_{\rm qr}^\natural$. Choisissons $x\in G$ tel que 
$\gamma'= x^{-1}\cdot \gamma \cdot x$ est en position standard. Puisque 
$P_{[\gamma]} = xP_{[\gamma']}x^{-1}$, si $P_{[\gamma]}\neq G$, alors $P=P_{[\gamma']}\neq G$ et 
d'aprs le lemme 2, on a $\Theta_\Pi(\gamma')= 0$. On conclut gr\^ace \`a la relation $(**)$ de \ref{les oprateurs T}.
\end{proof}

\subsection{Commentaire}\label{commentaire} On peut bien s\^ur traduire ces r\'esultats 
dans le langage classique des caract\`eres $(\theta,\omega)$--tordus de $G$. 
Un \'el\'ement $g\in G$ est dit:
\begin{itemize}
\item {\it $\theta$--r\'egulier} si 
l'\'el\'ement $g\theta\in G^\natural=G\theta$ est r\'egulier;
\item {\it $\theta$--quasi--r\'egulier} si $g\theta$ est quasi--r\'egulier.
\end{itemize}
Notons $G_{\theta-{\rm reg}}$ et $G_{\theta-{\rm qr}}$ les sous--ensembles de 
$G$ form\'es des \'el\'ements qui sont respectivement $\theta$--r\'eguliers et 
$\theta$--quasi--r\'eguliers. D'apr\`es la proposition de \ref{lments rguliers et quasi-rguliers}, on a l'inclusion
$$
G_{\theta-{\rm reg}}\subset G_{\theta-{\rm qr}};
$$
et ces deux ensembles sont ouverts denses dans $G$ (corollaire de \ref{lments rguliers et quasi-rguliers}).

Soit $(\pi,V)$ une repr\'esentation admissible de $G$ telle que $\omega\pi\simeq \pi^\theta$, et soit 
$A\in {\rm Isom}_G(\omega\pi,\pi^\theta)$. On suppose que $\pi$ est de type fini. 
Soit $K$ un sous--groupe 
ouvert compact de $G$. Pour $g\in G$, on note ${\rm T}_g={\rm T}_{(\pi,A),g}^K\in {\rm End}_\Bbb{C}(V)$ 
l'op\'erateur d\'efini par
\begin{align*}
{\rm T}_g &= {\rm vol}(K,dk)^{-1}\int_K\omega(k^{-1})\pi(k^{-1}g\theta(k))\circ A\h{0.2}dk\\
& ={\rm vol}(K,dk)^{-1}\int_K\pi(k^{-1})\circ \pi(g)\circ A\circ \pi(k)dk.\end{align*}
D'apr\`es le thorme de \ref{les oprateurs T}, pour $g\in G_{\theta-{\rm qr}}$, on a ${\rm T}_g\in {\rm End}_0(V)$, 
et l'application
$$G_{\theta-{\rm qr}}\rightarrow {\rm End}_0(V),\,g\mapsto {\rm T}_g$$ est localement 
constante. D'apr\`es le corollaire de \ref{les oprateurs T}, le caract\`ere tordu $\Theta_\pi^A$ est repr\'esent\'e sur 
$G_{\theta-{\rm qr}}$ par la fonction localement constante $g\mapsto \Theta_\pi^A(g)=
 {\rm tr}({\rm T}_g)$. 
En d'autres termes, pour toute fonction $f\in C^\infty_{\rm c}(G_{\theta-{\rm qr}})$, on a 
l'\'egalit\'e
$$
\Theta_\pi^A(f)=\int_{G_{\theta-{\rm qr}}}f(g)\Theta_\pi^A(g)dg.
$$
Puisque pour $f\in C^\infty_{\rm c}(G)$ et $x\in G$ on a 
$\Theta_\pi^A({^{x,\theta}f})=\omega(x^{-1})\Theta_\pi^A(f)$ (relation $(**)$ de \ref{caractres tordus (bis)}), on a l'\'egalit\'e
$$
\Theta_\pi^A(x^{-1}g\theta(x))=\omega(x)\Theta_\pi^A(g)\quad (g\in G_{\theta-{\rm qr}},\,
x\in G).
$$

La traduction des num\'eros \ref{induction parabolique et caractres} et \ref{restriction de Jacquet et caractres} est laiss\'ee au lecteur.

\section{S\'eries discr\`etes et repr\'esentations cuspidales}
 
Continuons avec les notations du ch.~5.

\subsection{Caract\`eres des reprsentations irrductibles essentiellement de carr intgrable}\label{caractres des sries ed}
Notons $Z=Z(G)$ le 
centre de $G$, et fixons une mesure de Haar $dz$ sur $Z$. Soit $d\bar{g}={dg\over dz}$ 
la mesure quotient sur le groupe $G/Z$.

Soit $(\Pi,V)$ une $\omega$--repr\'esentation lisse 
de $G^\natural$. Posons $\pi=\Pi^\circ$ et 
$A=\Pi(\delta_1)\in {\rm Isom}_G(\omega\pi,\pi^\theta)$. Pour $(v,\check{v})\in 
V\times \check{V}$, on 
note $\pi_{v,\check{v}}:G\rightarrow \Bbb{C}$ le coefficient de $\pi$ d\'efini par
$$
\pi_{v,\check{v}}(g)=
\langle \pi(g)(v),\check{v}\rangle \quad (g\in G),
$$
et $\Pi_{v,\check{v}}$ le coefficient de $\Pi$ d\'efini par
$$
\Pi_{v,\check{v}}(\gamma)=
\langle \Pi(\gamma)(v),\check{v}\rangle \quad (\gamma\in G^\natural).
$$
Pour $\gamma=g\cdot \delta_1$, on a donc $\Pi_{v,\check{v}}(\gamma)=\pi_{A(v),\check{v}}(g)$. 
Soit $\ES{A}(\Pi)$ le $\Bbb{C}$--espace vectoriel engendr\'e par les 
coefficients de $\Pi$. Pour $\varphi\in \ES{A}(\Pi)$, puisque $\Pi(z\cdot \gamma\cdot z)= \omega(z)\Pi(\gamma)$, 
on a
$$
\varphi(z^{-1}\cdot \gamma\cdot z)=\omega(z)\varphi(g)\quad (z\in Z,\,\gamma\in G^\natural).\leqno{(*)}
$$
Soit aussi $\ES{A}(\pi)$ le $\Bbb{C}$--espace vectoriel engendr\'e par les coefficients de $\pi$, et soit 
$$\ES{A}(\Pi)\rightarrow \ES{A}(\pi),\,\varphi\mapsto \varphi^\circ$$ l'application lin\'eaire d\'efinie par 
$\Pi_{v,\check{v}}^\circ= \pi_{v,\check{v}}$ pour tout $(v,\check{v})\in V\times \check{V}$.

On suppose $\pi$ irr\'eductible, et {\it essentiellement 
de carr\'e int\'egrable modulo $Z$}; i.e. il existe un caract\`ere $\psi$ de $G$ tel que 
les coefficients de $\pi\otimes \psi$ sont de carr\'e int\'egrable modulo $Z$. 
Soit $d(\pi)=d(\pi,d\bar{g})>0$ le degr\'e formel de $\pi$ d\'efini par $d\bar{g}$: 
pour tous $(v,\check{v}),\,(v',\check{v}')\in V\times \check{V}$, on a 
$$
\int_{G/Z} \pi_{v,\check{v}}(g) \pi_{v'\!,\check{v}'}(g^{-1})d\bar{g}
= d(\pi)^{-1}\langle v,\check{v}'\rangle \langle v',\check{v}\rangle.
$$

Fixons un sous--groupe ouvert compact $K$ de $G$. Pour 
$\varphi\in \ES{A}(\Pi)$ et $\gamma\in G^\natural$, on note 
$\Theta_{\varphi,\gamma}^K: G\rightarrow \Bbb{C}$ la fonction d\'efinie par
$$
\Theta_{\varphi,\gamma}^K(g)={\rm vol}(K,dk)^{-1}\int_K\omega(k^{-1})\varphi(g^{-1}k^{-1}\cdot 
\gamma\cdot kg)dk.
$$
D'apr\`es $(*)$, on a l'galit
$$
\Theta^K_{\varphi,\gamma}(zg)=\omega(z)\Theta^K_{\varphi,\gamma}(g)\quad (z\in Z,\,g\in G).\leqno{(**)}
$$
Pour $\theta={\rm id}_G$, 
$\omega=1$ et $\gamma\in G$ (semisimple) r\'egulier, le th\'eor\`eme suivant est d\^u \`a 
Rader et Silberger \cite[theo. 2]{RS} (cf. aussi \cite[theo.~A.12]{BH1}).  

\begin{montheo}
Pour toute fonction $\varphi\in \ES{A}(\Pi)$, on a
$$
\varphi^\circ(1) \Theta_\Pi(\gamma)=d(\pi)\int_{G/Z}\omega(g^{-1})\Theta_{\varphi,\gamma}^K(g)
 d\bar{g}\quad (\gamma \in G^\natural_{\rm qr});
$$
l'int\'egrale converge absolument et uniform\'ement\footnote{Prcisment, pour toute partie compacte $C$ de $G^\natural_{\rm qr}$, 
les fonctions $\omega^{-1}\Theta_{\varphi,\gamma}^K$ sur $G/Z$ ($\gamma\in C$) vrifient le ``critre de M--test'' de Weierstrass.} sur les 
parties compactes de $G^\natural_{\rm qr}$. Si de plus $\pi$ est cuspidale, 
alors pour toute fonction $\varphi\in \ES{A}(\Pi)$ et 
toute partie compacte $C$ de $G^\natural_{\rm qr}$, il existe une partie 
$\Omega=\Omega(\varphi,C)$ de $G$ compacte modulo $Z$ 
telle que
$$\Theta^K_{\varphi,\gamma}(g)=0\quad  (\gamma\in C,\,g\in G \smallsetminus \Omega).
$$
\end{montheo}

\begin{proof}
Puisque pour 
$\gamma\in G^\natural_{\rm qr}$, on a
$$\Theta_\Pi(g^{-1}\cdot\gamma\cdot g)=\omega(g)\Theta_\Pi(\gamma)\quad (g\in G),
$$ 
le th\'eor\`eme ne d\'epend pas du choix de $K$. On peut donc, comme dans la d\'emonstration 
du th\'eor\`eme de \ref{les oprateurs T}, supposer que 
$K$ est un sous--groupe ouvert distingu\'e de $K_\circ$ tel que $\omega\vert_K =1$. 
D'autre part, il suffit de traiter les fonctions $\varphi\in \ES{A}(\Pi)$ qui sont des 
coefficients de $\Pi$. Soit donc $\varphi=\Pi_{v,\check{v}}$ 
pour un couple $(v,\check{v})\in V\times \check{V}$. Soit aussi une partie compacte 
$C$ de $G^\natural_{\rm qr}$.

Pour $\gamma\in G^\natural$, reprenons l'op\'erateur ${\rm T}_\gamma=
{\rm T}^K_{\Pi,\gamma}\in {\rm End}_0(V)$ d\'efini en \ref{les oprateurs T}. D'apr\`es la d\'emons\-tration 
du thorme de loc.~cit., pour chaque $\gamma\in C$, il existe un voisinage $\ES{V}(\gamma)$ de $\gamma$ dans 
$G^\natural_{\rm qr}$ et un sous--groupe ouvert distingu\'e 
$K_\gamma$ de $K$, tels que
$$
{\rm T}_{x\cdot \delta\cdot y}={\rm T}_\gamma\quad (\delta\in \ES{V}(\gamma);\, x,\,y\in K_\gamma).
$$
On peut recouvrir $C$ par un nombre fini de tels voisinages $\ES{V}_\gamma$. 
On en d\'eduit l'existence d'un sous--groupe ouvert compact $J$ de $G$ 
tel que 
$$
{\rm T}_{x\cdot \gamma\cdot y}={\rm T}_\gamma\quad (\gamma\in C;\, x,\,y\in J).
$$
Par cons\'equent, pour $\gamma\in C$, l'op\'erateur ${\rm T}_\gamma$ appartient au sous--espace 
$V^J\times \check{V}^J$ de ${\rm End}_0(V)$. Posons $W=V^J$ et $W^*=\check{V}^J\;(=
{\rm Hom}_\Bbb{C}(W,\Bbb{C}))$. Choisissons une 
base $\{w_1,\ldots ,w_r\}$ de $W$ sur $\Bbb{C}$, et notons  $\{w_1^*,
\ldots ,w_r^*\}$ la base de $W^*$ duale de $\{w_1,\ldots ,w_r\}$. 
Pour $\gamma\in C$ et $v'\in V$, on a donc
$$
{\rm T}_\gamma(v')= \sum_{i=1}^r\sum_{j=1}^r \langle {\rm T}_\gamma(w_i),w_j^*\rangle \langle v'\!, w^*_i\rangle w_j.
$$
Par suite, pour $\gamma\in C$ et $g\in G$, on a
\begin{align*}
\Theta_{\varphi,\gamma}^K(g)
&={\rm vol}(K,dk)^{-1}\int_K \omega(k^{-1})\langle \omega(g)\pi(g^{-1})\circ \Pi(k^{-1}\cdot \gamma\cdot k)\circ \pi(g)(v), \check{v}\rangle dk\\
&=\omega(g)\langle {\rm T}_\gamma\circ \pi(g)(v),\check{\pi}(g)(\check{v})\rangle\\
&= \omega(g) 
\sum_{i=1}^r\sum_{j=1}^r \langle {\rm T}_\gamma(w_i),w_j^*\rangle \langle \pi(g)(v), w^*_i\rangle \langle w_j,\check{\pi}(g)(\check{v})\rangle\\
&= \omega(g) 
\sum_{i=1}^r\sum_{j=1}^r \langle {\rm T}_\gamma(w_i),w_j^*\rangle  
\pi_{v,w^*_i}(g)\pi_{w_j,\check{v}}(g^{-1}).
\end{align*}
Si $\pi$ est cuspidale, alors les coefficients de $\pi$ sont \`a support compact modulo $Z$, 
et d'apr\`es le calcul ci-dessus, il existe une partie 
$\Omega$ de $G$ compacte modulo $Z$ telle que pour $\gamma\in C$ et $g\in G\smallsetminus \Omega$, on a 
$\Theta_{\varphi,\gamma}^K(g)=0$.

Reprenons la d\'emonstration dans le cas g\'en\'eral. 
Puisque la fonction
$$
G^\natural_{\rm qr}\rightarrow {\rm End}_0(V),\,\gamma\mapsto {\rm T}_\gamma
$$
est localement constante (thorme de \ref{les oprateurs T}), 
il existe une constante $d_C>0$ telle que
$$
\vert \omega(g^{-1})\Theta_{\varphi,\gamma}^K(g)\vert \leq d_C
\sum_{i=1}^r\sum_{j=1}^r \vert 
\pi_{v,w^*_i}(g)\vert \vert \pi_{w_j,\check{v}}(g^{-1})\vert \quad (\gamma\in C,\,g\in G).
$$
Par suite, pour $\gamma\in C$ on a
\begin{align*}
d(\pi)\int_{G/Z}\omega(g^{-1})\Theta_{\varphi,\gamma}^K(g)
d\bar{g} & = \sum_{i=1}^r\sum_{j=1}^r  \langle {\rm T}_\gamma(w_i),w_j^*\rangle 
d(\pi)\int_{G/Z}\pi_{v,w^*_i}(g)\pi_{w_j,\check{v}}(g^{-1})d\bar{g}\\
&= \sum_{i=1}^r\sum_{j=1}^r  \langle {\rm T}_\gamma(w_i),w_j^*\rangle \langle v,\check{v}\rangle 
\langle w_i,w^*_j\rangle\\
&=\varphi^\circ(1)\sum_{i=1}^r \langle {\rm T}_\gamma(w_i),w_i^*\rangle,
\end{align*}
et
$$\sum_{i=1}^r \langle {\rm T}_\gamma(w_i),w_i\rangle={\rm tr}({\rm T}_\gamma)=\Theta_\Pi(\gamma).
$$ 

Choisissons un caract\`ere $\psi$ de $G$ tel que les coefficients de 
$\pi'=\pi\otimes \psi$ sont de carr\'e int\'egrable modulo $Z$. 
Pour $(u,\check{u})\in V\times \check{V}$, on a 
$$
\pi_{u,\check{u}}(g)= \psi(g^{-1})\pi'_{u,\check{u}}(g)\quad (g\in G),
$$
et notant 
$\vert\vert\;\vert\vert_2$ la norme ${\rm L}^2$ sur $G/Z$, on a 
$\vert\vert\pi'_{u,\check{u}}\vert\vert_2<+\infty$. Ainsi, la majoration pour 
$\vert \omega^{-1}\Theta_{\varphi,\gamma}^K\vert$ donn\'ee plus haut jointe \`a l'in\'egalit\'e 
de Schwartz, entra\^{\i}nent que
$$
\int_{G/Z}\vert \omega(g^{-1})\Theta^K_{\varphi,\gamma}(g)\vert d\bar{g}
\leq d_C\sum_{i=1}^r\sum_{j=1}^r \vert\vert\pi'_{v,w^*_i}\vert\vert_2\vert\vert\pi'_{w_j,\check{v}}\vert\vert_2\quad (\gamma\in C).
$$
D'o\`u l'assertion de convergence uniforme (d'aprs le ``critre de M--test'' de Weierstrass).
\end{proof}

\subsection{Caract\`eres des repr\'esentations irr\'eductibles cuspidales}\label{caractres des reprsentations cusp}
Soit $H$ 
un sous--groupe  ouvert (donc ferm\'e) de $G$, contenant $Z$ et compact modulo $Z$. 
Soit $H^\natural$ un $H$--espace tordu qui soit un sous--espace topologique tordu de $G^\natural$. 
On note encore $\omega$ le caract\`ere $\omega\vert_H$ de $H$. 
Soit $(\Sigma,W)$ une $\omega$--repr\'esentation lisse de $H^\natural$ telle que la 
repr\'esentation $\sigma=\Sigma^\circ$ de $H$ est irrductible. Puisque $H$ est compact modulo $Z$, 
l'espace $W$ est de dimension finie, et l'on pose $\dim(\sigma)=\dim_\Bbb{C}(W)$. 
Soit $(\Pi,V)= {^\omega{\rm ind}}_{H^\natural}^{G^\natural}(\Sigma,W)$ comme en \ref{induction compacte}, et 
posons $\pi=\Pi^\circ$. D'aprs \cite[theo.~1]{Bu}, la repr\'esentation $\pi$ est admissible si et seulement si 
elle est somme directe finie de reprsentations cuspidales. En particulier si la reprsentation $\pi$ est irrductible, 
alors elle est cuspidale et 
on peut lui appliquer le th\'eor\`eme de \ref{caractres des sries ed}.

On note $\Theta_\Sigma:G^\natural\rightarrow \Bbb{C}$ la fonction d\'efinie par
$$
\Theta_\Sigma(\gamma)= \left\{
\begin{array}{ll}{\rm tr}(\Sigma(\gamma)) &\hbox{si $\gamma\in H^\natural$}\\
0 & \hbox{sinon}
\end{array}
\right..
$$
On d\'efinit de la m\^eme mani\`ere la fonction $\Theta_\sigma:G\rightarrow \Bbb{C}$.

Fixons un sous--groupe ouvert compact $K$ de $H$. Le th\'eor\`eme suivant est une 
g\'en\'eralisation de \cite[appendix, theo. A.14]{BH1}. 

\begin{montheo}
Supposons que $\pi$ est irr\'eductible. 
\begin{enumerate}
\item[(1)] On a $d(\pi)= \dim(\sigma){\rm vol}(H/Z,d\bar{g})^{-1}$.
\item[(2)]Pour $\gamma\in G^\natural_{\rm qr}$, on a
$$
\Theta_\Pi(\gamma)={\rm vol}(K,dg)^{-1}\sum_{g\in G/H}
\omega(g^{-1})\int_K\omega(k^{-1})
\Theta_\Sigma(g^{-1}k^{-1}\cdot\gamma\cdot kg)dk;
$$
la somme est finie. Mieux: pour toute partie compacte $C$ de $G^\natural_{\rm qr}$, il 
existe un sous--ensemble fini $\Omega= \Omega(C)$ de $G/H$ tel que
$$
\int_K\omega(k^{-1})
\Theta_\Sigma(g^{-1}k^{-1}\cdot\gamma\cdot kg)dk=0 \quad (\gamma\in C,\, g\in 
(G/H)\smallsetminus \Omega).
$$
\end{enumerate}
\end{montheo}

\begin{proof}
Choisissons une base $\{w_1,\ldots ,w_n\}$ de $W$, et 
notons $\{w_1^*,\ldots ,w_n^*\}$ la base de $W^*={\rm Hom}_\Bbb{C}(W,\Bbb{C})$ duale de 
$\{w_1,\ldots ,w_n\}$. Pour $i=1,\ldots ,n$, on note $v_i\in V$ et $\check{v}_i\in {\rm ind}_H^G(W^*)$ 
les fonctions telles ${\rm Supp}(v_i)=H={\rm Supp}(\check{v}_i)$, $v_i(1)=w_i$ et $\check{v}_i(1)=w^*_i$. 
Via l'identification canonique $\check{V}= {\rm ind}_H^G(W^*)$, 
on a
$$\Theta_\Sigma(\gamma)= 
\sum_{i=1}^n\Pi_{v_i,\check{v}_i}(\gamma)\quad (\gamma\in G^\natural).
$$
En particulier, la fonction $\Theta_\Sigma$ appartient \`a l'espace 
$\ES{A}(\Pi)$ des coefficients de $\Pi$, et pour $g\in G$, on a 
$\Theta_\Sigma^\circ(g)=\Theta_\sigma(g)$.

D'apr\`es la d\'efinition de $d(\pi)$, on a
\begin{align*}
d(\pi)\int_{G/Z}\Theta_\sigma(g)\Theta_\sigma(g^{-1})d\bar{g} &=
d(\pi)\sum_{i=1}^n\sum_{j=1}^n\int_{G/Z}\pi_{v_i,\tilde{v_i}}(g)\pi_{v_j,\tilde{v_j}}(g^{-1}) d\bar{g}\\
& = \sum_{i=1}^n\sum_{j=1}^n \langle v_i, \check{v}_j\rangle 
\langle v_j,\check{v}_i\rangle\\
&=n=\dim(\sigma).
\end{align*}
D'autre part, notant $d\bar{h}$ la restriction de $d\bar{g}$ \`a $H/Z$, on a aussi
\begin{align*}
\int_{G/Z}\Theta_\sigma(g)\Theta_\sigma(g^{-1})d\bar{g} &=
\int_{H/Z}\Theta_\sigma(h)\Theta_\sigma(h^{-1})d\bar{h} \\
& =\sum_{i=1}^n\sum_{j=1}^n \int_{H/Z}\langle \sigma(h)(w_i),w_i^*\rangle
\langle \sigma(h^{-1})(w_j),w_j^*\rangle d\bar{h}\cr
&= \sum_{i=1}^n\sum_{j=1}^n{{\rm vol}(H/Z,d\bar{g})\over \dim(\sigma)}\langle w_i, w_j^*\rangle 
\langle w_j, w_i^*\rangle\\
&= {\rm vol}(H/Z,d\bar{g}).
\end{align*}
D'o\`u le point (1).

Soit $K$ un sous--groupe ouvert compact de $H$. 
Appliquons le thorme de \ref{caractres des sries ed} \`a la fonction $\varphi=\Theta_\Sigma$. Pour 
$\gamma\in G^\natural_{\rm qr}$, on a
$$
\dim(\sigma)\Theta_\Pi(\gamma)=d(\pi)\int_{G/Z}\omega(g^{-1})\Theta^K_{\varphi,\gamma}(g)d\bar{g}.
$$
Pour $\delta\in G^\natural$ et $h\in H$, si $\delta\not\in H^\natural$ on a $\varphi(h^{-1}\cdot\delta \cdot h)=\varphi(\delta)=0$, et si $\delta\in H^\natural$ on a 
$$
\varphi(h^{-1}\cdot\delta\cdot h)={\rm tr}(\Sigma(h^{-1}\cdot \delta\cdot h)=
\omega(h)\varphi(\delta).
$$
Par cons\'equent pour $\gamma\in G^\natural_{\rm qr}$, la fonction $G/Z \rightarrow \Bbb{C},\,
g\mapsto \omega(g^{-1})\Theta^K_{\varphi,\gamma}(g)$ se factorise \`a travers $G/H$, et d'apr\`es 
(1), on a
$$
\Theta_\Pi(\gamma)=\sum_{g\in G/H}\omega(g^{-1})\Theta^K_{\varphi,\gamma}(g).
$$
D'o\`u le point (2), l'assertion de finitude se d\'eduisant 
directement de loc.~cit.
\end{proof}

\section{Int\'egrales orbitales et caract\`eres}

Continuons avec les notations des ch.~5 et 6.

\subsection{Int\'egrales orbitales tordues}\label{IO tordues}
Soit un \'el\'ement $\gamma\in G^\natural$ tel que:
\begin{itemize}
\item la $G$--orbite
$\ES{O}_G(\gamma)=\{g^{-1}\cdot \gamma\cdot g:g\in G\}$ est 
ferm\'ee dans $G$;
\item le centralisateur $G_\gamma=\{g\in G: g^{-1}\cdot \gamma\cdot g
= \gamma\}$ est unimodulaire.
\end{itemize}
Le choix d'une mesure de Haar $dg_\gamma$ sur $G_\gamma$ d\'efinit une distribution 
$\Lambda^G_\omega(\cdot , \gamma)= \Lambda^G_\omega (\cdot ,\gamma ,dg_\gamma)$ sur 
$G^\natural$: pour $\phi\in C^\infty_{\rm c}(G^\natural)$, on pose
$$
\Lambda^G_\omega(\phi,\gamma)=\left\{
\begin{array}{ll}
0 &\hbox{si $\omega\vert_{G_\gamma}\neq 1$} \\
\int_{G_\gamma \backslash G}\omega(g)\phi(g^{-1}\cdot \gamma \cdot g) {dg\over dg_\gamma}
&\hbox{sinon}
\end{array}\right..
$$
Puisque la $G$--orbite de $\gamma$ est ferm\'ee dans $G$, l'int\'egrale 
est absolument convergente (d'ailleurs c'est m\^eme une somme finie). On appelle 
$\Lambda^G_\omega(\cdot ,\gamma)$ la {\it $\omega$--int\'egrale orbitale de $\gamma$}. Rappelons que l'on a 
pos ${^x{\phi}}=\phi\in {\rm Int}_{\bf G}(x^{-1})$ ($\phi\in C^\infty_{\rm c}(G)$, $x\in G$). 
Pour $\phi\in C^\infty_{\rm c}(G^\natural)$ 
et $x\in G$, on a
$$
\Lambda^G_\omega({^x\phi},\gamma)= 
\omega(x^{-1})\Lambda^G_\omega(\phi,\gamma)
=\Lambda^G_\omega(\phi,x^{-1}\cdot\gamma\cdot x).\leqno{(*)}
$$

\begin{marema}
{\rm Il est fort probable que le centralisateur $G_\gamma$ d'un lment $\gamma\in G^\natural$ dont la 
$G$--orbite est ferm\'ee dans $G^\natural$, soit 
automatiquement unimodulaire; mais nous n'essaierons pas de le d\'emontrer ici. Notons que si $p=1$ (i.e. si le corps de base $F$ est de caractristique nulle) et si ${\bf G}^\natural$ est localement de type fini, \cad si le $F$--automorphisme de $Z(G)$ dfini par ${\bf G}^\natural$ est d'ordre fini, on sait que pour tout lment $\gamma\in G^\natural$, le centralisateur $G_\gamma$ est unimodulaire.\hfill $\blacksquare$}
\end{marema}

Soit $\gamma\in G^\natural$ quasi--semisimple. D'apr\`es la proposition 1 de \ref{la topo p-adique}, la $G$--orbite 
$\ES{O}_G(\gamma)$ est ferm\'ee dans $G^\natural$. Puisque le 
groupe ${\bf G}_\gamma^\circ$ est r\'eductif (thorme 2 de \ref{automorphismes qss}) et d\'efini sur $F$ 
(thorme de \ref{automorphismes stabilisant (B,T)}), le groupe 
$G_\gamma^\circ = {\bf G}_\gamma^\circ(F)$ est unimodulaire. Comme le groupe 
quotient $G_\gamma/G_\gamma^\circ$ est fini, le choix d'une mesure de Haar 
$dg_\gamma^\circ$ sur $G_\gamma^\circ$ d\'efinit une distribution 
$\Lambda^G_\omega(\cdot , \gamma)= \Lambda^G_\omega (\cdot ,\gamma ,dg_\gamma^\circ)$ sur 
$G^\natural$: pour $\phi\in C^\infty_{\rm c}(G^\natural)$, on pose
$$
\Lambda^G_\omega(\phi,\gamma)=\left\{
\begin{array}{ll}0 &\hbox{si $\omega\vert_{G_\gamma^\circ}\neq 1$}\\ 
\int_{G_\gamma^\circ \backslash G}\omega(g)\phi(g^{-1}\cdot \gamma \cdot g) 
{dg\over dg_\gamma^\circ}
&\hbox{sinon}
\end{array}
\right..
$$
Notons que si $\omega\vert_{G_\gamma}\neq 1$, on a $\Lambda^G_\omega(\phi,\gamma)=0$. 
Par ailleurs, le groupe $G_\gamma$ est lui aussi unimodulaire, et 
si $dg_\gamma$ est une mesure de Haar sur $G_\gamma$, alors 
notant $dg_\gamma^\circ$ la restriction \`a $G_\gamma^\circ$ de la 
mesure de Haar $\vert G_\gamma/G_\gamma^\circ\vert^{-1}dg_\gamma$ sur $G_\gamma$, on a
$$
\Lambda^G_\omega(\cdot,\gamma,dg_\gamma)=\Lambda^G_\omega(\cdot,\gamma,dg_\gamma^\circ).
$$

\subsection{Descente parabolique}\label{descente parabolique} Soit $P\in \ES{P}_\circ$. Rappelons que les mesures de Haar 
$dg$, $dm=dm_P$ et $du=du_P$ sur $G$, $M_P$ et $U_P$ sont celles normalis\'ees par $K_\circ$ (cf. \ref{mesures normalises}). 
Pour toute fonction $f\in C^\infty_{\rm c}(G)$, 
on a la formule d'int\'egration
$$
\int_Gf(g)dg= \int\!\!\!\int\!\!\!\int_{M_P\times U_P\times K_\circ}f(muk)dmdudk.\leqno{(*)}
$$

Soit maintenant $P^\natural\in \ES{P}_\circ^\natural$. Posons $P=N_G(P^\natural)$, $M=M_P$, $M^\natural=
M_P^\natural$ et $P^{-,\natural}= M^\natural\cdot U_{P^-}$. Posons aussi $\frak{m}=\frak{m}_P$. Soit un \'el\'ement $\gamma\in M^\natural$ tel que:
\begin{itemize}
\item la $M$--orbite 
$\ES{O}_M(\gamma)=\{m\cdot\gamma \cdot m^{-1}:m\in M\}$ est ferm\'ee dans $M^\natural$;
\item le centralisateur $M_\gamma=\{m\in M:m^{-1}\cdot \gamma\cdot m=m\}$ est unimodulaire.
\end{itemize}
Alors le choix 
d'une mesure de Haar $dm_\gamma$ sur $M_\gamma$ d\'efinit comme en \ref{IO tordues} une distribution 
$\Lambda^M_\omega(\cdot,\gamma)= \Lambda^M_\omega(\cdot,\gamma,dm_\gamma)$ sur 
$M^\natural$.

\begin{monlem}
\begin{enumerate}
\item[(1)] Pour $\gamma\in P^\natural$, on a $\delta_{P^\natural}(\gamma)=
\vert \textstyle{\det_F}({\rm Ad}_{P^\natural}(\gamma);\frak{u}_P)\vert_F$.
\item[(2)] Pour $\gamma\in M^\natural$, on a 
$\vert \textstyle{\det_F}({\rm Ad}_{P^\natural}(\gamma)^{-1};\frak{u}_P)\vert_F
=\vert \textstyle{\det_F}({\rm Ad}_{P^{-,\natural}}(\gamma);\frak{u}_{P^-})\vert_F$.
\end{enumerate}
\end{monlem}

\begin{proof}
Soit $\gamma\in P^\natural$, et montrons (1). L'application 
${\rm Int}_{P^\natural}(\gamma):U_P\rightarrow U_P$ 
est un automorphisme de vari\'et\'e $\varpi$--adique, de Jacobien constant 
$J(\gamma)=\vert \det_F({\rm Ad}_{P^\natural}(\gamma);\frak{u}_P)\vert_F$. 
D'autre part, d'apr\`es la d\'efinition du module d'un automorphisme 
de $U_P$ 
(cf. \ref{modules}), on a $J(\gamma)
=\Delta_{U_P}({\rm Int}_{P^\natural}(\gamma)\vert_{U_P})$. Posons 
$\overline{P}=P/U_P$, $\smash{\overline{P}}^\natural =P^\natural/U_P$ (c'est 
un $\overline{P}$--espace tordu) et 
$\bar{\gamma}=\gamma\cdot U_P \in \smash{\overline{P}}^\natural$. D'après la remarque de \ref{sous-espaces paraboliques}, on a 
$\delta_{P^\natural}(\gamma)= \Delta_{P^\natural}(\gamma)^{-1}\Delta_{\smash{\overline{P}}^\natural}(\bar{\gamma})$. 
En rempla\c{c}ant $U_P$ par $P$ (resp. $\overline{P}$) dans le raisonnement ci-dessus, on 
obtient
\begin{align*}
\delta_{P^\natural}(\gamma)&= \Delta_P({\rm Int}_{P^\natural}(\gamma))
\Delta_{\overline{P}}({\rm Int}_{\smash{\overline{P}}^\natural}(\bar{\gamma})^{-1})\\
& =\vert \textstyle{\det_F}({\rm Ad}_{P^\natural}(\gamma);\frak{p})\vert_F
\vert \textstyle{\det_F}({\rm Ad}_{\smash{\overline{P}}^\natural}(\bar{\gamma});\frak{p}/\frak{u}_P\vert_F^{-1}\\
&= \vert \textstyle{\det_F}({\rm Ad}_{P^\natural}(\gamma); \frak{u}_P)\vert_F.
\end{align*}

Supposons que $\gamma\in M^\natural$, et montrons (2). D'apr\`es la 
d\'emonstration du point (1), on a
\begin{align*}
\Delta_{G^\natural}(\gamma)^{-1} & = \vert \textstyle{\det_F}({\rm Ad}_{G^\natural}(\gamma);\frak{g})\vert_F\\
&= \vert \textstyle{\det_F}({\rm Ad}_{P^{-,\natural}}(\gamma);\frak{u}_{P^-})\vert_F
\vert \textstyle{\det_F}({\rm Ad}_{M^\natural}(\gamma);\frak{m})\vert_F
\vert \textstyle{\det_F}({\rm Ad}_{P^\natural}(\gamma);\frak{u}_P)\vert_F\\
&= \vert \textstyle{\det_F}({\rm Ad}_{P^{-,\natural}}(\gamma);\frak{u}_{P^-})\vert_F\Delta_{M^\natural}(\gamma)^{-1}
\vert \textstyle{\det_F}({\rm Ad}_{P^\natural}(\gamma);\frak{u}_P)\vert_F.
\end{align*}
On en d\'eduit que
$$
\vert \textstyle{\det_F}({\rm Ad}_{P^\natural}(\gamma)^{-1};\frak{u}_P)\vert_F
=\Delta_{G^\natural}(\gamma)\Delta_{M^\natural}(\gamma)^{-1}\vert \textstyle{\det_F}({\rm Ad}_{P^{-,\natural}}(\gamma);\frak{u}_{P^-})\vert_F.
$$
Or d'après la remarque de \ref{sous-espaces paraboliques}, on a $\Delta_{G^\natural}=1$ et $\Delta_{M^\natural}=1$. 
D'où le point (2).
\end{proof}

Pour $\gamma\in M^\natural$, l'automorphisme ${\rm Ad}_{G^\natural}(\gamma)$ de $\frak{g}$ 
stabilise $\frak{m}$; en fait on a l'galit
$$
{\rm Ad}_{G^\natural}(\gamma)\vert_\frak{m}={\rm Ad}_{M^\natural}(\gamma).
$$
On 
peut donc poser (pour la d\'efinition de $G^\natural/M$, cf. \ref{espaces topo tordus})
$$
D_{ G^\natural/M}(\gamma)=\textstyle{\det_F}({\rm id}-{\rm Ad}_{G^\natural}(\gamma); \frak{g}/\frak{m})
\quad (\gamma\in M^\natural). 
$$

\begin{mapropo}
Soit $\gamma\in M^\natural$ tel que:
\begin{itemize}
\item la $M$--orbite $\ES{O}_M(\gamma)$ 
est ferm\'ee dans $M^\natural$;
\item le centralisateur $M_\gamma$ est 
unimodulaire;
\item on a les inclusions $G_\gamma\subset M$ et $\frak{g}_\gamma\subset \frak{m}$.
\end{itemize}
Alors $D_{G^\natural/M}(\gamma)\neq 0$ et la $G$--orbite $\ES{O}_G(\gamma)$ est ferm\'ee dans $G^\natural$. 
De plus, pour 
toute fonction $\phi\in C^\infty_{\rm c}(G^\natural)$, on a la formule de 
descente
$$
\vert D_{G^\natural/M}(\gamma)\vert_F^{1\over 2} \Lambda^G_\omega(\phi,\gamma)= 
\Lambda^M_\omega({^\omega\phi_{P^\natural,K_\circ}},\gamma);
$$
o\`u les distributions $\Lambda^G_\omega(\cdot,\gamma)$ et $\Lambda^M_\omega(\cdot,\gamma)$ sur $G^\natural$ et sur $M^\natural$ 
sont d\'efinies par la m\^eme mesure de Haar $dg_\gamma$ sur $G_\gamma=M_\gamma$ 
(pour la d\'efinition de ${^\omega\phi_{P^\natural,K_\circ}}$, cf. \ref{induction parabolique et caractres}).
\end{mapropo}

\begin{proof}
Pour $m\in M$, puisque $G_{m^{-1}\cdot \gamma \cdot m}\cap U_P=\{1\}$ et $\frak{g}_{m^{-1}\cdot\gamma \cdot m}\cap \frak{u}_P=\{0\}$, l'application
$$
U_P\rightarrow U_P,\, u\mapsto {\rm Int}_{P^\natural}(m^{-1}\cdot \gamma \cdot m)^{-1}(u^{-1})u
$$
est un automorphisme de vari\'et\'e $\varpi$--adique, de Jacobien
$$
J(m^{-1}\cdot \gamma \cdot m)= \vert \textstyle{\det_F} ( {\rm id}-{\rm Ad}_{P^\natural}(m^{-1}\cdot \gamma \cdot m)^{-1};\frak{u}_P)\vert_F\neq 0.
$$
Comme
$${\rm id}_\frak{p}-{\rm Ad}_{P^\natural}(m^{-1}\cdot \gamma \cdot m)^{-1}= 
{\rm Ad}_P(m^{-1})\circ({\rm id}_\frak{p}-
{\rm Ad}_{P^\natural}(\gamma)^{-1})\circ {\rm Ad}_P(m),$$
on a
$$
J(m^{-1}\cdot \gamma\cdot m)=J(\gamma).
$$
D'apr\`es le lemme, on a
\begin{align*}
\delta_{P^\natural}^{1/2}(\gamma)J(\gamma)&= 
\vert \textstyle{\det_F}({\rm Ad}_{P^\natural}(\gamma);\frak{u}_P)\vert_F^{1\over 2}\vert \textstyle{\det_F}({\rm id}-{\rm Ad}_{P^\natural}(\gamma)^{-1};\frak{u}_P)\vert_F\\
&=\vert \textstyle{\det_F}({\rm Ad}_{P^\natural}(\gamma)-{\rm id};\frak{u}_P)\vert_F^{1\over 2}\vert \textstyle{\det_F}({\rm id}-{\rm Ad}_{P^\natural}(\gamma)^{-1};\frak{u}_P)\vert_F^{1\over 2}\\
&=\vert \textstyle{\det_F}({\rm id}-{\rm Ad}_{P^\natural}(\gamma);\frak{u}_P)\vert_F^{1\over 2}\vert \textstyle{\det_F}({\rm id}-{\rm Ad}_{P^{-,\natural}}(\gamma);\frak{u}_{P^-})\vert_F^{1\over 2}\\
&=\vert \textstyle{\det_F}({\rm id}-{\rm Ad}_{G^\natural}(\gamma);\frak{g}/\frak{m})\vert_F^{1\over 2}.
\end{align*}
On a donc
$$
\delta_{P^\natural}^{1/2}(\gamma)J(\gamma)=\vert D_{G^\natural/M}(\gamma)\vert_F^{1\over 2}\neq 0.
$$

Pour $g\in G$, crivons $g=kmu$ ($k\in K_\circ$, $m\in M$, $u\in U_P$) et posons $\gamma'=m^{-1}\cdot\gamma \cdot m$. 
Alors on a
$$
g^{-1}\cdot \gamma \cdot g = {\rm Int}_G(k^{-1})(\gamma' \cdot {\rm Int}_{P^\natural}(\gamma')^{-1}(u^{-1})u).
$$
Par consquent
$$
\ES{O}_G(\gamma)= {\rm Int}_G(K)(\ES{O}_M(\gamma)\cdot U_P),
$$
et puisque la $M$--orbite $\ES{O}_M(\gamma)$ est ferme dans $M^\natural$, la $G$--orbite $\ES{O}_G(\gamma)$ est ferme dans 
$G^\natural$.

Soit une fonction $\phi\in C^\infty_{\rm c}(G^\natural)$. 
D'apr\`es la relation $(*)$, on a
$$
\Lambda^G_\omega(\phi,\gamma)=\int\!\!\!\int\!\!\!\int_{M_\gamma \backslash M
\times U_P \times K_\circ}\omega(muk)\phi(k^{-1}u^{-1}m^{-1}\cdot \gamma \cdot muk){dm\over dg_\gamma}dudk.
$$
Comme $\omega\vert_{U_P}=1$ (lemme de \ref{induction parabolique et caractres}), 
le changement de variables $u\mapsto {\rm Int}_G(m^{-1}\cdot\gamma\cdot m)^{-1}(u^{-1})u$ dans la 
formule pour $\Lambda^G_\omega(\phi,\gamma)$ donne
\begin{align*}
\Lambda^G_\omega(\phi,\gamma)&=\int_{M_\gamma\backslash M} \omega(m)\delta_{P^\natural}^{-1/ 2}(\gamma)
{^\omega\phi_{P^\natural,K_\circ}}(m^{-1}\cdot \gamma \cdot m)J(\gamma)^{-1}{dm\over dg_\gamma}\\
&= \vert D_{G^\natural/M}(\gamma)\vert_F^{-{1\over 2}}\Lambda^M_\omega({^\omega\phi_{P^\natural,K_\circ}},\gamma).
\end{align*}
La proposition est d\'emontr\'ee.
\end{proof}

\subsection{Formule d'int\'egration de Weyl}\label{formule d'intgration de Weyl}
Soit $(S,T,T^\natural)$ un triplet de Cartan de $G^\natural$ (\ref{tores maximaux et sous-espaces de Cartan}). Puisque $Z_G(T^\natural)=S$, l'application
$$
\pi:S\backslash G \times T^\natural \rightarrow G^\natural,\, (g,\gamma)\mapsto g^{-1}\cdot \gamma\cdot g
$$
est bien d\'efinie, et c'est un morphisme de vari\'et\'es $\varpi$--adiques. Posons $\frak{g}={\rm Lie}(G)$, 
$\frak{s}={\rm Lie}(S)$, $\overline{\frak{g}}=\frak{g}/\frak{s}$ et $\frak{t}={\rm Lie}(T)$. Soit $(\overline{g},\gamma)\in S\backslash G \times T^\natural$, et soit $g$ un rel\`evement de $\bar{g}$ dans $G$. Via les isomorphismes de vari\'et\'es 
$\varpi$--adiques
\begin{align*}
& S\backslash G \rightarrow S\backslash G,\, Sx\rightarrow Sxg,\\
& T\rightarrow T^\natural,\, t\mapsto t\cdot \gamma,\\
& G\rightarrow G^\natural, x \mapsto g^{-1}\cdot (x\cdot \gamma) \cdot g,
\end{align*}
identifions l'espace tangent \`a $S\backslash G$ (resp. $T^\natural$, $G^\natural$) en 
$\bar{g}$ (resp. $\gamma$, $g^{-1}\cdot \gamma \cdot g$) \`a $\overline{\frak{g}}$ 
(resp. $\frak{t}$, $\frak{g}$). Notons $p:\frak{g}\rightarrow \overline{\frak{g}}$ la projection canonique. Alors 
la diff\'erentielle ${\rm d}(\pi)_{\bar{g},\gamma}:\overline{\frak{g}}\times \frak{t}\rightarrow \frak{g}$ de 
$\pi$ au point $(\bar{g},\gamma)$, est donn\'ee par
$$
{\rm d}(\pi)_{\bar{g},\gamma}(p(X),Y)= {\rm Ad}_G(\gamma)(X)-X +Y
$$
D'apr\`es le d\'ebut de la d\'emonstration de la proposition de \ref{lments rguliers et quasi-rguliers}, si
$\gamma\in T^\natural \cap G^\natural_{\rm reg}$, on a l'\'egalit\'e $\frak{g}(1-\gamma)+\frak{t}=\frak{g}$. 
On en d\'eduit le

\begin{monlem1}
Pour tout $(\bar{g},\gamma)\in S\backslash G\times(T^\natural\cap G^\natural_{\rm reg})$, 
l'application $\pi$ est submersive en $(\bar{g},\gamma)$.
\end{monlem1}

Notons $W=W(G,T)$ le groupe de Weyl $N_G(T)/T$. Puisque 
$${\rm Int}_{G^\natural}(\gamma)(N_G(T))= N_G({\rm Int}_{G^\natural}(\gamma)(T))\quad 
(\gamma\in G^\natural),
$$ 
tout \'el\'ement $\gamma\in T^\natural$ induit un automorphisme 
$\tau_\gamma = {\rm Int}_{G^\natural}(\gamma)\vert_{N_G(T)}$ de $N_G(T)$. Pour $\gamma\in T^\natural$, $t\in T$ et $n\in N_G(T)$, on a
$$
\tau_{t\cdot \gamma}(n)= t\cdot \tau_\gamma(n)\cdot t^{-1}= t \tau_\gamma(n)t^{-1}\tau_\gamma(n)^{-1} \tau_\gamma(n)\in T\tau_\gamma(n).
$$
On en d\'eduit que $\tau_\gamma$ induit
\begin{itemize}
\item par passage au quotient un automorphisme 
$\tau_W$ de $W$,
\item par restriction un automorphisme $\tau=\tau_T$ de $T$.
\end{itemize}
Ces deux automorphismes 
{\it ne d\'ependent pas} du choix de $\gamma\in T^\natural$ (bien que $\tau_\gamma$ en 
d\'epende). Posons 
\begin{align*}
&N_G(S^\natural)=\{g\in G: g\cdot S^\natural \cdot g^{-1}= S^\natural\},\\
& N_G(T^\natural)=\{g\in G: g\cdot T^\natural \cdot g^{-1}= T^\natural\}.
\end{align*}
On a les inclusions
$$
N_G(S^\natural)\subset N_G(T^\natural)\subset N_G(S) \subset N_G(T).
$$
L'inclusion du milieu et celle  droite rsultent des galits $Z_G(T^\natural)=S$ et 
$Z_G(S)=T$ (proposition de \ref{tores maximaux et sous-espaces de Cartan}). Quant  l'inclusion de gauche, 
puisque $Z_G(S^\natural)= S$ (loc.~cit.), on a l'inclusion $N_G(S^\natural)\subset N_G(S)$. Choisissons un lment 
$\delta_0\in S^\natural\cap G^\natural_{\rm reg}$. Pour $n\in N_G(S^\natural)$, on a ${\rm Int}_{{\bf G}^\natural}(n)(\delta_0)\in S^\natural 
=S\cdot \delta_0$, et comme $T^\natural = T\cdot \delta_0$, on obtient
$$
n\cdot T^\natural \cdot n^{-1} = (nTn^{-1})\cdot {\rm Int}_{{\bf G}^\natural}(n)(\delta_0)\subset T\cdot S^\natural = T^\natural. 
$$
Notons $W^\natural$ le sous--groupe de $W$ form\'e des \'el\'ements fix\'es par $\tau_W$, et posons
\begin{align*}
& W(G,S^\natural)= N_G(S^\natural)/S,\\
& W(G,T^\natural)=N_G(T^\natural)/T.
\end{align*}
%

\begin{monlem2}
\begin{enumerate}
\item[(1)] On a l'\'egalit\'e $W^\natural = W(G,T^\natural)$.
\item[(2)] Supposons que $T_\tau =S$, $T(1-\tau)\cap S =\{1\}$ et $T=T_\tau T(1-\tau)$. Alors l'inclusion 
$N_G(S^\natural)\subset N_G(T^\natural)$ induit par passage aux quotients 
une identification
$$W(G,S^\natural) = W(G,T^\natural).$$
\end{enumerate}
\end{monlem2}

\begin{proof}
Pour $n\in N_G(T)$ et $\gamma\in T^\natural$, on a 
$n\cdot \gamma \cdot n^{-1}= n\tau_\gamma(n^{-1})\cdot \gamma$, et 
$n\cdot \gamma \cdot n^{-1}\in T^\natural$ si et seulement si 
$n\tau_\gamma(n^{-1})\in T$. D'o\`u le point (1).

Supposons que $T_\tau=S$, $T(1-\tau)\cap S =\{1\}$ et $T=T_\tau T(1-\tau)$, et montrons (2). 
Il s'agit de montrer que 
$N_G(T^\natural)= N_G(S^\natural)T\;(=TN_G(S^\natural))$ et $N_G(S^\natural)\cap T=S$. 
Soit $n\in N_G(T^\natural)$ et $\gamma\in T^\natural$. Posons $t= n\tau_\gamma(n^{-1})\in T$, 
et \'ecrivons $t= s \tau(x)x^{-1}= x^{-1}s\tau(x)$ avec $s\in S$ et $x\in T$. On a donc
$$
n\cdot \gamma \cdot n^{-1} = n\tau_\gamma(n^{-1})\cdot \gamma = x^{-1}s\tau_\gamma(x)\cdot \gamma = x^{-1}s\cdot \gamma \cdot x,
$$
d'o\`u $xn\cdot \gamma \cdot (xn)^{-1}=s\cdot \gamma$. Comme $xn\in N_G(S)$, 
on en d\'eduit que $xn\in N_G(S^\natural)$ et $n\in TN_G(S^\natural)$. Soit maintenant 
$y\in N_G(S^\natural)\cap T$. Puisque $y \cdot S^\natural \cdot y^{-1}= y\tau(y^{-1})\cdot S^\natural$,  on a $y\tau(y^{-1})\in T(1-\tau)\cap S =\{1\}$, 
i.e. $y\in T_\tau= S$.
\end{proof}

\begin{marema1}
{\rm Les hypothses du point (2) sont rarement satisfaites. D'apr\`es la proposition de \ref{automorphismes ss et u}, elles le sont 
si $\tau$ est unipotent. \hfill $\blacksquare$}
\end{marema1}

On d\'efinit comme suit une op\'eration \`a gauche, continue et libre, de $S\backslash T$ sur 
$S\backslash G \times T^\natural$: pour $\bar{t}\in S\backslash T$ et $(\bar{g},\gamma)\in 
S\backslash G\times T^\natural$, on choisit un repr\'esentant $t$ de $\bar{t}$ dans $T$ et 
un repr\'esentant $g$ de $\bar{g}$ dans $G$, et l'on pose
\begin{align*}
& \bar{t}\cdot \bar{g} = S tg,\\
& \bar{t}\cdot (\bar{g},\gamma)= (\bar{t}\cdot\bar{g},
t\cdot \gamma\cdot t^{-1})=(\bar{t}\cdot\bar{g}, t\tau(t^{-1})\cdot \gamma).
\end{align*}
Puisque $T=Z_G(S)$ et $\tau\vert_S={\rm id}$, l'\'el\'ement $\bar{t}\cdot (\bar{g},\gamma)$ est 
bien d\'efini. Soit
$$
X=S\backslash G \times_{S\backslash T} T^\natural
$$
le quotient de $S\backslash G\times T^\natural$ par la relation d'\'equivalence $\sim$ 
d\'efinie par:
$$
(\bar{g},\gamma)\sim (\bar{g}',\gamma') \Leftrightarrow  \hbox{il existe un $\bar{t}\in S\backslash T$ tel 
que $(\bar{g}',\gamma')= \bar{t}\cdot (\bar{g},\gamma)$}. 
$$
On note $q:S\backslash G \times T^\natural \rightarrow X$ la projection canonique. 
L'application $\pi:S\backslash G\times T^\natural\rightarrow G^\natural$ se factorise \`a travers $q$: on note $\bar{\pi}: X\rightarrow G^\natural$ 
l'unique application telle que $\bar{\pi}\circ q = \pi$. Par construction, $X$ est une 
vari\'et\'e $\varpi$--adique (pour la topologie quotient) de dimension $\dim (G^\natural)=\dim (G)$, 
et $q$ et $\bar{\pi}$ 
sont des morphismes de vari\'et\'es $\varpi$--adiques.

On d\'efinit comme suit une op\'eration \`a gauche, continue et 
libre, de $W^\natural$ sur 
$X$: pour $w\in W^\natural$ et 
$(\bar{g},\gamma)\in S\backslash G \times T^\natural$, on choisit un repr\'esentant 
$g$ de $\bar{g}$ dans $G$ et un repr\'esentant $n$ de $w$ dans $N_G(T^\natural)$, et l'on pose
$$
w\cdot q(\bar{g},\gamma)=q(Sng, n\tau_\gamma(n)\cdot \gamma).
$$
L'op\'eration est bien d\'efinie: si $g'=xg$ et $n'=yn$ pour des $x,\,y\in T$, alors posant $\gamma'=x\tau(x^{-1})\cdot \gamma$ et $t= ynxn^{-1}\in T$, on a
\begin{align*}
(Sn'g', n'\tau_{\gamma'}(n'^{-1})\cdot \gamma')&= (Stg,ynx\tau(x^{-1})\cdot \gamma \cdot n^{-1}y^{-1})\\
& = (Stg, tn\tau_\gamma(n^{-1})\tau_\gamma(nx^{-1})\tau_\gamma(n^{-1}y^{-1})\cdot \gamma )\\
& = (Stg, tn\tau_\gamma (n^{-1})\tau_\gamma(t^{-1})\cdot \gamma)\\
& = (Stg, t\tau(t^{-1}) n\tau_\gamma (n^{-1})\cdot \gamma).
\end{align*}
Pour $(\bar{g},\gamma)\in S\backslash G\times T^\natural$ et $w\in W^\natural$, on a
$$
\bar{\pi}(w\cdot q(\bar{g},\gamma))= \pi(\bar{g},\gamma).\leqno{(*)}
$$ 

Pour toute partie $Y$ de $G^\natural$, notons ${^GY}$ l'ensemble des 
$g^{-1}\cdot \delta\cdot g$ pour $g\in G$ et $\delta\in Y$. 
Ainsi on a ${^G{T^\natural}}={\rm Im}(\pi)$, et ${^G{T^\natural}}\cap G^\natural_{\rm reg}=
{^G(T^\natural \cap G^\natural_{\rm reg})}=\pi (S\backslash G\times (T^\natural \cap G^\natural_{\rm reg}))$.

\begin{monlem3}
Pour $\delta\in {^G{T^\natural}}\cap G^\natural_{\rm reg}$, 
la fibre
$$\bar{\pi}^{-1}(\delta)=
q(\{(\bar{g},\gamma)\in S\backslash G\times T^\natural: \pi(\bar{g},\gamma)=\delta\})$$
de $\bar{\pi}$ au dessus de $\delta$, est un espace principal homog\`ene sous 
$W^\natural$.
\end{monlem3}

\begin{proof}
Soit $\delta\in {^G{T^\natural}}\cap G^\natural_{\rm reg}$. \'Ecrivons 
$\delta =\pi(\bar{g}'\!,\gamma')$ pour un $(\bar{g}'\!,\gamma')\in S\backslash G\times 
T^\natural$, et soit $(\bar{g},\gamma)\in S\backslash G\times 
T^\natural$ tel que $\pi(\bar{g},\gamma)= \delta$. Choisissons des rel\`evements 
$g$ et $g'$ de $\bar{g}$ et $\bar{g}'$ dans $G$. Puisque $\delta\in 
{^G{T^\natural}}\cap G^\natural_{\rm reg}$, on a $\gamma,\, \gamma'\in 
T^\natural \cap G^\natural_{\rm reg}$, et comme $g'g^{-1}\cdot \gamma \cdot gg'^{-1}=\gamma'$ et 
${\bf G}_\gamma^\circ (F) = S = {\bf G}_{\gamma'}^\circ(F)$, on a 
$n=g'g^{-1}\in N_G(S)$. Alors $n\cdot \gamma\cdot n^{-1}= n\tau_\gamma(n^{-1})\cdot \gamma = 
\gamma'\in T^\natural$. Par cons\'equent $n\tau_\gamma(n)\in T$, i.e. $n\in N_G(T^\natural)$. 
Posons $w\in nT\in N_G(T^\natural)/T=W^\natural$. On a $w\cdot q(\bar{g},\gamma)= q(\bar{g'},\gamma')$. 
D'o\`u le lemme, puisque $W^\natural$ op\`ere librement sur $X$.
\end{proof}

Posons
$$
X_{\rm reg}= \bar{\pi}^{-1}({^G{T^\natural}}\cap G^\natural_{\rm reg})=
q(S\backslash G\times (T^\natural\cap G^\natural_{\rm reg})).
$$
D'apr\`es le lemme 1, l'ensemble
$$
{^GT^\natural}\cap G^\natural_{\rm reg}=\pi(S\backslash G\times (T^\natural\cap 
G^\natural_{\rm reg}))
$$
est ouvert dans $G^\natural$, par suite $X_{\rm reg}$ est ouvert dans $X$. 
Rappelons que $\dim (X)=\dim (G^\natural)$. D'apr\`es les lemmes 1 et 3, 
l'application $\bar{\pi}:X\rightarrow G^\natural$ induit par restriction  
une application surjective 
$\bar{\pi}_{\rm reg}:X_{\rm reg}\rightarrow {^GT^\natural}\cap G^\natural_{\rm reg}$ 
v\'erifiant:
\begin{itemize}
\item pour $x\in X_{\rm reg}$, la diff\'erentielle ${\rm d}(\bar{\pi})_x$ de $\bar{\pi}$ 
au point $x$ est un isomorphisme;
\item pour $\delta\in {^GT^\natural}\cap G^\natural_{\rm reg}$, la fibre $\bar{\pi}^{-1}(\delta)$ est un espace principal homog\`ene sous $W^\natural$.
\end{itemize}
En d'autres termes, $\bar{\pi}_{\rm reg}$ est un rev\^etement galoisien principal de 
groupe $W^\natural$.

Notons $\overline{T}$ le groupe quotient $T(1-\tau)\backslash T$, et 
$\smash{\overline{T}}^\natural$ le 
$\overline{T}$--espace tordu $T(1-\tau)\backslash T^\natural$. 
Pour $\gamma\in T^\natural$, on pose
$$
D_{G^\natural/T}(\gamma)={\rm \det}_F({\rm id}-{\rm Ad}_{G^\natural}(\gamma);\frak{g}/\frak{t}). 
$$
Pour $\gamma\in T^\natural$ et 
$t\in T$, on a $t^{-1}\tau(t)\cdot\gamma = t^{-1}\cdot \gamma \cdot t$, d'o\`u
\begin{align*}
D_{G^\natural/T}(t^{-1}\tau(t)\cdot\gamma)&= {\rm \det}_F({\rm id}-{\rm Ad}_{G^\natural}(t^{-1}\cdot \gamma\cdot t);\frak{g}/\frak{t})\\
&= {\rm det}_F({\rm Ad}_G(t^{-1})\circ ({\rm id}-{\rm Ad}_{G^\natural}(\gamma))\circ
{\rm Ad}_G(t);\frak{g}/\frak{t})\\
&= D_{G^\natural/T}(\gamma).
\end{align*}
Par cons\'equent l'application $D_{G^\natural/T}:T^\natural\rightarrow F$ se factorise \`a travers 
$\smash{\overline{T}}^\natural$.

Toute mesure de Haar 
$ds$ sur $S$ d\'efinit 
comme suit une mesure de Haar $d\bar{\gamma}=d\bar{\gamma}(ds)$ 
sur $\smash{\overline{T}}^\natural$. Puisque le groupe $T_\tau/S$ est fini, il existe une unique 
mesure de Haar sur $T_\tau$ prolongeant $ds$, que l'on note encore 
$ds$. Via la suite exacte longue de groupes topologiques
$$
1\rightarrow T_\tau \rightarrow T\buildrel 1-\tau\over{\longrightarrow} T \rightarrow \overline{T}
\rightarrow 1,
$$
$ds$ d\'efinit une mesure de Haar $d\bar{t}=d\bar{t}(ds)$ sur le groupe 
quotient $\overline{T}$. Pr\'ecis\'ement, choisissons une mesure de Haar 
$dt$ sur $T$, et notons $d\tilde{t}$ l'image de la mesure quotient ${dt\over ds}$ sur 
$T_\tau\backslash T$ par l'isomorphisme de groupes 
topologiques de $T_\tau\backslash T\buildrel 1-\tau \over{\longrightarrow} T(1-\tau)$. 
La mesure quotient $d\bar{t}= {dt\over d\tilde{t}}$ sur $\overline{T}$ d\'epend 
seulement de $ds$, et pas du choix de $dt$. Alors $d\bar{\gamma}$ est la mesure de 
Haar gauche sur $\smash{\overline{T}}^\natural$ associ\'ee \`a $d\bar{t}$ comme en \ref{module d'un espace tordu} 
(c'est aussi une mesure de Haar  droite); i.e. 
on pose $d\bar{\gamma}= \bar{\delta}\cdot d\bar{t}$ pour un (i.e. pour tout) $\bar{\delta}\in \smash{\overline{T}}^\natural$. 
De mani\`ere \'equivalente, notant $d\gamma$ la mesure de Haar  gauche sur 
$T^\natural$ associ\'ee \`a $dt$, la mesure quotient ${d\gamma \over d\tilde{t}}$ sur $\smash{\overline{T}}^\natural$ 
d\'epend seulement de $ds$, et pas du choix de $dt$, et elle 
co\"{\i}ncide avec $d\bar{\gamma}$.

Notons $d\delta$ la mesure de Haar $\delta_1\cdot dg$ sur $G^\natural$.

\begin{marema2}
{\rm 
Posons $Y=T^\natural\cap G^\natural_{\rm reg}$. 
Puisque ${^GY}$
est ouvert dans $G^\natural$ et que $G^\natural_{\rm reg}$ est dense 
dans $G^\natural$, l'ensemble ${^GT^\natural}\smallsetminus {^GY}$ est n\'egligeable dans 
$G^\natural$ (par rapport \`a $d\delta$). En particulier,  
toute fonction int\'egrable $\phi$ sur $G^\natural$ est int\'egrable sur ${^GT^\natural}$, et 
l'on a $\int_{{^GT^\natural}}\phi(\delta)d\delta= \int_{^GY}\phi(\delta)d\delta$. 
D'autre part, l'image $\overline{Y}$ de $Y$ dans 
$\smash{\overline{T}}^\natural$ est ouverte 
dans $\smash{\overline{T}}^\natural$, et l'ensemble 
$\smash{\overline{T}}^\natural \smallsetminus 
\overline{Y}$ est n\'egligeable dans $\smash{\overline{T}}^\natural$ 
(par rapport \`a $d\bar{\gamma}$).\hfill $\blacksquare$}
\end{marema2}

On en d\'eduit la \og formule d'intgration de H.~Weyl \fg{} suivante:

\begin{mapropo}
Soit $ds$ une mesure de Haar sur $S$. Posons $d\bar{g}={dg\over ds}$ et $d\bar{\gamma}=d\bar{\gamma}(ds)$. 
Soit $\phi$ une 
fonction int\'egrable sur ${^GT^\natural}$ (par rapport \`a $d\delta$). Pour presque tout 
$\gamma\in T^\natural$, 
la fonction $S\backslash G\rightarrow \Bbb{C},\, \bar{g}\mapsto \phi\circ \pi(\bar{g},\gamma)$ 
est int\'egrable sur $S\backslash G$ (par rapport \`a $d\bar{g}$), et l'on a
$$
\int_{^GT^\natural}\phi(\delta)d\delta = {1\over \vert W^\natural\vert}
\int_{\smash{\overline{T}}^\natural}\vert D_{G^\natural/T}(\gamma)\vert_F\left\{\int_{S\backslash G}\phi(g^{-1}\cdot \gamma\cdot g)d\bar{g}\right\}d\bar{\gamma}.
$$
\end{mapropo}

Fixons un syst\`eme de repr\'esentants $\ES{C}$ des classes de $G$--conjugaison de 
sous--espaces de Cartan $T^\natural$ de $G^\natural$ ---  ou, ce qui revient au m\^eme, de 
triplets de Cartan $(S,T,T^\natural)$ de $G^\natural$. Pour chaque 
\'el\'ement $(S,T,T^\natural)$ de $\ES{C}$, le choix d'une mesure de Haar $ds$ sur 
$S$ d\'efinit une mesure de Haar $d\bar{\gamma}= d\bar{\gamma}(ds)$ sur 
$\smash{\overline{T}}^\natural$ et, pour tout $\gamma \in T^\natural \cap G^\natural_{\rm reg}$, une distribution $\Lambda^G_\omega(\cdot,\gamma)$ sur $G^\natural$: pour $\phi\in C^\infty_{\rm c}(G^\natural)$, on pose
$$
\Lambda^G_\omega(\phi,\gamma)=\int_{S\backslash G}\omega(g)\phi(g^{-1}\cdot \gamma\cdot g){dg\over ds}.
$$ 

\begin{moncoro1}
Soit $\Theta$ une 
fonction localement int\'egrable sur $G^\natural$ (par rapport \`a $d\delta$), 
d\'efinie sur $G^\natural_{\rm reg}$ et telle que
$$
\Theta(g^{-1}\cdot \gamma \cdot g)=\omega(g)\Theta(\gamma)\quad (g\in G,\, \gamma\in G^\natural_{\rm reg}).
$$
Pour toute fonction $\phi\in C^\infty_{\rm c}(G^\natural)$, on a 
$$
\int_{G^\natural}\Theta(\delta)\phi(\delta)d\delta =\sum_{T^\natural}
{1\over \vert W(G,T^\natural)\vert}
\int_{\smash{\overline{T}}^\natural}\vert D_{G^\natural/T}(\gamma)\vert_F 
\Theta(\gamma) \Lambda^G_\omega(\phi, \gamma)d\bar{\gamma},
$$
o\`u $T^\natural$ parcourt les \'el\'ements de $\ES{C}$.
\end{moncoro1}

\begin{proof}
Remarquons tout d'abord que pour $T^\natural\in \ES{C}$, $\gamma\in 
T^\natural\cap G^\natural_{\rm reg}$ et $t\in T$, on a
$$\Lambda^G_\omega(\phi, t^{-1}\cdot\gamma
\cdot t)=\omega(t^{-1})\Lambda^G_\omega(\phi,\gamma).$$
Par cons\'equent la fonction 
$\gamma\mapsto \Theta(\gamma)\Lambda^G_\omega(\phi,\gamma)$ sur $T^\natural\cap G^\natural$ se factorise \`a travers l'image de $T^\natural\cap G^\natural$ dans $\smash{\overline{T}}^\natural$, et l'\'enonc\'e a bien un sens.

La fonction $\delta\mapsto \Theta(\delta)\phi(\delta)$ 
est int\'egrable sur $G^\natural$ (et d\'efinie sur $G^\natural_{\rm reg}$). On a donc
$$
\int_{G^\natural}\Theta(\delta)\phi(\delta)d\delta =\sum_{T^\natural}
\int_{^GT^\natural}\Theta(\delta)\phi(\delta)d\delta
$$
o\`u $T^\natural$ parcourt les \'el\'ements de $\ES{C}$. Pour 
$(S,T,T^\natural)\in \ES{C}$, 
d'apr\`es la proposition on a
\begin{align*}
\int_{^GT^\natural}\Theta(\delta)\phi(\delta)d\delta &= {1\over \vert W^\natural\vert}
\int_{\smash{\overline{T}}^\natural}\vert D_{G^\natural/T}(\gamma)\vert_F\left\{\int_{S\backslash G}
\Theta(g^{-1}\cdot\gamma\cdot g)\phi(g^{-1}\cdot \gamma\cdot g){dg\over ds}\right\}d\bar{\gamma}\\
 & = {1\over \vert W^\natural\vert}
\int_{\smash{\overline{T}}^\natural}\vert D_{G^\natural/T}(\gamma)\vert_F
\Theta(\gamma)\Lambda^G_\omega(\phi,\gamma)d\bar{\gamma};
\end{align*}
o\`u l'on a pos\'e $W^\natural=W(G,T^\natural)$. D'o\`u le corollaire.
\end{proof}

\begin{moncoro2}
Soit $\Pi$ une $\omega$--repr\'esentation admissible de $G^\natural$ telle que $\Pi^\circ$ est 
de type fini. Pour toute fonction $\phi\in C^\infty_{\rm c}(G^\natural_{\rm reg})$, on a
$$
\Theta_\Pi(\phi) =\sum_{T^\natural}
{1\over \vert W(G,T^\natural)\vert}
\int_{\smash{\overline{T}}^\natural}\vert D_{G^\natural/T}(\gamma)\vert_F 
\Theta_\Pi(\gamma) \Lambda^G_\omega(\phi, \gamma)d\bar{\gamma},
$$
o\`u $T^\natural$ parcourt les \'el\'ements de $\ES{C}$.
\end{moncoro2}

\begin{moncoro3}
Soit $P^\natural\in \ES{P}^\natural_\circ$, $\Sigma$ une $\omega$--repr\'esentation 
admissible de $M_P^\natural$ telle que $\Sigma^\circ$ est de type fini, 
et $\Pi={^\omega\iota}_{P^\natural}^{G^\natural}(\Sigma)$.
\begin{enumerate}
\item[(1)] Le support de la fonction caract\`ere $\Theta_\Pi: G^\natural_{\rm qr}\rightarrow {\Bbb C}$ (\cad l'ensemble des $\gamma\in G^\natural_{\rm qr}$ tels que $\Theta_\Pi(\gamma)\neq 0$) est contenu dans ${^G(M^\natural_P\cap G_{\rm qr}^\natural)}
={^G{M_P^\natural}}\cap G_{\rm qr}^\natural$.
\item[(2)] Soit $T^\natural$ un sous--espace de Cartan de $M_P^\natural$, et soit  
$\{T_1^\natural,\ldots ,T_s^\natural\}$ un syst\`eme de repr\'esentants des classes de $M_P$--conjugaison 
de sous--espaces de Cartan de $M_P^\natural$ conjugu\'es \`a $T^\natural$ dans $G$: 
\begin{itemize}
\item $T^\natural_i$ est un sous--espace de Cartan de $M^\natural_P$ de la forme $g_i^{-1}\cdot T^\natural \cdot g_i$ pour un \'el\'ement $g_i\in G$;
\item tout sous--espace de Cartan de $M_P^\natural$ conjugu\'e \`a $T^\natural$ dans 
$G$ est conjugu\'e \`a l'un des $T^\natural_i$ dans $M_P$;
\item les $T^\natural_i$ sont deux--\`a--deux non conjugu\'es dans $M_P$.
\end{itemize}
Pour $i=1,\ldots s$, posons $\smash{\overline{W}}_i^\natural=W(M_P,T^\natural_i)\backslash W(G,T^\natural_i)$, 
et soit $\{n_w: w\in \smash{\overline{W}}^\natural_i\}$ 
un syst\`eme de repr\'esentants des \'el\'ements 
de $\smash{\overline{W}}^\natural_i\;(=N_{M_P}(T^\natural_i)\backslash N_G(T^\natural_i))$ dans $N_G(T^\natural_i)$. Pour $\gamma \in T^\natural\cap G_{\rm reg}^\natural$, 
on a l'\'egalit\'e
$$
\Theta_\Pi(\gamma)= \sum_{i=1}^s\sum_{w\in \smash{\overline{W}}^\natural_i}\omega(g_in_w)^{-1}\vert D_{G^\natural/M_P}(\gamma^{g_i n_w})\vert_F^{-{1\over 2}}\Theta_\Sigma(\gamma^{g_i n_w})
$$
o\`u l'on a pos\'e $\gamma^{g_in_w}={\rm Int}'_G(g_i n_w)^{-1}(\gamma)\;(=n_w^{-1}g_i^{-1}\cdot \gamma\cdot g_in_w)$.
\end{enumerate}
\end{moncoro3}

\begin{proof}
Posons $M=M_P$ et $M^\natural=M^\natural_P$. Puisque $\Sigma^\circ$ est admissible 
et de type fini, $\Pi^\circ$ l'est aussi. Les fonctions caract\`eres $\Theta_\Sigma$ sur 
$M^\natural_{\rm qr}$ et $\Theta_\Pi$ sur $G^\natural_{\rm qr}$ sont donc bien d\'efinies. 

Montrons (1). D'apr\`es le thorme de \ref{induction parabolique et caractres}, il suffit de montrer que si $\phi$ est une fonction 
dans $C^\infty_{\rm c}(G^\natural_{\rm qr})$ dont le support ne rencontre pas ${^GM^\natural}\cap G^\natural_{\rm qr}$, alors 
${^\omega\phi_{P^\natural,K_\circ}}=0$. Soit une fonction $\phi\in C^\infty_{\rm c}(G^\natural_{\rm qr})$ 
telle que ${^\omega \phi_{P^\natural,K_\circ}}\neq 0$. Puisque $M^\natural\cap G^\natural_{\rm qr}$ est ouvert dense dans $M^\natural$, il existe un 
$\delta\in M^\natural\cap G^\natural_{\rm qr}$ tel que ${^\omega\phi_{P^\natural,K_\circ}}(\delta)\neq 0$. D'apr\`es la d\'efinition de 
${^\omega\phi_{P^\natural,K_\circ}}$, il existe un $k\in K_\circ$ et un $u\in U_P$ tels que $\phi(k^{-1}\cdot\delta \cdot uk)\neq 0$, 
et d'apr\`es la d\'emonstration de la proposition de \ref{descente parabolique}, il existe un (unique) $u'\in U_P$ tel que 
$u= {\rm Int}_{G^\natural}(\delta)^{-1}(u'^{-1})u'$. On a donc $k^{-1}\cdot\delta \cdot uk=k^{-1}u'^{-1}\cdot\delta
\cdot u'k\in {^GM^\natural}\cap G^\natural_{\rm qr}$, d'o\`u le point (1).

Montrons (2). Soit $\Theta$ la fonction sur $T^\natural\cap G^\natural_{\rm reg}$ d\'efinie par
$$
\Theta(\gamma)=\sum_{i=1}^s\sum_{w\in \smash{\overline{W}}^\natural_i}\omega(g_in_w)^{-1}\vert D_{G^\natural/M}(\gamma^{g_i n_w})\vert_F^{-{1\over 2}}\Theta_\Sigma(\gamma^{g_i n_w }).\leqno{(**)}
$$
Puisque pour $\delta\in M^\natural\cap G^\natural_{\rm reg}$ et $m\in M$, on a
$$\omega(m)^{-1}\vert D_{G^\natural/M}(\delta^m)\vert_F^{-{1\over 2}}\Theta_\Sigma(\delta^m)=
\vert D_{G^\natural/M}(\delta)\vert_F^{-{1\over 2}}\Theta_\Sigma(\delta),
$$ l'expression \`a droite de l'\'egalit\'e dans $(**)$ ne d\'epend pas du choix de $\{g_1,\ldots ,g_s\}$ ni de celui de 
$\{n_w:w\in \smash{\overline{W}}^\natural_i\}$; d'ailleurs elle ne d\'epend pas non plus du choix de 
$\{T_1^\natural,\ldots ,T_s^\natural\}$. On en d\'eduit que si $\gamma\in T^\natural 
\cap G^\natural_{\rm reg}$ et $g\in G$ sont tels que 
$g^{-1}\cdot \gamma \cdot g\in T^\natural$, puisque $g\in N_G(T^\natural)$, on a 
$$\Theta(g^{-1}\cdot \gamma \cdot g)=\omega(g)\Theta(\gamma).$$
On peut donc prolonger $\Theta$ \`a ${^GT^\natural}\cap G^\natural_{\rm reg}$ en posant
$$
\Theta(g^{-1}\cdot \gamma \cdot g)=\omega(g)\Theta(\gamma)\quad (\gamma\in T^\natural\cap G^\natural_{\rm reg},\, g\in G).
$$
Pour $\phi\in C^\infty_{\rm c}({^GT^\natural}\cap G^\natural_{\rm reg})$, $\phi\geq 0$, la fonction 
$\vert \Theta\vert \phi$ sur ${^GT^\natural}\cap G^\natural_{\rm reg}$, prolong\'ee par 
$0$ \`a $G^\natural$ tout entier, est mesurable (par rapport \`a $d\delta$), 
et d'apr\`es le corollaire 1, on a
\begin{align*}
\lefteqn{
\int_{G^\natural}\vert\Theta(\delta)\vert\phi(\delta)d\delta =}\\
 &= {1\over \vert W^\natural \vert}\int_{\smash{\overline{T}}^\natural}\vert D_{G^\natural/T}(\gamma)\vert_F \vert\Theta(\gamma)\vert 
\Lambda^G_{\vert \omega\vert}(\phi,\gamma)d\bar{\gamma}\\
& \leq 
{1\over \vert W^\natural\vert}\sum_{i=1}^s
\sum_{w\in\smash{\overline{W}}^\natural_i}\int_{\smash{\overline{T}}^\natural}\vert \omega(\gamma^{g_i n_w})\vert^{-1}{\vert D_{G^\natural/T}(\gamma)\vert_F\over \vert D_{G^\natural/M}(\gamma^{g_i n_w})\vert_F^{1\over 2}}\vert\Theta_\Sigma(\gamma^{g_i n_w})\vert \Lambda^G(\phi,\gamma)d\bar{\gamma};
\end{align*}
o\`u l'on a pos\'e $W^\natural=W(G,T^\natural)$. Pour $i=1,\ldots ,s$, $w\in \smash{\overline{W}}^\natural_i$ et 
$\gamma\in T^\natural\cap G_{\rm reg}$, on a
$$
\vert \omega (\gamma^{g_i n_w})\vert^{-1}\Lambda^G_{\vert\omega\vert}(\phi,\gamma)=
\Lambda^G_{\vert\omega\vert}(\phi,\gamma^{g_i n_w}),
$$
pourvu que les mesures de Haar $ds_i$ sur $S_i=Z(T^\natural_i)$ aient \'et\'e choisies de 
mani\`ere compatible, ce que l'on suppose. Alors d'apr\`es la proposition de \ref{descente parabolique}, on a
$$
\Lambda^G_{\vert\omega\vert}(\phi,\gamma^{g_i n_w})= \vert D_{G^\natural/M}(\gamma^{g_i n_w})\vert_F^{-{1\over 2}}\Lambda^M_{\vert \omega\vert}({^\omega\phi_{P^\natural,K_\circ}},\gamma^{g_i n_w}).
$$
Comme
$$
\vert D_{G^\natural/T}(\gamma)\vert_F = \vert D_{G^\natural/T}(\gamma^{g_i n_w})\vert_F= 
\vert D_{G^\natural/M}(\gamma^{g_i n_w})\vert_F \vert
D_{M^\natural/T}(\gamma^{g_i n_w})\vert_F
$$
et
$$
\vert W^\natural\vert = \vert W(G,T^\natural_i)\vert,
$$ on obtient
\begin{align*}
\lefteqn{
\int_{G^\natural}\vert\Theta(\delta)\vert\phi(\delta)d\delta }\\
 & \leq {1\over \vert W(G,T^\natural)\vert}\sum_{i=1}^s
\int_{\smash{\overline{T}}^\natural}\vert D_{M^\natural/T}(\gamma^{g_i n_w})\vert_F
\vert\Theta_\Sigma(\gamma^{g_i n_w})\vert \Lambda^M_{\vert\omega\vert}({^\omega\phi_{P^\natural,K_\circ}},\gamma^{g_i n_w})d\bar{\gamma}\\
& =\sum_{i=1}^s {1\over \vert W(G,T^\natural_i)\vert}\sum_{w\in\smash{\overline{W}}^\natural_i}
\int_{\smash{\overline{T}}^\natural_i}\vert D_{M^\natural/T}(\gamma_i^{n_w})\vert_F
\vert\Theta_\Sigma(\gamma_i^{n_w})\vert \Lambda^M_{\vert\omega\vert}({^\omega\phi_{P^\natural,K_\circ}},\gamma_i^{n_w})d\bar{\gamma}_i\\
& =\sum_{i=1}^s {1\over \vert W(M,T^\natural_i)\vert}
\int_{\smash{\overline{T}}^\natural_i}\vert D_{M^\natural/T}(\gamma_i)\vert_F
\vert\Theta_\Sigma(\gamma_i)\vert \Lambda^M_{\vert \omega\vert}({^\omega\phi_{P^\natural,K_\circ}},\gamma_i)d\bar{\gamma}_i,
\end{align*}
d'o\`u (\`a nouveau d'apr\`es le corollaire 1)
$$
\int_{G^\natural}\vert\Theta(\delta)\vert\phi(\delta)d\delta\leq \int_{{^GT}^\natural \cap M^\natural_{\rm reg}}
\vert\Theta_\Sigma(\delta_M)\vert{^\omega\phi_{P^\natural,K_\circ}}(\delta_M)d\delta_M.
$$
La fonction $\Theta \phi$ est donc int\'egrable sur $G^\natural$, et v\'erifie
$$
\int_{G^\natural}\Theta(\delta)\phi(\delta)d\delta= \int_{{^GT}^\natural\cap M^\natural}\Theta_\Sigma(\delta_M){^\omega\phi_{P^\natural,K_\circ}}(\delta_M)d\delta_M=\Theta_\Sigma({^\omega\phi_{P^\natural,K_\circ}}).
$$
On conclut gr\^ace au th\'eor\`eme \ref{induction parabolique et caractres}.
\end{proof}

\begin{marema3}
{\rm 
Pour les caract\`eres non tordus des repr\'esentations admissibles, de type fini et {\it unitarisables}, 
le corollaire 2 est d\^u \`a Van Dijk \cite{VD}. L'hypoth\`ese d'unitarisabilit\'e a ensuite \'et\'e supprim\'ee 
par Clozel \cite{Cl1}. 
D'ailleurs dans loc. cit., le corollaire 2 est aussi dmontr dans le cas tordu 
suivant (changement de base pour le groupe lin\'eaire): $G$ est le groupe des points $F$--rationnels 
de ${\bf G}={\rm Res}_{E/F}(\Bbb{GL}_{n/F}\times_F E)$ pour un entier $n\geq 1$ et une extension finie cyclique 
$E$ de $F$ --- on a donc $G={\rm GL}_n(E)$ ---, 
$G^\natural = G\theta$ o\`u $\theta$ est le $F$--automorphisme de ${\bf G}$ donn\'e par un g\'en\'erateur 
du groupe ${\rm Gal}(E/F)$, et $\omega=1$.\hfill $\blacksquare$
}
\end{marema3}


\appendix

\section{Reprsentations irrductibles d'un $G$--espace tordu}

Dans cette annexe A, on fixe un groupe topologique 
localement profini $G$, un $G$--espace tordu $G^\natural$ et un caract\`ere $\omega$ de $G$. 
On fixe aussi une mesure de Haar \`a gauche $d_lg$ sur $G$, et l'on note 
$d_l\delta\;(=\delta \cdot d_lg)$ la mesure de Haar 
\`a gauche sur $G^\natural$ associ\'ee \`a $d_lg$ (cf. \ref{module d'un espace tordu}). Sauf pr\'ecision contraire, les 
modules consid\'er\'es sont des modules {\it \`a gauche}.

\subsection{Rappels sur les repr\'esentations (lisses) irr\'eductibles de $G$}\label{rappels sur les reprsentations irrductibles de G}
Soit $\ES{H}=\ES{H}(G)$ 
l'alg\`ebre de Hecke de $G$, 
i.e. l'espace $C^\infty_{\rm c}(G)$ muni du produit de convolution d\'efini par
$$
f*f'(x)=\int_Gf(g)f'(g^{-1}x)d_lg\quad (f,\, f'\in C^\infty_{\rm c}(G);\, x\in G).
$$
C'est une $\Bbb{C}$--alg\`ebre \`a idempotents, en g\'en\'eral sans unit\'e. 
Si $(\pi,V)$ est une repr\'esentation lisse de $G$, l'espace $V$ est naturellement muni d'une structure 
de $\ES{H}$--module: pour $f\in \ES{H}$ et $v\in V$, on pose $f\cdot v= \pi(f)(v)$, cf. \ref{caractres}. 
On obtient ainsi un isomorphisme entre la cat\'egorie des repr\'esentations 
lisses de $G$ et une sous--cat\'egorie pleine de la cat\'egorie des $\ES{H}$--modules: 
celle form\'ee des $\ES{H}$--modules $V$ {\it non d\'eg\'en\'er\'es}, c'est--\`a--dire tels que $\ES{H}\cdot V=V$. 

Soit $K$ un sous--groupe ouvert compact de $G$. On note $e_K$ la fonction caract\'eristique de $K$ 
divis\'ee par ${\rm vol}(K,d_lg)$; c'est un idempotent de $\ES{H}$. 
Notons $\ES{H}_K=\ES{H}(G,K)$ l'alg\`ebre de Hecke $e_K*\ES{H}*e_K$, 
i.e. l'espace $C_{\rm c}(K\backslash G/K)$ muni du produit de convolution d\'efini comme ci-dessus. 
C'est une $\Bbb{C}$--alg\`ebre \`a unit\'e ($e_K$ est l'unit\'e). 
Si $(\pi,V)$ est une repr\'esentation lisse de $G$, le sous--espace $V^K$ de $V$ form\'e des vecteurs 
$K$--invariants co\"{\i}ncide avec $\pi(e_K)(V)$. De plus la structure de $\ES{H}$--module sur $V$ induit 
une structure de $\ES{H}_K$--module sur $V^K$.

Si $(\pi,V)$ est une repr\'esentation lisse de $G$, on a $V=\bigcup_K V^K$ o\`u 
$K$ parcourt les sous--groupes ouverts compacts de $G$ (resp. un syst\`eme fondamental de voisinages de 
$1$ dans $G$ form\'e de tels sous--groupes). D'ailleurs, puisque $\ES{H}=\bigcup_K e_K* \ES{H}$, un $\ES{H}$--module $V'$ est non 
d\'eg\'en\'er\'e si et seulement si $V'=\bigcup_K e_K\cdot V'$. Si $(\pi,V)$ est une repr\'esentation lisse irr\'eductible de 
$G$, alors pour tout sous--groupe ouvert compact $K$ de $G$, le 
$\ES{H}_K$--module $V^K$ est nul ou simple \cite[prop.~2.10]{BZ}. D'autre part (loc.~cit.), fix\'e $K$, l'application $(\pi,V)\mapsto V^K$ 
induit une bijection entre:
\begin{itemize}
\item  l'ensemble des classes d'isomorphisme de repr\'esentations lisses irr\'eductibles 
de $G$ ayant un vecteur non nul fix\'e par $K$;
\item l'ensemble des classes d'isomorphisme de 
$\ES{H}_K$--modules simples.
\end{itemize}
Cette bijection induit par restriction une bijection entre:
\begin{itemize}
\item  l'ensemble des classes d'isomorphisme de repr\'esentations admissibles irr\'eductibles 
de $G$ ayant un vecteur non nul fix\'e par $K$;
\item l'ensemble des classes d'isomorphisme de 
$\ES{H}_K$--modules simples de dimension finie.
\end{itemize}
Tout comme pour les groupes finis, 
 le caract\`ere d'une repr\'esentation admissible irr\'eductible de $G$ d\'etermine cette repr\'esentation 
 \`a isomorphisme pr\`es \cite[prop.~2.19]{BZ}:
 
\begin{mapropo}
Soit $\pi_1,\,\pi_2,\,\ldots,\pi_n$ des repr\'esentations 
admissibles irr\'eductibles de $G$ deux--\`a--deux non isomorphes. Alors les distributions 
$\Theta_{\pi_1},\, \Theta_{\pi_2},\,\ldots,\, \Theta_{\pi_n}$ sont lin\'eairement ind\'ependantes.
\end{mapropo}

\begin{moncoro}
Deux repr\'esentations 
admissibles irr\'eductibles $\pi_1$ et $\pi_2$ de $G$ sont isomor\-phes si et seulement si 
$\Theta_{\pi_1}=\Theta_{\pi_2}$.
\end{moncoro}

\subsection{$\omega$--reprsentations $G$--irrductibles de $G^\natural$.}\label{représentations G-irréductibles} 
Cette annexe A a pour principal objectif la généralisation de la proposition de \ref{rappels sur les reprsentations irrductibles de G} aux $\omega$--représentations admissibles de $G^\natural$, précisément l'indépendance linéaire des distributions $\Theta_\Pi$ pour $\Pi$ parcourant les $\omega$--représentations admissibles $G$--irréductibles de $G^\natural$ (proposition de \ref{ind\'ependance lin\'eaire des caract\`eres tordus}); où par {\it $G$--irréductible} on entend une $\omega$--représentation (lisse) $\Pi$ de $G^\natural$ telle que la représentation $\Pi^\circ$ de $G$ est irréductible. On pourrait donc se limiter à l'étude des $\omega$--représentations admissibles $G$--irréductibles de $G^\natural$. On va voir qu'on peut facilement décrire une classe beaucoup plus large d'objets simples de la catégorie $\mathfrak{R}(G^\natural,\omega)$ en termes des représentations irréductibles de $G$. 

Une $\omega$--repr\'esentation lisse non nulle $(\Pi,V)$ de $G^\natural$ est dite {\it irr\'eductible} si c'est un objet simple dans la catégorie  $\mathfrak{R}(G^\natural,\omega)$, autrement si le seul sous--espace 
non nul $V'$ de $V$ tel que $\Pi(\delta)(V')=V'$ pour tout élément $\delta\in G^\natural$ --- ou, ce qui revient au m\^eme, tel que 
$V'$ est $G$--stable et $\Pi(\delta)(V')=V'$ pour un \'el\'ement $\delta\in G^\natural$ --- est $V$ lui--m\^eme. Bien sûr, si une $\omega$--représentation $\Pi$ de $G^\natural$ est $G$--irréductible, alors elle est automatiquement irréductible; mais l'inverse n'est en général pas vrai. 

Fixons un \'el\'ement $\delta_1\in G^\natural$ 
et posons $\theta ={\rm Int}_{G^\natural}(\delta_1)$. 

Si $\pi$ est une repr\'esentation lisse 
de $G$, pour $k\in \Bbb{Z}_{\geq 1}$, on pose $\frak{N}_{\theta,k}(x)=x \theta(x)\cdots \theta^{k-1}(x)$ ($x\in G$), 
$\omega_k=\omega\circ \frak{N}_{\theta,k}$ --- c'est un caract\`ere de $G$, ind\'ependant du choix de $\delta_1$ ---, 
et l'on note $\pi(k)$ la repr\'esentation 
$\omega_k^{-1}\pi^{\theta^k}$ de $G$. \`A isomorphisme pr\`es, $\pi(k)$ ne 
d\'epend pas du choix de $\delta_1$ dans $G^\natural$. Pour $k\in \Bbb{Z}_{\leq 1}$, 
on note $\pi(-k)$ la repr\'esentation lisse de $G$ telle que $\pi(-k)(k)=\pi(0)$. Posant 
$\pi(0)=\pi$, on a $\pi(k)(k')=\pi(k+k')$ ($k,\,k'\in \Bbb{Z}$). Si $\pi$ est irr\'eductible et s'il existe un entier 
$k\geq 1$ tel que $\pi(k)\simeq \pi$, on note $s(\pi)$ le plus petit entier $k\geq 1$ v\'erifiant cette propri\'et\'e; 
sinon on pose $s(\pi)=+\infty$.

Soit $(\Pi,V)$ une $\omega$--repr\'esentation 
lisse irr\'eductible de $G^\natural$. La reprsentation $\pi=\Pi^\circ$ de $G$ n'est en général pas irrductible (on l'a dit plus haut), ni mme de type fini. Supposons qu'il existe une sous--représentation irréductible $(\pi_0,V_0)$ de $(\pi,V)$. Pour $k\in \Bbb{Z}$, posons $V_k=\Pi(\delta_1)^{-k}(V_0)$; 
c'est un sous--espace $G$--stable de $V$. 
Puisque $\pi_0$ est irr\'eductible, si $V_0\cap V_k\neq 0$, alors $V_k=V_0$ 
et $\Pi(\delta_1)^k\vert_{V_0}$ est un isomorphisme de $\pi_0$ sur $\pi_0(k)$. 
Comme d'autre part $\Pi$ est irr\'eductible, les $V_k$ ($k\in \Bbb{Z}$) engendrent 
$V$ sur $\Bbb{C}$. On a donc deux cas possibles: ou bien $s(\pi_0)=+\infty$ et $V=\oplus_{k\in \Bbb{Z}}V_k$; ou bien 
$s(\pi_0)\neq +\infty$ et $V=\oplus_{k=0}^{s(\pi_0)-1}V_k $. Pour $k\in \Bbb{Z}$, notons 
$\pi_k$ la restriction de $\pi$ \`a $V_k$. La restriction de $\Pi(\delta_1)^{-k}$ \`a $V_0$ induit un 
isomorphisme de $\pi_0(k)$ sur $\pi_k$. On en d\'eduit que $s(\pi_k)=s(\pi_0)$. L'invariant $s(\pi_0)$ ne d\'epend donc 
pas du choix de $\pi_0$ --- ni du choix de $\delta_1$ comme on l'a dit plus haut; 
on le note $s(\Pi)$. D'apr\`es ce qui pr\'ec\`ede, 
la repr\'esentation $\pi$ de $G$ est de type fini si et seulement si $s(\Pi)<+\infty$, et elle est est irrductible si et seulement si $s(\Pi)=0$. 

Supposons $s(\Pi)=+\infty$. Soit $(\Pi',V')$ la $\omega$--repr\'esentation lisse de $G^\natural$ d\'efinie 
par
$$\Pi'^\circ=\oplus_{k\in \Bbb{Z}}\pi_0(k)
$$ et 
$$\Pi'(\delta_1)((v_k)_{k\in \Bbb{Z}})= (v'_k)_{k\in \Bbb{Z}},\quad v'_k= v_{k+1}.
$$
Alors l'application
$$(V_0)^\Bbb{Z}\rightarrow V,\,(v_k)_{k\in \Bbb{Z}}\mapsto \sum_{k\in \Bbb{Z}}\Pi(\delta_1)^{-k}(v_k)$$ 
est un isomorphisme de $(\Pi',V')$ sur $(\Pi,V)$.

Supposons maintenant $s(\Pi)<+\infty$. Posons $s=s(\Pi)$. 
Soit $(\Pi',V')$ la $\omega$--repr\'esentation lisse de $G^\natural$ 
d\'efinie par
$$\Pi'^\circ = \oplus_{k=0}^{s-1}\pi_0(k)$$ et
$$\Pi'(\delta_1)(v_0,\ldots ,v_{s-1})=(v_1,\ldots ,
v_{s-1},\Pi(\delta_1)^s(v_0)).
$$ Alors l'application
$$(V_0)^s\rightarrow \oplus_{k=0}^{s-1}V_i,\,(v_0,\ldots ,v_{s-1})\mapsto \sum_{k=0}^{s-1}\Pi(\delta_1)^{-k}(v_k)$$ 
est un isomorphisme de $(\Pi',V')$ sur $(\Pi,V)$.

\v1 Récapitulons. Soit $(\Pi,V)$ une $\omega$--représentation lisse irréductible de $G$ telle qu'il existe une sous--représentation irréductible $\pi_0$ de $\Pi^\circ$. Alors l'invariant $s(\Pi)= s(\pi_0)\in {\Bbb Z}_{\geq 1} \cup \{+\infty\}$ ne dépend pas de $\pi_0$ (il dépend seulement de la classe d'isomorphisme de $\Pi$), et on a
$$
\Pi^\circ \simeq \left\{\begin{array}{ll}\oplus_{k\in {\Bbb Z}}\pi_0(k)\hfill  & \hbox{si $s(\Pi)=+\infty$}\\
\oplus_{k=0}^{s(\Pi)-1}\pi_0(k) & \hbox{sinon}
\end{array}\right..
$$ 
En particulier la représentation $\Pi^\circ$ est semisimple, et elle est de type fini (i.e. de longueur finie) si et seulement si $s(\Pi)<+\infty$. 

\begin{marema}
{\rm 
On vient de voir que pour une $\omega$--représentation lisse irréductible $(\Pi,V)$ de $G^\natural$, 
les deux conditions suivantes sont équivalentes:
\begin{itemize}
\item il existe une sous--représentation irréductible de $\Pi^\circ$;
\item la représentation $\Pi^\circ$ de $G$ est semisimple.
\end{itemize}
On verra en \ref{reprsentations irrductibles et H-modules simples} que si une certaine propriété $({\rm P}_2)$ du $G$--espace tordu $G^\natural$ est vérifiée\footnote{On verra aussi (cf. \ref{la condition P2}) que cette propriété $({\rm P}_2)$ est toujours vérifiée si $G={\bf G}(F)$ et $G^\natural={\bf G}(F)$ pour un groupe réductif connexe $\bf{G}$ défini sur $F$ et un $\bf{G}$--espace tordu défini sur $F$ et possédant un point $F$--rationnel; où $F$ est un corps commutatif localement compact non archimédien.}, alors 
les deux conditions ci--dessus sont automatiquement satisfaites si $\Pi^\circ$ est admissible. Il est possible qu'elles le soient 
pour une classe beaucoup plus large de $\omega$--représentations lisses irréductibles de $G^\natural$, au moins 
sous certaines hypothèses de finitude (par exemple s'il existe un entier $l\geq 1$ tel que $\theta^l={\rm Int}_G(g)$ pour un élément 
$g\in G$, auquel cas on peut identifier $G^\natural$ à un sous--groupe de $(G\rtimes \langle \theta \rangle)/C$ comme en \ref{automorphismes ss et u}; où $C$ est le sous--groupe distingué de $G\rtimes \langle \theta \rangle$ engendré par $g^{-1} \rtimes \theta^l$). Mais il semble vain d'espérer qu'elles le soient en général\footnote{On verra (\ref{H-modules}, et remarque 1 de \ref{reprsentations irrductibles et H-modules simples}) qu'une $\omega$--représentation lisse de $G^\natural$ n'est autre qu'un module (non dégénéré) $V$ sur l'agèbre $\widetilde{\ES{H}}^\natural$ des polynômes de Laurent sur $\ES{H}$ tordus par un automorphisme $\theta_{\ES{H}}$ de $\ES{H}$, la représentation $\Pi^\circ$ de $G$ correspondant au $\ES{H}$--module $V$. \`A moins d'imposer certaines conditions sur l'automorphisme $\theta_{\ES{H}}$, le fait que le $\widetilde{\ES{H}}^\natural$--module $V$ soit simple n'implique pas qu'il soit semisimple comme $\ES{H}$--module.}.\hfill $\blacksquare$ }
\end{marema}

\subsection{$(\ES{H}^\natural,\omega)$--modules et $(\ES{H}^\natural_K,\omega)$--modules.}
\label{H-modules}
Notons $\ES{H}^\natural=\ES{H}(G^\natural)$ 
l'espace $C^\infty_{\rm c}(G^\natural)$ muni de la structure de $\ES{H}$--bimodule donn\'ee par 
$(f\in \ES{H},\,\phi\in C^\infty_{\rm c}(G^\natural),\,\delta\in G^\natural)$:
$$
f*\phi(\delta) =\int_Gf(g)\phi(g^{-1}\cdot \delta)d_lg,\quad
\phi* f (\delta)=\int_G\phi(\delta\cdot g)f(g^{-1})d_lg;
$$
puisque $d_l(g^{-1})= \Delta_G(g^{-1})d_lg$, on a aussi
$$
\phi* f(\delta)=\int_G\phi(\delta\cdot g^{-1})f(g)\Delta_G(g^{-1})d_lg.
$$ 
On appelle {\it $(\ES{H}^\natural,\omega)$--module} un $\ES{H}$--module $V$ muni d'une application 
$$\ES{H}^\natural\rightarrow {\rm End}_\Bbb{C}(V),\,\phi\mapsto (v\mapsto \phi\cdot v)$$ 
telle que 
$$
(f*\phi* f') \cdot v= f\cdot (\phi\cdot (\omega f'\cdot v))\quad (\phi\in \ES{H}^\natural;\, f,\,f'\in \ES{H};\, v\in V).
$$
Les $(\ES{H}^\natural,\omega)$--modules forment une sous--cat\'egorie (non pleine) de la cat\'egorie des 
$\ES{H}$--modules: un morphisme entre deux $(\ES{H}^\natural,\omega)$--modules $V_1$ et $V_2$ est 
simplement un morphisme de $\ES{H}$--modules $u:V_1\rightarrow V_2$ tel que $u(\phi\cdot v)=\phi\cdot u(v)$ 
pour tout $\phi\in \ES{H}^\natural$ et tout $v\in V_1$. 

Un $(\ES{H}^\natural,\omega)$--module $V$ 
est dit {\it non dgnr} si $\ES{H}^\natural\cdot V=V$ (puisque $\ES{H}^\natural= \ES{H}*\ES{H}^\natural$, le $\ES{H}$--module 
sous--jacent est lui aussi non d\'eg\'en\'er\'e). Les $(\ES{H}^\natural,\omega)$--modules non dgnrs forment une sous--catgorie pleine 
de la catgorie des $(\ES{H}^\natural,\omega)$--modules. Notons qu'un morphisme entre deux $(\ES{H}^\natural,\omega)$--modules non 
d\'eg\'en\'er\'es $V_1$ et $V_2$ est une application $\Bbb{C}$--lin\'eaire $u:V_1\rightarrow V_2$ 
telle que $u(\phi\cdot v)=\phi\cdot u(v)$ pour tout $\phi\in \ES{H}^\natural$ et tout $v\in V_1$ (une telle 
application est automatiquement $\ES{H}$--lin\'eaire). Si $(\Pi,V)$ 
est une $\omega$--repr\'esentation lisse de $G^\natural$, l'espace $V$ est naturellement muni d'une 
structure de $(\ES{H}^\natural,\omega)$--module: pour $\phi\in \ES{H}^\natural$ et $v\in V$, on pose 
$\phi\cdot v= \Pi(\phi)(v)$, cf. \ref{caractres tordus}. On v\'erifie que pour $f,\, f'\in \ES{H}$, on a 
$\Pi(f*\phi*f')=\pi(f)\circ \Pi(\phi)\circ \pi(\omega f')$, o\`u l'on a pos\'e $\pi=\Pi^\circ$. 
On obtient ainsi un isomorphisme entre:
\begin{itemize}
\item la cat\'egorie des $\omega$--repr\'esentations lisses 
de $G^\natural$;
\item la catgorie des $(\ES{H}^\natural,\omega)$--modules non d\'eg\'en\'er\'es.
\end{itemize}

\v1
Soit $K$ un sous--groupe ouvert compact de $G$. On note 
$\ES{H}^\natural_K$ l'espace $C_{\rm c}(K\backslash G^\natural /K)$ 
muni de la structure 
de $\ES{H}_K$--bimodule d\'efinie comme plus haut. En d'autres termes, on pose $\ES{H}^\natural_K=e_K*\ES{H}^\natural * e_K$. 
On a $\ES{H}^\natural_K=\ES{H}_K*\ES{H}_K^\natural =\ES{H}^\natural_K* \ES{H}_K$. Si de plus il existe un \'el\'ement $\delta\in G^\natural$ 
normalisant $K$, i.e. tel que ${\rm Int}_{G^\natural}(\delta)(K)=K$, alors $\ES{H}^\natural_K$ est un $\ES{H}_K$--module 
(\`a gauche ou \`a droite) libre de rang $1$, engendr\'e par la fonction caract\'eristique de $K\cdot \delta=\delta\cdot K$. 

Supposons de plus que le caractre $\omega$ est trivial sur $K$. 
On appelle encore {\it $(\ES{H}^\natural_K,\omega)$--module} un $\ES{H}_K$--module $W$ muni 
d'une application
$$\ES{H}^\natural_K\rightarrow {\rm End}_\Bbb{C}(W),\,\phi\mapsto (w\mapsto \phi\cdot w)$$
telle que 
$$
(f*\phi *f')\cdot w=f\cdot(\phi\cdot (\omega f' \cdot w))\quad (\phi\in \ES{H}^\natural_K;\, f,\,f'\in \ES{H}_K;\,w\in W).
$$
Les notions de morphismes entre deux $(\ES{H}^\natural_K,\omega)$--modules et de $(\ES{H}^\natural_K,\omega)$--module 
non d\'eg\'en\'er\'e sont d\'efinies comme plus haut. 
Si $(\Pi,V)$ est une $\omega$--repr\'esentation lisse de $G^\natural$, l'espace $V^K$ est 
naturellement muni d'une structure de $(\ES{H}^\natural_K,\omega)$--module (celle d\'eduite par restriction du 
$(\ES{H}^\natural,\omega)$--module $V$). Puisque $\ES{H}_K^\natural = \ES{H}_K^\natural*e_K$, on a 
$\ES{H}^\natural_K\cdot V=\ES{H}^\natural_K\cdot V^K$. Posons $\pi=\Pi^\circ$, et 
pour $\gamma\in G$, notons $\phi^K_\gamma$ la fonction caract\'eristique de $K\cdot \gamma \cdot K$ 
divis\'ee par ${\rm vol}(K\cdot \gamma \cdot K,d_l\delta)$. 
Pour $v\in V$, on a 
$$\phi^K_\gamma \cdot v = \pi (e_K)\circ \Pi(\gamma)\circ \pi(e_K)(v).
$$

\begin{marema}
{\rm Si de plus il existe un \'el\'ement 
$\delta\in G^\natural$ normalisant $K$, alors 
$\Pi(\delta)$ induit par restriction un automorphisme de $V^K$, qui co\"{\i}ncide 
avec $w\mapsto \phi^K_{\delta}\cdot w$. En ce cas on a l'égalité $\ES{H}^\natural_K\cdot V^K=V^K$, 
i.e. le $(\ES{H}^\natural_K,\omega)$--module 
$V^K$ est automatiquement non d\'eg\'en\'er\'e. 
En g\'en\'eral, on a seulement l'inclusion $\ES{H}^\natural_K\cdot V^K\subset V^K$. Notons que puisque $\ES{H}^\natural\cdot V=V$ et 
$\ES{H}^\natural = \bigcup_{K'}\ES{H}^\natural_{K'}$ 
o\`u $K'$ parcourt les sous--groupes ouverts compacts de $G^\natural$ tels que $\omega\vert_{K'}=1$, 
on a toujours l'égalité $V=\bigcup_{K'}\ES{H}^\natural_{K'}\cdot V^{K'}$.\hfill $\blacksquare$
}
\end{marema}

Continuons avec le sous--groupe ouvert compact $K$ de $G$. On suppose 
toujours $\omega\vert_K=1$. On suppose aussi qu'il existe un \'el\'ement $\delta_1\in G^\natural$ normalisant $K$, et on pose 
$\theta ={\rm Int}_{G^\natural}(\delta_1)$. Pour $f\in \ES{H}$, on note ${^{\theta\!}f}\in \ES{H}$ la fonction $f\circ \theta^{-1}$. 
Puisque $d_l\theta(g)=d_lg$, pour $f,\,f'\in \ES{H}$, on a ${^\theta(f*f')}={^{\theta\!}f}*{^{\theta\!}f'}$, i.e. 
l'application $\ES{H}\rightarrow \ES{H},\,f\mapsto{^{\theta\!}f}$ est un automorphisme d'alg\`ebres. Par restriction, on 
obtient un automorphisme d'alg\`ebres $\ES{H}_K\rightarrow \ES{H}_K,\,f\mapsto {^{\theta\!}f}$. Pour $g\in G$, notons 
$f^K_g$ la fonction caract\'eristique de $KgK$ divis\'ee par ${\rm vol}(K,d_lg)$. Puisque
$$
\delta_1 \cdot KgK = 
K\cdot \delta_1 \cdot gK= K\theta(g)\cdot \delta_1 \cdot K= K\theta(g)K\cdot \delta_1,
$$
on a
$$
\phi^K_{\delta_1} * f^K_g = \phi^K_{\delta_1\cdot g}=\phi^K_{\theta(g)\cdot \delta_1}=f^K_{\theta(g)}*\phi^K_{\delta_1}.
$$
Par lin\'earit\'e on obtient l'\'egalit\'e
$$
\phi^K_{\delta_1}* f ={^{\theta\!}f} * \phi^K_{\delta_1}\quad (f\in \ES{H}_K).
$$
On en d\'eduit qu'un $(\ES{H}^\natural_K,\omega)$--module $W$ est non d\'eg\'en\'er\'e si et seulement si 
l'application $w\mapsto \phi^K_{\delta_1}\cdot w$ est un automorphisme de $W$. 

On suppose toujours $\omega\vert_K=1$ et $\theta(K)=K$, o\`u $\theta={\rm Ind}_{G^\natural}(\delta_1)$. 
D'apr\`es ce qui pr\'ec\`ede, la donn\'ee d'un $(\ES{H}^\natural_K,\omega)$--module $W$ \'equivaut \`a 
celle d'un {\it $(\ES{H}_K,\theta,\omega)$--module} $(W,\theta_W)$, c'est--\`a--dire un $\ES{H}_K$--module 
$W$ muni d'un $\Bbb{C}$--endomorphisme $\theta_W$ v\'erifiant l'\'egalit\'e
$$
\theta_W(\omega f\cdot w)= {^{\theta\!}f}\cdot \theta_W(v)\quad (f\in \ES{H}_W,\,w\in W).
$$
On passe d'un point de vue \`a l'autre gr\^ace \`a l'\'egalit\'e $\theta_W(w)=\phi^K_{\delta_1}\cdot w$ ($w\in W$), 
d'o\`u les notions (\'evidentes) de morphismes entre deux $(\ES{H}_K,\theta,\omega)$--modules et de 
$(\ES{H}_K,\theta,\omega)$--module non d\'eg\'en\'er\'e.

\subsection{$\omega$--repr\'esentations irr\'eductibles de $G^\natural$ et $(\ES{H}^\natural_K,\omega)$--modules simples}
\label{reprsentations irrductibles et H-modules simples} 
Si $K$ est un sous--groupe ouvert compact 
de $G$ tel que $\omega\vert_K=1$, un $(\ES{H}^\natural_K,\omega)$--module non nul et non d\'eg\'en\'er\'e 
$W$ est dit {\it simple} si c'est un objet simple dans la catégorie des $\ES{H}^\natural_K$--modules non dégénérés, \cad si le seul 
sous--$\ES{H}_K$--module non nul $W'$ de $W$ tel que $\ES{H}_K^\natural\cdot W'=W'$ 
est $W$ lui--m\^eme --- bien sûr, s'il existe un élément $\delta\in G^\natural$ normalisant $K$, alors la condition 
$\ES{H}^\natural_K\cdot W'=W'$ peut être remplacée par $\phi^K_\delta\cdot W'=W'$. 

Consid\'erons les propri\'et\'es $({\rm P}_1)$ et $({\rm P}_2)$ suivantes:

\begin{enumerate}
\item[$({\rm P}_1)$]Il existe un syst\`eme fondamental de voisinages de $1$ dans $G$ form\'e 
de sous--groupes ouverts compacts $K$ de $G$ 
tels que ${\rm Int}_{G^\natural}(\delta_K)(K)=K$ pour un \'el\'ement $\delta_K\in G^\natural$.
\item[$({\rm P}_2)$]Il existe un \'el\'ement $\delta_1\in G^\natural$ et un 
syst\`eme fondamental de voisinages de $1$ dans $G$ form\'e 
de sous--groupes ouverts compacts $K$ de $G$ tels que ${\rm Int}_{G^\natural}(\delta_1)(K)=K$. 
\end{enumerate}

On a clairement l'implication $({\rm P}_2)\Rightarrow ({\rm P}_1)$. 

{\bf On suppose que la propri\'et\'e 
$({\rm P}_1)$ est v\'erifi\'ee}. Fixons un syst\`eme fondamental de voisinages de $1$ dans $G$, 
disons $\ES{K}$, form\'e de sous--groupes ouverts compacts $K$ de $G$ tels que $\omega\vert_K=1$ et 
${\rm Int}_{G^\natural}(\delta_K)(K)=K$ pour un $\delta_K\in G^\natural$.

\begin{mapropo1}
(En supposant $({\rm P}_1)$.)
\begin{enumerate}
\item[(1)]Une $\omega$--repr\'esentation lisse $(\Pi,V)$ de $G^\natural$ est irr\'eductible si et seulement si 
$V\neq 0$ et pour tout $K\in \ES{K}$, le 
$(\ES{H}^\natural_K,\omega)$--module $V^K$ est nul ou simple.
\item[(2)]Soit $(\Pi_1,V_1)$ et $(\Pi_2,V_2)$ deux $\omega$--repr\'esentations lisses irr\'eductibles de $G^\natural$, 
et $K\in \ES{K}$ tel que $V_i^K\neq 0$ pour $i=1,\,2$. 
Alors les $\omega$--repr\'esentations $\Pi_1$ et $\Pi_2$ sont isomorphes si et seulement si les 
$(\ES{H}^\natural_K,\omega)$--modules 
$V_1^K$ et $V_2^K$ le sont.
\end{enumerate}
\end{mapropo1}

\begin{proof}
Il s'agit de recopier celle 
de \cite[prop.~2.10.(a)--(b)]{BZ}. Soit $(\Pi,V)$ une 
$\omega$--repr\'esen\-tation lisse de $G^\natural$, 
et soit $\pi=\Pi^\circ$. 
Soit $K\in \ES{K}$, et soit $W'$ un sous--$\ES{H}_K$--module de $V^K$ tel que $\ES{H}_K^\natural\cdot W'=W'$. Notons $V'=G^\natural\cdot W'$ le sous--espace de $V$ engendr\'e par les $\Pi(\delta)(w)$ pour 
$\delta\in G^\natural$ et $w\in W'$. Puisque $G^\natural = G\cdot \delta_K$ et $\Pi(\delta_K)(W')= W'$, $V'$ est aussi le sous--espace $G\cdot W$ de $V$ engendré par les $\pi(g)(w)$ pour $g\in G$ et $w\in W$, et on a $\Pi(\delta)(V')=V'$ pour tout élément $\delta\in G^\natural$. En d'autres termes, $V'$ définit une sous--$\omega$--représentation $\Pi'$ de $\Pi$. De plus, comme $\ES{H}^\natural_K\cdot V'=V'^K$ est engendr\'e par les $\phi^K_\delta\cdot w$ pour $\delta\in G^\natural$ 
et $w\in W$, on a l'égalité $V'^K=W'$.    

Montrons (1). Supposons $\Pi$ irr\'eductible, et soit $K\in \ES{K}$ tel que $V^K\neq 0$. 
Soit $W$ un sous--$\ES{H}_K$--module non nul de $V^K$ tel que $\ES{H}_K^\natural\cdot W=W$. 
Le sous--espace $G^\natural\cdot W$ de $V$ est non nul, et puisqu'il définit une sous--$\omega$--représentation de $\Pi$, 
on a $G^\natural\cdot W=V$. Comme $(G^\natural\cdot W)^K=W$, 
on obtient $W=V^K$. Donc le $(\ES{H}^\natural_K,\omega)$--module $V^K$ est simple. 
R\'eciproquement, supposons qu'il existe un $K\in \ES{K}$ tel que le $(\ES{H}^\natural_K,\omega)$--module 
$V^K$ n'est ni nul ni simple. Alors il existe un sous--$\ES{H}_K$--module non nul $W'$ de $V^K$, distinct de $V^K$, tel $\ES{H}^\natural_K\cdot W'= W'$. Le sous--espace $V'=G^\natural\cdot W$ de $V$ définit une sous--$\omega$--représentation non nulle de $\Pi$. Comme $V'^K=W'$, on a $V'\neq V$ et $\Pi$ n'est pas irr\'eductible. Cela ach\`eve la d\'emonstration 
du point (1).

Montrons (2). Si les $\omega$--repr\'esentations $\Pi_1$ et $\Pi_2$ sont isomorphes, alors les $(\ES{H}^\natural_K,\omega)$--modules 
$V_1^K$ et $V_2^K$ le sont aussi. R\'eciproquement, supposons qu'il existe un isomorphisme de 
$(\ES{H}^\natural_K,\omega)$--modules $u:V^K_1\rightarrow V^K_2$. Alors $W'=\{(w,u(w):w\in V_1^K\}$ est un sous--$\ES{H}_K$--module non nul de $V_1^K\times V_2^K$ tel que $\ES{H}^\natural_K\cdot W' =W'$, et $V'=G^\natural\cdot W'$ est un sous--espace non nul de $V_1\times V_2$ qui définit une sous--$\omega$--représentation de $\Pi_1\times \Pi_2$. Comme $V'^K=W'$, pour $i=1,\,2$, $V'$ ne peut pas contenir $V_i$, ni \^etre contenu dans 
$V_i$ (car sinon $W'$ contiendrait $V_i^K$, ou serait contenu dans $V_i^K$, ce qui est impossible). Puisque $\Pi_1$ 
et $\Pi_2$ sont irr\'eductibles, pour $i=1,\,2$, la projection $V'\rightarrow V_i$ est un isomorphisme de 
$(\ES{H}^\natural,\omega)$--modules. Les $\omega$--repr\'esentations $\Pi_1$ et $\Pi_2$ sont donc isomorphes.
\end{proof}

{\bf On suppose maintenant que la propri\'et\'e $({\rm P}_2)$ est v\'erifi\'ee}. Fixons un \'el\'ement 
$\delta_1\in G^\natural$ et un syst\`eme fondamental de voisinages de $1$ dans $G$, 
disons $\ES{K}_1$, form\'e de sous--groupes ouverts compacts $K$ de $G$ tels que 
$\omega\vert_K=1$ et $\theta(K)=K$, o\`u l'on a pos\'e
$$\theta={\rm Int}_{G^\natural}(\delta_1).$$ 

Tout $(\ES{H}^\natural,\omega)$--module $V$ d\'efinit, pour chaque $K\in \ES{K}_1$, un $(\ES{H}_K,\theta,\omega)$--module 
$(V^K,\theta_{V^K})$: pour $w\in V^K$, on a $\theta_{V^K}(w)=\phi^K_{\delta_1}\cdot w$. De plus, la famille 
$\{(V^K,\theta_{V^K}):K\in \ES{K}_1\}$ est compatible, 
au sens o\`u si $v\in V^K\cap V^{K'}$ pour des $K,\,K'\in \ES{K}_1$, on a $\theta_{V^K}(v)=\theta_{V^{K'}}(v)$. Il suffit pour cela 
de choisir un $K''\in \ES{K}_1$ tel que $K''\subset K\cap K'$, et de remarquer que $\phi^{K''}_{\delta_1}*e_K=\phi^K_{\delta_1}$ 
et $\phi^{K''}_{\delta_1}*e_{K'}=\phi^{K'}_{\delta_1}$. 

On d\'efinit comme en \ref{H-modules} la notion de {\it $(\ES{H},\theta,\omega)$--module non d\'eg\'en\'er\'e}. D'apr\`es ce qui pr\'ec\`ede, 
la donn\'ee d'un $(\ES{H}^\natural,\omega)$--module non d\'eg\'en\'er\'e $V$ \'equivaut \`a celle 
d'un $(\ES{H},\theta,\omega)$--module non d\'eg\'en\'er\'e $(V,\theta_V)$: pour $v\in V$, on a 
$\theta_V(v)=\phi^K_{\delta_1}\cdot v$ pour tout $K\in \ES{K}_1$ tel que $v\in e_K\cdot V=V^K$. 

\begin{marema1}
{\rm Notons $\theta_\ES{H}$ le $\Bbb{C}$--automorphisme de $\ES{H}$ donn\'e par 
$$
\theta_\ES{H}(f)= {^\theta(\omega^{-1} f)}\quad (f\in \ES{H}).
$$  
On v\'erifie que pour $f,\,f'\in \ES{H}$, on a
$$
\theta_\ES{H}(\omega f * f')={^{\theta\!}f}*\theta_\ES{H}(f').
$$
En d'autres termes, $(\ES{H},\theta_\ES{H})$ est un $(\ES{H},\theta,\omega)$--module non d\'eg\'en\'er\'e. D'ailleurs pour 
$K\in \ES{K}_1$, la restriction de $\theta_\ES{H}$ \`a $\ES{H}_K\;(\subset \ES{H}^K=e_K* \ES{H})$ induit un $\Bbb{C}$--automorphisme 
de $\ES{H}_K$, disons $\theta_{\ES{H}_K}$, et $(\ES{H}_K,\theta_{\ES{H}_K})$ est un $(\ES{H}_K,\theta,\omega)$--module 
non d\'eg\'en\'er\'e.

Notons aussi qu'un $\ES{H}^\natural$--module (non dégénéré) n'est autre qu'un module (non dégénéré) sur l'algèbre $\widetilde{\ES{H}}^\natural$ des polynômes de Laurent sur $\ES{H}$ tordus par l'automorphisme $\theta_{\ES{H}}$. \hfill $\blacksquare$}
\end{marema1}

\begin{marema2}
{\rm 
Si $(\Pi,V)$ est une $\omega$--représentation admissible de $G^\natural$, et si $K$ est un sous--groupe ouvert compact de $G$ tel que $V^K \neq 0$, alors puisque $\dim_{\Bbb C}(V^K)<+\infty$, il existe un sous--$\ES{H}_K$--module simple de $V^K$ (i.e. le socle du $\ES{H}_K$--module $V^K$ n'est pas nul).\hfill $\blacksquare$}
\end{marema2}

\begin{mapropo2}
(En supposant $({\rm P}_2)$.) Soit $K\in \ES{K}_1$, et soit 
$W$ un $(\ES{H}^\natural_K,\omega)$--module non d\'eg\'en\'er\'e simple tel que le socle du $\ES{H}_K$--module 
$W$ n'est pas nul. Alors le $\ES{H}_K$--module $W$ est semisimple, et 
il existe une $\omega$--repr\'esentation lisse irr\'eductible $(\Pi,V)$ de 
$G^\natural$ telle que le $(\ES{H}^\natural_K,\omega)$--module $V^K$ est isomorphe \`a $W$. De plus, la représentation 
$\Pi^\circ$ de $G$ est semisimple, et elle est de longueur finie (resp. irréductible) si et seulement si le $\ES{H}_K$--module 
$W$ est de longueur finie (resp. simple).
\end{mapropo2}

\begin{proof}
Il s'agit d'adapter celle de \cite[prop. 2.10.(c)]{BZ}. Posons $\theta_W(w)=\phi^K_{\delta_1}\cdot w$ 
($w\in W$) comme en \ref{H-modules}. Puisque $W$ est simple, 
le seul sous--$\ES{H}_K$--module $\theta_W$--invariant de $W$ est $W$ 
lui-m\^eme. Fixons un sous--espace non nul $W_0$ de $W$ qui soit un $\ES{H}_K$--module 
simple, et pour $i\in \Bbb{Z}$, posons $W_i=\theta_W^i(W_0)$, o\`u $\theta^i_W=\theta_W\circ\cdots \circ \theta_W$ ($i$ fois) si $i\geq 0$ et 
$\theta_W^i =(\theta_W^{-1})^{-i}$ sinon. Chaque $W_i$ est un sous--$\ES{H}_K$--module simple de $W$, 
et comme le sous--espace de $W$ engendr\'e par les $W_i$ pour $i\in \Bbb{Z}$ est \`a la fois 
$\ES{H}_K$--stable et $\theta_W$--invariant, c'est $W$ tout entier. On distingue deux cas: ou bien $W_0\cap W_i=\{0\}$ pour tout 
$i\in \Bbb{Z}\smallsetminus \{0\}$; ou bien il existe un plus petit entier $d\geq 1$ tel que $W_0\cap W_d\neq \{0\}$, auquel 
cas $W_d=W_0$. Fixons un \'el\'ement $w_0\in W_0\smallsetminus \{0\}$ et 
posons $\frak{J}_{K,0}=\{f\in \ES{H}_K:f\cdot w_0=0\}$. C'est un id\'eal \`a gauche dans $\ES{H}_K$, 
et $W_0$ est isomorphe (comme $\ES{H}_K$--module) \`a $\ES{H}_K/\frak{J}_{K,0}$. Pour $i\in \Bbb{Z}$, 
$\frak{J}_{K,i}=\theta_\ES{H}^i(\frak{J}_{K,0})$ est encore un id\'eal \`a gauche dans $\ES{H}_K$, et 
$W_i$ est isomorphe \`a $\ES{H}_K/\frak{J}_{K,i}$. Notons $\ES{A}_i=\ES{H}*\ES{H}_K$ 
(resp. $\frak{J}_i=\ES{H}*\frak{J}_{K,i}$) l'id\'eal 
\`a gauche dans $\ES{H}$ engendr\'e par $\ES{H}_K$ (resp. par $\frak{J}_{K,i}$), et posons $\ES{A}'_i= \ES{A}_i/\frak{J}_i$.

{\it Pla\c{c}ons--nous dans le premier cas: $W=\oplus_{i\in \Bbb{Z}}W_i$}. Notons $\ES{A}$ le $\ES{H}$--module $\oplus_{i\in \Bbb{Z}}\ES{A}_i$ 
(pour l'action diagonale de $\ES{H}$), et $\theta_\ES{A}$ le $\Bbb{C}$--automorphisme de 
$\ES{A}$ d\'efini 
par $\theta_\ES{A}({\bf f})_i = \theta_\ES{H}({\bf f}_{i-1})$ pour tout ${\bf f}=({\bf f}_i)_{i\in \Bbb{Z}}\in \ES{A}$. 
Pour $f\in \ES{H}$ et ${\bf f}=({\bf f}_i)_{i\in \Bbb{Z}}\in \ES{A}$, on a
$$\theta_\ES{A}(\omega f\cdot {\bf f})_i= \theta_\ES{H}(\omega f\cdot {\bf f}_{i-1})={^{\theta\!}f}\cdot \theta_\ES{A}({\bf f})_i. 
$$
En d'autres termes, $(\ES{A},\theta_\ES{A})$ est un $(\ES{H},\theta,\omega)$--module non d\'eg\'en\'er\'e.  
Posons $\frak{J}=\oplus_{i\in \Bbb{Z}}\frak{J}_i$. C'est un sous--$\ES{H}$--module 
$\theta_\ES{A}$--invariant de $\ES{A}$, et 
d'apr\`es le d\'ebut de la d\'emonstration de la proposition 1, on a $e_K\cdot \ES{A}=\oplus_{i\in \Bbb{Z}}\ES{H}_K$ et 
$e_K\cdot \frak{J}=\oplus_{i\in \Bbb{Z}}\frak{J}_{K,i}$. Par passage au quotient, $\theta_\ES{A}$ induit un 
$\Bbb{C}$--automorphisme $\theta_{\ES{A}'}$ de $\ES{A}'=\ES{A}/\frak{J}$, 
et $(\ES{A}',\theta_{\ES{A}'})$ est un $(\ES{H},\theta,\omega)$--module non d\'eg\'en\'er\'e tel que 
$e_K\cdot \ES{A}'\simeq W$. On a 
$\ES{A}'=\oplus_{i\in \Bbb{Z}}\ES{A}'_i$ et $\theta_{\ES{A}'}(\ES{A}'_{i-1})=\ES{A}'_i$ ($i\in {\Bbb Z}$). 
Soit $V'_0$ un sous--$\ES{H}$--module de $\ES{A}'_0$ tel que le $\ES{H}$--module 
quotient $V_0=\ES{A}'_0/V'_0$ est simple (puisque $\ES{A}'_0= \ES{H}\cdot f_0$ pour tout $f_0\in e_K\cdot \ES{A}'_0$ 
tel que $f_0\neq 0$, le $\ES{H}$--module $\ES{A}'_0$ est de type fini, et un tel $V'_0$ existe d'après le lemme de Zorn). 
Puisque le $\ES{H}_K$--module $W_0$ est simple, d'apr\`es la fin de la d\'emonstration de \cite[prop.~2.10]{BZ}, le 
$\ES{H}$--module (non dégénéré, simple) $V_0$ vérifie $e_K\cdot V'_0= \{0\}$ et $e_K\cdot V_0\simeq W_0$. Pour $i\in \Bbb{Z}$, on d\'efinit par r\'ecurrence 
un sous--$\ES{H}$--module $V'_i$ de $\ES{A}'_i$ en posant 
$\theta_{\ES{A}'}(V'_{i-1})=V'_i$. 
Soient $V'=\oplus_{i\in \Bbb{Z}}V'_i $ et $V=\ES{A}'/V'$. Par construction, 
$\theta_{\ES{A}'}$ induit par passage au quotient un $\Bbb{C}$--automorphisme 
$\theta_V$ de $V$, et $(V,\theta_V)$ est un $(\ES{H},\theta,\omega)$--module non dégénéré simple tel que 
$e_K\cdot V= W$.

{\it Pla\c{c}ons--nous maintenant dans le second cas: $W=\oplus_{i=0}^{d-1}W_i$ et $\theta_W^d(W_0)=W_0$}. 
Notons $\theta_{W_0}$ la restriction de $\theta_W^d$ \`a 
$W_0$. On a défini en \ref{représentations G-irréductibles} un caractère $\omega_d= \omega \circ \frak{N}_d$ de $G$, où 
(rappel) $\frak{N}_d:G\rightarrow G$ est l'application $x\mapsto x\theta(x)\cdots \theta^{d-1}(x)$. 
Ce caract\`ere est trivial sur $K$, et l'on a
$$
\theta_{W_0}(\omega_df\cdot w)={^{\theta^d\!\!}f}\cdot \theta_{W_0}(w)\quad (f\in \ES{H}_K,\,w\in W_0).
$$
En d'autres termes, $(W_0,\theta_{W_0})$ est un $(\ES{H}_K,\theta^d,\omega_d)$--module non dégénéré simple 
(puisque le $\ES{H}_K$--module $W_0$ est simple). Notons $\ES{A}$ le $\ES{H}$--module $\oplus_{i=0}^{d-1}\ES{A}_i$ 
(pour l'action diagonale de $\ES{H}$), 
et $\theta_\ES{A}$ le $\Bbb{C}$--automorphisme de $\ES{A}$ d\'efini par $\theta_\ES{A}({\bf f})_i= \theta_\ES{H}({\bf f}_{i-1})$ 
($i=0,\ldots, d-1$) pour tout ${\bf f}=({\bf f}_i)_{i=0}^{d-1}$; où l'on a posé 
${\bf f}_{-1}= {\bf f}_{d-1}$. Pour ${\bf f}\in 
\ES{A}$, on a
$$
\theta_\ES{A}(\omega f\cdot {\bf f})_i= \theta_\ES{H}(\omega f\cdot {\bf f}_{i-1})={^{\theta\!}f}\cdot \theta_\ES{A}({\bf f})_i
\quad (i=0,\ldots ,d-1),
$$
i.e. $(\ES{A},\theta_\ES{A})$ est un $(\ES{H},\theta,\omega)$--module non d\'eg\'en\'er\'e. De plus, $\theta_\ES{A}^d$ 
induit par restriction un $\Bbb{C}$--automor\-phisme $\theta_{\ES{A}_0}$ de $\ES{A}_0$ v\'erifiant
$$
\theta_{\ES{A}_0}(\omega_d f * f')= {^{\theta^d\!\!}f}* \theta_{\ES{A}_0}(f')\quad (f,\,f'\in \ES{H}),
$$
i.e. $(\ES{A}_0,\theta_{\ES{A}_0})$ est un $(\ES{H},\theta^d,\omega_d)$--module non d\'eg\'en\'er\'e. 
Posons $\frak{J}=\oplus_{i=0}^{d-1}\frak{J}_i$. C'est un sous--$\ES{H}$--module $\theta_\ES{A}$--invariant 
de $\ES{A}$, et $\frak{J}_0$ est un id\'eal $\theta_{\ES{A}_0}$--invariant de 
$\ES{A}_0\;(=\ES{H}* \ES{H}_K)$. Soit $(\ES{A}'=\ES{A}/\frak{J},\theta_{\ES{A}'})$ le 
$(\ES{H},\theta,\omega)$--module non d\'eg\'en\'er\'e d\'eduit de $(\ES{A},\theta_\ES{A})$ 
par passage au quotient. On a $\ES{A}'= \oplus_{i=0}^{d-1} \ES{A}'_i$ et $\theta_{\ES{A}'}(\ES{A}'_{i-1})=\ES{A}'_i$ ($i=0,\ldots d-1$), où l'on a posé $\ES{A}'_{-1}= \ES{A}'_{d-1}$. De m\^eme, soit 
$(\ES{A}'_0=\ES{A}_0/\frak{J}_0,\theta_{\ES{A}'_0})$ le $(\ES{H},\theta^d,\omega_d)$--module non d\'eg\'en\'er\'e 
d\'eduit de $(\ES{A}_0,\theta_{\ES{A}_0})$ 
par passage au quotient. Comme dans le premier cas, on a $e_K\cdot \ES{A}'\simeq W$ et 
$e_K\cdot \ES{A}'_0\simeq W_0$, et il suffit de montrer qu'il existe un sous--$\ES{H}$--module 
{\it $\theta_{\ES{A}'_0}$--invariant} $V'_0$ de $\ES{A}'_0$ tel que le $(\ES{H},\theta^d,\omega_d)$--module 
$(V_0=\ES{A}'_0/V'_0,\theta_{V_0})$ d\'eduit de $(\ES{A}'_0,\theta_{\ES{A}'_0})$ par passage 
au quotient est $\ES{H}$--simple et v\'erifie $e_K\cdot V_0\simeq W_0$. Soit $V'_0$ un sous--$\ES{H}$--module de $\ES{A}'_0$ tel que le $\ES{H}$--module quotient $V_0= \ES{A}'_0/V'_0$ 
est simple (cf. le premier cas); on a $e_K \cdot V'_0 =\{0\}$ et $e_K\cdot V_0\simeq W_0$. Pour $i\in {\Bbb Z}$, posons $V'_i = \theta_{\ES{A}'_0}^i(V'_0)$. C'est un sous--$\ES{H}$--module de $\ES{A}'_0$ qui vérifie $e_K\cdot V'_i= \theta_{W_0}^i(e_K\cdot V'_0)=\{0\}$. 
Soit $V''_0=\sum_{i\in {\Bbb Z}}V'_i$ le sous--$\ES{H}$--module de $\ES{A}'_0$ engendré par les espaces $V'_i$ pour $i\in {\Bbb Z}$. 
Il est $\theta_{\ES{A}'_0}$--invariant, et vérifie $e_K\cdot V''_0= \{0\}$. On a donc $V''_0 \neq \ES{A}'_0$. Comme $V'_0\subset V''_0 \subset \ES{A}'_0$ et $\ES{A}'_0/V'_0$ est $\ES{H}$--simple, on obtient $V''_0=V'_0$. Par conséquent le $\ES{H}$--module $V'_0$ est $\theta_{\ES{A}'_0}$--invariant, et $\theta_{\ES{A}'_0}$ induit par passage au quotient un ${\Bbb C}$--automorphisme $\theta_{V_0}$ de $V_0$ qui fait de $(V_0,\theta_{V_0})$ un $(\ES{H},\omega, \theta)$--module non dégénéré. Reste à poser $V'= \oplus_{i=0}^{d-1}
\theta_{\ES{A}'}^i(V'_0)\subset \ES{A}'$ et $V= \ES{A}'/V'$. Par construction, 
$\theta_{\ES{A}'}$ induit par passage au quotient un $\Bbb{C}$--automorphisme 
$\theta_V$ de $V$, et $(V,\theta_V)$ est un $(\ES{H},\theta,\omega)$--module non dégénéré simple tel que 
$e_K\cdot V= W$.
\end{proof}

\begin{marema3}
{\rm 
(Sans supposer $({\rm P}_2)$). D'après la preuve de la proposition 2, si $K$ est un sous--groupe ouvert compact 
de $G$ normalisé par un élément de $G^\natural$, alors pour un 
$(\ES{H}^\natural_K,\omega)$--module non dégénéré simple $W$ (tout comme pour un $(\ES{H},\omega)$--module non dégénéré simple, cf. la remarque de \ref{représentations G-irréductibles}), les deux conditions suivantes sont équivalentes:
\begin{itemize}
\item le socle du $\ES{H}_K$--module $W$ n'est pas nul;
\item le $\ES{H}_K$--module $W$ est semisimple.\hfill $\blacksquare$
\end{itemize} 
}
\end{marema3}

\begin{moncoro}(En supposant $({\rm P}_2)$.)
Soit $K\in \ES{K}_1$. 
L'application $V\mapsto V^K$ induit une bijection entre:
\begin{itemize}
\item l'ensemble des classes d'isomorphisme de 
$\omega$--repr\'esentations lisses irr\'eductibles $(\Pi,V)$ de $G^\natural$ telles que $V^K\neq 0$ et la représentation $\Pi^\circ$ de $G$ est semisimple;
\item l'ensemble des 
classes d'isomorphisme de $(\ES{H}_K^\natural,\omega)$--modules non d\'eg\'en\'er\'es 
simples et $\ES{H}_K$--semisimples.
\end{itemize}
Cette bijection se restreint en une bijection entre:
\begin{itemize}
\item l'ensemble des classes d'isomorphisme de 
$\omega$--repr\'esentations lisses $G$--irréductibles $(\Pi,V)$ de $G^\natural$ telles que $V^K\neq 0$;
\item l'ensemble des 
classes d'isomorphisme de $(\ES{H}_K^\natural,\omega)$--modules non d\'eg\'en\'er\'es $\ES{H}_K$--simples.
\end{itemize}
Elle se restreint aussi en une bijection entre:
\begin{itemize}
\item l'ensemble des classes d'isomorphisme de 
$\omega$--repr\'esentations (lisses) admissibles irr\'eductibles $(\Pi,V)$ de $G^\natural$ telles que $V^K\neq 0$;
\item l'ensemble des 
classes d'isomorphisme de $(\ES{H}_K^\natural,\omega)$--modules non d\'eg\'en\'er\'es 
simples et de dimension finie (sur ${\Bbb C}$).
\end{itemize}
\end{moncoro}
 
\begin{marema4}{\rm La d\'emonstration du lemme de Schur donn\'ee dans \cite[2.11]{BZ} fonctionne 
aussi pour les $\omega$--repr\'esentations lisses irr\'eductibles de $G^\natural$ (sans supposer $({\rm P}_2)$ ni $({\rm P}_1)$): 
soit $(\Pi,V)$ une $\omega$--repr\'esentation lisse irr\'eductible de 
$G^\natural$. Si le groupe $G$ est d\'enombrable \`a l'infini (c'est--\`a--dire r\'eunion d\'enombrable de sous--ensembles 
compacts) alors on a
$${\rm End}_{G^\natural}(\Pi)=\Bbb{C}{\rm id}_V.
$$
Si la propriété $({\rm P}_1)$ est vérifiée, alors d'après la proposition 1, la conclusion reste 
vraie pour toutes les $\omega$--représentations admissibles irréductibles de $G^\natural$ (sans supposer $G$ dénombrable à l'infini). 
\hfill 
$\blacksquare$}
\end{marema4}

\subsection{Ind\'ependance lin\'eaire des caract\`eres tordus}\label{ind\'ependance lin\'eaire des caract\`eres tordus}
Dans ce num\'ero, on suppose que la 
propri\'et\'e la propri\'et\'e $({\rm P}_2)$ est v\'erifi\'ee. Comme en \ref{reprsentations irrductibles et H-modules simples}, 
on fixe un \'el\'ement $\delta_1\in G^\natural$ et 
un syst\`eme fondamental de voisinages de $1$ dans $G$, disons $\ES{K}_1$,  form\'e de sous--groupes 
ouverts compacts $K$ de $G$ tels que $\omega\vert_K=1$ et ${\rm Int}_{G^\natural}(\delta_1)(K)=K$. 

Notons que si $\Pi_1$, $\Pi_2$ sont deux $\omega$--repr\'esentations lisses de 
$G^\natural $ telles que les repr\'esentations $\Pi_1^\circ$, $\Pi_2^\circ$ de $G$ sont irr\'eductibles, alors 
on a ${\rm Hom}_{G^\natural}(\Pi_1,\Pi_2)={\rm Hom}_G(\Pi_1^\circ,\Pi_2^\circ)$. La proposition de 
\ref{rappels sur les reprsentations irrductibles de G} et son corollaire se g\'en\'eralisent de la mani\`ere suivante:

\begin{mapropo}(En supposant $({\rm P}_2)$.)
Soit $\Pi_1,\,\Pi_2,\,\ldots,\Pi_n$ des 
$\omega$--repr\'esentations admissibles $G$--irrductibles de $G^\natural$ 
telles que les repr\'esentations $\Pi_1^\circ,\,\Pi_2^\circ,\,\ldots,\Pi_n^\circ$ de $G$ 
sont deux--\`a--deux non isomorphes. Alors les distributions 
$\Theta_{\Pi_1},\, \Theta_{\Pi_2},\,\ldots,\, \Theta_{\Pi_n}$ sur $G^\natural$ sont lin\'eairement ind\'ependantes.
\end{mapropo}

\begin{proof}
On reprend celle de \cite[prop.~2.19]{BZ}. Pour $i=1,\,\ldots ,n$, notons $V_i$ l'espace de $\Pi_i$. 
Choisissons un groupe $K\in \ES{K}_1$ tel que pour $i=1,\ldots ,n$, on a $V_i^K\neq 0$. 
Les $\ES{H}_K$--modules 
$V_i^K$ sont simples, de dimension finie, et deux--\`a--deux non isomorphes. Pour $i=1,\ldots, n$, 
l'automorphisme $\Pi(\delta_1)$ de $V_i^K$ induit par restriction un automorphisme $A_i$ de $V_i^K$, qui 
co\"{\i}ncide avec la restriction de $\Pi_i(\phi^K_{\delta_1})$ \`a $V_i^K$; on a donc
$$
\Theta_{\Pi_i}(\phi)={\rm tr}(\Pi_i(\phi)\circ A_i)\quad (\phi\in \ES{H}_K).
$$
On conclut gr\^ace au th\'eor\`eme de Frobenius--Schur \cite[ch.~VIII, \S13, prop.~2]{Bou}: pour $i=1,\ldots ,n$, le choix d'une base de $V_i$ sur 
$\Bbb{C}$ permet d'identifier ${\rm End}_\Bbb{C}(V_i^K)$ \`a ${\rm M}_{d_i}(\Bbb{C})$. Pour 
$1\leq k,l\leq d_i$, notons $u^i_{k,l}:\ES{H}_K\rightarrow \Bbb{C}$ l'application qui \`a 
$\phi\in \ES{H}_K$ associe le coefficient en la place $(k,l)$ de l'endomorphisme $w\mapsto \phi\cdot w$ 
de $V_i$ (vu comme un \'el\'ement de ${\rm M}_{d_i}(\Bbb{C})$). Le th\'eor\`eme de Frobenius--Schur dit 
que les fonctions 
$u^i_{k,l}$ ($i=1,\ldots ,n$, $1\leq k,l\leq d_i$) sont lin\'eairement ind\'ependantes sur $\Bbb{C}$.
\end{proof}

\begin{moncoro}
(En supposant $({\rm P}_2)$.)
Soit $\Pi_1$, $\Pi_2$ deux 
$\omega$--repr\'esentations admissibles $G$--irrductibles de $G^\natural $. Alors 
$\Pi_1$ et $\Pi_2$ sont isomorphes si et seulement si 
$\Theta_{\Pi_1}=\Theta_{\Pi_2}$.
\end{moncoro}

\subsection{La condition $({\rm P}_2)$ pour $G^\natural= {\bf G}^\natural(F)$}\label{la condition P2}
On reprend maintenant les hypoth\`eses du 
chapitre 5: $G={\bf G}(F)$ et $G^\natural={\bf G}^\natural(F)$ pour un groupe r\'eductif connexe ${\bf G}$ d\'efini sur $F$ et un ${\bf G}$--espace 
tordu ${\bf G}^\natural$ d\'efini sur $F$ et poss\'edant un point $F$-rationnel; où $F$ est un corps commutatif localement compact non archimédien.

\begin{monlem}
Le $G$--espace tordu $G^\natural$ 
v\'erifie ($P_2$): il existe un \'el\'ement $\delta_1\in G^\natural$ et une 
base de voisinages de $1$ dans $G$ form\'ee de sous--groupes ouverts compacts de $G$ 
normalis\'es par $\delta_1$.
\end{monlem}

\begin{proof}
Soit $F^{\rm nr}$ une extension non ramifi\'ee maximale de $F$. Posons $\frak{o}=\frak{o}_F$ et notons 
$\frak{o}^{\rm nr}$ l'anneau des entiers de $F^{\rm nr}$. Rappelons que $\varpi$ d\'esigne une uniformisante de $F$. 
Fixons un \'element $\delta\in G^\natural$. Notons $\Sigma$ le groupe de Galois de l'extension $F^{\rm nr}/F$, 
et fixons une chambre $\Sigma$--stable $\ES{C}$ de l'immeuble \'etendu $\frak{I}^{\rm nr}$ de ${\bf G}(F^{\rm nr})$. 
Notons $I^{\rm nr}$ le stabilisateur de $\ES{C}$ dans ${\bf G}(F^{\rm nr})$. Par transport 
de structure, le $F$--automorphisme ${\rm Int}_{{\bf G}^\natural}(\delta)$ de ${\bf G}$ induit un automorphisme de 
$\frak{I}^{\rm nr}$, qui commute \`a l'action de $\Sigma$; il envoie donc $\ES{C}$ sur une autre chambre $\Sigma$--stable 
de $\frak{I}^{\rm nr}$, disons $\ES{C}'$. On sait que $G$ op\`ere transitivement sur l'ensemble des 
chambres $\Sigma$--stables de $\frak{I}^{\rm nr}$. Soit donc un \'el\'ement $g\in G$ tel que $g\cdot \ES{C}'=\ES{C}$. 
Posons $\delta_1=g\cdot\delta$ et $\theta={\rm Int}_{{\bf G}^\natural}(\delta_1)$. Puisque $\theta(\ES{C})=\ES{C}$, 
on a $\theta(I^{\rm nr})=I^{\rm nr}$. On sait \cite[4.6.30, 5.1.30]{BT2} que $I=(I^{\rm nr})^\Sigma$ est le 
groupe des points $\frak{o}$--rationnels d'un $\frak{o}$-sch\'ema en groupes affine lisse (pas n\'ecessairement connexe) 
$\frak{G}$ de fibre 
g\'en\'erique $\frak{G}\times_\frak{o}F={\bf G}$, caract\'eris\'e \`a isomorphisme unique pr\`es par l'\'egalit\'e 
$\frak{G}(\frak{o}^{\rm nr})=I^{\rm nr}$. Puisque $\theta(I^{\rm nr})=I^{\rm nr}$, $\theta$ se prolonge de mani\`ere unique en un morphisme de $\frak{o}$-sch\'emas 
$u:\frak{G}\rightarrow \frak{G}$, qui est un isomorphisme de $\frak{o}$-sch\'ema en groupes.

Pour chaque entier $n\geq 1$, notons $\frak{G}(\frak{o}^{\rm nr})^n$ le {\it $n$--i\`eme sous--groupe de 
congruence de $\frak{G}(\frak{o}^{\rm nr})$}, d\'efini par
$$
\frak{G}(\frak{o}^{\rm nr})^n=\ker\{\pi_n:\frak{G}(\frak{o}^{\rm nr})\rightarrow \frak{G}(\frak{o}^{\rm nr}/\varpi^n\frak{o}^{\rm nr})\},
$$
o\`u la $\pi_n$ d\'esigne l'application canonique (r\'eduction modulo $\varpi^n$). D'apr\`es \cite[2.8]{Y}, il existe un 
$\frak{o}$--sch\'ema en groupes affine lisse $\frak{G}^n$ de fibre g\'en\'erique 
$\frak{G}^n\times_\frak{o}F={\bf G}$ tel que $\frak{G}^n(\frak{o}^{\rm nr})=\frak{G}(\frak{o}^{\rm nr})^n$. 
De plus, l'\'egalit\'e $\frak{G}^n(\frak{o}^{\rm nr})=\frak{G}(\frak{o}^{\rm nr})^n$ caract\'erise 
$\frak{G}^n$ \`a isomorphisme unique pr\`es, et pour $m\in \Bbb{Z}_{\geq 1}$, on a $(\frak{G}^n)^m=\frak{G}^{n+m}$. 
Posons $u_n=u\times_\frak{o}\frak{o}_n$, $\frak{o}_n=\frak{o}/\varpi^n\frak{o}$. Pour 
$x\in \frak{G}(\frak{o}^{\rm nr})$, on a $\pi_n(u(x))=u_n(\pi_n(x))$. Par suite on a 
$\theta(\frak{G}(\frak{o}^{\rm nr})^n)=\frak{G}(\frak{o}^{\rm nr})^n$ et 
$\theta(\frak{G}^n(\frak{o}))=\frak{G}^n(\frak{o})$. La famille $\ES{K}_1=\{\frak{G}^n(\frak{o}):n\in \Bbb{Z}_{\geq 1}\}$ est 
un syst\`eme fondamental de voisinages de $1$ dans $G$ v\'erifiant les propri\'et\'es voulues.
\end{proof}

\begin{marema}{\rm 
Pour $\Pi_1,\ldots ,\Pi_2$ comme dans l'\'enonc\'e de la proposition de \ref{ind\'ependance lin\'eaire des caract\`eres tordus}, on ne sait pas 
a priori si les fonctions 
caract\`eres $\Theta_{\Pi_i}: G^\natural_{\rm qr}\rightarrow \Bbb{C}$ sont lin\'eairement ind\'ependantes sur 
$\Bbb{C}$ (sauf bien s\^ur si elles sont localement int\'egrables sur $G^\natural$).\hfill $\blacksquare$}
\end{marema}

\section{Repr\'esentations $l$-modulaires}

Dans cette annexe B, on d\'ecrit bri\`evement comment les r\'esultats des ch.~2 et 5 s'\'etendent au cas des repr\'esentations \`a valeurs dans le groupe des 
automorphismes d'un espace vectoriel sur un corps de caract\'eristique $l$ 
diff\'erente de la caract\'eristique r\'esiduelle de $F$.

\subsection{G\'en\'eralit\'es \cite[ch.~1]{V}}  Soit $R$ un anneau commutatif, poss\'edant une unit\'e $1_R$. 
On note $R^\times$ le groupe des \'el\'ements inversibles de $R$, et $d\Bbb{Z}$ l'id\'eal de 
$\Bbb{Z}$ noyau du morphisme canonique $\Bbb{Z}\rightarrow  R,\,a\mapsto a1_R$. Ce morphisme 
se prolonge au sous--anneau $\ES{A}$ de $\Bbb{Q}$ engendr\'e par les inverses des entiers 
premiers \`a $d$.

Si $X$ est un td--espace, on note $C^\infty_{\rm c}(X,R)$ l'espace des fonctions sur 
$X$ \`a valeurs dans $R$, qui sont localement constantes et \`a support compact, et 
l'on pose
$$C^\infty_{\rm c}(X,R)^*={\rm Hom}_R(C^\infty_{\rm c}(X,R),R).
$$ Les 
\'el\'ements de $C^\infty_{\rm c}(X,R)^*$ sont appel\'es {\it $R$--distributions} sur $X$. 

Soit $G$ un groupe localement profini. On appelle {\it $R$--caract\`ere} de $G$ 
un morphisme de groupes $G\rightarrow R^\times$ dont le noyau contient un 
sous--groupe ouvert.

Une {\it $R$--mesure de Haar} \`a gauche sur $G$ est 
par d\'efinition une distribution non nulle sur $G$ invariante pour l'action de $G$ op\'erant 
sur lui--m\^eme par translations \`a gauche. On sait qu'une telle mesure existe si et seulement s'il 
existe un sous--groupe ouvert compact $K$ de $G$ dont le pro--ordre est premier \`a $d$, 
auquel cas cette mesure est unique \`a multiplication pr\`es par un \'el\'ement de $R^\times$. 
On suppose d\'esormais qu'il existe une $R$--mesure de Haar \`a gauche sur $G$, et l'on en fixe une 
$\mu_G$. On note $\Delta_{G,R}:G\rightarrow R^\times$ le {\it $R$--module de $G$}, c'est--\`a--dire 
le $R$--caract\`ere d\'efini comme en \ref{modules} par
$$
\int_Gf(gx^{-1})d\mu_G(g)=\Delta_{G,R}(x)\int_Gf(g)d\mu_G(g)\quad (f\in C^\infty_{\rm c}(G,R),\,x\in G).
$$
Alors $\Delta_{G,R}^{-1}\mu_G$ est une $R$--mesure de Haar \`a droite sur $G$, et on les obtient toutes de 
cette mani\`ere.

\begin{marema}{\rm 
Le module $\Delta_G(x)$ d'un \'el\'ement $x\in G$ est un indice g\'en\'eralis\'e: pour tout sous--groupe ouvert compact $K$ de $G$, on a 
$$
\Delta_G(x)= {[K:K\cap xKx^{-1}]\over [xKx^{-1}:K\cap xKx^{-1}]}={[K:K\cap xKx^{-1}]\over [K:x^{-1}Kx\cap K]}.
$$
En choisissant $K$ de pro--ordre premier \`a $d$, on voit que $\Delta_G(x)$ appartient \`a l'anneau 
$\ES{A}$. Son image dans $R$ appartient \`a $R^\times$: 
c'est le $R$--module $\Delta_{G,R}(x)$.\hfill $\blacksquare$}
\end{marema}

Soit $H$ un sous--groupe ferm\'e de $G$. Alors il existe une mesure de Haar 
\`a gauche sur $H$, et une mesure de Haar \`a droite sur l'espace quotient $H\backslash G$, 
c'est--\`a--dire une forme lin\'eaire non nulle sur $\boldsymbol{S}(H\backslash G,R)$, invariante pour l'action de $G$ sur $H\backslash G$ par translations \`a droite; o\`u $\boldsymbol{S}(H\backslash G,R)$ d\'esigne, comme en \ref{modules}, l'espace des fonctions 
$f:H\backslash G\rightarrow R$ qui sont uniform\'ement localement constantes \`a droite, \`a support compact 
modulo $H$, et v\'erifient
$$
f(hg)=\Delta_{G,R}(h)\Delta_{H,R}(h)^{-1}f(h)\quad (h\in H,\,g\in G).
$$

Soit $\theta$ un automorphisme de $G$. On d\'efinit le $R$--module $\Delta_{G,R}(\theta)\in R^\times$ comme en \ref{modules}. 
D'apr\`es le lemme de \ref{modules}, on a $\Delta_{G,R}=\Delta_{G,R}\circ \theta$, et 
s'il existe une partie ouverte compacte $\theta$--stable $\Omega$ de $G$ telle que ${\rm vol}(\Omega,\mu_G)\neq 0$, alors $\Delta_{G,R}(\theta)=1$. Si $H$ est un sous--groupe ferm\'e 
$\theta$--stable de $G$, on d\'efinit le $R$--module $\Delta_{H\backslash G,R}(\theta)$ comme en \ref{modules}. La 
relation $(*)$ de \ref{modules} est vraie pour les $R$--modules.

\subsection{$R$--repr\'esentations lisses}
On appelle {\it $R$--repr\'esentation lisse} de $G$ la 
donn\'ee d'un $R$--module $V$ et d'un morphisme de groupes $\pi:G\rightarrow {\rm Aut}_R(V)$ 
tel que le stabilisateur de $v$ dans $G$ est ouvert, pour tout $v\in V$. Comme pour les 
repr\'esentations complexes, on a les notions de $R$--repr\'esentation (lisse) admissible, irr\'eductible, semisimple, 
de type fini, de longueur finie. Si $R$ est un corps, alors pour toute $R$--repr\'esentation admissible 
$\pi$, on d\'efinit comme en \ref{caractres} le caract\`ere--distribution $\Theta_\pi ={\rm tr}(\pi):C^\infty_{\rm c}(G,R)\rightarrow R$ 
(il d\'epend du choix de la mesure de Haar $\mu$ sur $G$). Si $R$ est un 
corps alg\'ebriquement clos, la proposition de \ref{rappels sur les reprsentations irrductibles de G} 
reste vraie: les caract\`eres--distributions des repr\'esentations admissibles 
irr\'eductibles de $G$ sont lin\'eairement ind\'ependants sur $R$. 

Soit $G^\natural$ un $G$--espace tordu, et soit $\omega$ un $R$--caract\`ere de $G$. On d\'efinit 
le $R$--module de $G^\natural$ comme en \ref{module d'un espace tordu} : pour $\gamma\in G^\natural$, on pose
$$
\Delta_{G^\natural\!,R}(\gamma)=\Delta_{G,R}({\rm Int}_{G^\natural}(\gamma)^{-1}).
$$ 
Le lemme de \ref{module d'un espace tordu} est vrai pour les $R$--modules. On d\'efinit la notion de 
$(\omega,R)$--repr\'esentation 
lisse de $G^\natural$ 
comme en \ref{caractres tordus (bis)}, en rempla\c{c}ant $\Bbb{C}$ par $R$. La 
cat\'egorie des $(\omega,R)$--repr\'esentations lisses de $G^\natural$ est not\'ee 
$\frak{R}(G^\natural,\omega,R)$.

Soit $H$ un sous--groupe ferm\'e de $G$, et $H^\natural$ un $H$--espace tordu qui soit 
un sous--espace topologique tordu de $G^\natural$. On d\'efinit comme en \ref{induction compacte} un 
foncteur induction compacte (lisse)
$$
{^\omega{\rm ind}}^{G^\natural}_{H^\natural}:\frak{R}(H^\natural,\omega,R)
\rightarrow \frak{R}(G^\natural,\omega,R).
$$

{\bf On suppose d\'esormais que l'anneau $R$ est un corps}, de caract\'eristique $l$. On a donc $d=l$. 
Pour toute partie ouverte compacte $\Omega$ de $G$ telle que l'ensemble $H\cap \Omega$ est non vide et 
de pro--ordre premier \`a $l$, 
on d\'efinit l'endomorphisme $\phi\mapsto \phi\vert_\Omega$ de 
$C^\infty_{\rm c}(G^\natural\!,R)$ comme en \ref{caractres des induites compactes}. La 
proposition de \ref{caractres des induites compactes} reste vraie pour les $(\omega,R)$--repr\'esentations de $G^\natural$ 
pourvu que $K$ soit de pro--ordre premier à $l$. En effet, fix\'e un 
sous--groupe ouvert compact $K$ de $G$ de pro--ordre premier \`a $l$, tel que $HK=KH$, et un système de représentants $\{x_1=1,x_2,\ldots ,x_n\}$ de $HK\backslash G$ dans $G$, la partie ouverte compacte $\Omega'= \bigcup_{i=1}^nKx_i$ de $G$ vérifie 
$H\cap \Omega' = H\cap K$. Puisque le pro--ordre de $H\cap \Omega'$ divise celui de $K$, 
il est premier \`a $l$.

\subsection{Le principe de submersion d'Harish--Chandra}\label{le principe de submersion d'Harish--Chandra (bis)}
Reprenons les hypoth\`eses et les 
notations du chapitre 5. Soit $p$ la caract\'eristique r\'esiduelle de $F$. 
Puisqu'on a suppos\'e qu'il existe un sous--groupe ouvert compact de 
$G$ dont le pro--ordre est premier \`a $l$, on a $l\neq p$. Puisque $l$ peut diviser 
le pro--ordre du sous--groupe 
compact maximal $K_\circ$ de $G$, on ne peut normaliser les mesures de 
Haar comme en \ref{mesures normalises}.  
Le th\'eor\`eme de \ref{les oprateurs T} et son corollaire restent vrais, pourvu que le sous--groupe 
ouvert compact d\'efinissant les op\'erateurs ${\rm T}_\gamma$ soit de pro--ordre premier 
\`a $l$:

\begin{montheo}Soit $\Pi$ une $(\omega,R)$--repr\'esentation admissible de $G^\natural$ 
telle que la $R$--reprsen\-tation $\Pi^\circ$ de $G$ est de type fini, et soit $K$ 
un sous--groupe ouvert compact de $G$ dont le pro--ordre est premier \`a $l$. 
Le caract\`ere-distribution $\Theta_{\Pi}={\rm tr}(\Pi)$ de $\Pi$ est repr\'esent\'e sur $G^\natural_{\rm qr}$ par la fonction 
localement constante $\gamma\mapsto \Theta_\Pi(\gamma)={\rm tr}({\rm T}_\gamma)$, o\`u l'on 
a pos\'e
$${\rm T}_\gamma = {\rm vol}(K,\mu_G)^{-1}\int_K\omega(k)^{-1}\Pi(k^{-1}\cdot\gamma\cdot k)d\mu_G(k).$$
\end{montheo}

\subsection{Induction parabolique et restriction de Jacquet}
Pour $P^\natural\in \ES{P}^\natural_\circ$, on 
d\'efinit le foncteur induction parabolique {\it non normalis\'e}
$$
{^\omega\underline{\iota}_{P^\natural}^{G^\natural}}:\frak{R}(M_P^\natural,\omega,R)
\rightarrow \frak{R}(G^\natural,\omega ,R)
$$
comme en \ref{induction parabolique et caractres}, 
en supprimant le facteur $\delta_{P^\natural}^{1/2}$. On d\'efinit aussi le foncteur de Jacquet {\it non normalis\'e}
$$
{^\omega\underline{r}_{G^\natural}^{P^\natural}}:\frak{R}(G^\natural,\omega,R)
\rightarrow \frak{R}(M_P^\natural,\omega,R)
$$
comme en \ref{restriction de Jacquet et caractres}, en supprimant le facteur $\delta_{P^\natural}^{-1/2}$.

\begin{marema1}
{\rm Pour $\delta\in P^\natural$, le module 
$\Delta_{P^\natural}(\delta)=\delta_{P^\natural}(\delta)^{-1}$ appartient \`a $\Bbb{Z}[1/q]$, o\`u 
$q$ est le cardinal du corps r\'esiduel de $F$. Si l'image de $q$ dans $R$ a une racine carr\'ee dans 
$R$ (par exemple si $R$ est alg\'ebriquement clos), on en choisit une; cela permet de d\'efinir le $R$--caract\`ere 
$\delta_{P^\natural}^{1/2}$ de $P^\natural$, et les foncteurs induction parabolique et restriction de Jacquet 
{\it normalis\'es} pour les $(\omega,R)$--repr\'esentations lisses.\hfill $\blacksquare$
}\end{marema1}

Choisissons un sous--groupe ouvert $J_\circ$ de $K_\circ$, de pro--ordre premier \`a $l$, tel que pour tout $P\in \ES{P}_\circ$, on a la 
d\'ecomposition triangulaire
$$
J_\circ = (J_\circ \cap U_{P^-})(J_\circ \cap M_P)(J_\circ \cap U_P).
$$
Un tel $J_\circ$ existe, d'apr\`es la remarque 4 de \ref{la paire parabolique associe  gamma}. On suppose d\'esormais que 
les $R$--mesures de Haar \`a gauche $\mu_G$, $\mu_{M_P}$, 
$\mu_{U_P}$, $\mu_P$ sur $G$, $M_P$, $U_P$, $P$ sont celles 
normalis\'ees par $J_\circ$, c'est--\`a--dire 
telles que ${\rm vol}(J_\circ,\mu_G)=1$, ${\rm vol}(M_P\cap J_\circ,\mu_{M_P})=1$, etc.\h{0.2}.

Soit $P^\natural\in \ES{P}_\circ^\natural$. Choisissons un syst\`eme de repr\'esentants $\{x_1,\ldots ,x_n\}$ dans $K_\circ$ de l'espace quotient $PJ_\circ\backslash G$ tel que $1\in \{x_1,\ldots ,x_n\}$, et posons $\Omega=\coprod_{i=1}^nJ_\circ x_i$. Alors $\Omega \cap P= J_\circ \cap P$ est de pro--ordre premier \`a $l$. On note
$$
C^\infty_{\rm c}(G^\natural\!,R)\rightarrow C^\infty_{\rm c}(M^\natural_P,R),\,
\phi\mapsto {^\omega{\phi}_{P^\natural,J_\circ}}
$$
l'application lin\'eaire d\'efinie par
$$
{^\omega\phi_{P^\natural,J_\circ}}(\delta)= \int\!\!\!\int_{U_P\times J_\circ}\omega(k)\phi(k^{-1}\cdot \delta\cdot
uk)d\mu_{U_P}(u)d\mu_G(k)\quad (\delta\in M_P^\natural).
$$

Soit $\Sigma$ une $(\omega,R)$--repr\'esentation admissible de $M_P^\natural$, et 
$\Pi={^\omega\underline{\iota}^{G^\natural}_{P^\natural}}(\Sigma)$. Alors $\Pi$ est une $(\omega,R)$--reprsentation admissible de $G^\natural$, et 
pour toute fonction $\phi\in C^\infty_{\rm c}(G^\natural\!,R)$, on a la formule de descente (thorme de \ref{induction parabolique et caractres})
$$
\Theta_\Pi(\phi)=\Theta_\Sigma({^\omega\phi_{P^\natural,J_\circ}}).
$$

\begin{marema2}
{\rm Nous n'essaierons pas de traduire ici l'\'egalit\'e ci-dessus 
en termes de fonctions caract\`eres (corollaire 3 de \ref{formule d'intgration de Weyl}). Notons d'ailleurs que si $\Sigma^\circ$ est de type fini, on ne sait 
pas si $\Pi^\circ$ l'est aussi, m\^eme si $R$ est alg\'ebriquement clos. On sait en revanche, si $R$ 
est alg\'ebriquement clos, qu'une repr\'esentation admissible de type fini est de 
longueur finie \cite[ch.~II, 5.10]{V}, et que les foncteurs induction parabolique 
$\underline{\iota}_P^G:\frak{R}(M_P,R)\rightarrow \frak{R}(G,R)$ et 
restriction de Jacquet $r_G^P:\frak{R}(G,R)\rightarrow \frak{R}(M_P,R)$ pr\'eservent la 
propri\'et\'e d'\^etre de longueur finie \cite[ch.~II, 5.14]{V}.\hfill $\blacksquare$}
\end{marema2}

Soit $\Pi$ une $(\omega,R)$-repr\'esentation admissible de $G^\natural$ telle que 
$\Pi^\circ$ est de type fini, et soit $\Sigma={^\omega\underline{r}^{P^\natural}_{G^\natural}}(\Pi)$. 
D'apr\`es \cite[3.2, 3.3]{V}, 
le premier lemme de Jacquet est vrai (cf. la d\'emonstration du lemme 2 de \ref{restriction de Jacquet et caractres}), et $\Sigma$ est admissible. 
De plus $\Sigma^\circ$ est encore de type fini, et 
pour $\gamma\in M^\natural_P\cap G^\natural_{\rm qr}$ tel que $P_{[\gamma]}=P$ et 
$A(\gamma)=A_P$, on a l'\'egalit\'e (thorme de \ref{restriction de Jacquet et caractres})
$$
\Theta_\Pi(\gamma)=\Theta_\Sigma(\gamma).
$$

\subsection{Commentaires}
Pour $\theta={\rm id}$ et $\omega=1$, Meyer et Solleveld \cite{MS} ont 
r\'ecemment obtenu le th\'eor\`eme de \ref{le principe de submersion d'Harish--Chandra (bis)} 
par une m\'ethode diff\'erente de celle d'Harish--Chandra, 
utilisant des syst\`emes de coefficients sur l'immeuble de Bruhat--Tits. Leur r\'esultat est d'ailleurs plus 
fort puisqu'ils contr\^{o}lent le voisinage d'un \'el\'ement semisimple r\'egulier de $G$ sur lequel la 
fonction caract\`ere d'une $R$--repr\'esentation admissible de longueur finie de $G$ est constant. 

Il est probablement possible d'\'etendre une partie du ch.~6 aux $(\omega,R)$--repr\'esentations 
de $G^\natural$, du moins pour celles dont les coefficients sont \`a support compact modulo le centre, et 
en supposant $R$ alg\'ebriquement clos (lemme de Schur, degr\'e formel, etc.).

\section{Action d'un groupe algbrique et points rationnels}

Dans cette annexe C, on s'int\'eresse \`a l'action d'un groupe localement profini $G$ 
sur un espace totalement discontinu $X$, et plus particuli\`erement au cas o\`u cette action provient par 
passage aux points $F$--rationnels d'une action alg\'ebrique d\'efinie sur un corps commutatif localement 
compact non archim\'edien $F$. On reprend ici les rsultats de Bernstein--Zelevinski 
\cite[\S6, Appendix]{BZ}, en particulier le thorme de constructibilit. Gr\^ace aux techniques de loc.~cit., 
on obtient un rsultat nouveau (proposition de \ref{une consquence du lemme cl}). 
On rappelle aussi certains rsultats plus rcents de 
Moret--Bailly, Gabber et Gille \cite{MB2, GGMB}. 

\subsection{Rappels topologiques}\label{rappels topo}
Soit $X$ un espace topologique, et $Y$ un sous--ensemble de $X$. 
Notons $\overline{Y}$ la fermeture de $Y$ (dans $X$). Rappelons que $Y$ est dit {\it localement ferm\'e (dans $X$)} si 
les conditions \'equivalentes suivantes sont v\'erifi\'ees:
\begin{itemize}
\item $Y$ est l'intersection d'un ouvert et d'un ferm\'e de $X$;
\item $Y$ est ouvert dans $\overline{Y}$;
\item $Y$ est ferm\'e au voisinage de chacun de ses points (dans $X$).
\end{itemize}
Notons $\ES{U}(Y)=\ES{U}_X(Y)$ l'ensemble des $y\in Y$ tels que $Y$ est ferm\'e au voisinage de $y$ (dans $X$). 
C'est un sous--ensemble localement ferm\'e de $X$, ouvert dans $Y$, et $Y$ est localement ferm\'e si et seulement 
si $\ES{U}(Y)=Y$. 

On d\'efinit par r\'ecurrence une suite
$$Y=Y_0\supset Y_1\supset Y_2\supset \cdots \supset Y_k \supset \cdots $$
de sous--ensembles ferm\'es de $Y$: on pose $Y_1= Y\smallsetminus \ES{U}(Y)$ et $Y_k=(Y_{k-1})_1$. 
On dit que $Y$ est {\it constructible (dans X)} si $Y$ est union d'un nombre fini de sous--ensembles 
localement ferm\'es. D'apr\`es 
\cite[6.7]{BZ}, si $Y$ est constructible, alors il existe un plus petit entier $k\geq 1$ tel que 
$Y_k=\emptyset$, et $Y$ 
s'\'ecrit
$$
Y=\ES{U}(Y)\cup\ES{U}(Y_1)\cup\cdots \cup \ES{U}(Y_{k-1}).
$$
De plus, $\ES{U}(Y)$ est dense dans $Y$. Inversement, s'il existe un entier $k\geq 1$ tel que $Y_k=\emptyset$, alors $Y$ est constructible.

\subsection{Actions r\'eguli\`eres, localement r\'eguli\`eres, et constructibles}\label{actions rgulires}
Soit $X$ un td--espace muni d'une action continue 
d'un groupe topologique localement profini $G$. On note 
$$\alpha:G\times X\rightarrow X,\,(g,x)\mapsto 
\alpha(g,x)=g\cdot x$$ cette action, $G\backslash X$ l'ensemble des $G$--orbites dans $X$, et $\rho:X\rightarrow G\backslash X$ 
la projection canonique. On munit $G\backslash X$ de la topologie quotient: un sous--ensemble $\widetilde{Y}$ de $G\backslash X$ est ouvert si 
et seulement si $\rho^{-1}(\widetilde{Y})$ est ouvert dans $X$. Cela fait de $\rho$ une application 
continue ouverte. Notons que $G\backslash X$ n'est pas n\'ecessairement {\it s\'epar\'e} (au sens de Hausdorff), mais puisque 
$\rho$ envoie tout compact de $X$ sur un quasi--compact de $G\backslash X$, chaque point de $G\backslash X$ poss\`ede une base de 
voisinages ouverts quasi--compacts. 

D'apr\`es \cite[lemma 6.4]{BZ}, les 
conditions suivantes sont \'equivalentes:
\begin{itemize}
\item l'action $\alpha$ est {\it r\'eguli\`ere} au sens o\`u son graphe
$$
\ES{G}_X^\alpha=\{(x,g\cdot x):x\in X,\, g\in G\}\subset X\times X
$$
est ferm\'e dans $X\times X$;
\item la diagonale
$$
\Delta_{G\backslash X}=\{(\tilde{x},\tilde{x}):\tilde{x}\in G\backslash X\}\subset G\backslash X\times G\backslash X
$$
est ferm\'ee dans $G\backslash X\times G\backslash X$;
\item l'espace topologique $G\backslash X$ est s\'epar\'e.
\end{itemize}
%

\begin{marema1}
{\rm Si l'action $\alpha$ est r\'eguli\`ere, toutes les $G$--orbites dans 
$X$ sont ferm\'ees, et $G\backslash X$ est un td--espace. Notons que pour que l'action $\alpha$ soit r\'eguli\`ere, 
il faut et il suffit que chaque point de $x$ poss\`ede un voisinage ouvert ferm\'e et $G$--invariant sur lequel 
elle est r\'eguli\`ere.\hfill $\blacksquare$}
\end{marema1}

\begin{exemple}
{\rm Si $G$ est un sous--groupe ferm\'e d'un 
groupe topologique localement profini $H$, l'ensemble des classes $G\backslash H$ 
muni de la topologie quotient est un td--espace, et la projection canonique 
$H\rightarrow G\backslash H$ est ouverte \cite[cor.~6.5]{BZ}. En effet, on pose $X=H$ et pour $\alpha$ on prend 
l'action de $G$ sur $H$ donn\'ee par la multiplication \`a gauche. Le graphe
$$\ES{G}_H^\alpha=\{(h,h')\in H\times H: h'h^{-1}\in G\}$$ 
est ferm\'e dans $H\times H$, d'o\`u le r\'esultat. \hfill $\blacksquare$}
\end{exemple}

Soit un lment $\tilde{x}\in G\backslash X$ tel que la diagonale $\Delta_{G\backslash X}$ 
est ferm\'ee au voisinage de $(\tilde{x},\tilde{x})$ dans $G\backslash X\times G\backslash X$. Choisissons un voisinage ouvert 
$\widetilde{U}$ de $\tilde{x}$ dans $G\backslash X$ tel que $\Delta_{G\backslash X}\cap (\widetilde{U}\times \widetilde{U})$ est 
ferm\'e dans $\widetilde{U}\times \widetilde{U}$. Posons $U=\rho^{-1}(\widetilde{U})$. C'est un ouvert 
(non vide) $G$--invariant de $X$ sur lequel l'action $\alpha$ 
est r\'eguli\`ere: le graphe
$$\ES{G}_U^\alpha =\{(x,g\cdot x):x\in U,\,g\in G\}=\ES{G}_X^\alpha\cap (U\times U)
$$ est ferm\'e dans $U\times U$. R\'eciproquement, si la restriction de $\alpha$ \`a un 
ouvert $G$--invariant $U$ de $X$ est r\'eguli\`ere, alors la diagonale $\Delta_{G\backslash U}= \Delta_{G\backslash X}\cap (\rho(U)\times \rho(U))$ 
est ferm\'ee dans $\rho(U)\times \rho(U)$. On en d\'eduit que les conditions suivantes sont \'equivalentes:
\begin{itemize}
\item l'action $\alpha$ est {\it localement r\'eguli\`ere} au sens o\`u son graphe $\ES{G}_X^\alpha$ 
est localement ferm\'e dans $X\times X$;
\item la diagonale $\Delta_{G\backslash X}$ 
est localement ferm\'ee dans $G\backslash X\times G\backslash X$.
\end{itemize}

\begin{marema2}
{\rm Si l'action $\alpha$ est localement r\'eguli\`ere, alors toutes les $G$--orbites dans $X$ sont ferm\'ees, i.e. $G\backslash X$ est un espace {\it de Fr\'echet}. 
En effet, 
si $Y$ est une $G$--orbite dans $X$, et si $x$ appartient \`a $\overline{Y}\smallsetminus Y$ o\`u $\overline{Y}$ 
d\'esigne la fermeture de $Y$ dans $X$, alors $\alpha$ ne peut pas \^etre r\'eguli\`ere au voisinage de $x$.\hfill $\blacksquare$}
\end{marema2}

D'apr\`es \cite[prop.~6.8]{BZ}, les conditions suivantes sont \'equivalentes:
\begin{itemize}
\item l'action $\alpha$ est {\it constructible} au sens o\`u son graphe $\ES{G}_X^\alpha$ 
est constructible dans $X\times X$;
\item la diagonale $\Delta_{G\backslash X}$ 
est constructible dans $G\backslash X\times G\backslash X$.
\end{itemize}
Si l'action $\alpha$ est constructible, on a (loc.~cit.):
\begin{itemize}
\item si $X$ est non vide, alors il existe un sous--ensemble ouvert non vide et $G$--invariant $U$ de $X$ 
sur lequel l'action $\alpha$ est r\'eguli\`ere;
\item toutes les $G$--orbites dans $X$ sont localement ferm\'ees.
\end{itemize}
%

\subsection{Rappels sur la topologie dfinie par $F$}\label{rappels sur la topo dfinie par F} 
Soit $F$ un corps commutatif localement compact non archimdien\footnote{Le cas d'un corps localement compact non archimdien 
est celui qui nous intresse ici. Signalons cependant que la plupart des noncs ci--dessous s'tendent au cas 
d'un corps valu henslien $F$ tel que le compl\'et\'e $\widehat{F}$ de $F$ est une extension sparable de $F$ (cf. 
\cite{MB1, MB2}).}, et soit ${\bf X}$ une vari\'et\'e 
alg\'ebrique (a priori ni affine ni lisse) d\'efinie sur $F$, identifi\'ee comme dans le ch.~3 \`a 
l'ensemble de ses points $\overline{F}$--rationnels pour une cl\^oture alg\'ebrique 
$\overline{F}$ de $F$. Comme en \ref{gnralits} on note $F^{\rm sep}\!/F$ la sous--extension s\'eparable maximale de 
$F$ dans $\overline{F}$, et $\Sigma=\Sigma(F^{\rm sep}\!/F)$ son groupe de Galois. 

La vari\'et\'e ${\bf X}$ \'etant d\'efinie sur $F$, elle est munie d'une {\it $F$--structure}\footnote{Rappelons que si 
${\bf X}$ est affine, la $F$--structure sur ${\bf X}$ est donn\'ee par une sous--$F$--alg\`ebre 
$F[{\bf X}]$ de l'alg\`ebre affine $\overline{F}[{\bf X}]$ de ${\bf X}$ telle que $\overline{F}\otimes_FF[{\bf X}]=\overline{F}[{\bf X}]$. 
En ce cas, les parties $F$--ferm\'ees de ${\bf X}$ correspondent aux id\'eaux $I$ de $F[{\bf X}]$ tels que la $F$--alg\`ebre $F[{\bf X}]/I$ est 
r\'eduite. Soit ${\bf Z}$ est une partie $F$--ferm\'ee de ${\bf X}$. Notons $F[{\bf Z}]$ 
la restriction de $F[{\bf X}]$ \`a ${\bf Z}$ --- c'est une $F$--alg\`ebre r\'eduite --- et $F({\bf Z})$ son anneau des fractions. 
Ce dernier est un produit (fini) d'extensions du corps $F$, et on a 
le crit\`ere suivant \cite[ch. AG, 12.1]{Bor}: la sous--vari\'et\'e ${\bf Z}$ de ${\bf X}$ est d\'efinie sur $F$ si et seulement si la $\overline{F}$--alg\`ebre 
$\overline{F}\otimes_FF({\bf Z})$ est r\'eduite --- ou, ce qui revient au m\^eme d'apr\`es \cite[ch. AG, 2.2]{Bor}, si et seulement $F({\bf Z})$ 
est un produit d'extensions s\'eparables de $F$. Ce crit\`ere s'\'etend naturellement au cas o\`u la vari\'et\'e ${\bf X}$ 
n'est pas affine \cite[ch. AG, 12.2]{Bor}.} et en particulier d'une {\it $F$--topologie} moins fine que la topologie 
de Zariski, cf. \cite[ ch. AG, 11.3]{Bor}. Rappelons qu'une partie ferm\'ee de ${\bf X}$ est {\it $F$--ferm\'ee} si et seulement si 
elle est d\'efinie sur $F^{p^{-\infty}}$. En particulier, la $F$--topologie et la $F^{p^{-\infty}}$--topologie sur ${\bf X}$ co\"{\i}ncident. 
On munit l'ensemble $X={\bf X}(F)$ des points 
$F$--rationnels de ${\bf X}$ de la topologie d\'efinie par $F$:
\begin{itemize}
\item pour tout $F$--ouvert ${\bf U}$ de ${\bf X}$,  
l'ensemble ${\bf U}(F)={\bf U}\cap X$ des points $F$--rationnels de ${\bf U}$ est ouvert dans $X$;
\item si ${\bf X}$ est affine, la topologie sur $X$ est celle d\'efinie en \ref{gnralits}.
\end{itemize}
Cela fait de $X$ un td--espace.

\v1
Pour d\'eterminer si une partie de ${\bf X}$ 
est $F$--ouverte (resp. $F$--ferm\'ee), on dispose du crit\`ere galoisien habituel. 
L'action de $\Sigma$ sur $F^{\rm sep}$ s'\'etend de mani\`ere unique en une action sur 
$\overline{F}$, et on d\'efinit comme dans \cite[ch.~AG, 14.3]{Bor} 
une action de $\Sigma$ sur ${\bf X}={\bf X}(\overline{F})$ --- qui d'ailleurs prolonge celle d\'efinie dans loc. cit. 
sur ${\bf X}(F^{\rm sep})$. 
D'apr\`es \cite[lemma A.4]{BZ}, une partie ouverte (resp. ferm\'ee) de ${\bf X}$ est $F$--ouverte (resp. $F$--ferm\'ee) si et seulement si 
elle est $\Sigma$--stable. On en d\'eduit (loc.~cit.) que pour toute partie $\ES{X}$ de ${\bf X}(F)$, 
la fermeture $\overline{\ES{X}}$ de $\ES{X}$ dans ${\bf X}$ pour la topologie de Zariski est une sous--vari\'et\'e ferm\'ee 
de ${\bf X}$ d\'efinie sur $F$.

\v1
Soit ${\bf X}$ et ${\bf Y}$ deux vari\'et\'es alg\'ebriques d\'efinies sur $F$. Le produit ${\bf X}\times {\bf Y}$ est 
une vari\'et\'e 
alg\'ebrique d\'efinie sur $F$ \cite[ch.~AG, 12.4]{Bor}, et $({\bf X}\times {\bf Y})(F)$ est naturellement hom\'eomorphe 
\`a $X\times Y$, o\`u l'on a pos\'e $X={\bf X}(F)$ et $Y= {\bf Y}(F)$. D'autre part, 
tout $F$--morphisme --- i.e. morphisme d\'efini sur $F$ --- de vari\'et\'es alg\'ebriques $\alpha :{\bf Y}\rightarrow {\bf X}$ 
induit, par passage aux points 
$F$--rationnels, une application continue $\alpha_F:Y\rightarrow X$. 

\begin{exemple}{\rm 
Soit ${\bf X}$ une varit algbrique dfinie sur $F$, et ${\bf H}$ un groupe algbrique affine dfini sur $F$. On suppose que 
${\bf X}$ est muni d'une action algbrique ( gauche) de ${\bf H}$ elle aussi dfinie sur $F$, \cad d'un $F$--morphisme
$$\alpha:{\bf H}\times {\bf X}\rightarrow {\bf X},\, (h,x)\mapsto \alpha(h,x)=h\cdot x$$
Posons $X={\bf X}(F)$ et $H={\bf H}(F)$. 
Rappelons 
que le groupe $H$ (muni de la topologie d\'efinie par $F$) est localement profini. 
Par passage aux points $F$--rationnels, on obtient une action continue
$$
\alpha_F: ({\bf H}\times {\bf X})(F)= H\times X\rightarrow X.
$$
dont le graphe $\ES{G}_X^{\alpha_F}$ est naturellement 
hom\'eomorphe \`a l'image de $({\bf H}\times {\bf X})(F)$ par l'application $\gamma_F$, o\`u 
$\gamma$ d\'esigne le $F$--morphisme ${\bf H}\times {\bf X}\rightarrow {\bf X}\times {\bf X},\,(h,x)\mapsto (x,h\cdot x)$. 
\hfill $\blacksquare$
}
\end{exemple}

\subsection{Le thorme de constructibilit}\label{le thorme de constructibilit}
Le rsultat suivant \cite[theorem A.2]{BZ} est la version $\varpi$--adique du thorme de constructibilit de Chevalley 
\cite[ch.~AG, cor.~10.2]{Bor}.

\begin{montheo}
Soit $\alpha:{\bf Y}\rightarrow {\bf X}$ un $F$--morphisme 
de vari\'et\'es alg\'ebriques d\'efinies sur $F$. L'image $\alpha_F({\bf Y}(F))$ est constructible dans ${\bf X}(F)$.
\end{montheo}

\begin{mesrems}
{\rm 
\begin{enumerate}
\item[(1)]
D'apr\`es \cite[A.5]{BZ}, le th\'eor\`eme est impliqu\'e par le r\'esultat plus faible 
suivant: {\it il existe une sous--vari\'et\'e $F$--ouverte non vide ${\bf U}$ de ${\bf Y}$ telle que l'image $\alpha_F({\bf U}(F))$ est constructible 
dans ${\bf X}(F)$}. En effet, posons $\ES{Y}={\bf Y}(F)\smallsetminus {\bf U}(F)$. 
D'apr\`es \ref{rappels sur la topo dfinie par F}, la fermeture ${\bf Z}=\overline{\ES{Y}}$ de $\ES{Y}$ dans ${\bf Y}$ pour la topologie de Zariski 
est une sous--vari\'et\'e ferm\'ee de ${\bf Y}$ d\'efinie sur $F$, v\'erifiant ${\bf Z}(F)=\ES{Y}$. 
Puisque ${\bf Y}$ est un espace topologique n\oe therien, on peut par induction supposer 
que $\alpha_F({\bf Z}(F))$ est constructible dans ${\bf X}(F)$. Comme ${\bf Y}(F)\smallsetminus 
{\bf Z}(F)={\bf U}(F)$ et que par hypoth\`ese $\alpha_F({\bf U}(F))$ est constructible dans ${\bf X}(F)$, on obtient que 
$\alpha_F({\bf Y}(F))$ est constructible dans ${\bf X}(F)$.
\item[(2)]Soit $\alpha:{\bf Y}\rightarrow {\bf X}$ un morphisme dominant de vari\'et\'es alg\'ebriques 
affines irr\'eductibles. Le comorphisme $\alpha^\sharp:\overline{F}[{\bf X}]\rightarrow \overline{F}[{\bf Y}]$ 
induit par passage aux quotients un morphisme injectif de corps $\overline{F}({\bf X})\hookrightarrow \overline{F}({\bf Y})$, 
qui fait de $\overline{F}({\bf Y})$ une extension de $\overline{F}({\bf X})$, 
et $\alpha$ est s\'eparable si et seulement si cette extension est s\'eparable. Pour 
$x\in \alpha({\bf Y})$, la fibre $\alpha^{-1}(x)$ est une sous--vari\'et\'e ferm\'ee de ${\bf Y}$ --- m\^eme si ${\bf Y}$, ${\bf X}$ et $\alpha$ 
sont d\'efinis sur $F$ et si $x$ appartient \`a ${\bf X}(F)$, cette sous--varit n'est en g\'en\'eral pas d\'efinie sur $F$ --- et si ${\bf Z}$ est une composante 
irr\'eductible de $\alpha^{-1}(x)$, on a l'in\'egalit\'e \cite[ch.~AG, theo.~10.1]{Bor} 
$$\dim{\bf Z}\geq \dim{\bf Y}-\dim{\bf X}.$$ De plus (loc.~cit.)  
il existe un ouvert non vide ${\bf U}$ de ${\bf X}$ contenu dans 
$\alpha({\bf Y})$ tel que si $x$ appartient \`a ${\bf U}$, alors pour toute composante irr\'eductible ${\bf Z}$ de 
$\alpha^{-1}(x)$, l'ingalit ci--dessus est une galit.
\item[(3)]La d\'emonstration du th\'eor\`eme donn\'ee dans \cite[Appendix]{BZ} consiste \`a remplacer le $F$--morphisme 
$\alpha:{\bf Y}\rightarrow {\bf X}$ par un morphisme s\'eparable, apr\`es s'\^etre ramen\'e --- gr\^ace aux points (1) et (2) --- 
au cas 
particulier $(*)$ suivant:
\begin{itemize}
\item les vari\'et\'es alg\'ebriques ${\bf Y}$ et ${\bf X}$ sont affines, irr\'eductibles et lisses;
\item $\overline{\alpha({\bf Y})}={\bf X}$, i.e. le morphisme $\alpha$ est dominant;
\item toutes les composantes irr\'eductibles de toutes les fibres $\alpha^{-1}(\alpha(y))$ ont 
m\^eme dimension.\hfill $\blacksquare$
\end{itemize}
\end{enumerate}
}\end{mesrems}

\begin{exemple}
{\rm Continuons avec l'exemple de \ref{rappels sur la topo dfinie par F}. D'apr\`es le th\'eor\`eme, le graphe $\ES{G}_{\alpha_F}^X$ de $\alpha_F$ 
est constructible dans $X\times X$, i.e. l'action $\alpha_F$ est constructible. 
En particulier, toutes les $H$--orbites dans $X$ sont localement ferm\'ees. \hfill $\blacksquare$}
\end{exemple}

\begin{moncoro}
Soit $\alpha:{\bf H}\rightarrow {\bf G}$ un $F$--morphisme de groupes 
alg\'ebriques d\'efinis sur $F$. L'image $\alpha_F({\bf H}(F))$ est un sous--groupe localement ferm de ${\bf G}(F)$.
\end{moncoro}

\begin{proof}
Le $F$--morphisme $\alpha$ munit ${\bf G}$ d'une action alg\'ebrique de ${\bf H}$ 
d\'efinie sur $F$, et l'image $\alpha_F({\bf H}(F))$ n'est autre que la ${\bf H}(F)$--orbite de l'\'el\'ement neutre de ${\bf G}(F)$ 
pour cette action. D'o le corollaire.
\end{proof}

\subsection{Quelques cas particuliers utiles}\label{cas particuliers}
Le cas particulier suivant  \cite[lemma A.3]{BZ} du thorme de \ref{le thorme de constructibilit} est aussi un outil  essentiel 
 sa dmonstration:

\begin{monlem1}
Soit $\alpha:{\bf Y}\rightarrow {\bf X}$ un $F$--morphisme de vari\'et\'es 
alg\'ebriques affines lisses d\'efinies sur $F$, tel que pour tout $y\in {\bf Y}$ la diff\'erentielle ${\rm d}(\alpha)_y:
{\rm T}({\bf Y})_y\rightarrow {\rm T}({\bf X})_{\alpha(y)}$ 
est surjective (i.e. tel que $\alpha$ est lisse de dimension relative $\dim({\bf Y})-\dim({\bf X})$). 
L'application $\alpha_F:{\bf Y}(F)\rightarrow {\bf X}(F)$ est ouverte.
\end{monlem1}

On a aussi (cf. \cite[remark 1.3.1]{MB2}):

\begin{monlem2}
Soit $\alpha: {\bf Y}\rightarrow {\bf X}$ un $F$--morphisme {\rm propre} de varits algbriques dfinies sur $F$. Alors l'application 
$\alpha_F:{\bf Y}(F)\rightarrow {\bf X}(F)$ est propre (donc ferme).
\end{monlem2}

\begin{mesrems}{\rm Soit $\alpha: {\bf Y}\rightarrow {\bf X}$ un $F$--morphisme de varits algbriques dfinies sur $F$. On a aussi les variantes suivantes des 
lemmes 1 et 2. Elles sont valables pour tout corps valu henslien $F$ tel que le complt $\widehat{F}$ de $F$ est une 
extension sparable de $F$ \cite{MB1, MB2}:
\begin{enumerate}
\item[(1)]Si le morphisme $\alpha$ est {\it tale}, alors l'application $\alpha_F$ est un homomorphisme local.
\item[(2)]Si le morphisme $\alpha$ est {\it lisse}, alors l'application $\alpha_F$ est ouverte.
\item[(3)]Si le morphisme $\alpha$ est {\it fini}, alors l'application $\alpha_F$ est ferm\'ee.
\item[(4)]Si le morphisme $\alpha$ est {\it propre}, alors l'image $\alpha_F({\bf Y}(F))$ 
est ferme dans ${\bf X}(F)$.
\end{enumerate}
Le point (1) est une simple version du {\it thorme des fonctions implicites}. Le point (2) est une gnralisation du lemme 1 (les varits ${\bf X}$ et ${\bf Y}$ n'tant plus 
supposes ni affines ni lisses) et rsulte du fait que $\alpha_F$ admet des {\it sections locales} en chaque point de ${\bf Y}(F)$. Le point (3) est une consquence de la proprit de {\it continuit des racines} 
d'une expression polynomiale. Le point (4) est une consquence du {\it principe de Hasse infinitsimal} (cf. \cite[cor.~1.2.2]{MB2}). 
Notons que si le morphisme $\alpha$ est propre mais n'est pas fini, et bien s\^ur si le corps $F$ n'est pas localement compact, 
alors l'application $\alpha_F$ n'est en gnral pas ferme. \hfill $\blacksquare$
}\end{mesrems}

\begin{exemples}
{\rm Reprenons l'exemple de \ref{rappels sur la topo dfinie par F}, et supposons de plus que la varit ${\bf X}$ est un ${\bf H}$--espace homogne, 
\cad que le $F$--morphisme $\gamma:{\bf H}\times {\bf X}\rightarrow {\bf X}\times {\bf X}$ est surjectif.
\begin{enumerate}
\item[(1)]Puisque la $\overline{F}$--algbre $\overline{F}[t,t^{-1}]$ est un $\overline{F}[t^2,t^{-2}]$--module de type fini, 
le morphisme $\alpha:\Bbb{G}_{\rm m}\rightarrow \Bbb{G}_m,\, t\mapsto t^2$ est fini (et purement insparable si $p=2$). 
En particulier, $(F^\times)^2= \alpha_F(F^\times)$ est un sous--groupe ferm de $F^\times$ --- ce que l'on savait dj! Ainsi 
dans l'exemple (3) de \ref{la topo p-adique}, toutes les ${\bf H}(F)$--orbites dans $\ES{O}_{\bf H}(\delta)(F)$ sont fermes. 
\item[(2)]L'exemple suivant est donn dans \cite[6.1]{GGMB}. Supposons $p>1$, et prenons pour ${\bf H}$ le groupe algbrique affine $\Bbb{G}_{\rm a}\times \Bbb{G}_{\rm m}$ 
oprant sur ${\bf X}= \Bbb{A}^1$ via le morphisme
$$
\alpha: {\bf H}\times {\bf X}\rightarrow {\bf X},\, ((a,b),x)\mapsto a^p+ b^px.
$$
Pour $x\in {\bf X}(F)= F$, la ${\bf H}(F)$--orbite de $x$ est ${\bf H}(F)\cdot x=F^p + (F^\times)^px$, par suite toute ${\bf H}(F)$--orbite 
dans $F$ contient $0$ dans sa fermeture. On a donc deux cas possibles: ou bien $x\in F^p$, auquel cas 
${\bf H}(F)\cdot x =F^p$ est ferm dans $F$; ou bien $x\not\in F^p$, auquel cas ${\bf H}(F)\cdot x$ n'est pas ferm dans $F$. \hfill $\blacksquare$
\end{enumerate}
}
\end{exemples}

\subsection{Un crit\`ere local de sparabilit (rappels)}\label{un critre local} 
Rappelons qu'un morphisme dominant de varits algbriques affines irrductibles 
$\alpha:{\bf Y}\rightarrow {\bf X}$ est dit {\it sparable} si $\overline{F}({\bf Y})$ est une extension sparable de 
$\overline{F}({\bf X})$. Plus gnralement, un morphisme de varits affines $\alpha:{\bf Y}\rightarrow {\bf X}$ tel que ${\bf Y}$ 
est irrductible, est dit sparable si $\overline{F}({\bf Y})$ est une extension sparable de 
$\overline{F}(\overline{\alpha({\bf Y})})$; 
o $\overline{\alpha({\bf Y})}$ dsigne la fermeture de Zariski de l'image $\alpha({\bf Y})$ dans ${\bf X}$. 

Comme pour les morphismes orbites $\pi_v:{\bf H}\rightarrow {\bf H}\cdot v$ (cf. \ref{groupes algbriques}), on a 
un crit\`ere local de s\'eparabilit\'e pour tout morphisme de vari\'et\'es alg\'ebriques $\alpha: {\bf Y}\rightarrow {\bf X}$. 
Pour $x\in {\bf X}$, on note:
\begin{itemize}
\item $\frak{o}_x=\frak{o}_{{\bf X},x}$ l'anneau local de ${\bf X}$ en $x$; 
\item $\frak{p}_x=\frak{p}_{{\bf X},x}$ l'ideal maximal de $\frak{o}_x$; 
\item $\kappa(x) =\kappa_{\bf X}(x)$ le corps r\'esiduel 
$\frak{o}_x/\frak{p}_x$ (il est isomorphe \`a $\overline{F}$).
\end{itemize}
Pour $y\in {\bf Y}$ et $x=\alpha(y)$, on note:
\begin{itemize}
\item $\alpha^\sharp_y: \frak{o}_x\rightarrow \frak{o}_y$ le 
comorphisme de $\alpha$ en $y$;
\item ${\bf Y}_x = {\bf Y}\times_{\bf X}{\rm Spec}\,\kappa(x)$ la fibre g\'eom\'etrique de $\alpha$ en $x$ --- un $\kappa(x)$--schma;
\item $\frak{o}_{\alpha,y}= \frak{o}_y/\alpha^\sharp_y(\frak{p}_x)\frak{o}_y$ l'anneau local de ${\bf Y}_x$ en $y$;
\item $\frak{p}_{\alpha,y}$ l'id\'eal maximal de $\frak{o}_{\alpha,y}$;
\item $\kappa_\alpha(y)$ le corps r\'esiduel $\frak{o}_{\alpha,y}/\frak{p}_{\alpha,y}$.
\end{itemize}
L'espace topologique sous--jacent \`a 
${\bf Y}_x$ est hom\'eomorphe \`a $\alpha^{-1}(x)$, mais le $\kappa(x)$--schma ${\bf Y}_x$ n'est en g\'en\'eral pas r\'eduit. 
On dit que {\it $\alpha$ est s\'eparable en $y$} si l'anneau local $\frak{o}_{\alpha,y}$ 
est {\it r\'egulier} \cite[ch.~AG, 3.9]{Bor}, \cad si la dimension du $\kappa_\alpha(y)$--espace vectoriel $\frak{p}_{\alpha,y}/\frak{p}_{\alpha,y}^2$ 
est \'egale \`a la dimension de Krull de $\frak{o}_{\alpha,y}$, not\'ee $\dim (\frak{o}_{\alpha,y})$. 
Soit
$$
{\rm T}({\bf Y}_x)_y={\rm Hom}_{\kappa_\alpha(x)}(\frak{p}_{\alpha,y}/\frak{p}_{\alpha,y}^2,\kappa_\alpha(x))
$$
l'espace tangent de Zariski de ${\bf Y}_x$ en $y$. On a l'ingalit (loc.~cit.)
$$
\dim (\frak{o}_{\alpha,y})\leq \dim_{\kappa_\alpha(y)}({\rm T}({\bf Y}_x)_y)
$$
avec \'egalit\'e si et seulement si $\alpha$ est s\'eparable en $y$. 

\begin{mesrems}{\rm Soit $\alpha:{\bf Y}\rightarrow {\bf X}$ un morphisme de varits algbriques. 
\begin{enumerate}
\item[(1)]Pour $y\in {\bf Y}$ et $x= \alpha(x)$, 
l'anneau local $\frak{o}_{\alpha^{-1}(x),y}$ 
de la fibre $\alpha^{-1}(x)$ en $y$ est isomorphe au quotient $\frak{o}_{\alpha,y,{\rm red}}$ de $\frak{o}_{\alpha,y}$ par 
l'id\'eal form\'e par les \'el\'ements nilpotents. Les anneaux locaux $\frak{o}_{\alpha,y,{\rm red}}$ sont r\'eduits (par d\'efinition), et 
pour $x\in \alpha({\bf Y})$, les $y\in \alpha^{-1}(x)$ tels que $\frak{o}_{\alpha,y,{\rm red}}$ est r\'egulier forment un 
ouvert dense de $\alpha^{-1}(x)$ \cite[ch.~AG, cor.~17.2]{Bor}. En particulier, la fibre $\alpha^{-1}(x)$ est une 
vari\'et\'e lisse 
si et seulement si {\it tous} les anneaux locaux r\'eduits $\frak{o}_{\alpha,y,{\rm red}}$ ($y\in \alpha^{-1}(x)$) 
sont r\'eguliers. La dimension de Krull de $\frak{o}_{\alpha,y}$ co\"{\i}ncide (par d\'efinition) 
avec la dimension de Krull de $\frak{o}_{\alpha,y,{\rm red}}$, et aussi avec 
$\dim_y(\alpha^{-1}(x))$; o\`u $\dim_y(\alpha^{-1}(x))$ d\'esigne l'inf des $\dim ({\bf U}')$ 
pour ${\bf U'}$ parcourant les ouverts de $\alpha^{-1}(x)$ contenant $y$, c'est--\`a--dire l'inf des $\dim({\bf Z})$ 
pour ${\bf Z}$ parcourant les composantes irr\'eductibles de $\alpha^{-1}(x)$ contenant $y$. 
\item[(2)]
Si de plus ${\bf Y}$ et ${\bf X}$ sont affines et irr\'eductibles, et si $\overline{\alpha({\bf Y})}={\bf X}$, 
alors d'apr\`es \cite[ch.~AG, theo.~10.1]{Bor} on a l'in\'egalit\'e
$$
\dim (\frak{o}_{\alpha,y})\geq \dim {\bf Y}-\dim {\bf X}\quad (y\in {\bf Y}),
$$
et il existe un ouvert non vide ${\bf U}\subset {\bf X}$ tel que pour tout $x\in {\bf U}$ et tout $y\in \alpha^{-1}(x)$, l'ingalit ci--dessus 
est une galit (cf. la reamrque (2) de \ref{le thorme de constructibilit}).\hfill $\blacksquare$
\end{enumerate}}
\end{mesrems}

On a le critre local de sparabilit suivant \cite[lemma A.9]{BZ}:

\begin{monlem1}
Soit $\alpha:{\bf Y}\rightarrow {\bf X}$ un morphisme de varits algbriques, et soit $y$ un point lisse de ${\bf Y}$
tel que $x=\alpha(y)$ est un point lisse de ${\bf X}$. 
Les conditions suivantes sont \'equivalentes:
\begin{itemize}
\item la diff\'erentielle ${\rm d}(\alpha)_y:{\rm T}({\bf Y})_y\rightarrow {\rm T}({\bf X})_x$ de $\alpha$ en $y$ est surjective;
\item $\alpha$ est s\'eparable en $y$ et on a l'\'egalit\'e entre dimensions de Krull
$$
\dim \frak{o}_{\alpha,y}= \dim \frak{o}_{{\bf Y},y} - \dim \frak{o}_{{\bf X},x}.
$$
\end{itemize}
\end{monlem1}

Si $\alpha:{\bf Y}\rightarrow {\bf X}$ v\'erifie les hypoth\`eses du cas $(*)$ de la remarque (3) de \ref{le thorme de constructibilit}, 
alors d'aprs le lemme 1 et la remarque (2) ci--dessus, pour tout $y\in {\bf Y}$, le morphisme $\alpha$ est sparable 
en $y$ si et seulement si la diffrentielle ${\rm d}(\alpha)_y:{\rm T}({\bf Y})_y\rightarrow {\rm T}({\bf X})_x$ de $\alpha$ en $y$ est surjective. 
D'aprs \cite[ch.~AG, theo.~17.3]{Bor}, on en dduit le

\begin{monlem2}
Soit $\alpha:{\bf Y}\rightarrow {\bf X}$ un morphisme de varits algbriques vrifiant les hypothses du cas $(*)$ de la remarque 
(3) de \ref{le thorme de constructibilit}. Les conditions suivantes sont \'equivalentes:
\begin{itemize}
\item le morphisme $\alpha$ est sparable, i.e. $\overline{F}({\bf Y})$ est une extension sparable de $\overline{F}({\bf X})$;
\item il existe un point $y\in {\bf Y}$ tel que $\alpha$ est s\'eparable en $y$;
\item il existe un ouvert non vide ${\bf U}\subset {\bf Y}$ tel que 
$\alpha$ est s\'eparable en tout point de ${\bf U}$.
\end{itemize}
\end{monlem2}

\subsection{Produit fibr\'e (rappels)}\label{produit fibr}
Si $\alpha:{\bf Y}\rightarrow {\bf X}$ et $\beta:{\bf X}'\rightarrow {\bf X}$ sont 
deux morphismes de vari\'et\'es alg\'ebriques affines, le produit fibr\'e $\frak{Y}'={\bf Y}\times_{\bf X}{\bf X}'$ est par d\'efinition 
le $\overline{F}$--sch\'ema affine correspondant \`a la $\overline{F}$--alg\`ebre $\overline{F}[{\bf Y}]\otimes_{\overline{F}[{\bf X}]}\overline{F}[{\bf X}']$, 
muni des morphismes de $\overline{F}$--sch\'emas 
$p_1:\frak{Y}'\rightarrow {\bf Y}$ et $p_2:\frak{Y}'\rightarrow {\bf X}'$ correspondant aux injections 
naturelles $\overline{F}[{\bf Y}]\hookrightarrow \overline{F}[\frak{Y}']$ et $\overline{F}[{\bf X}']\hookrightarrow 
\overline{F}[\frak{Y}']$. Son espace topologique sous--jacent est
$$
\{(y,x')\in {\bf Y}\times {\bf X}':\alpha(y)=\beta(x')\}.
$$
Notons que $\frak{Y}'$ est un sous--sch\'ema ferm\'e de ${\bf Y}\times {\bf X}'$, 
mais n'est en g\'en\'eral pas une vari\'et\'e car son alg\`ebre affine 
$\overline{F}[\frak{Y}']=\overline{F}[{\bf Y}]\otimes_{\overline{F}[{\bf X}]}\overline{F}[{\bf X}']$ peut avoir des 
\'el\'ements nilpotents. Pour $y'=(y,x')\in \frak{Y}'$, les anneaux locaux $\frak{o}_{p_2,y'}$ et $\frak{o}_{\alpha,y}$ 
sont naturellement isomorphes. Si $\alpha$ et $\beta$ sont des $F$--morphismes de vari\'et\'es 
alg\'ebriques d\'efinies sur $F$, alors $\frak{Y}'$ est un $F$-sch\'ema et les projections $p_1$ et $p_2$ sont 
des morphismes de $F$--sch\'emas.

\`A l'inclusion $\frak{Y}'\subset {\bf Y}\times {\bf X}'$ correspond un sous--$\overline{F}$--sch\'ema 
ferm\'e r\'eduit ${\bf Y}'=\frak{Y}'_{\rm red}$ de ${\bf Y}\times {\bf X}'$ d'alg\`ebre affine le quotient $\overline{F}[\frak{Y}']_{\rm red}$ 
de $\overline{F}[\frak{Y}']$ par l'id\'eal form\'e par les \'el\'ements nilpotents, qui lui 
\og est\fg{} une vari\'et\'e alg\'ebrique \cite[ch.~AG, 6.2 et 6.3]{Bor}. Le morphisme $\pi:{\bf Y}'\rightarrow \frak{Y}'$ correspondant 
\`a la projection canonique $\overline{F}[\frak{Y}']\rightarrow \overline{F}[\frak{Y}']_{\rm red}$ est un hom\'eomorphisme 
sur les espaces topologiques sous--jacents. 
Notons $\alpha':{\bf Y}'\rightarrow {\bf X}'$ le morphisme de vari\'et\'es alg\'ebriques $p_2\circ \pi$. 
Pour $y'=(y,x')\in {\bf Y}'$, le morphisme $p_1\circ \pi: {\bf Y}'\rightarrow {\bf Y}$ induit 
un hom\'eomorphisme sur les fibres
$$\alpha'^{-1}(x')\buildrel \sim\over{\longrightarrow} \alpha^{-1}(\beta(x')).
$$ 

Supposons que $\alpha$ et $\beta$ sont des $F$--morphismes de vari\'et\'es alg\'ebriques affines 
d\'efinies sur $F$. Alors $\frak{Y}'$ et ${\bf Y}'$ sont 
des $F$--sch\'emas et $\pi:{\bf Y}'\rightarrow \frak{Y}'$ est un morphisme de $F$--sch\'emas. De plus le $F$--sch\'ema ${\bf Y}'$, identifi\'e \`a une partie $F$--ferm\'ee 
de ${\bf Y}\times {\bf X}'$, \og est\fg{} une vari\'et\'e alg\'ebrique d\'efinie sur $F^{p^{-\infty}}$, et $p_1\circ \pi:{\bf Y}'\rightarrow {\bf Y}$ et 
$p_2\circ \pi:{\bf Y}'\rightarrow {\bf X}'$ sont des $F^{p^{-\infty}}$--morphismes de vari\'et\'es alg\'ebriques d\'efinies sur $F^{p^{-\infty}}$. 

\subsection{Restriction \`a la Weil et morphisme de Frobenius (rappels)}\label{restriction et frobenius}
Soit ${\bf X}$ une vari\'et\'e alg\'ebrique affine 
d\'efinie sur une extension finie $L$ de $F$. On note ${\bf X}'={\rm Res}_{L/F}({\bf X})$ la vari\'et\'e alg\'ebrique affine 
d\'efinie sur $F$, obtenue par restriction des scalaires de $L$ \`a $F$. Concr\^etement, \'ecrivons
$$
{\bf X}={\rm Spec}\; \overline{F}[x_1,\ldots ,x_n]/(f_1,\ldots ,f_m), \quad f_k\in L[x_1,\ldots ,x_n],
$$
et choisissons une base $e_1,\ldots ,e_d$ 
de $L$ sur $F$. Pour $k=1,\ldots ,m$, on pose
$$
{\bf X}'= {\rm Spec}\;\overline{F}[(x'_{i,j})_{i=1,\ldots ,n; j=1,\ldots ,d}]/(f'_{k,j}: k=1,\ldots ,m; j=1,\ldots ,d)
$$
o\`u les $f'_{k,j}\in F[(x'_{i,j})]$ sont donn\'es par
$$
(f'_{k,1}e_1+\cdots +f'_{k,d}e_d)((x'_{i,j}))= f_k(x_1,\ldots ,x_n),\quad x_i = x'_{i,1}e_1+ \cdots + x'_{i,d}e_d.
$$
Le morphisme $\pi_{L/F,{\bf X}}:{\bf X}'\rightarrow {\bf X}$ qui \`a $(x'_{i,j})$ associe $(x_1,\ldots ,x_n)$ comme 
ci-dessus, est d\'efini sur $L$, et il induit un hom\'eomorphisme ${\bf X}'(F)\rightarrow {\bf X}(L)$, o\`u l'on 
munit ${\bf X}'(F)$ de la topologie d\'efinie par $F$ et ${\bf X}(L)$ de la topologie d\'efinie par $L$. Si 
${\bf X}$ est irr\'eductible (resp. lisse), alors ${\bf X}'$ est irr\'eductible (resp. lisse). D'autre part le foncteur 
${\bf X}\mapsto {\bf X}'$ commute au produit: si ${\bf X}_1$, ${\bf X}_2$ sont deux vari\'et\'es 
alg\'ebriques affines d\'efinies sur $L$, la vari\'et\'e ${\rm Res}_{L/F}({\bf X}_1\times {\bf X}_2)$ est $F$--isomorphe 
\`a ${\bf X}'_1\times {\bf X}'_2$, o\`u l'on a pos\'e ${\bf X}'_i= {\rm Res}_{L/F}({\bf X}_i)$, $i=1,\,2$. 
On en d\'eduit en particulier que si 
${\bf X}$ est un groupe alg\'ebrique, alors ${\bf X}'$ est un groupe alg\'ebrique, 
et $\pi_{L/K,{\bf X}}$ est un morphisme de groupes alg\'ebriques.

Supposons $p>1$, et soit $q=p^s$ pour un entier $s\geq 1$. Posons 
$L= F^{1/q}$. C'est une extension finie de $F$, purement ins\'eparable de degr\'e $q$. 
Si ${\bf X}$ est une vari\'et\'e alg\'ebrique affine d\'efinie sur $F$, on note ${^q{\bf X}}$ la vari\'et\'e 
alg\'ebrique affine d\'efinie sur $L$, obtenue en appliquant le morphisme $x\mapsto x^q$ sur les 
coordonn\'ees. Pr\'ecis\'ement, on \'ecrit
$$
{\bf X}={\rm Spec}\; \overline{F}[x_1,\ldots ,x_n]/(f_1,\ldots ,f_m), \quad f_k\in F[x_1,\ldots ,x_n]
$$
et l'on pose
$$
{^q{\bf X}}= {\rm Spec}\; \overline{F}[y_1,\ldots ,y_n]/(h_1,\ldots ,h_m),
$$
o\`u les $h_k\in L[x_1,\ldots ,x_n]$ sont donn\'es par
$$f_k(x_1,\ldots ,x_n)=h_k(y_1,\ldots ,y_n)^q,\quad x_i=y_i^q.
$$ 
Le morphisme $\pi_{q,{\bf X}}:{^q{\bf X}}\rightarrow {\bf X}$ qui \`a $(y_1,\ldots ,y_n)$ associe $(y_1^q,\ldots ,y_n^q)$ 
est d\'efini sur $L$, et il induit un hom\'eomorphisme ${^q{\bf X}}(L)\rightarrow {\bf X}(F)$, o\`u l'on 
munit ${^q{\bf X}}(L)$ de la topologie d\'efinie par $L$ et ${\bf X}(F)$ de la topologie d\'efinie par $F$. 
Si ${\bf X}$ est irr\'eductible (resp. lisse), alors ${^q{\bf X}}$ est irr\'eductible (resp. lisse), et 
comme pour le foncteur restriction \`a la Weil, le foncteur ${\bf X}\mapsto {^q{\bf X}}$ commute au produit. 
En particulier si ${\bf X}$ est un groupe alg\'ebrique, alors ${^q{\bf X}}$ est un groupe alg\'ebrique, et 
$\pi_{q,{\bf X}}$ est un morphisme de groupes alg\'ebriques. 

Continuons avec les hypoth\`eses du paragraphe pr\'ec\'edent ($p>1$, $L=F^{1/q}$), et appliquons le foncteur ${\rm Res}_{L/F}$ \`a 
la vari\'et\'e ${^q{\bf X}}$. On obtient une vari\'et\'e alg\'ebrique ${\bf X}_q={\rm Res}_{L/F}({^q{\bf X}})$ d\'efinie sur $F$, munie 
d'un morphisme
$$\beta_{{\bf X},q}= \pi_{L/F,{^q{\bf X}}}\circ \pi_{q,{\bf X}}:{\bf X}_q\rightarrow {\bf X}$$ lui aussi est d\'efini sur $F$. 
En effet, choisissons une base $\eta_1,\ldots ,\eta_q$ de $F$ sur $F^q$, et posons $e_i=\eta_i^{1/q}$. 
Alors $e_1,\ldots ,e_q$ est une base de $L=F^{1/q}$ sur $F$, et comme plus haut on pose
$${\bf X}_q= {\rm Spec}\;\overline{F}[(y'_{i,j})_{i=1,\ldots ,n; j=1,\ldots ,q}]/(h'_{k,j}: k=1,\ldots ,m; j=1,\ldots ,q),$$
o\`u les $h'_{i,j}\in F[(y'_{i,j})]$ sont donn\'es par
$$
(h'_{k,1}e_1+\cdots +h'_{k,q}e_q)((y'_{i,j}))= h_k(y_1,\ldots ,y_n),\quad y_i = y'_{i,1}e_1+ \cdots + y'_{i,q}e_q.
$$
Alors $\beta=\beta_{{\bf X},q}$ est donn\'e par
$$\beta((y'_{i,j}))=(x_1,\ldots ,x_n),\quad x_i= y'^q_{i,1}\eta_1+ \cdots + y'^q_{i,q}\eta_q, 
$$ 
et d'apr\`es ce qui pr\'ec\`ede, $\beta_F:{\bf X}_q(F)\rightarrow {\bf X}(F)$ est un hom\'eomorphisme. 
D'ailleurs pour toute sous--extension $E/F$ de $F^{\rm sep}\!/F$, puisque $E^q= E\otimes_F F^q$, l'application $\beta_E:{\bf X}_q(E)\rightarrow 
{\bf X}(E)$ induite par $\beta$ sur les points $E$--rationnels, est un hom\'eomorphisme; o\`u l'on munit ${\bf X}_q(E)$ et ${\bf X}(E)$ de 
la topologie d\'efinie par $E$.

\subsection{Le lemme cl\'e}\label{le lemme cl} Le lemme suivant est une simple variante de \cite[lemma A.13]{BZ}, implicitement d\'emontr\'ee dans \cite[A.14]{BZ}. 
On reprend ici les arguments de 
loc.~cit.

\begin{monlem}
Soit $\alpha:{\bf Y}\rightarrow {\bf X}$ 
un $F$--morphisme de vari\'et\'es alg\'ebriques affines d\'efinies sur $F$, tel que 
${\bf Y}$ est irr\'eductible. Il existe un $F$--ouvert non vide ${\bf U}$ de ${\bf Y}$, 
une vari\'et\'e alg\'ebrique affine ${\bf X}'$ d\'efinie sur $F$, et un 
$F$--morphisme de vari\'et\'es alg\'ebriques $\beta:{\bf X}'\rightarrow {\bf X}$ tels que, notant 
${\bf U}'=({\bf U}\times_{\bf X}{\bf X}')_{\rm red}$ le sous--$F$--sch\'ema ferm\'e r\'eduit 
de ${\bf U}\times_{\bf X}{\bf X}'$ correspondant \`a l'inclusion de ${\bf U}\times_{\bf X}{\bf X}'$ dans ${\bf U}\times {\bf X}'$ (cf. \ref{produit fibr}) et 
$\eta':{\bf U}'\rightarrow {\bf X}'$ la projection sur le second facteur, on a: 
\begin{enumerate}
\item[(1)] $\beta_F:{\bf X}'(F)\rightarrow {\bf X}(F)$ est un hom\'eomorphisme;
\item[(2)] $\eta'$ est s\'eparable en tout point de ${\bf U}'$.
\end{enumerate}
Plus prcisment: si $\alpha$ est sparable (e.g. si $p=1$) on peut prendre ${\bf X}'={\bf X}$ et $\beta={\rm id}_{\bf X}$, 
et si $\alpha$ n'est pas sparable (auquel cas $p>1$) on peut prendre ${\bf X}'={\bf X}_q$ et 
$\beta =\beta_{{\bf X},q}$ pour un entier $q=p^s$ ($s\geq 1$) suffisamment grand.
\end{monlem}

\begin{proof}
Commen\c{c}ons par le cas le plus simple: supposons que le $F$--morphisme $\alpha:{\bf Y}\rightarrow {\bf X}$ 
est s\'eparable. Alors d'apr\`es \ref{un critre local}, il existe 
un ouvert non vide ${\bf U}$ de ${\bf Y}$ tel que $\alpha$ est s\'eparable en tout point de ${\bf U}$. 
Quitte \`a remplacer 
l'ouvert ${\bf U}$ par l'union de ses translat\'es sous $\Sigma$, on peut 
le supposer $\Sigma$--stable. Alors ${\bf U}$ est $F$--ouvert et le lemme est d\'emontr\'e. 
En particulier si $p=1$, le lemme est d\'emontr\'e. 

Passons au cas g\'en\'eral. Supposons $p>1$. Notons $L$ le corps 
des fractions de l'anneau $\alpha^\sharp(\overline{F}[{\bf X}])\subset \overline{F}[{\bf Y}]$. C'est un sous--corps du corps 
$\overline{F}({\bf Y})$ des fonctions rationnelles sur ${\bf Y}$\footnote{Rappelons qu'on n'impose pas aux morphismes sparables d'\^etre dominants: 
la fermeture de Zariski ${\bf X}_1=\overline{\alpha({\bf Y})}$ de l'image $\alpha({\bf Y})$ dans ${\bf X}$ est une sous--varit ferme irrductible 
de ${\bf X}$ dfinie sur $F$, et $L$ est isomorphe au corps $\overline{F}({\bf X}_1)$ des fonctions rationnelles sur ${\bf X}_1$.}. 
D'apr\`es \cite[lemma A.14]{BZ}, il existe un entier $q=p^s$ 
($s\geq 1$) tel que $K=(\overline{F}({\bf Y})\otimes_L L^{1/ q})_{\rm red}$ est une extension {\it s\'eparable} de $L^{1/q}$; 
o\`u $(\overline{F}({\bf Y})\otimes_L L^{1/ q})_{\rm red}$ d\'esigne le quotient de l'anneau 
$\overline{F}({\bf Y})\otimes_L L^{1/q}$ par l'id\'eal form\'e par les \'el\'ements nilpotents. \`A cette extension 
correspond un morphisme sparable  de vari\'et\'es 
alg\'ebriques affines irr\'eductibles $\alpha'':{\bf Y}''\rightarrow {\bf X}''$ d\'efini comme suit. On pose ${\bf X}''= {^q{\bf X}}$, 
$\pi= \pi_{q,{\bf X}}:{\bf X}''\rightarrow {\bf X}$ et ${\bf Y}'' = ({\bf Y}\times_{\bf X}{\bf X}'')_{\rm red}$, 
et l'on note $\alpha'':{\bf Y}''\rightarrow {\bf X}''$ le morphisme de vari\'et\'es algbriques donn\'e par  
la projection sur le second facteur. Les vari\'et\'es 
${\bf X}''$ et ${\bf Y}''$ ainsi que le morphisme $\alpha''$ sont 
d\'efinis sur $F^{p^{-\infty}}$, et par construction $\alpha''$ est sparable. 
D'apr\`es \ref{un critre local}, il existe un ouvert ${\bf V}''$ de ${\bf Y}''$, que l'on peut supposer 
d\'efini sur $F^{p^{-\infty}}$ (c'est--\`a--dire $F^{p^{-\infty}}$--ouvert), tel que $\alpha''$ est s\'eparable en tout point de ${\bf V}''$. 
Comme $\pi: {\bf X}''\rightarrow {\bf X}$ est un hom\'eomorphisme pour la 
$F^{p^{-\infty}}$--topologie, la projection sur le premier facteur $\delta:{\bf Y}''\rightarrow {\bf Y}$ est aussi un 
hom\'eomorphisme pour la $F^{p^{-\infty}}$--topologie. Par suite ${\bf U}=\delta({\bf V}'')$ est $F^{p^{-\infty}}$--ouvert, 
et donc $F$--ouvert \cite[ch.~AG, 12.1]{Bor}, dans ${\bf Y}$. 
Par construction, le morphisme de vari\'et\'es donn\'e par la projection sur le 
second facteur
$$
\eta'': {\bf U}''=({\bf U}\times_{\bf X}{\bf X}'')_{\rm red}\rightarrow {\bf X}''
$$
est s\'eparable en tout point de ${\bf U}''$. 
Pour que la condition (1) soit v\'erifi\'ee, il suffit de remplacer ${\bf X}''$ par ${\bf X}'={\bf X}_q\;(={\rm Res}_{F^{1/q}/F}({\bf X}''))$, et 
$\pi$ par $\beta= \rho\circ \pi$, o\`u $\rho =\pi_{F^{1/q}/F,{\bf X}''}:{\bf X}'\rightarrow {\bf X}''$. 
D'apr\`es \ref{restriction et frobenius}, l'application $\beta_F:{\bf X}'(F)\rightarrow {\bf X}(F)$ est un hom\'eomorphisme. Notons 
$$
\eta': {\bf U}'=({\bf U}\times_{\bf X}{\bf X}')_{\rm red}\rightarrow {\bf X}'
$$
le morphisme de vari\'et\'es donn\'e par la projection sur le second facteur. 
Puisque le morphisme $\rho$ est s\'eparable (pour la d\'efinition de $\pi_{F^{1/q}/F,{\bf X}''}$, cf. \ref{restriction et frobenius}), 
le $\overline{F}$--sch\'ema ${\bf U}''\times_{{\bf X}''}{\bf X}'$ est 
r\'eduit. Il est donc isomorphe \`a ${\bf U}'$, et 
pour $y'= (u,x')\in {\bf U}\times_{\bf X}{\bf X}'$ et $y''= (u,\rho(x'))\in {\bf U}\times_{\bf X}{\bf X}''$, 
les anneaux locaux $\frak{o}_{\alpha',y'}$ et $\frak{o}_{\alpha'',y''}$ sont isomorphes. Par cons\'equent $\eta'$ 
est s\'eparable en tout point de ${\bf U}'$, et le lemme est d\'emontr\'e.
\end{proof}

\begin{marema1}{\rm 
D'apr\`es la d\'emonstration du lemme, on peut choisir le $F^{p^{-\infty}}$--ouvert 
${\bf V}''$ de ${\bf Y}''$ tel que (en fait il est implicitement choisi ainsi): 
\begin{itemize}
\item ${\bf V}''$ est lisse;
\item $\alpha''({\bf V}'')$ est contenu dans l'ouvert dense de $\overline{\alpha''({\bf Y}'')}$ form des points simples;
\item les anneaux locaux $\frak{o}_{\alpha''\!,y''}$ pour $y''\in {\bf V}''$ ont tous la 
m\^eme dimension.
\end{itemize}
En ce cas ${\bf U}'$ est lisse et 
les anneaux locaux $\frak{o}_{\alpha'\!,y'}$ pour $y'\in {\bf U}''$ ont tous la 
m\^eme dimension. Si de plus on suppose que ${\bf X}$ est lisse et que $\alpha$ est dominant (ce qui 
implique que ${\bf X}$ est irr\'eductible), 
alors le morphisme $\eta':{\bf U}'\rightarrow {\bf X}'$ v\'erifie les conditions du cas $(*)$.\hfill $\blacksquare$
}\end{marema1}

\begin{marema2}{\rm La vari\'et\'e ${\bf U}'$ n'est en g\'en\'eral pas d\'efinie sur $F$, 
mais seulement sur une extension finie purement ins\'eparable $F'$ de $F$. De plus, puisque $\alpha'$ 
est un morphisme de $F$--sch\'emas, c'est un $F'$--morphisme de vari\'et\'es alg\'ebriques d\'efinies sur 
$F'$. 
Comme dans \cite[A.13]{BZ} on peut remplacer ${\bf U}'$ par une vari\'et\'e ${\bf Z}$ d\'efinie sur $F$: 
l'ensemble
$${\bf U}'(F)={\bf U}'\cap ({\bf Y}\times {\bf X}')(F)$$ est \'egal \`a
$$
\{(y,x')\in {\bf U}(F)\times {\bf X}'(F): \alpha_F(u)=\beta_F(x')\},
$$
et puisque $\beta_F:{\bf X}'(F)\rightarrow {\bf X}(F)$ est un hom\'eomorphisme, ${\bf U}'(F)$ est hom\'eomorphe 
\`a ${\bf U}(F)$, et $\beta_F$ induit un hom\'eomorphisme de $\alpha'({\bf U}'(F))$ sur $\alpha_F({\bf U}(F))$. 
Soit ${\bf Z}=\overline{{\bf U}'(F)}$ la fermeture de ${\bf U}'(F)$ dans ${\bf U}\times {\bf X}'$ pour la topologie de Zariski. 
C'est une sous--vari\'et\'e ferm\'ee de ${\bf U}\times {\bf X}'$, contenue dans ${\bf U}'$ et 
d\'efinie sur $F$ (\ref{rappels sur la topo dfinie par F}). On a ${\bf U}'(F)={\bf Z}(F)$, et $\alpha'$ 
induit un $F$--morphisme $\alpha_{\bf Z}:{\bf Z}\rightarrow {\bf X}'$ de vari\'et\'es alg\'ebriques d\'efinies sur $F$, qui 
est s\'eparable en tout point $y'$ de ${\bf Z}$ tel que $\dim \frak{o}_{\alpha_{\bf Z},y'}=\dim \frak{o}_{\alpha'\!,y'}$. 
En effet, pour $y'\in {\bf Z}$, l'anneau local $\frak{o}_{\alpha_{\bf Z},y'}$ est un quotient de 
$\frak{o}_{\alpha'\!,y'}$, et puisque $\frak{o}_{\alpha'\!,y'}$ est r\'egulier, on a 
$$\dim \frak{o}_{\alpha_{\bf Z},y'}\leq \dim \frak{o}_{\alpha'\!,y'}$$ 
avec \'egalit\'e si et seulement si 
$\frak{o}_{\alpha_{\bf Z},y'}=\frak{o}_{\alpha'\!,y'}$. En particulier, si le morphisme $\alpha\vert_{\bf U}$ est {\it quasi--fini}, 
i.e. si pour tout $y\in {\bf U}$ la fibre $\alpha^{-1}(\alpha(y))$ 
est finie, alors le morphisme $\alpha'$ est quasi--fini, et le morphisme $\alpha_{\bf Z}$ (lui aussi 
quasi--fini) est sparable en tout point de ${\bf Z}$.  \hfill $\blacksquare$ }
\end{marema2}

\subsection{Un rsultat bien connu}\label{un rsultat bien connu} 
Le lemme suivant, bien connu des spcialistes, est valable pour n'importe quel corps commutatif algbriquement clos $\overline{F}$. 

Rappelons \cite[ch.~AG, theo.~10.1]{Bor} que si $f: {\bf Y}\rightarrow {\bf X}$ est un morphisme dominant de varits algbriques irrductibles, 
on a $\dim({\bf Y})=\dim({\bf X})$ si et seulement si l'ensemble des $x\in {\bf X}$ tels que la fibre $f^{-1}(x)\subset {\bf Y}$ est finie, contient un 
ouvert dense de ${\bf X}$.

\begin{monlem}Soit $f: {\bf Y}\rightarrow {\bf X}$ un morphisme de varits algbriques affines irrduc\-tibles. On suppose 
que $f$ est dominant, et 
que $\dim({\bf Y})=\dim({\bf X})$. Alors il existe un ouvert affine ${\bf U}\subset {\bf X}$ 
tel que le morphisme 
$$f\vert_{f^{-1}({\bf U})}:f^{-1}({\bf U})\rightarrow {\bf U}$$
est fini. 
\end{monlem} 

\begin{proof}Notons $A=\overline{F}[{\bf X}]$ et $B= \overline{F}[{\bf Y}]$ les algbres affines de ${\bf X}$ et ${\bf Y}$, et soit 
$K=\overline{F}({\bf X})$ et $L= \overline{F}({\bf Y})$ leurs corps des fractions. 
Le comorphisme
$$
f^\sharp: A\rightarrow B
$$
est injectif, et fait de $B$ une $A$--algbre de type fini. Il induit par passage aux corps des fractions 
un morphisme injectif de corps
$K\hookrightarrow L$,
qui fait de $L$ une extension (de type fini) de $K$. Commen\c{c}ons par montrer que cette extension est {\it finie}. 
La $K$--algbre (de type fini) 
$K\otimes_AB$ est l'agbre affine $K[{\bf Y}_\eta]$ de la {\it fibre gnrique} 
${\bf Y}_\eta= {\bf Y}\times_{\bf X} {\rm Spec}(K)$ --- un $K$--schma rduit, mais pas gomtriquement rduit --- 
de $f$; o $\eta$ est le point gnrique de ${\bf X}$. Puisque $K\otimes_AB$ est isomorphe au localis $S^{-1}B$, $S= 
f^\sharp(A)\smallsetminus \{0\}$, c'est un anneau intgre. D'aprs le \og lemme de normalisation de 
N\oe ther \fg{}, il existe des lments $x_1,\ldots ,x_n$ dans $K\otimes_AB$, algbriquement indpendants 
sur $K$, faisant de $K\otimes_AB$ un $K[x_1,\ldots ,x_n]$--module de type fini. 
Comme $n=\dim({\bf Y}_\eta)=\dim({\bf Y})-\dim({\bf X})$, on a $n=0$. Par consquent 
$K\otimes_AB$ est un $K$--espace vectoriel de dimension finie. C'est 
donc un corps (puisque c'est un anneau intgre), et $L= K\otimes_AB$ est une extension finie 
de $K$.

Choisissons un sous--ensemble fini $\{b_1,\dots ,b_m\}\subset B$ engendrant $B$ sur $A$. D'aprs le paragraphe prcdent, chaque 
$b_i$ est algbrique sur $K$, i.e. $Q_i(b_i)=0$ pour un polyn\^ome non nul $Q_i(t)\in K[t]$. Puisque 
$K$ est le corps des fractions de $A$, il existe un lment $a\in A$, $a\neq 0$, tel que 
$aQ_i(t)\in A[t]$ pour $i=1,\ldots ,m$. Posant $b= f^\sharp(a)\in B\smallsetminus\{0\}$, l'anneau $B[b^{-1}]$ est entier 
sur $A[a^{-1}]$. Mais comme $B[b^{-1}]$ est une $A[a^{-1}]$--algbre de type fini, c'est un $A[a^{-1}]$--module de type fini, et le morphisme
$$
f\vert_{{\rm spec}(B[b^{-1}]}:{\rm Spec}(B[b^{-1}])\rightarrow {\rm Spec}(A[a^{-1}])
$$
est fini.
\end{proof}

\subsection{Une cons\'equence du lemme cl\'e}\label{une consquence du lemme cl}
Soit ${\bf H}$ un groupe algbrique, et $\alpha:{\bf Y}\rightarrow {\bf X}$ est un morphisme de vari\'et\'es 
alg\'ebriques. On dit que $\alpha$ est un {\it ${\bf H}$--morphisme} si les varits ${\bf Y}$ et ${\bf X}$ sont munies d'une action 
algbrique ( gauche) de ${\bf H}$, et si l'on a
$$
\alpha(h\cdot y)= h\cdot \alpha(y)\quad (h\in {\bf H},\, y\in {\bf Y}). 
$$
Plus gnralement, si $\beta:{\bf G}\rightarrow {\bf H}$ est un morphisme de groupes algbriques, on dit que $\alpha$ 
est un {\it $\beta$--morphisme} si la varit ${\bf Y}$ est muni d'une action algbrique de ${\bf G}$, la varit ${\bf X}$ est muni d'une action 
algbrique de ${\bf H}$, et si l'on a
$$
\alpha(g\cdot y)= \beta(g)\cdot \alpha(y)\quad (g\in {\bf G},\,y\in {\bf Y}).
$$ 

La proposition suivante est une consquence da la dmonstration du lemme de \ref{le lemme cl}. 

\begin{mapropo}
Soit ${\bf H}$ un groupe algbrique affine dfini sur $F$, et soit 
$\alpha:{\bf Y}\rightarrow {\bf X}$ un $F$--morphisme {\rm fini} de varits algbriques affines 
dfinies sur $F$, tel que ${\bf Y}$ est irrductible. On suppose que les varits ${\bf Y}$ et ${\bf X}$ 
sont munies d'une action algbrique de ${\bf H}$ dfinie sur $F$, et que 
$\alpha$ est un ${\bf H}$--morphisme. On suppose aussi que ${\bf Y}$ est un ${\bf H}$--espace homogne. 
Si le morphisme $\alpha$ est sparable, 
alors il est tale. Sinon, 
il existe un $F$--morphisme $\alpha_1: {\bf Y}_1\rightarrow {\bf X}_1$ de varits algbriques affines dfinies sur $F$, et des 
$F$--morphismes de varits algbriques $\gamma: {\bf Y}_1\rightarrow {\bf Y}$ et $\zeta: {\bf X}_1\rightarrow {\bf X}$, tels que:
\begin{enumerate}
\item[(1)]le morphisme $\alpha_1$ est fini et tale;
\item[(2)]l'application $\gamma_F:{\bf Y}_1(F)\rightarrow {\bf Y}(F)$ est un homomorphisme;
\item[(3)]l'application $\zeta_F:{\bf X}_1(F)\rightarrow {\bf X}(F)$ induit par restriction un homomorphisme
$$\alpha_{1,F}({\bf Y}_1(F))\rightarrow \alpha_F({\bf Y}(F));$$
\item[(4)]on a l'galit $\alpha_F= \zeta_F\circ \alpha_{1,F}\circ \gamma_F^{-1}$.
\end{enumerate}
\end{mapropo}

\begin{proof} 
Puisque le morphisme $\alpha$ est fini, il est surjectif. Par suite la varit ${\bf X}$ est irrductible, et c'est un ${\bf H}$--espace homogne. 
En particulier (par homognit), le morphisme $\alpha$ vrifie les hypothses du cas $(*)$ de la remarque (3) 
de \ref{le thorme de constructibilit}. Notons aussi que puisque $\alpha$ est fini, on a $\dim({\bf Y})=\dim({\bf X})$. 

Soit $L\subset \overline{F}({\bf Y})$ le corps des fractions de l'anneau 
$\alpha^\sharp(\overline{F}[{\bf X}])$.

Si le morphisme $\alpha$ est sparable, i.e. si l'extension $\overline{F}[Y]/L$ est sparable, alors 
d'aprs le lemme 2 de \ref{un critre local} (par homognit), le morphisme $\alpha$ est sparable en tout point de ${\bf Y}$. 
Par suite (lemme 1 de \ref{un critre local}), pour tout 
$y\in {\bf Y}$, la diffrentielle ${\rm d}(\alpha)_y:{\rm T}({\bf Y})_y\rightarrow {\rm T}({\bf X})_{\alpha(y)}$ est surjective, donc bijective. 
En d'autres termes, le morphisme $\alpha$ est lisse de dimension relative $0$. Il est donc tale. 

Supposons maintenant que l'extension $\overline{F}({\bf Y)}/L$ n'est pas sparable (on a donc $p>1$), et 
choisisssons comme dans la dmonstration du lemme de \ref{le lemme cl} un entier $q=p^s$, $s\geq 1$, tel que 
$K= (\overline{F}({\bf Y})\otimes_L L^{1/q})_{\rm red}$ est une extension sparable de $L^{1/q}$. 
Reprenons la d\'emonstration de loc.~cit., en tenant compte de l'action de ${\bf H}$. 

Posons ${\bf H}''= {^q{\bf H}}$ et notons $\pi_{\bf H}$ le $F^{p^{-\infty}}$--morphisme $\pi_{q,{\bf H}}: {\bf H}''\rightarrow {\bf H}$ (cf. \ref{restriction et frobenius}). C'est un morphisme de groupes algébriques, et par fonctorialit\'e, la vari\'et\'e ${\bf X}''={^q{\bf X}}$ est munie d'une action alg\'ebrique de 
${\bf H}''$ d\'efinie sur $F^{p^{-\infty}}$, qui fait de $\pi=\pi_{q,{\bf X}}:{\bf X}''\rightarrow {\bf X}$ un $\pi_{\bf H}$--morphisme. 
On en d\'eduit une action alg\'ebrique de ${\bf H}''$ sur la vari\'et\'e ${\bf Y}\times {\bf X}''$, d\'efinie sur $F^{p^{-\infty}}$: pour 
$h''\in {\bf H}''$ et $(y,x'')\in {\bf Y}\times {\bf X}''$, on pose
$$
h''\cdot (y,x'')= (\pi_{\bf H}(h)\cdot y, h''\cdot x'').
$$ 
Cette action se restreint en une action algbrique sur ${\bf Y}''=({\bf Y}\times_{\bf X}{\bf X}'')_{\rm red}$, elle aussi d\'efinie sur $F^{p^{-\infty}}$.
La projection sur le second facteur $\alpha'':{\bf Y}'\rightarrow {\bf X}''$ est un ${\bf H}''$--morphisme surjectif, et la projection 
sur le premier facteur ${\bf Y}''\rightarrow {\bf Y}$ est un $\pi_{\bf H}$--morphisme surjectif. Les morphismes 
$\pi: {\bf X}''\rightarrow {\bf X}$ et $\pi_{\bf H}: {\bf H}''\rightarrow {\bf H}$ sont des homéomorphismes pour la $F^{p^{-\infty}}$--topologie. On en déduit que la variété ${\bf X}''$ est un ${\bf H}''$--espace homogène, et que la variété ${\bf Y}''$ est elle aussi un ${\bf H}''$--espace homogène. Puisque ${\bf H}''$ est lisse, les variétés ${\bf X}''$ et ${\bf Y}''$ le sont aussi --- pour ${\bf X}''={^q{\bf X}}$ on le savait déjà, puisque ${\bf X}$ est lisse ---, et comme le ${\bf H}''$--morphisme $\alpha''$ est séparable, par homogénéité il est séparable en tout point de ${\bf Y}''$ (i.e. on peut, dans la dmonstration du lemme de \ref{le lemme cl}, prendre 
${\bf U}={\bf Y}$ et ${\bf U}''={\bf Y}''$). De plus  $\alpha''$ est un morphisme quasi--fini. Il est donc fini (par homognit, d'aprs le 
lemme de \ref{un rsultat bien connu}), et lisse de dimension relative $0$, donc tale. 

Il faut ensuite (comme dans la démonstration du lemme de \ref{un rsultat bien connu})) remplacer ${\bf X}''$ par ${\bf X}'= {\bf X}_q\;(= {\rm Res}_{F^{1/q}/F}({\bf X}''))$, et $\pi$ par $\beta = \rho \circ \pi$, où $\rho= \pi_{F^{1/q}/F,{\bf X}'}: {\bf X}'\rightarrow {\bf X}''$. La variété ${\bf X}'$ est définie sur $F$, et $\beta$ est un $F$--morphisme qui, par passage aux point $F$--rationnels, donne un homéomorphisme $\beta_F: {\bf X}'(F)\rightarrow {\bf X}(F)$. Posons ${\bf Y}'=({\bf Y}\times_{\bf X}{\bf X}')_{\rm red}$ et notons
$$
\alpha':{\bf Y}'\rightarrow {\bf X}'
$$
la projection sur le second facteur. C'est un $F^{p^{-\infty}}$--morphisme surjectif. Montrons qu'il est fini et étale. 
Le morphisme $\rho$ est séparable, par conséquent le $\overline{F}$--schéma ${\bf Y}''\times_{{\bf X}''}{\bf X}'$ est réduit, et il est isomorphe à ${\bf Y}'$. Puisque le morphisme $\alpha'': {\bf Y}''\rightarrow {\bf X}''$ est fini, le morphisme $\alpha' \simeq \alpha''\times_{{\bf X}''}{\bf X}'$ l'est aussi. Pour $y'=(y,x')\in {\bf Y}\times_{\bf X}{\bf X}'$ 
et $y''=(y,\rho(x'))\in {\bf Y}\times_{\bf X}{\bf X}''$, les anneaux locaux $\frak{o}_{y',\alpha'}$ et $\frak{o}_{\alpha'',y''}$ sont isomorphes. 
Le morphisme $\alpha'$ est donc séparable en tout point de ${\bf Y}'$, et même lisse de dimension relative $0$, donc étale\footnote{En particulier, si $F'/F$ est une sous--extension finie 
de $F^{p^{-\infty}}/F$ telle que ${\bf Y}'$ est dfini sur $F'$, alors le morphisme $\alpha'$ est dfini sur $F'$, et l'application 
$\alpha'_{F'}: {\bf Y}'(F')\rightarrow {\bf X}'(F')$ est ferme, et c'est un homomorphisme local.}. 

Notons ${\bf Y}_1$ la fermeture de ${\bf Y}'(F)={\bf Y}'\cap ({\bf Y}\times {\bf X}')(F)$ 
dans ${\bf Y}\times {\bf X}'$ pour la topologie de Zariski. C'est une sous--vari\'et\'e ferm\'ee de ${\bf Y}'$ dfinie sur $F$ 
telle que ${\bf Y}_1(F)= {\bf Y}'(F)$, et 
$\alpha'$ induit par restriction un $F$--morphisme $\alpha'_1:{\bf Y}_1\rightarrow {\bf X}'$ qui, d'aprs la 
remarque 2 de \ref{le lemme cl}, est sparable en tout point de ${\bf Y}_1$. Puisque $\alpha'_1$ est le compos d'une immersion ferme 
${\bf Y}_1\hookrightarrow {\bf Y}'$ et du morphisme fini $\alpha':{\bf Y}'\rightarrow {\bf X}'$, c'est lui--m\^eme un morphisme fini. L'image 
${\bf X}_1=\alpha'_1({\bf Y}_1)$ est donc une sous--varit ferme de ${\bf X}'$ dfinie sur $F$, et le $F$--morphisme $\alpha_1:{\bf Y}_1\rightarrow {\bf X}_1$ 
dduit  de $\alpha'_1$ est surjectif, fini, et sparable en tout point de ${\bf Y}_1$. Il est donc lisse, et m\^eme tale.

En dfinitive on a le diagramme commutatif suivant 
$$
\xymatrix{
{\bf Y}'(F) = {\bf Y}_1(F) \ar[d] \ar[r]^-{\alpha_{1,F}} & {\bf X}_1(F) \ar[r]& {\bf X}'(F)\ar[d]^{\beta_F}\\
{\bf Y}(F)\ar[rr]_{\alpha_F}&& {\bf X}(F)
}
$$
o la fl\`eche verticale de gauche est celle dduite par restriction de la projection sur le premier facteur 
${\bf Y}(F)\times {\bf X}'(F)\rightarrow {\bf Y}(F)$, et la fl\`eche horizontale du haut  droite est l'immersion ferme canonique 
(inclusion). On note $\gamma:{\bf Y}_1\rightarrow {\bf Y}$ le $F$--morphisme dduit de la projection sur le premier facteur, et 
$\zeta: {\bf X}_1\rightarrow {\bf X}$ le $F$--morphisme dduit de $\beta$. Puisque
$$
{\bf Y}'(F)=\{(y,x')\in {\bf Y}(F)\times {\bf X}'(F): \alpha_F(y)=\beta_F(x')\}
$$
et que $\beta_F:{\bf X}'(F)\rightarrow {\bf X}(F)$ est un homomorphisme, on obtient que:
\begin{itemize}
\item l'application $\gamma_F:{\bf Y}_1(F)\rightarrow {\bf Y}(F)$ 
est un homomorphisme;
\item l'application $\zeta_F:{\bf X}_1(F)\rightarrow {\bf X}(F)$ induit par restriction un homomorphisme
$$\alpha_{1,F}({\bf Y}_1(F))\rightarrow \alpha_F({\bf Y}(F)).$$
\end{itemize}
Cela achve la dmonstration de la proposition. 
\end{proof}

\begin{moncoro}
L'application (ferme) $\alpha_F: {\bf Y}(F)\rightarrow {\bf X}(F)$ est un 
homomorphisme local sur son image $\alpha_F({\bf Y}(F))$. 
\end{moncoro}

\begin{proof}
Les applications $\alpha_F:{\bf Y}(F)\rightarrow {\bf X}(F)$ et $\alpha_{1,F}:{\bf Y}_1(F)\rightarrow {\bf X}_1(F)$ sont 
fermes (lemme 2 de \ref{cas particuliers}), et 
l'application $\alpha_{1,F}:{\bf Y}_1(F)\rightarrow {\bf X}_1(F)$ 
est un homomorphisme local (remarque (3) de \ref{cas particuliers}). D'o le corollaire.
\end{proof}

\begin{mesrems}
{\rm 
\begin{enumerate}
\item[(1)]
Si l'extension $\overline{F}[Y]/L$ est sparable, puisque le morphisme $\alpha$ est tale 
(et fini), l'application $\alpha_F:{\bf Y}(F)\rightarrow {\bf X}(F)$ 
est un homomorphisme local (et elle est ferme). En particulier l'image 
$\alpha_F({\bf Y}(F))$ est une sous--varit $\varpi$--adique ouverte et ferme de ${\bf X}(F)$.

\item[(2)]
Si l'extension $\overline{F}[Y]/L$ n'est pas sparable, comme l'application 
$$\alpha_F= \beta_F\circ \alpha_{1,F}\circ \gamma_F^{-1}: {\bf Y}(F)\rightarrow {\bf X}(F)$$
est un homomorphisme local sur son image $\alpha_F({\bf Y}(F))$, elle la munit d'une structure de varit 
$\varpi$--adique, mais cette dernire {\it n'est pas} une sous--varit $\varpi$--adique de ${\bf X}(F)$.

\item[(3)]
Soit $L'/L$ la sous--extension sparable maximale de $\overline{F}({\bf Y})/L$. Si $L'\neq \overline{F}({\bf Y})$, 
alors on peut prendre pour $q$ le degr de l'extension (purement insparable) $\overline{F}({\bf Y})/L'$. 
En effet on a l'galit $L'^{1/q}=\overline{F}({\bf Y})$, d'o l'inclusion 
$L^{1/q}\subset \overline{F}({\bf Y})$. Comme $L^{1/q}/L$ est une extension purement insparable 
de degr $q$, on a l'galit $\overline{F}({\bf Y})=L^{1/q}L'$; o $L^{1/q}L'$ est l'extension compose de $L^{1/q}/L$ et $L'/L$ dans $\overline{F}/L$. 
Comme d'autre part $K= (\overline{F}({\bf Y})\otimes_L L^{1/q})_{\rm red}$ s'identifie  $\overline{F}({\bf Y})L^{1/q}=\overline{F}({\bf Y})$, 
l'extension $K/L^{1/q}$ est sparable (de degr $[L':L]$).  \hfill$\blacksquare$
\end{enumerate}
}
\end{mesrems}


\end{document}